%% file: main.tex
\documentclass[10pt,a4paper]{article}
\usepackage[utf8]{inputenc}

\usepackage{
    mathrsfs,
    xcolor,
    verbatim,
    amsmath,
    amsthm,
    amsfonts,
    amssymb,
    dsfont,
    nicefrac,
    enumitem,
    bm,
    mathtools,
    xparse,
    etoolbox
}
\usepackage[normalem]{ulem}

\numberwithin{equation}{section}

\usepackage{thmtools}
\usepackage[margin=1in]{geometry}

\usepackage{tikz}
\usetikzlibrary{
    matrix,
    chains,
    positioning,
    decorations.pathreplacing,
    arrows,
    shapes,
    math,
    arrows.meta,
    calc,
    shapes.geometric
}

\usepackage[most]{tcolorbox}
\tcbset{
    theorembox/.style={
        colback=white!95!gray,
        colframe=black,
        boxrule=0.5pt,
        left=5pt,
        right=5pt,
        top=3pt,
        bottom=3pt,
        before skip=1\baselineskip,
        after skip=1\baselineskip,
        sharp corners,
        enhanced,
        breakable=false
    }
}
\usepackage{capt-of}
\usepackage{xurl}

\newcommand{\Power}[1]{\mathcal{P}(#1)}
\newcommand{\Ran}{\mathbb{R}_{\mathrm{an}}}

\newif\ifrevisionnotes
\revisionnotestrue
\NewDocumentCommand{\note}{+m}{%
    \ifrevisionnotes
        \begingroup
        \normalfont\color{red!70!black}%
        \textbf{[Note:} #1\textbf{]}%
        \endgroup
    \fi
}

\NewDocumentCommand{\replace}{O{} m m}{%
    \ifrevisionnotes
        \textcolor{red}{\sout{#2}}\,\textcolor{green!50!black}{#3}%
        \ifstrempty{#1}{}{\space\note{#1}}%
    \else
        #3%
    \fi
}

\newcommand{\R}{\mathbb R}
\newcommand{\dd}{\,\mathrm d}
\newcommand{\ddd}{\mathrm d}
\newcommand{\Id}{\mathrm{Id}}
\newcommand{\grad}{\nabla}

\newcommand{\bfX}{\mathbf{X}}
\newcommand{\bfY}{\mathbf{Y}}
\newcommand{\indicator}[1]{\mathds{1}_{\smash{#1}}}

\newcommand{\rel}[2]{\mathcal{N}^{#1}_{#2}}
\newcommand{\rell}[3]{\mathcal{N}^{#1}_{#2,#3}}

\newcommand{\N}{\mathbb{N}}
\DeclarePairedDelimiter{\abs}{\lvert}{\rvert}
\DeclarePairedDelimiter{\braces}{\lbrace}{\rbrace}
\DeclarePairedDelimiter{\brackets}{[}{]}
\DeclarePairedDelimiter{\norm}{\lvert}{\rvert}
\DeclarePairedDelimiter{\opnorm}{\lVert}{\rVert}

\NewDocumentCommand{\signedint}{}{\mathop{\mathrlap{\bm{\int}}}\!\int}
\NewDocumentCommand{\TVmeasure}{m}{%
  \mathopen{\pmb{\pmb{\boldsymbol{\lvert}}}}#1\mathclose{\pmb{\pmb{\boldsymbol{\rvert}}}}\cfadd{def:signed-radon-measure}%
}
\newcommand{\dsubmanifold}[3]{$#3$-dimensional $C^{#1}$-submanifold of $\R^{#2}$}

\newcommand{\EMPH}[1]{#1}

\usepackage{graphicx} %

\newcommand{\midcup}{\mathop{\mathchoice
  {\vcenter{\hbox{\scalebox{1.1}{$\textstyle\bigcup$}}}}%
  {\bigcup}%
  {\bigcup}%
  {\bigcup}%
}\nolimits}
\newcommand{\midcap}{\mathop{\mathchoice
  {\vcenter{\hbox{\scalebox{1.1}{$\textstyle\bigcap$}}}}%
  {\bigcap}%
  {\bigcap}%
  {\bigcap}%
}\nolimits}
\usepackage[sort,nocompress]{cite}

\usepackage{hyperref}
\usepackage[sort,capitalise,nameinlink]{cleveref}

\theoremstyle{plain}
\newtheorem{theorem}{Theorem}[section]

\newtheorem{definition}[theorem]{Definition}
\theoremstyle{definition}

\theoremstyle{remark}

\hypersetup{
    colorlinks=true,
    linkcolor=blue!80!black,
    citecolor=green,
    urlcolor=blue!80!black
}

\crefname{enumi}{item}{items}
\crefname{equation}{}{}
\crefname{subsection}{Subsection}{Subsections}
\crefname{figure}{Figure}{Figures}
\crefname{setting}{Setting}{Settings}

\makeatletter
\AtBeginEnvironment{tcolorbox}{%
   \def\@mpfn{footnote}%
}
\makeatother

\input{commands.tex}

\makeatother
\ExplSyntaxOn
\bool_gset_false:N \g_kuckuck_debug_colors_bool
\bool_gset_false:N \g_kuckuck_debug_logs_bool
\ExplSyntaxOff
\abbr{GD}{gradient descent}{GDs}{gradient descents}{GD}
\abbr{DE}{differential equation}{DEs}{differential equations}{DE}
\abbr{KL}{Kurdyka-\L ojasiewicz}{KLs}{Kurdyka-\L ojasiewiczs}{KL}
\abbr{DTGC}{definable Thom gradient conjecture}{DTGCs}{definable Thom gradient conjectures}{DTGC}
\abbr{GF}{gradient flow}{GFs}{gradient flows}{GF}
\abbr{ANN}{artificial neural network}{ANNs}{artificial neural networks}{ANN}
\abbr{DNN}{deep neural network}{DNNs}{deep neural networks}{DNN}
\abbr{RMSprop}{root mean square propagation}{RMSprops}{the root mean square propagation}{RMSprop}
\abbr{SGD}{stochastic gradient descent}{SGDs}{stochastic gradient descent}{SGD}
\abbr{iid}{independent and identically distributed}{i.i.d}{independent and identically distributed}{i.i.d.}
\abbr{SOP}{stochastic optimization problem}{SOPs}{stochastic optimization problems}{SOP}
\abbr{AL}{deep learning}{ALs}{deep learning}{AL}
\abbr{DL}{deep learning}{DLs}{deep learning}{DL}
\abbr{Adam}{adaptive moment estimation \SGD}{Adams}{adaptive moment estimation \SGD}{Adam}
\abbr{Muon}{momentum orthogonalized by Newton-Schulz \SGD}{Muons}{momentum orthogonalized by Newton-Schulz \SGD}{Muon}
\abbr{NS}{Newton-Schulz}{NSs}{Newton-Schulz}{NS}
\abbr{AdamW}{\Adam\ with decoupled weight decay}{AdamWs}{\Adam\ with decoupled weight decay}{AdamW}
\abbr{AI}{artificial intelligence}{AIs}{artificial intelligence}{AI}
\abbr{PDE}{partial differential equation}{PDEs}{partial differential equations}{PDE}
\abbr{as}{almost surely}{ass}{almost surely}{a.s.}
\abbr{pas}{$\P$-almost surely}{pass}{almost surely}{$\P$-a.s.}
\makeatother

\title{On the gradient conjecture at infinity:\\
Finite spherical length for gradient\\
trajectories of definable potentials}
\author{Thang Do$^{1}$, Arnulf Jentzen$^{2,3}$, and Benno Kuckuck$^{4}$
	\bigskip
	\\
    	\small{$^1$ School of Data Science, The Chinese University of Hong Kong, Shenzhen}
	\vspace{-0.1cm}\\
	\small{ (CUHK-Shenzhen), China, e-mail: \texttt{minhthangdo@link.cuhk.edu.cn}}
 \smallskip
	\\
	\small{$^2$ School of Data Science, The Chinese University of Hong Kong, Shenzhen}
	\vspace{-0.1cm}\\
	\small{(CUHK-Shenzhen), China, e-mail: \texttt{ajentzen@cuhk.edu.cn}}
	\smallskip
	\\
 \small{$^3$ Applied Mathematics: Institute for Analysis and Numerics, Faculty of Mathematics and}
	\vspace{-0.1cm}\\
	\small{Computer Science, University of M{\"u}nster, Germany, e-mail: \texttt{ajentzen@uni-muenster.de}}
	\smallskip
	\\
    \small{$^4$ Applied Mathematics: Institute for Analysis and Numerics, Faculty of Mathematics and}
	\vspace{-0.1cm}\\
	\small{Computer Science, University of M{\"u}nster, Germany, e-mail: \texttt{bkuckuck@uni-muenster.de}}
	\smallskip
	\\
}

\date{\today}
\begin{document}
\maketitle

\begin{abstract}
We prove that for every open $ U \subseteq \R^n $ the spherical projection
of every gradient trajectory $ \theta(t) $, $ t \geq 0 $,
of a \emph{definable} $ C^1 $-potential $ F \colon U \to \R $
has finite length and, in particular, converges (directional convergence).
In other words, every gradient trajectory of a definable $ C^1 $-potential has
\emph{finite total angular motion}.
Here, \emph{definable} means definable in an arbitrary o-minimal structure expanding the real field.

Notably, the gradient trajectory in this result is an arbitrary solution of the gradient flow differential equation (whose solutions may not be unique) with no further assumptions:
The gradient trajectory
may be unbounded and diverge to infinity 
or it may be convergent.
The established theorem therefore provides affirmative answers to both
the \emph{gradient conjecture at infinity} for definable potentials
and the \emph{definable Thom gradient conjecture}.

We extend our conclusions from the Euclidean metric to a general class of ambient metrics
including Riemannian metrics and we establish uniform local upper bounds for the lengths
of the spherical projections of gradient flow trajectories.
We finally apply our findings to the training of deep neural networks (DNNs)
to conclude that, in the case of definable and continuously differentiable activation and loss functions,
the spherical projection of the DNN training flow must always converge (directional convergence).
\end{abstract}

\pagebreak

\setcounter{tocdepth}{2}
\tableofcontents
\clearpage

\section{Introduction and main result}
\label{sec:introduction}

The study of solutions of \GF\ \DEs\
has a long history in mathematics
(see, \eg, \cite{Loj84_MR771152,MR713046,Thom89,K98,KMP00,MR2274510} and \cref{sec:literature_review} below).

In the main result of this work, \cref{maintheorem} below, we prove
for every natural number
$ n \in \N = \{ 1, 2, 3, \dots \} $
and every open $ U \subseteq \R^n $
that the spherical projection of every \GF\ trajectory
$ \theta( t ) $, $ t \geq 0 $, of a $ C^1 $-potential $ F \colon U \to \R $ definable in
an arbitrary o-minimal structure over the reals
has finite length and must, in particular, be convergent (\emph{directional convergence}).

In the following, definability of $ F $ means that the graph of the function $ F $ is a definable set
in an arbitrary fixed o-minimal structure over the reals (see \cref{sec:basics} below).
A \GF\ trajectory is an arbitrary solution of the \GF\ \DE.
It does not need to be bounded or have a convergent subsequence but may \emph{diverge to infinity}.

Our findings have strong connections to
the \emph{gradient conjecture}~\cite{MR713046,Thom89,KMP00},
the \emph{definable Thom gradient conjecture}~\cite{Sanz2024_MR4789467,KP06_MR2325137},
the \emph{gradient conjecture at infinity} \cite{DHL22_MR4453693,G07_MR2290270},
and also to \DL\ \cite{Cai25,JT20,KGDJ25}
and we refer to \cref{sec:literature_review} below
for a short literature review.

\subsection{Finite spherical length of gradient flow (GF) trajectories for the Euclidean metric}

We now present \cref{maintheorem} in a self-contained way and, thereafter,
we explain the conclusion and consequences of \cref{maintheorem} and also relate
it to the state-of-the-art in the current literature.

\begin{tcolorbox}[theorembox]
\begin{athm}{theorem}{maintheorem}
Let $ n \in \mathbb{N} $,
$ \vartheta \in \R^n $,
let $ U \subseteq \R^n $ be open,
let $ F \in C^1( U, \R ) $ be \EMPH{definable},
let $ I \subseteq \R $ be convex with nonempty interior,
let
$ \theta \in C^1( I, U ) $ satisfy for all $ t \in I $ that
\begin{equation}
\label{eq:first_theorem_ODE}
\EMPH{
  \theta'(t) = - (\nabla F)( \theta( t ) )
}
  ,
\end{equation}
and let
$ P \colon U \to \R^n $
satisfy\footnotemark\ for all $ v \in U \backslash \{ \vartheta \} $
that
$
  P( v ) = | v - \vartheta |^{ - 1 } ( v - \vartheta )
$.
Then \EMPH{$ P \circ \theta $} has finite length.
\end{athm}
\end{tcolorbox}
\footnotetext{Note that for all $ d \in \N $, $ v = ( v_1, \dots, v_d ) \in \R^d $
it holds that $ | v | = ( \sum_{ i = 1 }^d | v_i |^2 )^{ 1 / 2 } $ (standard norm).}

\cref{maintheorem} is proved as \cref{maintheorem3} in \cref{subsec:general-time-intervals} below.
The conclusion of \cref{maintheorem} can be viewed as follows.
The point $ \vartheta $ in \cref{maintheorem} serves as the
reference point of polar coordinates,
with every point $ v \in \R^n \backslash \{ \vartheta \} $
having a unique representation
by the polar coordinate pair consisting
of the distance $ | v - \vartheta | \in (0,\infty) $ to the reference point of the polar coordinates
and its angular component\footnote{Note that the case
where the \GF\ trajectory hits the reference point $ \vartheta $ (cf., \eg, \cref{lem:fibre_is_convex})
is covered in \cref{maintheorem} by allowing the spherical projection $ P $ at $ \vartheta $ to be an arbitrary element of $ \R^n $.},
i.e., its spherical projection
$
  | v - \vartheta |^{ - 1 } ( v - \vartheta )
$.
\cref{maintheorem} then proves for every definable $ C^1 $-potential $ F \colon U \to \R $
that every \GF\ trajectory $ \theta \colon I \to U $ of $ F $
has \emph{finite total angular motion} over its entire time interval $ I $, that is,
the length of $ P \circ \theta $ is finite.
Clearly, the fact that $ P \circ \theta $ has finite length implies that there exists a unique
$ v \in \R^n $ such that
$ \limsup_{ t \to \sup( I ) } | ( P \circ \theta )( t ) - v | = 0 $
(\emph{directional convergence}); see, \eg, \cref{lem:finite-length-endpoint-limit}.

Note that the time interval $ I $ in \cref{maintheorem} may be bounded
or unbounded, \eg, $ I = [0,\infty) $
and the endpoints of $ I $ may or may not be elements of $ I $.
\cref{maintheorem} covers the situation where
$ \liminf_{ t \to \sup( I ) } | \theta( t ) | = \infty $ (divergence to infinity)
as well as the situation where
$ \liminf_{ t \to \sup( I ) } | \theta( t ) | < \infty $
(cf., \eg, \cref{lem:locally-definable-endpoint} below).

Concerning smoothness, \cref{maintheorem} only assumes that the potential $ F \colon U \to \R $
is continuously differentiable (an element of $ C^1( U, \R ) $) and
thus only imposes that the gradient $ \nabla F $ is continuous
but not, e.g., that it is locally Lipschitz continuous.
Solutions of \cref{eq:first_theorem_ODE} may thus not be unique. For example,
in the case of the time interval $ I = [0,\infty) $
and the potential $ F = ( \R \ni v \mapsto - \frac{ 2 }{ 3 } | v |^{ 3 / 2 } \in \R ) $
we have for every $ c \in I $ that
$ I \ni t \mapsto - ( \max\{ 0, (t - c) / 2 \} )^2 \in \R $
is a solution of \cref{eq:first_theorem_ODE}.

\cref{maintheorem} assumes that the potential is definable
with respect to some \emph{arbitrary} given o-minimal structure.
In \cref{definition:definable-function} in \cref{sec:basics} below, we briefly recall the notion of a definable function.
Note that all polynomials and all semialgebraic functions (see, \eg, \cite[Section~2.2]{MR1659509})
are definable functions in every o-minimal structure but not every analytic function is a definable function
(\eg, the function
$ \R \ni x \mapsto \sin( x ) \in \R $ is not definable in any o-minimal structure).
The conclusion of \cref{maintheorem} even fails to hold in general if the definability assumption is replaced by analyticity (see \cref{cor:analytic-graph-counterexample} and \cref{prop:analytic-counterexample} below
and cf.\ also \cite[Theorem~1.1, page 592]{DHL22_MR4453693}).

However, the regularity assumption does allow the following weakening: In the case where the \GF\ trajectory $ \theta(t) $, $ t \in I $, in \cref{maintheorem}
has accumulation points in $ U $ at both ends of the time interval $ I $,
the conclusion of \cref{maintheorem} also holds when the potential is only locally definable.
This statement is the subject of the following \cref{maincor},
which is a direct consequence of \cref{maintheorem} (see \cref{cor:spherical}
in \cref{sec:locally-definable-corollary} below).
In \cref{definition:locally-definable-function} in \cref{sec:basics} below
we recall the precise definition of a locally definable function.

\begin{tcolorbox}[theorembox]
\begin{athm}{corollary}{maincor}
Let $ n \in \mathbb{N} $,
$ \vartheta \in \R^n $,
let $ U \subseteq \R^n $ be open,
let $ F \in C^1( U, \R ) $ be \EMPH{locally definable},
let $ I \subseteq \R $ be convex with nonempty interior,
let
$ \theta \in C^1( I, U ) $ satisfy for all $ t \in I $ that
\begin{equation}
\label{eq:first_theorem_ODE_BC}
\EMPH{
  \theta'(t) = - (\nabla F)( \theta( t ) )
}
  ,
\end{equation}
let $ a, b \in U $ satisfy
\EMPH{$
  \liminf_{ t \to \inf( I ) } | \theta(t) - a |
  =
  \liminf_{ t \to \sup( I ) } | \theta(t) - b |
  = 0
$},
and let
$ P \colon U \to \R^n $
satisfy for all $ v \in U \backslash \{ \vartheta \} $
that
$
  P( v ) = | v - \vartheta |^{ - 1 } ( v - \vartheta )
$.
Then \EMPH{$ P \circ \theta $} has finite length.
\end{athm}
\end{tcolorbox}

\cref{maincor} is proved as
\cref{cor:spherical} in \cref{sec:locally-definable-corollary} below.
\cref{maincor} assumes that the potential is a \emph{locally definable} function.
Analytic, locally bounded semianalytic, and
locally bounded subanalytic functions on open sets
(cf., \eg, \cite[Definitions~2.1, 3.1, and 3.2]{MR972342})
are all examples
of locally $ \Ran $-definable functions (see \cref{definition:locally-definable-function} below).
Here, $ \Ran $ is the structure generated by all restricted real analytic 
functions,
whose definable sets are precisely the globally subanalytic sets 
(see, \eg, \cite[Section~3.1, page 137]{Wilkie2009_MR2605320}).
O-minimality of this structure was established in \cite{Gabrielov1968,vanDenDries1986}.
Moreover, every definable function is also locally definable.

\subsection{Related work and key novelty of this article}
\label{sec:literature_review}

In the case where the potential is not only locally definable (as in \cref{maincor}) but analytic,
the \emph{gradient conjecture} (also known as the \emph{gradient conjecture of R.\ Thom})
asked whether the spherical projection of a \GF\ trajectory converging
to the reference point is convergent;
see, \eg, \cite[Section~III.B, page 204]{Thom89}, \cite[Problem 18, page 64]{MR713046}, and \cite[Section IV.9, page 1592]{MR1275210}.

Finite length of the spherical projection of the \GF\ trajectories is a stronger conclusion
than mere convergence as asserted by the gradient conjecture:
finite length controls the total angular motion, not only the limit of the angular motion.
In the case of analytic potentials, finite spherical length of \GF\ trajectories that converge
to the reference point in the domain $ U $ (and thereby a complete solution to the classical gradient conjecture of R.\ Thom)
was established in Kurdyka, Mostowski, \& Parusi\'nski \cite[Theorem~7.1]{KMP00}.
Before the complete solution of the classical gradient conjecture in \cite{KMP00},
special cases of it (such as the two-dimensional case $ n = 2 $)
had been established in \cite[Theorem~3, page 60]{Kuiper1991},
\cite[Section 1, page 134]{Ichikawa92_MR1156895},
and \cite[Section~3, page 450]{Moussu1997_MR1437048}.

It is natural to ask whether the conclusion of the gradient conjecture also holds
for (locally) definable potential functions $ F \colon U \to \R $
and this extension of the gradient conjecture has sometimes been referred
to as the \emph{definable Thom gradient conjecture} (cf., \eg, \cite{P05} and \cite[Section 7.8.3.2]{Sanz2024_MR4789467}).

The definable Thom gradient conjecture and finite spherical length for convergent \GF\ trajectories
have been established in the two-dimensional case $ n = 2 $ \cite[Item (2) (page 139) and \S 8]{KP06_MR2325137}
as well as for polynomially bounded o-minimal structures \cite[Item (3) (page 139) and \S 9]{KP06_MR2325137}
(where an o-minimal structure is called
polynomially bounded if every function from $ \R $ to $ \R $ that is definable with respect to 
this o-minimal structure grows at most polynomially).

In analogy to R.\ Thom's gradient conjecture that deals with gradient trajectories that are bounded and converge to a point,
there is also the \emph{gradient conjecture at infinity} which asserts that the spherical
projection $ | \theta(t) |^{ - 1 } \theta(t) $ converges as $ t $ tends to $ T $
for every \GF\ trajectory $ \theta(t) $, $ t \in [0,T) $, that escapes to infinity in the sense that
$ \liminf_{ t \to T } | \theta( t ) | = \infty $ where $ T \in (0,\infty] $; cf., \eg, \cite[Section~6]{G07_MR2290270}.
In particular,
in the case where the potential is $ C^2 $ and semialgebraic,
\cite[Theorem~5.1]{G07_MR2290270} not only solves the gradient conjecture at infinity
but also establishes finite length of the spherical projection
of every \GF\ trajectory that escapes to infinity.

The work \cite[Theorem~1.1, page 592]{DHL22_MR4453693} explicitly specifies for every $ k \in \N $
a convex $ C^k $-potential $ F \colon \R^2 \to \R $ (whose restriction to $ \R^2 \backslash \{ 0 \} $ is analytic)
for which the gradient conjecture at infinity \emph{fails to hold} for the positive gradient flow of $ F $
(or, equivalently, for the negative gradient flow of the concave potential $ - F $),
with escape to infinity occurring at a finite maximal time.
This does not contradict \cref{maintheorem} above
as the potentials in \cite{DHL22_MR4453693} are not definable.
Moreover, in the setting of suitable non-Euclidean Riemannian metrics,
the work \cite[Proposition 5.1, page 76]{G14MR3200186} specifies a smooth semialgebraic
potential for which the gradient conjecture at infinity also fails to hold.

For bounded \GF\ trajectories, \KL-type inequalities near critical points are powerful tools
for establishing finite length and convergence of \GF\ trajectories \cite{K98}.
In a related direction, Bochnak--\L ojasiewicz-type gradient inequalities
at infinity, near asymptotic critical values, are studied in
\cite[Proposition~2.7, page~497]{MR2166522},
\cite[Proposition~3.2, page~42]{MR2208534}, and
\cite[Proposition~6, page~69]{Loi_MR3483169}.
In particular, \cite[Proposition~2.7, page~497]{MR2166522}
in conjunction with \cite[Remark~9, page~69]{Loi_MR3483169}
establishes such an inequality
for polynomially bounded o-minimal structures.

\cref{maintheorem} provides affirmative answers to the definable Thom gradient conjecture and to the definable gradient conjecture at infinity
and it not only establishes \emph{directional convergence} but also
\emph{finite spherical length} for both regimes
(at a critical point and at infinity)
for any arbitrary $ C^1 $-potential that is definable with respect to an arbitrary
o-minimal structure over the reals.

In contrast to the classical Thom gradient conjecture
and the definable Thom gradient conjecture,
the \emph{strong gradient conjecture} due to R.\ Moussu (which is also known as the \emph{non-oscillation conjecture})
(see,
\eg,
\cite[Section~1]{Sanz1998_MR1656007},
\cite[Theorem~5.1, page 99]{FS01_MR1844525},
\cite[Theorem 2.1, item 5, pages 2519--2520]{Goldstein2009_MR2568063},
\cite[Conjecture (F), page 781]{K98},
and
\cite[Theorem~1, page 1685]{GS13_MR3072668})
remains an open problem of research in most cases,
even for analytic potential functions $ F \colon U \to \R $.

For further reading and further references on the topic of the dynamics
of the gradient vector field of an analytic or definable function (potential)
we refer to the survey article \cite{Sanz2024_MR4789467}.

\subsection{Finite spherical length of GF trajectories for ambient metrics}

In \cref{maintheorem} gradients are taken with respect to the standard Euclidean metric.
In the following result, \cref{thm:ambient-conditional-reference} below, we extend \cref{maintheorem}
to the setting of more general ambient metrics. Specifically,
if $ n \in \N $,
if $ U \subseteq \R^n $ is an open set,
if $ F \in C^1( U, \R ) $,
if $ \mathbb{S} = \{ A \in \R^{ n \times n } \colon A \text{ is symmetric positive-definite} \} $,
if
$ G \colon U \to \mathbb{S} $
and
$ B \colon U \to \mathbb{S} $
satisfy for all $ x \in U $ that $ B(x) = ( G(x) )^{ - 1 } $ (ambient metric and inverse ambient metric),
if $ g_x \colon \R^n \times \R^n \to \R $, $ x \in U $,
satisfy\footnote{Note that for all $ d \in \N $, $ v = ( v_1, \dots, v_d ) $, $ w = ( w_1, \dots, w_d ) \in \R^d $
it holds that $ \langle v, w \rangle = \sum_{ i = 1 }^d v_i w_i $ (standard scalar product).}
for all $ x \in U $, $ v, w \in \R^n $ that
$ g_x( v, w ) = \langle v, G(x) w \rangle $ (scalar products associated with $ G $),
then for every $ x \in U $ we have that
$ B(x) \nabla F( x ) $
is the metric gradient in the sense that
for all $ v \in \R^n $ it holds that
$
  g_x( B(x) \nabla F( x ), v ) = F'( x ) v = \langle \nabla F( x ), v \rangle
$.

\begin{tcolorbox}[theorembox]
\begin{theorem}
\label{thm:ambient-conditional-reference}
Let $n\in\N$, $\vartheta_0,\vartheta_1,\vartheta_2\in\R^n$, $c,\varepsilon\in(0,\infty)$,
let $I\subseteq\R$ be an interval with nonempty interior,
let $U\subseteq\R^n$ be open,
let $F\in C^1(U,\R)$,
let $V\subseteq U$ be a set, assume that $F|_V$ is definable,
let $B\colon\R^n\to\R^{n\times n}$ and
$\lambda\colon\R\to[0,\infty)$ be measurable,
assume for all $v\in\R^n$ that $B(v)$ is symmetric positive-definite,
let $\theta\colon I\to U$ be locally absolutely continuous,
assume for Lebesgue-almost all $t\in I$ that
\begin{equation}
  \label{eq:ambient-conditional-flow}
  \EMPH{\theta'(t)=-\lambda(t)B(\theta(t))(\nabla F)(\theta(t))},
\end{equation}
assume\footnotemark\ for all $i\in\{0,1,2\}$, $v\in\R^n$ with
$0<\norm{v-\vartheta_i}
  (\indicator{\theta(I)}(\vartheta_1)+|i-1|)
  <\varepsilon$ that
$v\in V$ and $\opnorm{B(v)-B(\vartheta_i)}\leq c\norm{v-\vartheta_i}^{\varepsilon}$,
assume $\sum_{z\in\{-1,1\}}\liminf_{t\to z\sup(zI)}
  \bigl([1-\indicator{I}(z\sup(zI))]\,
    \norm{\theta(t)-\vartheta_{z+1}}\bigr)=0$,
and let $P\colon U\to\R^n$ satisfy for all
$v\in U\backslash\{\vartheta_1\}$ that
$P(v)=\norm{v-\vartheta_1}^{-1}(v-\vartheta_1)$.
Then \EMPH{$P\circ\theta$} has finite length.
\end{theorem}
\end{tcolorbox}
\footnotetext{Note that for all $m,n\in\N$, $A\in\R^{m\times n}$ it holds that
$\opnorm{A}=\sup_{v\in\R^n,\,|v|=1}|Av|$ (spectral norm).}

\Cref{thm:ambient-conditional-reference} is proved as
\cref{maintheorem:ambient1} in \cref{subsec:ambient-finite-length} below.
According to the condition
$\sum_{z\in\{-1,1\}}\liminf_{t\to z\sup(zI)}
  ([1-\indicator{I}(z\sup(zI))]\,
    \norm{\theta(t)-\vartheta_{z+1}})=0$
in \cref{thm:ambient-conditional-reference}, we assume in
\cref{thm:ambient-conditional-reference} that, if $\inf(I)\notin I$,
then $\vartheta_0$ must be an accumulation point of $\theta(t)$ as
$t\to\inf(I)$ and that, if $\sup(I)\notin I$, then $\vartheta_2$ must be an
accumulation point of $\theta(t)$ as $t\to\sup(I)$.
\Cref{thm:ambient-conditional-reference} applies, \eg, to smooth Riemannian metrics
as well as to ambient metrics satisfying the stated local H\"{o}lder estimates
around $\vartheta_0$, $\vartheta_1$, and $\vartheta_2$.
Analytic potentials in conjunction with analytic Riemannian metrics have been treated in \cite[Section~9]{KMP00}.

It is also possible to establish a uniform upper bound
for the local length of the spherical projections of \GF\ trajectories counting
all visits in a sufficiently small ball around $ \vartheta $.
This is the subject of the following \cref{thm:ambient-uniform-general-interval}.

\begin{tcolorbox}[theorembox]
\begin{theorem}
\label{thm:ambient-uniform-general-interval}
Let $n\in\N$, $\vartheta\in\R^n$, $c,\varepsilon\in(0,\infty)$,
let $V=\{v\in\R^n\colon0<|v-\vartheta|<\varepsilon\}$,
let $\mathbb S=\{A\in\R^{n\times n}\colon A\text{ is symmetric positive-definite}\}$,
let $\mathbb B\in\mathbb S$,
and let $f\in C(V,\R)$ be definable.
Then there exist $\delta,L\in(0,\infty)$ such that
for every interval with nonempty interior $I\subseteq\R$,
every open $U\subseteq\R^n$,
every $F\in C^1(U,\R)$,
every locally absolutely continuous $\theta\colon I\to U$,
and all measurable $P\colon U\to\R^n$, $\lambda\colon\R\to[0,\infty)$, and
$B\colon\R^n\to\mathbb S$
with $F|_{V\cap U}=f|_{V\cap U}$,
with the property that for all $v\in U\backslash\{\vartheta\}$ it holds that
$P(v)=|v-\vartheta|^{-1}(v-\vartheta)$ and
$\opnorm{B(v)-\mathbb B}\indicator{V}(v)\leq c|v-\vartheta|^\varepsilon$,
and with the property that for Lebesgue-almost all $t\in I$ it holds that
\begin{equation}
  \label{eq:ambient-uniform-general-flow}
  \EMPH{\theta'(t)=-\lambda(t)B(\theta(t))(\nabla F)(\theta(t))}
\end{equation}
we have that
\begin{enumerate}[label=(\roman*)]
\item it holds that $P\circ\theta$ is Lebesgue-almost everywhere differentiable and
\item it holds that
\begin{equation}
  \label{eq:ambient-uniform-general-length}
  \EMPH{\int_I |(P\circ\theta)'(s)|\,
    \indicator{[0,\delta]}(|\theta(s)-\vartheta|)\,\mathrm{d}s\leq L}.
\end{equation}
\end{enumerate}
\end{theorem}
\end{tcolorbox}

\Cref{thm:ambient-uniform-general-interval} is proved as
\cref{maintheorem:ambient2} in \cref{subsec:ambient-uniform-length} below.
We emphasize that the constants $ \delta $ and $ L $
in \cref{eq:ambient-uniform-general-length}
in \cref{thm:ambient-uniform-general-interval} are
independent of the considered specific \GF\ trajectory $ \theta \colon I \to U $
and the considered ambient metric $ \R^n \ni v \mapsto [ B( v ) ]^{ - 1 } \in \R^{ n \times n } $.

\subsection{Definability, o-minimal structures, and rectifiable curves}
\label{sec:basics}

For the convenience of the reader, in this subsection, we briefly recall
those concepts from analysis and o-minimal geometry
that are required to accurately retrace the statements of
\cref{maintheorem,maincor}
above.
For further reading in the area of o-minimal geometry we refer, \eg, to
\cite[Chapter~1]{vdD_MR1633348}
and
\cite[Chapter~1]{Coste}.

\begin{athm}{definition}{definition:real-structure}[Structures on the real numbers]
We say that $ \mathcal{S} $ is a structure if and only if
it holds that
$ \mathcal S = ( \mathcal S_n )_{ n \in \N } \colon \N \to ( \midcup_{ n = 1 }^{ \infty } \Power{ \Power{ \R^n } } ) $
is a function from $ \N $ to $ ( \midcup_{ n = 1 }^{ \infty } \Power{ \Power{ \R^n } } ) $
which satisfies that
\begin{enumerate}[label=(\roman*)]

\item it holds for every $ n \in \N $, every polynomial
$ p \colon \R^n \to \R $, and every\footnote{Note that
for every $ I \subseteq \R $ it holds that $ I $ is an interval if and only if $ I $ is convex.}
interval $ I \subseteq \R $ that
$ p^{ - 1 }( I ) \in \mathcal{S}_n $,

\item it holds for all $ n \in \N $, $ A, B \in \mathcal{S}_n $ that
$ \{ A \cap B, A \cup B, A \backslash B \} \subseteq \mathcal{S}_n \subseteq \Power{ \R^n } $,

\item it holds for all $ n, m \in \N $, $ A \in \mathcal{S}_n $,
$ B \in \mathcal{S}_m $ that $ A \times B \in \mathcal{S}_{ n + m } $,
and

\item it holds for all $ n \in \N $, $ A \in \mathcal{S}_{ n + 1 } $ that
$ \midcup_{ x = (x_1, \dots, x_{n+1}) \in A } \{ ( x_1, \dots, x_n ) \} \in \mathcal{S}_n $.
\end{enumerate}
\end{athm}

\begin{athm}{definition}{definition:o-minimal-structure}[O-minimal structures]
Let $\mathcal S $
be a structure\cfadd{definition:real-structure}\dott \cfload[.]
Then we say that $\mathcal S$ is o-minimal if and only if it
holds for every $ A \in \mathcal S_1 $ that
$ A $ is a finite union of intervals\dott \cfout[.]
\end{athm}

\begin{athm}{definition}{definition:definable-set}[Definable sets]
Let $ \mathcal S $
be an o-minimal structure\cfadd{definition:real-structure}\cfadd{definition:o-minimal-structure}\dott \cfload[.]
Then we say that $ A $ is an $ \mathcal{S} $-definable set
(we say that $ A $ is definable) if and only if there exists $ n \in \N $
such that $ A \in \mathcal{S}_n $\dott \cfout[.]
\end{athm}

\begin{athm}{definition}{definition:definable-function}[Definable functions]
Let $ \mathcal S $
be an o-minimal structure\cfadd{definition:real-structure}\cfadd{definition:o-minimal-structure}\dott \cfload[.]
Then we say that $ f $ is an $ \mathcal{S} $-definable function
(we say that $ f $ is definable) if and only if
it holds\footnote{Note that for every $ X $, every $ Y $, and every $ f \colon X \to Y $ it holds that
$ \operatorname{Graph}(f) = \{ (x,y) \in X \times Y \colon y = f(x) \} $.} that
$ \operatorname{Graph}(f) $
is an $ \mathcal{S} $-definable set\dott \cfout[.]
\end{athm}

In our formulation of \cref{maincor} above we also employ the concept of a locally definable function,
which we present in \cref{definition:locally-definable-function} below.
For further details on locally definable functions we refer, \eg, to
\cite[Section~2]{MR3569073}.

\begin{athm}{definition}{definition:locally-definable-function}[Locally definable functions]
Let $ \mathcal S $
be an o-minimal structure\cfadd{definition:real-structure}\cfadd{definition:o-minimal-structure}\dott \cfload[.]
Let $ n \in \N $, let $ U \subseteq \R^n $ be open,
and let $ f \colon U \to \R $ be a function.
Then we say that $ f $ is a locally definable function with respect to $ \mathcal{S} $
(we say that $ f $ is locally definable) if and only if for all $ x \in U $
there exists $ \varepsilon \in (0,\infty) $ such that
$ f|_{ \{ y \in U \colon | x - y | < \varepsilon \} } $ is
definable with respect to $ \mathcal{S} $\dott \cfout[.]
\end{athm}

In the next step we briefly recall the concept
of the length of a function from an interval to a subset of $ \R^n $. For further reading
on curve lengths we refer, \eg, to \cite[Section~2.3.1]{BBI01}.

\NewDocumentCommand{\curvelength}{m}{%
  \operatorname{len}(#1)\cfadd{definition:curve-length}%
}
\begin{athm}{definition}{definition:curve-length}[Length of a curve]
Let $ I \subseteq \R $ be an interval,
let $ n \in \N $,
let $ U \subseteq \R^n $ be a set,
and let $ \gamma \colon I \to U $ be a function\dott \cfload[.]
\cfconsiderloaded{definition:curve-length}%
Then we denote by $ \curvelength{ \gamma } \in [0,\infty] $ the extended real number
given by
\begin{multline}
\llabel{eq1}
  \curvelength{ \gamma }
  =
  \sup\Bigl(
    \bigl\{ 0 \bigr\}
    \cup
    \Bigl\{
      r \in \R
      \colon \,
      \exists \, m \in \N, t_0, t_1, \ldots, t_m \in I
      \colon
      \Bigl(
        t_0 < t_1 < \ldots < t_m
\\
\textstyle
        \text{ and }
        r = \sum\limits_{ k = 1 }^m \norm{ \gamma( t_k ) - \gamma( t_{ k - 1 } ) }
      \Bigr)
    \Bigr\}
  \Bigr) ,
\end{multline}
we call $ \curvelength{ \gamma } $ the length of $ \gamma $,
and we say that $ \gamma $ has finite length if and only if $ \curvelength{ \gamma } < \infty $\dott \cfout[.]
\end{athm}

\subsection{Connections to the training of deep neural networks (DNNs)}
\label{sec:connections_to_deep_learning}

Even though \cref{maintheorem}
is a result in the area of pure mathematics,
it also has strong connections
to applied mathematics, optimization, and computer science, specifically, to \DL, that is, to the training
of \DNNs\ by means of gradient descent based optimization algorithms.

\GF\ trajectories are the basis for many optimization procedures
in the sense that, under suitable assumptions, the solutions $ \theta(t) $, $ t \geq 0 $, of \cref{eq:first_theorem_ODE}
(with the time interval being $ I = [0,\infty) $) converge to minimizers of the potential $ F \colon U \to \R $
as $ t $ tends to infinity. \KL-type inequalities and related techniques
are also commonly used to
establish convergence
of discretizations of \cref{eq:first_theorem_ODE} (gradient based optimization methods)
to critical points of $ F $ under, among other hypotheses, the assumption that the considered optimization process
is bounded (remains in a compactum).

In the training of \DNNs\ by means of \SGD\ based optimization methods
(such as standard \SGD, \Adam~\cite{KingmaBa15}, and \Muon~\cite{JordanEtAl24Muon})
the situation is, however, more subtle
as extensive numerical simulations \cite{JKR25,DJR25,KGDJ25},
complemented by analytical results in specific settings
\cite{Cai25,ChizatBach20,GJL22,JT20,KGDJ25,LyuLi20,Soudry18},
strongly suggest that, in most cases of \DNN\ training,
the \GF\ trajectories and their discretizations by means of \SGD\ methods
neither converge to a critical point nor remain bounded,
but diverge to infinity (though note that the divergence to infinity of \GF\ trajectories can only 
occur at a rather slow speed of
divergence; see, \eg, \cite[Proposition~2.54]{GJL22}).

In such a diverging scenario it is, of course, no longer possible to conclude convergence of the \GF\ trajectory
but \cref{maintheorem} allows us
in broad generality (definable and continuously differentiable activation and loss)
to conclude
\emph{directional convergence} (convergence of the spherical projection)
for the \DNN\ training flow.
This is precisely the subject of the next result, \cref{maincor2},
which is an immediate consequence of \cref{maintheorem}.

\begin{tcolorbox}[theorembox]
\begin{athm}{corollary}{maincor2}
\label{cor:DNN_training}
Let $ m, L \in \N $, $ \ell = ( \ell_0, \ell_1, \dots, \ell_L ) \in \N^{ L + 1 } $,
let
$
  n =
  \sum_{ i = 1 }^L
  \ell_i ( \ell_{ i - 1 } + 1 )
$,
let $ A \in C^1( \R, \R ) $, $ \mathscr{L} \in C^1( \R^{ \ell_L } \times \R^{ \ell_L } \times \R^n , [ 0, \infty ) ) $ be definable,
for every $ \theta = ( \theta_1, \dots, \theta_n ) \in \R^n $ let
$
  \rel{ v }{ \theta }
  = ( \rell{ v }{ \theta }{ 1 }, \dots, \rell{ v }{ \theta }{ \ell_v } )
  \colon \R^{ \ell_0 } \to \R^{ \ell_v }
$,
$
  v \in \{ 0, 1, \dots, L \}
$,
satisfy for all $ v \in \{ 0, 1, \dots, L - 1 \} $,
$ x = ( x_1, \dots, x_{ \ell_0 } ) \in \R^{ \ell_0 } $,
$ i \in \{ 1, 2, \dots, \ell_{ v + 1 } \} $
that
\begin{equation}
\label{eq:neural_network_def}
\begin{split}
  \rell{ v + 1 }{ \theta }{ i }( x )
&
  =
  \theta_{ \ell_{ v + 1 } \ell_{ v } + i + \sum_{ h = 1 }^v \ell_h ( \ell_{ h - 1 } + 1 ) }
\\ & \quad
  + \textstyle\sum\limits_{ j = 1 }^{ \ell_v }
  \theta_{ ( i - 1 ) \ell_v + j + \sum_{ h = 1 }^v \ell_h ( \ell_{ h - 1 } + 1 ) }
  \Bigl(
    x_{ \min\{ j, \ell_0 \} }
    \indicator{ \{ 0 \} }( v )
    +
    A( \rell{ v }{ \theta }{ j }( x ) )
    \indicator{ \N }( v )
  \Bigr)
  ,
\end{split}
\end{equation}
let
$ \vartheta \in \R^n $,
$ x_1, x_2, \dots, x_m \in \R^{ \ell_0 } $,
$ y_1, y_2, \dots, y_m \in \R^{ \ell_L } $,
let $ F \in C( \R^n, \R ) $ and $ P \colon \R^n \to \R^n $
satisfy for all $ \theta \in \R^n \backslash \{ \vartheta \} $ that
\begin{equation}
  F( \theta ) =
  \frac{ 1 }{ m }
  \sum_{ i = 1 }^m \mathscr{L}\bigl( \rel{ L }{ \theta }( x_i ) , y_i, \theta \bigr)
\end{equation}
and \EMPH{$ P( \theta ) = \norm{ \theta - \vartheta }^{ - 1 } ( \theta - \vartheta ) $},
and let $ \theta \in C^1( [0,\infty), \R^n ) $
satisfy for all $ t \in [0,\infty) $ that
\EMPH{$ \theta'(t) = - \bigl(\nabla F\bigr)\bigl( \theta(t) \bigr) $}.
Then there exists a unique $ v \in \R^n $ such that
$
  \limsup_{ t \to \infty } | ( P \circ \theta )( t ) - v | = 0
$.
\end{athm}
\end{tcolorbox}

\cref{maincor2} is proved
as \cref{cor:neural-network-training-full.item4} of \cref{cor:neural-network-training-full}
in \cref{sec:directional-convergence} below.
In order to state the conclusion of \cref{cor:DNN_training} we need
to specify realizations of \DNNs, which we do in \cref{eq:neural_network_def}
(as, \eg, in \cite[Definition~1.1]{DoJR25}).
Many of the activation and loss functions used in \DL\ are definable and continuously differentiable
and thus satisfy the hypotheses of \cref{cor:DNN_training}; cf., \eg, \cite[Example~2.4]{KGDJ25} and \cite[Sections~1.2 and 5.4]{JKW25Book}.

\subsection{Outline of the proof and structure of this article}
\label{subsec:proof-outline}

In this subsection we provide a brief outline of the proof of
\cref{maintheorem} and describe the structure of the
remainder of this article.
\Cref{maintheorem} is proved as \cref{maintheorem3} in
\cref{sec:applications} below.
In the first step of the proof of \cref{maintheorem3} we reduce
\cref{maintheorem3} (by separately considering the gradient trajectories
$\theta(t)$ and $\theta(-t)$ of the potentials $F$ and $-F$, respectively)
to the special case where the time interval
$I$ is equal to the
half-open interval $I=[\tau,T)$ for some $T\in(\tau,\infty]$, $\tau\in\R$.
This special case is proved as \cref{maintheorem2} in
\cref{subsec:open-domain-length} below.

\Cref{maintheorem2}, in turn, is proved by reducing it to the special case
where the \GF\ trajectory $\theta(t)$, $t\in[\tau,T)$, never hits the
reference point $\vartheta$ for the spherical projection in the sense that
$\{t\in[\tau,T)\colon\theta(t)=\vartheta\}=\emptyset$.
This special case is proved as \cref{prop:spherical-length-away} in
\cref{subsec:open-domain-length}.
The challenging cases in the proof of \cref{prop:spherical-length-away}
are the case in which
$\limsup_{t\to T}\norm{\theta(t)-\vartheta}=0$
(convergence to the reference point $\vartheta\in\R^n$) -- the case of the
\emph{definable Thom gradient conjecture} -- and the case in which
$\liminf_{t\to T}\norm{\theta(t)}=\infty$
(divergence to infinity) -- the case of the
\emph{gradient conjecture at infinity}.
All other cases can be handled by employing the well-established classical
\KL\ theory \cite[Theorem~1]{K98} (see \cref{maintheorem1} in
\cref{subsec:finite-and-infinite-limits} below for details).

We treat these two challenging cases in the proof of
\cref{prop:spherical-length-away} jointly by applying an appropriate
transformation to the \GF\ \DE\ in
\eqref{eq:first_theorem_ODE} above.
Specifically, we transform \eqref{eq:first_theorem_ODE} to \DE\ systems on
definable $C^2$-submanifolds of $\R^{n+1}$ of the form
\begin{eqn}
  \label{eq:intro-radial-angular-system}
  r'(t)=-\alpha(t)(r(t))^2(\partial_r h)(r(t),x(t))
  \qquad\text{and}\qquad x'(t)=-\alpha(t)(\nabla_Mh)(r(t),x(t))
\end{eqn}
for Lebesgue-almost all $t\in I=[\tau,T)$ and it is precisely the subject of
\cref{sec:radial-integration} to study the dynamics of \DE\ systems of the form
\eqref{eq:intro-radial-angular-system}.
Here, $\alpha\colon I\to[0,\infty)$ is a Lebesgue measurable function,
$h\colon\mathcal O\to\R$ is a definable $C^1$-function on a relatively
open subset $\mathcal O$ of $(0,\mathfrak R)\times M$ where
$\mathfrak R\in(0,\infty)$ is a positive real number and where
$M\subseteq\R^n$ is a compact and definable $C^2$-submanifold of $\R^n$,
and $r\colon I\to(0,\mathfrak R)$ and $x\colon I\to M$ are locally
absolutely continuous functions such that the solution pair
$(r(t),x(t))$, $t\in I$,
takes values in $\mathcal O$.

In the setting of \cref{prop:spherical-length-away} in \cref{sec:applications}
below and \cref{maintheorem} above we
have that the submanifold $M=\{v\in\R^n\colon\norm{v}=1\}$ is the unit
sphere in $\R^n$, we have that $x(t)=P(\theta(t))$ is the spherical
projection of $\theta(t)$, we have that
$r(t)=\norm{\theta(t)-\vartheta}^\nu$ with $\nu=1$ in the case where
$\limsup_{t\to T}\norm{\theta(t)-\vartheta}=0$ and with $\nu=-1$ in the
case where $\liminf_{t\to T}\norm{\theta(t)}=\infty$, we have that
$\alpha(t)=(r(t))^{-2\nu}(1+(F(\theta(t)))^2)^{3/2}$ for $t\in I$,
we have that
$\mathcal O=\{(s,u)\in(0,\mathfrak R)\times M\colon\vartheta+s^\nu u\in U\}$,
and we have that
\begin{eqn}
  \label{eq:intro-normalized-potential}
  h(s,u)=\varphi(F(\vartheta+s^\nu u))
\end{eqn}
for $(s,u)\in\mathcal O$,
where $\varphi\colon\R\to(-1,1)$ satisfies for all $e\in\R$ that
$\varphi(e)=(1+e^2)^{-1/2}e$.
We note that in both cases, $\nu=1$ and $\nu=-1$, we have that
$\limsup_{t\to T}\abs{r(t)}=0$.

As $r(t)$ represents the radius $\norm{\theta(t)-\vartheta}$ (in the case $\nu=1$)
and its inverse (in the case $\nu=-1$) and as $x(t)$ represents the angular component
of $\theta(t)$, we refer to \eqref{eq:intro-radial-angular-system} as a
\emph{radial-angular \DE\ system} and it is the key objective of
\cref{sec:radial-integration} to establish that $x(t)$ in
\eqref{eq:intro-radial-angular-system} has finite length.
Specifically, \itref{cor:angular-length}{item2} proves under suitable
assumptions that $\int_\tau^T\norm{x'(t)}\dd t<\infty$ and this is used in
the proof of \cref{prop:spherical-length-away} to establish that
$P\circ\theta$ has finite length.
Our proof of \cref{cor:angular-length} is based on applications of
\cref{prop:buildQ} in \cref{subsec:bounded-state} below,
\cref{lem:radial-integrability} in \cref{subsec:radial-integrability} below,
and \cref{prop:cone-annular-length} in \cref{subsec:cone-annular-occupation} below.

On a technical level, the most fundamental contribution of this work is
\cref{prop:buildQ}, which establishes a certain \emph{Lyapunov-type stability}
for solutions $(r(t),x(t))\in\mathcal O$, $t\in I$, of the radial-angular
\DE\ systems.
Specifically, \eqref{prop:buildQ.concl2} in
\cref{prop:buildQ} shows that there exists
$\eta\in(0,\mathfrak R)$ such that for all $a,b\in I$ with $a\leq b$
and $r([a,b])\subseteq(0,\eta)$ it holds that
\begin{eqn}
  \label{eq:intro-radial-angular-lyapunov}
  Q(r(b),x(b))-Q(r(a),x(a))
    \leq-\int_a^b\norm{x'(t)}\dd t
      +\int_a^b m(r(t))\abs{r'(t)}\dd t
  \dott
\end{eqn}
The function $m\colon(0,\eta)\to[0,\infty)$ in
\eqref{eq:intro-radial-angular-lyapunov} is a continuous, integrable, and
definable function and the function
$Q\colon\mathcal O\cap((0,\eta)\times M)\to[0,\infty)$ in
\eqref{eq:intro-radial-angular-lyapunov} is continuous and bounded.
Note that \eqref{eq:intro-radial-angular-lyapunov} shows that $Q$ serves
as a Lyapunov-type
function for the radial-angular \DE\ system in
\eqref{eq:intro-radial-angular-system}
up to the radial error term in the second integral on the right hand
side of \eqref{eq:intro-radial-angular-lyapunov}.
Furthermore, we observe that rearranging
\eqref{eq:intro-radial-angular-lyapunov} shows that for all $a,b\in I$ with
$a\leq b$ and $r([a,b])\subseteq(0,\eta)$ it holds that
\begin{eqn}
  \label{eq:intro-angular-length-from-lyapunov}
  \int_a^b\norm{x'(t)}\dd t
    \leq \int_a^b m(r(t))\abs{r'(t)}\dd t
      +\biggl[{\textstyle\sup\limits_{v\in\mathcal O\cap((0,\eta)\times M)}Q(v)}\biggr]
  \dott
\end{eqn}
Note that \eqref{eq:intro-angular-length-from-lyapunov} is only
applicable to subintervals $[a,b]\subseteq I$ of $I$ on which the radial
component $r$ is smaller than $\eta$ in the sense that
$r([a,b])\subseteq(0,\eta)$.
However, as $r$ converges to zero in the sense that
$\limsup_{t\to T}\abs{r(t)}=0$, we obtain that there exists
$\sigma\in I$ which satisfies
$r(I\cap[\sigma,\infty))\subseteq(0,\eta)$.
In the proof of \cref{cor:angular-length} we also employ
\cref{lem:radial-integrability} and \cref{prop:cone-annular-length} to show,
under suitable assumptions, that the integral on the right hand side of
\eqref{eq:intro-angular-length-from-lyapunov} can be uniformly bounded.
Combining this with \eqref{eq:intro-angular-length-from-lyapunov}
and the fact that $Q$ is bounded then proves the finite length estimate
for the tail $\int_\sigma^T\norm{x'(t)}\dd t<\infty$.
This and the fact that for every locally absolutely continuous
$y\colon I\to\R^n$ we have that
$\int_\tau^\sigma\norm{y'(t)}\dd t<\infty$ then allow us to establish, in
\itref{cor:angular-length}{item2}, the finite length estimate
$\int_\tau^T\norm{x'(t)}\dd t
=\int_\tau^\sigma\norm{x'(t)}\dd t
+\int_\sigma^T\norm{x'(t)}\dd t<\infty$.

The central task in the proof of \cref{prop:buildQ} is thus
to construct the Lyapunov-type function $Q$ so that
\eqref{eq:intro-radial-angular-lyapunov} holds.
One could try to directly employ the function $h$ in
\eqref{eq:intro-radial-angular-system} as a Lyapunov-type function and
applying the fundamental theorem of calculus and the chain rule to the
Lyapunov-type energy $h(r(t),x(t))$, $t\in I$, shows that for all
$a,b\in I$ with $a<b$ it holds that
\begin{eqn}[a]
  \label{eq:intro-original-energy-dissipation}
  h(r(b),x(b))-h(r(a),x(a))&=\int_a^b\bigl((\partial_r h)(r(t),x(t))r'(t)
    +\langle(\nabla_Mh)(r(t),x(t)),x'(t)\rangle\bigr)\dd t\\
  &=-\int_a^b\alpha(t)\bigl((r(t))^2
    \abs{(\partial_r h)(r(t),x(t))}^2
    +\norm{(\nabla_Mh)(r(t),x(t))}^2\bigr)\dd t\\
  &\leq-\int_a^b\alpha(t)\norm{(\nabla_Mh)(r(t),x(t))}^2\dd t
  \dott
\end{eqn}
This shows that the Lyapunov energy $h(r(t),x(t))$, $t\in I$,
is non-increasing and that $h$ is a Lyapunov function for
\eqref{eq:intro-radial-angular-system} but the dissipation
$\int_a^b\alpha(t)\norm{(\nabla_Mh)(r(t),x(t))}^2\dd t<\infty$
in \eqref{eq:intro-original-energy-dissipation} associated to the
Lyapunov function $h$ is not large enough to directly establish
finiteness of the length of $x$, i.e.,
$\int_\tau^T\norm{x'(t)}\dd t
=\int_\tau^T\alpha(t)\norm{(\nabla_Mh)(r(t),x(t))}\dd t$
(see the second equality in \eqref{eq:intro-radial-angular-system}).
The problem is that the norm of the gradient
$(\nabla_Mh)(r(t),x(t))$ appears quadratically inside the integral in
\eqref{eq:intro-original-energy-dissipation} but in the length the norm
of the gradient $(\nabla_Mh)(r(t),x(t))$ appears only linearly inside
the integral and if the norm of the gradient
$\norm{(\nabla_Mh)(r(t),x(t))}$ is small, then
$\norm{(\nabla_Mh)(r(t),x(t))}^2$ is much smaller than
$\norm{(\nabla_Mh)(r(t),x(t))}$.

We want to make use of \eqref{eq:intro-original-energy-dissipation} and
our objective is to transform the energy $h(r(t),x(t))$ so that the
dissipation of the transformed energy/the transformed Lyapunov-type
function is large enough to control -- up to a radial error
(see \eqref{eq:intro-radial-angular-lyapunov} and
\eqref{eq:intro-angular-length-from-lyapunov} above) -- the length
$\int_\tau^T\norm{x'(t)}\dd t$.
The transformation may, however, depend not only on the Lyapunov energy
$h(r(t),x(t))$ but also on the value of the radial component $r(t)$ and
we consider Lyapunov-type functions of the form
$Q(r,x)=\Phi(r,h(r,x))$ for some suitable function $\Phi$.
Specifically, the approach initiated by \cref{lem:divide} in
\cref{sec:atlas} below is to construct a natural number $N\in\N$,
strictly positive reals $\rho\in(0,\mathfrak R)$, $H\in(0,\infty)$,
definable $[0,H]$-valued transformation functions
$\mathcal F_0,\mathcal F_1,\dots,\mathcal F_{N-1}$,
definable back-transformation functions $E_1,E_2,\dots,E_N$, and definable
energy boundary functions $\mathfrak e_0,\mathfrak e_1,\dots,\mathfrak e_N$
that specify pairwise disjoint regions in the radial-energy plane
$\mathbf{H}_{0},\mathbf{H}_{1},\dots,\mathbf{H}_{N-1}$
through
$\forall\,i\in\{0,1,\dots,N-1\}\colon
\mathbf{H}_{i}=\{(\mathfrak r,\mathfrak h)\in(0,\rho)\times\R\colon
\mathfrak e_i(\mathfrak r)<\mathfrak h<\mathfrak e_{i+1}(\mathfrak r)\}$
and to consider the transformed energy
\begin{eqn}
  \label{eq:intro-transformed-lyapunov-energy}
  \mathcal F_i\bigl(r(t),h\bigl(r(t),x(t)\bigr)\bigr)
\end{eqn}
in the regime where $(r(t),h(r(t),x(t)))\in\mathbf{H}_{i}$
for $t\in I$, $i\in\{0,1,\dots,N-1\}$.
In the following we denote for every $i\in\{0,1,\dots,N-1\}$ by
$Z_i(r,x)=\mathcal F_i(r,h(r,x))$ the transformed energy associated to
radial component $r$ and angular component $x$, with
$Z_i\colon\mathcal O\cap((0,\rho)\times M)\to[0,H]$ mapping from
$\mathcal O\cap((0,\rho)\times M)$ to $[0,H]$.
In \cref{lem:divide} we also supply crucial (lower and) upper bounds for
the partial derivative with respect to the transformed energy variable
$(\partial_z E_i)(r,z)$ of the back-transformations $E_i$,
$i\in\{1,2,\dots,N\}$.
Specifically, the fundamental theorem of calculus, the chain rule,
the inverse function differentiation rule,
the second equation in \eqref{eq:intro-radial-angular-system}, and
\cref{lem:divide} show that for all
$i\in\{0,1,\dots,N-1\}$, $a,b\in I$ with $a<b$ and
$\midcup_{t\in[a,b]}\bigl\{\bigl(r(t),h(r(t),x(t))\bigr)\bigr\}\subseteq\mathbf{H}_{i}$
it holds that
\begin{eqn}[a]
  \label{eq:intro-transformed-energy-dissipation}
  &Z_i(r(b),x(b))-Z_i(r(a),x(a))\\
  &\quad=\int_a^b\bigl((\partial_r Z_i)(r(t),x(t))r'(t)
    +\langle(\nabla_M Z_i)(r(t),x(t)),x'(t)\rangle\bigr)\dd t\\
  &\quad=\int_a^b(\partial_r Z_i)(r(t),x(t))r'(t)\dd t
    -\int_a^b\frac{\norm{(\nabla_Mh)(r(t),x(t))}\,\norm{x'(t)}}
      {(\partial_z E_{i+1})(r(t),Z_i(r(t),x(t)))}\dd t\\
  &\quad\leq\int_a^b(\partial_r Z_i)(r(t),x(t))r'(t)\dd t
    -\frac1c\int_a^b\norm{x'(t)}\dd t
  \dott
\end{eqn}
Here the final inequality uses
$0<(\partial_z E_{i+1})(r(t),Z_i(r(t),x(t)))
\leq c\norm{(\nabla_Mh)(r(t),x(t))}$ with a strictly positive real constant
$c\in(0,\infty)$ according to
\itref{lem:divide}{item5,item6}.
The inequality in \eqref{eq:intro-transformed-energy-dissipation} does not,
however, directly provide us with
\eqref{eq:intro-radial-angular-lyapunov} as
\eqref{eq:intro-transformed-energy-dissipation} only applies to solutions
with the radial component and the energy staying in one fixed
radial-energy region $\mathbf{H}_i$ and as the dependence of the
transformation functions $\mathcal F_0,\mathcal F_1,\dots,\mathcal F_{N-1}$
not only on the energy but also on the radial component produces an
additional radial error term (see the first integral on the right hand
side of \eqref{eq:intro-transformed-energy-dissipation}).

To this end, we establish in \cref{lem:residual} in \cref{sec:atlas} below
upper bounds for the absolute difference between the radial partial
derivative of the energy $\partial_r h$ and the radial partial derivative
of the back-transformations $\partial_r E_{i+1}$ involving nonnegative,
definable, and integrable functions $\mathbf m_i(r)$.
We use these estimates to explicitly specify suitable definable $C^1$
weight functions
$\omega_i\colon(0,\eta)\times(0,H)\to(0,1]$,
$i\in\{0,1,\dots,N-1\}$,
(see \eqref{prop:buildQ.setting4} in the proof of
\cref{prop:buildQ} below) that balance the radial error term
(see the first integral on the right hand side of
\eqref{eq:intro-transformed-energy-dissipation}) against the angular
dissipation term with the right sign
(see the second integral on the right hand side of
\eqref{eq:intro-transformed-energy-dissipation}), where
$\eta\in(0,\rho)$ is a sufficiently small strictly positive real number.
For every $i\in\{0,1,\dots,N-1\}$ we denote by
$P_i\colon(0,\eta)\times[0,H]\to[0,H]$ the one-parameter primitive
(with respect to the transformed energy variable) of the weight function
$\omega_i$ in the sense that for all $r\in(0,\eta)$, $z\in[0,H]$
we have that $P_i(r,z)=\int_0^z\omega_i(r,u)\dd u$
and the Lyapunov-type function in
\eqref{eq:intro-radial-angular-lyapunov} is then, for every
$(r,x)\in\mathcal O\cap((0,\eta)\times M)$, given by
\begin{eqn}
  \label{eq:intro-lyapunov-construction}
  Q(r,x)=\kappa\sum_{i=0}^{N-1}P_i(r,Z_i(r,x))
    =\kappa\sum_{i=0}^{N-1}P_i\bigl(r,\mathcal F_i(r,h(r,x))\bigr)
\end{eqn}
where $\kappa\in(0,\infty)$ is a sufficiently large real number.
\Cref{lem:primitive} in \cref{sec:radial-integration} provides regularity
properties for $Q$ and, in particular, allows us to conclude that $Q$ is
continuous.
Moreover, the fact that for all $i\in\{0,1,\dots,N-1\}$ we have that
$\operatorname{im}(P_i)\subseteq[0,H]$ shows for all
$(r,x)\in\mathcal O\cap((0,\eta)\times M)$ that
$0\leq Q(r,x)\leq\kappa NH$.
Hence, $Q$ is a bounded and continuous Lyapunov-type function.

To prove \eqref{eq:intro-radial-angular-lyapunov}, we also need to handle
the boundaries between the radial-energy regions, at which the energy
transformation functions $\mathcal F_0,\mathcal F_1,\dots,\mathcal F_{N-1}$
need not be differentiable.
To solve this issue we employ a \emph{truncation and approximation procedure} for the
transformed energy $Z_i(r,x)$.
Specifically, we approximate $Z_i(r,x)$ by
$\max\{\varepsilon,\min\{H-\varepsilon,Z_i(r,x)\}\}$ for
$\varepsilon\in(0,H/2)$ and, thereby, also obtain approximations for the
Lyapunov-type function $Q$ in \eqref{eq:intro-lyapunov-construction} and
in \cref{lem:truncated} in
\cref{sec:truncation} we establish regularity properties and estimates for the
approximated quantities and functions.
These regularity properties and estimates allow us to apply the
fundamental theorem of calculus for absolutely continuous functions to
the approximated quantities and, thereafter, to pass to the limit as the
truncation parameter $\varepsilon$ tends to zero.
For the times $t\in I$ spent exactly at the boundaries of the
radial-energy regions in the sense that
$h(r(t),x(t))=\mathfrak e_i(r(t))$, we apply the elementary result in
\cref{lem:AC-level-set} (which shows for every $c\in\R$ that the
derivative of an absolutely continuous function vanishes almost
everywhere on its $c$-level set) to the real valued function
$t\mapsto h(r(t),x(t))-\mathfrak e_i(r(t))$.

In the proof of \cref{prop:buildQ} we actually perform
the arguments and transformation procedures
\eqref{eq:intro-transformed-lyapunov-energy}--\eqref{eq:intro-lyapunov-construction}
outlined above
not directly with the energy function $h\colon\mathcal O\to\R$ from
\eqref{eq:intro-radial-angular-system} but with the normalized energy function
$\mathbf h=\varphi\circ h\colon\mathcal O\to(-1,1)$.
(This is possible as a solution $(r(t),x(t))\in\mathcal O$, $t\in I$, of
\eqref{eq:intro-radial-angular-system} is also a solution of
\eqref{eq:intro-radial-angular-system} but with
$h\colon\mathcal O\to\R$ replaced by
$\mathbf h\colon\mathcal O\to(-1,1)$ and
$\alpha\colon I\to[0,\infty)$ replaced by a modified
$\widetilde\alpha\colon I\to[0,\infty)$; see
\eqref{prop:buildQ.eq0.3} in the proof of
\cref{prop:buildQ} for details.)
However, if $h$ emerges from the polar coordinate type transformation in the proof of
\cref{prop:spherical-length-away}, in which case we have for all
$(s,u)\in\mathcal O$ that
$h(s,u)=\varphi(F(\vartheta+s^\nu u))$ so that $h$ already takes values
in $(-1,1)$ (see \eqref{eq:intro-normalized-potential} above), then this
additional normalization step is not needed and the
procedure
\eqref{eq:intro-transformed-lyapunov-energy}--\eqref{eq:intro-lyapunov-construction}
is directly applicable as described above.

\section{Regularity properties of radial upper bounds for definable \texorpdfstring{$C^1$}{C1}-scalar fields}
\label{sec:bounds}

In the main results of this section, \cref{lem:normalized} and
\cref{lem:contact-bound} in
\cref{subsec:normalized-bound} below, we establish integrability,
smoothness, and definability properties of radial upper bounds for
definable $C^1$-scalar fields on definable $C^r$-submanifolds of $\R^n$.

In particular, in \cref{lem:normalized} we consider a definable
$C^1$-scalar field $Z=(Z(r,x))_{(r,x)\in U}\colon U\to\R$ on a
relatively open subset $U\subseteq(0,\mathfrak r_0)\times M$ of the
Cartesian product of an open interval $(0,\mathfrak r_0)$ and a bounded
and definable $C^1$-submanifold $M$ of $\R^n$ and establish an upper bound
for the partial derivative with respect to the radial variable
$(\partial_r Z)(r,x)$ using a definable and (locally) integrable $C^1$-function
$\eta(r)$.
We employ \cref{lem:normalized} in the proof of \cref{lem:residual}
(which establishes lower and upper bounds for the partial derivative
with respect to the radial variable of the energy back-transformations)
in \cref{sec:atlas} below and (through its consequence in
\cref{lem:normalized-interval}) in
the proof of \cref{lem:primitive} (which establishes regularity properties
for one-parameter primitives of definable two-parameter functions,
in particular, for the weight functions
$\omega_i\colon(0,\tau)\times(0,H)\to(0,1]$ for
$i\in\{0,1,\dots,N-1\}$; cf.\ \eqref{eq:intro-lyapunov-construction})
in \cref{sec:radial-integration} below.
\Cref{lem:residual} and \cref{lem:primitive} are both used in the proof
of \cref{prop:buildQ} to establish Lyapunov-type stability
for the radial-angular \DE\ system.

In \cref{lem:contact-bound} we consider a definable $C^1$-function
$a\colon(0,\mathfrak r_0)\to\R$ and a definable $C^1$-scalar field
$h=(h(r,x))_{(r,x)\in\mathcal O}\colon\mathcal O\to\R$ on a relatively
open subset $\mathcal O\subseteq(0,\mathfrak r_0)\times M$ of
$(0,\mathfrak r_0)\times M$, where
$M\subseteq\R^n$ is a bounded and definable $C^1$-submanifold of $\R^n$, and we
establish an upper bound of the form
$|(\partial_r h)(r,x)-a'(r)|\leq\beta(r)|(\nabla_Mh)(r,x)|$
with the constraint $h(r,x)=a(r)$ for sufficiently small $r$ using
a definable and (locally) integrable $C^1$-function
$\beta(r)$.
We employ \cref{lem:contact-bound} in the proof of
\cref{prop:buildQ} (which establishes Lyapunov-type
stability for the radial-angular \DE\ system) to control the radial
error (cf.\ the second integral on the right hand side of
\eqref{eq:intro-radial-angular-lyapunov}) at times $t$ at which the energy
$h(r(t),x(t))$ coincides with an
energy boundary value in the sense that
$h(r(t),x(t))\in\midcup_{i=0}^N\{\mathfrak e_i(r(t))\}$.

Before we formulate and prove \cref{lem:normalized} and
\cref{lem:contact-bound}, we first recall several well-known
differential geometric preliminaries in
\cref{subsec:unary-variation,subsec:smooth-functions-manifolds} below
(cf., \eg, \cite[Chapter 3]{Boumal23}).
We include proofs for the statements in
\cref{subsec:unary-variation,subsec:smooth-functions-manifolds} for completeness.
In \cref{subsec:normalized-bound} we then use some of the auxiliary
results from \cref{subsec:smooth-functions-manifolds} to prove
\cref{lem:normalized} and
\cref{lem:contact-bound}.

\subsection{Embedded submanifolds}
\label{subsec:unary-variation}

In the following notion, \cref{definition:embedded-manifold}, we recall the standard
concept of an embedded $m$-dimensional $C^r$-submanifold of $\R^n$.
For further reading we refer, \eg, to
\cite[Definition~3.10 and Theorem~3.12]{Boumal23}.

\begin{athm}{definition}{definition:embedded-manifold}
  [Embedded submanifolds with specified dimension]
   Let $r,n\in\N$, $m\in\N_0$ and let $M\subseteq\R^n$ be a set\dott \cfload[.]
  Then we say that $M$ is an \dsubmanifold{r}{n}{m} if
  and only if for every $x\in M$ there exist an open
  $U\subseteq\R^m$, an open neighbourhood $V\subseteq\R^n$ of $x$,
  and a homeomorphism $\phi\in C^r(U,V\cap M)$ such that for all
  $u\in U$ it holds that
  \begin{eqn}
    \operatorname{rank}(\phi'(u))=m
    \dott
  \end{eqn}
  \cfout[.]
\end{athm}

\newcommand{\submanifold}[2]{%
  $C^{#1}$-submanifold of $\R^{#2}$\cfadd{definition:submanifold}%
}
\begin{athm}{definition}{definition:submanifold}[Embedded submanifolds with unspecified dimension]
  Let $r,n\in\N$ and let $M\subseteq\R^n$ be a set\dott \cfload[.]
  \cfconsiderloaded{definition:submanifold}%
  Then we say that $M$ is a \submanifold{r}{n} if and only if
  there exists $m\in\N_0$ such that $M$ is an
  \dsubmanifold{r}{n}{m}\cfadd{definition:embedded-manifold}\dott \cfout[.]
\end{athm}

In the next notion we recall the concept of the tangent space of a
submanifold at a point (cf., \eg, \cite[Definition~3.14]{Boumal23} for
further reading).

\begin{athm}{definition}{definition:tangent-space}[Tangent spaces]
  Let $n\in\N$, let $M\subseteq\R^n$ be a \submanifold{1}{n}, and
  let $x\in M$\dott \cfload[.]
  Then we denote by $T_xM$ the set given by
  \begin{eqn}
    T_xM=
    \bigl\{v\in\R^n\colon
      \bigl[\exists\,\varepsilon\in(0,\infty), \gamma\in
        C^1((-\varepsilon,\varepsilon),M)\colon
        \bigl((\gamma(0)=x)\wedge(\gamma'(0)=v)\bigr)\bigr]
    \bigr\}
    \dott
  \end{eqn}
  \cfout[.]
\end{athm}

\NewDocumentCommand{\CManifold}{m m m}{%
  \mathbf{C}^{#1}(#2,#3)\cfadd{definition:c1-on-manifold}%
}
\NewDocumentCommand{\COne}{m m}{%
  \CManifold{1}{#1}{#2}%
}

In the next lemma we recall a well-known local representation of a
$C^1$-submanifold of $\R^n$ by a diffeomorphism
(cf., \eg, \cite[Theorem~3.12]{Boumal23}).

\begin{athm}{lemma}{lem:embedded-manifold-local-straightening}
  [Local normal forms of embedded submanifolds]
  Let $n\in\N$, $m\in\N_0$, let $M\subseteq\R^n$ be an
  \dsubmanifold{1}{n}{m}\cfadd{definition:embedded-manifold}, and let
  $x\in M$\dott \cfload[.]
  Then there exist an open neighbourhood $U\subseteq\R^n$ of $x$,
  an open $V\subseteq\R^n$, and
  a $C^1$-diffeomorphism $\Phi\colon U\to V$ such that
  \begin{eqn}
    \Phi(x)=0\qquad\text{and}\qquad
    \Phi(M\cap U)=V\cap(\R^m\times\{0\}^{n-m})
    \dott
  \end{eqn}
  \cfout[.]
\end{athm}
\begin{aproof}
  Throughout this proof let $P\colon \R^n\to\R^{n-m}$ and
  $Q\colon \R^n\to\R^m$ satisfy for all $y=(y_1,\dots,y_n)\in \R^n$
 that
  \begin{eqn}
    P(y_1,\ldots,y_n)=(y_{m+1},\ldots,y_n)
    \qquad\text{and}\qquad
    Q(y_1,\ldots,y_n)=(y_1,\ldots,y_m)
    \dott
  \end{eqn}
  \cfload[.]
  \startnewargseq
  \argument{
    the assumption that $M$ is an \dsubmanifold{1}{n}{m};
  }{
    that there exist an open 
    $A\subseteq\R^m$, a point $\xi\in A$,
    a set $O\subseteq M$ open in $M$,
    and an injective $\phi\in C^1(A,\R^n)$ which satisfy that
    $\phi(\xi)=x\in M$,
    $\phi(A)=O$,
    $\phi\colon A\to O$ is a homeomorphism,
    and
    \begin{eqn}
      \llabel{eq0.0a.1}
     \forall\, u\in A\colon \operatorname{rank}\phi'(u)=m
      \dott
    \end{eqn}
  }
  \startnewargseq
  \argument{
    \lref{eq0.0a.1};
  }{
    that there exist distinct $i_1,i_2,\ldots,i_m\in\{1,2,\ldots,n\}$
    and $\pi\colon\R^n\to\R^m$
    which satisfy that for all $y=(y_1,y_2,\ldots,y_n)\in\R^n$ it holds that
    \begin{eqn}
      \pi(y)=(y_{i_1},y_{i_2},\ldots,y_{i_m})
    \end{eqn}
    and that
    \begin{eqn}
      \llabel{eq0.0a.minor}
      \pi\circ \phi'(\xi)
    \end{eqn}
    is a linear isomorphism\dott
  }
  \startnewargseq
  \argument{
    \lref{eq0.0a.minor};
  }{
    that there exist distinct
    $j_1,j_2,\ldots,j_{n-m}\in\{1,2,\ldots,n\}\backslash\{i_1,i_2,\ldots,i_m\}$
    and a linear isomorphism $L\colon\R^n\to\R^m\times\R^{n-m}$
    which satisfy that for all $y=(y_1,y_2,\ldots,y_n)\in\R^n$ it holds that
    \begin{eqn}
      \llabel{def: L}
      L(y)=\bigl((y_{i_1},y_{i_2},\ldots,y_{i_m}),
        (y_{j_1},y_{j_2},\ldots,y_{j_{n-m}})\bigr)
    \end{eqn}
    and that
    \begin{eqn}
      \llabel{eq0.0a.2}
      (Q\circ L\circ\phi)'(\xi)
      =\pi\circ \phi'(\xi)
    \end{eqn}
    is a linear isomorphism\dott
  }
  In the following let $h=Q\circ L\circ\phi$.
  \startnewargseq
  \argument{
    \lref{eq0.0a.2};
    the fact that
    $h=Q\circ L\circ\phi$;
  }{
    that $h\in C^1(A,\R^m)$ and its derivative
    \begin{eqn}
      \llabel{eq0.0a.2a}
      h'(\xi)
      =(Q\circ L\circ\phi)'(\xi)
      \colon \R^m\to\R^m
    \end{eqn}
    is a linear isomorphism\dott
  }
  \argument{
    \lref{eq0.0a.2a};
    the fact that $A$ is open and contains
    $\xi$;
    the $C^1$ inverse-function theorem applied to $h$ at
    $\xi$;
  }{
    that there exist open $W,B\subseteq\R^m$ which satisfy that
    $\xi\in W\subseteq A$,
    $h(\xi)\in B$, $h(W)=B$, and the restriction
    \begin{eqn}
      \llabel{eq0.0a.2b}
      \mathbf h=h|_{W}\colon W\to B
    \end{eqn}
    is bijective with $\mathbf h\in C^1(W,B)$ and
    $\mathbf h^{-1}\in C^1(B,W)$\dott
  }
  \startnewargseq
  In the following let $a=\mathbf h(\xi)\in B$, $b=P(L(x))\in\R^{n-m}$, and let $f\colon B\to\R^{n-m}$ satisfy $f=P\circ L\circ\phi\circ
    \mathbf h^{-1}$.
  \argument{
    \lref{def: L};\lref{eq0.0a.2b};
    the fact that $\phi(\xi)=x$;
    the fact that $h=Q\circ L\circ\phi$
  }{
    that
    \begin{eqn}\llabel{eqt1}
      L(x)=(\mathbf h(\xi),P(L(x)))=(a,b)
      \dott
    \end{eqn}
  }
  \argument{
    the fact that $a=\mathbf h(\xi)$;
    the fact that
    $\mathbf h^{-1}\circ\mathbf h=\Id_{W}$;
  }{
    that
    \begin{eqn}
      \llabel{eq0.0a.2c}
      \mathbf h^{-1}(a)
      =\mathbf h^{-1}(\mathbf h(\xi))
      =\xi
      \dott
    \end{eqn}
  }
  \argument{
    \lref{eq0.0a.2c};
    the fact that $P\in C^1(\R^n,\R^{n-m})$,
    $L\in C^1(\R^n,\R^m\times\R^{n-m})$,
    $\phi\in C^1(A,\R^n)$, and $\mathbf h^{-1}\in C^1(B,W)$;
    the fact that $\phi(\xi)=x$;
    the fact that $b=P(L(x))$; the fact that $f=P\circ L\circ\phi\circ
    \mathbf h^{-1}$
  }{
    that $f\in C^1(B,\R^{n-m})$ and that
    \begin{eqn}
      \llabel{eq0.0a.2d}
      f(a)
      =P(L(\phi(\mathbf h^{-1}(a))))
      =P(L(x))
      =b
      \dott
    \end{eqn}
  }
  \argument{
    \lref{eq0.0a.2b};
    the fact that
    $\mathbf h\circ\mathbf h^{-1}=\Id_{B}$;
    the fact that for all $u\in A$ it holds that
    $h(u)=Q(L(\phi(u)))$;
    the fact that
    $f=P\circ L\circ\phi\circ\mathbf h^{-1}$;
    the fact that for all $y\in\R^n$ it holds that $(Q(y),P(y))=y$;
  }{
    that for all $u\in B$ it holds that
    \begin{eqn}
      \llabel{eq0.0a.3}
      (L\circ\phi\circ\mathbf h^{-1})(u)
      =\bigl(\mathbf h(\mathbf h^{-1}(u)),f(u)\bigr)
      =(u,f(u))
      \dott
    \end{eqn}
  }
  \argument{
    the fact that $W$ is open in $A$;
    the fact that $\phi\colon A\to O$ is a homeomorphism;
    the fact that $O$ is open in $M$;
    the fact that $\xi\in W$; the fact that $\phi(\xi)=x$;
  }{
    that there exists an open $C\subseteq\R^n$ which satisfies
    \begin{eqn}
      \llabel{eq0.0a.4}
      x\in C
      \qquad\text{and}\qquad
      M\cap C=\phi(W)
      \dott
    \end{eqn}
  }
  \startnewargseq
  \argument{
    \lref{eqt1};\lref{eq0.0a.4};
    the fact that $L$ is a homeomorphism;
    the fact that $a\in B$ ; the fact that $B$ is open ;
  }{
    that there exist open $K\subseteq B$, $H\subseteq\R^{n-m}$
    which satisfy
    \begin{eqn}
      \llabel{eq0.0a.5}
      (a,b)\in K\times H\subseteq L(C)
      \dott
    \end{eqn}
  }
  In the following let $G\subseteq B$, $U\subseteq\R^n$ satisfy
  \begin{eqn}
    \llabel{eq0.0a.6}
    G=K\cap f^{-1}(H)
    \qquad\text{and}\qquad
    U=L^{-1}(G\times H)
    \dott
  \end{eqn}
  \startnewargseq
  \argument{  \lref{eqt1};  \lref{eq0.0a.2d};  \lref{eq0.0a.5};
    \lref{eq0.0a.6};
    the fact that $B$ and $K$ are open;
    the fact that $H$ is open;
    the fact that $f$ is continuous;
    the fact that $L\colon\R^n\to\R^m\times\R^{n-m}$ is a homeomorphism;
  }{
    that $G$ is open, $U$ is open, and
    \begin{eqn}
      \llabel{eq0.0a.7}
      a\in G\subseteq B,
      \qquad x\in U\subseteq C,
      \qquad\text{and}\qquad f(G)\subseteq H
      \dott
    \end{eqn}
  }
  \argument{
    \lref{eq0.0a.7};
    \lref{eq0.0a.3};
    \lref{eq0.0a.4};
    \lref{eq0.0a.6};
    the fact that $\mathbf h^{-1}(B)=W$;
    the fact that $L$ is bijective;
  }{
    that
    \begin{eqn}
      \llabel{eq0.0a}
      L(M\cap U)
        =L(\phi(W))\cap(G\times H)
        =\{(u,f(u))\in G\times H\colon u\in G\}
      \dott
    \end{eqn}
  }
In the following let
  \begin{eqn}
    \llabel{eq0.0b}
    V=\{(u-a,v-f(u))\in\R^m\times\R^{n-m}
      \colon (u,v)\in G\times H\}
  \end{eqn}
  and let $\Phi\colon U\to V$ satisfy for all
  $(u,v)\in G\times H$ that
  \begin{eqn}
    \llabel{eq0.0c}
    \Phi(L^{-1}(u,v))=(u-a,v-f(u))
    \dott
  \end{eqn}
  \startnewargseq
  \argument{\lref{eqt1};
    \lref{eq0.0b};
    \lref{eq0.0c};
    the fact that $b=f(a)$;
    the fact that $a\in G\subseteq B$;
    the fact that $G$ is open and $H$ is open;
    the fact that $f\in C^1(B,\R^{n-m})$;
  }{
    that
    \begin{eqn}
      \llabel{eq0.0b.1}
      V=\bigl\{(s,t)\in\R^m\times\R^{n-m}\colon
        \bigl[(s+a\in G)\wedge(t+f(s+a)\in H)\bigr]\bigr\},
    \end{eqn}
    $V$ is open, $0\in V$, $\Phi(x)=0$, and $\Phi\in C^1(U,V)$\dott
  }
  In the following let $\Theta\colon V\to U$ satisfy for all $(s,t)\in V$
  that
  \begin{eqn}
    \llabel{eq0.0b.2}
    \Theta(s,t)=L^{-1}\bigl(s+a,t+f(s+a)\bigr)
    \dott
  \end{eqn}
  \startnewargseq
  \argument{
    \lref{eq0.0b.1};
    the fact that $U=L^{-1}(G\times H)$;
    the fact that $G\subseteq B$;
    the fact that $f\in C^1(B,\R^{n-m})$;
  }{
    that $\Theta\in C^1(V,U)$\dott
  }
  \argument{
    \lref{eq0.0c};
    \lref{eq0.0b.2};
  }{
    that for all $(u,v)\in G\times H$ and $(s,t)\in V$ it holds that
    \begin{eqn}
      \llabel{eq0.0c.1}
      \Theta\bigl(\Phi(L^{-1}(u,v))\bigr)
      =L^{-1}(u,v)
      \qquad\text{and}\qquad
      \Phi(\Theta(s,t))=(s,t)
      \dott
    \end{eqn}
  }
  \argument{
    \lref{eq0.0c.1};
    the fact that $\Phi\in C^1(U,V)$ and $\Theta\in C^1(V,U)$;
  }{
    that $\Phi\colon U\to V$ is a $C^1$ diffeomorphism and it holds
    for all $(s,t)\in V$ that
    \begin{eqn}
      \llabel{eq0.0d}
      \Phi^{-1}(s,t)
      =\Theta(s,t)
      =L^{-1}\bigl(s+a,t+f(s+a)\bigr)
      \dott
    \end{eqn}
  }
  \argument{
    the fact that $f(G)\subseteq H$;
    the fact that $U=L^{-1}(G\times H)$;
    \lref{eq0.0c};
  }{
    that for all $u\in G$ it holds that $(u,f(u))\in G\times H$,
    $L^{-1}(u,f(u))\in U$, and
    \begin{eqn}
      \llabel{eq0.0c.2}
      \Phi(L^{-1}(u,f(u)))
      =(u-a,f(u)-f(u))
      =(u-a,0)
      \dott
    \end{eqn}
  }
  \argument{
    \lref{eq0.0c.2};
    \lref{eq0.0a};
    \lref{eq0.0b};
  }{
    that
    \begin{eqn}
      \llabel{eq0.0}
      \Phi(M\cap U)
      =V\cap\bigl(\R^m\times\{0\}^{n-m}\bigr)
      \dott
    \end{eqn}
  }
\end{aproof}

\subsection{Smooth functions on embedded submanifolds}
\label{subsec:smooth-functions-manifolds}

In the following notion we recall the concept of the set of all
$S$-valued $r$-times continuously differentiable functions on an
$m$-dimensional $C^r$-submanifold $M$ of $\R^n$
(cf., \eg, \cite[Definition~8.5 and Remark~8.6]{Boumal23}).

\begin{athm}{definition}{definition:c1-on-manifold}
  [Smooth functions on embedded submanifolds]
  Let $r,k,n\in\N$, $m\in\N_0$, let $S\subseteq\R^k$ be a set, and let
  $M\subseteq\R^n$ be an \dsubmanifold{r}{n}{m}\cfadd{definition:embedded-manifold}\dott
  \cfload[.]
  \cfconsiderloaded{definition:c1-on-manifold}%
  Then we denote by $\CManifold{r}{M}{S}$ the set of all functions
  $f\colon M\to S$ such that for every open $U\subseteq\R^m$,
  every open $V\subseteq\R^n$, and every
  $\phi\in C^r(U,V\cap M)$ with the property that for all $u\in U$
  it holds that $\operatorname{rank}(\phi'(u))=m$ we have that
  \begin{eqn}
    \llabel{eq1}
    f\circ\phi\in C^r(U,S)
    \dott
  \end{eqn}
  \cfout[.]
\end{athm}

\begin{athm}{lemma}{lem:angular-gradient-existence}
  [Existence and uniqueness of the angular gradient]
  Let $n\in\N$, $\mathfrak r\in(0,\infty)$, let
  $M\subseteq\R^n$ be a \submanifold{1}{n}, let
  $U\subseteq(0,\mathfrak r)\times M$ be open in
  $(0,\mathfrak r)\times M$, and let $Z\in\COne{U}{\R}$\dott \cfload[.]
  Then for all $(r,x)\in U$ there exists a unique
  $v\in T_xM$\cfadd{definition:tangent-space} such that for all
  $\varepsilon\in(0,\infty)$,
  $\gamma\in C^1((-\varepsilon,\varepsilon),M)$,
  $\sigma\in C^1((-\varepsilon,\varepsilon),\R)$ with $\gamma(0)=x$,
  $\midcup_{t\in(-\varepsilon,\varepsilon)}\{(r,\gamma(t))\}\subseteq U$,
  and $\forall\,t\in(-\varepsilon,\varepsilon)\colon\sigma(t)=Z(r,\gamma(t))$
  it holds that
  \begin{eqn}
    \llabel{concl1}
    \langle v,\gamma'(0)\rangle=\sigma'(0)
    \dott
  \end{eqn}
  \cfout[.]
\end{athm}
\begin{aproof}
  We first prove that for all $(r,x)\in U$ it holds that there exists $v\in T_xM$
  such that for all $\varepsilon\in(0,\infty)$,
  $\gamma\in C^1((-\varepsilon,\varepsilon),M)$,
  $\sigma\in C^1((-\varepsilon,\varepsilon),\R)$ with $\gamma(0)=x$,
  $\{r\}\times\gamma((-\varepsilon,\varepsilon))\subseteq U$, and
  $\forall\,t\in(-\varepsilon,\varepsilon)\colon\sigma(t)=Z(r,\gamma(t))$
  that
  \begin{eqn}
    \llabel{eq1}
    \langle v,\gamma'(0)\rangle=\sigma'(0)
    \dott
  \end{eqn}
  \startnewargseq
  To prove \lref{eq1}, let $(r,x)\in U$\dott
  \startnewargseq
  In the following let $m\in\N_0$ satisfy that $M$ is an
  \dsubmanifold{1}{n}{m}, and let
  $Q\colon\R^n\to\R^m$ satisfy for all $y=(y_1,\ldots,y_n)\in\R^n$ that
  \begin{eqn}
    \llabel{eq3}
    Q(y)=(y_1,\ldots,y_m)
    \dott
  \end{eqn}
  \startnewargseq
  \argument{
    \cref{lem:embedded-manifold-local-straightening};
  }{
    that there exist open $O,W\subseteq\R^n$ and a
    $C^1$-diffeomorphism $\Phi\colon O\to W$ which satisfy that 
    \begin{eqn}
      \llabel{eq4}
      x\in O,\qquad \Phi(x)=0,
      \qquad\text{and}\qquad
      \Phi(M\cap O)=W\cap(\R^m\times\{0\}^{n-m})
      \dott
    \end{eqn}
  }
  In the following let
  \begin{eqn}
    \llabel{eq5}
    D=\{a\in\R^m\colon(a,0)\in W,\ (r,\Phi^{-1}(a,0))\in U\}
  \end{eqn}
  and let $\phi\colon D\to M$ and $f\colon D\to\R$ satisfy for all $a\in D$ that
  \begin{eqn}
    \llabel{eq6}
    \phi(a)=\Phi^{-1}(a,0)
    \qquad\text{and}\qquad
    f(a)=Z(r,\phi(a))
    \dott
  \end{eqn}
  \startnewargseq
  \argument{
    \lref{eq3};
    \lref{eq4};
    \lref{eq5};
    \lref{eq6};
    the fact that $U$ is open in $(0,\mathfrak r)\times M$;
    the chain rule;
  }{
    that $D$ is open, that $0\in D$, that $\phi\in C^1(D,M)$,
    that $\phi(D)$ is open in $M$, that
    \begin{eqn}
      \llabel{eq7}
      \phi(0)=x
      \qquad\text{and}\qquad
      \{r\}\times\phi(D)\subseteq U,
    \end{eqn}
    and that for all $a\in D$, $y\in\phi(D)$ it holds that
    \begin{eqn}
      \llabel{eq8}
      Q(\Phi(\phi(a)))=a
      \qquad\text{and}\qquad
      \phi(Q(\Phi(y)))=y
      \dott
    \end{eqn}
  }
  \argument{
    \lref{eq8};
    \lref{eq3};
    \lref{eq4};
    \lref{eq6};
    the fact that $D$ is open and $\phi(D)$ is open in $M$;
    the chain rule;
  }{
    that for all $a\in D$ it holds that
    \begin{eqn}
      \llabel{eq9}
      \phi'(a)\colon\R^m\to T_{\phi(a)}M
      \text{ is a linear isomorphism}
      \dott
    \end{eqn}
  }
  In the following let
  \begin{eqn}
    \llabel{eq10}
    A=\{(s,a)\in(0,\mathfrak r)\times D\colon(s,\phi(a))\in U\}
  \end{eqn}
  and let $\Psi\colon A\to U$ satisfy for all $(s,a)\in A$ that
  \begin{eqn}
    \llabel{eq11}
    \Psi(s,a)=(s,\phi(a))
    \dott
  \end{eqn}
  \startnewargseq
  \argument{
    \lref{eq7};
    \lref{eq9};
    \lref{eq10};
    \lref{eq11};
    the fact that $U$ is open in $(0,\mathfrak r)\times M$;
  }{
    that $A$ is open, that $\{r\}\times D\subseteq A$,
    that $\Psi\in C^1(A,U)$, and that for all $(s,a)\in A$ it holds that
    \begin{eqn}
      \llabel{eq12}
      \operatorname{rank}(\Psi'(s,a))=m+1
      \dott
    \end{eqn}
  }
  \argument{
    \lref{eq12};
    \lref{eq6};
    \lref{eq11};
    the fact that $Z\in\COne{U}{\R}$;
    the chain rule;
  }{
    that $Z\circ\Psi\in C^1(A,\R)$ and that $f\in C^1(D,\R)$\dott
  }
  \argument{
    \lref{eq7};
    \lref{eq9};
    the fact that $f\in C^1(D,\R)$;
    the fact that $f'(0)\circ(\phi'(0))^{-1}$ is a linear functional on $T_xM$;
    the finite-dimensional Riesz representation theorem;
  }{
    that there exists a unique $v\in T_xM$ which satisfies for all $a\in\R^m$ that
    \begin{eqn}
      \llabel{eq13}
      \langle v,\phi'(0)a\rangle=f'(0)a
      \dott
    \end{eqn}
  }
  We claim that for all $\varepsilon\in(0,\infty)$,
  $\gamma\in C^1((-\varepsilon,\varepsilon),M)$,
  $\sigma\in C^1((-\varepsilon,\varepsilon),\R)$ with $\gamma(0)=x$,
  $\{r\}\times\gamma((-\varepsilon,\varepsilon))\subseteq U$, and
  $\forall\,t\in(-\varepsilon,\varepsilon)\colon\sigma(t)=Z(r,\gamma(t))$
  it holds that
  \begin{eqn}
    \llabel{eq13.1}
    \langle v,\gamma'(0)\rangle=\sigma'(0)
    \dott
  \end{eqn}
  \startnewargseq
  To prove \lref{eq13.1}, let $\varepsilon\in(0,\infty)$,
  $\gamma\in C^1((-\varepsilon,\varepsilon),M)$,
  $\sigma\in C^1((-\varepsilon,\varepsilon),\R)$ satisfy for all
  $t\in(-\varepsilon,\varepsilon)$ that
  \begin{eqn}
    \llabel{eq2}
    \gamma(0)=x,
    \qquad
    (r,\gamma(t))\in U,
    \qquad\text{and}\qquad
    \sigma(t)=Z(r,\gamma(t))
    \dott
  \end{eqn}
  \startnewargseq
  \argument{\lref{eq2};
    \lref{eq7};
    the fact that $\phi(D)$ is open in $M$;
    the fact that $\gamma$ is continuous;
  }{
    that there exists $\eta\in(0,\varepsilon)$ which satisfies that
    \begin{eqn}
      \llabel{eq14}
      \gamma((-\eta,\eta))\subseteq\phi(D)
      \dott
    \end{eqn}
  }
  In the following let $\alpha\colon(-\eta,\eta)\to D$ satisfy for all
  $t\in(-\eta,\eta)$ that
  \begin{eqn}
    \llabel{eq15}
    \alpha(t)=Q(\Phi(\gamma(t)))
    \dott
  \end{eqn}
  \startnewargseq
  \argument{
    \lref{eq2};
    \lref{eq3};
    \lref{eq4};
    \lref{eq6};
    \lref{eq7};
    \lref{eq8};
    \lref{eq14};
    \lref{eq15};
    the chain rule;
  }{
    that $\alpha\in C^1((-\eta,\eta),D)$, that $\alpha(0)=0$,
    and that for all $t\in(-\eta,\eta)$ it holds that
    \begin{eqn}
      \llabel{eq16}
      \gamma(t)=\phi(\alpha(t))
      \qquad\text{and}\qquad
      \sigma(t)=f(\alpha(t))
      \dott
    \end{eqn}
  }
  \argument{
    \lref{eq16};
    the chain rule;
  }{
    that
    \begin{eqn}
      \llabel{eq17}
      \gamma'(0)=\phi'(0)\alpha'(0)
      \qquad\text{and}\qquad
      \sigma'(0)=f'(0)\alpha'(0)
      \dott
    \end{eqn}
  }
  \argument{
    \lref{eq17};
    \lref{eq13};
    the fact that the choices of $\varepsilon$, $\gamma$, and $\sigma$
    were arbitrary;
  }[verbs=ep]{
    \lref{eq13.1}\dott
  }
  \startnewargseq
  \argument{
    \lref{eq13};
    \lref{eq13.1};
    the fact that the choice of $(r,x)$ was arbitrary;
  }[verbs=ep]{
    \lref{eq1}\dott
  }

   Next we prove that for all $(\mathbf r,\mathbf x)\in U$,
  $w,\mathbf w\in T_{\mathbf x}M$ with the property that for all
  $\varepsilon\in(0,\infty)$,
  $\gamma\in C^1((-\varepsilon,\varepsilon),M)$,
  $\sigma\in C^1((-\varepsilon,\varepsilon),\R)$ with $\gamma(0)=\mathbf x$,
  $\{\mathbf r\}\times\gamma((-\varepsilon,\varepsilon))\subseteq U$, and
  $\forall\,t\in(-\varepsilon,\varepsilon)\colon\sigma(t)=Z(\mathbf r,\gamma(t))$
  it holds that $\langle w,\gamma'(0)\rangle
    =\langle\mathbf w,\gamma'(0)\rangle=\sigma'(0)$ we have that
  \begin{eqn}
    \llabel{eq18}
    w=\mathbf w
    \dott
  \end{eqn}
  \startnewargseq
  To prove \lref{eq18}, let $(\mathbf r,\mathbf x)\in U$,
  $w,\mathbf w\in T_{\mathbf x}M$
  satisfy for all $\varepsilon\in(0,\infty)$,
  $\gamma\in C^1((-\varepsilon,\varepsilon),M)$,
  $\sigma\in C^1((-\varepsilon,\varepsilon),\R)$ with $\gamma(0)=\mathbf x$,
  $\{\mathbf r\}\times\gamma((-\varepsilon,\varepsilon))\subseteq U$, and
  $\forall\,t\in(-\varepsilon,\varepsilon)\colon\sigma(t)=Z(\mathbf r,\gamma(t))$ that
  \begin{eqn}
    \llabel{eq19}
    \langle w,\gamma'(0)\rangle
      =\langle\mathbf w,\gamma'(0)\rangle=\sigma'(0)
    \dott
  \end{eqn}
  \startnewargseq
  In the following let $m\in\N_0$ satisfy that $M$ is an
  \dsubmanifold{1}{n}{m}\dott
  \startnewargseq
  \argument{
    \cref{lem:embedded-manifold-local-straightening};
  }{
    that there exist open $O,W\subseteq\R^n$ and a
    $C^1$-diffeomorphism $\Phi\colon O\to W$ which satisfy that
    \begin{eqn}
      \llabel{eq19.1}
      \mathbf x\in O,\qquad \Phi(\mathbf x)=0,
      \qquad\text{and}\qquad
      \Phi(M\cap O)=W\cap(\R^m\times\{0\}^{n-m})
      \dott
    \end{eqn}
  }
  In the following let
  \begin{eqn}
    \llabel{eq19.2}
    A=\{(s,a)\in(0,\mathfrak r)\times\R^m\colon
      (a,0)\in W,\ (s,\Phi^{-1}(a,0))\in U\}\quad\text{and}\quad
    D=\{a\in\R^m\colon(\mathbf r,a)\in A\},
  \end{eqn}
  let $\Psi\colon A\to U$ satisfy for all $(s,a)\in A$ that
  \begin{eqn}
    \llabel{eq19.3}
    \Psi(s,a)=(s,\Phi^{-1}(a,0)),
  \end{eqn}
  and let $\phi\colon D\to M$ and $f\colon D\to\R$ satisfy for all $a\in D$ that
  \begin{eqn}
    \llabel{eq19.4}
    \phi(a)=\Phi^{-1}(a,0)
    \qquad\text{and}\qquad
    f(a)=Z(\mathbf r,\phi(a))
    \dott
  \end{eqn}
  \startnewargseq
  \argument{
    \lref{eq19.1};
    \lref{eq19.2};
    \lref{eq19.3};
    \lref{eq19.4};
    the fact that $U$ is open in $(0,\mathfrak r)\times M$;
    the chain rule;
  }{
    that $A$ and $D$ are open, that $0\in D$,
    that $\Psi\in C^1(A,U)$, that $\phi\in C^1(D,M)$,
    that $\phi'(0)\colon\R^m\to T_{\mathbf x}M$ is a linear isomorphism, that
    \begin{eqn}
      \llabel{eq19.5}
      \phi(0)=\mathbf x
      \qquad\text{and}\qquad
      \{\mathbf r\}\times\phi(D)\subseteq U,
    \end{eqn}
    and that for all $(s,a)\in A$ it holds that
    \begin{eqn}
      \llabel{eq19.7}
      \operatorname{rank}(\Psi'(s,a))=m+1
      \dott
    \end{eqn}
  }
  \argument{
    \lref{eq19.7};
    the fact that $Z\in\COne{U}{\R}$;
  }{
    that $Z\circ\Psi\in C^1(A,\R)$\dott
  }
  \argument{
    \lref{eq19.2};
    \lref{eq19.3};
    \lref{eq19.4};
  }{
    that for all $a\in D$ it holds that
    \begin{eqn}
      \llabel{eq19.8}
      (\mathbf r,a)\in A
      \qquad\text{and}\qquad
      f(a)=Z(\mathbf r,\phi(a))=(Z\circ\Psi)(\mathbf r,a)
      \dott
    \end{eqn}
  }
  \argument{
    \lref{eq19.8};
    the fact that $Z\circ\Psi\in C^1(A,\R)$;
    the fact that $(D\ni a\mapsto(\mathbf r,a)\in A)\in C^1(D,A)$;
    the chain rule;
  }{
    that $f\in C^1(D,\R)$\dott
  }
  We claim that for all $a\in\R^m$ it holds that
  \begin{eqn}
    \llabel{eq20}
    \langle w,\phi'(0)a\rangle
      =\langle\mathbf w,\phi'(0)a\rangle=f'(0)a
    \dott
  \end{eqn}
  \startnewargseq
  To prove \lref{eq20}, let $a\in\R^m$\dott
  \startnewargseq
  \argument{
    the fact that $D$ is open and $0\in D$;
  }{
    that there exists $\tau\in(0,\infty)$ which satisfies for all
    $t\in(-\tau,\tau)$ that
    \begin{eqn}
      \llabel{eq21}
      ta\in D
      \dott
    \end{eqn}
  }
  In the following let $\beta\colon(-\tau,\tau)\to M$ and
  $\chi\colon(-\tau,\tau)\to\R$ satisfy for all $t\in(-\tau,\tau)$ that
  \begin{eqn}
    \llabel{eq22}
    \beta(t)=\phi(ta)
    \qquad\text{and}\qquad
    \chi(t)=f(ta)
    \dott
  \end{eqn}
  \startnewargseq
  \argument{
    \lref{eq19.4};
    \lref{eq19.5};
    \lref{eq21};
    \lref{eq22};
    the fact that $\phi\in C^1(D,M)$ and $f\in C^1(D,\R)$;
    the chain rule;
  }{
    that $\beta\in C^1((-\tau,\tau),M)$,
    that $\chi\in C^1((-\tau,\tau),\R)$,
    that $\beta(0)=\mathbf x$, that for all $t\in(-\tau,\tau)$ it holds that
    \begin{eqn}
      \llabel{eq23}
      (\mathbf r,\beta(t))\in U,
      \qquad
      \chi(t)=Z(\mathbf r,\beta(t)),
      \qquad
      \beta'(0)=\phi'(0)a
      \qquad\text{and}\qquad
      \chi'(0)=f'(0)a
      \dott
    \end{eqn}
  }
  \argument{
    \lref{eq19};
    \lref{eq23};
    the fact that $a\in\R^m$ was arbitrary;
  }[verbs=ep]{
    \lref{eq20}\dott
  }
  \startnewargseq
  \argument{
    \lref{eq19.5};
    \lref{eq20};
    the fact that the choices of $(\mathbf r,\mathbf x)$, $w$, and $\mathbf w$
    were arbitrary;
  }[verbs=ep]{
    \lref{eq18}\dott
  }
  \startnewargseq
  \argument{
    \lref{eq1};
    \lref{eq18};
  }[verbs=ep]{
    \lref{concl1}\dott
  }
\end{aproof}

\begin{athm}{definition}{definition:angular-gradient}[Gradients with respect to embedded submanifolds]
  Let $n\in\N$, $\mathfrak{r}_0\in(0,\infty)$, let
  $M\subseteq\R^n$ be a definable \submanifold{1}{n},
  let $U\subseteq(0,\mathfrak{r}_0)\times M$ be
  definable and open in $(0,\mathfrak{r}_0)\times M$, and let
  $Z\in\COne{U}{\R}$ be definable\dott \cfload[.]
  Then we denote by $\nabla_MZ\colon U\to\R^n$ the function which satisfies
  for all $(r,x)\in U$ that
  \begin{enumerate}[label=(\roman*)]
    \item it holds that $(\nabla_MZ)(r,x)\in T_xM$\cfadd{definition:tangent-space} and
    \item it holds for all $\varepsilon\in(0,\infty)$,
    $\gamma\in C^1((-\varepsilon,\varepsilon),M)$,
    $\sigma\in C^1((-\varepsilon,\varepsilon),\R)$ with $\gamma(0)=x$,\\
    $\midcup_{t\in(-\varepsilon,\varepsilon)}\{(r,\gamma(t))\}\subseteq U$,
    and $\forall\,t\in(-\varepsilon,\varepsilon)\colon\sigma(t)=Z(r,\gamma(t))$ that
    \begin{eqn}
      \langle(\nabla_MZ)(r,x),\gamma'(0)\rangle
      =\sigma'(0)
      \dott
    \end{eqn}
  \end{enumerate}
  \cfout[.]
\end{athm}
\begin{athm}{lemma}{lem:radial-derivative-well-defined}
  [Continuity and definability of the partial derivative with respect to the radial variable]
  Let $n\in\N$, $\mathfrak{r}\in(0,\infty)$, let $M\subseteq\R^n$
  be a definable\cfadd{definition:definable-set} \submanifold{1}{n}, let
  $U\subseteq(0,\mathfrak{r})\times M$ be definable and open in
  $(0,\mathfrak{r})\times M$, and let
  $f=(f(r,x))_{(r,x)\in U}\in \COne{U}{\R}$ be definable%
  \cfadd{definition:definable-function}\cfadd{definition:c1-on-manifold}\dott
  \cfload[.]
  Then
  \begin{enumerate}[label=(\roman*)]
    \item \llabel{item1} it holds for all $x\in M$ that
    $\{s\in(0,\mathfrak{r})\colon(s,x)\in U\}\ni r\mapsto f(r,x)\in\R$
    is differentiable and
    \item \llabel{item2} it holds that
    $U\ni(r,x)\mapsto(\partial_r f)(r,x)\in\R$
    is continuous and definable\dott
  \end{enumerate}
  \cfout[.]
\end{athm}
\begin{aproof}
  Throughout this proof let $m\in\N_0$ satisfy that $M$ is an
  \dsubmanifold{1}{n}{m}\cfadd{definition:embedded-manifold}\dott \cfload[.]
  \startnewargseq
  \argument{
    the fact that $M$ is an $m$-dimensional $C^1$ submanifold of $\R^n$;
  }{
    that for every $p\in M$ there exist an open
    $W_p\subseteq\R^m$, a set $O_p\subseteq M$ open in $M$,
    and a homeomorphism $\phi_p\in C^1(W_p,O_p)$ which
    satisfy that $p\in O_p$ and that for all $u\in W_p$ it holds that
    $\operatorname{rank}(\phi_p)'(u)=m$%
    \cfadd{definition:embedded-manifold}\dott
  }
  In the following for every $p\in M$ let
  $A_p=\{(r,u)\in(0,\mathfrak{r})\times W_p\colon(r,\phi_p(u))\in U\}$
  and let
  $\Phi_p\colon A_p\to U\cap((0,\mathfrak{r})\times O_p)$
  and
  $g_p=(g_p(r,u))_{(r,u)\in A_p}\colon A_p\to\R$ satisfy for all
  $(r,u)\in A_p$ that
  \begin{eqn}
    \llabel{eq1}
    \Phi_p(r,u)=(r,\phi_p(u))
    \qquad\text{and}\qquad
    g_p(r,u)=f(r,\phi_p(u))
    \dott
  \end{eqn}
  \startnewargseq
  \argument{
    \lref{eq1};
    the assumption that $U$ is open in $(0,\mathfrak{r})\times M$;
    the fact that for all $p\in M$ it holds that $A_p$ is open;
    the fact that for all $p\in M$ it holds that
    $\Phi_p\in C^1(A_p,U\cap((0,\mathfrak{r})\times O_p))$
    is a homeomorphism;
    the fact that for all $p\in M$, $(r,u)\in A_p$
    it holds that
    $\operatorname{rank}(\Phi_p)'(r,u)=m+1$;
    the fact that $f\in \COne{U}{\R}$;
  }{
    that for all $p\in M$ it holds that
    $g_p\in C^1(A_p,\R)$\dott
    \llabel{arg1}
  }
  \argument{
    \lref{arg1};
    \lref{eq1};
    the fact that for all $p\in M$ it holds that $p\in O_p$;
  }[verbs=ep]{
    \lref{item1}\dott
  }
  \startnewargseq
  \argument{
    \lref{eq1};
    \lref{item1};
    the fact that for all $p\in M$ it holds that
    $g_p\in C^1(A_p,\R)$;
  }{
    that for all $p\in M$, $(r,u)\in A_p$
    it holds that
    \begin{eqn}[a]
      \llabel{eq2}
      (\partial_r f)(r,\phi_p(u))
      =\lim_{t\to0}
        \frac{f(r+t,\phi_p(u))-f(r,\phi_p(u))}{t}=\lim_{t\to0}
        \frac{g_p(r+t,u)-g_p(r,u)}{t}
        =(\partial_r g_p)(r,u)
      \dott
    \end{eqn}
  }
  \argument{
    \lref{eq2};
    the fact that for all $p\in M$ it holds that
    $\phi_p$ is a homeomorphism;
  }{
    that for all $p\in M$, $(r,y)\in U\cap((0,\mathfrak{r})\times O_p)$
    it holds that
    \begin{eqn}
      \llabel{eq3}
      (\partial_r f)(r,y)
      =(\partial_r g_p)(r,\phi_p^{-1}(y))
      \dott
    \end{eqn}
  }
  \argument{
    \lref{eq3};
    the fact that for all $p\in M$ it holds that
    $\partial_r g_p$ is continuous;
    the fact that for all $p\in M$ it holds that
    $U\cap((0,\mathfrak{r})\times O_p)\ni(r,y)\mapsto(r,\phi_p^{-1}(y))
      \in A_p$
    is continuous;
    the fact that $(U\cap((0,\mathfrak{r})\times O_p))_{p\in M}$ is an open
    cover of $U$;
  }{
    that $\partial_r f$ is continuous\dott
  }
  \argument{\lref{item1};
  }{
    that
    \begin{equation}
      \llabel{eq4}
      \begin{aligned}
      &\bigl\{(a,r,y)\in\R\times U\colon
        a=(\partial_r f)(r,y)\bigr\}\\
      &=\braces*{(a,r,y)\in\R\times U\colon\brackets*{
        \begin{aligned}
          &\forall\,\varepsilon\in(0,\infty)\ \exists\,\delta\in(0,\infty)\colon\bigl(
            \forall\,s\in(0,\mathfrak{r}),\, (s,y)\in U,\, 0<\abs{s-r}<\delta\colon\\
          &\qquad\abs{f(s,y)-f(r,y)-a(s-r)}
            <\varepsilon\abs{s-r}\bigr)
        \end{aligned}
        }}
      \dott
      \end{aligned}
    \end{equation}
  }
  \argument{
    \lref{eq4};
    the fact that $f$ is definable;
    the fact that $U$ is definable;
    \cite[Theorem~1.13]{Coste};
  }{
    that the graph of $\partial_r f$ is definable\dott
    \llabel{arg3}
  }
  \argument{
    \lref{arg3};
    the fact that $\partial_r f$ is continuous;
  }[verbs=ep]{
    \lref{item2}\dott
  }
\end{aproof}

\begin{athm}{lemma}{lem:angular-gradient-well-defined}
  [Continuity and definability of the gradient with respect to the angular variable]
  Let $n\in\N$, $\mathfrak r\in(0,\infty)$, let $M\subseteq\R^n$ be a
  definable \submanifold{1}{n}, let
  $U\subseteq(0,\mathfrak r)\times M$ be definable and open in
  $(0,\mathfrak r)\times M$, and let $Z\in\COne{U}{\R}$ be definable\dott
  \cfload[.]
  Then $\nabla_MZ\colon U\to\R^n$ is continuous and
  definable\cfadd{definition:angular-gradient}\dott \cfout[.]
\end{athm}
\begin{aproof}
  Throughout this proof let $m\in\N_0$ satisfy that $M$ is an
  \dsubmanifold{1}{n}{m}\dott\cfload[.]
  \startnewargseq
  In our proof we distinguish between the case $M=\emptyset$, the case
  ($M\neq\emptyset$ and $m=0$), and the case
  ($M\neq\emptyset$ and $m\geq1$)\dott
  \startnewargseq
  We first prove that $\nabla_MZ$ is continuous and definable in the case that
  \begin{eqn}
    \llabel{case1}
    M=\emptyset\dott
  \end{eqn}
  \startnewargseq
  \argument{
    \lref{case1};
    the fact that $U\subseteq(0,\mathfrak r)\times M$;
  }{
    that $U=\emptyset$\dott
    \llabel{arg:empty-domain}
  }
  \argument{
    \lref{arg:empty-domain};
  }[verbs=ep]{
    that $\nabla_MZ$ is continuous and definable in the case
    $M=\emptyset$\dott
  }
  We now prove that $\nabla_MZ$ is continuous and definable in the case that
  \begin{eqn}
    \llabel{case2}
    M\neq\emptyset\qquad\text{and}\qquad m=0\dott
  \end{eqn}
  \startnewargseq
  \argument{
    \lref{case2};
    the fact that $M$ is an \dsubmanifold{1}{n}{m};
    the fact that for all $(r,x)\in U$ it holds that
    $(\nabla_MZ)(r,x)\in T_xM$;
  }{
    that for all $(r,x)\in U$ it holds that $T_xM=\{0\}$ and
    $(\nabla_MZ)(r,x)=0$\dott
    \llabel{arg:zero-gradient}
  }
  \argument{
    \lref{arg:zero-gradient};
    the fact that $U$ is definable;
  }[verbs=ep]{
    that $\nabla_MZ$ is continuous and definable in the case
    $M\neq\emptyset$ and $m=0$\dott
  }
  We now prove that $\nabla_MZ$ is continuous and definable in the case that
  \begin{eqn}
    \llabel{case3}
    M\neq\emptyset\qquad\text{and}\qquad m\geq1\dott
  \end{eqn}
  \startnewargseq
  We claim that for all $x\in M$ there exist a definable set $V\subseteq M$
  open in $M$ with $x\in V$ and a coordinate projection
  $\pi\colon\R^n\to\R^m$ such that $\pi(V)$ is definable and open and that
  \begin{eqn}
    \llabel{eq:charts}
    \pi|_V\colon V\to\pi(V)\text{ is a definable }C^1\text{-diffeomorphism}
    \dott
  \end{eqn}
  \startnewargseq
  To prove \lref{eq:charts}, let $x\in M$\dott\cfload[.]
  \startnewargseq
  \argument{
    \cref{lem:embedded-manifold-local-straightening};
    the inverse function theorem;
  }{
    that there exist a set $W\subseteq M$ open in $M$ with $x\in W$
    and a coordinate projection $\pi\colon\R^n\to\R^m$ which satisfy that
    $\pi|_W$ is a $C^1$-diffeomorphism onto an open subset of $\R^m$\dott
  }
  \startnewargseq
  \argument{
    the fact that $W$ is open in $M$ and $x\in W$;
  }{
    that there exists $\delta\in(0,\infty)$ which satisfies that
    \begin{eqn}
      \llabel{eq:chart-radius}
      \{y\in M\colon|y-x|<\delta\}\subseteq W
      \dott
    \end{eqn}
  }
  In the following let
  \begin{eqn}
    \llabel{eq:chart-domain}
    V=\{y\in M\colon|y-x|<\delta\}
    \dott
  \end{eqn}
  \startnewargseq
  \argument{
    \lref{eq:chart-radius};
    \lref{eq:chart-domain};
    the fact that $M$ is definable;
  }{
    that $V$ is open in $M$,
    that $\pi(V)$ is definable and open, and that
    $\pi|_V\colon V\to\pi(V)$ is a $C^1$-diffeomorphism\dott
  }
  \argument{
    the fact that $V$ is definable;
    the fact that $\pi|_V\colon V\to\pi(V)$ is bijective;
  }{
    that
    \begin{eqn}
      \operatorname{graph}((\pi|_V)^{-1})
        =\{(\pi(y),y)\colon y\in V\}
    \end{eqn}
    is definable and that $(\pi|_V)^{-1}$ is definable\dott
  }
  \argument{
    the fact that $V$ is definable and open in $M$ and $x\in V$;
    the fact that $\pi$ is a coordinate projection and $\pi(V)$ is definable and open;
    the fact that $(\pi|_V)^{-1}\colon\pi(V)\to V$ is a definable
    $C^1$-diffeomorphism;
    the fact that $x\in M$ was arbitrary;
  }[verbs=ep]{
    \lref{eq:charts}\dott
  }
\startnewargseq

  \argument{
    \lref{eq:charts};
    the fact that $M$ is definable;
    the existence of a finite definable $C^1$-atlas for every definable
    $C^1$-submanifold (see \cite[Section~5.2]{Fischer2010Peano});
  }{
    that there exist $N\in\N$, definable sets $\mathcal O_j\subseteq M$
    open in $M$, definable open $\mathcal D_j\subseteq\R^m$, and definable
    $C^1$-diffeomorphisms $\phi_j\colon\mathcal D_j\to\mathcal O_j$,
    $j\in\{1,\ldots,N\}$, which satisfy that
    \begin{eqn}
      \llabel{eq:cover}
      M=\midcup_{j=1}^N\mathcal O_j
    \end{eqn}
    and which satisfy for all $j\in\{1,2,\ldots,N\}$, $u\in\mathcal D_j$ that
    $(\phi_j)'(u)\colon\R^m\to T_{\phi_j(u)}M$ is a linear isomorphism\dott
  }
  In the following for every $j\in\{1,2,\ldots,N\}$ let
  \begin{eqn}
    \mathcal U_j=\{(r,u)\in(0,\mathfrak r)\times\mathcal D_j\colon
      (r,\phi_j(u))\in U\}
  \end{eqn}
  and let $z_j\colon\mathcal U_j\to\R$ satisfy for all
  $(r,u)\in\mathcal U_j$ that $z_j(r,u)=Z(r,\phi_j(u))$\dott
  \startnewargseq
  \argument{
    the fact that $U$ is definable and open in $(0,\mathfrak r)\times M$;
    the fact that for all $j\in\{1,2,\dots,N\}$ it holds that
    $\phi_j$ is definable and continuously differentiable;
    the fact that $Z\in\COne{U}{\R}$ is definable;
  }{
    that for all $j\in\{1,\ldots,N\}$ it holds that $\mathcal U_j$ is
    definable and open and that $z_j\in C^1(\mathcal U_j,\R)$ is definable\dott
  }
  In the following for all $j\in\{1,\ldots,N\}$ let
  $A_j\colon\mathcal D_j\to\R^{n\times m}$ and
  $G_j\colon\mathcal D_j\to\R^{m\times m}$ satisfy for all
  $u\in\mathcal D_j$ that
  \begin{eqn}
    \llabel{eq:chart-matrices}
    A_j(u)=(\phi_j)'(u)
    \qquad\text{and}\qquad
    G_j(u)=(A_j(u))^{\mathsf T}A_j(u)
    \dott
  \end{eqn}
  \startnewargseq
  \argument{
    \lref{eq:chart-matrices};
    the fact that for all $j\in\{1,2,\ldots,N\}$, $u\in\mathcal D_j$
    it holds that $(\phi_j)'(u)\colon\R^m\to T_{\phi_j(u)}M$
    is a linear isomorphism;
  }{
    that for all $j\in\{1,2,\ldots,N\}$, $u\in\mathcal D_j$ it holds that
    $G_j(u)$ is positive-definite and invertible\dott
  }
  In the following for all $j\in\{1,2,\ldots,N\}$ let
  $g_j\colon\mathcal U_j\to\R^n$
  satisfy for all $(r,u)\in\mathcal U_j$ that
  \begin{eqn}
    \llabel{eq:local-gradient}
    g_j(r,u)=A_j(u)(G_j(u))^{-1}(\nabla_u z_j)(r,u)
    \dott
  \end{eqn}
  \startnewargseq
  \argument{
    the fact that for all $j\in\{1,2,\ldots,N\}$ it holds that
    $\phi_j$ and $z_j$ are continuously differentiable and definable;
    \itref{lem:radial-derivative-well-defined}{item2};
    \lref{eq:local-gradient};
  }{
    that for all $j\in\{1,2,\ldots,N\}$ it holds that $g_j$ is continuous and
    definable and that for all $(r,u)\in\mathcal U_j$ it holds that
    $g_j(r,u)\in T_{\phi_j(u)}M$\dott
  }
  We claim that for all $j\in\{1,\ldots,N\}$, $(r,u)\in\mathcal U_j$
  it holds that
  \begin{eqn}
    \llabel{eq:identification}
    (\nabla_MZ)(r,\phi_j(u))=g_j(r,u)
    \dott
  \end{eqn}
  \startnewargseq
  To prove \lref{eq:identification}, let $j\in\{1,2,\ldots,N\}$,
  $(r,u)\in\mathcal U_j$, $\varepsilon\in(0,\infty)$,
  $\gamma\in C^1((-\varepsilon,\varepsilon),M)$, and
  $\sigma\in C^1((-\varepsilon,\varepsilon),\R)$ satisfy for all
  $t\in(-\varepsilon,\varepsilon)$ that
  \begin{eqn}
    \llabel{eq:curve-data}
   \gamma(0)=\phi_j(u),\qquad (r,\gamma(t))\in U,
    \qquad\text{and}\qquad
    \sigma(t)=Z(r,\gamma(t))
    \dott
  \end{eqn}
  \startnewargseq
  \argument{
    \lref{eq:curve-data};
    the fact that $\gamma$ is continuous;
    the fact that $\phi_j(u)\in\mathcal O_j$; the fact that $\mathcal O_j$ is open in $M$;
  }{
    that there exists $\eta\in(0,\varepsilon)$ which satisfies that
    \begin{eqn}
      \llabel{eq:curve-interval}
      \gamma((-\eta,\eta))\subseteq\mathcal O_j
      \dott
    \end{eqn}
  }
  \startnewargseq
  \argument{
    \lref{eq:curve-data};
    \lref{eq:curve-interval};
    the fact that $\gamma\in C^1((-\varepsilon,\varepsilon),M)$;
    the fact that $\phi_j\colon\mathcal D_j\to\mathcal O_j$ is a
    $C^1$-diffeomorphism;
    the fact that for all $(s,v)\in\mathcal U_j$ it holds that
    $z_j(s,v)=Z(s,\phi_j(v))$;
    the chain rule;
  }{
    that there exists $\alpha\in C^1((-\eta,\eta),\mathcal D_j)$ which
    satisfies $\alpha(0)=u$ and satisfies for all $t\in(-\eta,\eta)$ that
    \begin{eqn}
      \llabel{eq:curve-coordinates}
      \gamma(t)=\phi_j(\alpha(t))
      \qquad\text{and}\qquad
      \sigma(t)=z_j(r,\alpha(t))
      \dott
    \end{eqn}
  }
  \startnewargseq
  \argument{
    \lref{eq:chart-matrices};
    \lref{eq:local-gradient};
    \lref{eq:curve-coordinates};
    the chain rule;
  }{
    that
    \begin{eqn}[a]
      \llabel{eq:curve-identity}
      \langle g_j(r,u),\gamma'(0)\rangle
        &=\langle A_j(u)(G_j(u))^{-1}(\nabla_u z_j)(r,u),
          A_j(u)\alpha'(0)\rangle\\
        &=\langle(\nabla_u z_j)(r,u),\alpha'(0)\rangle=\sigma'(0)
      \dott
    \end{eqn}
  }
  \argument{
    \lref{eq:curve-identity};
    the fact that $g_j(r,u)\in T_{\phi_j(u)}M$;
    the fact that the choices of $j$, $(r,u)$, $\varepsilon$, $\gamma$,
    and $\sigma$ were arbitrary;
    \cref{lem:angular-gradient-existence};
  }[verbs=ep]{
    \lref{eq:identification}\dott
  }

  \startnewargseq
  \argument{\lref{eq:cover};
    \lref{eq:identification};
    the fact that for all $j\in\{1,2,\dots,N\}$ it holds that
    $g_j$ and $\phi_j^{-1}$ are continuous;
  }{
    that $\nabla_MZ$ is continuous\dott
  }
  \startnewargseq
  \argument{  \lref{eq:cover};
    \lref{eq:identification};
    the fact that for all $j\in\{1,\ldots,N\}$ it holds that
    $\mathcal U_j$, $\phi_j$, and $g_j$ are definable;
  }{
    that
    \begin{eqn}
      \operatorname{graph}(\nabla_MZ)
        =\midcup_{j=1}^N
          \{((r,\phi_j(u)),g_j(r,u))\colon(r,u)\in\mathcal U_j\}
    \end{eqn}
    is definable\dott
    \llabel{arg:gradient-definable}
  }
  \argument{
    \lref{arg:gradient-definable};
    the fact that $\nabla_MZ$ is continuous;
  }[verbs=ep]{
    that $\nabla_MZ$ is continuous and definable in the case
    $M\neq\emptyset$ and $m\geq1$\dott
  }
\end{aproof}

\subsection{Integrability, smoothness, and definability properties of radial upper bounds}
\label{subsec:normalized-bound}

In the following elementary \cref{lem:support1} we recall
 the well-known fact that for every $r_0\in (0,\infty)$ and every bounded and definable curve $u\colon (0,r_0)\to\R^n$ we have that $u$ has finite length
 near zero. Only for completeness, we include a detailed proof for \cref{lem:support1} here.

\begin{athm}{lemma}{lem:support1}
  Let
  $n,p\in\N$,
  $\mathfrak{r}_0\in(0,\infty)$
  and let $u\colon (0,\mathfrak{r}_0)\to\R^n$ be bounded and definable\dott \cfload[.]
  Then there exists $\mathfrak{r}_1\in(0,\mathfrak{r}_0)$ such that
  \begin{enumerate}[label=(\roman*)]
    \item \llabel{item1} it holds that 
    $u|_{(0,\mathfrak{r}_1)}\in C^p((0,\mathfrak{r}_1),\R^n)$,
    \item \llabel{item2} there exists $v\in\R^n$ such that
    $\limsup_{s\to 0}|u(s)-v|=0$, and
    \item \llabel{item3} it holds that
    $\int_0^{\mathfrak{r}_1}|u'(r)|\dd r<\infty$\dott
  \end{enumerate}
  \cfout[.]
\end{athm}
\begin{aproof}
  Throughout this proof, for every $j\in\{1,\ldots,n\}$, let
  $\pi_j\colon \R^n\to\R$ be the $j$th coordinate projection\dott
  \argument{
    the assumption that $u$ is definable;
    the fact that for all $j\in\{1,\ldots,n\}$ it holds that
    $\pi_j$ is definable;
    \cite[Theorem~2.1 and Corollary~6.5]{Coste};
  }{
    that for all $j\in\{1,\ldots,n\}$ there exists
    $\mathfrak{R}_j\in(0,\mathfrak{r}_0)$ which satisfies that
    $(\pi_j\circ u)|_{(0,\mathfrak{R}_j)}\in C^p((0,\mathfrak{R}_j),\R)$ and that
    $(\pi_j\circ u)|_{(0,\mathfrak{R}_j)}$ is strictly increasing, strictly decreasing, or constant\dott
  }
  \startnewargseq
  In the following let $\rho\in \R$ satisfy
  \begin{eqn}\llabel{def: R}
      \rho=\min\{\mathfrak{R}_1,\mathfrak{R}_2,\dots,\mathfrak{R}_n\}.
  \end{eqn}
  \startnewargseq
  \argument{\lref{def: R};
  }{
    that for all $j\in\{1,\ldots,n\}$ it holds that
    $(\pi_j\circ u)|_{(0,\rho)}\in C^p((0,\rho),\R)$ and that
    $(\pi_j\circ u)|_{(0,\rho)}$ is strictly increasing, strictly decreasing, or constant\dott
  }
In the following let $\mathfrak{r}_1\in(0,\rho)$.
  \argument{
    the fact that for all $j\in\{1,\ldots,n\}$ it holds that
    $(\pi_j\circ u)|_{(0,\mathfrak{r}_1)}\in C^p((0,\mathfrak{r}_1),\R)$;
  }{
    that \begin{eqn}\llabel{eq1}
        u|_{(0,\mathfrak{r}_1)}\in C^p((0,\mathfrak{r}_1),\R^n).
    \end{eqn}
  }
  \argument{
    \lref{eq1};
  }[verbs=ep]{
    \lref{item1}\dott
  }
  \startnewargseq
  \argument{
    the fact that $u$ is bounded;
    the fact that for all $j\in\{1,\ldots,n\}$ it holds that
    $(\pi_j\circ u)|_{(0,\mathfrak{r}_1)}$ is strictly increasing, strictly decreasing, or constant;
  }{
    that for all $j\in\{1,\ldots,n\}$ there exists $v_j\in\R$ which satisfies that
    \begin{eqn}\llabel{def:uj}
       \lim_{r\to 0}(\pi_j\circ u)(r)=v_j. 
    \end{eqn}
  }
  \startnewargseq
  In the following let
  $v=(v_1,\ldots,v_n)\in\R^n$.
  \argument{
   \lref{def:uj}
  }{
    that
    \begin{eqn}
      \llabel{eq2}
      \lim_{r\to 0}u(r)=v
      \dott
    \end{eqn}
  }
  \argument{
    \lref{eq2};
  }[verbs=ep]{
    \lref{item2}\dott
  }
  \startnewargseq
  \argument{
    the fact that $p\in\N$;
    the fact that for all $j\in\{1,\ldots,n\}$ it holds that
    $(\pi_j\circ u)|_{(0,\mathfrak{r}_1)}\in C^1((0,\mathfrak{r}_1),\R)$ and that
    $(\pi_j\circ u)|_{(0,\mathfrak{r}_1)}$ is strictly increasing, strictly decreasing, or constant;
    the fundamental theorem of calculus;
  }{
    that for all $j\in\{1,\ldots,n\}$, $a,b\in(0,\mathfrak{r}_1)$ with $a<b$ it holds that
    \begin{eqn}
      \llabel{eq3}
        \int_a^b|(\pi_j\circ u)'(r)|\dd r
      =
        |(\pi_j\circ u)(b)-(\pi_j\circ u)(a)|
      \dott
    \end{eqn}
  }
  \argument{
    \lref{eq3};
    the fact that for all $j\in\{1,\ldots,n\}$ it holds that
    $\pi_j\circ u$ is bounded on $(0,\mathfrak{r}_1)$;
    the triangle inequality;
  }{
    that for all $j\in\{1,\ldots,n\}$, $a,b\in(0,\mathfrak{r}_1)$ with $a<b$ it holds that
    \begin{eqn}
      \llabel{eq4}
        \int_a^b|(\pi_j\circ u)'(r)|\dd r
      \leq
        2\sup_{s\in(0,\mathfrak{r}_1)}|(\pi_j\circ u)(s)|
      <
        \infty
      \dott
    \end{eqn}
  }
  \argument{
    \lref{eq4};
    the monotone convergence theorem;
  }{
    that for all $j\in\{1,\ldots,n\}$ it holds that
    \begin{eqn}
      \llabel{eq5}
      \int_0^{\mathfrak{r}_1}|(\pi_j\circ u)'(r)|\dd r<\infty
      \dott
    \end{eqn}
  }
  \argument{
    \lref{eq5};
    the fact that for all $r\in(0,\mathfrak{r}_1)$ it holds that
    $|u'(r)|\leq\sum_{j=1}^n|(\pi_j\circ u)'(r)|$;
  }{
    that
    \begin{eqn}
      \llabel{eq6}
        \int_0^{\mathfrak{r}_1}|u'(r)|\dd r
      \leq
        \sum_{j=1}^n\int_0^{\mathfrak{r}_1}|(\pi_j\circ u)'(r)|\dd r
      <
        \infty
      \dott
    \end{eqn}
  }
  \argument{
    \lref{eq6};
  }[verbs=ep]{
    \lref{item3}\dott
  }
\end{aproof}

\begin{athm}{lemma}{lem:normalized}
  Let $n\in\N$, $\mathfrak{r}_0\in(0,\infty)$, let
  $M\subseteq\R^n$ be a bounded and definable \submanifold{1}{n}, let
  $U\subseteq(0,\mathfrak{r}_0)\times M$ be definable,
  assume that $U$ is open in $(0,\mathfrak{r}_0)\times M$,
  let $\Omega\subseteq U$ be definable, let
  $Z=(Z(r,x))_{(r,x)\in U}\in\COne{U}{\R}$ be definable,
  assume $\sup_{v\in\Omega}|Z(v)|<\infty$,
  for every $r\in(0,\mathfrak{r}_0)$ let 
    $\mathcal{F}_r=\{x\in M\colon (r,x)\in\Omega\}$,
  and let $\eta\colon (0,\mathfrak{r}_0)\to[0,\infty]$ satisfy for all $r\in (0,\mathfrak{r}_0)$ that
  \begin{equation}\llabel{eq:eta}
    \eta(r)=\sup\Biggl(\{0\}\cup\biggl(\,
      \mathop{\textstyle\bigcup}\limits_{x\in\mathcal{F}_r}\biggl\{
        \frac{|(\partial_r Z)(r,x)|}{1+|(\nabla_MZ)(r,x)|}
      \biggr\}\biggr)\Biggr)\cfadd{definition:angular-gradient}\dott
  \end{equation}
  \cfload[.]
  Then there exists $\mathfrak{r}_1\in(0,\mathfrak{r}_0)$ such that
  \begin{enumerate}[label=(\roman*)]
    \item \llabel{item1} it holds that $\eta((0,\mathfrak{r}_1))\subseteq\R$ and
    \item \llabel{item2} it holds that $(0,\mathfrak{r}_1)\ni r\mapsto\eta(r)\in\R$
    is continuously differentiable, integrable, and definable\dott
  \end{enumerate}
\end{athm}
\begin{aproof}
  Throughout this proof let $q\colon \Omega\to[0,\infty)$ satisfy for all
  $(r,x)\in\Omega$ that
  \begin{equation}\llabel{def:q}
    q(r,x)=\frac{|(\partial_r Z)(r,x)|}{1+|(\nabla_MZ)(r,x)|}
  \end{equation}
  and let
  $I=\{r\in(0,\mathfrak{r}_0)\colon \eta(r)=+\infty\}$\dott \cfload[.]
  \startnewargseq
  \argument{
    the fact that $Z\in\COne{U}{\R}$ is definable;
    \cref{lem:radial-derivative-well-defined};
    \cref{lem:angular-gradient-well-defined};
  }{
    that $\partial_r Z\colon U\to\R$ and $\nabla_MZ\colon U\to\R^n$ are definable\dott
    \llabel{eq0a}
  }
  \argument{\lref{def:q};
    \lref{eq0a};
    the fact that $\Omega$ is definable;
    the fact that for all $(r,x)\in\Omega$ it holds that $1+|(\nabla_MZ)(r,x)|>0$ ;
  }{
    that $q\colon \Omega\to[0,\infty)$ is definable\dott
    \llabel{eq0b}
  }
  \argument{
    \lref{eq0b};
    the fact that $\Omega$ is definable;
  }{
    that the set
    \begin{eqn}
      \llabel{eq:unbounded-witnesses}
      H=\bigl\{(r,\lambda,x)\in(0,\mathfrak{r}_0)\times[0,\infty)\times M\colon
        \bigl[((r,x)\in\Omega)\wedge(q(r,x)>\lambda)\bigr]\bigr\}
    \end{eqn}
    is definable\dott
  }
  \argument{
    \lref{eq:unbounded-witnesses};
  }{
    that the set
    \begin{eqn}
      \llabel{eq:unbounded-pairs}
      E=\bigl\{(r,\lambda)\in(0,\mathfrak{r}_0)\times[0,\infty)\colon
        \bigl[\exists\,x\in M\colon(r,\lambda,x)\in H\bigr]\bigr\}
    \end{eqn}
    is definable\dott
  }
  \argument{
    \lref{eq:unbounded-pairs};
  }{
    that the set
    \begin{eqn}
      \llabel{eq:bounded-radii}
      J=\bigl\{r\in(0,\mathfrak{r}_0)\colon
        \bigl[\exists\,\lambda\in[0,\infty)\colon(r,\lambda)\notin E\bigr]\bigr\}
    \end{eqn}
    is definable\dott
  }
  \argument{\lref{eq:eta};\lref{eq:unbounded-pairs};
    \lref{eq:bounded-radii};
    \lref{eq:unbounded-witnesses};
    the fact that $I=\{r\in(0,\mathfrak{r}_0)\colon \eta(r)=+\infty\}$;
  }{
    that
    \begin{equation}
      \llabel{eq0c}
      \begin{split}
        I
        &=\bigl\{r\in(0,\mathfrak{r}_0)\colon
            \bigl[\forall\,\lambda\in[0,\infty)\ \exists\,x\in M\colon
            ((r,x)\in\Omega)\wedge(q(r,x)>\lambda)\bigr]\bigr\}\\
        &=\bigl\{r\in(0,\mathfrak{r}_0)\colon
            \bigl[\forall\,\lambda\in[0,\infty)\ \exists\,x\in M\colon
            (r,\lambda,x)\in H\bigr]\bigr\}\\
        &=\bigl\{r\in(0,\mathfrak{r}_0)\colon
            \bigl[\forall\,\lambda\in[0,\infty)\colon(r,\lambda)\in E\bigr]\bigr\}
          =(0,\mathfrak{r}_0)\backslash J
      \end{split}
    \end{equation}
    is definable\dott
  }
  In the following let
  \begin{equation}
    \begin{split}
    G&=\bigl\{(r,0)\colon\bigl[(r\in J)\wedge(\mathcal F_r=\emptyset)\bigr]\bigr\}\\
      &\cup\braces*{(r,y)\in J\times[0,\infty)\colon
        \brackets*{\bigl(\forall\,x\in\mathcal F_r\colon q(r,x)\leq y\bigr)
          \wedge\bigl(\forall\,\varepsilon\in(0,\infty)\ \exists\,x\in\mathcal F_r\colon
            q(r,x)>y-\varepsilon\bigr)}}.
    \end{split}
    \llabel{eq:finite-bound-graph}
  \end{equation}
  \startnewargseq
  \argument{
    \lref{eq0c};\lref{eq:finite-bound-graph};
    the fact that $\Omega$ and $q$ are definable;
  }{
    that $G$ is definable\dott
  }
  \argument{
    \lref{eq:eta};
    \lref{eq0c};
    \lref{eq:finite-bound-graph};
  }{
    that
    \begin{eqn}
      \llabel{eq0d}
      G=\operatorname{graph}\bigl(\eta|_{(0,\mathfrak{r}_0)\backslash I}\bigr)
      \dott
    \end{eqn}
  }
  \argument{
    \lref{eq0d};
    the fact that $G$ is definable;
  }{
    that $\eta|_{(0,\mathfrak{r}_0)\backslash I}$ is definable\dott
  }
  \argument{
    \lref{eq0c};
    the fact that for all definable $A\subseteq\R$ it holds that $A$ is
    a finite union of singletons and open intervals;
  }{
    that
    \begin{eqn}
      \llabel{eq1}
      \bigl(\exists\,\rho\in(0,\mathfrak{r}_0)\colon I\cap(0,\rho)=\emptyset\bigr)
      \quad\text{or}\quad
      \bigl(\exists\,\rho\in(0,\mathfrak{r}_0)\colon (0,\rho)\subseteq I\bigr)
      \dott
    \end{eqn}
  }
  In the following we claim that for all $\rho\in (0,\mathfrak{r}_0)$ it holds that
  \begin{eqn}\llabel{needtoprove}
      (0,\rho)\nsubseteq I.
  \end{eqn}
  We prove \lref{needtoprove} by contradiction. In the following assume for contradiction that $\bigl(\exists\,\rho\in(0,\mathfrak{r}_0)\colon (0,\rho)\subseteq I\bigr)$ and let $\boldsymbol{\rho}_2\in(0,\mathfrak{r}_0)$ satisfy
  $(0,\boldsymbol{\rho}_2)\subseteq I$.
  \startnewargseq
  \argument{
    the assumption that $(0,\boldsymbol{\rho}_2)\subseteq I$;
    the fact that $I=\{r\in(0,\mathfrak{r}_0)\colon \eta(r)=+\infty\}$ ;
    \lref{eq:eta};
  }{
    that for every $r\in(0,\boldsymbol{\rho}_2)$ there exists $x\in\mathcal{F}_r$ such that
    \begin{eqn}\llabel{def:x}
      \frac{|(\partial_r Z)(r,x)|}{1+|(\nabla_MZ)(r,x)|}\geq r^{-2}
      \dott
    \end{eqn}
  }
  In the following let
  \begin{equation}
    \llabel{eq:selector-witness-set}
    W=\bigl\{(r,x)\in\Omega\colon
      \bigl[(r\in(0,\boldsymbol{\rho}_2))\wedge(q(r,x)\geq r^{-2})\bigr]\bigr\}.
  \end{equation}
  \startnewargseq
  \argument{\lref{eq:selector-witness-set};
    the fact that $q$ is definable;
    the fact that $\Omega$ is definable;
  }{
    that $W$ is definable\dott
  }
  \argument{
   \lref{def:x}; \lref{eq:selector-witness-set};
    the definition of $q$;
  }{
    that for all $r\in(0,\boldsymbol{\rho}_2)$ it holds that
    $\mathcal{A}_r=\{x\in M\colon (r,x)\in W\}$ is nonempty\dott
  }
  \argument{
    the fact that for all $r\in(0,\boldsymbol{\rho}_2)$ it holds that $\mathcal{A}_r$ is nonempty;
    the fact that $W\subseteq(0,\boldsymbol{\rho}_2)\times M$;
  }{
    that
    \begin{eqn}
      \llabel{eq:selector-witness-projection}
      \operatorname{proj}_1(W)=(0,\boldsymbol{\rho}_2)
      \dott
    \end{eqn}
  }
  \argument{
    \lref{eq:selector-witness-projection};
    the fact that $W$ is definable;
    the definable choice theorem \cite[Theorem~3.1]{Coste}
  }{
    that there exists a definable function $\chi\colon (0,\boldsymbol{\rho}_2)\to M$ which
    satisfies for every $r\in(0,\boldsymbol{\rho}_2)$ that $(r,\chi(r))\in W$\dott
  }
  \startnewargseq
  \argument{\lref{def:q};    \lref{eq:selector-witness-set};
    the fact that for all $r\in(0,\boldsymbol{\rho}_2)$ it holds  that $(r,\chi(r))\in W$;
  }{
    that for all $r\in(0,\boldsymbol{\rho}_2)$ it holds that $\chi(r)\in\mathcal{F}_r$ and
    \begin{eqn}
      \llabel{eq:infinite-bound-test}
      \frac{|(\partial_r Z)(r,\chi(r))|}{1+|(\nabla_MZ)(r,\chi(r))|}\geq r^{-2}
      \dott
    \end{eqn}
  }
  In the following let $z_*\colon (0,\boldsymbol{\rho}_2)\to\R$ satisfy for all $r\in(0,\boldsymbol{\rho}_2)$ that
  $z_*(r)=Z(r,\chi(r))$\dott
  \startnewargseq
  \argument{
    the fact that $M$ is bounded;
    the fact that $\chi((0,\boldsymbol{\rho}_2))\subseteq M$;
    the fact that $Z$ is bounded on $\Omega$;
    the fact that for all $r\in(0,\boldsymbol{\rho}_2)$ it holds that $(r,\chi(r))\in\Omega$;
  }{
    that $\chi$ and $z_*$ are bounded and definable\dott
    \llabel{eq2}
  }
  \argument{
    \lref{eq2};
    \cref{lem:support1};
  }{
    that there exists $\boldsymbol{\rho}_3\in(0,\boldsymbol{\rho}_2)$ which satisfies that
    \begin{eqn}
      \llabel{eq3}
      \chi|_{(0,\boldsymbol{\rho}_3)}\in C^1((0,\boldsymbol{\rho}_3),M),
      \quad
      z_*|_{(0,\boldsymbol{\rho}_3)}\in C^1((0,\boldsymbol{\rho}_3),\R),
      \quad\text{and}\quad
      \int_0^{\boldsymbol{\rho}_3}\bigl(|\chi'(r)|+|(z_*)'(r)|\bigr)\dd r<\infty
      \dott
    \end{eqn}
  }
  \startnewargseq
  \argument{
    \lref{eq3};
    the fact that for all
    $r\in(0,\boldsymbol{\rho}_3)$, $\varepsilon\in(0,\infty)$ with
    $(r-\varepsilon,r+\varepsilon)\subseteq(0,\boldsymbol{\rho}_3)$
    it holds that
    $\chi|_{(r-\varepsilon,r+\varepsilon)}
      \in C^1((r-\varepsilon,r+\varepsilon),M)$ and
    $(\chi|_{(r-\varepsilon,r+\varepsilon)})'(r)=\chi'(r)$;
  }{
    that for all $r\in(0,\boldsymbol{\rho}_3)$ it holds that
    \begin{eqn}
      \llabel{eq:selector-tangent}
      \chi'(r)\in T_{\chi(r)}M
      \dott
    \end{eqn}
  }
  \argument{
    \lref{eq:selector-tangent};
    the fact that for all $r\in(0,\boldsymbol{\rho}_3)$ it holds that $z_*(r)=Z(r,\chi(r))$;
    the fact that for all $r\in(0,\boldsymbol{\rho}_3)$ it holds that
    $(r,\chi(r))\in U$;
    the fact that $Z\in\COne{U}{\R}$;
    the chain rule;
  }{
    that for all $r\in(0,\boldsymbol{\rho}_3)$ it holds that
    \begin{eqn}
      \llabel{eq4}
      (z_*)'(r)=(\partial_r Z)(r,\chi(r))+
      \langle(\nabla_MZ)(r,\chi(r)),\chi'(r)\rangle
      \dott
    \end{eqn}
  }
  \argument{
    \lref{eq4};
    the triangle inequality;
    the Cauchy--Schwarz inequality;
  }{
    that for all $r\in(0,\boldsymbol{\rho}_3)$ it holds that
    \begin{equation}\label{eq:normalized-trace}
      \frac{|(\partial_r Z)(r,\chi(r))|}{1+|(\nabla_MZ)(r,\chi(r))|}
      \leq |(z_*)'(r)|+|\chi'(r)|
      \dott
    \end{equation}
  }
  \argument{
    \lref{eq:infinite-bound-test};
    \eqref{eq:normalized-trace};
  }{
    that for all $r\in(0,\boldsymbol{\rho}_3)$ it holds that
    \begin{eqn}
      \llabel{eq5}
      r^{-2}\leq |(z_*)'(r)|+|\chi'(r)|
      \dott
    \end{eqn}
  }
  \argument{
    \lref{eq5};
    the fact that $\int_0^{\boldsymbol{\rho}_3} r^{-2}\dd r=+\infty$;
  }{
    that
    \begin{eqn}
      \llabel{eq6}
      \int_0^{\boldsymbol{\rho}_3}\bigl(|\chi'(r)|+|(z_*)'(r)|\bigr)\dd r=+\infty
      \dott
    \end{eqn}
  }
  \argument{
    \lref{eq6};
    \lref{eq3};
  }[verbs=p]{
    that
    \begin{equation}
      \llabel{eq6.1}
      \infty
        =\int_0^{\boldsymbol{\rho}_3}\bigl(|\chi'(r)|+|(z_*)'(r)|\bigr)\dd r
        <\infty
      \dott
    \end{equation}
  }
  \startnewargseq
  This contradiction proves \lref{needtoprove}\dott
  \startnewargseq

\argument{\lref{eq1};\lref{needtoprove};}{ that \begin{eqn}
      \llabel{eqt1}
     \exists\,\rho\in(0,\mathfrak{r}_0)\colon I\cap(0,\rho)=\emptyset.
    \end{eqn}}
  \argument{\lref{eqt1};
    the fact that $I=\{r\in(0,\mathfrak{r}_0)\colon \eta(r)=+\infty\}$ 
  }{
    that there exists $\rho\in (0,\mathfrak{r}_0)$ which satisfies for all $r\in(0,\rho)$ that 
    \begin{eqn}
      \llabel{eq7}
      \eta(r)<\infty
      \dott
    \end{eqn}
  }
  \startnewargseq
  \argument{
    \lref{eq7};
    the fact that $\eta|_{(0,\mathfrak{r}_0)\backslash I}$ is definable;
  }{
    that $\eta|_{(0,\rho)}\colon (0,\rho)\to[0,\infty)$ is definable\dott
    \llabel{eq8}
  }
  \argument{
    \lref{eq8};
    the fact that the zero set of $\eta|_{(0,\rho)}$ is definable;
  }{
    that there exists $\boldsymbol{\rho}_1\in (0,\rho)$ which satisfies that
    \begin{eqn}
      \llabel{eq9}
      \bigl(\forall\, r\in(0,\boldsymbol{\rho}_1)\colon \eta(r)=0\bigr)
      \quad\text{or}\quad
      \bigl(\forall\, r\in(0,\boldsymbol{\rho}_1)\colon \eta(r)>0\bigr)
      \dott
    \end{eqn}
  }
  In the following we prove that there exists $\varrho\in (0,\boldsymbol{\rho}_1)$ such that
  \begin{eqn}
    \llabel{eq17}
    \int_0^{\varrho}\eta(r)\dd r<\infty
    \dott
  \end{eqn}
  \startnewargseq
  In our proof of \lref{eq17} we distinguish between the case that for
  all $r\in(0,\boldsymbol{\rho}_1)$ it holds that $\eta(r)=0$ and the case that for
  all $r\in(0,\boldsymbol{\rho}_1)$ it holds that $\eta(r)>0$.
  \startnewargseq
  We first prove \lref{eq17} in the case that
  \begin{eqn}\llabel{case1}
    \forall\, r\in(0,\boldsymbol{\rho}_1)\colon \eta(r)=0
    \dott
  \end{eqn}
  \startnewargseq
  \argument{
    \lref{case1}
  }{
    that
    \begin{eqn}
      \llabel{eq10}
      \int_0^{\boldsymbol{\rho}_1}\eta(r)\dd r=0<\infty
      \dott
    \end{eqn}
  }
  This proves \lref{eq17} in the case that for all $r\in(0,\boldsymbol{\rho}_1)$ it
  holds that $\eta(r)=0$.
  \startnewargseq
  We now prove \lref{eq17} in the case that
  \begin{eqn}\llabel{case2}
    \forall\, r\in(0,\boldsymbol{\rho}_1)\colon \eta(r)>0
    \dott
  \end{eqn}
  \startnewargseq
  \argument{\lref{eq:eta};
    \lref{case2};
  }{
    that for every $r\in(0,\boldsymbol{\rho}_1)$ there exists $x\in\mathcal{F}_r$ such that
    \begin{eqn}
      \llabel{eq11}
      \frac{|(\partial_r Z)(r,x)|}{1+|(\nabla_MZ)(r,x)|}\geq\frac12\eta(r)
      \dott
    \end{eqn}
  }
 In the following let
  \begin{eqn}
    \llabel{eq:half-bound-witness-set}
    K
    =\bigl\{(r,x)\in\Omega\colon
      \bigl[(r\in(0,\boldsymbol{\rho}_1))\wedge(q(r,x)\geq\tfrac12\eta(r))\bigr]\bigr\}.
  \end{eqn}
  \startnewargseq
  \argument{
    the fact that $q$, $\eta|_{(0,\boldsymbol{\rho}_1)}$, $\Omega$, and $(0,\boldsymbol{\rho}_1)$
    are definable;
    \lref{eq:half-bound-witness-set};
  }{
    that $K$ is definable\dott
  }
  \argument{\lref{def:q};
    \lref{eq11};
    \lref{eq:half-bound-witness-set};
  }{
    that for all $r\in(0,\boldsymbol{\rho}_1)$ it holds that
    $\mathcal{Q}_r=\{x\in M\colon(r,x)\in K\}$ is
    nonempty\dott\llabel{argt5}
  }
  \argument{\lref{argt5};
    the fact that $K\subseteq(0,\boldsymbol{\rho}_1)\times M$;
  }{
    that
    \begin{eqn}
      \llabel{eq:half-bound-witness-projection}
      \operatorname{proj}_1(K)=(0,\boldsymbol{\rho}_1)
      \dott
    \end{eqn}
  }
  \argument{
    \lref{eq:half-bound-witness-projection};
    the fact that $K$ is definable;
    the definable choice theorem \cite[Theorem~3.1]{Coste};
  }{
    that there exists a definable function $\xi\colon (0,\boldsymbol{\rho}_1)\to M$
    which satisfies for all $r\in(0,\boldsymbol{\rho}_1)$ that
    $(r,\xi(r))\in K$\dott
  }
  \startnewargseq
  \argument{\lref{def:q};  \lref{eq:half-bound-witness-set};
    the fact that for all $r\in(0,\boldsymbol{\rho}_1)$ it holds that
    $(r,\xi(r))\in K$;
  }{
    that for all $r\in(0,\boldsymbol{\rho}_1)$ it holds that
    \begin{eqn}
      \llabel{eq12}
      \xi(r)\in\mathcal{F}_r
      \qquad\text{and}\qquad
      \frac{|(\partial_r Z)(r,\xi(r))|}
      {1+|(\nabla_MZ)(r,\xi(r))|}\geq\frac12\eta(r)
      \dott
    \end{eqn}
  }
  In the following let $\zeta\colon (0,\boldsymbol{\rho}_1)\to\R$ satisfy for all $r\in(0,\boldsymbol{\rho}_1)$ that
  $\zeta(r)=Z(r,\xi(r))$\dott
  \startnewargseq
  \argument{
    the fact that $M$ is bounded;
    the fact that $Z$ is bounded on $\Omega$;
    the fact that for all $r\in(0,\boldsymbol{\rho}_1)$ it holds that
    $(r,\xi(r))\in\Omega$;
    \cref{lem:support1};
  }{
    that there exists $\boldsymbol{\rho}_4\in (0,\boldsymbol{\rho}_1)$ which satisfies that
    \begin{eqn}
      \llabel{eq13}
     \xi|_{(0,\boldsymbol{\rho}_4)}\in C^1((0,\boldsymbol{\rho}_4),M),
      \quad
      \zeta|_{(0,\boldsymbol{\rho}_4)}\in C^1((0,\boldsymbol{\rho}_4),\R),
      \quad\text{and}\quad  \int_0^{\boldsymbol{\rho}_4}
      \bigl(|\xi'(r)|+|\zeta'(r)|\bigr)\dd r<\infty
      \dott
    \end{eqn}
  }
  \startnewargseq
  \argument{
    \lref{eq12};
    \lref{eq13};
    the fact that $Z\in\COne{U}{\R}$;
    the fact that for all
    $r\in(0,\boldsymbol{\rho}_4)$ it holds that $\zeta(r)=Z(r,\xi(r))$;
    the chain rule;
    the triangle inequality;
    the Cauchy--Schwarz inequality;
  }{
    that for all $r\in(0,\boldsymbol{\rho}_4)$ it holds that
    \begin{eqn}
      \llabel{eq14}
      \frac{|(\partial_r Z)(r,\xi(r))|}
      {1+|(\nabla_MZ)(r,\xi(r))|}
      \leq |\zeta'(r)|+|\xi'(r)|
      \dott
    \end{eqn}
  }
  \argument{
    \lref{eq12};
    \lref{eq14};
  }{
    that for all $r\in(0,\boldsymbol{\rho}_4)$ it holds that
    \begin{eqn}
      \llabel{eq15}
      \eta(r)\leq
      2|\zeta'(r)|+2|\xi'(r)|
      \dott
    \end{eqn}
  }
  \argument{
    \lref{eq15};
    \lref{eq13};
    the fact that $\eta$ is nonnegative;
  }{
    that
    \begin{eqn}
      \llabel{eq16}
      \int_0^{\boldsymbol{\rho}_4}\eta(r)\dd r<\infty
      \dott
    \end{eqn}
  }
  This proves \lref{eq17} in the case that for all $r\in(0,\boldsymbol{\rho}_1)$ it
  holds that $\eta(r)>0$.
  \startnewargseq

  \argument{
    the fact that $\eta|_{(0,\rho)}\colon (0,\rho)\to[0,\infty)$ is definable;
    the fact that $0<\varrho<\rho$;
    \cite[Corollary~6.5]{Coste};
  }{
    that there exists $\mathfrak{r}_1\in(0,\varrho)$ such that
    $\eta|_{(0,\mathfrak{r}_1)}\in C^1((0,\mathfrak{r}_1),\R)$ is
    locally bounded\dott
    \llabel{eq18}
  }
  \argument{
    \lref{eq18};
      \lref{eq7};
    \lref{eq17};
    the fact that $\eta|_{(0,\rho)}$ is definable;
  }[verbs=ep,plural]{
    \lref{item1} and \lref{item2}\dott
  }
\end{aproof}

In the following result we specialize \cref{lem:normalized} to the
case where $M=(0,H)$ is the $1$-dimensional $C^1$-submanifold $(0,H)$
of $\R$ for some strictly positive real number $H\in(0,\infty)$.

\begin{athm}{lemma}{lem:normalized-interval}
  Let $\mathfrak r,H\in(0,\infty)$, let
  $U\subseteq(0,\mathfrak r)\times(0,H)$ be definable and open,
  let $\Omega\subseteq U$ be definable, let
  $Z=(Z(r,x))_{(r,x)\in U}\in C^1(U,\R)$ be definable,
  assume $\sup_{v\in\Omega}|Z(v)|<\infty$, and let
  $\eta\colon(0,\mathfrak r)\to[0,\infty]$
  satisfy that for all $r\in(0,\mathfrak r)$ it holds that
  \begin{eqn}
    \llabel{ass1}
    \eta(r)=\sup\Biggl(\{0\}\cup\Biggl\{z\in\R\colon
      \Biggl[\exists\,x\in(0,H)\colon
      \biggl(\bigl[(r,x)\in\Omega\bigr]\wedge
        \biggl[z=\frac{|(\partial_r Z)(r,x)|}{1+|(\partial_x Z)(r,x)|}\biggr]\biggr)
    \Biggr]\Biggr\}\Biggr)
    \dott
  \end{eqn}
  Then there exists $\rho\in(0,\mathfrak r)$ such that
  \begin{enumerate}[label=(\roman*)]
    \item \llabel{item1} it holds that $\eta((0,\rho))\subseteq\R$ and
    \item \llabel{item2} it holds that $(0,\rho)\ni r\mapsto\eta(r)\in\R$
    is continuously differentiable, integrable, and definable\dott
  \end{enumerate}
\end{athm}
\begin{aproof}
  \argument{
    the fact that $(0,H)$ is a bounded and definable
    $1$-dimensional $C^1$-submanifold of $\R$;
    the fact that for all $(r,x)\in U$ it holds that
    $(\nabla_{(0,H)}Z)(r,x)=(\partial_x Z)(r,x)$;
    \lref{ass1};
    \cref{lem:normalized} (applied with $n\curvearrowleft1$,
    $\mathfrak r_0\curvearrowleft\mathfrak r$,
    $M\curvearrowleft(0,H)$,
    $\mathcal F_r\curvearrowleft\{x\in(0,H)\colon(r,x)\in\Omega\}$
    for every $r\in(0,\mathfrak r)$ in the notation of
    \cref{lem:normalized});
  }[verbs=ep,plural]{
    \lref{item1} and \lref{item2}\dott
  }
\end{aproof}

\begin{athm}{lemma}{lem:contact-bound}
  Let $n\in\N$, $\mathfrak{r}_0\in(0,\infty)$, let
  $M\subseteq\R^n$ be a bounded and definable \submanifold{1}{n},
  let $\mathcal O\subseteq(0,\mathfrak{r}_0)\times M$ be definable
  and open in $(0,\mathfrak{r}_0)\times M$, let
  $h=(h(r,x))_{(r,x)\in\mathcal O}\in \COne{\mathcal O}{\R}$,
  $a\in C^1((0,\mathfrak{r}_0),\R)$ be definable, and let
  $\beta\colon(0,\mathfrak{r}_0)\to[0,\infty]$ satisfy for all
  $r\in(0,\mathfrak{r}_0)$ that
  \begin{multline}
    \llabel{ass2}
    \beta(r)=\sup\Bigl(\{0\}\cup\Bigl\{z\in\R\colon
      \Bigl[\exists\,x\in
      \bigl\{y\in M\colon(r,y)\in\mathcal O\cap h^{-1}(\{a(r)\})\bigr\}\colon\\
      z|(\nabla_Mh)(r,x)|=|(\partial_r h)(r,x)-a'(r)|>0
    \Bigr]\Bigr\}\Bigr)\cfadd{definition:c1-on-manifold}\dott
  \end{multline}
  \cfload[.]
  Then there exists $\mathfrak{r}_1\in(0,\mathfrak{r}_0)$ such that
  \begin{enumerate}[label=(\roman*)]
    \item \llabel{item1} it holds that $\beta((0,\mathfrak{r}_1))\subseteq\R$ and
    \item \llabel{item2} it holds that $(0,\mathfrak{r}_1)\ni r\mapsto\beta(r)\in\R$
    is continuously differentiable, integrable, and definable\dott
  \end{enumerate}
  \cfout[.]
\end{athm}
\begin{aproof}
  Throughout this proof let $G\colon\mathcal O\to\R$ satisfy for all
  $v\in\mathcal O$ that $G(v)=|(\nabla_Mh)(v)|$, let
  \begin{eqn}\llabel{def:S}
    S=\bigl\{(r,x)\in\mathcal O\colon
       \bigl[(h(r,x)=a(r))\wedge(G(r,x)>0)\bigr]\bigr\},
  \end{eqn}
  let $q\colon S\to[0,\infty)$ satisfy that for all $(r,x)\in S$ it holds that
  \begin{eqn}\llabel{def:q}
    q(r,x)=\frac{|(\partial_r h)(r,x)-a'(r)|}{G(r,x)},
  \end{eqn}
  and let
  $I=\{r\in(0,\mathfrak{r}_0)\colon\beta(r)=+\infty\}$\dott
  \cfload[.]
  \startnewargseq
  \argument{
    the fact that $h\in \COne{\mathcal O}{\R}$ is definable;
    the fact that $\mathcal O$ is definable and open in
    $(0,\mathfrak{r}_0)\times M$;
    \cref{lem:radial-derivative-well-defined};
    \cref{lem:angular-gradient-well-defined};
  }{
    that $\partial_r h\colon\mathcal O\to\R$ and
    $\nabla_Mh\colon\mathcal O\to\R^n$ are \llabel{arg1} definable\dott
  }
  \argument{
   \lref{arg1};
    the fact that for all
    $(r,x)\in\mathcal O$ it holds that $G(r,x)=|(\nabla_Mh)(r,x)|$;
  }{
    that $G\colon\mathcal O\to\R$ is definable\dott
  }
  \argument{
    the fact that $a\in C^1((0,\mathfrak{r}_0),\R)$ is definable;
    \cref{lem:radial-derivative-well-defined} (applied with
    $n\curvearrowleft1$, $\mathfrak{r}\curvearrowleft\mathfrak{r}_0$,
    $M\curvearrowleft\{0\}$,
    $U\curvearrowleft(0,\mathfrak{r}_0)\times\{0\}$,
    $f\curvearrowleft\bigl((0,\mathfrak{r}_0)\times\{0\}\ni(r,0)
      \mapsto a(r)\in\R\bigr)$ in the notation of
    \cref{lem:radial-derivative-well-defined});
  }{
    that $a'\colon(0,\mathfrak{r}_0)\to\R$ is definable\dott
  }
  \argument{
    the fact that $h$, $a$, and $G$ are definable;
  }{
    that $S$ is \llabel{arg3} definable\dott
  }
  \argument{\lref{arg3};
    the fact that $\partial_r h$, $a'$, and $G$ are definable;
    the fact that for all $(r,x)\in S$ it holds that $G(r,x)>0$;
    \lref{def:q}
  }{
    that $q\colon S\to[0,\infty)$ is well-defined and definable\dott
    \llabel{eq1}
  }
  \argument{
    \lref{eq1};
    \lref{ass2};
  }{
    that $I$ is definable\dott
    \llabel{eq2}
  }
  \argument{
    \lref{eq2};
    \cite[Theorem~2.10]{Coste};
  }{
    that there exists $R\in(0,\mathfrak{r}_0)$ which satisfies
    \begin{eqn}
      \llabel{eq3}
      (0,R)\subseteq I
      \quad\text{or}\quad
      (0,R)\cap I=\emptyset\dott
    \end{eqn}
  }
  \startnewargseq
  We claim that
  \begin{eqn}
    \llabel{eq3.1}
    (0,R)\not\subseteq I
    \dott
  \end{eqn}
  \startnewargseq
  To prove this claim, assume for contradiction that
  \begin{eqn}
    \llabel{eq3.2}
    (0,R)\subseteq I
    \dott
  \end{eqn}
  \startnewargseq
  \argument{\lref{ass2};
    \lref{eq3.2};
the fact that $I=\{r\in(0,\mathfrak{r}_0)\colon\beta(r)=+\infty\}$;
  }{
    that for all $r\in(0,R)$ there exists $x\in M$ such that \llabel{arg5}
    $(r,x)\in S$ and $q(r,x)\geq r^{-2}$\dott
  }
  \argument{\lref{arg5};the definable choice theorem \cite[Theorem~3.1]{Coste}}{that there exists a definable function $\chi\colon(0,R)\to M$ 
  which satisfies for all $r\in(0,R)$ that
  \begin{eqn}
    \llabel{eq4}
    (r,\chi(r))\in S
    \quad\text{and}\quad
    q(r,\chi(r))\geq r^{-2}.
  \end{eqn}}
  \startnewargseq
  \argument{
    \lref{eq4};
    the fact that $M$ is bounded;
    \cref{lem:support1};
  }{that there exists $\mathfrak{R}_x\in(0,R)$ which satisfies that
  \begin{eqn}
    \llabel{eq5}
    \chi|_{(0,\mathfrak{R}_x)}\in C^1((0,\mathfrak{R}_x),M)
    \quad\text{and}\quad
    \int_0^{\mathfrak{R}_x}|\chi'(r)|\dd r<\infty.
  \end{eqn}}
  \startnewargseq
  \argument{\lref{def:S};\lref{eq4};
  }{
    that for all $r\in(0,\mathfrak{R}_x)$ it holds that
    \begin{eqn}
      \llabel{eq6}
      h(r,\chi(r))=a(r)
      \dott
    \end{eqn}
  }
  \argument{
    \lref{eq6};
    the fact that $h\in \COne{\mathcal O}{\R}$,
    $a\in C^1((0,\mathfrak{r}_0),\R)$, and
    $\chi|_{(0,\mathfrak{R}_x)}\in C^1((0,\mathfrak{R}_x),M)$;
    the chain rule;
  }{
    that for all $r\in(0,\mathfrak{R}_x)$ it holds that
    \begin{eqn}[a]
      \llabel{eq7}
      0
      &=
      \frac{\ddd}{\ddd r}(h(r,\chi(r))-a(r))
      =
      (\partial_r h)(r,\chi(r))
      +\langle(\nabla_Mh)(r,\chi(r)),\chi'(r)\rangle
      -a'(r)
      \dott
    \end{eqn}
  }
  \argument{
    \lref{eq7}
  }{
    that for all $r\in(0,\mathfrak{R}_x)$ it holds that
    \begin{eqn}
      \llabel{eq7.1}
      (\partial_r h)(r,\chi(r))-a'(r)
      =-\langle(\nabla_Mh)(r,\chi(r)),\chi'(r)\rangle
      \dott
    \end{eqn}
  }
  \argument{\lref{def:q};
    \lref{eq7.1};
    the Cauchy--Schwarz inequality;
    the fact that $G(r,x)=|(\nabla_Mh)(r,x)|$;
    the fact that for all $r\in(0,R)$ it holds that
    $G(r,\chi(r))>0$, by \lref{eq4};
  }{
    that for all $r\in(0,\mathfrak{R}_x)$ it holds that
    \begin{eqn}[a]
      \llabel{eq8}
      q(r,\chi(r))
      &=\frac{|(\partial_r h)(r,\chi(r))-a'(r)|}{G(r,\chi(r))}
      =\frac{|\langle(\nabla_Mh)(r,\chi(r)),\chi'(r)\rangle|}
      {|(\nabla_Mh)(r,\chi(r))|}
      \\
      &\leq\frac{|(\nabla_Mh)(r,\chi(r))|\,|\chi'(r)|}
      {|(\nabla_Mh)(r,\chi(r))|}
      =|\chi'(r)|
      \dott
    \end{eqn}
  }
  \argument{
   \lref{eq4};
    \lref{eq8};
  }{
    that for all $r\in(0,\mathfrak{R}_x)$ it holds that
    \begin{eqn}
      \llabel{eq9}
      r^{-2}\leq|\chi'(r)|
      \dott
    \end{eqn}
  }
  \argument{
    \lref{eq9};
    the fact that $\int_0^{\mathfrak{R}_x}r^{-2}\dd r=+\infty$;
  }{
    that
    \begin{eqn}
      \llabel{eq10}
      \int_0^{\mathfrak{R}_x}|\chi'(r)|\dd r=+\infty
      \dott
    \end{eqn}
  }
  \argument{
    \lref{eq10};
    \lref{eq5};
  }{
    that the assumption in \lref{eq3.2} leads to a contradiction\dott
    \llabel{eq11}
  }
  \argument{
    \lref{eq11};
  }[verbs=ep]{
    \lref{eq3.1}\dott
  }

  \startnewargseq
  \argument{\lref{eq3};
    \lref{eq3.1};the fact that $I=\{r\in(0,\mathfrak{r}_0)\colon\beta(r)=+\infty\}$
  }{
    that for all $r\in(0,R)$ it holds that
    $\beta(r)<+\infty$\dott
    \llabel{eq12}
  }
  \argument{
    \lref{ass2};
    \lref{def:q};
    \lref{eq12};
    the fact that $q$ is definable;
  }{
    that $\beta|_{(0,R)}\colon(0,R)\to[0,\infty)$ is definable\dott
  }
  \argument{
    the fact that $\beta|_{(0,R)}$ is definable;
    \cite[Theorem~2.10]{Coste};
    the fact that $\beta$ is nonnegative;
  }{that there exists $\mathfrak{R}_\beta\in(0,R)$ which satisfies
  \begin{eqn}
    \llabel{eq13}
    \beta|_{(0,\mathfrak{R}_\beta)}=0
    \quad\text{or}\quad
    \forall\, r\in(0,\mathfrak{R}_\beta)\colon\beta(r)>0;
  \end{eqn}}
  \startnewargseq
  We claim that
  \begin{eqn}
    \llabel{eq13.1}
    \exists\, \mathfrak{R}_{\mathrm{int}}\in(0,\mathfrak{R}_\beta)\colon
    \int_0^{\mathfrak{R}_{\mathrm{int}}}\beta(r)\dd r<\infty
    \dott
  \end{eqn}
  \startnewargseq
  We prove \lref{eq13.1} by contradiction\dott
  \startnewargseq
  In the following we thus assume that for all $\rho\in(0,\mathfrak{R}_\beta)$
  it holds that
  \begin{eqn}
    \llabel{eq13.1a}
    \int_0^\rho\beta(r)\dd r=\infty
    \dott
  \end{eqn}
  \startnewargseq
  \argument{
    \lref{eq13};
    \lref{eq13.1a};
  }{
    that for all $r\in(0,\mathfrak{R}_\beta)$ it holds that
    \begin{eqn}
      \llabel{eq13.2}
      \beta(r)>0
      \dott
    \end{eqn}
  }
  \argument{\lref{ass2};
    \lref{eq13.2};
  }{
    that for all $r\in(0,\mathfrak{R}_\beta)$ there exists $y\in M$ such that \llabel{arg6}
    $(r,y)\in S$ and $q(r,y)\geq\frac12\beta(r)$\dott
  }
  \argument{\lref{arg6};the definable choice theorem \cite[Theorem~3.1]{Coste}}{that there exists a
  definable function $\psi\colon(0,\mathfrak{R}_\beta)\to M$
  which satisfies that for all $r\in(0,\mathfrak{R}_\beta)$ it holds that
  \begin{eqn}
    \llabel{eq14}
    (r,\psi(r))\in S
    \quad\text{and}\quad
    q(r,\psi(r))\geq\frac12\beta(r);
  \end{eqn}}
  \startnewargseq
  \argument{
    \lref{eq14};
    the fact that $M$ is bounded;
    \cref{lem:support1};
  }{
    that there exists $\mathfrak{R}_y\in(0,\mathfrak{R}_\beta)$ which satisfies that
  \begin{eqn}
    \llabel{eq15}
    \psi|_{(0,\mathfrak{R}_y)}\in C^1((0,\mathfrak{R}_y),M)
    \quad\text{and}\quad
    \int_0^{\mathfrak{R}_y}|\psi'(r)|\dd r<\infty.
  \end{eqn}}
  \startnewargseq
  \argument{\lref{def:S};\lref{eq14};
  }{
    that for all $r\in(0,\mathfrak{R}_y)$ it holds that
    \begin{eqn}
      \llabel{eq15.1}
      h(r,\psi(r))=a(r)
      \dott
    \end{eqn}
  }
  \argument{
    \lref{eq15.1};
    the fact that $h\in \COne{\mathcal O}{\R}$,
    $a\in C^1((0,\mathfrak{r}_0),\R)$, and
    $\psi|_{(0,\mathfrak{R}_y)}\in C^1((0,\mathfrak{R}_y),M)$;
    the chain rule;
  }{
    that for all $r\in(0,\mathfrak{R}_y)$ it holds that
    \begin{eqn}[a]
      \llabel{eq15.2}
      0
      &=
      \frac{\ddd}{\ddd r}(h(r,\psi(r))-a(r))
      =
      (\partial_r h)(r,\psi(r))
      +\langle(\nabla_Mh)(r,\psi(r)),\psi'(r)\rangle
      -a'(r)
      \dott
    \end{eqn}
  }
  \argument{
    \lref{eq15.2};
  }{
    that for all $r\in(0,\mathfrak{R}_y)$ it holds that
    \begin{eqn}
      \llabel{eq15.3}
      (\partial_r h)(r,\psi(r))-a'(r)
      =-\langle(\nabla_Mh)(r,\psi(r)),\psi'(r)\rangle
      \dott
    \end{eqn}
  }
  \argument{
    \lref{def:S};
    \lref{eq14};
  }{
    that for all $r\in(0,\mathfrak{R}_\beta)$ it holds that
    \begin{eqn}
      \llabel{eq15.4}
      G(r,\psi(r))>0
      \dott
    \end{eqn}
  }
  \argument{
    \lref{eq15.4};
    \lref{def:q};
    \lref{eq15.3};
    the Cauchy--Schwarz inequality;
    the fact that for all $(r,x)\in\mathcal O$ it holds that
    $G(r,x)=|(\nabla_Mh)(r,x)|$;
  }{
    that for all $r\in(0,\mathfrak{R}_y)$ it holds that
    \begin{eqn}[a]
      \llabel{eq16}
      q(r,\psi(r))
      &=\frac{|(\partial_r h)(r,\psi(r))-a'(r)|}{G(r,\psi(r))}
      =\frac{|\langle(\nabla_Mh)(r,\psi(r)),\psi'(r)\rangle|}
      {|(\nabla_Mh)(r,\psi(r))|}
      \\
      &\leq\frac{|(\nabla_Mh)(r,\psi(r))|\,|\psi'(r)|}
      {|(\nabla_Mh)(r,\psi(r))|}
      =|\psi'(r)|
      \dott
    \end{eqn}
  }
  \argument{
 \lref{eq14};
    \lref{eq16};
  }{
    that for all $r\in(0,\mathfrak{R}_y)$ it holds that
    \begin{eqn}
      \llabel{eq17}
      \beta(r)\leq2|\psi'(r)|
      \dott
    \end{eqn}
  }
  \argument{
    \lref{eq17};
    \lref{eq13.1a};
    \lref{eq15};
    the fact that $\beta$ is nonnegative;
  }{
    that
    \begin{eqn}
      \llabel{eq18}
      \infty
        =\int_0^{\mathfrak{R}_y}\beta(r)\dd r
        \leq2\int_0^{\mathfrak{R}_y}|\psi'(r)|\dd r
        <\infty
      \dott
    \end{eqn}
  }
  This contradiction shows
    \lref{eq13.1}\dott
  \startnewargseq

\argument{\lref{eq13.1}}{that there exists $\mathfrak{R}_{\mathrm{int}}\in(0,\mathfrak{R}_\beta)$ which satisfies that
  \begin{eqn}
    \llabel{eq19}
    \int_0^{\mathfrak{R}_{\mathrm{int}}}\beta(r)\dd r<\infty.
  \end{eqn}}
  \startnewargseq
  \argument{
    \lref{eq19};
    the fact that $\beta|_{(0,\mathfrak{R}_{\mathrm{int}})}$ has finite values
    and is definable;
    \cite[Corollary~6.5]{Coste};
    the fact that $\beta$ is nonnegative;
  }[verbs=ep,plural]{ \lref{item1,item2}\dott
  }
\end{aproof}

\section{Construction of the radial-energy regions and the energy transformations}
\label{sec:atlas}

As outlined in \cref{subsec:proof-outline} above, the key issue in the
proof of \cref{maintheorem} is to establish Lyapunov-type stability
(in the precise sense of \eqref{eq:intro-radial-angular-lyapunov} in
\cref{subsec:proof-outline}) for the radial-angular \DE\ system
(see \eqref{eq:intro-radial-angular-system} in \cref{subsec:proof-outline})
as asserted in \cref{prop:buildQ} in
\cref{sec:radial-integration} below.
The system in \eqref{eq:intro-radial-angular-system} already admits a
Lyapunov function, namely, the scalar field $h$ itself.
However, the dissipation associated to this Lyapunov function and its
Lyapunov energy $h(r(t),x(t))$ is not large enough to directly
allow us to conclude finiteness of the length of the angular component
(see \eqref{eq:intro-original-energy-dissipation} in
\cref{subsec:proof-outline}).

It is the key subject of this section to construct and study suitable
transformation functions that we apply to the energy $h(r(t),x(t))$ to
obtain new Lyapunov-type functions with larger dissipation.
In \cref{subsec:cell-decomposition} we briefly recall some well-known
concepts from the area of o-minimal geometry that we employ in the
construction of the energy transformations.
In \cref{lem:divide} in \cref{subsec: energy coordinates} we construct the
desired transformations $\mathcal F_0,\mathcal F_1,\dots,\mathcal F_{N-1}$
and their associated back-transformations $E_1,E_2,\dots,E_N$ and establish
lower and upper bounds for the partial derivative with respect to the
energy variable of the back-transformations $(\partial_z E_i)(r,z)$.
In \cref{lem:residual} in \cref{subsec:residual-bound} we establish
lower and upper bounds for the partial derivative with respect to the
radial variable of the back-transformations $(\partial_r E_i)(r,z)$.
\Cref{lem:radial-decay} in \cref{subsec:residual-bound} establishes in
item~\ref{lem:radial-decay.item2} a basic growth bound near the origin for
definable and integrable functions, which we employ in the proof of
\cref{lem:residual}.

\subsection{Smooth cylindrical definable cell decompositions}
\label{subsec:cell-decomposition}

In this subsection we recall some standard concepts from o-minimal geometry.
We refer, \eg, to \cite[Definition~2.4 and Section~6.2]{Coste} for further reading.

\begin{athm}{definition}{definition:definable-partition}[Definable partitions]
  Let $d\in\N$, let $A\subseteq\R^d$ be a set, and let
  $\mathscr P$ be a set\dott \cfload[.]
  Then we say that $\mathscr P$ is a definable partition of $A$ if and
  only if
  \begin{enumerate}[label={(\roman*)}]
    \item it holds that $\mathscr P$ is a finite set,
    \item it holds that $\midcup_{B\in\mathscr P}B=A$,
    \item it holds for all $B\in\mathscr P$ that $B$ is nonempty and
    definable\cfadd{definition:definable-set}, and
    \item it holds for all $B,C\in\mathscr P$ with $B\ne C$ that
    $B\cap C=\emptyset$\dott
  \end{enumerate}
  \cfout[.]
\end{athm}

\begin{athm}{definition}{definition:planar-cell-decomposition}
  [Smooth cylindrical definable cell decompositions]
  Let $k\in\N$\dott \cfload[.]
  Then
  \begin{enumerate}[label={(\roman*)}]
    \item for every set $\mathscr P$ we say that $\mathscr P$ is a cylindrical definable $C^k$ cell
    decomposition of $\R$ if and only if
    \begin{enumerate}[label={(\alph*)}]
      \item it holds that $\mathscr P$ is a definable partition of
      $\R$\cfadd{definition:definable-partition} and
      \item it holds for all $A\in\mathscr P$ that $A$ is a singleton or a
      non-empty open interval,
    \end{enumerate}
    and
    \item for every set $\mathscr P$ and every $d\in\N$ we say that $\mathscr P$ is a
    cylindrical definable $C^k$ cell decomposition of $\R^{d+1}$ if and only if
    \begin{enumerate}[label={(\alph*)}]
      \item it holds that $\mathscr P$ is a definable partition of
      $\R^{d+1}$\cfadd{definition:definable-partition},
      \item it holds that
      $\midcup_{E\in\mathscr P}\bigl\{\{x=(x_1,\dots,x_d)\in\R^d\colon
        (\exists\,v\in\R\colon(x_1,\dots,x_d,v)\in E)\}\bigr\}$
      is a cylindrical definable $C^k$ cell decomposition of $\R^d$, and
      \item it holds for all
      $B\in\midcup_{E\in\mathscr P}\bigl\{\{x=(x_1,\dots,x_d)\in\R^d\colon
        (\exists\,v\in\R\colon(x_1,\dots,x_d,v)\in E)\}\bigr\}$
      that there exist $q\in\N_0$ and
      $\zeta_0,\zeta_1,\dots,\zeta_{q+1}\colon B\to[-\infty,\infty]$ such that
      \begin{enumerate}[label={(\Roman*)}]
        \item it holds for all $x\in B$ that
        $-\infty=\zeta_0(x)<\zeta_1(x)<\ldots<\zeta_q(x)<\zeta_{q+1}(x)=\infty$,
        \item it holds for all $i\in\{1,2,\dots,q\}$ that
        $B\ni x\mapsto\zeta_i(x)\in\R$ is definable\cfadd{definition:definable-function},
        \item it holds for all $i\in\{1,2,\dots,q\}$, $\xi\in B$ that there exists
        an open $U\subseteq\R^d$ and
        $F\in C^k(U,\R)$ such that $\xi\in U$ and
        $F|_{B\cap U}=\zeta_i|_{B\cap U}$, and
        \item it holds for all $V\in\mathscr P$ with $V\subseteq B\times\R$ that
        there exists $j\in\{0,1,\dots,q\}$ such that
        \begin{eqn}
          V\in\bigl\{\operatorname{Graph}(\zeta_j),
            \bigl\{(x_1,\dots,x_d,v)\in B\times\R\colon
              \zeta_j(x_1,\dots,x_d)<v<\zeta_{j+1}(x_1,\dots,x_d)\bigr\}\bigr\}\dott
        \end{eqn}
      \end{enumerate}
    \end{enumerate}
  \end{enumerate}
  \cfout[.]
\end{athm}

\begin{athm}{definition}{definition:partition-compatibility}[Compatible definable partitions]
  Let $d\in\N$, let $A\subseteq\R^d$ be a set, and let
  $\mathscr P$ be a definable partition of
  $\R^d$\cfadd{definition:definable-partition}\dott \cfload[.]
  Then we say that $\mathscr P$ is compatible with $A$ if and
  only if there exists $\mathscr A\subseteq\mathscr P$ such that
  $A=\midcup_{B\in\mathscr A}B$\dott
  \cfout[.]
\end{athm}

\begin{athm}{definition}{definition:partition-refinement}[Refinements of definable partitions]
  Let $d\in\N$ and let
  $\mathscr P$ and $\mathscr R$ be definable partitions of
  $\R^d$\cfadd{definition:definable-partition}\dott \cfload[.]
  Then we say that $\mathscr R$ is a refinement of
  $\mathscr P$ if and only if for all $A\in\mathscr R$
  there exists $B\in\mathscr P$ such that $A\subseteq B$\dott
  \cfout[.]
\end{athm}

\subsection{Upper and lower bounds for the energy-derivatives of the back-transformations}\label{subsec: energy coordinates}

The next result, \cref{lem:divide}, is a key tool for the construction of the
Lyapunov-type function $Q$ for the radial-angular \DE\ system
$(r(t),x(t))$, $t\in I$, in \cref{prop:buildQ} in
\cref{sec:radial-integration} below.
Naturally, the function $h$ in \cref{prop:buildQ} also serves as the Lyapunov function for
the radial-angular \DE\ system and the associated Lyapunov energy
$h(r(t),x(t))$, $t\in I$, is non-increasing over time but the dissipation
associated to the Lyapunov function $h$ is not large enough to directly
provide an estimate for the length of $x(t)$, $t\in I$.
Instead the approach initiated by \cref{lem:divide} is to employ transformations
using the functions $\mathcal F_0,\mathcal F_1,\dots,\mathcal F_{N-1}$ in
\cref{lem:divide} and to consider the transformed Lyapunov-type energy
\begin{eqn}
  \label{eq:transformed-lyapunov-energy}
  \mathcal F_i\bigl(r(t),(\varphi\circ h)(r(t),x(t))\bigr)
\end{eqn}
for $t\in I$, $i\in\{0,1,\dots,N-1\}$, where
$\varphi\colon\R\to\R$ is the normalization function which satisfies for
all $e\in\R$ that $\varphi(e)=(1+e^2)^{-1/2}e$.
The functions $E_1,E_2,\dots,E_N$ serve as the back-transformations in the
sense that for all $i\in\{1,2,\dots,N\}$, $r\in(0,\rho)$,
$z\in(0,H)$ it holds that $\mathcal F_{i-1}(r,E_i(r,z))=z$
(see \itref{lem:divide}{item5}).

\begin{athm}{lemma}{lem:divide}
  Let $D\in\N$, $R,c\in(0,\infty)$, let $M\subseteq\R^D$ be a
  nonempty, compact, and definable \submanifold{2}{D},
  let $\Omega\subseteq(0,R)\times M$ be definable and open in $(0,R)\times M$,
  let $h\in \COne{\Omega}{\R}$ be
  definable\cfadd{definition:c1-on-manifold}, let $a,b\colon(0,R)\to\R$
  be definable, assume for all $r\in(0,R)$
  that $a(r)<b(r)$, assume for all
  $(r,x)\in\Omega$ that
  \begin{eqn}
    \llabel{ass1}
    a(r)\leq h(r,x)\leq b(r),
  \end{eqn}
  and let
  $G\colon\Omega\to\R$ satisfy for all
  $(r,x)\in\Omega$ that
  $G(r,x)=\norm{(\nabla_Mh)(r,x)}$\cfadd{definition:angular-gradient}\dott
  \cfload[.]
  Then there exist $\rho\in(0,R)$, $N\in\N$, $H\in(0,\infty)$ and
  definable 
  $\mathfrak e_0,\mathfrak e_1,\dots,\mathfrak e_N\in C^2((0,\rho),\R)$,
  $\mathcal F_0,\mathcal F_1,\dots,\mathcal F_{N}
    \in C((0,\rho)\times\R,[0,H])$,
  $E_1,E_2,\dots,E_N\in C^2((0,\rho)\times(0,H),\R)$,
  $T_1,T_2,\dots,T_N\in C^1((0,\rho)\times(0,H),(0,\infty))$
  such that
  \begin{enumerate}[label=(\roman*)]
    \item \llabel{item1} for all $r\in(0,\rho)$ it holds that
    $a(r)=\mathfrak e_0(r)<\mathfrak e_1(r)<\ldots<\mathfrak e_N(r)=b(r)$,
    \item \llabel{item2} for all $r\in(0,\rho)$, $e\in\R$ with
    $\inf\bigl(\{\infty\}\cup
      \{z\in\R\colon\exists\,x\in M\colon
        (r,x)\in\Omega\cap h^{-1}(\{e\})\cap G^{-1}(\{z\})\}\bigr)=0$
    it holds that $e\in\{\mathfrak e_0(r),\mathfrak e_1(r),\dots,\mathfrak e_N(r)\}$,
    \item \llabel{item3} for all $i\in\{0,1,\dots,N-1\}$,
    $r\in(0,\rho)$, $u\in(-\infty,\mathfrak e_i(r)]$,
    $v\in[\mathfrak e_{i+1}(r),\infty)$,
    $z,w\in(\mathfrak e_i(r),\mathfrak e_{i+1}(r))$ with $z<w$ it holds that
    $0=\mathcal F_i(r,u)<\mathcal F_i(r,z)<\mathcal F_i(r,w)<\mathcal F_i(r,v)=H$,
    \item \llabel{item4} for all $i\in\{0,1,\dots,N-1\}$ it holds
    that
    $\mathcal F_i|_{\{(r,e)\in(0,\rho)\times\R\colon
      \mathfrak e_i(r)<e<\mathfrak e_{i+1}(r)\}}$
    is twice continuously differentiable,
    \item \llabel{item5} for all $i\in\{1,2,\dots,N\}$,
    $r\in(0,\rho)$, $z\in(0,H)$ it holds that
    \begin{eqn}
      \llabel{concl2}
      \mathfrak e_{i-1}(r)<E_i(r,z)<\mathfrak e_i(r),\quad
      \mathcal F_{i-1}(r,E_i(r,z))=z,\qquad\text{and}\qquad
      T_i(r,z)=\partial_z E_i(r,z)>0,
    \end{eqn}
    and
    \item \llabel{item6} for all $i\in\{1,2,\dots,N\}$,
    $(r,x)\in\Omega\cap((0,\rho)\times M)$ with
    $\mathfrak e_{i-1}(r)<h(r,x)<\mathfrak e_i(r)$ it holds that
    \begin{eqn}
      \llabel{concl3}
      T_i(r,\mathcal F_{i-1}(r,h(r,x)))\leq cG(r,x)\dott
    \end{eqn}
  \end{enumerate}
  \cfout[.]
\end{athm}
\begin{aproof}
  Throughout this proof let $\pi_j\colon\R^D\to\R$,
  $j\in\{1,\dots,D\}$, satisfy for all
  $j\in\{1,\dots,D\}$, $x=(x_1,\dots,x_D)\in\R^D$ that
  $\pi_j(x)=x_j$ and let
  $B\in[0,\infty)$, $H\in(0,\infty)$, $\mathcal D\subseteq(0,R)\times\R$
  satisfy
  \begin{equation}
    \llabel{setting1}
    \begin{split}
    &B=\sum_{j=1}^D\biggl(\max_{x\in M}\pi_j(x)-\min_{x\in M}\pi_j(x)\biggr),
      \qquad H=1+2B/c,\\
    \text{and}\qquad&\mathcal D
      =\bigl\{(r,e)\in(0,R)\times\R\colon
          \bigl[\exists\,x\in M\colon((r,x)\in\Omega)\wedge(h(r,x)=e)\bigr]\bigr\}\dott
    \end{split}
  \end{equation}
  \startnewargseq
  We claim that for all $u,v\in\R$ with $u<v$ and every definable
  function $\psi\colon(u,v)\times(0,1)\to M$ there exists
  $A\in(0,\infty)$ such that for all $\delta\in(0,1)$ there exists a finite set $\mathcal E\subseteq(u,v)$ such that
  $\psi(\cdot,\delta)|_{(u,v)\backslash\mathcal E}\in C^1((u,v)\backslash\mathcal E,M)$ and that
  \begin{eqn}
    \llabel{eq1}
    \int_{(u,v)\backslash\mathcal E}
      \norm[\big]{\bigl(\partial_e \psi\bigr)\bigl(e,\delta\bigr)}\dd e
      \leq A\dott
  \end{eqn}
  \startnewargseq
  To prove this claim, let $u,v\in\R$ satisfy $u<v$, and let
  $\psi\colon(u,v)\times(0,1)\to M$ be definable\dott
  \startnewargseq
  In the following let
  $\phi_j\colon(0,1)\times(u,v)\to\R$, $j\in\{1,\dots,D\}$,
  satisfy for all $j\in\{1,\dots,D\}$,
  $(\delta,e)\in(0,1)\times(u,v)$ that
  \begin{eqn}
    \llabel{setting1.1}
    \phi_j(\delta,e)=\pi_j(\psi(e,\delta))\dott
  \end{eqn}
  \startnewargseq
  \argument{
    \lref{setting1.1};
    the fact that $\psi$ is definable;
    \cite[Theorem~6.7]{Coste};
  }{
    that there exist definable partitions $\mathscr V_j$,
    $j\in\{1,\dots,D\}$, of $(0,1)\times(u,v)$\cfadd{definition:definable-partition}
    which satisfy that for all $j\in\{1,\dots,D\}$, $S\in\mathscr V_j$
    it holds that $S$ is a $C^1$ submanifold of $\R^2$\cfadd{definition:submanifold}
    and $\phi_j|_S\in \COne{S}{\R}$\dott\llabel{arg0.0}
  }
  \argument{
    \lref{arg0.0};
    \cite[Theorem~6.6]{Coste};
  }{
    that there exists a cylindrical definable $C^2$ cell decomposition
    $\mathscr P$ of $\R^2$\cfadd{definition:planar-cell-decomposition} which is
    compatible\cfadd{definition:partition-compatibility} with $(0,1)\times(u,v)$ and which satisfies for all
    $j\in\{1,\dots,D\}$, $S\in\mathscr V_j$ that
    \begin{eqn}
      \llabel{eq0.1}
      S=\midcup_{C\in\mathscr P,\,C\subseteq S}C\dott
    \end{eqn}
  }
  In the following let $\mathcal U\subseteq(0,1)\times(u,v)$ satisfy
  \begin{eqn}
    \llabel{setting1.2}
    \mathcal U=\bigl\{(\delta,e)\in(0,1)\times(u,v)\colon
      \bigl[\exists\,C\in\mathscr P,\ \eta\in(0,\infty)\colon
      \{\delta\}\times(e-\eta,e+\eta)\subseteq C\bigr]\bigr\}\dott
  \end{eqn}
  \startnewargseq
  \argument{
    \lref{eq0.1};
    the fact that for all $j\in\{1,\dots,D\}$ it holds that
    $\mathscr V_j$ is a definable partition of $(0,1)\times(u,v)$\cfadd{definition:definable-partition};
    the fact that for all $j\in\{1,\dots,D\}$, $S\in\mathscr V_j$
    it holds that $\phi_j|_S\in \COne{S}{\R}$;
  }{
    that for all $C\in\mathscr P$, $j\in\{1,\dots,D\}$ with
    $C\subseteq(0,1)\times(u,v)$ there exists
    $S\in\mathscr V_j$ such that $C\subseteq S$ and that $\phi_j|_C\in \COne{C}{\R}$\dott\llabel{arg0.1}
  }
  \argument{
    \lref{arg0.1};
    \lref{setting1.2};
    the fact that $\mathscr P$ is a definable partition of $\R^2$\cfadd{definition:definable-partition};
    the fact that $\mathscr P$ is compatible\cfadd{definition:partition-compatibility} with $(0,1)\times(u,v)$;
    the chain rule;
    \itref{lem:radial-derivative-well-defined}{item2};
  }{
    that $\mathcal U$ is definable and that there exist definable
     $d_j\colon\mathcal U\to\R$, $j\in\{1,\dots,D\}$,
    which satisfy for all $j\in\{1,\dots,D\}$,
    $(\delta,e)\in\mathcal U$ that
    \begin{eqn}
      \llabel{eq0.2}
      d_j(\delta,e)=\bigl(\partial_e \phi_j\bigr)\bigl(\delta,e\bigr)\dott
    \end{eqn}
  }
  In the following let $\mathcal L_j,\mathcal Z_j,\mathcal W_j
  \subseteq\mathcal U$, $j\in\{1,\dots,D\}$, satisfy for all
  $j\in\{1,\dots,D\}$ that
  \begin{gather}
    \llabel{setting1.3}
    \mathcal L_j=\{(\delta,e)\in\mathcal U\colon d_j(\delta,e)<0\},
      \qquad
      \mathcal Z_j=\{(\delta,e)\in\mathcal U\colon d_j(\delta,e)=0\},\\
    \text{and}\qquad\mathcal W_j
      =\{(\delta,e)\in\mathcal U\colon d_j(\delta,e)>0\}\dott
  \end{gather}
  \startnewargseq
  \argument{
    \lref{setting1.3};
    the fact that $\mathcal U$ is definable;
    the fact that for all $j\in\{1,\dots,D\}$ it holds that
    $d_j$ is definable;
    the fact that $\mathscr P$ is a definable partition of $\R^2$\cfadd{definition:definable-partition};
    \cite[Theorem~6.6]{Coste};
  }{
    that there exists a cylindrical definable $C^2$ cell decomposition
    $\mathscr S$ of $\R^2$ which is a refinement\cfadd{definition:partition-refinement} of $\mathscr P$,
    which is compatible\cfadd{definition:partition-compatibility} with $\mathcal U$, and which satisfies for all
    $j\in\{1,\dots,D\}$, $S\in\{\mathcal L_j,\mathcal Z_j,\mathcal W_j\}$ that
    \begin{eqn}
      \llabel{eq0.3}
      S=\midcup_{C\in\mathscr S,\,C\subseteq S}C\dott
    \end{eqn}
  }
  \startnewargseq
  \argument{  \lref{setting1.3};
    \lref{eq0.3};
    the fact that $\mathscr S$ is a definable partition of $\R^2$\cfadd{definition:definable-partition};
  }{
    that for all $C\in\mathscr S$, $j\in\{1,\dots,D\}$,
    $\xi,\eta\in C$ with $C\subseteq\mathcal U$ it holds that
    \begin{eqn}
      \llabel{eq0.4}
      \operatorname{sign}\bigl(d_j(\xi)\bigr)
        =\operatorname{sign}\bigl(d_j(\eta)\bigr)\dott
    \end{eqn}
  }
  \argument{
    \lref{setting1.2};
    the fact that $\mathscr S$ is a refinement\cfadd{definition:partition-refinement} of $\mathscr P$;
    the fact that $\mathscr S$ is compatible\cfadd{definition:partition-compatibility} with $\mathcal U$;
  }{
    that for all $C\in\mathscr S$, $\delta\in(0,1)$,
    $p,q\in(u,v)$ with $C\subseteq(0,1)\times(u,v)$,
    $\{\delta\}\times(p,q)\subseteq C$, and $p<q$ it holds that
    \begin{eqn}
      \llabel{eq0.5}
      C\subseteq\mathcal U\dott
    \end{eqn}
  }
  \argument{
    \lref{eq0.5};
    \lref{setting1.1};
    \lref{eq0.2};
    \lref{eq0.4};
    the fact that for all $C\in\mathscr P$, $j\in\{1,\dots,D\}$
    with $C\subseteq(0,1)\times(u,v)$ it holds that $\phi_j|_C\in \COne{C}{\R}$;
    the fact that $\mathscr S$ is a refinement\cfadd{definition:partition-refinement} of $\mathscr P$;
    the fact that for all $\delta\in(0,1)$, $C\in\mathscr S$ it holds
    that $\{e\in\R\colon(\delta,e)\in C\}$ is empty, a singleton,
    or an open interval;
    the fact that $\mathscr S$ is compatible\cfadd{definition:partition-compatibility} with
    $(0,1)\times(u,v)$;
    the mean value theorem;
  }{
    that for all $\delta\in(0,1)$ there exist
    $k\in\{1,\dots,\#\mathscr S\}$ and
    $t=(t_0,\dots,t_k)\in\R^{k+1}$ with $u=t_0<\cdots<t_k=v$
    such that for all $\ell\in\{1,\dots,k\}$,
    $j\in\{1,\dots,D\}$ it holds that
    $\psi(\cdot,\delta)|_{(t_{\ell-1},t_\ell)}\in C^1((t_{\ell-1},t_\ell),M)$ and that
    $\phi_j(\delta,\cdot)|_{(t_{\ell-1},t_\ell)}$ is non-decreasing
    or non-increasing\dott
  }
  In the following let $A\in(0,\infty)$ satisfy
  $A=1+B\,(\#\mathscr S)$\dott
  \startnewargseq
We claim that for all $\delta\in(0,1)$ there exists
  a finite set $\mathcal E\subseteq(u,v)$ such that
  $\psi(\cdot,\delta)|_{(u,v)\backslash\mathcal E}\in C^1((u,v)\backslash\mathcal E,M)$ and that
  \begin{eqn}
    \llabel{eq1.0}
    \int_{(u,v)\backslash\mathcal E}
      \norm[\big]{\bigl(\partial_e \psi\bigr)\bigl(e,\delta\bigr)}\dd e
      \leq A\dott
  \end{eqn}
  \startnewargseq
  To prove \lref{eq1.0}, let $\delta\in(0,1)$,
  $k\in\{1,\dots,\#\mathscr S\}$, and
  $t=(t_0,\dots,t_k)\in\R^{k+1}$, assume that
  $u=t_0<\cdots<t_k=v$ and that for all
  $\ell\in\{1,\dots,k\}$, $j\in\{1,\dots,D\}$ it holds that
  $\psi(\cdot,\delta)|_{(t_{\ell-1},t_\ell)}\in C^1((t_{\ell-1},t_\ell),M)$ and that
  $\phi_j(\delta,\cdot)|_{(t_{\ell-1},t_\ell)}$ is non-decreasing
  or non-increasing, and let $\mathcal E=\{t_1,\dots,t_{k-1}\}$\dott
  \startnewargseq
  \argument{
    the fact that for all $\ell\in\{1,\dots,k\}$,
    $j\in\{1,\dots,D\}$ it holds that
    $\phi_j(\delta,\cdot)|_{(t_{\ell-1},t_\ell)}\in C^1((t_{\ell-1},t_\ell),\R)$
    and that $\phi_j(\delta,\cdot)|_{(t_{\ell-1},t_\ell)}$ is
    either non-decreasing or non-increasing;
    the fundamental theorem of calculus;
    \lref{setting1.1};
    the fact that $\psi$ takes values in $M$;
  }{
    that for all $\ell\in\{1,\dots,k\}$, $j\in\{1,\dots,D\}$,
    $p,q\in(t_{\ell-1},t_\ell)$ with $p<q$ it holds that
    \begin{eqn}
      \llabel{eq1.1}
      \int_p^q\bigl|\bigl(\partial_e \phi_j\bigr)\bigl(\delta,e\bigr)\bigr|\dd e
        =\abs{\phi_j(\delta,q)-\phi_j(\delta,p)}
        \leq\max_{x\in M}\pi_j(x)-\min_{x\in M}\pi_j(x)\dott
    \end{eqn}
  }
  \argument{
    \lref{eq1.1};
    the monotone convergence theorem
  }{
    that for all $\ell\in\{1,\dots,k\}$, $j\in\{1,\dots,D\}$
    it holds that
    \begin{eqn}
      \llabel{eq1.2}
      \int_{t_{\ell-1}}^{t_\ell}
        \bigl|\bigl(\partial_e \phi_j\bigr)\bigl(\delta,e\bigr)\bigr|\dd e
        \leq\max_{x\in M}\pi_j(x)-\min_{x\in M}\pi_j(x)\dott
    \end{eqn}
  }
  \argument{
    \lref{setting1.1};
    the fact that for all $\ell\in\{1,\dots,k\}$ it holds that
    $\psi(\cdot,\delta)|_{(t_{\ell-1},t_\ell)}\in C^1((t_{\ell-1},t_\ell),M)$;
  }{
    that for all $e\in(u,v)\backslash\mathcal E$ it holds that
    \begin{eqn}
      \llabel{eq1.3}
      \norm[\big]{\bigl(\partial_e \psi\bigr)\bigl(e,\delta\bigr)}
        =\sqrt{\sum_{j=1}^D
          \bigl(\bigl(\partial_e \phi_j\bigr)\bigl(\delta,e\bigr)\bigr)^2}
        \leq\sum_{j=1}^D
          \bigl|\bigl(\partial_e \phi_j\bigr)\bigl(\delta,e\bigr)\bigr|\dott
    \end{eqn}
  }
  \argument{
    \lref{eq1.3};
    \lref{setting1};
    \lref{eq1.2};
    the fact that $(u,v)\backslash\mathcal E$
    is the disjoint union of $(t_{\ell-1},t_\ell)$,
    $\ell\in\{1,\dots,k\}$;
    the fact that $k\leq\#\mathscr S$;
  }{
    that
    \begin{eqn}[a]
      \llabel{eq2}
      \int_{(u,v)\backslash\mathcal E}
        \norm[\big]{\bigl(\partial_e \psi\bigr)\bigl(e,\delta\bigr)}\dd e
      &\leq\sum_{\ell=1}^k\sum_{j=1}^D
        \int_{t_{\ell-1}}^{t_\ell}
        \bigl|\bigl(\partial_e \phi_j\bigr)\bigl(\delta,e\bigr)\bigr|\dd e\\
      &\leq\sum_{\ell=1}^k\sum_{j=1}^D
        \biggl(\max_{x\in M}\pi_j(x)-\min_{x\in M}\pi_j(x)\biggr)
        =kB\leq B\,(\#\mathscr S)<A\dott
    \end{eqn}
  }
  \argument{
    \lref{eq2};
    the fact that $\mathcal E\subseteq(u,v)$ is finite;
    the fact that $\psi(\cdot,\delta)|_{(u,v)\backslash\mathcal E}\in C^1((u,v)\backslash\mathcal E,M)$;
    the fact that $\delta\in(0,1)$ was arbitrary;
  }[verbs=ep]{
    \lref{eq1.0}\dott
  }
  \startnewargseq
  \argument{
    \lref{eq1.0};
    the fact that $A\in(0,\infty)$;
    the fact that $u,v,\psi$ were arbitrary;
  }[verbs=ep]{
    \lref{eq1}\dott
  }

  In the following let $\mu\colon\mathcal D\to[0,\infty)$ and
  $\mathcal C\subseteq\mathcal D$ satisfy for all
  $(r,e)\in\mathcal D$ that
  \begin{eqn}
    \llabel{setting2}
    \mu(r,e)=\inf\bigl\{G(r,x)\colon\bigl[((r,x)\in\Omega)\wedge(h(r,x)=e)\bigr]\bigr\}
    \qquad\text{and}\qquad
    \mathcal C=\{(s,v)\in\mathcal D\colon\mu(s,v)=0\}\dott
  \end{eqn}
  \startnewargseq
  \argument{
    \lref{setting1};
    the fact that $h$ is definable;the fact that $\Omega$ is definable;
  }{
    that
    \begin{eqn}
      \llabel{eq2.1}
      \mathcal D=\{(r,h(r,x))\colon(r,x)\in\Omega\}
    \end{eqn}
    and that $\mathcal D$ is definable\dott
  }
  \argument{
    \lref{eq2.1};
    \lref{setting2};
  }{
    that for all $(r,e)\in\mathcal D$ there exists
    $x\in M$ such that $(r,x)\in\Omega$ and
    \begin{eqn}
      \llabel{eq2.2}
      h(r,x)=e\qquad\text{and}\qquad
      0\leq\mu(r,e)\leq G(r,x)<\infty\dott
    \end{eqn}
  }
  In the following let $\mathcal J\subseteq\mathcal D\times[0,\infty)$
  satisfy
  \begin{eqn}
    \llabel{eq2.3}
    \mathcal J=(\mathcal D\times[0,\infty))\backslash
      \bigl\{(r,h(r,x),z)\colon
        \bigl[((r,x)\in\Omega)\wedge(z\in[0,\infty))\wedge(G(r,x)<z)\bigr]\bigr\}
    \dott
  \end{eqn}
  \startnewargseq
  \argument{
    \lref{eq2.3};
    the fact that $\mathcal D$ and $\Omega$ are definable; the fact that $h$ and $G$ are definable;
  }{
    that $\mathcal J$ is definable\dott
  }
  \argument{
    \lref{eq2.3};
  }{
    that for all $(r,e,z)\in\mathcal J$, $x\in M$ with
    $(r,x)\in\Omega$ and $h(r,x)=e$ it holds that
    \begin{eqn}
      \llabel{eq2.4}
      z\leq G(r,x)\dott
    \end{eqn}
  }
  \argument{
    \lref{eq2.4};
    \lref{setting2};
  }{
    that for all $(r,e,z)\in\mathcal J$ it holds that
    \begin{eqn}
      \llabel{eq2.5}
      z\leq\mu(r,e)\dott
    \end{eqn}
  }
  \argument{
    \lref{setting2};
    \lref{eq2.2};
    \lref{eq2.3};
  }{
    that for all $(r,e)\in\mathcal D$ it holds that
    \begin{eqn}
      \llabel{eq2.6}
      (r,e,\mu(r,e))\in\mathcal J\dott
    \end{eqn}
  }
  In the following let $\Gamma\subseteq\mathcal J$ satisfy
  \begin{eqn}
    \llabel{eq2.7}
    \Gamma=\mathcal J\backslash
      \bigl\{(r,e,z)\in\mathcal J\colon
        \bigl[\exists\,v\in[0,\infty)\colon((r,e,v)\in\mathcal J)\wedge(z<v)\bigr]\bigr\}
    \dott
  \end{eqn}
  \startnewargseq
  \argument{
    \lref{eq2.7};
    the fact that $\mathcal J$ is definable;
  }{
    that $\Gamma$ is definable\dott
  }
  \argument{
    \lref{eq2.5};
    \lref{eq2.6};
    \lref{eq2.7};
  }{
    that for all $(r,e)\in\mathcal D$ it holds that
    \begin{eqn}
      \llabel{eq2.8}
      (r,e,\mu(r,e))\in\Gamma\dott
    \end{eqn}
  }
  \argument{
    \lref{eq2.6};
    \lref{eq2.7};
    the fact that for all $(r,e,z)\in\Gamma$, $v\in[0,\infty)$
    with $(r,e,v)\in\mathcal J$ it holds that $v\leq z$;
  }{
    that for all $(r,e,z)\in\Gamma$ it holds that
    \begin{eqn}
      \llabel{eq2.9}
      \mu(r,e)\leq z\dott
    \end{eqn}
  }
  \argument{
    \lref{eq2.9};
    \lref{eq2.5};
    \lref{eq2.8};
    the fact that $\Gamma\subseteq\mathcal J$;
  }{
    that
    \begin{eqn}
      \llabel{eq2.10}
      \Gamma=\{(r,e,\mu(r,e))\colon(r,e)\in\mathcal D\}\dott
    \end{eqn}
  }
  \argument{
    \lref{eq2.10};
    the fact that $\Gamma$ is definable;
  }{
    that $\mu$ is definable\dott\llabel{arg2.1}
  }
  \argument{
    \lref{arg2.1};
    \lref{setting2};
  }{
    that $\mathcal C=\mu^{-1}(\{0\})$ is definable\dott
  }
  We claim that for all $r\in(0,R)$ it holds that
  \begin{eqn}\llabel{eq3}
    \#\{e\in\R\colon(r,e)\in\mathcal C\}<\infty\dott
  \end{eqn}
  \startnewargseq
  To prove \lref{eq3}, let $r\in(0,R)$\dott
  \startnewargseq
  We prove \lref{eq3} by contradiction\dott
  \startnewargseq
  We thus assume that
  $\#\{e\in\R\colon(r,e)\in\mathcal C\}=\infty$\dott
  \startnewargseq
  \argument{
    the fact that $\mathcal C$ is definable;
    the fact that $r\in(0,R)$;
  }{
    that $\{e\in\R\colon(r,e)\in\mathcal C\}$ is definable\dott
    \llabel{arg3.1}
  }
  \argument{
    \lref{arg3.1};
    the assumption that $\#\{e\in\R\colon(r,e)\in\mathcal C\}=\infty$;
    the fact that definable subsets of $\R$ are finite unions of points
    and intervals;
  }{
    that there exist $u,v\in\R$ with $u<v$ which satisfy
    \begin{eqn}
      \llabel{eq3.1}
      (u,v)\subseteq\{e\in\R\colon(r,e)\in\mathcal C\}\dott
    \end{eqn}
  }
  \startnewargseq
  \argument{ \lref{setting2};
    \lref{eq3.1};
  }{
    that for all $e\in(u,v)$ it holds that
    \begin{eqn}
      \llabel{eq3.2}
      (r,e)\in\mathcal D\qquad\text{and}\qquad\mu(r,e)=0\dott
    \end{eqn}
  }
  \argument{
    \lref{eq3.2};
    \lref{eq2.5};
  }{
    that for all $(e,\delta)\in(u,v)\times(0,1)$ it holds that
    \begin{eqn}
      \llabel{eq3.3}
      (r,e,\delta)\notin\mathcal J\dott
    \end{eqn}
  }
  \argument{
    \lref{eq3.3};
    \lref{eq2.3};
    \lref{eq3.2};
  }{
    that for all $(e,\delta)\in(u,v)\times(0,1)$ 
    there exists $x\in M$ such that
    \begin{eqn}
      \llabel{eq3.4}
      (r,x)\in\Omega,\qquad h(r,x)=e,\qquad\text{and}\qquad
      G(r,x)<\delta\dott
    \end{eqn}
  }
  In the following let $\mathcal H\subseteq(u,v)\times(0,1)\times M$
  satisfy
  \begin{eqn}
    \llabel{eq3.5}
    \mathcal H=\bigl\{(e,\delta,x)\in(u,v)\times(0,1)\times M\colon
      \bigl[((r,x)\in\Omega)\wedge(h(r,x)=e)\wedge(G(r,x)<\delta)\bigr]\bigr\}\dott
  \end{eqn}
  \startnewargseq
  \argument{
    \lref{eq3.4};
    \lref{eq3.5};
    the fact that $M$, $\Omega$, $h$, and $G$ are definable;
  }{
    that $\mathcal H$ is definable and that
    \begin{eqn}
      \llabel{eq3.6}
      \bigl\{(e,\delta)\in(u,v)\times(0,1)\colon
        \bigl[\exists\,x\in M\colon(e,\delta,x)\in\mathcal H\bigr]\bigr\}
        =(u,v)\times(0,1)\dott
    \end{eqn}
  }
  \argument{
    \lref{eq3.6};
    \lref{eq3.5};
    the fact that $\mathcal H$ is definable;
   \cite[Theorem~3.1]{Coste};
  }{
    that there exists a definable function
    $\chi\colon(u,v)\times(0,1)\to M$ which satisfies for all
    $(e,\delta)\in(u,v)\times(0,1)$ that
    \begin{eqn}\llabel{eq4}
      (r,\chi(e,\delta))\in\Omega,\quad
      h(r,\chi(e,\delta))=e,\qquad\text{and}\qquad
      G(r,\chi(e,\delta))<\delta\dott
    \end{eqn}
  }
  \startnewargseq
  \argument{
    \lref{eq1};
    the fact that $\chi$ is definable;
  }{
    that there exists $A\in(0,\infty)$ which satisfies that for all
    $\delta\in(0,1)$ there exists a finite set
    $\mathcal E\subseteq(u,v)$ such that
    $\chi(\cdot,\delta)|_{(u,v)\backslash\mathcal E}\in C^1((u,v)\backslash\mathcal E,M)$ and
    \begin{eqn}
      \llabel{eq4.1}
      \int_{(u,v)\backslash\mathcal E}
        \norm[\big]{\bigl(\partial_e \chi\bigr)\bigl(e,\delta\bigr)}\dd e
        \leq A\dott
    \end{eqn}
  }
  In the following let $\delta\in(0,1)$ satisfy
  \begin{eqn}
    \llabel{eq4.1.1}
    \delta=\min\biggl\{\frac12,\frac{v-u}{2A}\biggr\},
  \end{eqn}
  and let $\mathcal E\subseteq(u,v)$ be finite with the property that
  $\chi(\cdot,\delta)|_{(u,v)\backslash\mathcal E}\in C^1((u,v)\backslash\mathcal E,M)$
  and with the property that
  \begin{eqn}
    \llabel{eq4.1.2}
    \int_{(u,v)\backslash\mathcal E}
      \norm[\big]{\bigl(\partial_e \chi\bigr)\bigl(e,\delta\bigr)}\dd e
      \leq A\dott
  \end{eqn}
  \startnewargseq
  \argument{
    \lref{eq4};
    the fact that $\mathcal E$ is finite;
    the fact that $\chi(\cdot,\delta)|_{(u,v)\backslash\mathcal E}\in C^1((u,v)\backslash\mathcal E,M)$;
    the fact that $h\in \COne{\Omega}{\R}$;
    the fact that for all $(s,x)\in\Omega$ it holds that
    $G(s,x)=\norm{(\nabla_Mh)(s,x)}$;
    the chain rule;
    the Cauchy--Schwarz inequality;
  }{
    that for all $e\in(u,v)\backslash\mathcal E$ it holds that
    \begin{eqn}[a]
      \llabel{eq5.1}
      1&=\frac{\ddd}{\ddd e}h(r,\chi(e,\delta))
        =\bigl\langle(\nabla_Mh)(r,\chi(e,\delta)),
          \bigl(\partial_e \chi\bigr)\bigl(e,\delta\bigr)\bigr\rangle\\
       &\leq G(r,\chi(e,\delta))
        \norm[\big]{\bigl(\partial_e \chi\bigr)\bigl(e,\delta\bigr)}
        \leq\delta\norm[\big]{\bigl(\partial_e \chi\bigr)\bigl(e,\delta\bigr)}
      \dott
    \end{eqn}
  }
  \argument{
    \lref{eq5.1};
    \lref{eq4.1.2};
    the fact that $\mathcal E$ is finite;
  }{
    that
    \begin{eqn}
      \llabel{eq5.2}
      v-u=\int_{(u,v)\backslash\mathcal E}1\dd e
        \leq\delta\int_{(u,v)\backslash\mathcal E}
          \norm[\big]{\bigl(\partial_e \chi\bigr)\bigl(e,\delta\bigr)}\dd e
        \leq\delta A\dott
    \end{eqn}
  }
  \argument{
    \lref{eq5.2};
    \lref{eq4.1.1};
    the fact that $v-u>0$;
  }{
    that \begin{eqn}[i]v-u\leq\delta A\leq\frac{v-u}{2}<v-u\dott\end{eqn}
  }
  This contradiction and the fact that $r\in(0,R)$ was arbitrary
  prove \lref{eq3}\dott
  \startnewargseq

  \argument{
    \lref{setting2};
    \lref{eq2.2};
  }{
    that for all $(r,e)\in\mathcal D\backslash\mathcal C$ it holds that
    \begin{eqn}
      \llabel{eq5.3}
      0<\mu(r,e)<2\mu(r,e)<\infty\dott
    \end{eqn}
  }
  \argument{
    \lref{eq5.3};
    \lref{eq2.5};
  }{
    that for all $(r,e)\in\mathcal D\backslash\mathcal C$ it holds that
    \begin{eqn}
      \llabel{eq5.4}
      (r,e,2\mu(r,e))\notin\mathcal J\dott
    \end{eqn}
  }
  \argument{
    \lref{eq5.4};
    \lref{eq2.3};
    \lref{eq5.3};
  }{
    that for all $(r,e)\in\mathcal D\backslash\mathcal C$ 
    there exists $x\in M$ such that
    \begin{eqn}
      \llabel{eq5.5}
      (r,x)\in\Omega,\qquad h(r,x)=e,\qquad\text{and}\qquad
      G(r,x)<2\mu(r,e)\dott
    \end{eqn}
  }
  In the following let
  $\mathcal K\subseteq(\mathcal D\backslash\mathcal C)\times M$ satisfy
  \begin{eqn}
    \llabel{eq5.6}
    \mathcal K=\bigl\{(r,e,x)\in(\mathcal D\backslash\mathcal C)\times M\colon
      \bigl[((r,x)\in\Omega)\wedge(h(r,x)=e)\wedge(G(r,x)<2\mu(r,e))\bigr]\bigr\}\dott
  \end{eqn}
  \startnewargseq
  \argument{
    \lref{eq5.5};
    \lref{eq5.6};
    the fact that $\mathcal D$, $\mathcal C$, $M$, and $\Omega$ are definable;
    the fact that $h$, $G$, and $\mu$ are definable;
  }{
    that $\mathcal K$ is definable and that
    \begin{eqn}
      \llabel{eq5.7}
      \bigl\{(r,e)\in\mathcal D\backslash\mathcal C\colon
        \bigl[\exists\,x\in M\colon(r,e,x)\in\mathcal K\bigr]\bigr\}
        =\mathcal D\backslash\mathcal C\dott
    \end{eqn}
  }
  \argument{
    \lref{eq5.7};
    \lref{eq5.6};
    the fact that $\mathcal K$ is definable;
    \cite[Theorem~3.1]{Coste};
  }{
    that there exists a definable function
    $\zeta\colon\mathcal D\backslash\mathcal C\to M$ which satisfies for all
    $(r,e)\in\mathcal D\backslash\mathcal C$ that
    \begin{eqn}\llabel{eq6}
      (r,\zeta(r,e))\in\Omega,\quad h(r,\zeta(r,e))=e,
      \qquad\text{and}\qquad G(r,\zeta(r,e))<2\mu(r,e)\dott
    \end{eqn}
  }
  \startnewargseq
  \argument{
    \lref{eq6};
    the fact that $\zeta$ is definable;
    the fact that for all $j\in\{1,\dots,D\}$ it holds that $\pi_j$ is definable;
    \cite[Theorem~6.7]{Coste};
  }{
    that there exist definable partitions $\mathscr X_j$,
    $j\in\{1,\dots,D\}$, of $\mathcal D\backslash\mathcal C$\cfadd{definition:definable-partition}
    which satisfy that for all $j\in\{1,\dots,D\}$, $S\in\mathscr X_j$
    it holds that $S$ is a $C^2$ submanifold of $\R^2$\cfadd{definition:submanifold}
    and $(\pi_j\circ\zeta)|_S\in\CManifold{2}{S}{\R}$\dott
  }
  \startnewargseq
  \argument{
    the fact that $a$ and $b$ are definable;
     \cite[Corollary~6.5]{Coste};
  }{
    that there exists $\kappa\in(0,R)$ which satisfies that
    $a|_{(0,\kappa)},b|_{(0,\kappa)}\in C^2((0,\kappa),\R)$\dott\llabel{arg1}
  }
  \argument{
    \lref{arg1};
    the fact that $\mathcal D$ and $\mathcal C$ are definable;
    the fact that $a$ and $b$ are definable;
    the fact that for all $j\in\{1,\dots,D\}$ it holds that
    $\mathscr X_j$ is a definable partition of $\mathcal D\backslash\mathcal C$\cfadd{definition:definable-partition};
    \cite[Theorem~6.6]{Coste};
  }{
    that there exists a cylindrical definable $C^2$ cell decomposition
    $\mathscr R$ of $\R^2$\cfadd{definition:planar-cell-decomposition}
    which is compatible\cfadd{definition:partition-compatibility} with
    $(0,\kappa/2)\times\R$, $\mathcal D$, $\mathcal C$,
    $\{(r,a(r))\colon r\in(0,\kappa/2)\}$, and
    $\{(r,b(r))\colon r\in(0,\kappa/2)\}$, and which satisfies for all
    $j\in\{1,\dots,D\}$, $S\in\mathscr X_j$ that
    \begin{eqn}
      \llabel{eq6.1}
      S=\midcup_{C\in\mathscr R,\,C\subseteq S}C\dott
    \end{eqn}
  }
  In the following let $\mathcal O\subseteq\mathcal D\backslash\mathcal C$
  satisfy
  \begin{eqn}
    \llabel{eq6.2}
    \mathcal O=\midcup_{C\in\mathscr R,\,C\subseteq\mathcal D\backslash\mathcal C,\,
      C\text{ open}}C\dott
  \end{eqn}
  \startnewargseq
  \argument{
    \lref{eq6.1};
    \lref{eq6.2};
    the fact that for all $j\in\{1,\dots,D\}$ it holds that
    $\mathscr X_j$ is a definable partition of $\mathcal D\backslash\mathcal C$\cfadd{definition:definable-partition};
    the fact that for all $j\in\{1,\dots,D\}$, $S\in\mathscr X_j$
    it holds that $(\pi_j\circ\zeta)|_S\in\CManifold{2}{S}{\R}$;
  }{
    that $\mathcal O$ is definable and open and that
    $\zeta|_{\mathcal O}\in C^2(\mathcal O,M)$\dott\llabel{arg6.1}
  }
  \argument{
    \lref{arg6.1};
    the fact that $\zeta$ is definable;
  }{
    that there exist definable functions
    $\nu_j\in C^1(\mathcal O,\R)$, $j\in\{1,\dots,D\}$, which satisfy
    for all $j\in\{1,\dots,D\}$, $(r,e)\in\mathcal O$ that
    \begin{eqn}
      \llabel{eq6.3}
      \nu_j(r,e)=\bigl(\partial_e (\pi_j\circ\zeta)\bigr)\bigl(r,e\bigr)\dott
    \end{eqn}
  }
  In the following let $\mathcal Y_{j,s}\subseteq\mathcal O$,
  $j\in\{1,\dots,D\}$, $s\in\{-1,0,1\}$, satisfy for all
  $j\in\{1,\dots,D\}$, $s\in\{-1,0,1\}$ that
  \begin{eqn}
    \llabel{eq6.4}
    \mathcal Y_{j,s}=\{(r,e)\in\mathcal O\colon
      \operatorname{sign}(\nu_j(r,e))=s\}\dott
  \end{eqn}
  \startnewargseq
  \argument{
    \lref{eq6.4};
    the fact that $\mathcal O$ is definable;
    the fact that for all $j\in\{1,\dots,D\}$ it holds that $\nu_j$ is definable;
    the fact that $\mathscr R$ is a definable partition of $\R^2$\cfadd{definition:definable-partition};
    \cite[Theorem~6.6]{Coste};
  }{
    that there exists a cylindrical definable $C^2$ cell decomposition
    $\mathscr T$ of $\R^2$ which is
    a refinement\cfadd{definition:partition-refinement} of $\mathscr R$, which is compatible\cfadd{definition:partition-compatibility} with $\mathcal O$,
    and which satisfies for all $j\in\{1,\dots,D\}$,
    $s\in\{-1,0,1\}$ that
    \begin{eqn}
      \llabel{eq6.5}
      \mathcal Y_{j,s}=\midcup_{C\in\mathscr T,\,C\subseteq\mathcal Y_{j,s}}C\dott
    \end{eqn}
  }
  \argument{
    \lref{eq6.5};
    \lref{eq6.4};
    the fact that $\mathscr T$ is a definable partition of $\R^2$\cfadd{definition:definable-partition};
  }{
    that for all $C\in\mathscr T$, $j\in\{1,\dots,D\}$,
    $\xi,\eta\in C$ with $C\subseteq\mathcal O$ it holds that
    \begin{eqn}
      \llabel{eq6.6}
      \operatorname{sign}(\nu_j(\xi))
        =\operatorname{sign}(\nu_j(\eta))\dott
    \end{eqn}
  }
  \argument{
    \lref{eq6.2};
    the fact that $\mathscr T$ is a refinement\cfadd{definition:partition-refinement} of $\mathscr R$;
    the fact that $\mathscr R$ is compatible\cfadd{definition:partition-compatibility} with $\mathcal D$ and $\mathcal C$;
    the fact that for all $C\in\mathscr R$ with nonempty interior in $\R^2$
    it holds that $C$ is open;
  }{
    that for all $C\in\mathscr T$ with
    $C\subseteq\mathcal D\backslash\mathcal C$ and $C$ open it holds that
    \begin{eqn}
      \llabel{eq6.7}
      C\subseteq\mathcal O\dott
    \end{eqn}
  }
  \argument{
    \lref{eq3};
    the fact that $\mathscr T$ is compatible\cfadd{definition:partition-compatibility} with $\mathcal C$;
    the fact that for all $C\in\mathscr T$ with $C\subseteq\mathcal C$
    and $C$ open it holds that there exist $r\in(0,R)$, $p,q\in\R$
    with $p<q$ and $\{r\}\times(p,q)\subseteq C$;
  }{
    that for all open $C\in\mathscr T$ it holds that
    \begin{eqn}
      \llabel{eq6.8}
      C\cap\mathcal C=\emptyset\dott
    \end{eqn}
  }
  \argument{
    \lref{ass1};
    \lref{setting1};
    \lref{setting2};
  }{
    that
    \begin{eqn}
      \llabel{eq6.9}
      \mathcal C\subseteq\mathcal D\subseteq
        \{(r,e)\in(0,R)\times\R\colon a(r)\leq e\leq b(r)\}\dott
    \end{eqn}
  }
  \argument{
    the fact that $\mathscr T$ is a cylindrical definable $C^2$ cell decomposition of $\R^2$;
  }{
    that there exists a definable partition $\mathscr I$ of $\R$\cfadd{definition:definable-partition}
    which satisfies that for all $J\in\mathscr I$ it holds that
    $J$ is a singleton or an open interval and that
    \begin{eqn}
      \llabel{eq6.10}
      \mathscr I=\bigl\{\{r\in\R\colon
        \bigl[\exists\,e\in\R\colon(r,e)\in C\bigr]\}\colon C\in\mathscr T\bigr\}\dott
    \end{eqn}
  }
  \argument{
    \lref{eq6.10};
    the fact that $\mathscr R$ is compatible\cfadd{definition:partition-compatibility} with $(0,\kappa/2)\times\R$;
    the fact that $\mathscr T$ is a refinement\cfadd{definition:partition-refinement} of $\mathscr R$;
  }{
    that
    \begin{eqn}
      \llabel{eq6.11}
      \{0\}\in\mathscr I\qquad\text{and}\qquad\{\kappa/2\}\in\mathscr I\dott
    \end{eqn}
  }
  \argument{
    \lref{eq6.11};
    the fact that $\mathscr I$ is a definable partition of $\R$\cfadd{definition:definable-partition};
    the fact that for all $J\in\mathscr I$ it holds that
    $J$ is a singleton or an open interval;
  }{
    that there exists $\varepsilon\in(0,\kappa)$ which satisfies
    \begin{eqn}
      \llabel{eq6.12}
      \varepsilon=\min\{s\in(0,\infty)\colon\{s\}\in\mathscr I\}\leq\kappa/2
      \qquad\text{and}\qquad(0,\varepsilon)\in\mathscr I\dott
    \end{eqn}
  }
  \argument{
    \lref{eq6.12};
    \lref{eq6.10};
    the fact that $\mathscr T$ is a cylindrical definable $C^2$ cell decomposition of $\R^2$;
    the fact that $\mathscr T$ is compatible\cfadd{definition:partition-compatibility} with
    $\{(r,a(r))\colon r\in(0,\kappa/2)\}$ and
    $\{(r,b(r))\colon r\in(0,\kappa/2)\}$;
    the fact that for all $r\in(0,R)$ it holds that $a(r)<b(r)$;
  }{
    that there exist $m\in\N$ with $m\geq2$ and definable functions
    $\mathbf{f}_j\in C^2((0,\varepsilon),\R)$, $j\in\{1,\dots,m\}$,
    which satisfy for all $r\in(0,\varepsilon)$ that
    \begin{eqn}
      \llabel{eq6.13}
      \mathbf{f}_1(r)<\cdots<\mathbf{f}_m(r),
    \end{eqn}
    and which satisfy
    \begin{eqn}
      \llabel{eq6.14}
      \bigl\{C\in\mathscr T\colon
        \bigl[(C\subseteq(0,\varepsilon)\times\R)\wedge(C\text{ is not open})\bigr]\bigr\}
        =\bigl\{\{(r,\mathbf{f}_j(r))\colon r\in(0,\varepsilon)\}\colon
          j\in\{1,\dots,m\}\bigr\}\dott
    \end{eqn}
  }
  \argument{
    \lref{eq6.14};
    \lref{eq6.13};
    the fact that $\mathscr T$ is a cylindrical definable $C^2$ cell decomposition of $\R^2$;
  }{
    that for all $j\in\{1,\dots,m-1\}$ it holds that
    \begin{eqn}
      \llabel{eq6.15}
      \{(r,e)\in(0,\varepsilon)\times\R\colon
        \mathbf{f}_j(r)<e<\mathbf{f}_{j+1}(r)\}\in\mathscr T\dott
    \end{eqn}
  }
  \argument{
    \lref{eq6.13};
    \lref{eq6.14};
    the fact that $\mathscr T$ is compatible\cfadd{definition:partition-compatibility} with
    $\{(r,a(r))\colon r\in(0,\kappa/2)\}$ and
    $\{(r,b(r))\colon r\in(0,\kappa/2)\}$;
    the fact that for all open $C\in\mathscr T$, $r\in(0,\varepsilon)$
    with $C\subseteq(0,\varepsilon)\times\R$ it holds that
    $\#\{e\in\R\colon(r,e)\in C\}=\infty$;
    the fact that for all $r\in(0,R)$ it holds that $a(r)<b(r)$;
  }{
    that there exist $p,q\in\{1,\dots,m\}$ with $p<q$ which satisfy
    for all $r\in(0,\varepsilon)$ that
    \begin{eqn}
      \llabel{eq6.16}
      a(r)=\mathbf{f}_p(r)\qquad\text{and}\qquad b(r)=\mathbf{f}_q(r)\dott
    \end{eqn}
  }
  In the following let $N\in\N$ satisfy $N=q-p$, let
  $\mathfrak e_j\colon(0,\varepsilon)\to\R$, $j\in\{0,1,\dots,N\}$,
  satisfy for all $j\in\{0,1,\dots,N\}$, $r\in(0,\varepsilon)$ that
  \begin{eqn}
    \llabel{eq6.17}
    \mathfrak e_j(r)=\mathbf{f}_{p+j}(r),
  \end{eqn}
  and let $\mathcal B_i\subseteq(0,\varepsilon)\times\R$,
  $i\in\{0,1,\dots,N-1\}$, satisfy for all $i\in\{0,1,\dots,N-1\}$ that
  \begin{eqn}
    \llabel{eq6.18}
    \mathcal B_i=\{(r,e)\in(0,\varepsilon)\times\R\colon
      \mathfrak e_i(r)<e<\mathfrak e_{i+1}(r)\}\dott
  \end{eqn}
  \startnewargseq
  \argument{
    \lref{eq6.13};
    \lref{eq6.15};
    \lref{eq6.17};
    \lref{eq6.18};
    the fact that $N=q-p$;
    the fact that for all $j\in\{1,\dots,m\}$ it holds that
    $\mathbf{f}_j$ is definable and $\mathbf{f}_j\in C^2((0,\varepsilon),\R)$;
  }{
    that for all $j\in\{0,1,\dots,N\}$ it holds that
    $\mathfrak e_j$ is definable and $\mathfrak e_j\in C^2((0,\varepsilon),\R)$ and that for all
    $i\in\{0,1,\dots,N-1\}$ it holds that $\mathcal B_i$ is open and
    \begin{eqn}
      \llabel{eq6.19}
      \mathcal B_i\in\mathscr T\dott
    \end{eqn}
  }
  \argument{
    \lref{eq6.19};
    \lref{eq6.8};
    the fact that $\mathscr T$ is compatible\cfadd{definition:partition-compatibility} with $\mathcal D$;
  }{
    that for all $i\in\{0,1,\dots,N-1\}$ it holds that
    \begin{eqn}
      \llabel{eq6.20}
      \mathcal B_i\cap\mathcal C=\emptyset\qquad\text{and}\qquad
      \bigl(\mathcal B_i\cap\mathcal D=\emptyset
        \quad\text{or}\quad\mathcal B_i\subseteq\mathcal D\backslash\mathcal C\bigr)\dott
    \end{eqn}
  }
  \argument{
    \lref{eq6.13};
    \lref{eq6.16};
    \lref{eq6.17};
    \lref{eq6.18};
    the fact that $N=q-p$;
  }{
    that
    \begin{eqn}[a]
      \llabel{eq6.21}
      \{(r,e)\in(0,\varepsilon)\times\R\colon a(r)\leq e\leq b(r)\}
      &=\midcup_{j=0}^N\{(r,\mathfrak e_j(r))\colon r\in(0,\varepsilon)\}
        \cup\Bigl(\midcup_{i=0}^{N-1}\mathcal B_i\Bigr)\dott
    \end{eqn}
  }
  \argument{
    \lref{eq6.21};
    \lref{eq6.9};
    \lref{eq6.13};
    \lref{eq6.16};
    \lref{eq6.17};
    \lref{eq6.20};
  }{
    that for all $r\in(0,\varepsilon)$ it holds that
    \begin{eqn}\llabel{eq7}
      a(r)=\mathfrak e_0(r)<\cdots<\mathfrak e_N(r)=b(r)
      \qquad\text{and}\qquad
      \{e\in\R\colon(r,e)\in\mathcal C\}\subseteq
      \{\mathfrak e_0(r),\dots,\mathfrak e_N(r)\}\dott
    \end{eqn}
  }
  \argument{
    \lref{eq6.7};
    \lref{eq6.19};
    the fact that $\zeta|_{\mathcal O}\in C^2(\mathcal O,M)$;
  }{
    that for all $i\in\{0,1,\dots,N-1\}$ with
    $\mathcal B_i\subseteq\mathcal D\backslash\mathcal C$ it holds that
    $\mathcal B_i\subseteq\mathcal O$ and that
    $\zeta|_{\mathcal B_i}\in C^2(\mathcal B_i,M)$\dott\llabel{arg6.2}
  }
  \argument{
    \lref{arg6.2};
    \lref{eq6.3};
    \lref{eq6.6};
    \lref{eq6.19};
  }{
    that for all $i\in\{0,1,\dots,N-1\}$, $j\in\{1,\dots,D\}$,
    $\xi,\eta\in\mathcal B_i$ with
    $\mathcal B_i\subseteq\mathcal D\backslash\mathcal C$ it holds that
    \begin{eqn}
      \llabel{eq7.1}
      \operatorname{sign}\bigl(\bigl(\partial_e (\pi_j\circ\zeta)\bigr)\bigl(\xi\bigr)\bigr)
        =\operatorname{sign}\bigl(\bigl(\partial_e (\pi_j\circ\zeta)\bigr)\bigl(\eta\bigr)\bigr)
      \dott
    \end{eqn}
  }
  In the following let $\mathcal A=\{i\in\{0,1,\dots,N-1\}\colon
  \mathcal B_i\subseteq\mathcal D\backslash\mathcal C\}$ and for every
  $i\in\mathcal A$, $j\in\{1,\dots,D\}$ let
  $\sigma_{ij}\in\{-1,0,1\}$ satisfy that for all
  $(r,e)\in\mathcal B_i$ it holds that
  \begin{eqn}
    \sigma_{ij}=\operatorname{sign}\bigl(
      \bigl(\partial_e (\pi_j\circ\zeta)\bigr)\bigl(r,e\bigr)\bigr)
  \end{eqn}
  and let
  $\mathfrak v_i\colon\mathcal B_i\to\R$ satisfy for all
  $(r,e)\in\mathcal B_i$ that
  $\mathfrak v_i(r,e)=\sum_{j=1}^D\sigma_{ij}\pi_j(\zeta(r,e))$\dott
  \startnewargseq
  \argument{\lref{eq6}; the chain rule; the Cauchy--Schwarz inequality;
   } {
    that for all $i\in\mathcal A$, $(r,e)\in\mathcal B_i$ it holds that
    \begin{gather}\llabel{eq8}
      1=\langle(\nabla_Mh)(r,\zeta(r,e)),\bigl(\partial_e\zeta\bigr)\bigl(r,e\bigr)\rangle
        \leq2\mu(r,e)\norm{\bigl(\partial_e\zeta\bigr)\bigl(r,e\bigr)}\\
      \qquad\text{and}\qquad\bigl(\partial_e\mathfrak v_i\bigr)\bigl(r,e\bigr)
        =\sum_{j=1}^D\abs{\bigl(\partial_e(\pi_j\circ\zeta)\bigr)\bigl(r,e\bigr)}
          \geq\frac1{2\mu(r,e)}>0\dott
    \end{gather}
  }
  \argument{\lref{eq8}; \lref{setting1};
    \cite[Corollary~6.5]{Coste};} {
    that there exist $\rho\in(0,\varepsilon)$ and definable
    functions $\alpha_i,\beta_i,\Delta_i\in C^2((0,\rho),\R)$,
    $i\in\mathcal A$, which satisfy for all
    $i\in\mathcal A$, $r\in(0,\rho)$ that
    \begin{gather}\llabel{eq9}
      \alpha_i(r)=\lim_{e\searrow\mathfrak e_i(r)}\mathfrak v_i(r,e),
      \qquad \beta_i(r)=\lim_{e\nearrow\mathfrak e_{i+1}(r)}\mathfrak v_i(r,e),\\
      \qquad\text{and}\qquad0<\Delta_i(r)=\beta_i(r)-\alpha_i(r)\leq B\dott
    \end{gather}
  }
  \startnewargseq
  \argument{\lref{eq7}; the fact that $\rho<\varepsilon$;}{
    that for all $r\in(0,\rho)$ it holds that
    \begin{eqn}
      \llabel{eq9.1}
      a(r)=\mathfrak e_0(r)<\mathfrak e_1(r)<\cdots<\mathfrak e_N(r)=b(r)\dott
    \end{eqn}
  }
  \argument{\lref{setting1}; \lref{setting2}; \lref{eq7};
    the fact that $\rho<\varepsilon$;}{
    that for all $r\in(0,\rho)$, $e\in\R$ with
    $\inf\bigl(\{\infty\}\cup
      \{z\in\R\colon\exists\,x\in M\colon
        (r,x)\in\Omega\cap h^{-1}(\{e\})\cap G^{-1}(\{z\})\}\bigr)=0$
    it holds that
    \begin{eqn}
      \llabel{eq9.2}
      e\in\{\mathfrak e_0(r),\mathfrak e_1(r),\dots,\mathfrak e_N(r)\}\dott
    \end{eqn}
  }
  In the following let $\mathcal F_N\colon(0,\rho)\times\R\to[0,H]$
  satisfy $\mathcal F_N=0$ and for every $i\in\{0,1,\dots,N-1\}$ let
  $\mathcal F_i\colon(0,\rho)\times\R\to[0,H]$ satisfy for all
  $(r,e)\in(0,\rho)\times\R$ that
  \begin{eqn}\llabel{setting3}
    \mathcal F_i(r,e)=
    \begin{cases}
      0&\colon e\leq\mathfrak e_i(r),\\
      H(\mathfrak v_i(r,e)-\alpha_i(r))/\Delta_i(r)
        &\colon \mathfrak e_i(r)<e<\mathfrak e_{i+1}(r),\ i\in\mathcal A,\\
      H(e-\mathfrak e_i(r))/(\mathfrak e_{i+1}(r)-\mathfrak e_i(r))
        &\colon \mathfrak e_i(r)<e<\mathfrak e_{i+1}(r),\ i\notin\mathcal A,\\
      H&\colon e\geq\mathfrak e_{i+1}(r)
    \end{cases}\dott
  \end{eqn}
  \startnewargseq
  \argument{\lref{eq8}; \lref{eq9}; \lref{setting3};}{
    that for all $i\in\{0,1,\dots,N-1\}$ it holds that $\mathcal F_i$
    is definable and that
    $\mathcal F_i|_{\mathcal B_i\cap((0,\rho)\times\R)}
      \in C^2(\mathcal B_i\cap((0,\rho)\times\R),\R)$,
    that for all $i\in\{0,1,\dots,N-1\}$, $r\in(0,\rho)$
    it holds that $\R\ni e\mapsto\mathcal F_i(r,e)\in[0,H]$ is continuous and
    non-decreasing,
    and that for all $i\in\{0,1,\dots,N-1\}$, $r\in(0,\rho)$,
    $e\in(\mathfrak e_i(r),\mathfrak e_{i+1}(r))$ it holds that
    \begin{eqn}\llabel{eq9.3}
      \bigl(\partial_e\mathcal F_i\bigr)\bigl(r,e\bigr)>0\dott
    \end{eqn}
  }
  \argument{
    \lref{eq6.18};
    the fact that for all $i\in\{0,1,\dots,N-1\}$ it holds that
    $\mathcal F_i|_{\mathcal B_i\cap((0,\rho)\times\R)}
      \in C^2(\mathcal B_i\cap((0,\rho)\times\R),\R)$;
  }{
    that for all $i\in\{0,1,\dots,N-1\}$ it holds that
    \begin{eqn}
      \llabel{eq9.4}
      \mathcal F_i|_{\{(r,e)\in(0,\rho)\times\R\colon
        \mathfrak e_i(r)<e<\mathfrak e_{i+1}(r)\}}
        \quad\text{is twice continuously differentiable}\dott
    \end{eqn}
  }
  \argument{
    \lref{eq9.4};
    \lref{setting3};
    \lref{eq9.3};
    the mean value theorem;
  }{
    that for all $i\in\{0,1,\dots,N-1\}$,
    $r\in(0,\rho)$, $u\in(-\infty,\mathfrak e_i(r)]$,
    $v\in[\mathfrak e_{i+1}(r),\infty)$,
    $z,w\in(\mathfrak e_i(r),\mathfrak e_{i+1}(r))$ with $z<w$
    it holds that
    \begin{eqn}
      \llabel{eq9.5}
      0=\mathcal F_i(r,u)<\mathcal F_i(r,z)<\mathcal F_i(r,w)<\mathcal F_i(r,v)=H\dott
    \end{eqn}
  }
  \argument{\lref{eq9.5}; \lref{setting3};
    the fact that for all $i\in\{0,1,\dots,N-1\}$, $r\in(0,\rho)$
    it holds that $\R\ni e\mapsto\mathcal F_i(r,e)\in[0,H]$ is continuous;
    the intermediate value theorem;
    the fact that for all $j\in\{0,1,\dots,N\}$ it holds that
    $\mathfrak e_j$ is continuous;
    the fact that for all $i\in\{0,1,\dots,N-1\}$ it holds that
    $\mathcal F_i$ is continuous on $\mathcal B_i\cap((0,\rho)\times\R)$;
  }{
    that for all $i\in\{0,1,\dots,N-1\}$, $r\in(0,\rho)$,
    $\eta\in(0,H)$ there exist
    $u,v\in(\mathfrak e_i(r),\mathfrak e_{i+1}(r))$ with $u<v$
    and an open $J\subseteq(0,\rho)$ with $r\in J$ such that
    for all $s\in J$, $e,f\in\R$ with $e\leq u$ and $f\geq v$
    it holds that
    \begin{eqn}\llabel{eq10}
      0\leq\mathcal F_i(s,e)\leq\mathcal F_i(s,u)<\eta
      \qquad\text{and}\qquad
      0\leq H-\mathcal F_i(s,f)\leq H-\mathcal F_i(s,v)<\eta\dott
    \end{eqn}
  }
  \argument{\lref{eq10}; \lref{setting3};
    the fact that for all $j\in\{0,1,\dots,N\}$ it holds that
    $\mathfrak e_j$ is continuous;
    the fact that for all $i\in\{0,1,\dots,N-1\}$ it holds that
    $\mathcal F_i$ is continuous on $\mathcal B_i\cap((0,\rho)\times\R)$;
  }{
    that for all $i\in\{0,1,\dots,N-1\}$ it holds that
    $\mathcal F_i$ is continuous on $(0,\rho)\times\R$\dott
    \llabel{arg2}
  }
  \argument{\lref{arg2}; \lref{eq9.3}; \lref{eq9.4}; \lref{eq9.5};
    the fact that for all $i\in\{0,1,\dots,N-1\}$ it holds that
    $\mathcal F_i$ is definable;
    the intermediate value theorem;
    the inverse function theorem;} {
    that there exist definable
    $E_1,E_2,\dots,E_N\in C^2((0,\rho)\times(0,H),\R)$ and
    $T_1,T_2,\dots,T_N\in C^1((0,\rho)\times(0,H),(0,\infty))$
    which satisfy for all $i\in\{1,2,\dots,N\}$,
    $r\in(0,\rho)$, $z\in(0,H)$ that
    \begin{eqn}
      \llabel{eq10.1}
      \mathfrak e_{i-1}(r)<E_i(r,z)<\mathfrak e_i(r),\quad
      \mathcal F_{i-1}(r,E_i(r,z))=z,\qquad\text{and}\qquad
      T_i(r,z)=\partial_z E_i(r,z)>0,
    \end{eqn}
    and which satisfy for all $i\in\{0,1,\dots,N-1\}$,
    $(r,e)\in\mathcal B_i$ with $r<\rho$ that
    \begin{eqn}\llabel{eq11}
      T_{i+1}(r,\mathcal F_i(r,e))=\frac1{\bigl(\partial_e\mathcal F_i\bigr)\bigl(r,e\bigr)}\dott
    \end{eqn}
  }
  \argument{\lref{eq11}; \lref{setting1}; \lref{setting2};
    \lref{eq8}; \lref{eq9}; \lref{setting3};} {
    that for all $i\in\mathcal A$, $(r,x)\in\Omega$ with $r<\rho$
    and $\mathfrak e_i(r)<h(r,x)<\mathfrak e_{i+1}(r)$ it holds that
    \begin{eqn}\llabel{eq12}
      T_{i+1}(r,\mathcal F_i(r,h(r,x)))
        \leq\frac{2\Delta_i(r)}H\mu(r,h(r,x))
        \leq\frac{2B}{1+2B/c}G(r,x)\leq cG(r,x)\dott
    \end{eqn}
  }
  \argument{\lref{eq12}; \lref{setting1}; \lref{eq6.18}; \lref{eq11};
    the fact that for all $i\in\{0,1,\dots,N-1\}\backslash\mathcal A$ it holds that
    $\mathcal B_i\cap\mathcal D=\emptyset$;
  }{
    that for all $i\in\{1,2,\dots,N\}$,
    $(r,x)\in\Omega\cap((0,\rho)\times M)$ with
    $\mathfrak e_{i-1}(r)<h(r,x)<\mathfrak e_i(r)$ it holds that
    \begin{eqn}
      \llabel{eq12.1}
      T_i(r,\mathcal F_{i-1}(r,h(r,x)))\leq cG(r,x)\dott
    \end{eqn}
  }
  \argument{
    \lref{eq9.1};
    \lref{eq9.2};
    \lref{eq9.4};
    \lref{eq9.5};
    \lref{eq10.1};
    \lref{eq12.1};
  }[verbs=ep]{
    \lref{item1}, \lref{item2}, \lref{item3}, \lref{item4},
    \lref{item5}, and \lref{item6}\dott
  }
\end{aproof}

\subsection{Upper and lower bounds for the radial derivatives of the back-trans\-for\-mations}
\label{subsec:residual-bound}

\begin{athm}{lemma}{lem:radial-decay}
  Let $\mathfrak{r}_0\in(0,\infty)$\dott \cfload[.]
  Then
  \begin{enumerate}[label=(\roman*)]
    \item \llabel{item1} for every definable\cfadd{definition:definable-function}
    and integrable $n\colon(0,\mathfrak{r}_0)\to[0,\infty)$ there exists
    $\mathfrak{r}_1\in(0,\mathfrak{r}_0)$ such that
    $n|_{(0,\mathfrak{r}_1)}$ is continuous and
    \begin{equation}
      \llabel{concl1}
    \textstyle  \limsup_{s\to 0}|sn(s)|=0
    \end{equation}
    and
    \item \llabel{item2} for every $\varepsilon\in(0,\infty)$,
    $N\in\N$ and all definable and integrable 
    $\mathfrak{n}_1,\mathfrak{n}_2,\ldots,\mathfrak{n}_N\colon(0,\mathfrak{r}_0)\to[0,\infty)$ there exists
    $R\in(0,\mathfrak{r}_0)$ such that for all $i\in\{1,2,\ldots,N\}$,
    $s\in(0,R)$ it holds that
    \begin{eqn}
      \llabel{concl2}
      s^2|\mathfrak{n}_i(s)|^2\leq\varepsilon
      \dott
    \end{eqn}
  \end{enumerate}
  \cfout[.]
\end{athm}
\begin{aproof}
  In the following we prove \lref{item1}\dott
  \startnewargseq
  To prove \lref{item1}, let
  $n\colon(0,\mathfrak{r}_0)\to[0,\infty)$ be
  definable\cfadd{definition:definable-function} and integrable\dott \cfload[.]
  \startnewargseq
  \argument{
    the assumption that $n$ is definable;
    \cite[Theorem~2.1 and Corollary~6.5]{Coste};
  }{
    that there exists $R\in(0,\mathfrak{r}_0)$ which satisfies that $n|_{(0,R)}$ is continuous
    and is strictly increasing, strictly decreasing, or constant\dott
  }
 In the following let
  $\mathfrak{r}_1\in(0,R)$.
  \startnewargseq
  \argument{the fact that
     $n|_{(0,R)}$ is continuous;
    the fact that $\mathfrak{r}_1<R$;
  }{
    that $n|_{(0,\mathfrak{r}_1)}$ is continuous\dott
  }
  In our proof of \lref{item1} we distinguish between the case that $n|_{(0,R)}$ is non-decreasing and the case that $n|_{(0,R)}$ is non-increasing. In the following we prove \lref{item1} in the case that
  \begin{eqn}
      n|_{(0,R)}\quad \text{is non-decreasing}.
  \end{eqn}
  \startnewargseq
  \argument{
    the assumption that $n|_{(0,R)}$ is non-decreasing;
    the fact that $n$ is nonnegative;
  }{
    that for all $r\in(0,\mathfrak{r}_1)$ it holds that
    \begin{eqn}
      \llabel{eq2}
      0\leq rn(r)\leq r n(\mathfrak{r}_1)
      \dott
    \end{eqn}
  }
  \argument{
    \lref{eq2};
  }{
    that 
    \begin{eqn}
      \llabel{eq3}
      \lim_{r\to 0}rn(r)=0
      \dott
    \end{eqn}
  }
  This shows \lref{item1} in the case that $n|_{(0,R)}$ is non-decreasing. 
  We now prove \lref{item1} in the case that
  \begin{eqn}
      n|_{(0,R)}\quad \text{is non-increasing}.
  \end{eqn}
  \startnewargseq
  \argument{
    the assumption that $n|_{(0,R)}$ is non-increasing;
    the fact that $n$ is nonnegative;
  }{
    that for all $r\in(0,\mathfrak{r}_1)$ it holds that
    \begin{eqn}
      \llabel{eq4}
      0\leq \frac r2 n(r)
      \leq \int_{r/2}^{r}n(s)\dd s
      \leq \int_0^{r}n(s)\dd s
      \dott
    \end{eqn}
  }
  \argument{
    the assumption that $n$ is non-negative and integrable on $(0,\mathfrak{r}_0)$; the dominated convergence theorem;
  }{
    that
    \begin{eqn}
      \llabel{eq5}
      \lim_{r\to 0}\int_0^r n(s)\dd s=0
      \dott
    \end{eqn}
  }
  \argument{
    \lref{eq5};
    \lref{eq4};
  }{
    that 
    \begin{eqn}
      \llabel{eq6}
      \lim_{r\to 0}rn(r)=0
      \dott
    \end{eqn}
  }
This establishes \lref{item1} in the case that $n|_{(0,R)}$ is non-increasing.
  \startnewargseq
  This completes the proof of \lref{item1}\dott

  \startnewargseq
  In the following we prove \lref{item2}\dott
  \startnewargseq
  To prove \lref{item2}, let $N\in\N$, $\varepsilon\in(0,\infty)$,
  and let $\mathfrak{n}_1,\ldots,\mathfrak{n}_N\colon
  (0,\mathfrak{r}_0)\to[0,\infty)$ be definable and integrable\dott
  \startnewargseq
  \argument{
    \lref{item1};
  }{
    that for all $i\in\{1,\ldots,N\}$ it holds that
    \begin{eqn}
      \llabel{eq7}
      \lim_{r\to 0}r^2(\mathfrak{n}_i(r))^2=0
      \dott
    \end{eqn}
  }
  \argument{
    \lref{eq7};
    the assumption that $\varepsilon>0$;
  }{
    that for every $i\in\{1,\ldots,N\}$ there exists
    $\boldsymbol{\rho}_i\in(0,\mathfrak{r}_0)$ which satisfies for all $r\in(0,\boldsymbol{\rho}_i)$ that
    \begin{eqn}
      \llabel{eq8}
      r^2(\mathfrak{n}_i(r))^2\leq\varepsilon
      \dott
    \end{eqn}
  }
 In the following let
  $\mathfrak{r}_\varepsilon=\min_{i\in\{1,2,\dots,N\}}\boldsymbol{\rho}_i$.
  \startnewargseq
  \argument{
    the fact that $N\in\N$;
    the fact that for all
    $i\in\{1,\ldots,N\}$ it holds that $\boldsymbol{\rho}_i\in(0,\mathfrak{r}_0)$;
    the fact that $\mathfrak{r}_\varepsilon=\min_{1\leq i\leq N}\boldsymbol{\rho}_i$;
  }{
    that $\mathfrak{r}_\varepsilon\in(0,\mathfrak{r}_0)$ and that for all
    $i\in\{1,\ldots,N\}$ it holds that $\mathfrak{r}_\varepsilon\leq\boldsymbol{\rho}_i$\dott
    \llabel{eq9}
  }
  \argument{
    \lref{eq9};
    \lref{eq8};
  }[verbs=ep,plural]{
    \lref{item2}\dott
  }
\end{aproof}

\begin{athm}{lemma}{lem:residual}
  Let $D,N\in\N$, $\mathfrak{r}_0,\mathbf{c},H\in(0,\infty)$,
  $\mathfrak{r}_1\in(0,\mathfrak{r}_0)$, let
  $M\subseteq\R^D$ be a nonempty, compact, and definable
  \submanifold{2}{D}, let
  $\mathcal O\subseteq(0,\mathfrak{r}_0)\times M$ be definable
  and open in $(0,\mathfrak{r}_0)\times M$, let
  $h\in \COne{\mathcal O}{\R}$,
  $\mathfrak{e}_0,\mathfrak{e}_1,\dots,\mathfrak{e}_N
    \in C^2((0,\mathfrak{r}_1),\R)$ be definable\cfadd{definition:c1-on-manifold},
  assume for all $r\in(0,\mathfrak{r}_1)$ that
  \begin{eqn}
    \llabel{ass1}
    \mathfrak{e}_0(r)<\mathfrak{e}_1(r)<\ldots<\mathfrak{e}_N(r),
  \end{eqn}
  for every $i\in\{1,2,\dots,N\}$ let
  $\mathcal F_i\in C((0,\mathfrak{r}_1)\times\R,[0,H])$,
  $E_i=(E_i(r,z))_{(r,z)\in(0,\mathfrak{r}_1)\times(0,H)}
    \in C^2((0,\mathfrak{r}_1)\times(0,H),\R)$ be definable,
  and assume for all $i\in\{1,2,\dots,N\}$,
  $r\in(0,\mathfrak{r}_1)$,
  $u,v\in(\mathfrak e_{i-1}(r),\mathfrak e_i(r))$, $z\in(0,H)$
  with $u<v$ that
  \begin{eqn}
    \llabel{ass5}
    \mathcal F_i(r,u)<\mathcal F_i(r,v),\quad
    \mathfrak{e}_{i-1}(r)<E_i(r,z)<\mathfrak{e}_i(r),\quad
    \mathcal{F}_i(r,E_i(r,z))=z,\quad\text{and}\quad
    \partial_z E_i(r,z)>0
    \dott
  \end{eqn}
  \cfload[.]
  Then there exist $\mathfrak{r}_2\in(0,\mathfrak{r}_1)$ and definable
  $\mathbf{m}_1,\mathbf{m}_2,\dots,\mathbf{m}_N
    \in C^1((0,\mathfrak{r}_2),[1,\infty))$ such that
  \begin{enumerate}[label=(\roman*)]
    \item \llabel{item0} for all $i\in\{1,2,\dots,N\}$,
    $r\in(0,\mathfrak{r}_2)$ it holds that
    $\int_0^{\mathfrak{r}_2}\mathbf{m}_i(s)\dd s<\infty$ and
    $r^2(\mathbf{m}_i(r))^2\leq\mathbf c$ and
    \item \llabel{item1} for all $i\in\{1,2,\dots,N\}$,
    $(r,x)\in\mathcal O\cap((0,\mathfrak{r}_2)\times M)$ with
    $\mathfrak e_{i-1}(r)<h(r,x)<\mathfrak e_i(r)$ it holds that\cfadd{definition:angular-gradient}
    \begin{eqn}
      \llabel{eq:residual}
      \bigl|\partial_r h(r,x)
        -(\partial_r E_i)(r,\mathcal F_i(r,h(r,x)))\bigr|
      \leq \mathbf{m}_i(r)\Bigl((\partial_z E_i)(r,\mathcal F_i(r,h(r,x)))+\mathbf{c}|(\nabla_Mh)(r,x)|\Bigr)
      \dott
    \end{eqn}
  \end{enumerate}
  \cfout[.]
\end{athm}
\begin{aproof}
  Throughout this proof let $G\colon\mathcal O\to\R$ satisfy for all
  $(r,x)\in\mathcal O$ that
  $G(r,x)=|(\nabla_Mh)(r,x)|$\cfadd{definition:angular-gradient} and
  for every $i\in\{1,2,\dots,N\}$ let
  $\Omega_i\subseteq\R^{D+1}$ satisfy
  \begin{eqn}
    \llabel{setting0}
    \Omega_i
      =\{(r,x)\in\mathcal O\cap((0,\mathfrak{r}_1)\times M)\colon
          \mathfrak{e}_{i-1}(r)<h(r,x)<\mathfrak{e}_i(r)\},
  \end{eqn}
  let
  $c_i,T_i\colon(0,\mathfrak{r}_1)\times(0,H)\to\R$ satisfy for all
  $r\in(0,\mathfrak{r}_1)$, $z\in(0,H)$ that
  \begin{eqn}
    \llabel{setting1}
    c_i(r,z)=\partial_r E_i(r,z)\qquad\text{and}\qquad T_i(r,z)=\partial_z E_i(r,z),
  \end{eqn}
  and let $Z_i\colon\mathcal O\cap((0,\mathfrak{r}_1)\times M)\to[0,H]$ satisfy for all
  $(r,x)\in\mathcal O\cap((0,\mathfrak{r}_1)\times M)$ that
  \begin{eqn}
    \llabel{setting1.1}
    Z_i(r,x)=\mathcal F_i(r,h(r,x))
    \dott
  \end{eqn}
  \cfload[.]
  \startnewargseq
  \argument{
    \lref{ass5};
  }{
    that for all $i\in\{1,2,\dots,N\}$, $r\in(0,\mathfrak r_1)$,
    $e\in(\mathfrak e_{i-1}(r),\mathfrak e_i(r))$ it holds that
    \begin{eqn}
      \llabel{eq:interior-values}
      0<\mathcal F_i(r,e)<H\dott
    \end{eqn}
  }
  \argument{
    \lref{eq:interior-values};
    \lref{ass5};
  }{
    that for all $i\in\{1,2,\dots,N\}$, $r\in(0,\mathfrak r_1)$,
    $e\in(\mathfrak e_{i-1}(r),\mathfrak e_i(r))$ it holds that
    \begin{eqn}
      \llabel{eq:inverse-identity}
      E_i(r,\mathcal F_i(r,e))=e\dott
    \end{eqn}
  }
  \argument{
    \lref{setting0};
    the fact that $\mathcal O$ is definable and open in
    $(0,\mathfrak{r}_0)\times M$;
    the fact that for all $i\in\{1,2,\dots,N\}$ it holds that
    $h$, $\mathfrak{e}_{i-1}$, and $\mathfrak{e}_i$ are definable and continuous;
  }{
    that for all $i\in\{1,2,\dots,N\}$ it holds that $\Omega_i$
    is definable and is open in $(0,\mathfrak{r}_1)\times M$\dott
    \llabel{eq0}
  }
  \argument{
    \lref{eq0};
    \lref{ass5};
    \lref{setting0};
    \lref{setting1.1};
    \lref{eq:interior-values};
    \lref{eq:inverse-identity};
    the fact that for all $i\in\{1,2,\dots,N\}$ it holds that
    $E_i\in C^2((0,\mathfrak r_1)\times(0,H),\R)$;
    the inverse function theorem;
    the fact that $h$ is definable and $h\in \COne{\mathcal O}{\R}$;
    the chain rule;
  }{
    that for all $i\in\{1,2,\dots,N\}$ it holds that
    $\mathcal F_i|_{\{(r,e)\in(0,\mathfrak r_1)\times\R\colon
      \mathfrak e_{i-1}(r)<e<\mathfrak e_i(r)\}}$
    is twice continuously differentiable,
    $Z_i|_{\Omega_i}\colon\Omega_i\to\R$ is definable,
    $Z_i|_{\Omega_i}\in \COne{\Omega_i}{\R}$, and for all
    $(r,x)\in\Omega_i$ it holds that
    \begin{eqn}
      \llabel{eq1}
      \mathfrak e_{i-1}(r)<h(r,x)<\mathfrak e_i(r)
      \qquad\text{and}\qquad 0<Z_i(r,x)<H
      \dott
    \end{eqn}
  }
  In the following for every $i\in\{1,2,\dots,N\}$ let
  $\mathbf{Z}_i=(\mathbf{Z}_i(r,x))_{(r,x)\in\Omega_i}=Z_i|_{\Omega_i}$\dott
  \startnewargseq
  \argument{
    \lref{ass5};
    \lref{setting1};
    \lref{setting1.1};
    \lref{eq:inverse-identity};
    \lref{eq1};
  }{
    that for all $i\in\{1,2,\dots,N\}$,
    $(r,x)\in\Omega_i$ it holds that
    \begin{eqn}
      \llabel{eq2}
      h(r,x)=E_i(r,Z_i(r,x))\qquad \text{and}\qquad  T_i(r,Z_i(r,x))>0
      \dott
    \end{eqn}
  }
  \argument{
    \lref{eq2};
    the fact that for all $i\in\{1,2,\dots,N\}$ it holds that
    $\Omega_i$ is open in $(0,\mathfrak{r}_1)\times M$;
  }{
    that for all $i\in\{1,2,\dots,N\}$,
    $(r,x)\in\Omega_i$ there exists
    $\delta\in(0,\min\{r,\mathfrak{r}_1-r\})$ which satisfies that for all
    $s\in(r-\delta,r+\delta)$ it holds that
    \begin{eqn}
      \llabel{eq2.1}
      (s,x)\in\Omega_i\qquad \text{and}\qquad  h(s,x)=E_i(s,Z_i(s,x))
      \dott
    \end{eqn}
  }
  \argument{
    \lref{eq2.1};
    the fact that for all $i\in\{1,2,\dots,N\}$ it holds that
    $\mathbf{Z}_i\in \COne{\Omega_i}{\R}$ and
    $E_i\in C^2((0,\mathfrak{r}_1)\times(0,H),\R)$;
    the fact that $h\in \COne{\mathcal O}{\R}$;
    the chain rule;
  }{
    that for all $i\in\{1,2,\dots,N\}$,
    $(r,x)\in\Omega_i$ it holds that
    \begin{eqn}[a]
      \llabel{eq2.2}
      \partial_r h(r,x)
      &=(\partial_r E_i)(r,Z_i(r,x))
        +\bigl(\partial_z E_i\bigr)\bigl(r,Z_i(r,x)\bigr)
          (\partial_r \mathbf{Z}_i)(r,x)
      \dott
    \end{eqn}
  }
  \argument{
    \lref{eq2};
    the fact that for all $i\in\{1,2,\dots,N\}$ it holds that
    $\Omega_i$ is open in $(0,\mathfrak{r}_1)\times M$;
  }{
    that for every $i\in\{1,2,\dots,N\}$,
    $(r,x)\in\Omega_i$, $\varepsilon\in(0,\infty)$ and every
    $\gamma\in C^1((-\varepsilon,\varepsilon),M)$ with $\gamma(0)=x$
    there exists $\delta\in(0,\varepsilon)$ which satisfies that for all
    $t\in(-\delta,\delta)$ it holds that
    \begin{eqn}
      \llabel{eq2.3}
      (r,\gamma(t))\in\Omega_i\qquad\text{and}\qquad
      h(r,\gamma(t))=E_i(r,Z_i(r,\gamma(t)))
      \dott
    \end{eqn}
  }
  \argument{
    \lref{eq2.3};
    the fact that for all $i\in\{1,2,\dots,N\}$ it holds that
    $\mathbf{Z}_i\in \COne{\Omega_i}{\R}$ and
    $E_i\in C^2((0,\mathfrak{r}_1)\times(0,H),\R)$;
    the fact that $h\in \COne{\mathcal O}{\R}$;
    the chain rule;
  }{
    that for every $i\in\{1,2,\dots,N\}$,
    $(r,x)\in\Omega_i$, $\varepsilon\in(0,\infty)$ and every
    $\gamma\in C^1((-\varepsilon,\varepsilon),M)$ with $\gamma(0)=x$
    it holds that
    \begin{eqn}[a]
      \llabel{eq2.4}
      \langle(\nabla_Mh)(r,x),\gamma'(0)\rangle
      &=\frac{\ddd}{\ddd t}h(r,\gamma(t))\big|_{t=0}
        =\bigl(\partial_z E_i\bigr)\bigl(r,Z_i(r,x)\bigr)
          \frac{\ddd}{\ddd t}\mathbf{Z}_i(r,\gamma(t))\big|_{t=0}\\
      &=\bigl(\partial_z E_i\bigr)\bigl(r,Z_i(r,x)\bigr)
          \langle(\nabla_M\mathbf{Z}_i)(r,x),\gamma'(0)\rangle
      \cfadd{definition:angular-gradient}\dott
    \end{eqn}
  }
  \argument{
    the fact that for all $i\in\{1,2,\dots,N\}$,
    $(r,x)\in\Omega_i$ it holds that
    $(\nabla_M\mathbf{Z}_i)(r,x)\in T_xM$\cfadd{definition:tangent-space};
    the fact that for all $i\in\{1,2,\dots,N\}$,
    $(r,x)\in\Omega_i$ it holds that
    $\bigl(\partial_z E_i\bigr)\bigl(r,Z_i(r,x)\bigr)\in\R$;
    the fact that for all $x\in M$ it holds that $T_xM$
    is a linear subspace of $\R^D$;
  }{
    that for all $i\in\{1,2,\dots,N\}$,
    $(r,x)\in\Omega_i$ it holds that
    \begin{eqn}
      \llabel{eq2.4.1}
      \bigl(\partial_z E_i\bigr)\bigl(r,Z_i(r,x)\bigr)
        (\nabla_M\mathbf{Z}_i)(r,x)\in T_xM
      \dott
    \end{eqn}
  }
  \argument{
    \lref{eq2.4.1};
    \lref{eq2.4};
    \cref{lem:angular-gradient-existence};
  }{
    that for all $i\in\{1,2,\dots,N\}$,
    $(r,x)\in\Omega_i$ it holds that
    \begin{eqn}
      \llabel{eq2.5}
      (\nabla_Mh)(r,x)
        =\bigl(\partial_z E_i\bigr)\bigl(r,Z_i(r,x)\bigr)
          (\nabla_M\mathbf{Z}_i)(r,x)
      \dott
    \end{eqn}
  }
  \argument{
    \lref{eq2.5};
    \lref{setting1};
    \lref{eq2.2};
  }{
    that for all $i\in\{1,2,\dots,N\}$,
    $(r,x)\in\Omega_i$ it holds that $ (\partial_r h)(r,x)
        =c_i(r,Z_i(r,x))
          +T_i(r,Z_i(r,x))(\partial_r \mathbf{Z}_i)(r,x)$ and
    \begin{eqn}
      \llabel{eq3}
      (\nabla_Mh)(r,x)
        =T_i(r,Z_i(r,x))(\nabla_M\mathbf{Z}_i)(r,x)
      \dott
    \end{eqn}
  }
  \argument{
    \lref{eq3};
    \lref{eq2};
    the fact that for all $(r,x)\in\mathcal O$ it holds that
    $G(r,x)=|(\nabla_Mh)(r,x)|$;
  }{
    that for all $i\in\{1,2,\dots,N\}$,
    $(r,x)\in\Omega_i$ it holds that
    \begin{eqn}
      \llabel{eq4}
      \frac{\bigl|(\partial_r \mathbf{Z}_i)(r,x)\bigr|}
           {1+|(\nabla_M\mathbf{Z}_i)(r,x)|}
      =
      \frac{\bigl|(\partial_r h)(r,x)
                    -c_i(r,Z_i(r,x))\bigr|}
           {T_i(r,Z_i(r,x))+G(r,x)}
      \dott
    \end{eqn}
  }
  In the following for every $i\in\{1,2,\dots,N\}$,
  $r\in(0,\mathfrak{r}_1)$ let $\mathcal X_{i,r}\subseteq M$ satisfy
  \begin{eqn}
    \llabel{setting2}
    \mathcal X_{i,r}=\{x\in M\colon(r,x)\in\Omega_i\}
  \end{eqn}
  and for every $i\in\{1,2,\dots,N\}$ let
  $\eta_i\colon(0,\mathfrak{r}_1)\to[0,\infty]$ satisfy that for all
  $r\in(0,\mathfrak{r}_1)$ it holds that
  \begin{eqn}
    \llabel{setting3}
    \eta_i(r)=\sup_{x\in\mathcal X_{i,r}}
      \frac{\bigl|(\partial_r \mathbf{Z}_i)(r,x)\bigr|}
           {1+|(\nabla_M\mathbf{Z}_i)(r,x)|},
  \end{eqn}
  where the supremum of the empty set is defined to be zero\dott
  \startnewargseq
  \argument{
    \lref{eq1};
    \lref{setting2};
    \lref{setting3};
    the fact that $M$ is a compact definable \submanifold{2}{D};
    the fact that for all $i\in\{1,2,\dots,N\}$ it holds that
    $\Omega_i$ is definable and is open in $(0,\mathfrak{r}_1)\times M$;
    the fact that for all $i\in\{1,2,\dots,N\}$ it holds that
    $\mathbf{Z}_i$ is definable and
    $\mathbf{Z}_i\in \COne{\Omega_i}{\R}$;
    \cref{lem:normalized} (applied for every $i\in\{1,2,\dots,N\}$
    with $n\curvearrowleft D$, $\mathfrak{r}_0\curvearrowleft \mathfrak{r}_1$,
    $U\curvearrowleft\Omega_i$,
    $\Omega\curvearrowleft\Omega_i$,
    $Z\curvearrowleft \mathbf{Z}_i$,
    $\eta\curvearrowleft\eta_i$ in the notation of
    \cref{lem:normalized});
  }{
    that for all $i\in\{1,2,\dots,N\}$ there exists
    $\boldsymbol{\rho}_i\in(0,\mathfrak{r}_1)$ which satisfies that
    \begin{enumerate}[label=(\Roman*)]
      \item \llabel{item:bound-finite} it holds that
      $\eta_i((0,\boldsymbol{\rho}_i))\subseteq\R$ and
      \item \llabel{eq5} it holds that
      $(0,\boldsymbol{\rho}_i)\ni r\mapsto\eta_i(r)\in\R$
      is continuously differentiable, definable, and integrable\dott
    \end{enumerate}
  }
  In the following let $R\in(0,\mathfrak{r}_1)$ satisfy
  \begin{eqn}
    \llabel{setting4}
    R=\frac12\min\{\boldsymbol{\rho}_i\colon i\in\{1,2,\dots,N\}\}
    \dott
  \end{eqn}
  \startnewargseq
  \argument{ \lref{item:bound-finite};
    \lref{eq2};
    \lref{eq4};
    \lref{setting2};
    \lref{setting3};
    \lref{setting4};
    the fact that $G$ is nonnegative;
  }{
    that for all $i\in\{1,2,\dots,N\}$,
    $(r,x)\in\Omega_i\cap((0,R)\times M)$ it holds that
    \begin{eqn}
      \llabel{eq6}
      \bigl|(\partial_r h)(r,x)
         -c_i(r,Z_i(r,x))\bigr|
      \leq\eta_i(r)\bigl(T_i(r,Z_i(r,x))+G(r,x)\bigr)
      \dott
    \end{eqn}
  }
  In the following for every $i\in\{1,2,\dots,N\}$ let
  $\mathbf{n}_i\colon(0,R)\to[1,\infty)$ satisfy that for all
  $r\in(0,R)$ it holds that
  \begin{eqn}
    \llabel{setting5}
    \mathbf{n}_i(r)=1+\max\{1,\mathbf{c}^{-1}\}\eta_i(r)
    \dott
  \end{eqn}
  \startnewargseq
  \argument{
    \lref{eq5};
    \lref{setting4};
    \lref{setting5};
    the fact that for all $i\in\{1,2,\dots,N\}$ it holds that
    $\eta_i|_{(0,\boldsymbol{\rho}_i)}$ is definable;
    the fact that for all $i\in\{1,2,\dots,N\}$ it holds that
    $\eta_i$ is nonnegative;
  }{
    that for all $i\in\{1,2,\dots,N\}$ it holds that
    $\mathbf{n}_i$ is definable, that
    $\mathbf{n}_i\in C^1((0,R),[1,\infty))$,
    and that
    \begin{eqn}
      \llabel{eq7}
      \int_0^R\mathbf{n}_i(r)\dd r
      =R+\max\{1,\mathbf{c}^{-1}\}\int_0^R\eta_i(r)\dd r<\infty
      \dott
    \end{eqn}
  }
  \argument{
    \lref{eq7};
    the fact that for all $i\in\{1,2,\dots,N\}$ it holds that
    $\mathbf{n}_i$ is definable and nonnegative;
    \cref{lem:radial-decay} (applied with $\mathfrak{r}_0\curvearrowleft R$,
    $(\mathfrak{n}_1,\ldots,\mathfrak{n}_N)\curvearrowleft
      (\mathbf{n}_1,\ldots,\mathbf{n}_N)$,
    $\varepsilon\curvearrowleft\mathbf c$ in the notation of
    \cref{lem:radial-decay});
  }{
    that there exists $\mathfrak{r}_2\in(0,R)$ which satisfies for all
    $i\in\{1,2,\dots,N\}$, $r\in(0,\mathfrak{r}_2)$ that
    \begin{eqn}
      \llabel{eq8}
      r^2(\mathbf{n}_i(r))^2\leq\mathbf c
      \dott
    \end{eqn}
  }
  In the following for every $i\in\{1,2,\dots,N\}$ let
  $\mathbf{m}_i\colon(0,\mathfrak{r}_2)\to[1,\infty)$ satisfy for all
  $r\in(0,\mathfrak{r}_2)$ that
  \begin{eqn}
    \llabel{eq:ni}
    \mathbf{m}_i(r)=\mathbf{n}_i(r)
      =1+\max\{1,\mathbf{c}^{-1}\}\eta_i(r)
    \dott
  \end{eqn}
  \startnewargseq
  \argument{
    \lref{eq7};
    \lref{eq:ni};
    the fact that $\mathfrak{r}_2\in(0,R)\subseteq(0,\mathfrak{r}_1)$;
    the fact that for all $i\in\{1,2,\dots,N\}$ it holds that
    $\mathbf{n}_i$ is definable and $\mathbf{n}_i\in C^1((0,R),[1,\infty))$;
  }{
    that for all $i\in\{1,2,\dots,N\}$ it holds that
    $\mathbf{m}_i$ is definable, that
    $\mathbf{m}_i\in C^1((0,\mathfrak{r}_2),[1,\infty))$, and that
    \begin{eqn}
      \llabel{eq9}
      \int_0^{\mathfrak{r}_2}\mathbf{m}_i(r)\dd r
      \leq\int_0^R\mathbf{n}_i(r)\dd r<\infty
      \dott
    \end{eqn}
  }
  \argument{
    \lref{eq9};
    \lref{eq8};
    \lref{eq:ni};
  }[verbs=ep]{
    \lref{item0}\dott
  }

  \startnewargseq
  \argument{
    \lref{eq:ni};
    the fact that for all $i\in\{1,2,\dots,N\}$ it holds that
    $\eta_i$ is nonnegative;
    the fact that $\max\{1,\mathbf{c}^{-1}\}\geq1$ ;the fact that
    $\mathbf{c}\max\{1,\mathbf{c}^{-1}\}\geq1$;
  }{
    that for all $i\in\{1,2,\dots,N\}$, $r\in(0,\mathfrak{r}_2)$
    it holds that
    \begin{eqn}
      \llabel{eq10}
      \eta_i(r)\leq \mathbf{m}_i(r)\qquad\text{and}\qquad \eta_i(r)\leq \mathbf{c}\mathbf{m}_i(r)
      \dott
    \end{eqn}
  }
  \argument{
    \lref{eq10};
    \lref{eq2};
    \lref{eq6};
    the fact that $G$ is nonnegative;
    the fact that $(0,\mathfrak{r}_2)\subseteq(0,R)$;
  }{
    that for all $i\in\{1,2,\dots,N\}$,
    $(r,x)\in\Omega_i\cap((0,\mathfrak{r}_2)\times M)$ it holds that
    \begin{eqn}[a]
      \llabel{eq11}
      \bigl|(\partial_r h)(r,x)
        -c_i(r,Z_i(r,x))\bigr|
      &\leq\eta_i(r)T_i(r,Z_i(r,x))+\eta_i(r)G(r,x)\\
      &\leq \mathbf{m}_i(r)T_i(r,Z_i(r,x))+\mathbf{c}\mathbf{m}_i(r)G(r,x)\\
      &=\mathbf{m}_i(r)\bigl(T_i(r,Z_i(r,x))+\mathbf{c}G(r,x)\bigr)
      \dott
    \end{eqn}
  }
  \argument{
    \lref{eq11};
    \lref{setting0};
    \lref{setting1};
    \lref{setting1.1};
    the fact that for all $(r,x)\in\mathcal O$ it holds that
    $G(r,x)=|(\nabla_Mh)(r,x)|$;
  }[verbs=ep]{
    \lref{item1}\dott
  }
\end{aproof}
\section{Construction and regularity of the truncation of the transformed energy}
\label{sec:truncation}

As outlined in \cref{subsec:proof-outline} above, the key issue in the
proof of \cref{maintheorem} is to establish Lyapunov-type stability
(in the precise sense of \eqref{eq:intro-radial-angular-lyapunov} in
\cref{subsec:proof-outline}) for the radial-angular \DE\ system
(see \eqref{eq:intro-radial-angular-system} in \cref{subsec:proof-outline})
as asserted in \cref{prop:buildQ} in
\cref{sec:radial-integration} below.
A difficulty that arises in the construction of $Q$ in
\eqref{eq:intro-lyapunov-construction} in \cref{subsec:proof-outline}
is that the transformed energy function $Z_i$ is $C^1$ on the set where
$(r,h(r,x))\in\mathbf{H}_i$, but need not be continuously differentiable
on its entire domain and its composition with a locally absolutely
continuous trajectory need thus not be locally absolutely continuous.

We solve this issue by employing a \emph{truncation and approximation procedure}
for the transformed energy $Z_i(r,x)$.
Specifically, we approximate $Z_i(r,x)$ by
$\max\{\varepsilon,\min\{H-\varepsilon,Z_i(r,x)\}\}$ for
$\varepsilon\in(0,H/2)$ and it is precisely the subject of
\cref{lem:truncated} in this section to establish
regularity properties and estimates for the truncated transformed energy.

\subsection{Vanishing of the derivatives of absolutely continuous functions on their level sets}
\label{subsec:level-set-derivative}

In the following elementary \cref{lem:AC-level-set} we recall the well-known
fact that for every $c\in\R$ and every absolutely continuous function $f$
we have that $f'$ vanishes almost everywhere on the $c$-level set
$f^{-1}(\{c\})$.
Only for completeness, we include a detailed proof for
\cref{lem:AC-level-set} here.

\begin{athm}{lemma}{lem:AC-level-set}
  Let $I\subseteq\R$ be a compact interval, let $c\in\R$, and let
  $f\colon I\to\R$ be absolutely continuous\dott \cfload[.]
  Then 
  \begin{eqn}
    \llabel{concl1}
    \int_{f^{-1}(\{c\})}\lvert f'(t)\rvert\dd t=0
    \dott
  \end{eqn}
  \cfout[.]
\end{athm}
\begin{aproof}
  Throughout this proof, let $\mu$ denote Lebesgue measure on $\R$,
  let $E=f^{-1}(\{c\})$,
  let $A\subseteq E$ be the set of isolated points of $E$,
  let $B\subseteq\operatorname{int}(I)$ satisfy
  \begin{eqn}
    \llabel{setting1}
    B=\{t\in\operatorname{int}(I)\colon
      f\text{ is not differentiable at }t\},
  \end{eqn}
  and let 
  $Z=A\cup B\cup(I\backslash\operatorname{int}(I))$\dott \cfload[.]
  \startnewargseq
  \argument{
    the fact that $A$ is the set of isolated points of $E$;
  }{
    that for all $t\in A$ there exist $p,q\in\mathbb Q$ with $p<t<q$
    such that
    \begin{eqn}
      \llabel{eq1}
      E\cap(p,q)=\{t\}
      \dott
    \end{eqn}
  }
  \argument{
    \lref{eq1};
    the fact that $\mathbb Q^2$ is countable;
  }{
    that $A$ is at most countable\dott
    \llabel{eq2}
  }
  \argument{
    \lref{eq2};
  }{
    that $A$ is Lebesgue measurable and that
    \begin{eqn}
      \llabel{eq2.1}
      \mu(A)=0
      \dott
    \end{eqn}
  }
  \argument{
    \lref{setting1};
    the fact that $f$ is absolutely continuous;
    the fundamental theorem of calculus for absolutely continuous functions;
  }{
    that $B$ is Lebesgue measurable and that
    \begin{eqn}
      \llabel{eq2.2}
      \mu(B)=0
      \dott
    \end{eqn}
  }
  \argument{
    \lref{eq2.1};
    \lref{eq2.2};
    the fact that $\mu(I\backslash\operatorname{int}(I))=0$;
  }{
    that $Z$ is Lebesgue measurable and that
    \begin{eqn}
      \llabel{eq2.4}
      0\leq\mu(Z)
      \leq\mu(A)+\mu(B)
        +\mu(I\backslash\operatorname{int}(I))=0
      \dott
    \end{eqn}
  }
  We claim that for all $t\in E\backslash Z$ it holds that 
  \begin{eqn}
    \llabel{claim1}
    f'(t)=0
    \dott
  \end{eqn}
  \startnewargseq
  To prove this claim, let $t\in E\backslash Z$\dott \cfload[.]
  \startnewargseq
  \argument{
    the fact that $t\in E\backslash Z$;
  }{
    that $f$ is differentiable at $t$ and there exists 
    $(s_k)_{k\in\N}\subseteq E\backslash\{t\}$ which satisfies that
    \begin{eqn}
      \llabel{eq3}
      \lim_{k\to\infty}s_k=t
      \dott
    \end{eqn}
  }
  \startnewargseq
  \argument{
    \lref{eq3};
    the fact that $f$ is differentiable at $t$;
    the fact that $t\in E$;
  }{
    that
    \begin{eqn}
      \llabel{eq4}
      f'(t)=\lim_{k\to\infty}\frac{f(s_k)-f(t)}{s_k-t}
      =\lim_{k\to\infty}\frac{c-c}{s_k-t}=0
      \dott
    \end{eqn}
  }
  \argument{
    \lref{eq4};
    the fact that $t\in E\backslash Z$ was chosen arbitrarily;
  }[verbs=ep]{
    \lref{claim1}\dott
  }

  \startnewargseq
  \argument{
    the fact that $f$ is absolutely continuous;
    the fundamental theorem of calculus for absolutely continuous functions;
  }{
    that $f'$ is Lebesgue measurable on $I$ and that
    \begin{eqn}
      \llabel{eq5}
      \int_I\lvert f'(t)\rvert\dd t<\infty
      \dott
    \end{eqn}
  }
  \argument{
 the fact that $f$ is continuous;
  }{
    that $E$ is compact and Lebesgue measurable\dott
    \llabel{eq6}
  }
  \argument{
    \lref{eq6};
    \lref{eq2.4};
    the fact that $Z$ is Lebesgue measurable;
  }{
    that $E\backslash Z$ and $E\cap Z$ are Lebesgue measurable and that
    \begin{eqn}
      \llabel{eq7}
      0\leq\mu(E\cap Z)\leq\mu(Z)=0
      \dott
    \end{eqn}
  }
  \argument{
    \lref{claim1};\lref{eq7};
  }{
    that
    \begin{eqn}
      \llabel{eq8}
      \int_{E\backslash Z}\lvert f'(t)\rvert\dd t
      =\int_{E\backslash Z}0\dd t=0
      \dott
    \end{eqn}
  }
  \argument{
    \lref{eq5};
    \lref{eq7};
    the fact that $E\cap Z\subseteq I$ is Lebesgue measurable;
  }{
    that
    \begin{eqn}
      \llabel{eq9}
      \int_{E\cap Z}\lvert f'(t)\rvert\dd t=0
      \dott
    \end{eqn}
  }
  \argument{
    \lref{eq9};
    \lref{eq5};
    \lref{eq8};
  }{
    that
    \begin{eqn}
      \llabel{eq10}
      \int_E\lvert f'(t)\rvert\dd t
      =\int_{E\backslash Z}\lvert f'(t)\rvert\dd t
        +\int_{E\cap Z}\lvert f'(t)\rvert\dd t=0
      \dott
    \end{eqn}
  }
  \argument{
    \lref{eq10};
  }[verbs=ep]{
    \lref{concl1}\dott
  }
\end{aproof}
\subsection{Construction and uniform regularity of the truncation of the transformed energy}
\label{subsec:truncation-estimates}

In the following well-known \cref{lem:product} we
recall two elementary properties of $C^1$-submanifolds, namely, that a
relatively open subset of an $m$-dimensional $C^1$-submanifold is again an
$m$-dimensional $C^1$-submanifold and that the Cartesian product of an
$m$-dimensional $C^1$-submanifold and an $n$-dimensional $C^1$-submanifold
is an $(m+n)$-dimensional $C^1$-submanifold.
Only for completeness, we include a proof for
\cref{lem:product} here.

\begin{athm}{lemma}{lem:product}
  [Products and subsets open in embedded manifolds]
  Let $d,\mathfrak{d}\in\N$, $m,n\in\N_0$, let $X\subseteq\R^d$ be an \dsubmanifold{1}{d}{m}, 
  let $Y\subseteq\R^{\mathfrak d}$ be an \dsubmanifold{1}{\mathfrak{d}}{n},\cfadd{definition:embedded-manifold} and let $U\subseteq X$ be
  open in $X$\dott \cfload[.]
  Then
  \begin{enumerate}[label=(\roman*)]
    \item \llabel{item1} it holds that $U$ is an \dsubmanifold{1}{d}{m} and
    \item \llabel{item2} it holds that $X\times Y$ is an \dsubmanifold{1}{d+\mathfrak{d}}{(m+n)}\dott
  \end{enumerate}
  \cfout[.]
\end{athm}
\begin{aproof}
  To prove \lref{item1}, let $p\in U$\dott \cfload[.]
  \startnewargseq
  \argument{
    the fact that $p\in U\subseteq X$;
    the assumption that $X$ is an \dsubmanifold{1}{d}{m};
  }{
    that there exist an open $W\subseteq\R^m$, a set $O\subseteq X$ open in $X$, and a homeomorphism
    $\phi\in C^1( W,O)$ which satisfy that $p\in O$ and that for all $u\in W$ it holds that
    \begin{eqn}
      \llabel{eq0.1}
      \operatorname{rank}\phi'(u)=m
      \dott
    \end{eqn}
  }
In the following let $G=\phi^{-1}(U)$ and let $\chi\colon G\to\R^d$ satisfy
  $\chi=\phi|_G$\dott \cfload[.]
  \startnewargseq
  \argument{
    the fact that $U$ is open in $X$;
    the fact that $O\subseteq X$;
  }{
    that $O\cap U$ is open in $O$\dott
    \llabel{eq0.2}
  }
  \argument{
    \lref{eq0.2};
    \lref{eq0.1};
    the fact that $\phi\colon W\to O$ is a homeomorphism;
    the fact that $G=\phi^{-1}(U)$;
  }{
    that $G=\phi^{-1}(O\cap U)$ is open in $W$\dott
    \llabel{eq0.3}
  }
  \argument{
    \lref{eq0.3};
    the fact that $W$ is open;
  }{
    that $G$ is open\dott
    \llabel{eq1}
  }
  \argument{
    \lref{eq1};
    \lref{eq0.1};
    the fact that $\chi=\phi|_G$;
    the fact that $\phi\colon W\to O$ is a homeomorphism;
  }{
    that $\chi\in C^1(G,O\cap U)$ is a homeomorphism and that for all
    $u\in G$ it holds that
    \begin{eqn}
      \llabel{eq2}
      \operatorname{rank}\chi'(u)
        =\operatorname{rank}\phi'(u)=m
      \dott
    \end{eqn}
  }
  \argument{
    \lref{eq0.1};
    the fact that $p\in O\cap U$;
    the fact that $G=\phi^{-1}(U)$;
    the fact that $\chi=\phi|_G$;
  }{
    that there exists $w\in G$ which satisfies that
    \begin{eqn}
      \llabel{eq2.0}
      \chi(w)=\phi(w)=p
      \dott
    \end{eqn}
  }
  \startnewargseq
  \argument{
    \lref{eq2};\lref{eq2.0};
    the fact that $G$ is open;
    the fact that $\chi\in C^1(G,\R^d)$ is injective;
    the fact that $\chi(G)=O\cap U$; the fact that $G\ni g\mapsto\chi(g)\in O\cap U$ is a homeomorphism; the fact that $p\in O\cap U$;
    the fact that $O\cap U$ is open in $U$;
    the fact that $p\in U$ was chosen arbitrarily;
  }[verbs=ep]{
    \lref{item1}\dott
  }
  To prove \lref{item2}, let $x\in X$, $y\in Y$\dott \cfload[.]
  \startnewargseq
  \argument{
    the fact that $x\in X$ and $y\in Y$;
    the assumption that $X$ is an \dsubmanifold{1}{d}{m}; the assumption that $Y$ is an \dsubmanifold{1}{\mathfrak{d}}{n};
  }{
    that there exist open $A\subseteq\R^m$ and $B\subseteq\R^n$,
    a set $S\subseteq X$ open in $X$,
    a set $T\subseteq Y$ open in $Y$, and homeomorphisms 
    $\psi\in C^1(A,S)$, $\eta\in C^1( B,T)$ which satisfy that $x\in S$, $y\in T$
   and 
    that for all $a\in A$, $b\in B$ it holds that
    \begin{eqn}
      \llabel{eq2.1}
      \operatorname{rank}\psi'(a)=m
      \qquad\text{and}\qquad
      \operatorname{rank}\eta'(b)=n
      \dott
    \end{eqn}
  }
In the following let $\Theta\colon A\times B\to\R^{d+\mathfrak{d}}$ satisfy for all
  $a\in A$, $b\in B$ that
  \begin{eqn}
    \llabel{setting3}
    \Theta(a,b)=(\psi(a),\eta(b))
    \dott
  \end{eqn}
  \cfload[.]
  \startnewargseq
  \argument{
    \lref{setting3};
    the fact that $\psi\colon A\to S$ and $\eta\colon B\to T$ are
    homeomorphisms;
  }{
    that $\Theta\colon A\times B\to S\times T$ is a homeomorphism and
    that for all $s\in S$, $t\in T$ it holds that
    \begin{eqn}
      \llabel{eq3}
      \Theta^{-1}(s,t)=(\psi^{-1}(s),\eta^{-1}(t))
      \dott
    \end{eqn}
  }
  \argument{
    the fact that $A$ is open;
    the fact that $B$ is open;
    the fact that $S$ is open in $X$;
    the fact that $T$ is open in $Y$;
    the fact that $x\in S$ and $y\in T$;
  }{
    that $A\times B$ is open and that $S\times T$ is open
    in $X\times Y$ and contains $(x,y)$\dott
    \llabel{eq4}
  }
  \argument{
    \lref{eq4};
    \lref{eq2.1};
    \lref{setting3};
    the fact that $\psi\in C^1(A,S)$ and $\eta\in C^1(B,T)$;
  }{
    that $\Theta\in C^1(A\times B,\R^{d+\mathfrak d})$ and that for all $a\in A$, $b\in B$ it holds that
    \begin{eqn}
      \llabel{eq5}
      \Theta'(a,b)
        =\begin{pmatrix}\psi'(a)&0\\0&\eta'(b)\end{pmatrix}
      \qquad\text{and}\qquad
      \operatorname{rank}\Theta'(a,b)=m+n
      \dott
    \end{eqn}
  }
  \argument{
    \lref{eq5};
    \lref{eq3};
    the fact that $A\times B$ is open;
    the fact that $\Theta$ is a homeomorphism; the fact that $(x,y)\in S\times T$;
    the fact that $S\times T$ is open in $X\times Y$;
    the fact that $x\in X$ and $y\in Y$ were chosen arbitrarily;
  }[verbs=ep]{
    \lref{item2}\dott
  }
\end{aproof}

\begin{athm}{lemma}{lem:truncated}
  Let $d\in\N$, $\mathfrak{R},H\in(0,\infty)$, $\varepsilon\in(0,\nicefrac{H}{2})$,
  let $M\subseteq\R^d$ be a \submanifold{1}{d},
  let $\mathcal O\subseteq(0,\mathfrak{R})\times M$ be open
  in $(0,\mathfrak{R})\times M$,
  let $Z\in C(\mathcal O,[0,H])$,
  let $\Omega\subseteq\mathcal O$ satisfy
  \begin{eqn}
    \llabel{ass1}
    \Omega=\{(\rho,y)\in\mathcal O\colon
      0<Z(\rho,y)<H\},
  \end{eqn}
  assume $Z|_{\Omega}\in\COne{\Omega}{[0,H]}$,
  let $\omega\in C((0,\mathfrak{R})\times(0,H),[0,1])$,
  $\mathbf m,\mathfrak b\in C((0,\mathfrak{R}),[0,\infty))$,
  let $P\colon(0,\mathfrak{R})\times[0,H]\to\R$
  satisfy for all $s,t\in(0,\mathfrak{R})$, $z\in[0,H]$ that
  \begin{eqn}
    \llabel{ass2}
    P(s,z)=\int_{(0,z)}\omega(s,u)\dd u
    \qquad\text{and}\qquad
    |P(t,z)-P(s,z)|
    \leq\int_{\min\{s,t\}}^{\max\{s,t\}}\mathfrak b(u)\dd u,
  \end{eqn}
  let $J\subseteq\R$ be an interval,
  let $\varrho\colon J\to(0,\mathfrak{R})$ and $\xi\colon J\to M$
  be locally absolutely continuous, assume
  $\midcup_{t\in J}\{(\varrho(t),\xi(t))\}\subseteq\mathcal O$, let
  $\alpha,\gamma\colon J\to[0,\infty)$ be Lebesgue
  measurable, assume that $J\ni t\mapsto\alpha(t)\gamma(t)\in\R$
  is locally integrable,
  let $v\colon J\to[0,H]$ satisfy for all $t\in J$ that
  \begin{eqn}
    \llabel{ass4}
    v(t)=Z(\varrho(t),\xi(t)),
  \end{eqn}
  assume for Lebesgue-almost all $t\in\{s\in J\colon0<v(s)<H\}$ that
  \begin{eqn}
    \llabel{ass5}
    \omega(\varrho(t),v(t))v'(t)
    \leq\min\bigl\{2\mathbf m(\varrho(t))|\varrho'(t)|,
      5\mathbf m(\varrho(t))|\varrho'(t)|
      -\tfrac{1}{4}\alpha(t)\gamma(t)\bigr\},
  \end{eqn}
  and let
  $\mathbf{Z}\colon\mathcal O\to[\varepsilon,H-\varepsilon]$,
  $q\colon\mathcal O\to\R$,
  $\zeta\colon J\to[\varepsilon,H-\varepsilon]$, and
  $Q\colon J\to\R$ satisfy for all
  $(\rho,y)\in\mathcal O$, $t\in J$ that
  \begin{eqn}
    \llabel{ass6}
    \mathbf{Z}(\rho,y)
      =\max\{\varepsilon,\min\{Z(\rho,y),H-\varepsilon\}\},
    \qquad
    q(\rho,y)=P(\rho,\mathbf{Z}(\rho,y)),
  \end{eqn}
  \begin{eqn}
    \llabel{ass7}
    \zeta(t)=\mathbf{Z}(\varrho(t),\xi(t)),
    \qquad\text{and}\qquad
    Q(t)=q(\varrho(t),\xi(t))
    \dott
  \end{eqn}
  \cfload[.]
  Then 
  \begin{enumerate}[label=(\roman*)]
    \item \llabel{item1} it holds that $\mathbf{Z}$ and
    $q$ are locally Lipschitz continuous,
    \item \llabel{item2} it holds that $\zeta$ and $Q$
    are locally absolutely continuous,
    \item \llabel{item3} it holds for Lebesgue-almost all $t\in J$ that
    \begin{eqn}[i]
      Q'(t)
      \leq\bigl(2\mathbf m(\varrho(t))+\mathfrak b(\varrho(t))\bigr)|\varrho'(t)|,
    \end{eqn}
    \item \llabel{item4} it holds for Lebesgue-almost all
    $t\in\{s\in J\colon\varepsilon<v(s)<H-\varepsilon\}$ that
    \begin{eqn}
      \llabel{concl2}
      Q'(t)
      \leq\bigl(5\mathbf m(\varrho(t))+\mathfrak b(\varrho(t))\bigr)|\varrho'(t)|
        -\tfrac{1}{4}\alpha(t)\gamma(t),
    \end{eqn}
    and
    \item \llabel{item5} it holds for Lebesgue-almost all $t\in J$ that
    \begin{eqn}
      \llabel{concl3}
      \bigl|Q'(t)-\omega(\varrho(t),\zeta(t))\zeta'(t)\bigr|
        \leq\mathfrak b(\varrho(t))|\varrho'(t)|
      \dott
    \end{eqn}
  \end{enumerate}
  \cfout[.]
\end{athm}
\begin{aproof}
  Throughout this proof let $\mu$ denote Lebesgue measure on $\R$,
  let $\mathfrak B\colon(0,\mathfrak{R})\to\R$ satisfy that for all
  $\rho\in(0,\mathfrak{R})$ it holds that
  \begin{eqn}
    \llabel{setting1}
    \mathfrak B(\rho)=\int_{\mathfrak{R}/2}^{\rho}\mathfrak b(u)\dd u,
  \end{eqn}
  let $\mathcal A=\{t\in J\colon
    \varepsilon<v(t)<H-\varepsilon\}$, and let
  $\mathcal D\subseteq\operatorname{int}(J)$ satisfy
  \begin{eqn}
    \llabel{setting2}
    \mathcal D
      =\{t\in\operatorname{int}(J)\colon
        \varrho,\zeta,Q
        \text{ are differentiable at }t\}
    \dott
  \end{eqn}
  \cfload[.]
  \startnewargseq
  \argument{
    \lref{ass2};
    the fact that $0\leq\omega\leq1$;
    the triangle inequality;
  }{
    that for all $\rho\in(0,\mathfrak{R})$, $\mathfrak z_1,\mathfrak z_2\in[0,H]$ it holds that
    \begin{eqn}
      \llabel{eq1}
      |P(\rho,\mathfrak z_2)-P(\rho,\mathfrak z_1)|\leq|\mathfrak z_2-\mathfrak z_1|
      \dott
    \end{eqn}
  }
  \argument{
    \lref{eq1};
    \lref{ass2};
    the triangle inequality;
  }{
    that for all $s,t\in(0,\mathfrak{R})$, $\mathfrak z_1,\mathfrak z_2\in[0,H]$ it holds that
    \begin{eqn}[a]
      \llabel{eq1.1}
      |P(t,\mathfrak z_2)-P(s,\mathfrak z_1)|
      &\leq|P(t,\mathfrak z_2)-P(s,\mathfrak z_2)|+|P(s,\mathfrak z_2)-P(s,\mathfrak z_1)|\\
      &\leq \biggl(\int_{\min\{s,t\}}^{\max\{s,t\}}\mathfrak b(w)\dd w\biggr)+|\mathfrak z_2-\mathfrak z_1|
      \dott
    \end{eqn}
  }
  \argument{
    \lref{eq1.1};
    the fact that $\mathfrak b$ is continuous;
  }{
    that for all $\rho_1,\rho_2\in(0,\mathfrak{R})$ with $\rho_1<\rho_2$
    there exists $L\in[1,\infty)$ which satisfies for all
    $s,t\in[\rho_1,\rho_2]$, $\mathfrak z_1,\mathfrak z_2\in[0,H]$ that
    \begin{eqn}
      \llabel{eq2}
      |P(t,\mathfrak z_2)-P(s,\mathfrak z_1)|
      \leq L\sqrt{(t-s)^2+(\mathfrak z_2-\mathfrak z_1)^2}
      \dott
    \end{eqn}
  }
  \startnewargseq
  \argument{
    \lref{ass1};
    the fact that $Z$ is continuous;
    the fact that $\Omega=Z^{-1}((0,H))$;
    the fact that $\mathcal O$ is open in $(0,\mathfrak{R})\times M$;
    the fact that $(0,\mathfrak{R})$ is a \dsubmanifold{1}{1}{1};
    the assumption that $M$ is a \submanifold{1}{d};
    \cref{lem:product};
  }{
    that $\Omega$ is open in $\mathcal O$
    and is a \submanifold{1}{d+1}\dott
    \llabel{eq2.1}
  }
  \argument{
    \lref{eq2.1};
    \cref{lem:embedded-manifold-local-straightening};
  }{
    that for every $p\in\Omega$ there exist
    $n\in\{0,1,\dots,d+1\}$, open  $A,B\subseteq\R^{d+1}$, and a $C^1$ diffeomorphism
    $\Psi\colon A\to B$ such that
    \begin{eqn}
      \llabel{eq2.1a}
      p\in A,
      \qquad 0\in B,
      \qquad \Psi(p)=0,
      \qquad\text{and}\qquad
      \Psi(\Omega\cap A)
        =B\cap\bigl(\R^n\times\{0\}^{d+1-n}\bigr)
      \dott
    \end{eqn}
  }
  We claim that for all $p\in\Omega$ there exist an open 
  $U\subseteq\R^{d+1}$ 
  and $F\in C^1(U,\R)$ which satisfy that
  \begin{eqn}
    \llabel{eq2.2}
    p\in U\qquad\text{and}\qquad F|_{U\cap\Omega}=Z|_{U\cap\Omega}
    \dott
  \end{eqn}
  \startnewargseq
  To prove this claim, let $p\in\Omega$, let $n\in\{0,1,\dots,d+1\}$,
  let $A\subseteq\R^{d+1}$ be open, let $B\subseteq\R^{d+1}$ be open,
  let $\Psi\colon A\to B$ be a $C^1$ diffeomorphism, and assume
  \begin{eqn}
    \llabel{setting2.1}
    p\in A,
    \qquad 0\in B,
    \qquad \Psi(p)=0,
    \qquad\text{and}\qquad
    \Psi(\Omega\cap A)
      =B\cap\bigl(\R^n\times\{0\}^{d+1-n}\bigr)
    \dott
  \end{eqn}
  \cfload[.]
  \startnewargseq
  \argument{
    the fact that $0\in B$;
    the fact that $B$ is open;
  }{
    that there exists $\delta\in(0,\infty)$ which satisfies that
    \begin{eqn}
      \llabel{eq2.2a}
      (-\delta,\delta)^{d+1}\subseteq B
      \dott
    \end{eqn}
  }
  \startnewargseq
  \argument{
    \lref{setting2.1};
    \lref{eq2.2a};
  }{
    that for all $a\in(-\delta,\delta)^n$ it holds that
    \begin{eqn}
      \llabel{eq2.2b}
      (a,0)\in B\cap\bigl(\R^n\times\{0\}^{d+1-n}\bigr)
      \qquad\text{and}\qquad
      \Psi^{-1}(a,0)\in\Omega\cap A
      \dott
    \end{eqn}
  }
  In the following let $U=\Psi^{-1}((-\delta,\delta)^{d+1})$,
  let $\pi\colon\R^{d+1}\to\R^n$ satisfy for all
  $a\in\R^n$, $b\in\R^{d+1-n}$ that
  \begin{eqn}
    \llabel{setting2.2}
    \pi(a,b)=a,
  \end{eqn}
  and let $g\colon(-\delta,\delta)^n\to\R$ satisfy for all $a\in(-\delta,\delta)^n$ that
  \begin{eqn}
    \llabel{setting3}
    g(a)=Z\bigl(\Psi^{-1}(a,0)\bigr)
    \dott
  \end{eqn}
  \cfload[.]
  \startnewargseq
  \argument{
    \lref{setting2.1};
    \lref{eq2.2a}; the fact that $U=\Psi^{-1}((-\delta,\delta)^{d+1})$;
    the fact that $\Psi\colon A\to B$ is a homeomorphism;
    the fact that $A$ is open;
  }{
    that $U$ is open and that
    \begin{eqn}
      \llabel{eq2.2c}
      p\in U
      \qquad\text{and}\qquad
      \pi(\Psi(U))=(-\delta,\delta)^n
      \dott
    \end{eqn}
  }
  \argument{
    \lref{setting2.1};
    \lref{eq2.2a};
    the fact that $U=\Psi^{-1}((-\delta,\delta)^{d+1})$;
    the fact that $\Psi\colon A\to B$ is bijective;
  }{
    that
    \begin{eqn}[a]
      \llabel{eq2.2c.1}
      \Psi(\Omega\cap U)
        =\Psi(\Omega\cap A)\cap(-\delta,\delta)^{d+1}
        =(-\delta,\delta)^n\times\{0\}^{d+1-n}
      \dott
    \end{eqn}
  }
  In the following let $\mathbf{T}\colon(-\delta,\delta)^n\to\Omega\cap U$
  satisfy that for all $a\in(-\delta,\delta)^n$ it holds that
  \begin{eqn}
    \llabel{setting3.1}
    \mathbf{T}(a)=\Psi^{-1}(a,0)
    \dott
  \end{eqn}
  \cfload[.]
  \startnewargseq
  \argument{
    \lref{setting2.1};
    \lref{setting2.2};
    \lref{eq2.2c.1};
    \lref{setting3.1};
    the fact that $\Psi\colon A\to B$ is a homeomorphism;
    the fact that $\R^n\ni a\mapsto(a,0)\in\R^{d+1}$ is continuous;
    the fact that $\pi$ is continuous;
  }{
    that $\mathbf{T}\colon(-\delta,\delta)^n\to\Omega\cap U$ is a
    homeomorphism and that for all $y\in\Omega\cap U$ it holds that
    \begin{eqn}
      \llabel{eq2.2c.2}
      \mathbf{T}^{-1}(y)=\pi(\Psi(y))
      \qquad\text{and}\qquad
      \mathbf{T}(0)=\Psi^{-1}(0)=p
      \dott
    \end{eqn}
  }
  \argument{
    \lref{eq2.2a};
    \lref{setting3.1};
    the fact that $\Psi^{-1}\in C^1(B,A)$;
    the fact that for all $z\in B$ it holds that $(\Psi^{-1})'(z)$
    is invertible;
    the fact that $\R^n\ni a\mapsto(a,0)\in\R^{d+1}$ is an
    injective linear function;
    the chain rule;
  }{
    that $\mathbf{T}\in C^1((-\delta,\delta)^n,\R^{d+1})$ and that for all
    $a\in(-\delta,\delta)^n$, $h\in\R^n$ it holds that
    \begin{eqn}[a]
      \llabel{eq2.2c.3}
      \mathbf{T}'(a)h=(\Psi^{-1})'(a,0)(h,0)\qquad\text{and}\qquad
      \operatorname{rank}\mathbf{T}'(a)=n
      \dott
    \end{eqn}
  }
  \argument{
    \lref{eq2.2c.3};
    \lref{setting3};
    \lref{setting3.1};
    the fact that $U$ is open;
    the fact that $\Omega\cap U$ is open in $\Omega$;
    the fact that $\mathbf{T}\colon(-\delta,\delta)^n\to\Omega\cap U$
    is a homeomorphism;
    the assumption that $Z|_{\Omega}\in\COne{\Omega}{[0,H]}$;
    \cref{definition:c1-on-manifold.eq1};
  }{
    that
    \begin{eqn}
      \llabel{eq2.2d}
      g=(Z|_{\Omega})\circ\mathbf{T}
        \in C^1((-\delta,\delta)^n,\R)
      \dott
    \end{eqn}
  }
  \argument{
    \lref{eq2.2d};
    \lref{eq2.2c};
    the fact that $U$ is open;
    the fact that $\Psi\in C^1(A,B)$;
    the fact that $\pi$ is linear;
  }{
    that there exists $F\in C^1(U,\R)$ which
    satisfies for all $y\in U$ that
    \begin{eqn}
      \llabel{eq2.2e}
      F(y)=g\bigl(\pi(\Psi(y))\bigr)
      \dott
    \end{eqn}
  }
  \startnewargseq
  \argument{
    \lref{setting2.1};
    \lref{setting3};
     \lref{eq2.2e};
    the fact that $U\subseteq A$;
    the fact that for all $y\in\Omega\cap U$ it holds that
    $\Psi(y)=(\pi(\Psi(y)),0)$;
  }{
    that for all $y\in\Omega\cap U$ it holds that
    \begin{eqn}[a]
      \llabel{eq2.2f}
      F(y)&=Z\bigl(\Psi^{-1}(\pi(\Psi(y)),0)\bigr)=Z\bigl(\Psi^{-1}(\Psi(y))\bigr)=Z(y)
      \dott
    \end{eqn}
  }
  \argument{
    \lref{eq2.2f};
    \lref{eq2.2c};
    the fact that $U$ is open;
    the fact that $F\in C^1(U,\R)$;
    the fact that $p\in\Omega$ was chosen arbitrarily;
  }[verbs=ep]{
    \lref{eq2.2}\dott
  }

  \startnewargseq
  \argument{
    \lref{eq2.2};
  }{
    that $Z|_{\Omega}$ is locally Lipschitz continuous\dott
    \llabel{eq2.3}
  }
  \argument{
    \lref{eq2.3};
    \lref{ass6};
    the fact that 
    $\R\ni u\mapsto\max\{\varepsilon,\min\{u,H-\varepsilon\}\}\in\R$
    is Lipschitz continuous;
  }{
    that $\mathbf{Z}|_{\Omega}$ is locally Lipschitz continuous\dott
  }
  \argument{
    \lref{ass6};
    the assumption that $Z$ is continuous;
  }{
    that $Z^{-1}((-\infty,\varepsilon))$ and
    $Z^{-1}((H-\varepsilon,\infty))$ are open in
    $\mathcal O$ and that for all
    $(\rho,y)\in Z^{-1}((-\infty,\varepsilon))
      \cup Z^{-1}((H-\varepsilon,\infty))$ it holds that
    \begin{eqn}
      \llabel{eq3}
      \mathbf{Z}(\rho,y)=
      \begin{cases}
        \varepsilon & \colon Z(\rho,y)<\varepsilon,\\
        H-\varepsilon & \colon Z(\rho,y)>H-\varepsilon
      \end{cases}
      \dott
    \end{eqn}
  }
  \argument{
    \lref{eq3};
    the fact that $\mathbf{Z}|_{\Omega}$ is locally Lipschitz continuous;
    the fact that
    $\Omega\cup Z^{-1}((-\infty,\varepsilon))
      \cup Z^{-1}((H-\varepsilon,\infty))=\mathcal O$;
  }{
    that $\mathbf{Z}$ is locally Lipschitz continuous\dott
    \llabel{eq4}
  }
  \argument{
    \lref{eq4};
     \lref{ass6};
    \lref{eq2};
  }[verbs=ep]{
    \lref{item1}\dott
  }
  \startnewargseq
  To prove \lref{item2}, let $I\subseteq J$ be a compact interval and let
  $\mathcal K=\{(\varrho(t),\xi(t))\colon t\in I\}$\dott \cfload[.]
  \startnewargseq
  \argument{
    the fact that $\varrho|_I$ and $\xi|_I$ are continuous;
    the assumption that $\midcup_{t\in J}\{(\varrho(t),\xi(t))\}\subseteq\mathcal O$;
  }{
    that $\mathcal K$ is a compact subset of $\mathcal O$\dott
    \llabel{eq4.1}
  }
  \argument{
    \lref{eq4.1};
    \lref{item1};
    \cite[Proposition~A.48]{Lee13};
  }{
    that $\mathbf{Z}|_{\mathcal K}$ and $q|_{\mathcal K}$ are Lipschitz
    continuous and that there exists $\Lambda\in(0,\infty)$ which
    satisfies for all $x,y\in\mathcal K$ that
    \begin{eqn}
      \llabel{eq4.2}
      |\mathbf{Z}(x)-\mathbf{Z}(y)|+|q(x)-q(y)|
        \leq\Lambda|x-y|
      \dott
    \end{eqn}
  }
  \argument{
    \lref{eq4.2};
    \lref{ass7};
    the fact that for all $s,t\in I$ it holds that
    $|(\varrho(t),\xi(t))-(\varrho(s),\xi(s))|
      \leq|\varrho(t)-\varrho(s)|+|\xi(t)-\xi(s)|$;
  }{
    that for all $N\in\N$, $a_1,a_2,\dots,a_N,b_1,b_2,\dots,b_N\in I$
    with $a_1<b_1\leq a_2<b_2\leq\ldots\leq a_N<b_N$ it holds that
    \begin{eqn}[a]
      \llabel{eq4.3}
      \sum_{j=1}^N
        \bigl(|\zeta(b_j)-\zeta(a_j)|+|Q(b_j)-Q(a_j)|\bigr)\leq\Lambda\sum_{j=1}^N
        \bigl(|\varrho(b_j)-\varrho(a_j)|+|\xi(b_j)-\xi(a_j)|\bigr)
      \dott
    \end{eqn}
  }
  \argument{
    \lref{eq4.3};
    the assumption that $\varrho|_I$ and $\xi|_I$ are absolutely continuous;
  }{
    that $\zeta|_I$ and $Q|_I$ are absolutely continuous\dott
    \llabel{eq4.4}
  }
  \argument{
    \lref{eq4.4};
    the fact that the compact interval $I\subseteq J$ was chosen arbitrarily;
  }[verbs=ep]{
    \lref{item2}\dott
  }
  \startnewargseq
  \argument{
    \lref{setting1};
    the assumption that $\mathfrak b$ is continuous and nonnegative;
    the fundamental theorem of calculus;
  }{
    that $\mathfrak B\in C^1((0,\mathfrak R),\R)$ and that for all $\rho,s,t\in(0,\mathfrak{R})$ it holds that
    \begin{eqn}
      \llabel{eq5}
      \mathfrak B'(\rho)=\mathfrak b(\rho)
      \qquad\text{and}\qquad
      \int_{\min\{s,t\}}^{\max\{s,t\}}\mathfrak b(u)\dd u
      =|\mathfrak B(t)-\mathfrak B(s)|
      \dott
    \end{eqn}
  }
  \argument{
    \lref{item2};
    the fact that $\varrho$ is locally absolutely continuous;
    the fundamental theorem of calculus for absolutely continuous functions;
  }{
    that for every compact interval $I\subseteq J$ it holds that
    for Lebesgue-almost all $t\in\operatorname{int}(I)$ it holds that
    $\varrho$, $\zeta$, and $Q$ are differentiable at $t$\dott
    \llabel{eq5.1}
  }
  \argument{
    \lref{eq5.1};
    \lref{setting2};
    the fact that $\operatorname{int}(J)$ is covered by the interiors of
    countably many compact intervals contained in $J$;
    the fact that $J\backslash\operatorname{int}(J)$ has at most two points;
  }{
    that $\mathcal D$ is Lebesgue measurable and that
    \begin{eqn}
      \llabel{eq6}
      \mu(J\backslash\mathcal D)=0
      \dott
    \end{eqn}
  }
  We claim that for all $t\in\mathcal D$ it holds that
  \begin{eqn}
    \llabel{claim1}
    \bigl|Q'(t)
      -\omega(\varrho(t),\zeta(t))\zeta'(t)\bigr|
    \leq\mathfrak b(\varrho(t))|\varrho'(t)|
    \dott
  \end{eqn}
  \startnewargseq
  To prove this claim, let $t\in\mathcal D$\dott \cfload[.]
  \startnewargseq
  \argument{
    \lref{ass6};
    \lref{ass7};
  }{
    that for all $h\in(0,\infty)$ with $t+h\in J$ it holds that
    \begin{eqn}[a]
      \llabel{eq7}
      Q(t+h)-Q(t)
      &=P(\varrho(t+h),\zeta(t+h))
        -P(\varrho(t),\zeta(t+h))\\
      &\quad+P(\varrho(t),\zeta(t+h))
        -P(\varrho(t),\zeta(t))
      \dott
    \end{eqn}
  }
  \argument{
    \lref{eq7};
    \lref{ass2};
    \lref{eq5};
  }{
    that for all $h\in(0,\infty)$ with $t+h\in J$ it holds that
    \begin{eqn}[a]
      \llabel{eq8}
      \Biggl|\frac{Q(t+h)-Q(t)}h
        -\frac{P(\varrho(t),\zeta(t+h))
          -P(\varrho(t),\zeta(t))}h\Biggr|\leq\frac{|\mathfrak B(\varrho(t+h))-\mathfrak B(\varrho(t))|}h
      \dott
    \end{eqn}
  }
  \argument{
    \lref{ass2};
    the fact that $\omega$ is continuous;
    the fundamental theorem of calculus;
  }{
    that for all $(\rho,z)\in(0,\mathfrak{R})\times(0,H)$ it holds that
    \begin{eqn}
      \llabel{eq9}
      \bigl(\partial_z P\bigr)\bigl(\rho,z\bigr)=\omega(\rho,z)
      \dott
    \end{eqn}
  }
  \argument{
    \lref{eq9};
    the fact that $t\in\mathcal D$;
    the fact that $\zeta(t)\in[\varepsilon,H-\varepsilon]
      \subseteq(0,H)$;
    the chain rule;
  }{
    that
    \begin{eqn}
      \llabel{eq10}
      \lim_{h\to 0}
      \frac{P(\varrho(t),\zeta(t+h))
        -P(\varrho(t),\zeta(t))}h
      =\omega(\varrho(t),\zeta(t))\zeta'(t)
      \dott
    \end{eqn}
  }
  \argument{
    \lref{eq5};
    the fact that $t\in\mathcal D$; the fact that $\varrho$ is continuous;
    the chain rule;
  }{
    that
    \begin{eqn}
      \llabel{eq11}
      \lim_{h\searrow 0}\frac{|\mathfrak B(\varrho(t+h))-\mathfrak B(\varrho(t))|}h
      =\mathfrak b(\varrho(t))|\varrho'(t)|
      \dott
    \end{eqn}
  }
  \argument{
    \lref{eq11};
    \lref{eq8};
    \lref{eq10};
    the fact that $Q$ is differentiable at $t$;
  }{
    that
    \begin{eqn}
      \llabel{eq12}
      \bigl|Q'(t)
        -\omega(\varrho(t),\zeta(t))\zeta'(t)\bigr|
      \leq\mathfrak b(\varrho(t))|\varrho'(t)|
      \dott
    \end{eqn}
  }
  \argument{
    \lref{eq12};
    the fact that $t\in\mathcal D$ was chosen arbitrarily;
  }[verbs=ep]{
    \lref{claim1}\dott
  }

  \startnewargseq
  \argument{
    \lref{ass4};
    \lref{ass6};
    \lref{ass7};
    the fact that $Z$, $\varrho$, and $\xi$ are continuous;
    the fact that $\mathcal A=\{t\in J\colon
    \varepsilon<v(t)<H-\varepsilon\}$;
  }{
    that $\mathcal A$ is open in $J$ and that for all
    $t\in\mathcal A$ it holds that
    \begin{eqn}
      \llabel{eq13}
      \zeta(t)=v(t)
      \dott
    \end{eqn}
  }
  \argument{
    \lref{eq13};
    \lref{ass5};
    the fact that $\mathcal A$ is open in $J$;
  }{
    that for Lebesgue-almost all $t\in\mathcal A$ it holds that
    \begin{eqn}[a]
      \llabel{eq14}
      \omega(\varrho(t),\zeta(t))\zeta'(t)
      \leq2\mathbf m(\varrho(t))|\varrho'(t)|\qquad\text{and}\qquad
      \omega(\varrho(t),\zeta(t))\zeta'(t)
      \leq-\frac{\alpha(t)\gamma(t)}4+5\mathbf m(\varrho(t))|\varrho'(t)|
      \dott
    \end{eqn}
  }
  \argument{
    \lref{eq14};
    \lref{eq6};
    \lref{claim1};
  }{
    that for Lebesgue-almost all $t\in\mathcal A$ it holds that
    \begin{eqn}[a]
      \llabel{eq15}
      Q'(t)
      \leq\bigl(2\mathbf m(\varrho(t))+\mathfrak b(\varrho(t))\bigr)|\varrho'(t)|\qquad\text{and}\qquad
      Q'(t)\leq-\frac{\alpha(t)\gamma(t)}4
        +\bigl(5\mathbf m(\varrho(t))+\mathfrak b(\varrho(t))\bigr)|\varrho'(t)|
      \dott
    \end{eqn}
  }
  \argument{
    \lref{ass4};
    \lref{ass6};
    \lref{ass7};
     the fact that $\mathcal A=\{t\in J\colon
    \varepsilon<v(t)<H-\varepsilon\}$;
    the fact that $\varepsilon<H-\varepsilon$;
  }{
    that
    \begin{eqn}
      \llabel{eq16}
      J\backslash\mathcal A
      =\{t\in J\colon\zeta(t)=\varepsilon\}
        \cup\{t\in J\colon\zeta(t)=H-\varepsilon\}
      \dott
    \end{eqn}
  }
  \argument{
    \lref{eq16};
    \lref{item2};
    \cref{lem:AC-level-set}
  }{
    that for every compact interval $I\subseteq J$ it holds that
    \begin{eqn}[a]
      \llabel{eq16.1}
      \int_{\{t\in I\colon\zeta(t)=\varepsilon\}}
        |\zeta'(t)|\dd t=0
      \qquad\text{and}\qquad
      \int_{\{t\in I\colon\zeta(t)=H-\varepsilon\}}
        |\zeta'(t)|\dd t=0
      \dott
    \end{eqn}
  }
  \argument{
    \lref{eq16.1};
    \lref{eq16};
    the fact that $\zeta$ is continuous;
    the fundamental theorem of calculus for absolutely continuous functions;
    the assumption that $J$ is an interval;
  }{
    that for Lebesgue-almost all $t\in J\backslash\mathcal A$
    it holds that
    \begin{eqn}
      \llabel{eq17}
      \zeta'(t)=0
      \dott
    \end{eqn}
  }
  \argument{
    \lref{eq17};
    \lref{eq6};
    \lref{claim1};
    the fact that $\mathbf m$ is nonnegative;
  }{
    that for Lebesgue-almost all $t\in J\backslash\mathcal A$
    it holds that
    \begin{eqn}
      \llabel{eq18}
      Q'(t)
      \leq\mathfrak b(\varrho(t))|\varrho'(t)|
      \leq\bigl(2\mathbf m(\varrho(t))+\mathfrak b(\varrho(t))\bigr)|\varrho'(t)|
      \dott
    \end{eqn}
  }
  \argument{
    \lref{eq18};
    \lref{setting2};
    \lref{eq6};
    \lref{eq15};
    the fact that $J=\mathcal A\cup(J\backslash\mathcal A)$;
  }[verbs=ep]{
    \lref{item3}\dott
  }
  \startnewargseq
  \argument{
    \lref{setting2};
    \lref{eq6};
    \lref{eq15};
    the fact that $\mathcal A=\{t\in J\colon
    \varepsilon<v(t)<H-\varepsilon\}$;
  }[verbs=ep]{
    \lref{item4}\dott
  }
  \startnewargseq
  \argument{
    \lref{setting2};
    \lref{eq6};
    \lref{claim1};
  }[verbs=ep]{
    \lref{item5}\dott
  }
\end{aproof}

\section{Lyapunov-type stability for the radial-angular differential equation (DE) system}
\label{sec:radial-integration}

In \cref{maintheorem3} in \cref{sec:applications} below we establish that, on an
arbitrary time interval $I\subseteq\R$ with non-empty interior, the
spherical projection of every \GF\ trajectory of a definable $C^1$-potential
has finite length.
One step in our proof of \cref{maintheorem3} is to perform an
appropriate coordinate transform
(polar coordinates, at infinity combined with the inversion method $\R^n\backslash\{0\}\ni x\mapsto\norm{x}^{-2}x\in\R^n\backslash\{0\}$)
to the \GF\ \DE\ to obtain an appropriate transformed \DE\
system that consists of a radial component
(which is $ r(t) = \norm{\theta(t)-\vartheta}^\nu $ for $ \nu \in \{ -1, 1 \} $)
and the angular component (which is the spherical projection) that takes values in
the sphere of $\R^n$ (see \eqref{prop:spherical-length-away.eq4} and
\eqref{prop:spherical-length-away.eq5} in the proof of
\cref{prop:spherical-length-away} for details).

It is precisely the subject of this section to study the dynamics of such
radial-angular \DE\ systems in an abstract form, where a solution of the
\DE\ system $(r(t),x(t))\in\mathcal O$, $t\in I$, takes values in a definable
set $\mathcal O\subseteq(0,\mathfrak R)\times M$ which is open in
$(0,\mathfrak R)\times M$, where $\mathfrak R\in(0,\infty)$ is a sufficiently
small real number, where $M$ is a compact and definable $C^2$-submanifold of
$\R^D$, and where $D\in\N$.

In particular, in \itref{cor:angular-length}{item2} in
\cref{subsec:product-length} below we prove that the length of the angular
variable $x(t)$, $t\in I$, of such radial-angular \DE\ systems
$(r(t),x(t))\in\mathcal O$, $t\in I$, is finite
provided that the radial variable $r(t)$ converges to zero as $t$ tends to
$\sup(I)$ and provided that the time interval $I$ agrees with the half-open
interval $I=[\tau,T)$ for some $T\in(\tau,\infty]$, $\tau\in\R$.
\Cref{cor:angular-length.item2} in \cref{cor:angular-length} is applied in the
proof of \cref{prop:spherical-length-away} to establish finite length of the
spherical projection of the \GF\ trajectory and, thereby, to complete the
proof of \cref{maintheorem3}.
Our proof of \cref{cor:angular-length} is based on applications of
\cref{prop:buildQ} in \cref{subsec:bounded-state} below,
\cref{lem:radial-integrability} in \cref{subsec:radial-integrability} below,
and \cref{prop:cone-annular-length} in \cref{subsec:cone-annular-occupation} below.

On a technical level, the most fundamental contribution of this work is
\cref{prop:buildQ}, which establishes a certain Lyapunov-type stability
for solutions $(r(t),x(t))\in\mathcal O$, $t\in I$, of the radial-angular
\DE\ systems.
Specifically, \eqref{prop:buildQ.concl2} in
\cref{prop:buildQ} shows that for all $a,b\in I$ with
$a\leq b$ it holds that
\begin{eqn}
  \label{eq:radial-angular-lyapunov}
  Q(r(b),x(b))-Q(r(a),x(a))
    \leq-\int_a^b\norm{x'(t)}\dd t
      +\int_a^b m(r(t))\abs{r'(t)}\dd t
  \dott
\end{eqn}
The function $m$ in \eqref{eq:radial-angular-lyapunov} is a nonnegative,
continuous, integrable, and definable function and the function $Q$ in
\eqref{eq:radial-angular-lyapunov} is a nonnegative, continuous, and bounded
function.
We note that \eqref{eq:radial-angular-lyapunov} shows that $Q$ serves as a Lyapunov-type
function for the radial-angular \DE\ system.
If $c\in[0,\infty)$ is an upper bound for $Q$, then rearranging
\eqref{eq:radial-angular-lyapunov} shows that for all $a,b\in I$ with
$a\leq b$ it holds that
\begin{eqn}
  \label{eq:angular-length-from-lyapunov}
  \int_a^b\norm{x'(t)}\dd t
    \leq 2c+\int_a^b m(r(t))\abs{r'(t)}\dd t
  \dott
\end{eqn}
In the proof of \cref{cor:angular-length} we employ
\cref{lem:radial-integrability} and \cref{prop:cone-annular-length} to show,
under suitable assumptions, that the integral on the right hand side of
\eqref{eq:angular-length-from-lyapunov} can be uniformly bounded.
Then, combining this with the consequence of the Lyapunov-type estimate in
\cref{prop:buildQ} in
\eqref{eq:angular-length-from-lyapunov} allows us to establish the finite
length estimate $\int_\tau^T\norm{x'(t)}\dd t<\infty$ in
\itref{cor:angular-length}{item2}.

\subsection{Regularity properties for one-parameter primitives of definable two-parameter functions}
\label{sec:primitive}

In the next result, \cref{lem:primitive} below, we establish some basic
properties for the one-parameter primitive (see \eqref{lem:primitive.ass1} below)
of a definable and continuously
differentiable two-parameter function.
We employ \cref{lem:primitive} in our proof of the Lyapunov-type estimate
in \cref{prop:buildQ} in \cref{subsec:bounded-state} below.
\Cref{prop:buildQ}, in turn, is used to establish the finite
length estimate for the angular component in \cref{cor:angular-length} in
\cref{subsec:product-length} below.
\Cref{cor:angular-length}, in turn, is employed in
\cref{prop:spherical-length-away} in \cref{sec:applications} to prove the
main result of this work; see \cref{maintheorem3} in
\cref{sec:applications}.

\begin{athm}{lemma}{lem:primitive}
  Let $\mathfrak{R},H\in(0,\infty)$, let
  $\omega\in C^1((0,\mathfrak{R})\times(0,H),[0,1])$ be definable,
  and let $P\colon(0,\mathfrak{R})\times[0,H]\to\R$ satisfy for all $r\in (0,\mathfrak{R})$,
  $z\in[0,H]$ that
  \begin{equation}
    \llabel{ass1}
    P(r,z)=\int_0^z\omega(r,u)\dd u.
  \end{equation}
  \cfload[.]
  Then there exist $\varrho\in(0,\mathfrak{R})$ and a definable
  $b\in C((0,\varrho),[0,\infty))$ such that
  \begin{enumerate}[label=(\roman*)]
    \item \llabel{item1} for all $r\in (0,\mathfrak{R})$, $z\in[0,H]$
    it holds that $P(r,0)=0$ and $0\leq P(r,z)\leq z\leq H$,
    \item \llabel{item2} for all $r\in(0,\mathfrak{R})$,
    $ z_1, z_2\in[0,H]$ it holds that
    \begin{eqn}[i]
      |P(r,z_1)-P(r,z_2)|\leq|z_1-z_2|,
    \end{eqn}
    \item \llabel{item3} for all $r\in(0,\mathfrak{R})$ it holds that
    $(0,H)\ni z\mapsto P(r,z)\in\R$ is differentiable,
    \item \llabel{item4} for all $r\in(0,\mathfrak{R})$, $z\in(0,H)$
    it holds that
    \begin{eqn}[i]
      \partial_z P(r,z)=\omega(r,z),
    \end{eqn}
    \item \llabel{item5} for all $s,t\in(0,\varrho)$, $z\in[0,H]$
    it holds that
    \begin{equation}
      \llabel{concl4}
      |P(s,z)-P(t,z)|\leq
      \biggl|\int_{s}^{t}b(r)\dd r\biggr|\leq \int_0^{\varrho}b(r)\dd r<\infty,
    \end{equation}
    and
    \item \llabel{item6} for all $\rho_1,\rho_2\in(0,\varrho)$ with
    $\rho_1<\rho_2$ there exists $L\in\R$
    which satisfies for all $s_1,s_2\in [\rho_1,\rho_2]$,
    $z_1, z_2\in [0,H]$ that
    \begin{equation}
      \llabel{concl6}
      |P(s_1,z_1)-P(s_2,z_2)|
      \leq L\sqrt{(s_1-s_2)^2+(z_1-z_2)^2}
      \dott
    \end{equation}
  \end{enumerate}
  \cfout[.]
\end{athm}
\begin{aproof}
  \argument{
    \lref{ass1};
    the fact that $\omega$ is continuous;
    the fact that for all $(r,u)\in(0,\mathfrak{R})\times(0,H)$ it holds
    that $0\leq\omega(r,u)\leq1$;
  }[verbs=ep]{
    \lref{item1}\dott
  }
  \startnewargseq
  \argument{
    \lref{ass1};
    the fact that for all $(r,u)\in(0,\mathfrak{R})\times(0,H)$ it holds
    that $0\leq\omega(r,u)\leq1$;
    the triangle inequality;
  }{
    that for all $r\in(0,\mathfrak{R})$, $a,b\in[0,H]$ with $a<b$
    it holds that
    \begin{equation}
      \llabel{eq1}
      |P(r,a)-P(r,b)|
      =\biggl|\int_a^b\omega(r,u)\dd u\biggr|\leq b-a=|a-b|
      \dott
    \end{equation}
  }
  \argument{
    \lref{eq1};
  }[verbs=ep]{
    \lref{item2}\dott
  }
  \startnewargseq
  \argument{
    \lref{ass1};
    the fact that $\omega$ is continuous;
    the fundamental theorem of calculus;
  }[verbs=ep]{
    \lref{item3}\dott
  }
  \argument{
    \lref{item3};
    \lref{ass1};
    the fact that $\omega$ is continuous;
    the fundamental theorem of calculus;
  }[verbs=ep]{
    \lref{item4}\dott
  }
  In the following let $W=(0,\mathfrak{R})\times(0,H)$ and let
  $\sigma\colon(0,\mathfrak{R})\to[0,\infty]$ satisfy for all
  $r\in(0,\mathfrak{R})$ that
  \begin{equation}
    \llabel{setting1}
    \sigma(r)=\sup_{u\in(0,H)}
      \frac{\bigl|(\partial_r \omega)(r,u)\bigr|}
           {1+\bigl|\bigl(\partial_z \omega\bigr)\bigl(r,u\bigr)\bigr|}
    \dott
  \end{equation}
  \startnewargseq
  \argument{
    \lref{setting1};
    the fact that $W$ is definable and open;
    the fact that $\omega\in C^1(W,[0,1])$ is definable and bounded;
    \cref{lem:normalized-interval} (applied with
    $\mathfrak r\curvearrowleft\mathfrak R$, $U\curvearrowleft W$,
    $\Omega\curvearrowleft W$, $Z\curvearrowleft\omega$,
    $\eta\curvearrowleft\sigma$ in the notation of
    \cref{lem:normalized-interval});
  }{
    that there exists $\varrho\in(0,\mathfrak{R})$ which satisfies that
    \begin{enumerate}[label=(\Roman*)]
      \item \llabel{item:bound-finite} it holds that
      $\sigma((0,\varrho))\subseteq\R$ and
      \item \llabel{eq2} it holds that
      $(0,\varrho)\ni r\mapsto\sigma(r)\in\R$
      is continuously differentiable, definable, and integrable\dott
    \end{enumerate}
  }
  In the following for every $\delta\in\{-1,1\}$ let
  $S_\delta\subseteq W$ and for every $r\in(0,\mathfrak{R})$ let
  $\mathcal V_{\delta,r}\subseteq(0,H)$ satisfy
  \begin{equation}
    \llabel{setting2}
    S_\delta=\bigl\{(r,z)\in W\colon
      \delta\bigl(\partial_z \omega\bigr)\bigl(r,z\bigr)>0\bigr\}
    \quad\text{and}\qquad
    \mathcal V_{\delta,r}=\{z\in(0,H)\colon(r,z)\in S_\delta\}
    \dott
  \end{equation}
  \startnewargseq
  \argument{
    \lref{setting2};
    the fact that $\omega\in C^1(W,[0,1])$ is definable;
    \cref{lem:radial-derivative-well-defined} (applied with
    $n\curvearrowleft1$, $\mathfrak{r}\curvearrowleft H$,
    $M\curvearrowleft(0,\mathfrak{R})$,
    $U\curvearrowleft(0,H)\times(0,\mathfrak{R})$,
    $f\curvearrowleft\bigl((0,H)\times(0,\mathfrak{R})\ni(z,r)
      \mapsto\omega(r,z)\in\R\bigr)$ in the notation of
    \cref{lem:radial-derivative-well-defined});
  }{
    that for all $\delta\in\{-1,1\}$ it holds that $S_\delta$ is
    definable and open\dott
    \llabel{eq3}
  }
  \argument{
    \lref{eq3};
    \lref{setting2};
    the fact that for every open and definable $A\subseteq\R$ it holds that
    $A$ is a finite union of pairwise disjoint open intervals;
  }{
    that for all $r\in(0,\mathfrak{R})$, $\delta\in\{-1,1\}$
    there exist $m_{r,\delta}\in\N_0$ and pairwise disjoint
    nonempty open intervals $I_{r,\delta,j}\subseteq(0,H)$,
    $j\in\{1,2,\ldots,m_{r,\delta}\}$, which satisfy
    \begin{equation}
      \llabel{eq4}
      \mathcal V_{\delta,r}=\midcup_{j=1}^{m_{r,\delta}}I_{r,\delta,j}
      \dott
    \end{equation}
  }
  \startnewargseq
  \argument{ \lref{setting2};
    \lref{eq4};
    the fact that for all $\delta\in\{-1,1\}$ it holds that $S_\delta$
    is definable;
    \cite[Theorem~3.12]{Coste};
  }{
    that there exists $N\in\N$ which satisfies that for all
    $r\in(0,\mathfrak{R})$ it holds that
    \begin{equation}
      \llabel{eq5}
      m_{r,-1}+m_{r,1}\leq N
      \dott
    \end{equation}
  }
  \startnewargseq
  \argument{\lref{setting2};
    \lref{eq4};
    the fact that $\omega\in C^1(W,[0,1])$;
    the fundamental theorem of calculus;
  }{
    that for all $r\in(0,\mathfrak{R})$, $\delta\in\{-1,1\}$,
    $j\in\{1,\ldots,m_{r,\delta}\}$,
    $a,b\in I_{r,\delta,j}$ with $a<b$ it holds that
    \begin{equation}
      \llabel{eq6}
      \int_a^b\bigl|\bigl(\partial_z \omega\bigr)\bigl(r,u\bigr)\bigr|\dd u
      =\delta\bigl(\omega(r,b)-\omega(r,a)\bigr)\leq1
      \dott
    \end{equation}
  }
  \argument{
    \lref{eq6};
    the fact that for all $r\in(0,\mathfrak{R})$, $\delta\in\{-1,1\}$, 
    $j\in\{1,\ldots,m_{r,\delta}\}$ it holds that $I_{r,\delta,j}$
    is an open interval;
    the monotone convergence theorem;
  }{
    that for all $r\in(0,\mathfrak{R})$, $\delta\in\{-1,1\}$, 
    $j\in\{1,\ldots,m_{r,\delta}\}$ it holds that
    \begin{equation}
      \llabel{eq7}
      \int_{I_{r,\delta,j}}
        \bigl|\bigl(\partial_z \omega\bigr)\bigl(r,u\bigr)\bigr|\dd u\leq1
      \dott
    \end{equation}
  }
  \argument{
    \lref{eq7};
    \lref{setting2};
    \lref{eq4};
    \lref{eq5};
    the fact that for all $r\in(0,\mathfrak{R})$,
    $u\in(0,H)\backslash(\mathcal V_{-1,r}\cup\mathcal V_{1,r})$ it holds that
    $\bigl(\partial_z \omega\bigr)\bigl(r,u\bigr)=0$
  }{
    that for all $r\in(0,\mathfrak{R})$ it holds that
    \begin{equation}
      \llabel{eq8}
      \int_{(0,H)}
        \bigl|\bigl(\partial_z \omega\bigr)\bigl(r,u\bigr)\bigr|\dd u
      \leq m_{r,-1}+m_{r,1}\leq N
      \dott
    \end{equation}
  }
  \argument{
    \lref{setting1};
    \lref{item:bound-finite};
  }{
    that for all $(r,u)\in(0,\varrho)\times(0,H)$ it holds that
    \begin{equation}
      \llabel{eq9}
      \bigl|(\partial_r \omega)(r,u)\bigr|
      \leq\sigma(r)
        \bigl(1+\bigl|\bigl(\partial_z \omega\bigr)\bigl(r,u\bigr)\bigr|\bigr)
      \dott
    \end{equation}
  }
  \argument{
    \lref{eq9};
    \lref{eq8};
    the fact that $\sigma$ is nonnegative;
  }{
    that for all $r\in(0,\varrho)$ it holds that
    \begin{equation}
      \llabel{eq10}
      \int_{(0,H)}
        \bigl|(\partial_r \omega)(r,u)\bigr|\dd u
      \leq(H+N)\sigma(r)
      \dott
    \end{equation}
  }
  In the following let $b\colon(0,\varrho)\to[0,\infty)$ satisfy 
  for all $r\in(0,\varrho)$ that
  \begin{equation}
    \llabel{setting3}
    b(r)=(H+N)\sigma(r)
    \dott
  \end{equation}
  \startnewargseq
  \argument{
    \lref{eq2};
    \lref{setting3};
    the fact that $H+N>0$;
  }{
    that $b$ is definable and continuous and satisfies
    \begin{equation}
      \llabel{eq11}
      \int_0^{\varrho}b(r)\dd r
      =(H+N)\int_0^{\varrho}\sigma(r)\dd r<\infty
      \dott
    \end{equation}
  }
  \argument{
    the fact that $\omega\in C^1(W,[0,1])$;
    the fundamental theorem of calculus;}
    {that for all $s,t\in(0,\varrho)$,
    $u\in(0,H)$ with $s<t$ it holds that
    \begin{equation}
      \llabel{eq12}
      |\omega(s,u)-\omega(t,u)|
      =\biggl|\int_s^t(\partial_r \omega)(r,u)\dd r\biggr|
      \leq\int_s^t
        \bigl|(\partial_r \omega)(r,u)\bigr|\dd r
      \dott
    \end{equation}
  }
  \argument{
    \lref{eq12};
    \lref{ass1};
    \lref{eq10};
    \lref{setting3};
    Tonelli's theorem;
  }{
    that for all $s,t\in(0,\varrho)$,
    $z\in[0,H]$ with $s<t$ it holds that
    \begin{equation}
      \llabel{eq13}
      \begin{aligned}
        |P(s,z)-P(t,z)|
        &\leq\int_0^z|\omega(s,u)-\omega(t,u)|\dd u\leq\int_0^z\int_s^t
          \bigl|(\partial_r \omega)(r,u)\bigr|\dd r\dd u\\
        &=\int_s^t\int_0^z
          \bigl|(\partial_r \omega)(r,u)\bigr|\dd u\dd r
        \leq\int_s^t b(r)\dd r\dott
      \end{aligned}
    \end{equation}
  }
  \argument{
    \lref{eq13};
    \lref{eq11};
  }[verbs=ep]{
    \lref{item5}\dott
  }
  \startnewargseq
  \argument{
    \lref{eq2};
    \lref{setting3};
    the extreme value theorem
  }{
    that for all $\rho_1,\rho_2\in(0,\varrho)$ with $\rho_1<\rho_2$
    there exists $K_{\rho_1,\rho_2}\in[0,\infty)$
    which satisfies for all $r\in[\rho_1,\rho_2]$ that
    \begin{equation}
      \llabel{eq14}
      b(r)\leq K_{\rho_1,\rho_2}
      \dott
    \end{equation}
  }
  \argument{
    \lref{eq14};
    \lref{item2};
    \lref{item5};
    the triangle inequality;
    the Cauchy--Schwarz inequality ;
  }{
    that for all $\rho_1,\rho_2\in(0,\varrho)$,
    $(\mathfrak s_1,\mathfrak z_1),(\mathfrak s_2,\mathfrak z_2)\in[\rho_1,\rho_2]\times[0,H]$
    with $\rho_1<\rho_2$ it holds that
    \begin{equation}
      \llabel{eq15}
      \begin{aligned}
        |P(\mathfrak s_1,\mathfrak z_1)-P(\mathfrak s_2,\mathfrak z_2)|
        &\leq|P(\mathfrak s_1,\mathfrak z_2)-P(\mathfrak s_2,\mathfrak z_2)|
          +|P(\mathfrak s_1,\mathfrak z_1)-P(\mathfrak s_1,\mathfrak z_2)|\\
        &\leq K_{\rho_1,\rho_2}|\mathfrak s_1-\mathfrak s_2|+|\mathfrak z_1-\mathfrak z_2|\\
        &\leq\sqrt{(K_{\rho_1,\rho_2})^2+1}
          \sqrt{(\mathfrak s_1-\mathfrak s_2)^2+(\mathfrak z_1-\mathfrak z_2)^2}
          \dott
      \end{aligned}
    \end{equation}
  }
  \argument{
    \lref{eq15};
  }[verbs=ep]{
    \lref{item6}\dott
  }
\end{aproof}

\subsection{Lyapunov-type stability for the radial-angular dynamics}
\label{subsec:bounded-state}

In \eqref{prop:buildQ.concl2} in
\cref{prop:buildQ} below we establish a crucial Lyapunov-type
estimate for the radial-angular \DE\ system.
The function $Q\colon\mathcal O\cap((0,\eta)\times M)\to[0,\infty)$ in
\cref{prop:buildQ} serves as the Lyapunov-type function for
the radial-angular \DE\ system in \eqref{prop:buildQ.ass2}.
In our proof of \cref{prop:buildQ} we employ the elementary
estimates for real numbers in the following \cref{lem:attenuation}.

\begin{athm}{lemma}{lem:attenuation}
  Let $r,T\in(0,\infty)$, $g\in[0,\infty)$, $n\in[1,\infty)$,
  $v,c,D,\mathcal B\in\R$ satisfy
  \begin{equation}
    \llabel{ass1}
    |v-c|\leq n(T+g),\quad r^2n^2\leq\tfrac{1}{64},\quad
    D=T+8r^2n(c^2+T^2)^{1/2},\quad\text{and}\quad
    \mathcal B=-D^{-1}(g^2+r^2v(v-c))
    \dott
  \end{equation}
  \cfload[.]
  Then the following two statements hold:
  \begin{enumerate}[label=(\roman*)]
    \item \llabel{item1} It holds that
    \begin{eqn}[i]
      \mathcal B\leq2nr^2|v|
      \dott
    \end{eqn}
    \item \llabel{item2} If $g\geq T$, then
    \begin{eqn}[i]
      \mathcal B\leq5nr^2|v|-\frac{g}{4}
      \dott
    \end{eqn}
  \end{enumerate}
  \cfout[.]
\end{athm}
\begin{aproof}
  Throughout this proof, let $Y,P\in\R$ satisfy
  \begin{equation}
    \llabel{setting1}
    Y=r^2n|c|\qquad\text{and}\qquad P=g^2+r^2v(v-c)
    \dott
  \end{equation}
  \startnewargseq
  \argument{
    \lref{ass1};
    the fact that $r>0$ and $n\geq1$;
  }{
    that
    \begin{equation}
      \llabel{eq1}
      0<r^2n\leq r^2n^2\leq\frac1{64}
      \dott
    \end{equation}
  }
  \argument{
    \lref{eq1};
    \lref{ass1};
    \lref{setting1};
    the fact that $T>0$;
    the fact that $|c|\leq\sqrt{c^2+T^2}\leq|c|+T$;
  }{
    that
    \begin{equation}
      \llabel{eq3}
      D\geq T>0,\qquad D\geq8Y,
    \end{equation}
    and
    \begin{equation}
      \llabel{eq3.1}
      D\leq T+8r^2n(|c|+T)
      =(1+8r^2n)T+8Y\leq\frac98T+8Y
      \dott
    \end{equation}
  }
  \argument{
    \lref{eq3.1};
    \lref{ass1};
    the fact that $g^2\geq0$ and $r^2>0$;
    the fact that $-v(v-c)\leq|v|\,|v-c|$;
  }{
    that
    \begin{equation}
      \llabel{eq4}
      \mathcal B=-\frac{g^2}{D}-\frac{r^2v(v-c)}{D}
      \leq\frac{r^2|v|\,|v-c|}{D}
      \dott
    \end{equation}
  }
  \argument{
    \lref{ass1};
    the fact that $n\geq1$;
  }{
    that if $g<T$, then it holds that
    \begin{equation}
      \llabel{eq5}
      |v-c|\leq n(T+g)\leq2nT
      \dott
    \end{equation}
  }
  \argument{
    \lref{eq5};
    \lref{eq3};
    \lref{eq4};
    the fact that $nr^2|v|\geq0$;
  }{
    that if $g<T$, then it holds that
    \begin{equation}
      \llabel{eq6}
      \mathcal B\leq\frac{r^2|v|\,|v-c|}{D}
      \leq2nr^2|v|\frac{T}{D}\leq2nr^2|v|
      \dott
    \end{equation}
  }
  \argument{
    \lref{ass1};
    the fact that $T>0$ and $n\geq1$;
  }{
    that if $g\geq T$, then it holds that
    \begin{equation}
      \llabel{eq7}
      g\geq T>0\qquad\text{and}\qquad |v-c|\leq n(T+g)\leq2ng
      \dott
    \end{equation}
  }
  \argument{
    \lref{eq7};
    \lref{setting1};
    the fact that $c(v-c)\geq-|c|\,|v-c|$;
  }{
    that if $g\geq T$, then it holds that
    \begin{equation}
      \llabel{eq8}
        P=g^2+r^2(v-c)^2+r^2c(v-c)
         \geq g^2-r^2|c|\,|v-c|\geq g^2-2Yg
      \dott
    \end{equation}
  }
  \argument{
    \lref{eq8};
    \lref{eq3};
    \lref{eq3.1};
    the fact that $T>0$;
  }{
    that if $g\geq T$ and $g\geq8Y$, then it holds that
    \begin{equation}
      \llabel{eq9}
      P\geq g^2-2Yg\geq\frac34g^2>0\qquad\text{and}\qquad
      0<D\leq\frac98T+8Y\leq\frac{17}{8}g
      \dott
    \end{equation}
  }
  \argument{
    \lref{eq9};
    \lref{ass1};
    \lref{setting1};
  }{
    that if $g\geq T$ and $g\geq8Y$, then it holds that
    \begin{equation}
      \llabel{eq10}
      \mathcal B=-\frac PD
      \leq-\frac{3g^2}{4D}
      \leq-\frac6{17}g\leq-\frac g4\leq0
      \dott
    \end{equation}
  }
  \argument{
    \lref{eq3};
    \lref{eq4};
    \lref{eq7};
  }{
    that if $g\geq T$ and $g<8Y$, then it holds that
    \begin{equation}
      \llabel{eq11}
      D\geq8Y>g\geq T>0\qquad\text{and}\qquad
      \mathcal B\leq\frac{r^2|v|\,|v-c|}{D}
      \leq2nr^2|v|\frac gD\leq2nr^2|v|
      \dott
    \end{equation}
  }
  \argument{
    \lref{eq11};
    \lref{eq6};
    \lref{eq10};
  }[verbs=ep]{
    \lref{item1}\dott
  }
  \startnewargseq
  \argument{
    \lref{ass1};
    \lref{setting1};
    \lref{eq7};
    the fact that $n\geq1$;
  }{
    that if $g\geq T$ and $g<8Y$, then it holds that
    \begin{equation}
      \llabel{eq12}
      |v-c|\leq2ng<16nY
      =16r^2n^2|c|\leq\frac{|c|}{4}
      \dott
    \end{equation}
  }
  \argument{
    \lref{eq12};
    the triangle inequality
  }{
    that if $g\geq T$ and $g<8Y$, then it holds that
    \begin{equation}
      \llabel{eq13}
      |v|\geq|c|-|v-c|\geq\frac34|c|
      \dott
    \end{equation}
  }
  \argument{
    \lref{eq13};
    \lref{setting1};
    the fact that $nr^2>0$;
  }{
    that if $g\geq T$ and $g<8Y$, then it holds that
    \begin{equation}
      \llabel{eq14}
      \frac g4<2Y=2nr^2|c|\leq\frac83nr^2|v|
      \dott
    \end{equation}
  }
  \argument{
    \lref{eq14};
    \lref{item1};
  }{
    that if $g\geq T$ and $g<8Y$, then it holds that
    \begin{equation}
      \llabel{eq15}
        \mathcal B\leq2nr^2|v|
          =-\frac g4+2nr^2|v|+\frac g4\leq-\frac g4+\frac{14}{3}nr^2|v|
          \leq-\frac g4+5nr^2|v|
      \dott
    \end{equation}
  }
  \argument{
    \lref{eq15};
    \lref{eq10};
    the fact that $5nr^2|v|\geq0$;
  }[verbs=ep]{
    \lref{item2}\dott
  }
\end{aproof}

\begin{athm}{proposition}{prop:buildQ}
  [Lyapunov-type estimate]
  Let $D\in\N$, $\mathfrak R,\mathsf a\in(0,\infty)$,
  $\mathsf b\in[\mathsf a,\infty)$, let $M\subseteq\R^D$ be a
  nonempty, compact, and definable\cfadd{definition:definable-set}
  \submanifold{2}{D}, let
  $\mathcal O\subseteq(0,\mathfrak R)\times M$ be definable and open
  in $(0,\mathfrak R)\times M$, and let
  $h=(h(r,x))_{(r,x)\in\mathcal O}\in \COne{\mathcal O}{\R}$
  be definable\cfadd{definition:c1-on-manifold}%
  \cfadd{definition:definable-function}\dott \cfload[.]
  Then there exist $\eta\in(0,\mathfrak R)$, a bounded
  $Q\in C(\mathcal O\cap((0,\eta)\times M),[0,\infty))$, and a definable and integrable
  $m\in C((0,\eta),[0,\infty))$\cfadd{definition:definable-function}
  such that for every interval $I\subseteq\R$,
  all locally absolutely continuous $r\colon I\to(0,\eta)$ and
  $x\colon I\to M$, and
  all Lebesgue measurable $\alpha\colon I\to[0,\infty)$ and
  $A\colon I\to\R^{D\times D}$
  with $\midcup_{t\in I}\{(r(t),x(t))\}\subseteq\mathcal O$ and
  with the property that for Lebesgue-almost every $t\in I$ and
  every $v\in T_{x(t)}M$\cfadd{definition:tangent-space} it holds that
  \begin{eqn}
    \llabel{ass1}
    A(t)v\in T_{x(t)}M,
    \qquad \langle A(t)v,v\rangle\geq\mathsf a\norm{v}^2,
    \qquad
    \norm{A(t)v}\leq\mathsf b\norm{v},
  \end{eqn}
  \begin{eqn}
    \llabel{ass2}
    r'(t)=-\alpha(t)(r(t))^2
      (\partial_r h)(r(t),x(t)),
    \qquad\text{and}\qquad
    x'(t)=-\alpha(t)A(t)(\nabla_Mh)(r(t),x(t))
  \end{eqn}
  \cfadd{definition:angular-gradient}%
  we have that for all $a,b\in I$ with $a\leq b$ it holds that
  \begin{eqn}
    \llabel{concl2}
    Q(r(b),x(b))-Q(r(a),x(a))
      \leq-\int_a^b\norm{x'(t)}\dd t
        +\int_a^b m(r(t))\abs{r'(t)}\dd t
    \dott
  \end{eqn}
  \cfout[.]
\end{athm}
\begin{aproof}
  Throughout this proof let $\varphi\colon\R\to(-1,1)$ and
  $\mathbf{h}\colon\mathcal O\to(-1,1)$ satisfy for all
  $e\in\R$ that $\varphi(e)=e/\sqrt{1+e^2}$ and 
  $\mathbf{h}=\varphi\circ h$\dott \cfload[.]
  \startnewargseq
  \argument{
    the fact that $h\in \COne{\mathcal O}{\R}$ is definable;
    the chain rule;
  }{
    that $\varphi\in C^1(\R,(-1,1))$ and
    $\mathbf{h}\in \COne{\mathcal O}{(-1,1)}$ are definable,
    that for all $e\in\R$ it holds that
    \begin{eqn}
      \llabel{eq0.1}
      \varphi'(e)=(1+e^2)^{-3/2}>0,
    \end{eqn}
    and that for all $(s,y)\in\mathcal O$ it holds that
    \begin{eqn}[a]
      \llabel{eq0.2}
      (\partial_r \mathbf{h})(s,y)
        =\varphi'(h(s,y))(\partial_r h)(s,y)\qquad\text{and}\qquad
      (\nabla_M\mathbf{h})(s,y)
        =\varphi'(h(s,y))(\nabla_Mh)(s,y)
      \dott
    \end{eqn}
  }
  \argument{
    the fact that for all $e\in\R$ it holds that
    $\varphi(e)=e/\sqrt{1+e^2}$;
    the fact that for all $(s,y)\in\mathcal O$ it holds that
    $\mathbf{h}(s,y)=\varphi(h(s,y))$;
  }{
    that for all $(s,y)\in\mathcal O$ it holds that
    \begin{eqn}
      \llabel{setting0}
      -1<\mathbf{h}(s,y)<1
      \dott
    \end{eqn}
  }
  In the following let $\mathsf c,\Lambda,\kappa\in(0,\infty)$
  satisfy
  \begin{eqn}
    \llabel{setting1}
    \mathsf c=\sqrt{\mathsf a},\qquad
    \Lambda=\frac{\mathsf b}{\mathsf c},\qquad
    \kappa=4\Lambda,
  \end{eqn}
  and let $G\colon\mathcal O\to[0,\infty)$ satisfy for all
  $(s,y)\in\mathcal O$ that
  \begin{eqn}
    \llabel{setting1.1}
    G(s,y)=\norm{(\nabla_M\mathbf{h})(s,y)}
    \dott
  \end{eqn}
  \cfadd{definition:angular-gradient}\cfload[.]
  \startnewargseq
  \argument{
    \lref{setting0};
    \lref{setting1};
    \lref{setting1.1};
    \cref{lem:divide}
    (applied with $R\curvearrowleft\mathfrak R$,
    $\Omega\curvearrowleft\mathcal O$,
    $h\curvearrowleft\mathbf{h}$,
    $a\curvearrowleft((0,\mathfrak R)\ni s\mapsto-1\in\R)$,
    $b\curvearrowleft((0,\mathfrak R)\ni s\mapsto1\in\R)$,
    $c\curvearrowleft\mathsf c$ in the notation of \cref{lem:divide});
  }{
    that there exist $\theta\in(0,\mathfrak R)$,
    $N\in\N$, $H\in(0,\infty)$ and definable
    $\mathfrak e_0,\mathfrak e_1,\dots,\mathfrak e_N\in C^2((0,\theta),\R)$,
    $\mathcal F_0,\mathcal F_1,\dots,\mathcal F_N
      \in C((0,\theta)\times\R,[0,H])$,
    $E_1,E_2,\dots,E_N\in C^2((0,\theta)\times(0,H),\R)$,
    $T_1,T_2,\dots,T_N
      \in C^1((0,\theta)\times(0,H),(0,\infty))$
    such that
    \begin{enumerate}[label=(\Roman*)]
      \item for all $s\in(0,\theta)$ it holds that
      \begin{eqn}
        \llabel{eq1}
        -1=\mathfrak e_0(s)<\cdots<\mathfrak e_N(s)=1,
      \end{eqn}
      \item for all $s\in(0,\theta)$, $e\in\R$ with
      $\inf\bigl(\{\infty\}\cup
        \{z\in\R\colon\exists\,y\in M\colon
          (s,y)\in\mathcal O\cap\mathbf h^{-1}(\{e\})\cap G^{-1}(\{z\})\}\bigr)=0$
      it holds that $e\in\{\mathfrak e_0(s),\mathfrak e_1(s),\dots,\mathfrak e_N(s)\}$,
      \item for all $i\in\{0,1,\dots,N-1\}$, $s\in(0,\theta)$,
      $u\in(-\infty,\mathfrak e_i(s)]$,
      $v\in[\mathfrak e_{i+1}(s),\infty)$,
      $z,w\in(\mathfrak e_i(s),\mathfrak e_{i+1}(s))$ with $z<w$
      it holds that
      \begin{eqn}
        \llabel{eq1.2}
        0=\mathcal F_i(s,u)<\mathcal F_i(s,z)
          <\mathcal F_i(s,w)<\mathcal F_i(s,v)=H,
      \end{eqn}
      \item for all $i\in\{0,1,\dots,N-1\}$ it holds that
      $\mathcal F_i|_{\{(s,u)\in(0,\theta)\times\R\colon
        \mathfrak e_i(s)<u<\mathfrak e_{i+1}(s)\}}$
      is twice continuously differentiable,
      \item for all $i\in\{1,2,\dots,N\}$, $s\in(0,\theta)$,
      $z\in(0,H)$ it holds that
      \begin{eqn}
        \llabel{eq1.5}
        \mathfrak e_{i-1}(s)<E_i(s,z)<\mathfrak e_i(s),\qquad
        \mathcal F_{i-1}(s,E_i(s,z))=z,\qquad\text{and}\qquad
        T_i(s,z)=\partial_z E_i(s,z)>0,
      \end{eqn}
      and
      \item for all $i\in\{1,2,\dots,N\}$,
      $(s,y)\in\mathcal O\cap((0,\theta)\times M)$ with
      $\mathfrak e_{i-1}(s)<\mathbf h(s,y)<\mathfrak e_i(s)$ it holds that
      \begin{eqn}
        \llabel{eq1.6}
        T_i(s,\mathcal F_{i-1}(s,\mathbf h(s,y)))\leq\mathsf c G(s,y)
        \dott
      \end{eqn}
    \end{enumerate}
  }
  In the following for every $i\in\{0,1,\dots,N-1\}$ let
  $Z_i\colon\mathcal O\cap((0,\theta)\times M)\to[0,H]$,
  $S_i\colon(0,\theta)\times(0,H)\to\R$, and
  $\Omega_i\subseteq\mathcal O\cap((0,\theta)\times M)$ satisfy that
  for all $(s,y)\in\mathcal O\cap((0,\theta)\times M)$ it holds that
  \begin{eqn}
    \llabel{setting2}
    Z_i(s,y)=\mathcal F_i(s,\mathbf{h}(s,y)),
  \end{eqn}
  that for all $(s,z)\in(0,\theta)\times(0,H)$ it holds that
  \begin{eqn}
    \llabel{setting2.1}
    S_i(s,z)=(\partial_s E_{i+1})(s,z),
  \end{eqn}
  and that
  \begin{eqn}
    \llabel{setting3}
    \Omega_i=\{(s,y)\in\mathcal O\cap((0,\theta)\times M)\colon
      \mathfrak e_i(s)<\mathbf{h}(s,y)<\mathfrak e_{i+1}(s)\}
    \dott
  \end{eqn}
  \cfload[.]
  \startnewargseq
  \argument{
    \lref{setting2};
    \lref{setting2.1};
    the fact that for all $i\in\{0,1,\dots,N-1\}$ it holds that $\mathbf{h}$ and
    $\mathcal F_i$ are continuous and definable;
    the fact that $(\mathcal O\cap((0,\theta)\times M))\ni(s,y)\mapsto s\in(0,\theta)$
    is continuous and definable;
  }{
    that for all $i\in\{0,1,\dots,N-1\}$ it holds that $Z_i$ is continuous
    and definable\dott
  }
  \argument{
    \lref{eq1.2};
    \lref{setting2};
    \lref{setting2.1};
    \lref{setting3};
  }{
    that for all $i\in\{0,1,\dots,N-1\}$ it holds that
    \begin{eqn}
      \llabel{eq2}
      \Omega_i=\{(s,y)\in\mathcal O\cap((0,\theta)\times M)\colon0<Z_i(s,y)<H\}
      \dott
    \end{eqn}
  }
  \argument{
    \lref{eq2};
    the fact that for all $i\in\{0,1,\dots,N-1\}$ it holds that $Z_i$ is
    continuous and definable;
    the fact that $M$ is a \submanifold{2}{D};
    the fact that $\mathcal O$ is definable and open in
    $(0,\mathfrak R)\times M$;
    \cref{lem:product};
  }{
    that for all $i\in\{0,1,\dots,N-1\}$ it holds that $\Omega_i$ is
    definable, open in $(0,\theta)\times M$, and a
    \submanifold{1}{D+1}\dott
  }
  We claim that for all $i\in\{0,1,\dots,N-1\}$ it holds that
  \begin{equation}
    \llabel{eq2.0}
    Z_i|_{\Omega_i}\in\COne{\Omega_i}{[0,H]}.
  \end{equation}
  \startnewargseq
  To prove this claim, let $i\in\{0,1,\dots,N-1\}$, let
  $n\in\{0,1,\dots,D\}$ satisfy that $M$ is an
  \dsubmanifold{2}{D}{n}, let $W\subseteq\R^{n+1}$
  be open, and let
  $\phi=(\phi_0,\ldots,\phi_D)\in C^1(W,\Omega_i)$
  satisfy for all $w\in W$ that $\operatorname{rank}\phi'(w)=n+1$\dott\cfload[.]
  \startnewargseq
  \argument{
    \lref{setting3};
    the fact that $M$ is an \dsubmanifold{2}{D}{n};
    the fact that $\phi\in C^1(W,\Omega_i)$;
    the fact that $\Omega_i\subseteq\mathcal O$;
    the fact that for all $w\in W$ it holds that
    $\operatorname{rank}\phi'(w)=n+1$;
    the fact that $\mathbf{h}\in \COne{\mathcal O}{(-1,1)}$;
  }{
    that $\phi_0\in C^1(W,\R)$ and
    $\mathbf{h}\circ\phi\in C^1(W,(-1,1))$ and that for all $w\in W$
    it holds that
    \begin{eqn}
      \llabel{eq2.1}
      \mathfrak e_i(\phi_0(w))<\mathbf{h}(\phi(w))
        <\mathfrak e_{i+1}(\phi_0(w))
      \dott
    \end{eqn}
  }
  \argument{
    \lref{eq2.1};
    \lref{setting2};
    \lref{setting2.1};
    the fact that $\phi_0\in C^1(W,\R)$ and
    $\mathbf{h}\circ\phi\in C^1(W,(-1,1))$;
    the fact that $\mathcal F_i|_{\{(s,u)\in(0,\theta)\times\R\colon
      \mathfrak e_i(s)<u<\mathfrak e_{i+1}(s)\}}
      \in C^2(\{(s,u)\in(0,\theta)\times\R\colon
        \mathfrak e_i(s)<u<\mathfrak e_{i+1}(s)\},[0,H])$;
    the chain rule
  }{
    that $Z_i\circ\phi\colon W\to[0,H]$ is differentiable and that
    for all $w\in W$ it holds that
    \begin{eqn}[a]
      \llabel{eq2.2}
      (Z_i\circ\phi)'(w)
        &=\bigl(\partial_s \mathcal F_i\bigr)
          \bigl(\phi_0(w),\mathbf{h}(\phi(w))\bigr)\,(\phi_0)'(w)+\bigl(\partial_u \mathcal F_i\bigr)
          \bigl(\phi_0(w),\mathbf{h}(\phi(w))\bigr)\,(\mathbf{h}\circ\phi)'(w)
      \dott
    \end{eqn}
  }
  \argument{
    \lref{eq2.2};
    \lref{eq2.1};
    the fact that $\phi_0$, $\mathbf{h}\circ\phi$, $(\phi_0)'$, and
    $(\mathbf{h}\circ\phi)'$ are continuous on $W$;
    the fact that $\partial_s\mathcal F_i$ and
    $\partial_u\mathcal F_i$ are continuous on
    $\{(s,u)\in(0,\theta)\times\R\colon
      \mathfrak e_i(s)<u<\mathfrak e_{i+1}(s)\}$;
  }{
    that $(Z_i\circ\phi)'$ is continuous on $W$\dott
    \llabel{eq2.2.1}
  }
  \argument{
    \lref{eq2.2.1};
    the fact that $Z_i\circ\phi$ is differentiable on $W$;
  }{
    that $Z_i\circ\phi\in C^1(W,[0,H])$\dott
    \llabel{eq2.3}
  }
  \argument{
    \lref{eq2.3};
    the arbitrariness of $i$, $W$, and $\phi$;
  }{
    \lref{eq2.0}\dott
  }

  \startnewargseq
  \argument{
    \lref{eq1.2};
    \lref{setting2};
    \lref{setting2.1};
    \lref{setting3};
  }{
    that for all $i\in\{0,1,\dots,N-1\}$, $(s,y)\in\mathcal O\cap((0,\theta)\times M)$
    it holds that
    \begin{eqn}
      \llabel{eq3}
      Z_i(s,y)=
      \begin{cases}
        0&\colon \mathbf{h}(s,y)\leq\mathfrak e_i(s),\\
        H&\colon \mathbf{h}(s,y)\geq\mathfrak e_{i+1}(s),\\
        \mathcal F_i(s,\mathbf{h}(s,y))&\colon(s,y)\in\Omega_i.
      \end{cases}\relax
    \end{eqn}
  }
  \argument{
    \lref{eq1.2};
    \lref{eq1.5};
    \lref{eq1.6};
    \lref{setting2};
    \lref{setting2.1};
    \lref{setting3};
    \lref{eq2};
  }{
    that for all $i\in\{0,1,\dots,N-1\}$, $(s,y)\in\Omega_i$,
    $z\in(0,H)$ it holds that
    \begin{eqn}
      \llabel{eq4}
      \mathbf{h}(s,y)=E_{i+1}(s,Z_i(s,y)),\qquad
      0<T_{i+1}(s,Z_i(s,y))\leq\mathsf c G(s,y),
      \qquad\text{and}\qquad
      T_{i+1}(s,z)=\bigl(\partial_z E_{i+1}\bigr)\bigl(s,z\bigr)
      \dott
    \end{eqn}
  }
  \argument{
    \lref{setting1};
    \lref{setting1.1};
    \lref{eq1};
    \lref{eq1.2};
    \lref{eq1.5};
    \lref{setting2};
    \lref{setting2.1};
    \lref{setting3};
    the fact that for all $i\in\{0,1,\dots,N-1\}$ it holds that
    $\mathcal F_i|_{\{(s,u)\in(0,\theta)\times\R\colon
      \mathfrak e_i(s)<u<\mathfrak e_{i+1}(s)\}}
      \in C^2(\{(s,u)\in(0,\theta)\times\R\colon
        \mathfrak e_i(s)<u<\mathfrak e_{i+1}(s)\},[0,H])$;
    \cref{lem:residual} (applied with
    $\mathfrak r_0\curvearrowleft\mathfrak R$,
    $\mathfrak r_1\curvearrowleft\theta$,
    $h\curvearrowleft\mathbf h$,
    $\mathbf c\curvearrowleft\min\{\mathsf c,1/64\}$,
    $(\mathcal F_1,\dots,\mathcal F_N)\curvearrowleft
      (\mathcal F_0,\mathcal F_1,\dots,\mathcal F_{N-1})$,
    $(\mathbf m_1,\dots,\mathbf m_N)\curvearrowleft
      (\mathbf m_0,\mathbf m_1,\dots,\mathbf m_{N-1})$
    in the notation of \cref{lem:residual});
  }{
    that there exist $\tau\in(0,\theta)$ and definable
    $\mathbf m_0,\mathbf m_1,\dots,\mathbf m_{N-1}
      \in C^1((0,\tau),[1,\infty))$, which satisfy that
    for all $i\in\{0,1,\dots,N-1\}$ it holds that
    \begin{eqn}
      \llabel{eq5}
      \int_0^\tau\mathbf m_i(s)\dd s<\infty,
    \end{eqn}
    that for all $i\in\{0,1,\dots,N-1\}$,
    $(s,y)\in\Omega_i\cap((0,\tau)\times M)$ it holds that
    \begin{eqn}
      \llabel{eq6}
      \abs[\big]{(\partial_r \mathbf{h})(s,y)-S_i(s,Z_i(s,y))}
        \leq\mathbf m_i(s)
          \bigl(T_{i+1}(s,Z_i(s,y))+\mathsf c G(s,y)\bigr),
    \end{eqn}
    and that for all $i\in\{0,1,\dots,N-1\}$, $s\in(0,\tau)$ it holds that
    \begin{eqn}
      \llabel{eq7}
      s^2(\mathbf m_i(s))^2\leq\min\{\mathsf c,1/64\}\leq \frac1{64}
      \dott
    \end{eqn}
  }
  In the following for all $i\in\{0,1,\dots,N-1\}$ let
  $\omega_i\colon(0,\tau)\times(0,H)\to(0,1]$ and
  $P_i\colon(0,\tau)\times[0,H]\to\R$ satisfy that for all
  $s\in(0,\tau)$, $z\in(0,H)$, $v\in[0,H]$ it holds that
  \begin{eqn}
    \llabel{setting4}
    \omega_i(s,z)=
      \frac{T_{i+1}(s,z)}
        {T_{i+1}(s,z)+8s^2\mathbf m_i(s)
          \sqrt{(S_i(s,z))^2+(T_{i+1}(s,z))^2}}
    \qquad\text{and}\qquad
    P_i(s,v)=\int_0^v\omega_i(s,u)\dd u
    \dott
  \end{eqn}
  \cfload[.]
  \startnewargseq
  \argument{
    \lref{setting2};
    \lref{setting2.1};
    \lref{setting4};
    the fact that for all $i\in\{0,1,\dots,N-1\}$ it holds that
    $E_{i+1}\in C^2((0,\theta)\times(0,H),\R)$ is definable;
    the fact that for all $i\in\{0,1,\dots,N-1\}$ it holds that
    $T_{i+1}\in C^1((0,\theta)\times(0,H),(0,\infty))$ is definable and positive;
    the fact that for all $i\in\{0,1,\dots,N-1\}$ it holds that
    $\mathbf m_i\in C^1((0,\tau),[1,\infty))$ is definable;
    \cref{lem:radial-derivative-well-defined};
  }{
    that for all $i\in\{0,1,\dots,N-1\}$ it holds that
    $\omega_i\in C^1((0,\tau)\times(0,H),(0,1])$
    is definable and that for all $i\in\{0,1,\dots,N-1\}$,
    $(s,z)\in(0,\tau)\times(0,H)$ it holds that
    \begin{eqn}
      \llabel{eq8}
      0<\omega_i(s,z)\leq1
      \dott
    \end{eqn}
  }
  \argument{
    \lref{eq8};
    \lref{setting1};
    \lref{setting1.1};
    \lref{setting4};
    \itref{lem:primitive}{item1,item2,item5};
    \cref{lem:contact-bound}
    (applied for every $j\in\{0,1,\dots,N\}$ with
    $n\curvearrowleft D$, $\mathfrak{r}_0\curvearrowleft\tau$,
    $\mathcal O\curvearrowleft\mathcal O\cap((0,\tau)\times M)$,
    $h\curvearrowleft\mathbf{h}|_{\mathcal O\cap((0,\tau)\times M)}$,
    $a\curvearrowleft\mathfrak e_j|_{(0,\tau)}$ in the notation of \cref{lem:contact-bound});
    the fact that for all $j\in\{0,1,\dots,N\}$ it holds that
    $\mathfrak e_j\in C^2((0,\theta),\R)$ is definable;
  }{
    that there exist $\eta\in(0,\tau)$, definable
    $\mathfrak b_i\in C((0,\eta),[0,\infty))$, $i\in\{0,1,\dots,N-1\}$,
    and $\beta_j\colon(0,\eta)\to[0,\infty)$, $j\in\{0,1,\dots,N\}$,
    which satisfy that for all $j\in\{0,1,\dots,N\}$ it holds that
    $(0,\eta)\ni s\mapsto\beta_j(s)\in\R$
    is definable, continuously differentiable, and integrable, that
    \begin{eqn}
      \llabel{eq9}
      \sum_{i=0}^{N-1}\int_0^\eta\mathfrak b_i(s)\dd s
        +\sum_{j=0}^{N}\int_0^\eta\beta_j(s)\dd s<\infty,
    \end{eqn}
    that for all $i\in\{0,1,\dots,N-1\}$, $s,t\in(0,\eta)$,
    $z,v\in[0,H]$ it holds that
    \begin{eqn}
      \llabel{eq10}
      0\leq P_i(s,z)\leq H,\quad
      \abs{P_i(s,z)-P_i(s,v)}\leq\abs{z-v},\quad\text{and}\quad
      \abs{P_i(s,z)-P_i(t,z)}
        \leq\biggl|\int_s^t\mathfrak b_i(u)\dd u\biggr|,
    \end{eqn}
    and that for all $j\in\{0,1,\dots,N\}$, $(s,y)\in\mathcal O\cap((0,\eta)\times M)$
    with $\mathbf{h}(s,y)=\mathfrak e_j(s)$ and $G(s,y)>0$ it holds that
    \begin{eqn}
      \llabel{eq11}
      \abs[\big]{(\partial_r \mathbf{h})(s,y)-(\mathfrak e_j)'(s)}
        \leq\beta_j(s)G(s,y)
      \dott
    \end{eqn}
  }
  In the following let
  $\mathbb{M},m\colon(0,\eta)\to[0,\infty)$ satisfy that for all
  $s\in(0,\eta)$ it holds that
  \begin{eqn}
    \llabel{setting5.1}
    \mathbb{M}(s)=\kappa\sum_{i=0}^{N-1}
      \bigl(5\mathbf m_i(s)+\mathfrak b_i(s)\bigr)
    \qquad\text{and}\qquad
    m(s)=\mathbb{M}(s)+\frac{\Lambda}{\mathsf c}
      \sum_{j=0}^{N}\beta_j(s)
  \end{eqn}
  and let $Q\colon\mathcal O\cap((0,\eta)\times M)\to[0,\infty)$
  satisfy that for all $(s,y)\in\mathcal O\cap((0,\eta)\times M)$
  it holds that
  \begin{eqn}
    \llabel{setting5}
    Q(s,y)=\kappa\sum_{i=0}^{N-1}P_i(s,Z_i(s,y))
    \dott
  \end{eqn}
  \cfload[.]
  \startnewargseq
  \argument{
    \lref{eq10};
    the fact that for all $i\in\{0,1,\dots,N-1\}$ it holds that
    $\mathfrak b_i$ is nonnegative;
    the triangle inequality;
  }{
    that for all $i\in\{0,1,\dots,N-1\}$, $s,t\in(0,\eta)$,
    $z,v\in[0,H]$ it holds that
    \begin{eqn}[a]
      \llabel{eq11.1}
      \abs{P_i(t,v)-P_i(s,z)}
        \leq\abs{P_i(t,v)-P_i(t,z)}
          +\abs{P_i(t,z)-P_i(s,z)}\leq\abs{v-z}
          +\int_{\min\{s,t\}}^{\max\{s,t\}}\mathfrak b_i(u)\dd u
      \dott
    \end{eqn}
  }
  \argument{
    \lref{eq11.1};
    the fact that for all $i\in\{0,1,\dots,N-1\}$ it holds that
    $\mathfrak b_i$ is continuous on $(0,\eta)$;
    the fundamental theorem of calculus;
  }{
    that for all $i\in\{0,1,\dots,N-1\}$ it holds that
    $P_i|_{(0,\eta)\times[0,H]}$ is continuous\dott
    \llabel{arg11.1}
  }
  \argument{
    \lref{arg11.1};
    \lref{eq10};
    \lref{setting5};
    the fact that $\kappa>0$;
    the fact that for all $i\in\{0,1,\dots,N-1\}$ it holds that
    $Z_i$ is continuous;
  }{
    that $Q$ is continuous and satisfies that for all
    $(s,y)\in\mathcal O\cap((0,\eta)\times M)$ it holds that
    \begin{eqn}
      \llabel{eq12}
      0\leq Q(s,y)\leq\kappa NH
      \dott
    \end{eqn}
  }
  \argument{
    \lref{setting5.1};
    \lref{setting5};
    the fact that for all $i\in\{0,1,\dots,N-1\}$ it holds that
    $\mathbf m_i$ and $\mathfrak b_i$ are definable, continuous, and
    nonnegative;
    the fact that for all $j\in\{0,1,\dots,N\}$ it holds that
    $\beta_j$ is definable, continuous, and nonnegative;
  }{
    that $m$ is definable, continuous, nonnegative, and
    Lebesgue measurable\dott
    \llabel{arg1}
  }
  \argument{
    \lref{arg1};
    \lref{eq5};
    \lref{eq9};
    \lref{setting5.1};
    \lref{setting5};
  }[verbs=ep]{
    that
    \begin{eqn}[i]
      \int_0^\eta m(s)\dd s<\infty
      \dott
    \end{eqn}
  }
  We claim that for every interval $I\subseteq\R$,
  all locally absolutely continuous $r\colon I\to(0,\eta)$ and
  $x\colon I\to M$, and
  all Lebesgue measurable $\alpha\colon I\to[0,\infty)$ and
  $A\colon I\to\R^{D\times D}$
  with $\midcup_{t\in I}\{(r(t),x(t))\}\subseteq\mathcal O$ and
  with the property that for Lebesgue-almost all $t\in I$ and
  every $v\in T_{x(t)}M$\cfadd{definition:tangent-space} it holds that
  \begin{eqn}
    A(t)v\in T_{x(t)}M,
    \qquad \langle A(t)v,v\rangle\geq\mathsf a\norm{v}^2,
    \qquad
    \norm{A(t)v}\leq\mathsf b\norm{v},
  \end{eqn}
  \begin{eqn}
    r'(t)=-\alpha(t)(r(t))^2
      (\partial_r \mathbf{h})(r(t),x(t)),
    \qquad\text{and}\qquad
    x'(t)=-\alpha(t)A(t)(\nabla_M\mathbf{h})(r(t),x(t))
  \end{eqn}
  \cfadd{definition:angular-gradient}%
  we have that for all $a,b\in I$ with $a\leq b$ it holds that
  \begin{eqn}
    \llabel{eq12.1}
    Q(r(b),x(b))-Q(r(a),x(a))
      \leq-\int_a^b\norm{x'(t)}\dd t
        +\int_a^b m(r(t))\abs{r'(t)}\dd t
    \dott
  \end{eqn}
  \cfload[.]
  \startnewargseq

  To prove \lref{eq12.1}, let $I\subseteq\R$ be an interval,
  let $r\colon I\to(0,\eta)$ and $x\colon I\to M$
  be locally absolutely continuous,
  let $\alpha\colon I\to[0,\infty)$ and
  $A\colon I\to\R^{D\times D}$ be Lebesgue measurable,
  assume $\midcup_{t\in I}\{(r(t),x(t))\}\subseteq\mathcal O$, assume that for Lebesgue-almost all $t\in I$
  and every $v\in T_{x(t)}M$ it holds that
  \begin{eqn}
    \llabel{setting6.1}
    A(t)v\in T_{x(t)}M,\qquad
    \langle A(t)v,v\rangle\geq\mathsf a\norm{v}^2,
    \qquad
    \norm{A(t)v}\leq\mathsf b\norm{v},
  \end{eqn}
  \begin{eqn}
    \llabel{setting6.2}
    r'(t)=-\alpha(t)(r(t))^2
      (\partial_r \mathbf{h})(r(t),x(t)),
    \qquad\text{and}\qquad
    x'(t)=-\alpha(t)A(t)(\nabla_M\mathbf{h})(r(t),x(t)),
  \end{eqn}
  and let
  $p\colon I\to\R^D$, $w\colon I\to\R$,
  $\gamma,\ell\colon I\to[0,\infty)$, and
  $\xi_i\colon I\to[0,H]$, $i\in\{0,1,\dots,N-1\}$, satisfy for all $t\in I$,
  $i\in\{0,1,\dots,N-1\}$ that
  \begin{gather}
    p(t)=(\nabla_M\mathbf{h})(r(t),x(t)),\qquad
    w(t)=(\partial_r \mathbf{h})(r(t),x(t)),\qquad
    \gamma(t)=\sqrt{\max\{0,\langle A(t)p(t),p(t)\rangle\}},\notag\\
    \ell(t)=\alpha(t)\norm{A(t)p(t)},\qquad\text{and}\qquad
    \xi_i(t)=Z_i(r(t),x(t)) \llabel{setting7}
    \dott
  \end{gather}
  \cfadd{definition:tangent-space}\cfadd{definition:angular-gradient}\cfload[.]
  \startnewargseq
  \argument{
    \lref{setting1};
    \lref{setting1.1};
    \lref{setting6.1};
    \lref{setting7};
    the Cauchy--Schwarz inequality;
  }{
    that for Lebesgue-almost all $t\in I$ it holds that
    \begin{eqn}
      \llabel{eq13}
      \mathsf c\norm{p(t)}\leq\gamma(t),\qquad
      (\gamma(t))^2\leq\norm{A(t)p(t)}\norm{p(t)},
      \qquad\text{and}\qquad
      \norm{A(t)p(t)}\leq\mathsf b\norm{p(t)}
      \dott
    \end{eqn}
  }
  \argument{
    \lref{eq13};
    \lref{setting1};
    \lref{setting1.1};
    \lref{setting6.2};
    \lref{setting7};
    the fact that for Lebesgue-almost all $t\in I$ with
    $\gamma(t)=0$ it holds that $p(t)=0$;
  }{
    that for Lebesgue-almost all $t\in I$ it holds that
    \begin{eqn}
      \llabel{eq14}
      \mathsf c\alpha(t)\gamma(t)\leq\ell(t)
        =\norm{x'(t)}\leq\Lambda\alpha(t)\gamma(t)
      \dott
    \end{eqn}
  }
  \argument{
    \lref{eq14};
    the fact that $x$ is locally absolutely continuous;
  }{
    that for all $a,b\in I$ with $a<b$ it holds that
    \begin{eqn}
      \llabel{eq15}
      \int_a^b\alpha(t)\gamma(t)\dd t
        \leq\frac1{\mathsf c}\int_a^b\norm{x'(t)}\dd t<\infty
      \dott
    \end{eqn}
  }
  \argument{
    \lref{setting6.1};
    \lref{setting6.2};
    \lref{setting7};
    the fact that $\mathbf{h}\in \COne{\mathcal O}{(-1,1)}$;
    the fact that for all compact subintervals $J\subseteq I$ it holds
    that $\{(r(t),x(t))\colon t\in J\}$ is a compact subset of $\mathcal O$;
    the chain rule;
  }{
    that $I\ni t\mapsto\mathbf{h}(r(t),x(t))\in\R$
    is locally absolutely continuous
    and that for Lebesgue-almost all $t\in I$ it holds that
    \begin{eqn}[a]
      \llabel{eq16}
      \frac{\ddd}{\ddd t}\mathbf{h}(r(t),x(t))
        &=(\partial_r \mathbf{h})(r(t),x(t))r'(t)
          +\langle(\nabla_M\mathbf{h})(r(t),x(t)),x'(t)\rangle
          =w(t)r'(t)+\langle p(t),x'(t)\rangle\\
        &=-\alpha(t)(r(t))^2(w(t))^2
          -\alpha(t)\langle A(t)p(t),p(t)\rangle
          =-\alpha(t)\bigl((\gamma(t))^2+(r(t))^2(w(t))^2\bigr)
      \dott
    \end{eqn}
  }
  \argument{
    \lref{eq2};
    \lref{setting7};
    the fact that $r$ and $x$ are continuous;
    the fact that for all $i\in\{0,1,\dots,N-1\}$ it holds that
    $Z_i$ is continuous;
  }{
    that for all $i\in\{0,1,\dots,N-1\}$ it holds that $\xi_i$ is
    continuous, that $\{u\in I\colon0<\xi_i(u)<H\}$ is open in $I$,
    and that
    \begin{eqn}
      \llabel{eq16.1}
      \{u\in I\colon0<\xi_i(u)<H\}
        =\{u\in I\colon(r(u),x(u))\in\Omega_i\}
      \dott
    \end{eqn}
  }
  \argument{
    \lref{eq16.1};
    \lref{eq2.0};
    \lref{setting7};
    the fact that $r$ and $x$ are locally absolutely continuous;
  }{
    that for all $i\in\{0,1,\dots,N-1\}$ it holds that
    $\xi_i|_{\{u\in I\colon0<\xi_i(u)<H\}}$
    is locally absolutely continuous\dott
    \llabel{eq16.2}
  }
  \argument{
    \lref{eq16.2};
    \lref{setting2};
    \lref{setting2.1};
    \lref{eq4};
    \lref{setting7};\lref{eq16.1};
    the fact that for all $i\in\{0,1,\dots,N-1\}$ it holds that
    $E_{i+1}\in C^2((0,\theta)\times(0,H),\R)$;
    the fact that $r$ is locally absolutely continuous;
  }{
    that for all $i\in\{0,1,\dots,N-1\}$ it holds that
    for Lebesgue-almost all
    $t\in\{u\in I\colon0<\xi_i(u)<H\}$ it holds that $\xi_i$
    is differentiable at $t$ and that
    \begin{eqn}[a]
      \llabel{eq16.3}
      &\frac{\ddd}{\ddd t}\mathbf{h}(r(t),x(t))
        =\frac{\ddd}{\ddd t}E_{i+1}(r(t),\xi_i(t))
          =(E_{i+1})'(r(t),\xi_i(t))\bigl(r'(t),(\xi_i)'(t)\bigr)\\
      &\quad=\bigl(\partial_s E_{i+1}\bigr)\bigl(r(t),\xi_i(t)\bigr)r'(t)
          +\bigl(\partial_z E_{i+1}\bigr)\bigl(r(t),\xi_i(t)\bigr)(\xi_i)'(t)\\
      &\quad=S_i(r(t),\xi_i(t))r'(t)
          +T_{i+1}(r(t),\xi_i(t))(\xi_i)'(t)
      \dott
    \end{eqn}
  }
  \argument{
    \lref{eq16.3};
    \lref{eq4};
    \lref{setting6.2};
    \lref{setting7};
     \lref{eq16};
  }{
    that for all $i\in\{0,1,\dots,N-1\}$ it holds that
    for Lebesgue-almost all
    $t\in\{u\in I\colon0<\xi_i(u)<H\}$ it holds that
    \begin{eqn}[a]
      \llabel{eq17}
      (\xi_i)'(t)
        &=\frac{\frac{\ddd}{\ddd t}\mathbf{h}(r(t),x(t))
          -S_i(r(t),\xi_i(t))r'(t)}{T_{i+1}(r(t),\xi_i(t))}\\
        &=-\frac{\alpha(t)}{T_{i+1}(r(t),\xi_i(t))}
          \bigl((\gamma(t))^2+(r(t))^2w(t)
            [w(t)-S_i(r(t),\xi_i(t))]\bigr)
      \dott
    \end{eqn}
  }
  \argument{
    \lref{setting1};
    \lref{setting1.1};
    \lref{eq4};
    \lref{setting7};
    \lref{eq13};
    \lref{eq16.1};
  }{
    that for all $i\in\{0,1,\dots,N-1\}$ it holds that
    for Lebesgue-almost all
    $t\in\{u\in I\colon0<\xi_i(u)<H\}$ it holds that
    \begin{eqn}
      \llabel{eq17.1}
      0<T_{i+1}(r(t),\xi_i(t))\leq\mathsf c G(r(t),x(t))
        =\mathsf c\norm{p(t)}\leq\gamma(t)
      \dott
    \end{eqn}
  }
  \argument{
    \lref{eq17.1};
    \lref{eq6};
    \lref{setting7};
    \lref{eq16.1};
    the fact that for all $t\in I$ it holds that $0<r(t)<\eta<\tau$;
    the fact that for all $i\in\{0,1,\dots,N-1\}$, $s\in(0,\tau)$
    it holds that $\mathbf m_i(s)\geq1$;
  }{
    that for all $i\in\{0,1,\dots,N-1\}$ it holds that
    for Lebesgue-almost all
    $t\in\{u\in I\colon0<\xi_i(u)<H\}$ it holds that
    \begin{eqn}[a]
      \llabel{eq17.1.1}
      \abs[\big]{w(t)-S_i(r(t),\xi_i(t))}
       & \leq\mathbf m_i(r(t))
          \bigl(T_{i+1}(r(t),\xi_i(t))+\mathsf c G(r(t),x(t))\bigr)\\
          &
        \leq\mathbf m_i(r(t))
          \bigl(T_{i+1}(r(t),\xi_i(t))+\gamma(t)\bigr)
      \dott
    \end{eqn}
  }
  \argument{
    \lref{eq7};
    the fact that for all $t\in I$ it holds that $0<r(t)<\eta<\tau$;
    the fact that for all $i\in\{0,1,\dots,N-1\}$, $s\in(0,\tau)$
    it holds that $\mathbf m_i(s)\geq1$;
  }{
    that for all $i\in\{0,1,\dots,N-1\}$, $t\in I$ it holds that
    \begin{eqn}
      \llabel{eq17.1.2}
      r(t)>0,\qquad \mathbf m_i(r(t))\geq1,\qquad\text{and}\qquad
        (r(t))^2(\mathbf m_i(r(t)))^2\leq\frac1{64}
      \dott
    \end{eqn}
  }
  \argument{
    \lref{eq17.1.2};
    \lref{setting4};
    \lref{setting6.2};
    \lref{setting7};
    \lref{eq17};
    \lref{eq17.1};
    the fact that for all $t\in I$ it holds that $\alpha(t)\geq0$;
  }{
    for all $i\in\{0,1,\dots,N-1\}$ 
    and Lebesgue-almost all
    $t\in\{u\in I\colon0<\xi_i(u)<H\}$ that
    \begin{equation}
      \llabel{eq17.2}
      \begin{aligned}
      &\omega_i(r(t),\xi_i(t))(\xi_i)'(t)\\
        &=-\alpha(t)
          \frac{(\gamma(t))^2+(r(t))^2w(t)
            [w(t)-S_i(r(t),\xi_i(t))]}
          {T_{i+1}(r(t),\xi_i(t))+8(r(t))^2\mathbf m_i(r(t))
            \sqrt{(S_i(r(t),\xi_i(t)))^2+(T_{i+1}(r(t),\xi_i(t)))^2}}
      \end{aligned}
    \end{equation}
    and 
    \begin{eqn}
      \llabel{eq17.3}
      \abs{r'(t)}
        =\abs[\big]{-\alpha(t)(r(t))^2w(t)}
        =\alpha(t)(r(t))^2\abs{w(t)}
      \dott
    \end{eqn}
  }
  \argument{
    \lref{eq17.2}; \lref{eq17.3};
    \lref{eq17.1};
    \lref{eq17.1.1}; \lref{eq17.1.2};
    \itref{lem:attenuation}{item1,item2} (applied with
    $r\curvearrowleft r(t)$,
    $T\curvearrowleft T_{i+1}(r(t),\xi_i(t))$,
    $g\curvearrowleft\gamma(t)$,
    $n\curvearrowleft\mathbf m_i(r(t))$,
    $v\curvearrowleft w(t)$,
    $c\curvearrowleft S_i(r(t),\xi_i(t))$ in the notation of \cref{lem:attenuation});
    the fact that for all $t\in I$ it holds that $\alpha(t)\geq0$;
  }{
    that for all $i\in\{0,1,\dots,N-1\}$ it holds that
    for Lebesgue-almost all
    $t\in\{u\in I\colon0<\xi_i(u)<H\}$ it holds that
    \begin{equation}
      \llabel{eq18}
      \textstyle\omega_i(r(t),\xi_i(t))(\xi_i)'(t)
        \leq 2\alpha(t)\mathbf m_i(r(t))(r(t))^2\abs{w(t)}
        =2\mathbf m_i(r(t))\abs{r'(t)}
    \end{equation}
    and 
    \begin{eqn}
      \llabel{eq18.1}
     \omega_i(r(t),\xi_i(t))(\xi_i)'(t)
        \leq5\alpha(t)\mathbf m_i(r(t))(r(t))^2\abs{w(t)}
          -\frac{\alpha(t)\gamma(t)}{4}
        =5\mathbf m_i(r(t))\abs{r'(t)}
          -\frac{\alpha(t)\gamma(t)}{4}
      \dott
    \end{eqn}
  }
  In the following for every $i\in\{0,1,\dots,N-1\}$,
  $\varepsilon\in(0,\nicefrac{H}{2})$ let
  $Y_i^\varepsilon\colon\mathcal O\cap((0,\eta)\times M)\to[\varepsilon,H-\varepsilon]$ and
  $q_i^\varepsilon\colon\mathcal O\cap((0,\eta)\times M)\to\R$ satisfy for all 
  $(s,y)\in\mathcal O\cap((0,\eta)\times M)$ that
  \begin{equation}
    \llabel{setting8}
    \textstyle Y_i^\varepsilon(s,y)=\max\{\varepsilon,\min\{Z_i(s,y),H-\varepsilon\}\}\qquad\text{and}
    \qquad
    q_i^\varepsilon(s,y)=P_i(s,Y_i^\varepsilon(s,y)),
  \end{equation}
 for every $\varepsilon\in(0,\nicefrac{H}{2})$ let
  $V_\varepsilon\colon\mathcal O\cap((0,\eta)\times M)\to\R$ and
  $\mathcal U_\varepsilon\subseteq\mathcal O\cap((0,\eta)\times M)$ satisfy for all
  $(s,y)\in\mathcal O\cap((0,\eta)\times M)$ that
  \begin{eqn}
    \llabel{setting9}
    \textstyle V_\varepsilon(s,y)=\kappa\sum\limits_{i=0}^{N-1}q_i^\varepsilon(s,y)
    \qquad\text{and}\qquad
    \mathcal U_\varepsilon
      =\bigcup\limits_{i=0}^{N-1}\{(u,v)\in\mathcal O\cap((0,\eta)\times M)
        \colon\varepsilon<Z_i(u,v)<H-\varepsilon\},
  \end{eqn}
  \startnewargseq
  and let $\Gamma\subseteq\mathcal O\cap((0,\eta)\times M)$ satisfy
  \begin{eqn}
    \llabel{setting10}
    \textstyle\Gamma=\bigcup\limits_{j=0}^{N}\{(s,y)\in\mathcal O\cap((0,\eta)\times M)
      \colon \mathbf{h}(s,y)=\mathfrak e_j(s)\}
    \dott
  \end{eqn}
  \cfload[.]
  \startnewargseq
  \argument{
    \lref{eq10};
    \lref{setting5.1};
    \lref{setting5};
    \lref{setting8};
    \lref{setting9};
    the triangle inequality;
    the fact that $\kappa>0$;
    the fact that for all $z\in[0,H]$, $\varepsilon\in(0,\nicefrac{H}{2})$
    it holds that
    $\abs{\max\{\varepsilon,\min\{z,H-\varepsilon\}\}-z}\leq\varepsilon$;
  }{
    that for all $\varepsilon\in(0,\nicefrac{H}{2})$,
    $(s,y)\in\mathcal O\cap((0,\eta)\times M)$ it holds that
    \begin{eqn}[a]
      \llabel{eq19}
     & \abs{V_\varepsilon(s,y)-Q(s,y)}
        =\kappa\abs[\bigg]{\sum_{i=0}^{N-1}
          \bigl(P_i(s,Y_i^\varepsilon(s,y))-P_i(s,Z_i(s,y))\bigr)}\\
        &\leq\kappa\sum_{i=0}^{N-1}
          \abs[\big]{P_i(s,Y_i^\varepsilon(s,y))-P_i(s,Z_i(s,y))}\leq\kappa\sum_{i=0}^{N-1}
          \abs{Y_i^\varepsilon(s,y)-Z_i(s,y)}
        \leq\kappa\sum_{i=0}^{N-1}\varepsilon
        =\kappa N\varepsilon
      \dott
    \end{eqn}
  }
  \argument{
    \lref{setting0};
    \lref{eq1};
    \lref{eq2};
    \lref{eq3};
    \lref{setting9};
    \lref{setting10};
  }{
    that for all $\varepsilon,\varepsilon'\in(0,\nicefrac{H}{2})$ with
    $\varepsilon'<\varepsilon$ it holds that
    $\mathcal U_\varepsilon\subseteq\mathcal U_{\varepsilon'}$ and that
    \begin{eqn}
      \llabel{eq20}
      \midcup_{\varepsilon\in(0,\nicefrac{H}{2})}\mathcal U_\varepsilon
        =(\mathcal O\cap((0,\eta)\times M))\backslash\Gamma
      \dott
    \end{eqn}
  }
  \argument{
    \lref{eq2};
    \lref{eq10};
    \lref{eq15};
    \lref{eq18};
    \lref{eq18.1};
    \lref{setting8};
    \itref{lem:truncated}{item2}
    (applied to $Z_i$ on $\mathcal O\cap((0,\eta)\times M)$);
    the fact that for all $i\in\{0,1,\dots,N-1\}$ it holds that
    $\mathbf m_i$ and $\mathfrak b_i$ are continuous;
  }{
    that for all $i\in\{0,1,\dots,N-1\}$,
    $\varepsilon\in(0,\nicefrac{H}{2})$ it holds that
    $I\ni t\mapsto q_i^\varepsilon(r(t),x(t))\in\R$
    is locally absolutely continuous\dott
  }
  \argument{
    \lref{eq2};
    \lref{eq10};
    \lref{eq15};
    \lref{eq18};
    \lref{eq18.1};
    \lref{setting8};
    \itref{lem:truncated}{item3};
  }{
    that for all $i\in\{0,1,\dots,N-1\}$, $\varepsilon\in(0,\nicefrac{H}{2})$
    it holds that for Lebesgue-almost all $t\in I$ it holds that
    \begin{eqn}
      \llabel{eq21}
      \frac{\ddd}{\ddd t}q_i^\varepsilon(r(t),x(t))
        \leq\bigl(2\mathbf m_i(r(t))+\mathfrak b_i(r(t))\bigr)
          \abs{r'(t)}
      \dott
    \end{eqn}
  }
  \argument{
    \lref{eq2};
    \lref{eq10};
    \lref{eq15};
    \lref{eq18};
    \lref{eq18.1};
    \lref{setting8};
    \itref{lem:truncated}{item4};
  }{
    that for all $i\in\{0,1,\dots,N-1\}$, $\varepsilon\in(0,\nicefrac{H}{2})$
    it holds that for Lebesgue-almost all
    $t\in\{u\in I\colon\varepsilon<\xi_i(u)<H-\varepsilon\}$ it holds that
    \begin{eqn}
      \llabel{eq22}
      \frac{\ddd}{\ddd t}q_i^\varepsilon(r(t),x(t))
        \leq\bigl(5\mathbf m_i(r(t))+\mathfrak b_i(r(t))\bigr)
            \abs{r'(t)}
          -\tfrac{1}{4}\alpha(t)\gamma(t)
      \dott
    \end{eqn}
  }
  \argument{
    \lref{setting9};
    \lref{eq21};
  }{
    that for all $\varepsilon\in(0,\nicefrac{H}{2})$ and Lebesgue-almost
    all $t\in I$ it holds that
    \begin{eqn}
      \llabel{eq22.1}
      \frac{\ddd}{\ddd t}V_\varepsilon(r(t),x(t))
        =\kappa\sum_{i=0}^{N-1}
          \frac{\ddd}{\ddd t}q_i^\varepsilon(r(t),x(t))
      \dott
    \end{eqn}
  }
  \argument{
    \lref{eq22.1};
    \lref{eq21};
    \lref{eq22};
    the fact that for all $i\in\{0,1,\dots,N-1\}$ it holds that
    $\mathbf m_i$ is nonnegative;
    the fact that $\kappa>0$;
  }{
    that for all $\varepsilon\in(0,\nicefrac{H}{2})$ and Lebesgue-almost
    all $t\in I$ it holds that
    \begin{eqn}
      \llabel{eq22.2}
      \frac{\ddd}{\ddd t}V_\varepsilon(r(t),x(t))
        \leq-\frac{\kappa\alpha(t)\gamma(t)}4
          \sum_{i=0}^{N-1}\indicator{(\varepsilon,H-\varepsilon)}(\xi_i(t))
        +\kappa\sum_{i=0}^{N-1}
          \bigl(5\mathbf m_i(r(t))+\mathfrak b_i(r(t))\bigr)
          \abs{r'(t)}
      \dott
    \end{eqn}
  }
  \argument{
    \lref{setting7};
    \lref{setting9};
  }{
    that for all $\varepsilon\in(0,\nicefrac{H}{2})$, $t\in I$
    it holds that
    \begin{eqn}
      \llabel{eq22.2.1}
      \indicator{\mathcal U_\varepsilon}(r(t),x(t))
        =\min\biggl\{1,\sum_{i=0}^{N-1}
          \indicator{(\varepsilon,H-\varepsilon)}(\xi_i(t))\biggr\}
        \leq\sum_{i=0}^{N-1}
          \indicator{(\varepsilon,H-\varepsilon)}(\xi_i(t))
      \dott
    \end{eqn}
  }
  \argument{
    \lref{eq22.2.1};
    the fact that for all $t\in I$ it holds that $\alpha(t)\geq0$
    and $\gamma(t)\geq0$;
    the fact that $\kappa>0$;
  }{
    that for all $\varepsilon\in(0,\nicefrac{H}{2})$, $t\in I$
    it holds that
    \begin{eqn}
      \llabel{eq22.2.2}
      -\frac{\kappa\alpha(t)\gamma(t)}4
        \sum_{i=0}^{N-1}\indicator{(\varepsilon,H-\varepsilon)}(\xi_i(t))
        \leq-\frac{\kappa\alpha(t)\gamma(t)}4
          \indicator{\mathcal U_\varepsilon}(r(t),x(t))
      \dott
    \end{eqn}
  }
  \argument{
    \lref{setting1};
    \lref{setting1.1};
    \lref{setting5.1};
    \lref{setting5};
    the fact that for all $t\in I$ it holds that $r(t)\in(0,\eta)$;
  }{
    that for all $t\in I$ it holds that
    \begin{eqn}
      \llabel{eq22.2.3}
      \kappa\sum_{i=0}^{N-1}
        \bigl(5\mathbf m_i(r(t))+\mathfrak b_i(r(t))\bigr)
        =\mathbb{M}(r(t))
      \qquad\text{and}\qquad
      \frac\kappa4=\Lambda
      \dott
    \end{eqn}
  }
  \argument{
    \lref{eq22.2.3};
    \lref{eq22.2};
    \lref{eq22.2.2};
  }{
    that for all $\varepsilon\in(0,\nicefrac{H}{2})$ and Lebesgue-almost
    all $t\in I$ it holds that
    \begin{eqn}[a]
      \llabel{eq22.3}
      \frac{\ddd}{\ddd t}V_\varepsilon(r(t),x(t))
        &\leq-\frac{\kappa\alpha(t)\gamma(t)}4
          \indicator{\mathcal U_\varepsilon}(r(t),x(t))
        +\mathbb{M}(r(t))\abs{r'(t)}
        \\
        &=-\Lambda\alpha(t)\gamma(t)
          \indicator{\mathcal U_\varepsilon}(r(t),x(t))
          +\mathbb{M}(r(t))\abs{r'(t)}
      \dott
    \end{eqn}
  }
  \argument{
    \lref{eq22.3};
    \lref{eq14};
    the fact that for all $\varepsilon\in(0,\nicefrac{H}{2})$, $t\in I$
    it holds that $\indicator{\mathcal U_\varepsilon}(r(t),x(t))\geq0$;
  }{
    that for all $\varepsilon\in(0,\nicefrac{H}{2})$ and Lebesgue-almost
    all $t\in I$ it holds that
    \begin{eqn}
      \llabel{eq23}
      \frac{\ddd}{\ddd t}V_\varepsilon(r(t),x(t))
        \leq-\ell(t)\indicator{\mathcal U_\varepsilon}(r(t),x(t))
          +\mathbb{M}(r(t))\abs{r'(t)}
      \dott
    \end{eqn}
  }
  \argument{
    \lref{eq23};
    \lref{setting9};
    the fact that for all $i\in\{0,1,\dots,N-1\}$,
    $\varepsilon\in(0,\nicefrac{H}{2})$ it holds that
    $I\ni t\mapsto q_i^\varepsilon(r(t),x(t))\in\R$
    is locally absolutely continuous;
    the fundamental theorem of calculus;
  }{
    that for all $\varepsilon\in(0,\nicefrac{H}{2})$, $a,b\in I$ with $a<b$
    it holds that
    \begin{eqn}[a]
      \llabel{eq24}
      V_\varepsilon(r(b),x(b))-V_\varepsilon(r(a),x(a))
        \leq-\int_a^b\ell(t)\indicator{\mathcal U_\varepsilon}(r(t),x(t))\dd t+\int_a^b\mathbb{M}(r(t))\abs{r'(t)}\dd t
      \dott
    \end{eqn}
  }
  \argument{
    \lref{eq19};
  }{
    that for all $(s,y)\in\mathcal O\cap((0,\eta)\times M)$
    it holds that
    \begin{eqn}
      \llabel{eq24.2}
      \lim_{\varepsilon\to 0}V_\varepsilon(s,y)=Q(s,y)
      \dott
    \end{eqn}
  }
  \argument{
    \lref{eq24.2};
  }{
    that for all $a,b\in I$ with $a<b$ it holds that
    \begin{eqn}[a]
      \llabel{eq24.3}
      &\lim_{\varepsilon\to 0}
        \bigl(V_\varepsilon(r(b),x(b))-V_\varepsilon(r(a),x(a))\bigr)
        =Q(r(b),x(b))-Q(r(a),x(a))
      \dott
    \end{eqn}
  }
  \argument{
    \lref{eq14};
    the fact that $\ell$ is nonnegative and measurable;
    the fact that $\mathbb{M}$ is continuous and nonnegative;
    the fact that $r$ and $x$ are locally absolutely continuous;
    the fundamental theorem of calculus for absolutely continuous functions;
    the extreme value theorem;
  }{
    that for all $a,b\in I$ with $a<b$ it holds that
    \begin{eqn}[a]
      \llabel{eq24.4}
      0\leq\int_a^b\ell(t)\dd t
        &=\int_a^b\norm{x'(t)}\dd t<\infty\qquad\text{and}\qquad
      0\leq\int_a^b\mathbb{M}(r(t))\abs{r'(t)}\dd t
        <\infty
      \dott
    \end{eqn}
  }
  \argument{\lref{setting9};\lref{eq20};
    \lref{eq24.4};
    the fact that $\ell$ is nonnegative and measurable;
    the monotone convergence theorem
  }{
    that for all $a,b\in I$ with $a<b$ it holds that
    \begin{eqn}
      \llabel{eq24.6}
      \lim_{\varepsilon\to 0}\int_a^b
        \ell(t)\indicator{\mathcal U_\varepsilon}(r(t),x(t))\dd t
        =\int_a^b\ell(t)
          \indicator{(\mathcal O\cap((0,\eta)\times M))\backslash\Gamma}(r(t),x(t))\dd t
        <\infty
      \dott
    \end{eqn}
  }
  \argument{
    \lref{eq24.6};
    \lref{eq24};
    \lref{eq24.3};
    \lref{eq24.4};
  }{
    that for all $a,b\in I$ with $a<b$ it holds that
    \begin{eqn}[a]
      \llabel{eq25}
      Q(r(b),x(b))-Q(r(a),x(a))
        \leq-\int_a^b\ell(t)
          \indicator{(\mathcal O\cap((0,\eta)\times M))\backslash\Gamma}(r(t),x(t))\dd t+\int_a^b\mathbb{M}(r(t))\abs{r'(t)}\dd t
      \dott
    \end{eqn}
  }
  In the following for all $j\in\{0,1,\dots,N\}$ let
  $f_j\colon I\to\R$ satisfy that for all $t\in I$ it holds that
  \begin{eqn}
    \llabel{setting11}
    f_j(t)=\mathbf{h}(r(t),x(t))-\mathfrak e_j(r(t))
    \dott
  \end{eqn}
  \cfload[.]
  \startnewargseq
  \argument{
    \lref{setting6.2};
    \lref{setting7};
    \lref{eq16};
    \lref{setting11};
    the fact that for all $j\in\{0,1,\dots,N\}$ it holds that
    $\mathfrak e_j\in C^1((0,\theta),\R)$;
    \cref{lem:AC-level-set};
    the chain rule;
  }{
    that for all $j\in\{0,1,\dots,N\}$ it holds that
    for Lebesgue-almost all
    $t\in f_j^{-1}(\{0\})$ it holds that
    \begin{eqn}
      \llabel{eq26}
      0=(f_j)'(t)
        =-\alpha(t)\bigl((\gamma(t))^2
          +(r(t))^2w(t)[w(t)-(\mathfrak e_j)'(r(t))]\bigr)
      \dott
    \end{eqn}
  }
  \argument{
    \lref{eq26};
    \lref{setting6.2};
    \lref{setting7};
    the fact that for all $t\in I$ it holds that $\alpha(t)\geq0$;
  }{
    that for all $j\in\{0,1,\dots,N\}$ it holds that
    for Lebesgue-almost all
    $t\in f_j^{-1}(\{0\})$ it holds that
    \begin{eqn}[a]
      \llabel{eq26.1}
      \alpha(t)(\gamma(t))^2
        &=-\alpha(t)(r(t))^2w(t)
          [w(t)-(\mathfrak e_j)'(r(t))]
          =[w(t)-(\mathfrak e_j)'(r(t))]r'(t)\\
        &=\abs{[w(t)-(\mathfrak e_j)'(r(t))]r'(t)}
          =\abs{w(t)-(\mathfrak e_j)'(r(t))}\abs{r'(t)}
      \dott
    \end{eqn}
  }
  \argument{
    \lref{setting11};
  }{
    that for all $j\in\{0,1,\dots,N\}$, $t\in f_j^{-1}(\{0\})$
    it holds that
    \begin{eqn}
      \llabel{eq26.2}
      \mathbf{h}(r(t),x(t))=f_j(t)+\mathfrak e_j(r(t))=\mathfrak e_j(r(t))
      \dott
    \end{eqn}
  }
  \argument{
    \lref{eq26.2};
    \lref{eq11};
    \lref{setting7};
  }{
    that for all $j\in\{0,1,\dots,N\}$, $t\in f_j^{-1}(\{0\})$
    with $G(r(t),x(t))>0$ it holds that
    \begin{eqn}
      \llabel{eq26.3}
      \abs{w(t)-(\mathfrak e_j)'(r(t))}
        =\abs[\big]{(\partial_r \mathbf{h})(r(t),x(t))
          -(\mathfrak e_j)'(r(t))}
        \leq\beta_j(r(t))G(r(t),x(t))
      \dott
    \end{eqn}
  }
  \argument{
    \lref{setting1};
    \lref{setting1.1};
    \lref{setting7};
    \lref{eq13};
  }{
    that for Lebesgue-almost all $t\in I$ it holds that
    \begin{eqn}
      \llabel{eq26.4}
      \mathsf cG(r(t),x(t))=\mathsf c\norm{p(t)}\leq\gamma(t)
      \dott
    \end{eqn}
  }
  \argument{
    \lref{eq26.4};
    \lref{eq26.1};
    \lref{eq26.3};
    the fact that $\mathsf c>0$;
    the fact that for all $j\in\{0,1,\dots,N\}$ it holds that
    $\beta_j$ is nonnegative;
  }{
    that for all $j\in\{0,1,\dots,N\}$ it holds that
    for Lebesgue-almost all
    $t\in f_j^{-1}(\{0\})$ with $G(r(t),x(t))>0$ it holds that
    \begin{eqn}[a]
      \llabel{eq27}
      \alpha(t)(\gamma(t))^2
        \leq\beta_j(r(t))G(r(t),x(t))\abs{r'(t)}
        \leq\frac{\beta_j(r(t))}{\mathsf c}\gamma(t)\abs{r'(t)}
      \dott
    \end{eqn}
  }
  \argument{
    \lref{eq27};
    \lref{eq13};
    \lref{eq14};
  }{
    that for all $j\in\{0,1,\dots,N\}$ it holds that
    for Lebesgue-almost all
    $t\in f_j^{-1}(\{0\})$ it holds that
    \begin{eqn}
      \llabel{eq28}
      \ell(t)\leq\frac{\Lambda}{\mathsf c}
        \beta_j(r(t))\abs{r'(t)}
      \dott
    \end{eqn}
  }
  \argument{
    \lref{eq28};
    \lref{setting10};
    \lref{setting11};
    the fact that for all $j\in\{0,1,\dots,N\}$ it holds that
    $\beta_j$ is nonnegative;
  }{
    that for all $a,b\in I$ with $a<b$ it holds that
    \begin{eqn}
      \llabel{eq29}
      \int_a^b\ell(t)\indicator{\Gamma}(r(t),x(t))\dd t
        \leq\frac{\Lambda}{\mathsf c}
          \int_a^b\sum_{j=0}^{N}\beta_j(r(t))\abs{r'(t)}\dd t
      \dott
    \end{eqn}
  }
  \argument{
    \lref{eq29};
    \lref{setting5.1};
    \lref{setting5};
    \lref{eq14};
    \lref{eq24.4};
    \lref{eq25};
    the fact that $r$ is locally absolutely continuous;
    the fundamental theorem of calculus for absolutely continuous functions;
    the extreme value theorem;
  }{
    that for all $a,b\in I$ with $a<b$ it holds that
    \begin{eqn}[a]
      \llabel{eq30}
      &Q(r(b),x(b))-Q(r(a),x(a))\\
      &\quad\leq-\int_a^b\ell(t)
          \indicator{(\mathcal O\cap((0,\eta)\times M))\backslash\Gamma}(r(t),x(t))\dd t
          +\int_a^b\mathbb{M}(r(t))\abs{r'(t)}\dd t\\
      &\quad=-\int_a^b\ell(t)\dd t
          +\int_a^b\ell(t)\indicator{\Gamma}(r(t),x(t))\dd t
          +\int_a^b\mathbb{M}(r(t))\abs{r'(t)}\dd t\\
      &\quad\leq-\int_a^b\ell(t)\dd t
          +\frac{\Lambda}{\mathsf c}
            \int_a^b\sum_{j=0}^{N}\beta_j(r(t))\abs{r'(t)}\dd t
          +\int_a^b\mathbb{M}(r(t))\abs{r'(t)}\dd t\\
      &\quad=-\int_a^b\ell(t)\dd t
          +\int_a^b\biggl(\mathbb{M}(r(t))
            +\frac{\Lambda}{\mathsf c}\sum_{j=0}^{N}\beta_j(r(t))\biggr)
            \abs{r'(t)}\dd t\\
      &\quad=-\int_a^b\ell(t)\dd t
          +\int_a^b m(r(t))\abs{r'(t)}\dd t\\
      &\quad=-\int_a^b\norm{x'(t)}\dd t
          +\int_a^b m(r(t))\abs{r'(t)}\dd t
      \dott
    \end{eqn}
  }
  \argument{
    \lref{eq30};
    the arbitrariness of $I$, $r$, $x$, $\alpha$, and $A$;
  }[verbs=ep]{
    \lref{eq12.1}\dott
  }

  To prove \lref{concl2}, let $I\subseteq\R$ be an interval,
  let $r\colon I\to(0,\eta)$ and $x\colon I\to M$
  be locally absolutely continuous,
  let $\alpha\colon I\to[0,\infty)$ and
  $A\colon I\to\R^{D\times D}$ be Lebesgue measurable,
  assume $\midcup_{t\in I}\{(r(t),x(t))\}\subseteq\mathcal O$,
  and assume that for Lebesgue-almost all $t\in I$ and
  every $v\in T_{x(t)}M$\cfadd{definition:tangent-space} it holds that
  \begin{eqn}
    \llabel{eq0.2.1}
    A(t)v\in T_{x(t)}M,
    \qquad \langle A(t)v,v\rangle\geq\mathsf a\norm{v}^2,
    \qquad
    \norm{A(t)v}\leq\mathsf b\norm{v},
  \end{eqn}
  \begin{eqn}
    \llabel{eq0.2.2}
    r'(t)=-\alpha(t)(r(t))^2
      (\partial_r h)(r(t),x(t)),
    \qquad\text{and}\qquad
    x'(t)=-\alpha(t)A(t)(\nabla_Mh)(r(t),x(t))
    \dott
  \end{eqn}
  \cfadd{definition:angular-gradient}\cfload[.]
  \startnewargseq
  In the following let $\lambda\colon I\to[0,\infty)$ satisfy
  that for all $t\in I$ it holds that
  \begin{eqn}
    \llabel{eq0.3}
    \lambda(t)=\frac{\alpha(t)}{\varphi'(h(r(t),x(t)))}
      =\alpha(t)\bigl(1+(h(r(t),x(t)))^2\bigr)^{3/2}
    \dott
  \end{eqn}
  \cfload[.]
  \startnewargseq
  \argument{
    \lref{eq0.3};
    the fact that $h$ is continuous;
    the fact that $\alpha$ is nonnegative and Lebesgue measurable;
    the fact that $r$ and $x$ are locally absolutely continuous;
  }{
    that $r$ and $x$ are continuous and that $\lambda$ is
    nonnegative and Lebesgue measurable\dott
  }
  \argument{
    \lref{eq0.1};
    \lref{eq0.2};
    \lref{eq0.2.2};
    \lref{eq0.3};
  }{
    that for Lebesgue-almost all $t\in I$ it holds that
    \begin{eqn}[a]
      \llabel{eq0.4}
      r'(t)&=-\alpha(t)(r(t))^2
        (\partial_r h)(r(t),x(t))
        =-\lambda(t)(r(t))^2
          (\partial_r \mathbf{h})(r(t),x(t))\\
     \text{and}\qquad x'(t)&=-\alpha(t)A(t)(\nabla_Mh)(r(t),x(t))
        =-\lambda(t)A(t)(\nabla_M\mathbf{h})(r(t),x(t))
      \dott
    \end{eqn}
  }
  \argument{
    \lref{eq0.4};
    \lref{eq12.1} ;
    \lref{eq0.2.1};
    the fact that $\lambda$ is nonnegative and Lebesgue measurable;
    the fact that $A$ is Lebesgue measurable;
    the fact that $r$ and $x$ are locally absolutely continuous;
    the fact that $\midcup_{t\in I}\{(r(t),x(t))\}\subseteq\mathcal O$;
  }{
    that for all $a,b\in I$ with $a\leq b$ it holds that
    \begin{eqn}
      \llabel{eq31}
      Q(r(b),x(b))-Q(r(a),x(a))
        \leq-\int_a^b\norm{x'(t)}\dd t
          +\int_a^b m(r(t))\abs{r'(t)}\dd t
      \dott
    \end{eqn}
  }
  \argument{
    \lref{eq31};
    the arbitrariness of $I$, $r$, $x$, $\alpha$, and $A$;
  }[verbs=ep]{
    \lref{concl2}\dott
  }
\end{aproof}

\subsection{Upper bounds for the error in the Lyapunov-type estimate}
\label{subsec:radial-integrability}

\begin{athm}{lemma}{lem:v}
  Let $\varrho,\mathfrak C\in(0,\infty)$, let $I\subseteq\R$ be an interval,
  let $r\colon I\to(0,\infty)$ be locally absolutely continuous,
  for every $R\in(0,\varrho)$ let
  $\boldsymbol{\nu}_R\colon\mathcal B((0,R])\to[0,\infty]$
  satisfy for all $B\in\mathcal B((0,R])$ that
  \begin{eqn}
    \llabel{ass3}
    \boldsymbol{\nu}_R(B)
      =\int_{r^{-1}(B)} | r'(t) | \dd t,
  \end{eqn}
  let $\mathfrak v\colon I\to[0,\infty)$ be Lebesgue measurable,
  assume for Lebesgue-almost all $t\in I$ that
  $| r'(t) | \leq \mathfrak v(t)$,
  and assume for all $\rho\in(0,\varrho)$ that
  $
    \int_{ r^{ - 1 }( [\rho/2,\rho])}\mathfrak v(t)\dd t
      \leq\mathfrak C\rho
    \dott
  $
  \cfload[.]
  Then
  \begin{enumerate}[label=(\roman*)]
    \item \llabel{item1} for all $\rho\in(0,\varrho)$ it holds that
    \begin{eqn}
      \llabel{concl1}
      \int_{r^{-1}((0,\rho])}\mathfrak v(t)\dd t
        \leq 2\mathfrak C\rho,
    \end{eqn}
    \item \llabel{item2} for all $R\in(0,\varrho)$ it holds that
    $\boldsymbol{\nu}_R$ is a finite measure, and
    \item \llabel{item3} for all $R\in(0,\varrho)$, $\rho\in(0,R]$
    it holds that
    \begin{eqn}[i]
      \boldsymbol{\nu}_R((0,\rho])\leq 2\mathfrak C\rho
      \dott
    \end{eqn}
  \end{enumerate}
  \cfout[.]
\end{athm}
\begin{aproof}
  Throughout this proof for every $\rho\in(0,\varrho)$,
  $j\in\N_0$ let $\mathcal T_{\rho,j}\subseteq I$ satisfy
  \begin{eqn}
    \llabel{setting1}
    \mathcal T_{\rho,j}
      =r^{-1}((2^{-j-1}\rho,2^{-j}\rho])
    \dott
  \end{eqn}
  \cfload[.]
  \startnewargseq
  \argument{
    \lref{setting1};
    the fact that $r$ is continuous;
  }{
    that for all $\rho\in(0,\varrho)$ it holds that
    $\mathcal T_{\rho,j}$, $j\in\N_0$, are Lebesgue measurable,
    are pairwise disjoint, and satisfy that
    \begin{eqn}
      \llabel{eq1}
      \midcup_{j\in\N_0}\mathcal T_{\rho,j}=r^{-1}((0,\rho])
      \dott
    \end{eqn}
  }
  \argument{
    the assumption that for all $\rho\in(0,\varrho)$ it holds that
    $\int_{r^{-1}([\rho/2,\rho])}\mathfrak v(t)\dd t\leq\mathfrak C\rho$;
    \lref{setting1};
    the fact that $\mathfrak v$ is nonnegative;
  }{
    that for all $\rho\in(0,\varrho)$, $j\in\N_0$ it holds that
    \begin{eqn}
      \llabel{eq2}
      \int_{\mathcal T_{\rho,j}}\mathfrak v(t)\dd t
        \leq\int_{r^{-1}([2^{-j-1}\rho,2^{-j}\rho])}
          \mathfrak v(t)\dd t
        \leq\mathfrak C2^{-j}\rho
      \dott
    \end{eqn}
  }
  \argument{
    \lref{eq2};
    \lref{eq1};
    the fact that for all $\rho\in(0,\varrho)$ it holds that
    $\mathcal T_{\rho,j}$, $j\in\N_0$, are pairwise disjoint;
  }{
    that for all $\rho\in(0,\varrho)$ it holds that
    \begin{eqn}[a]
      \llabel{eq3}
      \int_{r^{-1}((0,\rho])}\mathfrak v(t)\dd t
        &=\sum_{j=0}^{\infty}
          \int_{\mathcal T_{\rho,j}}\mathfrak v(t)\dd t
      \leq\mathfrak C\rho\sum_{j=0}^{\infty}2^{-j}
          =2\mathfrak C\rho
      \dott
    \end{eqn}
  }
  \argument{
    \lref{eq3};
  }[verbs=ep]{
    \lref{item1}\dott
  }
  \startnewargseq
  \argument{
    \lref{ass3};
    the fact that $r$ is Borel measurable;
  }{
    that for every $R\in(0,\varrho)$ and every sequence
    $(B_j)_{j\in\N_0}$ of pairwise disjoint Borel subsets of
    $(0,R]$ it holds that $ \boldsymbol{\nu}_R(\emptyset)=0$ and
    \begin{equation}
      \llabel{eq5}
     \textstyle \boldsymbol{\nu}_R\bigl(\midcup_{j\in\N_0}B_j\bigr)
        =\int_{\midcup_{j\in\N_0}r^{-1}(B_j)}\abs{r'(t)}\dd t
      =\sum_{j=0}^{\infty}
          \int_{r^{-1}(B_j)}\abs{r'(t)}\dd t
          =\sum_{j=0}^{\infty}\boldsymbol{\nu}_R(B_j)
      \dott
    \end{equation}
  }
  \argument{
    \lref{eq5};
    the assumption that for Lebesgue-almost all $t\in I$ it holds that
    $|r'(t)|\leq\mathfrak v(t)$;
    \lref{ass3};
    \lref{concl1};
  }{
    that for all $R\in(0,\varrho)$ it holds that
    $\boldsymbol{\nu}_R$ is a nonnegative Borel measure on
    $(0,R]$ and satisfies for all $\rho\in(0,R]$ that
    \begin{eqn}
      \llabel{eq6}
      \boldsymbol{\nu}_R((0,\rho])
        =\int_{r^{-1}((0,\rho])}\abs{r'(t)}\dd t
        \leq\int_{r^{-1}((0,\rho])}\mathfrak v(t)\dd t
        \leq 2\mathfrak C\rho
      \dott
    \end{eqn}
  }
  \argument{
    \lref{eq6};
  }[verbs=ep]{
    \lref{item2}\dott
  }
  \argument{
    \lref{eq6};
  }[verbs=ep]{
    \lref{item3}\dott
  }
\end{aproof}

\begin{athm}{lemma}{lem:radial-integrability}
  Let $R\in(0,\infty)$ and let
  $m\colon(0,R]\to[0,\infty)$ be locally bounded,
  definable\cfadd{definition:definable-function}, and integrable\dott
  \cfload[.]
  Then
  \begin{enumerate}[label=(\roman*)]
    \item \llabel{item1} for every $S\in(0,R]$,
    $K\in(0,\infty)$ and every finite measure
    $\mu\colon\mathcal B((0,S])\to[0,\infty]$ with the property that for all
    $\rho\in(0,S]$ it holds that $\mu((0,\rho])\leq K\rho$
    we have for every non-increasing integrable
    $w\colon(0,S]\to[0,\infty)$ that
    \begin{eqn}
      \llabel{concl1}
      \int_0^S w(s)\,\mu(\mathrm{d}s)
        \leq K\int_0^S w(s)\dd s,
    \end{eqn}
    \item \llabel{item2} there exist $\delta\in(0,R)$ and a definable,
    non-increasing, and integrable
    $\psi\colon(0,\delta]\to[0,\infty)$ such that for all
    $s\in(0,\delta]$ it holds that
    \begin{eqn}
      \llabel{concl2}
      m(s)\leq\psi(s),
    \end{eqn}
    \item \llabel{item3} for every $\kappa\in(0,\infty)$, $\eta\in(0,R]$ and every finite measure
    $\nu\colon\mathcal B((0,R])\to[0,\infty]$ with the property that for all
    $\rho\in(0,\eta]$ it holds that
    $\nu((0,\rho])\leq\kappa\rho$ we have that
    \begin{eqn}
      \llabel{concl3}
      \int_0^R m(s)\,\nu(\mathrm{d}s)<\infty,
    \end{eqn}
    and
    \item \llabel{item4} for every $\delta\in(0,R)$, $\varrho,\mathfrak C\in(0,\infty)$,
    every definable, non-increasing, and integrable
    $\psi\colon(0,\delta]\to[0,\infty)$,
    every interval $I\subseteq\R$,
    and every locally absolutely continuous $r\colon I\to(0,\infty)$
    with the property that for all $s\in(0,\delta]$ it holds
    that $m(s)\leq\psi(s)$,
    and with the property that for all $u\in(0,\varrho)$ it holds that
    \begin{eqn}
      \llabel{ass2}
      \int_{r^{-1}([u/2,u])}\abs{r'(t)}\dd t
        \leq\mathfrak C u
    \end{eqn}
    we have for all $S\in(0,\delta]\cap(0,\varrho)$ that
    \begin{eqn}
      \llabel{concl4}
      \int_{r^{-1}((0,S])}m(r(t))\abs{r'(t)}\dd t
        \leq 2\mathfrak C\int_0^S\psi(s)\dd s<\infty
      \dott
    \end{eqn}
  \end{enumerate}
  \cfout[.]
\end{athm}
\begin{aproof}
  Throughout this proof let $\lambda$ be the Lebesgue measure on
  $\R$\dott \cfload[.]
  \startnewargseq
  We first claim that for every $S\in(0,R]$, $K\in(0,\infty)$,
  every finite nonnegative Borel measure $\mu$ on $(0,S]$,
  and every non-increasing integrable function
  $w\colon(0,S]\to[0,\infty)$ with the property that for all
  $\rho\in(0,S]$ it holds that $\mu((0,\rho])\leq K\rho$
  we have that
  \begin{eqn}
    \llabel{eq0}
    \int_{(0,S]}w(s)\dd\mu(s)
      \leq K\int_0^S w(s)\dd s
    \dott
  \end{eqn}
  \cfload[.]
  \startnewargseq
  To prove \lref{eq0}, let $S\in(0,R]$, $K\in(0,\infty)$,
  let $\mu$ be a finite
  nonnegative Borel measure on $(0,S]$ which satisfies for all
  $\rho\in(0,S]$ that
  \begin{eqn}
    \llabel{setting1}
    \mu((0,\rho])\leq K\rho,
  \end{eqn}
  let $w\colon(0,S]\to[0,\infty)$ be non-increasing and Lebesgue
  integrable, for every $a\in[0,\infty)$ let
  $A_a\subseteq(0,S]$ satisfy
  \begin{eqn}
    \llabel{setting2}
    A_a=\{s\in(0,S]\colon w(s)>a\},
  \end{eqn}
  and let $F\colon(0,S]\times[0,\infty)\to\{0,1\}$ satisfy
  for all $(s,a)\in(0,S]\times[0,\infty)$ that
  \begin{eqn}
    \llabel{setting2.1}
    F(s,a)=\indicator{[0,w(s))}(a)
    \dott
  \end{eqn}
  \cfload[.]
  \startnewargseq
  \argument{
    \lref{setting2};
    the fact that $w$ is non-increasing;
  }{
    that $w$ is Borel measurable and that for every
    $a\in[0,\infty)$ there exists $b\in[0,S]$ such that
    \begin{eqn}
      \llabel{eq1}
      (0,b)\subseteq A_a\subseteq(0,b]
      \qquad\text{and}\qquad
      \lambda(A_a)=b
      \dott
    \end{eqn}
  }
  \argument{
    \lref{eq1};
    \lref{setting1};
    the fact that $\mu(\emptyset)=0$;
  }{
    that for all $a\in[0,\infty)$ it holds that
    \begin{eqn}
      \llabel{eq2}
      \mu(A_a)\leq K\lambda(A_a)
      \dott
    \end{eqn}
  }
  \argument{
    \lref{setting2};
    \lref{setting2.1};
    the fact that $w$ is Borel measurable;
  }{
    that $F$ is nonnegative and measurable and that
    for all $(s,a)\in(0,S]\times[0,\infty)$ it holds that
    \begin{eqn}
      \llabel{eq2.1}
      F(s,a)=\indicator{A_a}(s)
      \dott
    \end{eqn}
  }
  \argument{
    \lref{setting2.1};
    the fact that for all $s\in(0,S]$ it holds that
    $w(s)\in[0,\infty)$;
  }{
    that for all $s\in(0,S]$ it holds that
    \begin{eqn}
      \llabel{eq2.2}
      \int_0^\infty F(s,a)\dd a
        =\lambda([0,w(s)))=w(s)
      \dott
    \end{eqn}
  }
  \argument{
    \lref{eq2.2};
    \lref{eq2.1};
    the fact that $F$ is nonnegative and measurable;
    the fact that $\mu$ is finite and
    $\lambda|_{[0,\infty)}$ is sigma-finite;
    Tonelli's theorem
  }{
    that
    \begin{eqn}[a]
      \llabel{eq2.3}
      \textstyle\int_{(0,S]}w(s)\dd\mu(s)
        =\int_{(0,S]}\bigl(\int_0^\infty F(s,a)\dd a\bigr)
          \dd\mu(s)=\int_0^\infty\bigl(\int_{(0,S]}F(s,a)\dd\mu(s)\bigr)
          \dd a
        =\int_0^\infty\mu(A_a)\dd a
      \dott
    \end{eqn}
  }
  \argument{
    \lref{eq2.1};
    \lref{eq2.2};
    the fact that $F$ is nonnegative and measurable;
    the fact that $\lambda|_{(0,S]}$ is finite and
    $\lambda|_{[0,\infty)}$ is sigma-finite;
    Tonelli's theorem
  }{
    that
    \begin{eqn}[a]
      \llabel{eq2.4}
      \int_0^\infty\lambda(A_a)\dd a
        =\int_0^\infty\biggl(\int_0^S F(s,a)\dd s\biggr)\dd a=\int_0^S\biggl(\int_0^\infty F(s,a)\dd a\biggr)\dd s
        =\int_0^S w(s)\dd s
      \dott
    \end{eqn}
  }
  \argument{
    \lref{eq2.4};
    \lref{eq2};
    \lref{eq2.3};
    the fact that $w$ is integrable;
  }{
    that
    \begin{eqn}[a]
      \llabel{eq3}
      \int_{(0,S]}w(s)\dd\mu(s)
        &=\int_0^\infty\mu(A_a)\dd a
      \leq K\int_0^\infty\lambda(A_a)\dd a
      =K\int_0^S w(s)\dd s<\infty
      \dott
    \end{eqn}
  }
  \argument{
    \lref{eq3};
    the arbitrariness of $S$, $K$, $\mu$, and $w$;
  }[verbs=ep]{
    \lref{item1}\dott
  }
  \startnewargseq
  \argument{
    the fact that $m$ is definable;
    \cite[Theorem~2.1]{Coste};
  }{
    that there exists $\varepsilon\in(0,R)$ which satisfies that
    $m|_{(0,\varepsilon)}$ is non-increasing or non-decreasing\dott
  }
  In the following let $\delta\in(0,\varepsilon)$ and let
  $\psi\colon(0,\delta]\to[0,\infty)$ satisfy for all
  $s\in(0,\delta]$ that
  \begin{eqn}
    \llabel{setting3}
    \psi(s)=
    \begin{cases}
      m(s)&\colon m|_{(0,\varepsilon)}\text{ is non-increasing},\\
      m(\delta)&\colon m|_{(0,\varepsilon)}\text{ is not non-increasing}
    \end{cases}
    \dott
  \end{eqn}
  \cfload[.]
  \startnewargseq
  \argument{
    \lref{setting3};
    the fact that $m$ is definable;
    the fact that $m|_{(0,\varepsilon)}$ is non-increasing or
    non-decreasing;
    the fact that $0<\delta<\varepsilon<R$;
  }{
    that $\psi$ is definable and non-increasing and that for all
    $s\in(0,\delta]$ it holds that
    \begin{eqn}
      \llabel{eq4}
      0\leq m(s)\leq\psi(s)\leq m(s)+m(\delta)
      \dott
    \end{eqn}
  }
  \argument{
    \lref{eq4};
    the assumption that $m$ is integrable;
    the fact that $\psi$ is nonnegative and non-increasing;
  }{
    that
    \begin{eqn}
      \llabel{eq5}
      \int_0^\delta\psi(s)\dd s
        \leq\biggl(\int_0^\delta m(s)\dd s\biggr)+\delta m(\delta)<\infty
      \dott
    \end{eqn}
  }
  \argument{
    \lref{eq5};
    \lref{eq4};
    the fact that $\delta\in(0,R)$ ; the fact that $\psi$ is definable and
    non-increasing;
  }[verbs=ep]{
    \lref{item2}\dott
  }
  We next prove \lref{item3}\dott
  \startnewargseq
  In the following let $\kappa\in(0,\infty)$, $\eta\in(0,R]$,
  let $\nu$ be a finite nonnegative Borel measure on $(0,R]$
  which satisfies for all $\rho\in(0,\eta]$ that
  \begin{eqn}
    \llabel{setting4}
    \nu((0,\rho])\leq\kappa\rho,
  \end{eqn}
  and let $\zeta\in(0,R)$, $B\in[0,\infty)$ satisfy
  \begin{eqn}
    \llabel{setting5}
   \textstyle \zeta=\min\{\delta,\eta\}
    \qquad\text{and}\qquad
    B=\sup_{s\in[\zeta,R]}m(s)
    \dott
  \end{eqn}
  \cfload[.]
  \startnewargseq
  \argument{
    \lref{item1};
    \lref{eq4};
    \lref{eq5};
    \lref{setting4};
    \lref{setting5};
    the fact that $\psi|_{(0,\zeta]}$ is non-increasing;
  }{
    that
    \begin{eqn}
      \llabel{eq6}
      \int_0^\zeta m(s)\,\nu(\mathrm{d}s)
        \leq\int_0^\zeta\psi(s)\,\nu(\mathrm{d}s)
        \leq\kappa\int_0^\zeta\psi(s)\dd s<\infty
      \dott
    \end{eqn}
  }
  \argument{
    the assumption that $m$ is locally bounded;
    the fact that $[\zeta,R]$ is a compact subset of $(0,R]$;
    \lref{setting5};
    the fact that $\nu$ is finite;
  }{
    that
    \begin{eqn}
      \llabel{eq7}
      \int_{(\zeta,R]}m(s)\dd\nu(s)
        \leq B\nu((0,R])<\infty
      \dott
    \end{eqn}
  }
  \argument{
    \lref{eq7};
    \lref{eq6};
    the arbitrariness of $\kappa$, $\eta$, and $\nu$;
  }[verbs=ep]{
    \lref{item3}\dott
  }
  We finally prove \lref{item4}\dott
  \startnewargseq
  In the following let $\theta\in(0,R)$,
  $\varrho,\mathfrak C\in(0,\infty)$, let
  $\varphi\colon(0,\theta]\to[0,\infty)$ be definable,
  non-increasing, and integrable and satisfy for all
  $s\in(0,\theta]$ that
  \begin{eqn}
    \llabel{setting6}
    m(s)\leq\varphi(s),
  \end{eqn}
  let $I\subseteq\R$ be an interval,
  let $r\colon I\to(0,\infty)$ be locally absolutely continuous,
  and assume for all $u\in(0,\varrho)$ it holds that
  \begin{eqn}
    \llabel{setting7}
    \int_{r^{-1}([u/2,u])}\abs{r'(t)}\dd t
      \leq\mathfrak C u,
  \end{eqn}
  let $\sigma\in(0,\theta]\cap(0,\varrho)$, and let
  $L\colon\mathcal B((0,\sigma])\to[0,\infty]$ satisfy
  for all $D\in\mathcal B((0,\sigma])$ that
  \begin{eqn}
    \llabel{setting8}
    L(D)=\int_{r^{-1}(D)}\abs{r'(t)}\dd t
    \dott
  \end{eqn}
  \cfload[.]
  \startnewargseq
  \argument{
    \lref{setting7};
    \lref{setting8};
    the fact that $r$ is locally absolutely continuous;
    \itref{lem:v}{item2,item3};
  }{
    that $L$ is a finite nonnegative Borel measure on
    $(0,\sigma]$ and that for all $\rho\in(0,\sigma]$ it holds that
    \begin{eqn}
      \llabel{eq8}
      L((0,\rho])\leq 2\mathfrak C\rho
      \dott
    \end{eqn}
  }
  \argument{
    \lref{setting8};
    the fact that $m$ is nonnegative and Borel measurable;
    the monotone convergence theorem;
  }{
    that
    \begin{eqn}
      \llabel{eq9}
      \int_{r^{-1}((0,\sigma])}m(r(t))\abs{r'(t)}\dd t
        =\int_{(0,\sigma]}m(s)\dd L(s)
      \dott
    \end{eqn}
  }
  \argument{
    \lref{eq9};
    \lref{item1};
    \lref{setting6};
    \lref{eq8};
    the fact that $\varphi|_{(0,\sigma]}$ is non-increasing and
    integrable;
  }{
    that
    \begin{eqn}[a]
      \llabel{eq10}
      \int_{r^{-1}((0,\sigma])}m(r(t))\abs{r'(t)}\dd t
        \leq\int_0^\sigma\varphi(s)\,L(\mathrm{d}s)
        \leq 2\mathfrak C\int_0^\sigma\varphi(s)\dd s<\infty
      \dott
    \end{eqn}
  }
  \argument{
    \lref{eq10};
    the arbitrariness of $\theta$, $\varphi$, $\varrho$,
    $\mathfrak C$, $I$, $r$, and $\sigma$;
  }[verbs=ep]{
    \lref{item4}\dott
  }
\end{aproof}

\subsection{Upper bounds for parts of the length of the radial component in the DE system}
\label{subsec:cone-annular-occupation}

\begin{athm}{proposition}{prop:cone-annular-length}
  Let $D\in\N$, $\mathfrak R,\mathsf a\in(0,\infty)$,
  $\mathsf b\in[\mathsf a,\infty)$, let $M\subseteq\R^D$ be a
  nonempty, compact, and definable\cfadd{definition:definable-set}
  \submanifold{2}{D}, let
  $\Omega\subseteq(0,\mathfrak R)\times M$ be definable and open in
  $(0,\mathfrak R)\times M$, and let
  $h=(h(r,x))_{(r,x)\in\Omega}\in \COne{\Omega}{\R}$ be
  definable\cfadd{definition:c1-on-manifold}\cfadd{definition:definable-function}\dott
  \cfload[.]
  Then there exists $\mathfrak C\in\R$ such that for every
  interval $I\subseteq\R$, all locally absolutely continuous
  $r\colon I\to(0,\mathfrak R)$ and $x\colon I\to M$, and
  all Lebesgue measurable
  $\alpha\colon I\to[0,\infty)$ and $A\colon I\to\R^{D\times D}$
  with $\midcup_{t\in I}\{(r(t),x(t))\}\subseteq\Omega$ and
  with the property that for Lebesgue-almost every $t\in I$ and
  every $v\in T_{x(t)}M$\cfadd{definition:tangent-space} it holds that
  \begin{equation}
    \llabel{ass1}
    A(t)v\in T_{x(t)}M,
    \qquad \langle A(t)v,v\rangle\geq\mathsf a\norm{v}^2,
    \qquad
    \norm{A(t)v}\leq\mathsf b\norm{v},
  \end{equation}
  \begin{eqn}
    \llabel{ass2}
    r'(t)=-\alpha(t)(r(t))^2
      (\partial_r h)(r(t),x(t)),
    \qquad\text{and}\qquad
    x'(t)=-\alpha(t)A(t)(\nabla_Mh)(r(t),x(t))
  \end{eqn}
  we have for all $\rho\in(0,\mathfrak R)$ that
  \begin{eqn}
    \llabel{concl1}
    \int_{r^{-1}([\rho/2,\rho])}
      \sqrt{\abs{r'(t)}^2+(r(t))^2\norm{x'(t)}^2}\dd t
      \leq \mathfrak C\rho
    \dott
  \end{eqn}
  \cfadd{definition:angular-gradient}\cfout[.]
\end{athm}
\begin{aproof}
  Throughout this proof let $K\subseteq\R^{D+1}$ and
  $\mathcal W\subseteq\R^2\times M$ satisfy
  \begin{eqn}
    \llabel{setting1}
    K=[1/2,1]\times M
    \qquad\text{and}\qquad
    \mathcal W=\{(\rho,s,y)\in(0,\mathfrak R)\times(0,\infty)\times M
      \colon(\rho s,y)\in\Omega\},
  \end{eqn}
  let $f\colon\mathcal W\to\R$ and
  $H\colon\mathcal W\to[0,\infty)$ satisfy that for all
  $(\rho,s,y)\in\mathcal W$ it holds that
  \begin{eqn}
    \llabel{setting2}
    f(\rho,s,y)=h(\rho s,y)
    \qquad\text{and}\qquad
    H(\rho,s,y)=\sqrt{
      \rho^2\bigl((\partial_r h)(\rho s,y)\bigr)^2
      +\norm{(\nabla_Mh)(\rho s,y)}^2},
  \end{eqn}
  let $\mathcal D\subseteq(0,\mathfrak R)\times K$ satisfy
  \begin{eqn}
    \llabel{setting3}
    \mathcal D=\mathcal W\cap((0,\mathfrak R)\times K),
  \end{eqn}
  let $\mathcal E\subseteq(0,\mathfrak R)\times\R$ satisfy
  \begin{eqn}
    \llabel{setting3.1}
    \mathcal E=\bigl\{(\rho,e)\in(0,\mathfrak R)\times\R\colon\bigl[
      \exists\,u\in K\colon((\rho,u)\in\mathcal D)\wedge(f(\rho,u)=e)
      \bigr]\bigr\},
  \end{eqn}
  let $\kappa,\Lambda,\mathbf B\in(0,\infty)$ satisfy
  $\kappa=\min\{1/4,\mathsf a\}$,
  $\Lambda=\max\{1,\mathsf b\}$, and
  \begin{eqn}
    \llabel{setting4}
    \mathbf B=\frac12+\sum_{k=1}^D
      \biggl(\sup_{y=(y_1,\dots,y_D)\in M}y_k
        -\inf_{y=(y_1,\dots,y_D)\in M}y_k\biggr),
  \end{eqn}
  let $\mu\colon\mathcal E\to[0,\infty)$ satisfy
   for all $(\rho,e)\in\mathcal E$ that
  \begin{eqn}
    \llabel{eq1}
    \mu(\rho,e)=\inf\bigl\{H(\rho,u)\colon\bigl[
      ((\rho,u)\in\mathcal D)\wedge(f(\rho,u)=e)\bigr]\bigr\},
  \end{eqn}
  and let $\mathcal V\subseteq\mathcal E$ satisfy
  \begin{eqn}
    \llabel{setting5}
    \mathcal V=\{(\rho,e)\in\mathcal E\colon\mu(\rho,e)=0\}
    \dott
  \end{eqn}
  \cfload[.]
  \startnewargseq
  \argument{
    \lref{setting1};
    the assumptions that $M$ is nonempty, compact, and definable;
    the assumptions that $\Omega$ is definable and open in $(0,\mathfrak R)\times M$;
  }{
    that $K$ is nonempty, compact, and definable and that
    $\mathcal W$ is definable and open in
    $(0,\mathfrak R)\times(0,\infty)\times M$\dott
    \llabel{eq0.1}
  }
  \argument{
    \lref{eq0.1};
    \lref{setting2};
    the assumption that $h\in \COne{\Omega}{\R}$ is definable;
  }{
    that $f\in \COne{\mathcal W}{\R}$ is definable\dott
  }
  \argument{
    the assumption that $h\in \COne{\Omega}{\R}$ is definable;
    the assumptions that $\Omega$ is definable and open in $(0,\mathfrak R)\times M$;
    \cref{lem:radial-derivative-well-defined};
    \cref{lem:angular-gradient-well-defined};
  }{
    that
    $\partial_r h\colon\Omega\to\R$ and
    $\nabla_Mh\colon\Omega\to\R^D$
    are continuous and definable\dott
    \llabel{arg0.2}
  }
  \argument{
    \lref{arg0.2};
    \lref{setting1};
    \lref{setting2};
  }{
    that $H$ is continuous and definable on $\mathcal W$\dott
  }
  \argument{
    \lref{setting3};
    \lref{setting3.1};
    the fact that $K$, $\mathcal W$, and $f$ are definable;
    \cite[Theorem~1.13]{Coste};
  }{
    that $\mathcal D$ and $\mathcal E$ are definable and that for all
    $(\rho,e)\in\mathcal E$ it holds that
    \begin{eqn}
      \llabel{eq0.2}
      \bigl\{H(\rho,u)\colon\bigl[
        ((\rho,u)\in\mathcal D)\wedge(f(\rho,u)=e)\bigr]\bigr\}
        \ne\emptyset
      \dott
    \end{eqn}
  }
  \argument{
    \lref{eq0.2};
    \lref{eq1};
    the fact that $H$ is finite and nonnegative;
    the fact that $\mathcal D$, $f$, and $H$ are definable;
    \cite[Theorem~1.13]{Coste};
  }{
    that $\mu$ is finite and definable\dott
  }
  \argument{
    \lref{eq1};
    \lref{eq0.2};
    the fact that $\mathcal D$, $f$, $H$, and $\mu$ are definable;
     \cite[Theorem~3.1]{Coste};
  }{
    that there exists a definable function
    $\psi=(\psi_0,\psi_1,\dots,\psi_D)\colon\mathcal E\to K$
    which satisfies that for all $(\rho,e)\in\mathcal E$ it holds that
    \begin{eqn}
      \llabel{eq2}
      (\rho,\psi(\rho,e))\in\mathcal D,
      \qquad f(\rho,\psi(\rho,e))=e,
      \qquad\text{and}\qquad \mu(\rho,e)\leq H(\rho,\psi(\rho,e)),
    \end{eqn}
    and which satisfies that for all $(\rho,e)\in\mathcal E$
    with $\mu(\rho,e)>0$ it holds that
    \begin{eqn}
      \llabel{eq2a}
      H(\rho,\psi(\rho,e))<2\mu(\rho,e)\dott
    \end{eqn}
  }
  \argument{
    \lref{eq1};
    \lref{setting5};
    the fact that $\mathcal D$, $f$, $H$, and $\mu$ are definable;
    \cite[Theorem~3.1]{Coste};
  }{
    that $\mathcal V$ is definable and that there exists a definable
    function
    $\chi=(\chi_0,\chi_1,\dots,\chi_D)\colon
      \{(\rho,\eta,e)\in(0,\mathfrak R)\times(0,\infty)\times\R
        \colon(\rho,e)\in\mathcal V\}\to K$
    which satisfies that for all
    $(\rho,\eta,e)\in(0,\mathfrak R)\times(0,\infty)\times\R$
    with $(\rho,e)\in\mathcal V$ it holds that
    \begin{eqn}
      \llabel{eq2.0}
      (\rho,\chi(\rho,\eta,e))\in\mathcal D,
      \qquad f(\rho,\chi(\rho,\eta,e))=e,
      \qquad\text{and}\qquad
      H(\rho,\chi(\rho,\eta,e))<\eta
      \dott
    \end{eqn}
  }
  \argument{
    \lref{eq2.0};
    the fact that $\chi$ is definable;
    \cite[Theorems~6.6 and~6.7]{Coste};
  }{
    that there exists $P\in\N$ which satisfies that for all
    $\rho\in(0,\mathfrak R)$, $\eta\in(0,\infty)$ 
    there exists a definable partition\cfadd{definition:definable-partition} $\mathscr P$ of
    $\{e\in\R\colon(\rho,e)\in\mathcal V\}$ such that
    \begin{enumerate}[label=(\Roman*)]
      \item it holds that $\#\mathscr P\leq P$ and for all
      $J\in\mathscr P$ it holds that $J$ is a singleton or an open
      interval,
      \item for all open intervals $J\in\mathscr P$ it holds that
      $\chi(\rho,\eta,\cdot)|_J\in C^1(J,K)$, and
      \item for all open intervals $J\in\mathscr P$ 
      there exists
      $\varsigma=(\varsigma_0,\varsigma_1,\dots,\varsigma_D)
        \in\{-1,1\}^{D+1}$ such that for all
      $j\in\{0,1,\dots,D\}$, $e\in J$ it holds that
      \begin{eqn}
        \varsigma_j\bigl(\partial_e \chi_j\bigr)\bigl(\rho,\eta,e\bigr)
          \geq0\dott
      \end{eqn}
    \end{enumerate}
    \llabel{arg2.0}
  }
  \argument{
    \lref{arg2.0};
    \lref{setting1};
    \lref{setting4};
    the fact that
    $\chi\bigl(\{(\rho,\eta,e)\in(0,\mathfrak R)\times(0,\infty)\times\R
      \colon(\rho,e)\in\mathcal V\}\bigr)\subseteq K$;
    the fundamental theorem of calculus;
    the monotone convergence theorem;
  }{
    that for every $\rho\in(0,\mathfrak R)$, $\eta\in(0,\infty)$,
    $a,b\in\R$ with $a<b$ and
    $\{\rho\}\times(a,b)\subseteq\mathcal V$ it holds that
    \begin{eqn}[a]
      \llabel{eq2.1}
      \int_a^b\bigl|\bigl(\partial_e \chi\bigr)
        \bigl(\rho,\eta,e\bigr)\bigr|\dd e
      &\leq\sum_{k=0}^D\int_a^b
        \bigl|\bigl(\partial_e \chi_k\bigr)\bigl(\rho,\eta,e\bigr)\bigr|\dd e
      \leq P\mathbf B
      \dott
    \end{eqn}
  }
  \argument{
    \lref{eq2.1};
    \lref{setting2};
    \lref{eq2.0};
    the chain rule;
    the Cauchy--Schwarz inequality;
  }{
    that for all $\rho\in(0,\mathfrak R)$, $\eta\in(0,\infty)$,
    $a,b\in\R$ with $a<b$ and
    $\{\rho\}\times(a,b)\subseteq\mathcal V$ it holds that for
    Lebesgue-almost all $e\in(a,b)$ it holds that
    \begin{eqn}
      \llabel{eq2.1.1}
      1=\partial_e\bigl(f(\rho,\chi(\rho,\eta,e))\bigr)
        \leq H(\rho,\chi(\rho,\eta,e))
          \bigl|\bigl(\partial_e \chi\bigr)\bigl(\rho,\eta,e\bigr)\bigr|
        \leq\eta\bigl|\bigl(\partial_e \chi\bigr)\bigl(\rho,\eta,e\bigr)\bigr|
      \dott
    \end{eqn}
  }
  \argument{
    \lref{eq2.1.1};
    \lref{eq2.1}
  }{
    that for all $\rho\in(0,\mathfrak R)$, $\eta\in(0,\infty)$,
    $a,b\in\R$ with $a<b$ and
    $\{\rho\}\times(a,b)\subseteq\mathcal V$ it holds that
    \begin{eqn}
      \llabel{eq2.1.2}
      b-a\leq\eta\int_a^b
        \bigl|\bigl(\partial_e \chi\bigr)\bigl(\rho,\eta,e\bigr)\bigr|\dd e
        \leq\eta P\mathbf B
      \dott
    \end{eqn}
  }
  \argument{
    \lref{eq2.1.2};
    the fact that $P\mathbf B\in(0,\infty)$;
    the fact that $\mathcal V$ is definable;
  }{
    that for all $\rho\in(0,\mathfrak R)$ it holds that
    \begin{eqn}
      \llabel{eq2.2}
      \#\{e\in\R\colon(\rho,e)\in\mathcal V\}<\infty
      \dott
    \end{eqn}
  }
  \argument{
    the fact that $\mathcal E$ is definable;
    the fact that for all $j\in\{0,1,\dots,D\}$ it holds that
    $\psi_j\colon\mathcal E\to\R$ is definable;
    \cite[Theorem~6.7]{Coste};
  }{
    that there exist definable partitions\cfadd{definition:definable-partition}
    $\mathscr T_j$, $j\in\{0,1,\dots,D\}$, of $\mathcal E$
    which satisfy that for all $j\in\{0,1,\dots,D\}$,
    $T\in\mathscr T_j$ it holds that $T$ is a $C^1$
    submanifold of $\R^2$\cfadd{definition:submanifold} and
    \begin{eqn}
      \llabel{eq2.2.0}
      \psi_j|_T\in \COne{T}{\R}\dott
    \end{eqn}
  }
  \argument{
    \lref{eq2.2.0};
    the fact that $\mathcal E$ and $\mathcal V$ are definable;
    \cite[Theorem~6.6]{Coste};
  }{
    that there exists a cylindrical definable $C^2$ cell decomposition
    $\mathscr P$ of $\R^2$ which is compatible\cfadd{definition:partition-compatibility} with
    $(0,\mathfrak R)\times\R$, $\mathcal E$, $\mathcal V$,
    and every element of $\mathscr T_j$ for every
    $j\in\{0,1,\dots,D\}$\cfadd{definition:planar-cell-decomposition}\dott
  }
  In the following let $\mathscr Q\subseteq\mathscr P$ be the set
  of all $Q\in\mathscr P$ with $Q\subseteq\mathcal E$
  such that for all $\rho\in(0,\mathfrak R)$ it holds that
  $\{e\in\R\colon(\rho,e)\in Q\}$ is empty or an open interval\dott
  \startnewargseq
  \argument{
    the fact that $\mathscr P$ is a definable partition\cfadd{definition:definable-partition} of $\R^2$;
    the fact that for all $Q\in\mathscr Q$ it holds that
    $Q\in\mathscr P$ and $Q\subseteq\mathcal E$;
    the fact that for all $j\in\{0,1,\dots,D\}$ it holds that
    $\mathscr T_j$ is a definable partition\cfadd{definition:definable-partition} of $\mathcal E$;
    the fact that for all $j\in\{0,1,\dots,D\}$,
    $T\in\mathscr T_j$ it holds that $T$ is a union of elements
    of $\mathscr P$;
  }{
    that for all $Q\in\mathscr Q$, $j\in\{0,1,\dots,D\}$
    it holds that
    \begin{eqn}
      \llabel{eq2.2.0.1}
      \exists\,T\in\mathscr T_j\colon Q\subseteq T
      \dott
    \end{eqn}
  }
  \argument{
    \lref{eq2.2.0.1};
    \lref{eq2.2.0};
    the fact that for all $Q\in\mathscr Q$ it holds that
    $Q\subseteq(0,\mathfrak R)\times\R$;
    the fact that for all $Q\in\mathscr Q$, $\rho\in(0,\mathfrak R)$
    it holds that $\{e\in\R\colon(\rho,e)\in Q\}$ is empty
    or an open interval;
    the fact that for all $\rho\in\R$ it holds that
    $\bigl(\R\ni e\mapsto(\rho,e)\in\R^2\bigr)\in C^1(\R,\R^2)$;
    the chain rule;
  }{
    that for all $Q\in\mathscr Q$, $\rho\in(0,\mathfrak R)$
    with $\{e\in\R\colon(\rho,e)\in Q\}\ne\emptyset$ it holds that
    \begin{eqn}[a]
      \llabel{eq2.2.0.2}
      \psi(\rho,\cdot)|_{\{e\in\R\colon(\rho,e)\in Q\}}\in C^1\bigl(\{e\in\R\colon(\rho,e)\in Q\},\R^{D+1}\bigr)\dott
    \end{eqn}
    \llabel{eq2.2.0.3}
  }
  \argument{
    \lref{eq2.2.0.3};
    the fact that for all $Q\in\mathscr Q$ it holds that $Q$ is definable;
    the fact that for all $j\in\{0,1,\dots,D\}$ it holds that
    $\psi_j\colon\mathcal E\to\R$ is definable;
    the fact that derivatives of definable functions are definable;
  }{
    that for all $Q\in\mathscr Q$, $j\in\{0,1,\dots,D\}$
    it holds that
    $\bigl(Q\ni(\rho,e)\mapsto
      \bigl(\partial_e \psi_j\bigr)\bigl(\rho,e\bigr)\in\R\bigr)$
    is definable\dott
  }
  In the following for all $Q\in\mathscr Q$,
  $j\in\{0,1,\dots,D\}$, $\tau\in\{-1,0,1\}$ let
  $\mathcal S_{Q,j}^{\tau}\subseteq Q$ satisfy
  \begin{eqn}
    \llabel{eq2.2.1}
    \mathcal S_{Q,j}^{\tau}
      =\bigl\{(\rho,e)\in Q\colon
        \operatorname{sign}\bigl(
          \bigl(\partial_e \psi_j\bigr)\bigl(\rho,e\bigr)\bigr)=\tau\bigr\}\dott
  \end{eqn}
  \startnewargseq
  \argument{
    \lref{eq2.2.1};
    the fact that for all $Q\in\mathscr Q$, $j\in\{0,1,\dots,D\}$
    it holds that
    $(\rho,e)\mapsto\bigl(\partial_e \psi_j\bigr)\bigl(\rho,e\bigr)$
    is definable on $Q$;
    the fact that $\mathscr P$ is finite;
    \cite[Theorem~6.6]{Coste};
  }{
    that there exists a cylindrical definable $C^2$ cell decomposition
    $\mathscr R$ of $\R^2$ which is a refinement\cfadd{definition:partition-refinement} of $\mathscr P$
    and which satisfies that for all $Q\in\mathscr Q$,
    $j\in\{0,1,\dots,D\}$, $\tau\in\{-1,0,1\}$ it holds that
    $\mathscr R$ is compatible\cfadd{definition:partition-compatibility}
    with $\mathcal S_{Q,j}^{\tau}$\dott
  }
  In the following let
  \begin{eqn}
    \llabel{setting5.0}
    \mathscr A=\bigl\{
      \bigl\{\rho\in\R\colon\bigl[\exists\,e\in\R\colon(\rho,e)\in U\bigr]\bigr\}
      \colon U\in\mathscr R\bigr\}
  \end{eqn}
  and let $L\in\N$, $\mathfrak C\in(0,\infty)$ satisfy
  \begin{eqn}
    \llabel{setting10}
    L=1+\#\mathscr R
    \qquad\text{and}\qquad
    \mathfrak C=\max\{1,2(\Lambda/\kappa)L\mathbf B\}\dott
  \end{eqn}
  \startnewargseq
  \argument{
    \lref{setting5.0};
    the fact that $\mathscr R$ is compatible\cfadd{definition:partition-compatibility} with $(0,\mathfrak R)\times\R$;
  }{
    that $\mathscr A$ is a definable partition\cfadd{definition:definable-partition} of $\R$
    and that for all $J\in\mathscr A$ it holds that $J$ is a singleton
    or an open interval and that either $J\subseteq(0,\mathfrak R)$ or
    $J\cap(0,\mathfrak R)=\emptyset$\dott
  }
  We claim that for every $J\in\mathscr A$ with
  $J\subseteq(0,\mathfrak R)$,
  every interval $I\subseteq\R$,
  all locally absolutely continuous $r\colon I\to(0,\mathfrak R)$
  and $x\colon I\to M$,
  and all Lebesgue measurable
  $\alpha\colon I\to[0,\infty)$,
  $A\colon I\to\R^{D\times D}$,
  and $\mathfrak v\colon I\to[0,\infty)$
  with $\midcup_{t\in I}\{(r(t),x(t))\}\subseteq\Omega$,
  with the property that for Lebesgue-almost all $t\in I$
  and every $v\in T_{x(t)}M$ it holds that
  \begin{eqn}
    A(t)v\in T_{x(t)}M,
    \qquad \langle A(t)v,v\rangle\geq\mathsf a\norm{v}^2,
    \qquad \text{and}\qquad
    \norm{A(t)v}\leq\mathsf b\norm{v},
  \end{eqn}
  with the property that for Lebesgue-almost all $t\in I$
  it holds that 
  \begin{eqn}
    r'(t)=-\alpha(t)(r(t))^2
      (\partial_r h)(r(t),x(t))
    \qquad\text{and}\qquad
    x'(t)=-\alpha(t)A(t)(\nabla_Mh)(r(t),x(t)),
  \end{eqn}
  and with the property that for Lebesgue-almost all $t\in I$
  it holds that
  \begin{eqn}
    \mathfrak v(t)
      =\sqrt{\abs{r'(t)}^2+(r(t))^2\norm{x'(t)}^2}{}
  \end{eqn}
  we have for all $\rho\in J$, $\tau,\theta\in I$
  with $\tau<\theta$ that
  \begin{eqn}
    \llabel{eq8.1}
    \int_{[\tau,\theta]\cap r^{-1}([\rho/2,\rho])}
      \mathfrak v(t)\dd t
      \leq\mathfrak C\rho
    \dott
  \end{eqn}
  \startnewargseq
  To prove \lref{eq8.1}, let $J\in\mathscr A$ satisfy
  $J\subseteq(0,\mathfrak R)$,
  let $I\subseteq\R$ be an interval,
  let $r\colon I\to(0,\mathfrak R)$ and $x\colon I\to M$
  be locally absolutely continuous,
  let $\rho\in J$, $\tau,\theta\in I$ satisfy $\tau<\theta$,
  and let $\alpha\colon I\to[0,\infty)$,
  $A\colon I\to\R^{D\times D}$, and
  $\mathfrak v\colon I\to[0,\infty)$ be Lebesgue measurable,
  assume that $\midcup_{t\in I}\{(r(t),x(t))\}\subseteq\Omega$,
  assume that for Lebesgue-almost all $t\in I$ and every
  $v\in T_{x(t)}M$ it holds that
  \begin{eqn}
    \llabel{setting10.1}
    A(t)v\in T_{x(t)}M,
    \qquad \langle A(t)v,v\rangle\geq\mathsf a\norm{v}^2,
    \qquad \text{and}\qquad
    \norm{A(t)v}\leq\mathsf b\norm{v},
  \end{eqn}
  assume that for Lebesgue-almost all $t\in I$ it holds that
  \begin{eqn}
    \llabel{setting10.2}
    r'(t)=-\alpha(t)(r(t))^2
      (\partial_r h)(r(t),x(t))
    \qquad\text{and}\qquad
    x'(t)=-\alpha(t)A(t)(\nabla_Mh)(r(t),x(t)),
  \end{eqn}
  and assume that for Lebesgue-almost all $t\in I$ it holds that
  \begin{eqn}
    \llabel{setting10.3}
    \mathfrak v(t)
      =\sqrt{\abs{r'(t)}^2+(r(t))^2\norm{x'(t)}^2}
    \dott
  \end{eqn}
  \startnewargseq
  \argument{
    \lref{eq2.2};
    \lref{eq2.2.0.2};
    \lref{eq2.2.1};
    \lref{setting5.0};
    \lref{setting10};
    the fact that $\mathscr R$ is
    compatible\cfadd{definition:partition-compatibility} with $\mathcal E$ and $\mathcal V$;
    the fact that for all $Q\in\mathscr Q$, $j\in\{0,1,\dots,D\}$,
    $\omega\in\{-1,0,1\}$ it holds that $\mathscr R$ is compatible\cfadd{definition:partition-compatibility}
    with $\mathcal S_{Q,j}^{\omega}$;
  }{
    that there exist $N\in\{1,\dots,L\}$, a set
    $\mathcal L\subseteq\{0,1,\dots,N-1\}$, definable functions
    $\mathfrak e_j\colon J\to\R$,
    $j\in\{1,\dots,N-1\}$, functions
    $\mathfrak e_0,\mathfrak e_N\colon J\to\R\cup\{-\infty,\infty\}$, and
    $\sigma_{ik}\in\{-1,0,1\}$, $i\in\mathcal L$,
    $k\in\{0,1,\dots,D\}$, which satisfy
    \begin{enumerate}[label=\alph*)]
      \item for all $\rho\in J$ it holds that
      \begin{gather}\llabel{eq3}
        -\infty=\mathfrak e_0(\rho)
        <\cdots<\mathfrak e_N(\rho)=\infty,\notag\\
        \{e\in\R\colon(\rho,e)\in\mathcal V\}
          \subseteq\{\mathfrak e_j(\rho)\colon j\in\{1,\dots,N-1\}\},\\
        \text{and}\qquad\{e\in\R\colon(\rho,e)\in\mathcal E\}
          \backslash\{\mathfrak e_j(\rho)\colon j\in\{1,\dots,N-1\}\}
          =\midcup_{i\in\mathcal L}
            (\mathfrak e_i(\rho),\mathfrak e_{i+1}(\rho))\notag,
      \end{gather}
      \item for all $j\in\{1,\dots,N-1\}$, $\rho\in J$
      it holds that $\mathfrak e_j(\rho)\in\R$,
      \item for all $i\in\mathcal L$, $\rho\in J$ it holds that
      $\psi(\rho,\cdot)|_{(\mathfrak e_i(\rho),\mathfrak e_{i+1}(\rho))}
        \in C^1((\mathfrak e_i(\rho),\mathfrak e_{i+1}(\rho)),K)$, and
      \item for all $i\in\mathcal L$, $k\in\{0,1,\dots,D\}$,
      $\rho\in J$, $e\in(\mathfrak e_i(\rho),\mathfrak e_{i+1}(\rho))$
      it holds that
      \begin{eqn}
        \llabel{eq4}
        (\rho,e)\in\mathcal E\backslash\mathcal V,
        \qquad\mu(\rho,e)>0,
        \qquad\text{and}\qquad
        \operatorname{sign}\bigl(\bigl(\partial_e \psi_k\bigr)
          \bigl(\rho,e\bigr)\bigr)=\sigma_{ik}\dott
      \end{eqn}
    \end{enumerate}
  }
  In the following let $\mathcal O_i\subseteq J\times\R$,
  $i\in\mathcal L$, satisfy
  \begin{eqn}
    \llabel{setting6}
    \mathcal O_i=\{(\rho,e)\in J\times\R
      \colon\mathfrak e_i(\rho)<e<\mathfrak e_{i+1}(\rho)\},
  \end{eqn}
  and let $z_i\colon\mathcal O_i\to\R$, $i\in\mathcal L$,
  satisfy that for all $i\in\mathcal L$, $(\rho,e)\in\mathcal O_i$
  it holds that
  \begin{eqn}
    \llabel{setting7}
    \textstyle z_i(\rho,e)=2\sum_{k=0}^D\sigma_{ik}\psi_k(\rho,e)
    \dott
  \end{eqn}
  \cfload[.]
  \startnewargseq
  \argument{
    \lref{setting2};
    \lref{eq2};
    \lref{eq2a};
    \lref{eq4};
    \lref{setting7};
    the chain rule;
    the fact that for all $i\in\mathcal L$, $(\rho,e)\in\mathcal O_i$
    it holds that
    $\bigl(\bigl(\partial_e\psi_j\bigr)\bigl(\rho,e\bigr)\bigr)_{j=1}^D
      \in T_{(\psi_j(\rho,e))_{j=1}^D}M$%
    \cfadd{definition:tangent-space};
    the Cauchy--Schwarz inequality;
  }{
    that for all $i\in\mathcal L$, $(\rho,e)\in\mathcal O_i$,
    $(s,y)\in K$, $(\nu_1,\dots,\nu_D)\in\R^D$ with
    $\psi(\rho,e)=(s,y)$ and
    $(\nu_1,\dots,\nu_D)=(\nabla_Mh)(\rho s,y)$ it holds that
    \begin{eqn}[a]
      \llabel{eq5}
      1
        &\textstyle=\partial_e\bigl(f(\rho,\psi(\rho,e))\bigr)
      =\rho\,(\partial_r h)(\rho s,y)
            \bigl(\partial_e \psi_0\bigr)\bigl(\rho,e\bigr)
          +\sum_{k=1}^D
            \nu_k
            \bigl(\partial_e \psi_k\bigr)\bigl(\rho,e\bigr)
      \\\textstyle
        &\textstyle\leq H(\rho,\psi(\rho,e))
             \norm{\bigl(\partial_e \psi\bigr)\bigl(\rho,e\bigr)}
          \leq2\mu(\rho,e)\sum_{k=0}^D
          \abs{\bigl(\partial_e \psi_k\bigr)\bigl(\rho,e\bigr)}
          =\mu(\rho,e)\bigl(\partial_e z_i\bigr)\bigl(\rho,e\bigr)
      \dott
    \end{eqn}
  }
  \argument{
    \lref{eq5};
    \lref{setting1};
    \lref{setting4};
    \lref{eq4};
    \lref{setting7};
    the fact that $\psi$ takes values in $K$;
    the mean value theorem;
  }{
    that for all $i\in\mathcal L$, $\rho\in J$,
    $\xi,\zeta\in(\mathfrak e_i(\rho),\mathfrak e_{i+1}(\rho))$
    with $\xi<\zeta$ it holds that
    \begin{eqn}[a]
      \llabel{eq6}
     \textstyle 0<z_i(\rho,\zeta)-z_i(\rho,\xi)
        =2\sum_{k=0}^D
          \abs{\psi_k(\rho,\zeta)-\psi_k(\rho,\xi)}
          \leq 2\mathbf{B}
      \dott
    \end{eqn}
  }
  \argument{
    \lref{eq6};
    \lref{setting7};
  }{
    that for all $i\in\mathcal L$, $\rho\in J$ it holds that
    \begin{eqn}
      \llabel{eq7}
      0<\lim_{e\nearrow\mathfrak e_{i+1}(\rho)}z_i(\rho,e)
          -\lim_{e\searrow\mathfrak e_i(\rho)}z_i(\rho,e)
        \leq 2\mathbf{B}
      \dott
    \end{eqn}
  }
  In the following let $\beta_i^-,\beta_i^+,d_i\colon
  J\to\R$, $i\in\mathcal L$, satisfy that for all
  $i\in\mathcal L$, $\rho\in J$ it holds that
  \begin{eqn}
    \llabel{setting8}
    \beta_i^-(\rho)=\lim_{e\searrow\mathfrak e_i(\rho)}z_i(\rho,e),
    \qquad
    \beta_i^+(\rho)=\lim_{e\nearrow\mathfrak e_{i+1}(\rho)}z_i(\rho,e),
    \qquad \text{and}\qquad
    d_i(\rho)=\beta_i^+(\rho)-\beta_i^-(\rho)
  \end{eqn}
  and let $q_i\colon J\times\R\to\R$,
  $i\in\mathcal L$, satisfy that for all $i\in\mathcal L$,
  $(\rho,e)\in J\times\R$ it holds that
  \begin{eqn}
    \llabel{setting9}
    q_i(\rho,e)=
    \begin{cases}
      0 & \colon e\leq\mathfrak e_i(\rho),\\
      z_i(\rho,e)-\beta_i^-(\rho)
        & \colon\mathfrak e_i(\rho)<e<\mathfrak e_{i+1}(\rho),\\
      d_i(\rho) & \colon e\geq\mathfrak e_{i+1}(\rho)
    \end{cases}
    \dott
  \end{eqn}
  \cfload[.]
  \startnewargseq
  \argument{
    \lref{eq1};
    \lref{eq2};
    \lref{eq4};
    \lref{eq5};
    \lref{eq6};
    \lref{eq7};
    \lref{setting8};
    \lref{setting9};
  }{
    that for all $i\in\mathcal L$, $\rho\in J$
    it holds that
    $q_i(\rho,\cdot)|_{(\mathfrak e_i(\rho),\mathfrak e_{i+1}(\rho))}
      \in C^1((\mathfrak e_i(\rho),\mathfrak e_{i+1}(\rho)),\R)$
    and that $q_i(\rho,\cdot)$ is continuous and non-decreasing
    on $\R$, is strictly increasing on
    $(\mathfrak e_i(\rho),\mathfrak e_{i+1}(\rho))$, satisfies
    \begin{eqn}
      \{q_i(\rho,e)\colon e\in\R\}\subseteq[0,d_i(\rho)]
      \qquad\text{and}\qquad
      \{q_i(\rho,e)\colon
        e\in(\mathfrak e_i(\rho),\mathfrak e_{i+1}(\rho))\}
        =(0,d_i(\rho)),
    \end{eqn}
    and satisfies that for all $u\in K$ with
    $(\rho,u)\in\mathcal D$ and $(\rho,f(\rho,u))\in\mathcal O_i$
    it holds that
    \begin{eqn}
      \llabel{eq8}
      \bigl(\partial_e q_i\bigr)\bigl(\rho,f(\rho,u)\bigr)
        \geq\frac1{\mu(\rho,f(\rho,u))}
        \geq\frac1{H(\rho,u)}>0
      \dott
    \end{eqn}
  }
  In the following let $E,w\colon I\to\R$ and
  $\mathbf{g}\colon I\to\R^D$ satisfy for all $t\in I$ that
  \begin{eqn}
    \llabel{setting11}
    E(t)=h(r(t),x(t)),
    \qquad w(t)=(\partial_r h)(r(t),x(t)),
    \qquad \text{and}\qquad \mathbf{g}(t)=(\nabla_Mh)(r(t),x(t)),
  \end{eqn}
  let $\delta\in(0,\infty)$ satisfy
  \begin{eqn}
    \llabel{setting12}
    \delta=\frac12\min\bigl(\{1\}\cup\{d_i(\rho)\colon i\in\mathcal L\}\bigr),
  \end{eqn}
  let $\mathbf{q}_i^\varepsilon\colon\R\to\R$,
  $i\in\mathcal L$, $\varepsilon\in(0,\delta)$, satisfy for all
  $i\in\mathcal L$, $\varepsilon\in(0,\delta)$, $e\in\R$ that
  \begin{eqn}
    \llabel{setting13}
    \mathbf{q}_i^\varepsilon(e)
      =\min\{d_i(\rho)-\varepsilon,
          \max\{\varepsilon,q_i(\rho,e)\}\},
  \end{eqn}
  let $\mathcal T_i^\varepsilon\subseteq[\tau,\theta]$,
  $i\in\mathcal L$, $\varepsilon\in(0,\delta)$, satisfy for all
  $i\in\mathcal L$, $\varepsilon\in(0,\delta)$ that
  \begin{eqn}
    \llabel{setting14}
    \mathcal T_i^\varepsilon
      =\bigl\{t\in[\tau,\theta]\colon\bigl[
        (\rho/2\leq r(t)\leq\rho)\wedge
        (\varepsilon<q_i(\rho,E(t))<d_i(\rho)-\varepsilon)\bigr]\bigr\},
  \end{eqn}
  and let $\mathcal U^\varepsilon\subseteq[\tau,\theta]$,
  $\varepsilon\in(0,\delta)$, satisfy for all
  $\varepsilon\in(0,\delta)$ that
  \begin{eqn}
    \llabel{setting15}
    \mathcal U^\varepsilon=\midcup_{i\in\mathcal L}\mathcal T_i^\varepsilon
    \dott
  \end{eqn}
  \cfload[.]
  \startnewargseq
  \argument{
    \lref{setting10.1};
    \lref{setting10.2};
    \lref{setting11};
    the fact that for all $t\in I$ it holds that $\mathbf{g}(t)\in T_{x(t)}M$;
    the fact that for all compact subintervals $[a,b]\subseteq I$
    it holds that $\{(r(t),x(t))\colon t\in[a,b]\}$ is a compact
    subset of $\Omega$;
    the assumption that $h\in \COne{\Omega}{\R}$;
    the chain rule;
  }{
    that $E$ is locally absolutely continuous and non-increasing on $I$,
    and satisfies that for
    Lebesgue-almost all $t\in I$ it holds that
    \begin{eqn}[a]
      \llabel{eq9}
      -E'(t)
        &=-(\partial_r h)(r(t),x(t))r'(t)
          -\langle(\nabla_Mh)(r(t),x(t)),x'(t)\rangle
      \\
        &=\alpha(t)\bigl((r(t))^2(w(t))^2+\langle A(t)\mathbf{g}(t),\mathbf{g}(t)\rangle\bigr)
      \geq\alpha(t)\bigl((r(t))^2(w(t))^2+\mathsf a\norm{\mathbf{g}(t)}^2\bigr)
          \geq0
      \dott
    \end{eqn}
  }
  \argument{
    \lref{setting2};
    \lref{setting11};
    the fact that $\kappa=\min\{1/4,\mathsf a\}$;
  }{
    that for all $t\in[\tau,\theta]$ with
    $\rho/2\leq r(t)\leq\rho$ it holds that
    \begin{eqn}[a]
      \llabel{eq9.1}
      (r(t))^2(w(t))^2+\norm{\mathbf{g}(t)}^2
        &\leq\bigl(H(\rho,r(t)/\rho,x(t))\bigr)^2
          =\rho^2(w(t))^2+\norm{\mathbf{g}(t)}^2
      \\
        &\leq4(r(t))^2(w(t))^2+\norm{\mathbf{g}(t)}^2
      \leq\frac{1}{\kappa}
          \bigl((r(t))^2(w(t))^2+\mathsf a\norm{\mathbf{g}(t)}^2\bigr)
      \dott
    \end{eqn}
  }
  \argument{
    \lref{eq9};
    \lref{eq9.1};
    the fact that for all $t\in I$ it holds that $\alpha(t)\geq0$;
  }{
    that for Lebesgue-almost all $t\in[\tau,\theta]$ with
    $\rho/2\leq r(t)\leq\rho$ it holds that
    \begin{eqn}[a]
      \llabel{eq9.2}
      -E'(t)
        &\geq\alpha(t)
          \bigl((r(t))^2(w(t))^2+\mathsf a\norm{\mathbf{g}(t)}^2\bigr)
      \geq\kappa\alpha(t)\bigl(H(\rho,r(t)/\rho,x(t))\bigr)^2
      \dott
    \end{eqn}
  }
  \argument{
    \lref{setting10.1};
    \lref{setting10.2};
    \lref{setting10.3};
    \lref{setting11};
    \lref{eq9.1};
    the fact that for all $t\in I$ it holds that
    $\mathbf{g}(t)\in T_{x(t)}M$ and $\alpha(t)\geq0$;
    the fact that $\Lambda=\max\{1,\mathsf b\}$;
  }{
    that for Lebesgue-almost all $t\in[\tau,\theta]$ with
    $\rho/2\leq r(t)\leq\rho$ it holds that
    \begin{eqn}[a]
      \llabel{eq10}
      \mathfrak v(t)
        &=\alpha(t)\sqrt{(r(t))^4(w(t))^2+(r(t))^2\norm{A(t)\mathbf{g}(t)}^2}
        =\alpha(t)r(t)\sqrt{(r(t))^2(w(t))^2+\norm{A(t)\mathbf{g}(t)}^2}
      \\
        &\leq\alpha(t)r(t)\sqrt{(r(t))^2(w(t))^2+(\mathsf b)^2\norm{\mathbf{g}(t)}^2}
        \leq\Lambda\alpha(t)r(t)\sqrt{(r(t))^2(w(t))^2+\norm{\mathbf{g}(t)}^2}
      \\
        &\leq\Lambda\alpha(t)\rho H(\rho,r(t)/\rho,x(t))
      \dott
    \end{eqn}
  }
  \argument{
    the fact that $E|_{[\tau,\theta]}$ is absolutely continuous;
    \cref{lem:AC-level-set} (applied for every $j\in\{1,\dots,N-1\}$
    with $I\curvearrowleft[\tau,\theta]$,
    $c\curvearrowleft\mathfrak e_j(\rho)$,
    $f\curvearrowleft E|_{[\tau,\theta]}$ in the notation of
    \cref{lem:AC-level-set});
  }{
    that for all $j\in\{1,\dots,N-1\}$ it holds that
    \begin{eqn}
      \llabel{eq10.0}
      \int_{\{t\in[\tau,\theta]\colon E(t)=\mathfrak e_j(\rho)\}}
        \abs{E'(t)}\dd t=0
      \dott
    \end{eqn}
  }
  \argument{
    \lref{eq9};
    \lref{eq10.0};
  }{
    that for all $j\in\{1,\dots,N-1\}$ it holds that there exists a Lebesgue
    measurable set $\mathcal Z_j\subseteq[\tau,\theta]$ of Lebesgue
    measure zero such that for all
    $t\in[\tau,\theta]\backslash\mathcal Z_j$ with
    $E(t)=\mathfrak e_j(\rho)$ it holds that
    \begin{eqn}
      \llabel{eq10.1}
      0=-E'(t)
        \geq\alpha(t)
          \bigl((r(t))^2(w(t))^2+\mathsf a\norm{\mathbf{g}(t)}^2\bigr)
          \geq0
      \dott
    \end{eqn}
  }
  \argument{
    \lref{setting10.2};
    \lref{setting11};
    \lref{eq10.1};
    the fact that $\mathsf a>0$;
    the fact that for all $t\in I$ it holds that
    $r(t)>0$ and $\alpha(t)\geq 0$;
  }{
    that for all $j\in\{1,\dots,N-1\}$ it holds that there exists a Lebesgue
    measurable set $\mathcal Z_j\subseteq[\tau,\theta]$ of Lebesgue
    measure zero such that for all
    $t\in[\tau,\theta]\backslash\mathcal Z_j$ with
    $E(t)=\mathfrak e_j(\rho)$ it holds that $r$ and $x$
    are differentiable at $t$ and that $\alpha(t)w(t)=0$,
    \begin{eqn}[a]
      \llabel{eq10.2}
        \qquad \alpha(t)\mathbf{g}(t)=0,\qquad r'(t)=-(r(t))^2\alpha(t)w(t)=0,\qquad\text{and}\qquad
   x'(t)=-A(t)\bigl(\alpha(t)\mathbf{g}(t)\bigr)=0
      \dott
    \end{eqn}
  }
  \argument{
    \lref{setting10.3};
    \lref{eq10.2};
    the fact that $E|_{[\tau,\theta]}$ is continuous ;the fact that
    $\mathfrak v$ is Lebesgue measurable;
  }{
    that for all $j\in\{1,\dots,N-1\}$ it holds that
    \begin{eqn}
      \llabel{eq11}
      \int_{\{t\in[\tau,\theta]\colon E(t)=\mathfrak e_j(\rho)\}}
        \mathfrak v(t)\dd t=0
      \dott
    \end{eqn}
  }
  \argument{
    \lref{eq7};
    \lref{setting8};
    \lref{setting9};
    \lref{setting12};
    \lref{setting13};
  }{
    that for all $i\in\mathcal L$, $\varepsilon\in(0,\delta)$
    it holds that
    $\{e\in\R\colon\varepsilon\leq q_i(\rho,e)
      \leq d_i(\rho)-\varepsilon\}\subseteq (\mathfrak e_i(\rho),\mathfrak e_{i+1}(\rho))$ and is a compact interval 
    and that
    $\mathbf{q}_i^\varepsilon$ is globally Lipschitz and non-decreasing\dott
    \llabel{arg2}
  }
  \argument{
    \lref{arg2};
    \lref{setting13};
    \lref{eq9};
  }{
    that for all $i\in\mathcal L$, $\varepsilon\in(0,\delta)$
    it holds that $\mathbf{q}_i^\varepsilon\circ E$ is absolutely continuous
    and non-increasing on $[\tau,\theta]$ and that for
    Lebesgue-almost all $t\in[\tau,\theta]$ it holds that
    \begin{eqn}
      \llabel{eq12}
      -\frac{\ddd}{\ddd t}\mathbf{q}_i^\varepsilon(E(t))\geq0
      \dott
    \end{eqn}
  }
  \argument{
    \lref{setting1};
    \lref{setting2};
    \lref{setting3};
    \lref{setting6};
    \lref{setting9};
    \lref{setting11};
    \lref{setting14};
    the assumption that $\midcup_{t\in I}\{(r(t),x(t))\}\subseteq\Omega$;
  }{
    that for all $i\in\mathcal L$, $\varepsilon\in(0,\delta)$,
    $t\in\mathcal T_i^\varepsilon$ it holds that
    \begin{eqn}[a]
      \llabel{eq12.1}
      (\rho,r(t)/\rho,x(t))\in\mathcal D,
        \quad f(\rho,r(t)/\rho,x(t))=E(t),\quad
      \mathfrak e_i(\rho)<E(t)<\mathfrak e_{i+1}(\rho),
        \quad \text{and}\quad(\rho,E(t))\in\mathcal O_i
      \dott
    \end{eqn}
  }
  \argument{
    \lref{eq8};
    \lref{eq12.1};
  }{
    that for all $i\in\mathcal L$, $\varepsilon\in(0,\delta)$,
    $t\in\mathcal T_i^\varepsilon$ it holds that
    \begin{eqn}
      \llabel{eq12.2}
      \bigl(\partial_e q_i\bigr)\bigl(\rho,E(t)\bigr)
        \geq\frac1{\mu(\rho,E(t))}
        \geq\frac1{H(\rho,r(t)/\rho,x(t))}>0
      \dott
    \end{eqn}
  }
  \argument{
    \lref{setting13};
    the fact that for all $i\in\mathcal L$ it holds that
    $q_i(\rho,\cdot)|_{(\mathfrak e_i(\rho),\mathfrak e_{i+1}(\rho))}
      \in C^1((\mathfrak e_i(\rho),\mathfrak e_{i+1}(\rho)),\R)$;
  }{
    that for all $i\in\mathcal L$, $\varepsilon\in(0,\delta)$,
    $e\in(\mathfrak e_i(\rho),\mathfrak e_{i+1}(\rho))$ with
    $\varepsilon<q_i(\rho,e)<d_i(\rho)-\varepsilon$ it holds that
    $\mathbf{q}_i^\varepsilon$ is differentiable at $e$ and satisfies
    \begin{eqn}
      \llabel{eq12.3}
      \mathbf{q}_i^\varepsilon(e)=q_i(\rho,e)
      \qquad\text{and}\qquad
      (\mathbf{q}_i^\varepsilon)'(e)
        =\bigl(\partial_e q_i\bigr)\bigl(\rho,e\bigr)
      \dott
    \end{eqn}
  }
  \argument{
    \lref{setting14};
    \lref{eq9};
    \lref{eq9.2};
    \lref{eq10};
    \lref{eq12.1};
    \lref{eq12.2};
    \lref{eq12.3};
    the chain rule;
  }{
    that for all $i\in\mathcal L$, $\varepsilon\in(0,\delta)$ it holds that
    for Lebesgue-almost all $t\in\mathcal T_i^\varepsilon$
    it holds that $E$ and $\mathbf{q}_i^\varepsilon\circ E$ are
    differentiable at $t$ and that
    \begin{eqn}[a]
      \llabel{eq13}
      -\frac{\ddd}{\ddd t}\mathbf{q}_i^\varepsilon(E(t))
        =\bigl(\partial_e q_i\bigr)\bigl(\rho,E(t)\bigr)(-E'(t))
      \geq\frac{-E'(t)}{H(\rho,r(t)/\rho,x(t))}
          \geq\kappa\alpha(t)H(\rho,r(t)/\rho,x(t))
      \geq\frac{\kappa}{\Lambda\rho}\mathfrak v(t)
      \dott
    \end{eqn}
  }
  \argument{
    \lref{eq13};
    \lref{setting14};
    \lref{setting15};
    \lref{eq12};
    the fact that for all $\varepsilon\in(0,\delta)$ it holds that
    $\mathcal T_i^\varepsilon$, $i\in\mathcal L$, are pairwise disjoint;
    the fundamental theorem of calculus;
  }{
    that for all $\varepsilon\in(0,\delta)$ it holds that
    \begin{eqn}[a]
      \llabel{eq14}
      \int_{\mathcal U^\varepsilon}\mathfrak v(t)\dd t
        &\leq\frac{\Lambda\rho}{\kappa}
          \sum_{i\in\mathcal L}\int_\tau^\theta
            -\frac{\ddd}{\ddd t}\mathbf{q}_i^\varepsilon(E(t))\dd t
     =\frac{\Lambda\rho}{\kappa}\sum_{i\in\mathcal L}
          \bigl(\mathbf{q}_i^\varepsilon(E(\tau))-\mathbf{q}_i^\varepsilon(E(\theta))\bigr)
      \dott
    \end{eqn}
  }
  \argument{
    \lref{eq14};\lref{setting10};
    \lref{eq7}; 
    \lref{setting8};
    \lref{setting13};
  }{
    that for all $\varepsilon\in(0,\delta)$ it holds that
    \begin{eqn}
      \llabel{eq15}
      \int_{\mathcal U^\varepsilon}\mathfrak v(t)\dd t
        \leq\frac{\Lambda\rho}{\kappa}\sum_{i\in\mathcal L}d_i(\rho)
        \leq\frac{2\Lambda}{\kappa}N\mathbf{B}\rho
        \leq\frac{2\Lambda}{\kappa}L\mathbf{B}\rho
        \leq \mathfrak C\rho
      \dott
    \end{eqn}
  }
  \argument{
    \lref{setting9};
    \lref{setting14};
    \lref{setting15};
  }{
    that for all $v,w\in (0,\delta)$ with $w<v$ it holds that $\mathcal U^v\subseteq \mathcal U^w$ and that
    \begin{eqn}
      \llabel{eq16}
      \midcup_{\varepsilon\in(0,\delta)}\mathcal U^\varepsilon
        =\bigl\{t\in[\tau,\theta]\colon\bigl[(\rho/2\leq r(t)\leq\rho)\wedge
          \bigl(E(t)\in\midcup_{i\in\mathcal L}
            (\mathfrak e_i(\rho),\mathfrak e_{i+1}(\rho))\bigr)\bigr]\bigr\}
      \dott
    \end{eqn}
  }
  \argument{
    \lref{eq16};
    \lref{eq3};
    \lref{eq11};
    \lref{eq15};
    the monotone convergence theorem;
    the fact that for all $t\in[\tau,\theta]$ with
    $\rho/2\leq r(t)\leq\rho$ it holds that
    $(\rho,r(t)/\rho,x(t))\in\mathcal D$ and
    $(\rho,E(t))\in\mathcal E$;
  }{
    that
    \begin{eqn}
      \llabel{eq17}
      \int_{[\tau,\theta]\cap r^{-1}([\rho/2,\rho])}
        \mathfrak v(t)\dd t
        \leq \mathfrak C\rho
      \dott
    \end{eqn}
  }
  \argument{
    \lref{eq17};
    the arbitrariness of $J$, $I$,
    $r$, $x$, $\alpha$, $A$, $\mathfrak v$, $\rho$, $\tau$, and $\theta$;
  }[verbs=ep]{
    \lref{eq8.1}\dott
  }

  \startnewargseq
  \argument{
    \lref{eq8.1};
    the fact that $\mathscr A$ is compatible\cfadd{definition:partition-compatibility} with $(0,\mathfrak R)$;
    the monotone convergence theorem;
  }[verbs=ep]{
    \lref{concl1}\dott
  }
\end{aproof}

\subsection{Finite length of the angular component in the DE system}
\label{subsec:product-length}

\begin{athm}{corollary}{cor:angular-length}
  Let $D\in\N$, $\mathfrak R,\mathsf a\in(0,\infty)$,
  $\mathsf b\in[\mathsf a,\infty)$, let $M\subseteq\R^D$ be a
  nonempty, compact, and definable\cfadd{definition:definable-set}
  \submanifold{2}{D}, let
  $\mathcal O\subseteq(0,\mathfrak R)\times M$ be definable and open
  in $(0,\mathfrak R)\times M$, and let
  $h=(h(r,x))_{(r,x)\in\mathcal O}\in \COne{\mathcal O}{\R}$
  be definable\cfadd{definition:c1-on-manifold}%
  \cfadd{definition:definable-function}\dott \cfload[.]
  Then there exist $\delta\in(0,\mathfrak R)$,
  $\mathfrak C\in\R$, and a definable, non-increasing, and integrable
  $\psi\colon(0,\delta]\to[0,\infty)$\cfadd{definition:definable-function}
  such that for every interval $I\subseteq\R$,
  all locally absolutely continuous $r\colon I\to(0,\mathfrak R)$ and
   $x\colon I\to M$, and
  all Lebesgue measurable $\alpha\colon I\to[0,\infty)$ and
  $A\colon I\to\R^{D\times D}$
  with  
  $\midcup_{t\in I}\{(r(t),x(t))\}\subseteq\mathcal O$ and
  with the property that for Lebesgue-almost every $t\in I$ and
  every $v\in T_{x(t)}M$\cfadd{definition:tangent-space} it holds that
  \begin{eqn}
    \llabel{ass1}
    A(t)v\in T_{x(t)}M,
    \qquad \langle A(t)v,v\rangle\geq\mathsf a\norm{v}^2,
    \qquad
    \norm{A(t)v}\leq\mathsf b\norm{v},
  \end{eqn}
  \begin{eqn}
    \llabel{ass2}
    r'(t)=-\alpha(t)(r(t))^2
      (\partial_r h)(r(t),x(t)),
    \qquad\text{and}\qquad
    x'(t)=-\alpha(t)A(t)(\nabla_Mh)(r(t),x(t))
  \end{eqn}
  \cfadd{definition:angular-gradient}%
  we have that
  \begin{enumerate}[label=(\roman*)]
    \item \llabel{item1} for all $R\in(0,\delta]$ with $r(I)\subseteq(0,R]$ it holds that
    \begin{eqn}
      \llabel{concl1}
      \int_I\sqrt{\abs{r'(t)}^2+\norm{x'(t)}^2}\dd t
        \leq \mathfrak C
          +\mathfrak C\int_0^R\psi(s)\dd s<\infty,
    \end{eqn}
    and
    \item \llabel{item2} for all $\tau\in\R$,
    $T\in(\tau,\infty]$ with $I=[\tau,T)$
    and $\limsup_{t\to T}|r(t)|=0$ it holds that
    \begin{eqn}
      \llabel{concl2}
      \int_{\tau}^T\norm{x'(t)}\dd t<\infty
      \dott
    \end{eqn}
  \end{enumerate}
  \cfout[.]
\end{athm}
\begin{aproof}
  \argument{
    \cref{prop:buildQ};
  }{
    that there exist $\eta\in(0,\mathfrak R)$,
    a bounded $Q\in C(\mathcal O\cap((0,\eta)\times M),[0,\infty))$, and a
    definable and integrable $m\in C((0,\eta),[0,\infty))$\cfadd{definition:definable-function}
    which satisfy that
    \begin{enumerate}[label=(\Roman*)]
      \item for all $a,b\in(0,\eta)$ with $a\leq b$ it holds that
      \begin{eqn}
        \llabel{eq1}
        \sup_{s\in[a,b]}m(s)<\infty
        \qquad\text{and}\qquad
        \int_0^\eta m(s)\dd s<\infty
      \end{eqn}
      and
      \item for every interval $I\subseteq\R$,
      all locally absolutely continuous $r\colon I\to(0,\eta)$
      and $x\colon I\to M$, and
      all Lebesgue measurable $\alpha\colon I\to[0,\infty)$ and
      $A\colon I\to\R^{D\times D}$
      with $\midcup_{t\in I}\{(r(t),x(t))\}\subseteq\mathcal O$ and
      with the property that for Lebesgue-almost all $t\in I$ and
      every $v\in T_{x(t)}M$\cfadd{definition:tangent-space} it holds that
      \begin{eqn}
        \llabel{eq0.1}
        A(t)v\in T_{x(t)}M,
        \qquad \langle A(t)v,v\rangle\geq\mathsf a\norm{v}^2,
        \qquad
        \norm{A(t)v}\leq\mathsf b\norm{v},
      \end{eqn}
      \begin{eqn}
        \llabel{eq0.2}
        r'(t)=-\alpha(t)(r(t))^2
          (\partial_r h)(r(t),x(t)),
        \qquad\text{and}\qquad
        x'(t)=-\alpha(t)A(t)(\nabla_Mh)(r(t),x(t))
      \end{eqn}
      \cfadd{definition:angular-gradient}%
      we have that for all $a,b\in I$ with $a\leq b$ it holds that
      \begin{eqn}
        \llabel{eq1.3}
        Q(r(b),x(b))-Q(r(a),x(a))
          \leq-\int_a^b\norm{x'(t)}\dd t
            +\int_a^b m(r(t))\abs{r'(t)}\dd t
        \dott
      \end{eqn}
    \end{enumerate}
  }
  \startnewargseq
  \argument{
    \cref{prop:cone-annular-length}
    (applied with $\Omega\curvearrowleft\mathcal O$ in the notation of \cref{prop:cone-annular-length});
  }{
    that there exist $\varrho\in(0,\mathfrak R)$ and
    $\mathfrak K\in(0,\infty)$ which satisfy that for every interval
    $I\subseteq\R$,
    all locally absolutely continuous $r\colon I\to(0,\mathfrak R)$
    and $x\colon I\to M$,
    and all Lebesgue measurable $\alpha\colon I\to[0,\infty)$ and
    $A\colon I\to\R^{D\times D}$
    with $\midcup_{t\in I}\{(r(t),x(t))\}\subseteq\mathcal O$,
    with the property that for Lebesgue-almost all $t\in I$ and
    every $v\in T_{x(t)}M$\cfadd{definition:tangent-space} it holds that
    \begin{eqn}
      \llabel{eq1.1}
      A(t)v\in T_{x(t)}M,\qquad
      \langle A(t)v,v\rangle\geq\mathsf a\norm{v}^2,
      \qquad\text{and}\qquad
      \norm{A(t)v}\leq\mathsf b\norm{v},
    \end{eqn}
    and with the property that for Lebesgue-almost all $t\in I$
    it holds that 
    \begin{eqn}
      \llabel{eq1.2}
      r'(t)=-\alpha(t)(r(t))^2
        (\partial_r h)(r(t),x(t))
      \qquad\text{and}\qquad
      x'(t)=-\alpha(t)A(t)(\nabla_Mh)(r(t),x(t))
    \end{eqn}
    \cfadd{definition:angular-gradient}%
    we have for all $u\in(0,\varrho)$ that
    \begin{eqn}
      \llabel{eq2}
      \int_{r^{-1}([u/2,u])}
        \sqrt{\abs{r'(t)}^2+(r(t))^2\norm{x'(t)}^2}\dd t
        \leq\mathfrak K u
      \dott
    \end{eqn}
  }
  \startnewargseq
  \argument{
    \lref{eq1.3};
  }{
    that for every interval $I\subseteq\R$,
    all locally absolutely continuous $r\colon I\to(0,\mathfrak R)$
    and $x\colon I\to M$,
    and all Lebesgue measurable $\alpha\colon I\to[0,\infty)$ and
    $A\colon I\to\R^{D\times D}$
    with $\midcup_{t\in I}\{(r(t),x(t))\}\subseteq\mathcal O$,
    with the property that for Lebesgue-almost all $t\in I$ and
    every $v\in T_{x(t)}M$ it holds that
    \begin{eqn}
      \llabel{eq2.1}
      A(t)v\in T_{x(t)}M,\qquad
      \langle A(t)v,v\rangle\geq\mathsf a\norm{v}^2,
      \qquad\text{and}\qquad
      \norm{A(t)v}\leq\mathsf b\norm{v},
    \end{eqn}
    and with the property that for Lebesgue-almost all $t\in I$
    it holds that
    \begin{eqn}
      \llabel{eq2.2}
      r'(t)=-\alpha(t)(r(t))^2
        (\partial_r h)(r(t),x(t))
      \qquad\text{and}\qquad
      x'(t)=-\alpha(t)A(t)(\nabla_Mh)(r(t),x(t))
    \end{eqn}
    we have that for all $a,b\in I$ with $a<b$ and
    $r([a,b])\subseteq(0,\eta)$ it holds that
    \begin{eqn}
      \llabel{eq3}
      Q(r(b),x(b))-Q(r(a),x(a))
        \leq-\int_a^b\norm{x'(t)}\dd t
          +\int_a^b m(r(t))\abs{r'(t)}\dd t
      \dott
    \end{eqn}
  }
  In the following let $\zeta,B\in\R$ satisfy
  \begin{eqn}[a]
    \llabel{setting1}
    \zeta=\frac12\min\{\eta,\varrho\}
    \qquad\text{and}\qquad
    B=\sup\bigl(\{0\}\cup Q(\mathcal O\cap((0,\eta)\times M))\bigr)
    \dott
  \end{eqn}
  \cfload[.]
  \startnewargseq
  \argument{
    \lref{eq1};
    \lref{setting1};
    the fact that $m|_{(0,\zeta]}$ is nonnegative, locally bounded, and
    definable;
    \itref{lem:radial-integrability}{item2} (applied with
    $R\curvearrowleft\zeta$, $m\curvearrowleft m|_{(0,\zeta]}$
    in the notation of \cref{lem:radial-integrability});
  }{
    that there exist $\delta\in(0,\zeta)$ and a definable,
    non-increasing, and integrable
    $\varphi\colon(0,\delta]\to[0,\infty)$ which satisfy that for all
    $s\in(0,\delta]$ it holds that
    \begin{eqn}
      \llabel{eq4}
      m(s)\leq\varphi(s)
      \dott
    \end{eqn}
  }
\startnewargseq
  \argument{
    \lref{setting1};
    the fact that $Q$ is nonnegative;
  }{
    that for all $(s,y)\in\mathcal O\cap((0,\eta)\times M)$
    it holds that
    \begin{eqn}
      \llabel{eq5}
      0\leq Q(s,y)\leq B
      \dott
    \end{eqn}
  }
  We first prove that for every interval $I\subseteq\R$,
  all locally absolutely continuous $r\colon I\to(0,\mathfrak R)$
  and $x\colon I\to M$,
  and all Lebesgue measurable $\alpha\colon I\to[0,\infty)$ and
  $A\colon I\to\R^{D\times D}$
  with $\midcup_{t\in I}\{(r(t),x(t))\}\subseteq\mathcal O$,
  with the property that for Lebesgue-almost all $t\in I$ and
  every $v\in T_{x(t)}M$\cfadd{definition:tangent-space} it holds that
  \begin{eqn}
    \llabel{eq5.1}
    A(t)v\in T_{x(t)}M,\qquad
    \langle A(t)v,v\rangle\geq\mathsf a\norm{v}^2,
    \qquad\text{and}\qquad
    \norm{A(t)v}\leq\mathsf b\norm{v},
  \end{eqn}
  and with the property that for Lebesgue-almost all $t\in I$
  it holds that 
  \begin{eqn}
    \llabel{eq5.2}
    r'(t)=-\alpha(t)(r(t))^2
      (\partial_r h)(r(t),x(t))
    \qquad\text{and}\qquad
    x'(t)=-\alpha(t)A(t)(\nabla_Mh)(r(t),x(t))
  \end{eqn}
  \cfadd{definition:angular-gradient}%
  we have that for every $R\in(0,\delta]$ with $\forall\,t\in I \colon r(t)\leq R$ it holds that
  \begin{eqn}
    \llabel{eq5.3}
    \int_I\sqrt{\abs{r'(t)}^2+\norm{x'(t)}^2}\dd t
      \leq 2\mathfrak K R+B
        +2\mathfrak K\int_0^R\varphi(s)\dd s<\infty
    \dott
  \end{eqn}
  \cfload[.]
  \startnewargseq
  To prove \lref{eq5.3}, let $I\subseteq\R$ be an interval,
  let $r\colon I\to(0,\mathfrak R)$ and $x\colon I\to M$
  be locally absolutely continuous,
  let $\alpha\colon I\to[0,\infty)$,
  $A\colon I\to\R^{D\times D}$, and
  $\mathfrak v\colon I\to[0,\infty)$ be Lebesgue measurable,
  assume that $\midcup_{t\in I}\{(r(t),x(t))\}\subseteq\mathcal O$,
  assume that for Lebesgue-almost all $t\in I$ and all
  $v\in T_{x(t)}M$ it holds that
  \begin{eqn}
    \llabel{setting2}
    A(t)v\in T_{x(t)}M,
    \qquad \langle A(t)v,v\rangle\geq\mathsf a\norm{v}^2,
    \qquad\text{and}\qquad
    \norm{A(t)v}\leq\mathsf b\norm{v},
  \end{eqn}
  and assume that for Lebesgue-almost all $t\in I$ it holds that
  \begin{eqn}
    \llabel{setting3}
    \begin{gathered}
      r'(t)=-\alpha(t)(r(t))^2(\partial_r h)(r(t),x(t)),
      \qquad
      x'(t)=-\alpha(t)A(t)(\nabla_Mh)(r(t),x(t)),\\
      \text{and}\qquad
      \mathfrak v(t)=\sqrt{\abs{r'(t)}^2+(r(t))^2\norm{x'(t)}^2}
      \dott
    \end{gathered}{}
  \end{eqn}
  \cfadd{definition:tangent-space}\cfadd{definition:angular-gradient}\cfload[.]
  \startnewargseq
  \argument{
    \lref{eq2};\lref{setting2};\lref{setting3};
    the fact that $\midcup_{t\in I}\{(r(t),x(t))\}\subseteq\mathcal O$;
  }{
    that for Lebesgue-almost all $t\in I$ it holds that
    $\abs{r'(t)}\leq\mathfrak v(t)$ and that for all
    $u\in(0,\varrho)$ it holds that
    \begin{eqn}
      \llabel{eq6}
      \int_{r^{-1}([u/2,u])}\mathfrak v(t)\dd t
        \leq\mathfrak K u
      \dott
    \end{eqn}
  }
  \argument{
    \lref{eq6};
    \lref{setting1};
    \lref{setting3};
    the fact that $\delta<\zeta<\varrho$;
    the fact that $r$ is locally absolutely continuous;
    \itref{lem:v}{item1} (applied with
    $\mathfrak C\curvearrowleft\mathfrak K$,
    $r\curvearrowleft\bigl(I\ni t\mapsto r(t)\in(0,\infty)\bigr)$,
    $\rho\curvearrowleft R$ in the notation of
    \cref{lem:v});
  }{
    that for all $R\in(0,\delta]$ it holds that
    \begin{eqn}
      \llabel{eq7}
      \int_{r^{-1}((0,R])}\abs{r'(t)}\dd t
        \leq\int_{r^{-1}((0,R])}\mathfrak v(t)\dd t
        \leq 2\mathfrak K R
      \dott
    \end{eqn}
  }
  \argument{
    \lref{eq1};
    \lref{setting1};
    \lref{eq4};
    \lref{setting3};
    \lref{eq6};
    the fact that $m|_{(0,\zeta]}$ is nonnegative, locally bounded, and
    definable;
    the fact that $\varphi$ is definable, non-increasing, and integrable;
    the fact that $r$ is locally absolutely continuous;
    the fact that $\delta<\zeta<\varrho$;
    the monotonicity of the Lebesgue integral;
    \itref{lem:radial-integrability}{item4} (applied with
    $R\curvearrowleft\zeta$, $m\curvearrowleft m|_{(0,\zeta]}$,
    $\mathfrak C\curvearrowleft\mathfrak K$, $\psi\curvearrowleft\varphi$,
    $r\curvearrowleft\bigl(I\ni t\mapsto r(t)\in(0,\infty)\bigr)$,
    $S\curvearrowleft R$ in the notation of
    \cref{lem:radial-integrability});
  }{
    that for all $R\in(0,\delta]$ it holds that
    \begin{eqn}
      \llabel{eq8}
      \int_{r^{-1}((0,R])}m(r(t))\abs{r'(t)}\dd t
        \leq 2\mathfrak K\int_0^R\varphi(s)\dd s<\infty
      \dott
    \end{eqn}
  }
  \argument{
    \lref{eq8};
    \lref{eq3};
    \lref{eq5};
    the fact that $\midcup_{t\in I}\{(r(t),x(t))\}\subseteq\mathcal O$;
    the fact that $\delta<\eta$;
    the fact that $m$ is nonnegative;
  }{
    that for every $R\in(0,\delta]$ with the property that
    for all $t\in I$ it holds that $r(t)\leq R$ we have for all
    $a,b\in I$ with $a<b$ that
    \begin{eqn}[a]
      \llabel{eq9}
      \int_a^b\norm{x'(t)}\dd t
        \leq Q(r(a),x(a))-Q(r(b),x(b))
          +\int_a^b m(r(t))\abs{r'(t)}\dd t\leq B+2\mathfrak K\int_0^R\varphi(s)\dd s
      \dott
    \end{eqn}
  }
  \argument{
    \lref{eq9}; the fact that $I$ is an interval;
    the monotone convergence theorem;
  }{
    that for every $R\in(0,\delta]$ with the property that
    for all $t\in I$ it holds that $r(t)\leq R$ we have that
    \begin{eqn}
      \llabel{eq10}
      \int_I\norm{x'(t)}\dd t
        \leq B+2\mathfrak K\int_0^R\varphi(s)\dd s<\infty
      \dott
    \end{eqn}
  }
  \argument{
    \lref{eq10};
    \lref{eq7};
    the fact that for all $a,b\in[0,\infty)$ it holds that
    $\sqrt{a^2+b^2}\leq a+b$;
  }{
    that for every $R\in(0,\delta]$ with the property that
    for all $t\in I$ it holds that $r(t)\leq R$ we have that
    \begin{eqn}[a]
      \llabel{eq11}
      \int_I\sqrt{\abs{r'(t)}^2+\norm{x'(t)}^2}\dd t
        \leq\int_I\abs{r'(t)}\dd t
          +\int_I\norm{x'(t)}\dd t\leq 2\mathfrak K R+B
          +2\mathfrak K\int_0^R\varphi(s)\dd s<\infty
      \dott
    \end{eqn}
  }
  \argument{
    \lref{eq11};
    the arbitrariness of $I$, $r$, $x$, $\alpha$, and $A$;
  }[verbs=ep]{
    \lref{eq5.3}\dott
  }
  We next prove that for every interval $I\subseteq\R$,
  all locally absolutely continuous $r\colon I\to(0,\mathfrak R)$
  and $x\colon I\to M$,
  and all Lebesgue measurable $\alpha\colon I\to[0,\infty)$ and
  $A\colon I\to\R^{D\times D}$
  with $\midcup_{t\in I}\{(r(t),x(t))\}\subseteq\mathcal O$,
  with the property that for Lebesgue-almost all $t\in I$ and
  every $v\in T_{x(t)}M$\cfadd{definition:tangent-space} it holds that
  \begin{eqn}
    \llabel{eq11.1}
    A(t)v\in T_{x(t)}M,\qquad
    \langle A(t)v,v\rangle\geq\mathsf a\norm{v}^2,
    \qquad\text{and}\qquad
    \norm{A(t)v}\leq\mathsf b\norm{v},
  \end{eqn}
  and with the property that for Lebesgue-almost all $t\in I$
  it holds that 
  \begin{eqn}
    \llabel{eq11.2}
    r'(t)=-\alpha(t)(r(t))^2
      (\partial_r h)(r(t),x(t))
    \qquad\text{and}\qquad
    x'(t)=-\alpha(t)A(t)(\nabla_Mh)(r(t),x(t))
  \end{eqn}
  we have that for every $\tau\in\R$,
  $T\in(\tau,\infty]$ with $I=[\tau,T)$
  and $\lim_{t\nearrow T}r(t)=0$ it holds that
  \begin{eqn}
    \llabel{eq11.3}
    \int_{\tau}^T\norm{x'(t)}\dd t<\infty
    \dott
  \end{eqn}
  \cfload[.]
  \startnewargseq
  To prove \lref{eq11.3}, let $I\subseteq\R$ be an interval,
  let $r\colon I\to(0,\mathfrak R)$ and $x\colon I\to M$
  be locally absolutely continuous,
  let $\alpha\colon I\to[0,\infty)$ and
  $A\colon I\to\R^{D\times D}$ be Lebesgue measurable,
  let $\tau\in\R$, $T\in(\tau,\infty]$,
  assume that $\midcup_{t\in I}\{(r(t),x(t))\}\subseteq\mathcal O$,
  assume that
  \begin{eqn}
    \llabel{setting5}
    I=[\tau,T)
    \qquad\text{and}\qquad
    \lim_{t\nearrow T}r(t)=0,
  \end{eqn}
  assume that for Lebesgue-almost all $t\in I$ and all
  $v\in T_{x(t)}M$ it holds that
  \begin{eqn}
    \llabel{setting6}
    A(t)v\in T_{x(t)}M,\qquad
    \langle A(t)v,v\rangle\geq\mathsf a\norm{v}^2,
    \qquad\text{and}\qquad
    \norm{A(t)v}\leq\mathsf b\norm{v},
  \end{eqn}
  and assume that for Lebesgue-almost all $t\in I$ it holds that
  \begin{eqn}
    \llabel{setting7}
    r'(t)=-\alpha(t)(r(t))^2
      (\partial_r h)(r(t),x(t))
    \qquad\text{and}\qquad
    x'(t)=-\alpha(t)A(t)(\nabla_Mh)(r(t),x(t))
    \dott
  \end{eqn}
  \cfadd{definition:tangent-space}\cfadd{definition:angular-gradient}\cfload[.]
  \startnewargseq
  \argument{
    \lref{setting5};
    the fact that $\lim_{t\nearrow T}r(t)=0$;
    the fact that $\delta>0$;
  }{
    that there exists
    $\sigma\in[\tau,T)$ such that for all
    $t\in[\sigma,T)$ it holds that
    \begin{eqn}
      \llabel{eq12}
      0<r(t)\leq\delta
      \dott
    \end{eqn}
  }
  \argument{
    \lref{eq12};
    \lref{eq5.3};
    \lref{setting5};
    \lref{setting6};
    \lref{setting7};
    the fact that $\midcup_{t\in I}\{(r(t),x(t))\}\subseteq\mathcal O$;
  }{
    that there exists
    $\sigma\in[\tau,T)$ which satisfies that
    \begin{eqn}[a]
      \llabel{eq13}
      \int_{\sigma}^T\norm{x'(t)}\dd t
        \leq\int_{\sigma}^T
          \sqrt{\abs{r'(t)}^2+\norm{x'(t)}^2}\dd t\leq 2\mathfrak K\delta+B
          +2\mathfrak K\int_0^\delta\varphi(s)\dd s<\infty
      \dott
    \end{eqn}
  }
  \startnewargseq
  \argument{
    \lref{eq13};
    the fact that $x$ is locally absolutely continuous;
  }{
    that
    \begin{eqn}
      \llabel{eq13.1}
      \int_\tau^T\norm{x'(t)}\dd t
        =\int_\tau^\sigma\norm{x'(t)}\dd t
          +\int_\sigma^T\norm{x'(t)}\dd t<\infty
      \dott
    \end{eqn}
  }
  \argument{
    \lref{eq13.1};
    the arbitrariness of $I$, $r$, $x$, $\alpha$, $A$, $\tau$, and $T$;
  }[verbs=ep]{
    \lref{eq11.3}\dott
  }
  \startnewargseq
  In the following let $\mathfrak C=\max\{2\mathfrak K,B\}$, and let
  $\psi\colon(0,\delta]\to[0,\infty)$ satisfy for all $s\in(0,\delta]$
  that $\psi(s)=1+\varphi(s)$\dott \cfload[.]
  \startnewargseq
  \argument{
    the fact that $\varphi$ is definable, non-increasing, and integrable;
    the fact that for all $s\in(0,\delta]$ it holds that $\psi(s)=1+\varphi(s)$;
    the fact that $\delta\in(0,\infty)$;
    the linearity of the integral;
  }{
    that $\psi$ is definable, non-increasing, and integrable and that
    for all $R\in(0,\delta]$ it holds that
    \begin{eqn}[i]
      \llabel{eq14}
      \int_0^R\psi(s)\dd s=R+\int_0^R\varphi(s)\dd s
      \dott
    \end{eqn}
  }
  \argument{
    \lref{eq14};
    \lref{eq5.3};
    \lref{eq11.3};
    the fact that $\mathfrak C=\max\{2\mathfrak K,B\}>0$;
    the fact that $\varphi$ is nonnegative;
  }[verbs=ep]{
    \lref{item1,item2}\dott}
\end{aproof}

\subsection{Uniform upper bounds for the length of the angular component in the DE system}
\label{subsec:angular-occupation}

\begin{athm}{proposition}{prop:angular-occupation}
  Let $D,n\in\N$, $ \mathfrak R, \mathsf a \in(0,\infty) $, $ \mathsf b \in [ \mathsf{a}, \infty ) $,
  let $M\subseteq\R^D$ be a
  nonempty, compact, connected, and definable\cfadd{definition:definable-set}
  \dsubmanifold{2}{D}{n}\cfadd{definition:embedded-manifold},
  and let $h=(h(r,x))_{(r,x)\in(0,\mathfrak R)\times M}
  \in \COne{(0,\mathfrak R)\times M}{\R}$ be
  definable\cfadd{definition:c1-on-manifold}\cfadd{definition:definable-function}\dott
  \cfload[.]
  Then there exist $R\in(0,\mathfrak R)$, $L\in\R$ such that
  for every interval $I\subseteq\R$,
  all locally absolutely continuous $r\colon I\to(0,\mathfrak R)$ and
   $x\colon I\to M$,
  and all Lebesgue measurable $\alpha\colon I\to[0,\infty)$
  and $A\colon I\to\R^{D\times D}$
  with the property that for Lebesgue-almost every $t\in I$ and
  every $v\in T_{x(t)}M$\cfadd{definition:tangent-space} it holds that
  \begin{eqn}
    \llabel{ass1}
    A(t)v\in T_{x(t)}M,\qquad
    \langle A(t)v,v\rangle\geq\mathsf a\norm{v}^2,
    \qquad\text{and}\qquad
    \norm{A(t)v}\leq\mathsf b\norm{v},
  \end{eqn}
  \begin{eqn}
    \llabel{ass2}
    r'(t)=-\alpha(t)(r(t))^2
      (\partial_r h)(r(t),x(t)),
    \qquad\text{and}\qquad
    x'(t)=-\alpha(t)A(t)(\nabla_Mh)(r(t),x(t))
  \end{eqn}
  \cfadd{definition:angular-gradient}%
  we have that
  \begin{eqn}
    \llabel{concl1}
    \int_{r^{-1}((0,R])}\norm{x'(t)}\dd t\leq L
    \dott
  \end{eqn}
  \cfout[.]
\end{athm}
\begin{aproof}
  \argument{
    \itref{cor:angular-length}{item1}
    (applied with $\mathcal O\curvearrowleft(0,\mathfrak R)\times M$,
    $R\curvearrowleft\delta$ in the notation of
    \cref{cor:angular-length});
    \cref{prop:cone-annular-length}
    (applied with $\Omega\curvearrowleft(0,\mathfrak R)\times M$
    in the notation of \cref{prop:cone-annular-length});
  }{
    that there exist $\delta,\varrho\in(0,\mathfrak R)$,
    $\mathfrak C\in(0,\infty)$, and $K\in[0,\infty)$ such that
    for every interval $I\subseteq\R$,
    all locally absolutely continuous $r\colon I\to(0,\mathfrak R)$
    and $x\colon I\to M$,
    and all Lebesgue measurable $\alpha\colon I\to[0,\infty)$ and
    $A\colon I\to\R^{D\times D}$
    with the property that for Lebesgue-almost all $t\in I$ and
    every $v\in T_{x(t)}M$\cfadd{definition:tangent-space} it holds that
    \begin{eqn}
      \llabel{eq0.1}
      A(t)v\in T_{x(t)}M,\qquad
      \langle A(t)v,v\rangle\geq\mathsf a\norm{v}^2,
      \qquad\text{and}\qquad
      \norm{A(t)v}\leq\mathsf b\norm{v},
    \end{eqn}
    and with the property that for Lebesgue-almost all $t\in I$
    it holds that 
    \begin{eqn}
      \llabel{eq0.2}
      r'(t)=-\alpha(t)(r(t))^2
        (\partial_r h)(r(t),x(t))
      \qquad\text{and}\qquad
      x'(t)=-\alpha(t)A(t)(\nabla_Mh)(r(t),x(t)),
    \end{eqn}
    \cfadd{definition:angular-gradient}%
  we have that
    \begin{enumerate}[label=(\Roman*)]
      \item if for all $t\in I$ it holds that $r(t)\leq\delta$, then
      \begin{eqn}
        \llabel{eq1}
        \int_I\norm{x'(t)}\dd t
          \leq\int_I\sqrt{\abs{r'(t)}^2+\norm{x'(t)}^2}\dd t
          \leq K,
      \end{eqn}
      and
      \item for all $u\in(0,\varrho)$ it holds that
      \begin{eqn}
        \llabel{eq2}
        \int_{r^{-1}([\nicefrac{u}{2},u])}\abs{r'(t)}\dd t
          \leq\int_{r^{-1}([\nicefrac{u}{2},u])}
            \sqrt{\abs{r'(t)}^2+(r(t))^2\norm{x'(t)}^2}\dd t
          \leq\mathfrak C u
        \dott
      \end{eqn}
    \end{enumerate}
  }
  In the following let $S\in(0,\mathfrak R)$ and $L\in[0,\infty)$ satisfy
  \begin{eqn}
    \llabel{setting1}
    S=\frac12\min\{\delta,\varrho\}
    \qquad\text{and}\qquad
    L=(1+2\mathfrak C)K
    \dott
  \end{eqn}
  \cfload[.]
  \startnewargseq
  We claim that for every interval $I\subseteq\R$,
  all locally absolutely continuous $r\colon I\to(0,\mathfrak R)$
  and $x\colon I\to M$,
  and all Lebesgue measurable $\alpha\colon I\to[0,\infty)$ and
  $A\colon I\to\R^{D\times D}$
  with the property that for Lebesgue-almost all $t\in I$ and every
  $v\in T_{x(t)}M$\cfadd{definition:tangent-space} it holds that
  \begin{eqn}
    A(t)v\in T_{x(t)}M,\qquad
    \langle A(t)v,v\rangle\geq\mathsf a\norm{v}^2,
    \qquad\text{and}\qquad
    \norm{A(t)v}\leq\mathsf b\norm{v},
  \end{eqn}
  and with the property that for Lebesgue-almost all $t\in I$
  it holds that 
  \begin{eqn}
    r'(t)=-\alpha(t)(r(t))^2
      (\partial_r h)(r(t),x(t))
    \qquad\text{and}\qquad
    x'(t)=-\alpha(t)A(t)(\nabla_Mh)(r(t),x(t))
  \end{eqn}
  we have for all $a,b\in I$ with $a<b$ that
  \begin{eqn}
    \llabel{eq3}
    \int_{[a,b]\cap r^{-1}((0,\nicefrac{S}{2}])}
      \norm{x'(t)}\dd t\leq L
    \dott
  \end{eqn}
  \cfload[.]
  \startnewargseq
  To prove this claim, let $I\subseteq\R$ be an interval,
  let $r\colon I\to(0,\mathfrak R)$ and $x\colon I\to M$
  be locally absolutely continuous,
  let $\alpha\colon I\to[0,\infty)$ and
  $A\colon I\to\R^{D\times D}$ be Lebesgue measurable,
  assume that for Lebesgue-almost all $t\in I$ and every
  $v\in T_{x(t)}M$ it holds that
  \begin{equation}
    \llabel{setting2}
    A(t)v\in T_{x(t)}M,\qquad
    \langle A(t)v,v\rangle\geq\mathsf a\norm{v}^2,
    \qquad\text{and}\qquad
    \norm{A(t)v}\leq\mathsf b\norm{v},
  \end{equation}
  assume that for Lebesgue-almost all $t\in I$ it holds that $r$ and $x$
  are differentiable at $t$ and that
  \begin{eqn}
    \llabel{setting3}
    r'(t)=-\alpha(t)(r(t))^2
      (\partial_r h)(r(t),x(t))
    \qquad\text{and}\qquad
    x'(t)=-\alpha(t)A(t)(\nabla_Mh)(r(t),x(t)),
  \end{eqn}
  and let $a,b\in I$ satisfy $a<b$\dott \cfload[.]
  \startnewargseq
  In the following let $U,E\subseteq[a,b]$ satisfy
  \begin{eqn}
    \llabel{setting4}
    U=\{t\in[a,b]\colon r(t)<S\}
    \qquad\text{and}\qquad
    E=[a,b]\cap r^{-1}((0,\nicefrac{S}{2}])
    \dott
  \end{eqn}
  \cfload[.]
  \startnewargseq
  \argument{
    the fact that $r|_{[a,b]}$ is absolutely continuous;
    the fact that for every interval $J\subseteq\R$ and every absolutely continuous
    $f\colon J\to\R$ it holds that $f$ is continuous;
  }{
    that $r|_{[a,b]}$ is continuous\dott
    \llabel{eq3.1}
  }
  \argument{
    \lref{eq3.1};
    \lref{setting4};
    the fact that $(-\infty,S)$ is open;
  }{
    that $U$ is open in $[a,b]$\dott
  }
  \argument{
    \lref{setting4};
    the fact that $S>0$;
  }{
    that
    \begin{eqn}[a]
      \llabel{eq3.4}
      E=\{t\in[a,b]\colon0<r(t)\leq\nicefrac{S}{2}\}\subseteq\{t\in[a,b]\colon r(t)<S\}=U
      \dott
    \end{eqn}
  }
  In the following let $W_t\subseteq U$, $t\in E$, satisfy for all
  $t\in E$ that
  \begin{eqn}
    \llabel{setting5.1}
    W_t=\bigcup\bigl\{V\subseteq U\colon\bigl[
      (V\text{ is an interval})\wedge(t\in V)\bigr]\bigr\}
    \dott
  \end{eqn}
  \cfload[.]
  \startnewargseq
  \argument{
    \lref{eq3.4};
    \lref{setting5.1};
  }{
    that for all $t\in E$ it holds that $W_t$ is an interval and that
    \begin{eqn}
      \llabel{eq3.4.2}
      \{t\}\subseteq W_t\subseteq U
      \dott
    \end{eqn}
  }
  \argument{
    \lref{eq3.4.2};
    \lref{setting5.1};
  }{
    that for every $t\in E$ and every interval $V$ with
    $W_t\subseteq V\subseteq U$ it holds that
    \begin{eqn}
      \llabel{eq3.4.3}
      V=W_t
      \dott
    \end{eqn}
  }
  \argument{
    \lref{setting5.1};
    \lref{eq3.4.2};
  }{
    that for every nonempty interval $J\subseteq U$ with the property that
    $J\cap E\neq\emptyset$ and with the property that for every
    interval $V$ with $J\subseteq V\subseteq U$ it holds that $V=J$
   we have for all $t\in J\cap E$ that
    \begin{eqn}
      \llabel{eq3.4.4}
      J=W_t
      \dott
    \end{eqn}
  }
  \argument{
    \lref{eq3.4.4};
    \lref{eq3.4.2};
    \lref{eq3.4.3};
    the fact that for all $t\in E$ it holds that $t\in W_t\cap E$;
  }{
    that there exists a set $\mathscr J$ which satisfies
    \begin{equation}
      \llabel{setting4.1}
      \begin{split}
      \mathscr J
        &=\{W_t\colon t\in E\}\\
        &=\bigl\{J\subseteq U\colon\bigl[
          (J\cap E\neq\emptyset)\wedge (J\text{ is an interval})\wedge
          \bigl(\forall\,V\text{ interval},\ J\subseteq V\subseteq U\colon V=J\bigr)
          \bigr]\bigr\}
      \dott
      \end{split}
    \end{equation}
  }
  In the following let $\mathcal H\subseteq\mathscr J$ satisfy
  \begin{eqn}
    \llabel{setting5}
    \mathcal H=\{J\in\mathscr J\colon a\notin J\}
    \dott
  \end{eqn}
  \cfload[.]
  \startnewargseq
  \argument{
    \lref{setting4.1};
  }{
    that for all $J\in\mathscr J$ it holds that $J$ is an
    interval\dott
    \llabel{eq3.2}
  }
  \argument{
    \lref{eq3.2};
    \lref{setting4.1};
  }{
    that for all $J,V\in\mathscr J$ with
    $J\cap V\neq\emptyset$ it holds that
    $J=J\cup V=V$\dott
    \llabel{eq3.2.1}
  }
  \argument{
    \lref{eq3.2.1};
  }{
    that for all
    $J,V\in\mathscr J$ with $J\neq V$ it holds that
    \begin{eqn}
      \llabel{eq3.3}
      J\cap V=\emptyset
      \dott
    \end{eqn}
  }
  \argument{
    \lref{setting4.1};
  }{
    that for all $t\in E$ it holds that
    \begin{eqn}
      \llabel{eq3.4.1}
      W_t\in\mathscr J
      \dott
    \end{eqn}
  }
  \argument{
    \lref{eq3.4.1};
    \lref{eq3.4.2};
  }{
    that
    \begin{eqn}
      \llabel{eq3.5}
      E\subseteq\midcup_{J\in\mathscr J}J
      \dott
    \end{eqn}
  }
  \argument{
    \lref{setting5};
  }{
    that
    \begin{eqn}
      \llabel{eq3.6}
      \mathscr J\backslash\mathcal H
        =\{J\in\mathscr J\colon a\in J\}
      \dott
    \end{eqn}
  }
  \argument{
    \lref{eq3.6};
    \lref{eq3.3};
    the fact that for all $J,V\in\mathscr J\backslash\mathcal H$
    it holds that $a\in J\cap V$;
  }{
    that for all $J,V\in\mathscr J\backslash\mathcal H$
    it holds that $J=V$\dott
    \llabel{eq3.7}
  }
  \argument{
    \lref{eq3.7};
    \lref{eq3.5};
  }{
    that
    \begin{eqn}
      \llabel{eq4}
      E\subseteq\midcup_{J\in\mathscr J}J
      \qquad\text{and}\qquad
      \#(\mathscr J\backslash\mathcal H)\leq1
      \dott
    \end{eqn}
  }
  In the following let
  $\sigma_J,\tau_J\in[a,b]$, $J\in\mathcal H$, satisfy for all
  $J\in\mathcal H$ that
  \begin{eqn}
    \llabel{setting6}
    \sigma_J=\inf J
    \qquad\text{and}\qquad
    \tau_J\in J\cap E
    \dott
  \end{eqn}
  \cfload[.]
  \startnewargseq
  \argument{
    \lref{setting4};
    \lref{setting4.1};
    \lref{setting5};
    \lref{setting6};
    the fact that $S>0$;
  }{
    that for all $J\in\mathcal H$ it holds that
    \begin{eqn}
      \llabel{eq4.1}
      a<\tau_J\leq b
      \qquad\text{and}\qquad
      r(\tau_J)\leq\frac{S}{2}<S
      \dott
    \end{eqn}
  }
  \argument{
    \lref{eq4.1};
    \lref{setting4};
    the fact that $r|_{[a,b]}$ is continuous;
  }{
    that for all $J\in\mathcal H$ there exists $s\in[a,\tau_J)$ such that
    \begin{eqn}
      \llabel{eq4.2}
      [s,\tau_J]\subseteq U
      \dott
    \end{eqn}
  }
  \argument{
    \lref{eq4.2};
    \lref{setting4.1};
    \lref{setting5};
    \lref{setting6};
  }{
    that for all $J\in\mathcal H$ there exists $s\in[a,\tau_J)$ such that
    \begin{eqn}
      \llabel{eq4.3}
      [s,\tau_J]\subseteq J
      \dott
    \end{eqn}
  }
  \argument{
    \lref{eq4.3};
    \lref{setting6};
  }{
    that for all $J\in\mathcal H$ it holds that
    \begin{eqn}
      \llabel{eq4.4}
      \sigma_J<\tau_J
      \dott
    \end{eqn}
  }
  \argument{
    \lref{eq4.4};
    \lref{setting4};
    \lref{setting4.1};
    \lref{setting5};
    \lref{setting6};
  }{
    that for all $J\in\mathcal H$ it holds that
    \begin{eqn}
      \llabel{eq5}
      a\leq\sigma_J<\tau_J\leq b,
      \qquad r(\sigma_J)=S,
      \qquad r(\tau_J)\leq\frac{S}{2},
      \qquad\text{and}\qquad (\sigma_J,\tau_J]\subseteq J
      \dott
    \end{eqn}
  }
  \argument{
    \lref{eq5};
    the fact that $r|_{[a,b]}$ is continuous;
    the intermediate value theorem;
  }{
    that for all $J\in\mathcal H$ it holds that
    $\{t\in[\sigma_J,\tau_J]\colon r(t)=\nicefrac{S}{2}\}$
    is nonempty and compact\dott
  }
  In the following let $\mathfrak b_J\in[a,b]$, $J\in\mathcal H$,
  satisfy for all $J\in\mathcal H$ that
  \begin{eqn}
    \llabel{setting7}
    \mathfrak b_J
      =\inf\{t\in[\sigma_J,\tau_J]
        \colon r(t)=\nicefrac{S}{2}\}
    \dott
  \end{eqn}
  \cfload[.]
  \startnewargseq
  \argument{
    \lref{eq5};
    \lref{setting7};
    the fact that $r|_{[a,b]}$ is continuous;
    the intermediate value theorem;
  }{
    that for all $J\in\mathcal H$ it holds that
    \begin{eqn}
      \llabel{eq6}
      \sigma_J<\mathfrak b_J\leq\tau_J,
      \qquad r(\mathfrak b_J)=\nicefrac{S}{2},
      \qquad
      [\sigma_J,\mathfrak b_J]
        \subseteq r^{-1}([\nicefrac{S}{2},S]),
      \qquad\text{and}\qquad (\sigma_J,\mathfrak b_J)\subseteq J
      \dott
    \end{eqn}
  }
  \argument{
    \lref{eq6};
    \lref{eq5};
    the fact that $r|_{[a,b]}$ is absolutely continuous;
    the fundamental theorem of calculus for absolutely continuous functions;
    the triangle inequality;
  }{
    that for all $J\in\mathcal H$ it holds that
    \begin{eqn}
      \llabel{eq7}
      \frac{S}{2}
        =\abs{r(\mathfrak b_J)-r(\sigma_J)}
        =\abs[\bigg]{\int_{\sigma_J}^{\mathfrak b_J}
          r'(t)\dd t}
        \leq\int_{\sigma_J}^{\mathfrak b_J}\abs{r'(t)}\dd t
      \dott
    \end{eqn}
  }
  \argument{
    \lref{eq7};
    \lref{eq2};
    \lref{setting1};
    \lref{eq6};
    the fact that for all $J,V\in\mathcal H$ with $J\neq V$
    it holds that $J\cap V=\emptyset$;
  }{
    that for every finite set $\mathcal F\subseteq\mathcal H$ it holds that
    \begin{eqn}[a]
      \llabel{eq8}
      \frac{S}{2}(\#\mathcal F)
        =\sum_{J\in\mathcal F}\frac{S}{2}
        \leq\sum_{J\in\mathcal F}
          \int_{\sigma_J}^{\mathfrak b_J}\abs{r'(t)}\dd t\leq\int_{r^{-1}([\nicefrac{S}{2},S])}
          \abs{r'(t)}\dd t
        \leq\mathfrak C S
      \dott
    \end{eqn}
  }
  \argument{
    \lref{eq8};
    \lref{eq4};
    the fact that $S>0$;
  }{
    that
    \begin{eqn}
      \llabel{eq9}
      \#\mathcal H\leq2\mathfrak C
      \qquad\text{and}\qquad
      \#\mathscr J
        =\#\mathcal H+\#(\mathscr J\backslash\mathcal H)
        \leq2\mathfrak C+1
      \dott
    \end{eqn}
  }
  \argument{
    \lref{setting4};
    \lref{setting4.1};
    the fact that $I$ is an interval;
  }{
    that for all $J\in\mathscr J$ it holds that $J$ is a nonempty
    interval and that
    \begin{eqn}
      \llabel{eq9.1}
      J\subseteq U\subseteq[a,b]\subseteq I
      \dott
    \end{eqn}
  }
  \argument{
    \lref{eq9.1};
    \lref{setting1};
    \lref{setting4};
  }{
    that for all $J\in\mathscr J$ and all $t\in J$ it holds that
    \begin{eqn}
      \llabel{eq9.2}
      0<(r|_J)(t)=r(t)<S\leq\delta
      \dott
    \end{eqn}
  }
  \argument{
    \lref{setting3};
    \lref{eq9.1};
  }{
    that for all $J\in\mathscr J$ it holds that
    for Lebesgue-almost all $t\in J$
    it holds that $r$, $x$, $r|_J$, and $x|_J$
    are differentiable at $t$ and that
    \begin{eqn}
      \llabel{eq9.3}
      (r|_J)'(t)=r'(t)
      \qquad\text{and}\qquad
      (x|_J)'(t)=x'(t)
      \dott
    \end{eqn}
  }
  \argument{
    \lref{setting2};
    \lref{eq9.1};
  }{
    that for all $J\in\mathscr J$ it holds that
    for Lebesgue-almost all $t\in J$ it holds that
    for all $v\in T_{(x|_J)(t)}M$ it holds that
    \begin{eqn}
      \llabel{eq9.4}
      (A|_J)(t)v\in T_{(x|_J)(t)}M,\qquad
      \langle (A|_J)(t)v,v\rangle\geq\mathsf a\norm{v}^2,
      \qquad\text{and}\qquad
      \norm{(A|_J)(t)v}\leq\mathsf b\norm{v}
      \dott
    \end{eqn}
  }
  \argument{
    \lref{setting3};
    \lref{eq9.3};
  }{
    that for all $J\in\mathscr J$ it holds that
    for Lebesgue-almost all $t\in J$
    it holds that 
    \begin{gather}
      \llabel{eq9.5}
      (r|_J)'(t)
        =-(\alpha|_J)(t)\bigl((r|_J)(t)\bigr)^2
          (\partial_r h)
            ((r|_J)(t),(x|_J)(t))\\
      \text{and}\qquad(x|_J)'(t)
        =-(\alpha|_J)(t)(A|_J)(t)
          (\nabla_Mh)((r|_J)(t),(x|_J)(t))
      \dott
    \end{gather}
  }
  \argument{
    \lref{eq9.5};
    \lref{eq1}
    (applied for every $J\in\mathscr J$ with
    $I\curvearrowleft J$, $r\curvearrowleft r|_J$,
    $x\curvearrowleft x|_J$, $\alpha\curvearrowleft\alpha|_J$,
    $A\curvearrowleft A|_J$);
    \lref{eq9.1};
    \lref{eq9.2};
    \lref{eq9.4};
    the fact that for all $J\in\mathscr J$ it holds that $r|_J$ and
    $x|_J$ are locally absolutely continuous;
    the fact that for all $J\in\mathscr J$ it holds that $\alpha|_J$
    and $A|_J$ are Lebesgue measurable;
  }{
    that for all $J\in\mathscr J$ it holds that
    \begin{eqn}
      \llabel{eq9.6}
      \int_J\norm{(x|_J)'(t)}\dd t
        \leq\int_J\sqrt{\abs{(r|_J)'(t)}^2
          +\norm{(x|_J)'(t)}^2}\dd t
        \leq K
      \dott
    \end{eqn}
  }
  \argument{
    \lref{eq9.6};
    \lref{eq9.3};
  }{
    that for all $J\in\mathscr J$ it holds that
    \begin{eqn}
      \llabel{eq10}
      \int_J\norm{x'(t)}\dd t\leq K
      \dott
    \end{eqn}
  }
  \argument{
    \lref{eq10};
    \lref{setting1};
    \lref{setting4};
    \lref{eq4};
    \lref{eq9};
  }{
    that
    \begin{eqn}[a]
      \llabel{eq11}
      \int_{[a,b]\cap r^{-1}((0,\nicefrac{S}{2}])}
          \norm{x'(t)}\dd t
        =\int_E\norm{x'(t)}\dd t
        \leq\sum_{J\in\mathscr J}\int_J\norm{x'(t)}\dd t\leq K(\#\mathscr J)
        \leq(1+2\mathfrak C)K=L
      \dott
    \end{eqn}
  }
  \argument{
    \lref{eq11};
    the fact that $I$, $r$, $x$, $\alpha$, $A$, $a$, and $b$ were
    arbitrary subject to the stated assumptions;
  }[verbs=ep]{
    \lref{eq3}\dott
  }

  \startnewargseq
  In the following let $R=S/2$\dott \cfload[.]
  \startnewargseq
  \argument{
    \lref{eq3};
    the fact that $R=S/2\in(0,\mathfrak R)$;
    the monotone convergence theorem;
  }[verbs=ep]{
    \lref{concl1}\dott
  }
\end{aproof}

\section{Spherical length at finite points and at infinity}
\label{sec:applications}

In this section we prove in \cref{maintheorem3} in
\cref{subsec:general-time-intervals} below that on an arbitrary time interval
with non-empty interior $I\subseteq\R$ we have that the spherical projection of every
\GF\ trajectory $\theta(t)$, $t\in I$, of a definable
$C^1$-potential $F\colon U\to\R$ defined on an open subset $U$ of
$\R^n$ has finite length.
\Cref{maintheorem} in the introduction essentially coincides with
\cref{maintheorem3}.

In the first step of our proof of \cref{maintheorem3} we reduce
\cref{maintheorem3} (by separately considering the \GF\ trajectories
$t\mapsto\theta(t)$ and $t\mapsto\theta(-t)$ of the potentials $F$ and
$-F$, respectively) to the special case where the time
interval $I$ is equal to the half-open interval $I=[\tau,T)$ for some
$T\in(\tau,\infty]$, $\tau\in\R$.
This special case is proved as \cref{maintheorem2} in
\cref{subsec:open-domain-length} below.

\Cref{maintheorem2} is proved by reducing it to the special case where the
\GF\ trajectory $\theta(t)$, $t\in[\tau,T)$, never hits the reference point
$\vartheta$ for the spherical projection in the sense that
$\{t\in[\tau,T)\colon\theta(t)=\vartheta\}=\emptyset$.
This special case is proved as \cref{prop:spherical-length-away} in
\cref{subsec:open-domain-length}.
Loosely speaking, the proof of \cref{prop:spherical-length-away} can be divided into three cases, namely,
\begin{enumerate}[label=(\Alph*)]
\item
\label{item:i}
the case in which there exists $ v \in \R^n \backslash \{ \vartheta \} $ such that
$
  \limsup_{ t \to T }
  | \theta( t ) - v |
  = 0
$
(convergence to a finite point $ v \in \R^n $ that is different from the reference point $ \vartheta $),
\item
\label{item:ii}
the case in which
$
  \limsup_{ t \to T }
  | \theta( t ) - \vartheta |
  = 0
$
(convergence to the reference point $ \vartheta \in \R^n $) -- the case of the \emph{definable Thom gradient conjecture},
and
\item
\label{item:iii}
the case in which $ \liminf_{ t \to T } | \theta(t) | = \infty $ (divergence to infinity) -- the case of the \emph{gradient conjecture at infinity}.
\end{enumerate}
The case in \cref{item:i} is straightforward as the spherical length can be controlled through the length of the \GF\ trajectory
and is handled here by \itref{maintheorem1}{item4} in \cref{subsec:finite-and-infinite-limits} below
(cf.\ also \cref{prop:annular-length} in \cref{subsec:euclidean-occupation} below).
The key difficulty lies in the cases in \cref{item:ii,item:iii}.

A central step in the proof of \cref{prop:spherical-length-away} in the case of \cref{item:ii,item:iii} is to
perform an appropriate coordinate transform to the \GF\ \DE\
to obtain a transformed \DE\ system that consists of a radial variable $r(t)$
(which is $r(t)=\norm{\theta(t)-\vartheta}^\nu$ with $\nu=1$ in the case of
\cref{item:ii} and $\nu=-1$ in the case of \cref{item:iii}) and the
angular variable $\omega(t)$ (which is the spherical
projection) that takes values in a regular submanifold of $\R^n$
(cf.\ \cref{definition:submanifold}), namely, in the setting of \cref{sec:applications}, the
sphere in $\R^n$.
The transformation is thus nothing but polar coordinates in the case
of \cref{item:ii} and polar coordinates combined with the inversion method
($\R^n\backslash\{0\}\ni x\mapsto\norm{x}^{-2}x\in\R^n\backslash\{0\}$)
in the case of \cref{item:iii}, even though the transformed \DE\ system appears
slightly more complicated.

The resulting \DE\ system is precisely in the form required by the machinery in
\cref{sec:radial-integration} -- specifically, \itref{cor:angular-length}{item2} --
to establish finite length of the angular variable on the submanifold and,
therefore, to establish finite spherical length of the original \GF\ trajectories.
A useful feature of the transformation and the resulting \DE\ system
is that after the transformation the radial variable always converges to
zero as both cases, \cref{item:ii} (convergence to the reference point
$\vartheta$) and \cref{item:iii} (divergence to infinity), are
handled in the regime where the transformed radial variable approaches zero.

\subsection{Uniform upper bounds for the spherical length of GF trajectories close to the origin}
\label{subsec:euclidean-occupation}

In the following \cref{prop:annular-length} we establish a \emph{uniform} upper
bound for the length of the spherical projection of a \GF\ trajectory
where we count only those time periods in which the \GF\ trajectory is in a sufficiently
small ball around the origin.
The potential $F\colon V\to\R$ in \cref{prop:annular-length} is assumed to be
continuously differentiable and definable.

\begin{athm}{corollary}{prop:annular-length}
  Let $n\in\N\backslash\{1\}$, $\mathfrak R\in(0,\infty)$,
  let $V=\{z\in\R^n\colon0<\norm{z}<\mathfrak R\}$, and let
  $F\in C^1(V,\R)$ be definable\cfadd{definition:definable-function}\dott
  \cfload[.]
  Then there exist $\varepsilon\in(0,\mathfrak R)$,
  $\mathfrak{C}\in\R$ such that for every interval $I\subseteq\R$,
  every locally absolutely continuous $\theta\colon I\to V$,
  every Lebesgue measurable $\lambda\colon I\to[0,\infty)$,
  and every $\omega\in C(I,\R^n)$ with the property that for all $t\in I$ it holds that
 $\omega(t)=\norm{\theta(t)}^{-1}\theta(t)$
  and 
    $\int_{I}|\theta'(s)+\lambda(s)(\grad F)(\theta(s))|\dd s=0$
  we have that
\begin{enumerate}[label=(\roman*)]
    \item it holds that $ \omega $ is locally absolutely continuous and
    \item it holds for all $ \rho \in (0,\varepsilon] $ that
    \begin{eqn}
    \llabel{concl1}
    \int_{\{t\in I\colon0<\norm{\theta(t)}\leq\rho\}}
      \norm{\theta'(t)}\dd t\leq\mathfrak{C}\rho
    \qquad\text{and}\qquad
    \int_{\{t\in I\colon0<\norm{\theta(t)}\leq\varepsilon\}}
      \norm{\omega'(t)}\dd t\leq\mathfrak{C}
    \dott
  \end{eqn}
 \end{enumerate}
  \cfout[.]
\end{athm}
\vspace{3pt}
\begin{aproof}
  Throughout this proof let $M\subseteq\R^n$ satisfy
  \begin{eqn}
    \llabel{setting1}
    M=\{u\in\R^n\colon\norm{u}=1\},
  \end{eqn}
  and let $h\colon(0,\mathfrak R)\times M\to\R$ satisfy for all
  $s\in(0,\mathfrak R)$, $u\in M$ that
  \begin{eqn}
    \llabel{setting2}
    h(s,u)=F(su)
    \dott
  \end{eqn}
  \cfload[.]
  \startnewargseq
  \argument{
    \lref{setting1};
    \lref{setting2};
    the fact that $n\geq2$;
    the assumption that $F\in C^1(V,\R)$ is definable;
  }{
    that $M$ is a nonempty\cfadd{definition:definable-set}, compact\cfadd{definition:embedded-manifold}, connected, and definable
    \dsubmanifold{2}{n}{(n-1)} and that\cfadd{definition:c1-on-manifold}
    $h\in \COne{(0,\mathfrak R)\times M}{\R}$ is definable\dott
    \llabel{eq1}
  }
  \argument{
    \lref{eq1};
    \cref{prop:cone-annular-length}
    (applied with $D\curvearrowleft n$,
    $\Omega\curvearrowleft(0,\mathfrak R)\times M$,
    $\mathsf a\curvearrowleft1$, $\mathsf b\curvearrowleft1$,
    $A\curvearrowleft(I\ni t\mapsto\mathrm{I}_n\in\R^{n\times n})$
    in the notation of \cref{prop:cone-annular-length});
    the monotonicity of the Lebesgue integral;
    \cref{prop:angular-occupation}
    (applied with $D\curvearrowleft n$, $n\curvearrowleft n-1$,
    $\mathsf a\curvearrowleft1$, $\mathsf b\curvearrowleft1$,
    $R\curvearrowleft\varepsilon$,
    $A\curvearrowleft(I\ni t\mapsto\mathrm{I}_n\in\R^{n\times n})$ in the notation of
    \cref{prop:angular-occupation});
  }{
    that there exist $\varepsilon\in(0,\mathfrak R)$,
    $K,L\in(0,\infty)$ which satisfy for every
    interval $I\subseteq\R$, all locally absolutely continuous
    $r\colon I\to(0,\mathfrak R)$ and $x\colon I\to M$,
    and every Lebesgue measurable 
    $\alpha\colon I\to[0,\infty)$
    with the property that for Lebesgue-almost all $t\in I$
    it holds that 
    \begin{eqn}
      \llabel{eq2}
      r'(t)=-\alpha(t)(r(t))^2
        (\partial_r h)(r(t),x(t))
      \qquad\text{and}\qquad
      x'(t)=-\alpha(t)(\nabla_Mh)(r(t),x(t)),
    \end{eqn}
    \cfadd{definition:angular-gradient}%
    we have for all $\rho\in(0,\mathfrak R)$ that
    \begin{eqn}
      \llabel{eq3}
      \int_{r^{-1}([\rho/2,\rho])}\abs{r'(t)}\dd t
        \leq K\rho
      \qquad\text{and}\qquad
      \int_{r^{-1}((0,\varepsilon])}\norm{x'(t)}\dd t\leq L
      \dott
    \end{eqn}
  }
  In the following let $\mathfrak{C}\in(0,\infty)$ satisfy
  \begin{eqn}
    \llabel{setting3}
    \mathfrak{C}=2K+L
    \dott
  \end{eqn}
  \startnewargseq
  We claim that for every interval $I\subseteq\R$,
  every locally absolutely continuous $y\colon I\to V$,
  every Lebesgue measurable $\lambda\colon I\to[0,\infty)$,
  and every $\omega\colon I\to\R^n$
  with the property that for Lebesgue-almost all $t\in I$ it
  holds that 
  \begin{eqn}
    y'(t)=-\lambda(t)(\grad F)(y(t)),
  \end{eqn}
  and with the property that for all $t\in I$ it holds that
  \begin{equation}
    \omega(t)=\frac{y(t)}{\norm{y(t)}}
  \end{equation}
  we have that $\omega$ is locally absolutely continuous and
  that for all $\rho\in(0,\varepsilon]$ it holds that
  \begin{eqn}
    \llabel{eq5.1}
    \int_{\{t\in I\colon0<\norm{y(t)}\leq\rho\}}
      \norm{y'(t)}\dd t\leq(2K+L)\rho
    \qquad\text{and}\qquad
    \int_{\{t\in I\colon0<\norm{y(t)}\leq\varepsilon\}}
      \norm{\omega'(t)}\dd t\leq L
    \dott
  \end{eqn}
  \startnewargseq
  To prove this claim, let $I\subseteq\R$ be an interval,
  let $\rho\in(0,\varepsilon]$,
  let $y\colon I\to V$ be locally absolutely continuous,
  let $\lambda\colon I\to[0,\infty)$ be Lebesgue measurable,
  assume that for Lebesgue-almost all $t\in I$ it holds that
  \begin{eqn}
    \llabel{setting3.1}
    y'(t)=-\lambda(t)(\grad F)(y(t)),
  \end{eqn}
  and let $\omega\colon I\to\R^n$ satisfy for all $t\in I$ that
  \begin{eqn}
    \llabel{setting3.2}
    \omega(t)=\frac{y(t)}{\norm{y(t)}}
    \dott
  \end{eqn}
  \startnewargseq
  In the following let $r\colon I\to(0,\mathfrak R)$ and
  $\alpha\colon I\to[0,\infty)$ satisfy for all $t\in I$ that
  \begin{eqn}
    \llabel{setting4}
    r(t)=\norm{y(t)}
    \qquad\text{and}\qquad
    \alpha(t)=\frac{\lambda(t)}{(r(t))^2}
    \dott
  \end{eqn}
  \startnewargseq
  \argument{
    \lref{setting1};
    \lref{setting3.2};
    \lref{setting4};
    the fact that $y$ is locally absolutely continuous;
    the assumption that $\lambda$ is Lebesgue measurable;
  }{
    that $\omega(I)\subseteq M$, that $\alpha$ is Lebesgue
    measurable, and that $r$ and $\omega$ are locally absolutely
    continuous\dott
  }
  \argument{
    \lref{setting2};
    the chain rule;
  }{
    that for all $s\in(0,\mathfrak R)$, $u\in M$ it holds that
    \begin{eqn}
      \llabel{eq6.1}
      (\partial_r h)(s,u)
        =\langle(\grad F)(su),u\rangle
      \dott
    \end{eqn}
  }
  \argument{
    \lref{setting1};
    \lref{setting2};
    the chain rule;
  }{
    that for all $s\in(0,\mathfrak R)$, $u\in M$,
    $\eta\in(0,\infty)$, $\gamma\in C^1((-\eta,\eta),M)$ with $\gamma(0)=u$ it holds that
    \begin{equation}
      \langle u,\gamma'(0)\rangle
        =\frac12\frac{\ddd}{\ddd t}\norm{\gamma(t)}^2\big|_{t=0}=0
    \end{equation}
    and
    \begin{equation}\llabel{eq6.2}
      \frac{\ddd}{\ddd t}h(s,\gamma(t))\big|_{t=0}
        =\frac{\ddd}{\ddd t}F(s\gamma(t))\big|_{t=0}
         =s\langle\bigl(\grad F\bigr)\bigl(su\bigr),\gamma'(0)\rangle
      \dott
    \end{equation}
  }
  \argument{
    \lref{eq6.2};
    \lref{setting1};
    the fact that for all $u\in M$ it holds that
    $T_uM=\{v\in\R^n\colon\langle u,v\rangle=0\}$\cfadd{definition:tangent-space};
    the fact that for all $u\in M$, $w\in\R^n$ it holds that
    $\langle w-\langle w,u\rangle u,u\rangle
      =\langle w,u\rangle(1-\norm{u}^2)=0$;
  }{
    that for all $s\in(0,\mathfrak R)$, $u\in M$ it holds that
    \begin{eqn}
      \llabel{eq6.3}
      s((\grad F)(su)-\langle(\grad F)(su),u\rangle u)
        \in T_uM
    \end{eqn}
    and that for all $s\in(0,\mathfrak R)$, $u\in M$,
    $\eta\in(0,\infty)$, $\gamma\in C^1((-\eta,\eta),M)$ with $\gamma(0)=u$ it holds that
    \begin{eqn}[a]
      \llabel{eq6.4}
      &\bigl\langle
        s\bigl(\bigl(\grad F\bigr)\bigl(su\bigr)-\langle\bigl(\grad F\bigr)\bigl(su\bigr),u\rangle u\bigr),
        \gamma'(0)\bigr\rangle\\
      &=s\langle\bigl(\grad F\bigr)\bigl(su\bigr),\gamma'(0)\rangle
        -s\langle\bigl(\grad F\bigr)\bigl(su\bigr),u\rangle
          \langle u,\gamma'(0)\rangle
        =\frac{\ddd}{\ddd t}h(s,\gamma(t))\big|_{t=0}
      \dott
    \end{eqn}
  }
  \argument{
    \lref{eq6.4};
    \lref{eq6.1};
    \cref{lem:angular-gradient-existence};
  }{
    that for all $s\in(0,\mathfrak R)$, $u\in M$ it holds that
    \begin{eqn}
      \llabel{eq7}
      (\partial_r h)(s,u)
        =\langle\bigl(\grad F\bigr)\bigl(su\bigr),u\rangle
      \qquad\text{and}\qquad
      (\nabla_Mh)(s,u)
        =s\bigl(\bigl(\grad F\bigr)\bigl(su\bigr)-\langle\bigl(\grad F\bigr)\bigl(su\bigr),u\rangle u\bigr)
      \dott
    \end{eqn}
  }
  \argument{
    \lref{setting3.2};
    \lref{setting4};
    the fact that $y$ is locally absolutely continuous;
    the chain rule;
  }{
    that for Lebesgue-almost all $t\in I$ it holds that
    \begin{eqn}
      \llabel{eq8}
      r'(t)=\langle\omega(t),y'(t)\rangle
      \qquad\text{and}\qquad
      \omega'(t)=\frac{y'(t)-r'(t)\omega(t)}{r(t)}
      \dott
    \end{eqn}
  }
  \argument{
    \lref{setting3.2};
    \lref{setting4};
  }{
    that for all $t\in I$ it holds that
    \begin{eqn}
      \llabel{eq8.1}
      y(t)=r(t)\omega(t)
      \qquad\text{and}\qquad
      \lambda(t)=\alpha(t)(r(t))^2
      \dott
    \end{eqn}
  }
  \argument{
    \lref{eq8.1};
    \lref{setting3.1};
    \lref{eq7};
    \lref{eq8};
  }{
    that for Lebesgue-almost all $t\in I$ it holds that
    \begin{eqn}[a]
      \llabel{eq9}
      r'(t)
        &=\langle\omega(t),y'(t)\rangle
         =-\lambda(t)\langle\bigl(\grad F\bigr)\bigl(y(t)\bigr),\omega(t)\rangle=-\lambda(t)
          \langle\bigl(\grad F\bigr)\bigl(r(t)\omega(t)\bigr),\omega(t)\rangle\\
          &=-\lambda(t)(\partial_r h)(r(t),\omega(t))
         =-\alpha(t)(r(t))^2
           (\partial_r h)(r(t),\omega(t))
      \dott
    \end{eqn}
  }
  \argument{
    \lref{eq9};
    \lref{setting3.1};
    \lref{eq7};
    \lref{eq8.1};
    the fact that for all $t\in I$ it holds that $r(t)>0$;
  }{
    that for Lebesgue-almost all $t\in I$ it holds that
    \begin{eqn}[a]
      \llabel{eq9.1}
      y'(t)-r'(t)\omega(t)
        &=-\lambda(t)\bigl(\bigl(\grad F\bigr)\bigl(y(t)\bigr)
          -\langle\bigl(\grad F\bigr)\bigl(y(t)\bigr),\omega(t)\rangle\omega(t)\bigr)=-\frac{\lambda(t)}{r(t)}
          (\nabla_Mh)(r(t),\omega(t))
      \dott
    \end{eqn}
  }
  \argument{
    \lref{eq9.1};
    \lref{setting4};
    \lref{eq8};
  }{
    that for Lebesgue-almost all $t\in I$ it holds that
    \begin{eqn}[a]
      \llabel{eq9.2}
      \omega'(t)
        =\frac{y'(t)-r'(t)\omega(t)}{r(t)}=-\frac{\lambda(t)}{(r(t))^2}
          (\nabla_Mh)(r(t),\omega(t))
         =-\alpha(t)(\nabla_Mh)(r(t),\omega(t))
      \dott
    \end{eqn}
  }
  \argument{
    \lref{eq8};
    the fact that for all $t\in I$ it holds that $r(t)>0$;
  }{
    that for Lebesgue-almost all $t\in I$ it holds that
    \begin{eqn}
      \llabel{eq10.1}
      y'(t)=r'(t)\omega(t)+r(t)\omega'(t)
      \dott
    \end{eqn}
  }
  \argument{
    \lref{eq3};
    the fact that $\omega(I)\subseteq M$;
    the fact that $\alpha$ is Lebesgue measurable;
    the fact that $r$ and $\omega$ are locally absolutely continuous;
    \lref{setting4};
    \lref{eq9};
    \lref{eq9.2};
    \itref{lem:v}{item1} (applied with
    $\varrho\curvearrowleft\mathfrak R$, $\mathfrak C\curvearrowleft K$);
  }{
    that
    \begin{eqn}
      \llabel{eq10.3}
      \int_{r^{-1}((0,\rho])}\abs{r'(t)}\dd t\leq2K\rho
      \qquad\text{and}\qquad
      \int_{r^{-1}((0,\varepsilon])}\norm{\omega'(t)}\dd t\leq L
      \dott
    \end{eqn}
  }
  \argument{
    \lref{eq10.3};
    \lref{setting1};
    \lref{setting4};
     \lref{eq10.1};
    the fact that $\omega(I)\subseteq M$;
    the fact that $0<\rho\leq\varepsilon$;
    the triangle inequality;
  }{
    that
    \begin{eqn}[a]
      \llabel{eq11}
      \int_{\{t\in I\colon0<\norm{y(t)}\leq\rho\}}\norm{y'(t)}\dd t
        &=\int_{r^{-1}((0,\rho])}\norm{y'(t)}\dd t\\
        &\leq\int_{r^{-1}((0,\rho])}\abs{r'(t)}\dd t
          +\rho\int_{r^{-1}((0,\varepsilon])}\norm{\omega'(t)}\dd t
        \leq(2K+L)\rho
      \dott
    \end{eqn}
  }
  \argument{
    \lref{eq11};
    \lref{setting4}; \lref{eq10.3};
    the fact that $\omega$ is locally absolutely continuous;
    the fact that $I$, $\rho$, $y$, $\lambda$, and
    $\omega$ were arbitrary subject to the stated assumptions;
  }[verbs=ep]{
    \lref{eq5.1}\dott
  }

  \startnewargseq
  \argument{\lref{eq3};
    \lref{setting3};
    \lref{eq5.1};
    the fact that $K\in(0,\infty)$;
    the fact that $L\in[0,\infty)$;
  }[verbs=ep]{
    \lref{concl1}\dott
  }
\end{aproof}

\subsection{Finite spherical length for GF trajectories converging to finite points}
\label{subsec:finite-and-infinite-limits}

In the proofs of our results on finite spherical length in
\cref{maintheorem1}, \cref{prop:spherical-length-away}, and \cref{maintheorem2} below,
we use
the well-known fact, which we recall in the following
\cref{lem:curve-length-speed}, that the length of a locally absolutely
continuous curve coincides with the integral of the norm of the (almost everywhere defined) derivative
of the curve.
Only for the sake of completeness, we include a detailed proof for
\cref{lem:curve-length-speed} here.

\begin{athm}{lemma}{lem:curve-length-speed}
  Let $n\in\N$, let $I\subseteq\R$ be an interval, and let
  $\gamma\colon I\to\R^n$ be locally absolutely continuous \cfload.
  Then
  \begin{eqn}
    \llabel{concl1}
    \curvelength{\gamma}=\int_I\norm{\gamma'(t)}\dd t
    \dott
  \end{eqn}
 \cfout[.]
\end{athm}

\begin{aproof}
  \argument{
    the fact that $\gamma$ is locally absolutely continuous;
    the fundamental theorem of calculus for absolutely continuous functions;
  }{
    that for all $k\in\N$, $t_0,\dots,t_k\in I$ with
    $t_0<t_1<\ldots<t_k$ it holds that
    \begin{eqn}[a]
      \llabel{eq1}
      \sum_{j=1}^k\norm{\gamma(t_j)-\gamma(t_{j-1})}
        &=\sum_{j=1}^k
          \norm[\bigg]{\int_{t_{j-1}}^{t_j}\gamma'(u)\dd u}
          \leq\sum_{j=1}^k
          \int_{t_{j-1}}^{t_j}\norm{\gamma'(u)}\dd u\\
          &
          =\int_{t_0}^{t_k}\norm{\gamma'(u)}\dd u
          \leq\int_I\norm{\gamma'(u)}\dd u
      \dott
    \end{eqn}
  }
  \argument{
    \lref{eq1};
  }{
    that
    \begin{eqn}
      \llabel{eq2}
      \curvelength{\gamma}\leq\int_I\norm{\gamma'(t)}\dd t
      \dott
    \end{eqn}
  }
  We claim that for all $a,b\in I$ with $a<b$ it holds that
  \begin{eqn}
    \llabel{claim1}
    \int_a^b\norm{\gamma'(t)}\dd t\leq\curvelength{\gamma}
    \dott
  \end{eqn}
  \startnewargseq
  To prove \lref{claim1}, let $a,b\in I$ satisfy $a<b$\dott
  \startnewargseq
  In the following let $q_{k,j}\in[a,b]$,
  $k\in\N$, $j\in\{0,1,\dots,k\}$, satisfy for all
  $k\in\N$, $j\in\{0,1,\dots,k\}$ that
  \begin{eqn}
    \llabel{setting1}
    q_{k,j}=a+\frac{j}{k}(b-a)
  \end{eqn}
  and let $g_k\colon[a,b)\to[0,\infty)$, $k\in\N$, satisfy
  for all $k\in\N$, $j\in\{0,1,\dots,k-1\}$,
  $t\in[q_{k,j},q_{k,j+1})$ that
  \begin{eqn}
    \llabel{setting2}
    g_k(t)=\frac{\norm{\gamma(q_{k,j+1})-\gamma(q_{k,j})}}
      {q_{k,j+1}-q_{k,j}}
    \dott
  \end{eqn}
  \startnewargseq
  \argument{
    \lref{setting1};
    \lref{setting2};
    the fact that $[a,b]\subseteq I$;
  }{
    that for all $k\in\N$ it holds that
    \begin{eqn}
      \llabel{eq3}
      \int_{[a,b)}g_k(t)\dd t
        =\sum_{j=0}^{k-1}
          \norm{\gamma(q_{k,j+1})-\gamma(q_{k,j})}
        \leq\curvelength{\gamma}
      \dott
    \end{eqn}
  }
  \argument{
    \lref{setting1};
    \lref{setting2};
    the fact that $\gamma|_{[a,b]}$ is absolutely continuous;
    the fact that for Lebesgue-almost all $t\in(a,b)$ it holds
    that $\gamma$ is differentiable at $t$;
  }{
    that for Lebesgue-almost all $t\in(a,b)$ it holds that
    \begin{eqn}
      \llabel{eq4}
      \lim_{k\to\infty}g_k(t)=\norm{\gamma'(t)}
      \dott
    \end{eqn}
  }
  \argument{
    \lref{eq4};
    \lref{eq3};
    the fact that for all $k\in\N$ it holds that $g_k$ is
    nonnegative and Lebesgue measurable;
    Fatou's lemma;
  }{
    that
    \begin{eqn}
      \llabel{eq5}
      \int_a^b\norm{\gamma'(t)}\dd t
        \leq\liminf_{k\to\infty}\int_{[a,b)}g_k(t)\dd t
        \leq\curvelength{\gamma}
      \dott
    \end{eqn}
  }
  \argument{
    \lref{eq5};
    the fact that $a,b\in I$ with $a<b$ were arbitrary;
  }[verbs=ep]{
    \lref{claim1}\dott
  }

  \startnewargseq
  \argument{
    \lref{claim1};
    the fact that $I$ is an interval;
    the fact that $\norm{\gamma'}$ is nonnegative and Lebesgue measurable;
    the monotone convergence theorem;
  }{
    that
    \begin{eqn}[a]
      \llabel{eq6}
      \int_I\norm{\gamma'(t)}\dd t
        =\sup\biggl(\{0\}\cup
          \biggl\{\int_c^d\norm{\gamma'(t)}\dd t\colon
            c,d\in I,\ c<d\biggr\}\biggr)
      \leq\curvelength{\gamma}
      \dott
    \end{eqn}
  }
  \argument{
    \lref{eq6};
    \lref{eq2};
  }[verbs=ep]{
    that
    \begin{eqn}
      \curvelength{\gamma}=\int_I\norm{\gamma'(t)}\dd t
      \dott
    \end{eqn}
  }
\end{aproof}

\begin{athm}{proposition}{maintheorem1}
Let $ n \in \mathbb{N} $,
$ \vartheta, p \in \R^n $,
$ \tau \in \R $, $ T \in ( \tau, \infty ] $,
let $ U \subseteq \R^n $ be definable and open,
let $ F \in C^1( U, \R ) $ be definable,
let
$ \theta \in C^1( [\tau,T), U ) $ satisfy
\begin{eqn}
  \llabel{ass1}
  \theta' = - (\nabla F)\circ\theta,
\end{eqn}
let
$ \omega \colon [ \tau, T ) \to \R^n $
satisfy for all $ t \in [ \tau, T ) $ with $ \theta( t ) \neq p $ that
$
  \omega( t ) = | \theta( t ) - p |^{ - 1 } ( \theta(t) - p )
$,
and let
$
  \mathfrak{t} =
  \sup( \{ -\infty \} \cup \theta^{ - 1 }( \{ p \} ) )
$\cfload.
Then the following statements hold:
\begin{enumerate}[label=(\roman*)]
\item \llabel{item2}
If $ \liminf_{ t \to T } | \theta(t) | < \infty $,
then there exists $ v \in \R^n $ such that
$
  \limsup_{ t \to T }
  | \theta(t) - v | = 0
$.

\item \llabel{item4}
If
$ \liminf_{ t \to T } | \theta(t) - \vartheta | = 0 $
and
there exists an open $ O \subseteq \R^n $
such that
$ ( \{ \vartheta \} \cap \{ p \} ) \subseteq O \subseteq ( U \cup \{ \vartheta \} ) $,
then
for every
$
  s \in [ \tau, T )
  \cap
  ( \mathfrak{t}, \infty )
$
the function $ \omega|_{ [s, T) } $
has finite length\cfadd{definition:curve-length}\dott

\end{enumerate}
\cfout[.]
\end{athm}

\begin{aproof}
  Throughout this proof let $J\subseteq[\tau,T)$ satisfy
  \begin{eqn}
    \llabel{setting1}
    J=[\tau,T)\cap(\mathfrak t,\infty)
    \dott
  \end{eqn}
  \cfload[.]
  \startnewargseq
  \argument{
    \lref{ass1};
    the chain rule;
  }{
    that $F\circ\theta$ is non-increasing and that for all
    $t\in(\tau,T)$ it holds that
    \begin{eqn}
      \llabel{eq1}
      \frac{\ddd}{\ddd t}F(\theta(t))
        =\langle\bigl(\grad F\bigr)\bigl(\theta(t)\bigr),\theta'(t)\rangle
        =-\norm{\bigl(\grad F\bigr)\bigl(\theta(t)\bigr)}^2
        =-\norm{\theta'(t)}^2
      \dott
    \end{eqn}
  }
  \argument{
    \lref{setting1};
    the fact that $
  \mathfrak{t} =
  \sup( \{ -\infty \} \cup \theta^{ - 1 }( \{ p \} ) )
$;
    the assumption that for all $t\in[\tau,T)$ with
    $\theta(t)\neq p$ it holds that
    $\omega(t)=(\theta(t)-p)/\norm{\theta(t)-p}$;
    the assumption that $\theta\in C^1([\tau,T),U)$;
    the chain rule;
  }{
    that for all $s\in J$ it holds that
    $\omega|_{[s,T)}\in C^1([s,T),\R^n)$ and that for all $t\in J$ it
    holds that
    \begin{eqn}
      \llabel{eq2}
      \theta(t)\neq p,
      \qquad \norm{\omega(t)}=1,
      \qquad\text{and}\qquad
      \omega'(t)
        =\frac{\theta'(t)-\langle\omega(t),\theta'(t)\rangle
          \omega(t)}{\norm{\theta(t)-p}}
      \dott
    \end{eqn}
  }
  \argument{
    \lref{eq2};
    \lref{setting1};
  }{
    that for all $t,s,b\in J$ with $s<b$ it holds that  $\min_{v\in[s,b]}\norm{\theta(v)-p}>0$,
    \begin{equation}
      \llabel{eq3}
      \begin{split}
      &\norm{\theta'(t)-\langle\omega(t),\theta'(t)\rangle
        \omega(t)}^2
        =\norm{\theta'(t)}^2
          -2\langle\omega(t),\theta'(t)\rangle
            \langle\theta'(t),\omega(t)\rangle
          +\abs{\langle\omega(t),\theta'(t)\rangle}^2
            \norm{\omega(t)}^2\\
        &=\norm{\theta'(t)}^2
          -2\abs{\langle\omega(t),\theta'(t)\rangle}^2
          +\abs{\langle\omega(t),\theta'(t)\rangle}^2=\norm{\theta'(t)}^2
          -\abs{\langle\omega(t),\theta'(t)\rangle}^2
          \leq\norm{\theta'(t)}^2,
      \end{split}
    \end{equation}
   and 
    \begin{eqn}[a]
      \llabel{eq4}
      \int_s^b\norm{\omega'(t)}\dd t
        &=\int_s^b
          \frac{\norm{\theta'(t)-\langle\omega(t),\theta'(t)\rangle
            \omega(t)}}{\norm{\theta(t)-p}}\dd t
        \leq\int_s^b\frac{\norm{\theta'(t)}}{\norm{\theta(t)-p}}\dd t\\
        &
          \leq\frac{\int_s^b\norm{\theta'(t)}\dd t}
            {\min_{t\in[s,b]}\norm{\theta(t)-p}}\leq\frac{(b-s)\max_{t\in[s,b]}\norm{\theta'(t)}}
          {\min_{t\in[s,b]}\norm{\theta(t)-p}}<\infty
      \dott
    \end{eqn}
  }
  We claim that if
  $\liminf_{t\to T}\norm{\theta(t)}<\infty$, then
  \begin{eqn}
    \llabel{eq12.2}
    \int_\tau^T\norm{\theta'(t)}\dd t<\infty
    \dott
  \end{eqn}
  \startnewargseq
  To prove this claim, assume that
  $\liminf_{t\to T}\norm{\theta(t)}<\infty$\dott
  \startnewargseq
  \argument{
    the assumption that
    $\liminf_{t\to T}\norm{\theta(t)}<\infty$;
  }{
    that there exists $\mathbf{v}\in\R^n$ which satisfies
    \begin{eqn}
      \llabel{eq13}
      \liminf_{t\to T}\norm{\theta(t)-\mathbf{v}}=0
      \dott
    \end{eqn}
  }
  In the following let
  $\chi\colon\R\to(-1,1)$, $H\colon U\to(-1,1)$, and
  $e\colon[\tau,T)\to(-1,1)$ satisfy for all
  $u\in\R$, $z\in U$, $t\in[\tau,T)$ that
  \begin{eqn}
    \llabel{setting4}
    \chi(u)=\frac{u}{\sqrt{1+u^2}},
    \qquad H(z)=\chi(F(z)),
    \qquad\text{and}\qquad e(t)=H(\theta(t))
    \dott
  \end{eqn}
  \startnewargseq
  \argument{
    \lref{setting4};
    the assumption that $F\in C^1(U,\R)$ is definable;
    the chain rule;
  }{
    that $H\in C^1(U,(-1,1))$ is definable and that for all
    $u\in\R$, $z\in U$ it holds that
    \begin{eqn}
      \llabel{eq14}
      \chi'(u)=(1+u^2)^{-3/2}>0
      \qquad\text{and}\qquad
      \bigl(\grad H\bigr)\bigl(z\bigr)=\chi'(F(z))\bigl(\grad F\bigr)\bigl(z\bigr)
      \dott
    \end{eqn}
  }
  \argument{
    \lref{eq14};
    \lref{ass1};
    \lref{eq1};
    \lref{setting4};
    the chain rule;
  }{
    that $e$ is non-increasing and that for all
    $t\in(\tau,T)$ it holds that
    \begin{eqn}
      \llabel{eq15}
      e'(t)
        =-\chi'(F(\theta(t)))\norm{\theta'(t)}^2
        =-\norm{\bigl(\grad H\bigr)\bigl(\theta(t)\bigr)}\norm{\theta'(t)}
        \leq0
      \dott
    \end{eqn}
  }
  \argument{
    \lref{eq15};
    the fact that $e([\tau,T))\subseteq(-1,1)$;
  }{
    that there exists $\ell\in[-1,1)$ which satisfies
    \begin{eqn}
      \llabel{eq16}
      \lim_{t\to T}e(t)=\ell
      \qquad\text{and}\qquad
      \inf_{t\in[\tau,T)}e(t)=\ell
      \dott
    \end{eqn}
  }
  In our proof of \lref{eq12.2} we distinguish between the case
  $e^{-1}(\{\ell\})\neq\emptyset$ and the case that for all
  $t\in[\tau,T)$ it holds that $e(t)>\ell$\dott
  \startnewargseq
  We first prove \lref{eq12.2} in the case
  $e^{-1}(\{\ell\})\neq\emptyset$\dott
  \startnewargseq
  \argument{
    \lref{eq14};
    \lref{eq15};
    \lref{eq16};
    the fundamental theorem of calculus;
  }{
    that for every $a\in[\tau,T)$ with $e(a)=\ell$ it holds
    that $e|_{[a,T)}$ and $\theta|_{[a,T)}$ are constant and that
    \begin{eqn}
      \llabel{eq17}
      \int_\tau^T\norm{\theta'(t)}\dd t
        =\int_\tau^a\norm{\theta'(t)}\dd t<\infty
      \dott
    \end{eqn}
  }
  \argument{
    \lref{eq17};
  }[verbs=ep]{
    \lref{eq12.2} in the case
    $e^{-1}(\{\ell\})\neq\emptyset$\dott
  }
  We now prove \lref{eq12.2} in the case that for all
  $t\in[\tau,T)$ it holds that $e(t)>\ell$\dott
  \startnewargseq
  In the following let
  $W=\bigl\{z\in U\colon\bigl[(\norm{z-\mathbf{v}}<1)
    \wedge(H(z)>\ell)\bigr]\bigr\}$, and let
  $f\colon W\to(0,\infty)$ satisfy for all $z\in W$ that
  \begin{eqn}
    \llabel{setting5}
    f(z)=H(z)-\ell
    \dott
  \end{eqn}
  \startnewargseq
  \argument{
    \lref{setting5};
    the fact that $W$ is open, bounded, and definable;
    the fact that $f\in C^1(W,(0,\infty))$ is positive and definable;
    \cite[Theorem~1]{K98};
  }{
    that there exist $c,\eta\in(0,\infty)$ and a strictly
    increasing $\varphi\in C^1((0,\eta),(0,\infty))$ which satisfy for all
    $z\in W$ with $f(z)<\eta$ that
    \begin{eqn}
      \llabel{eq18}
      \varphi'(f(z))\norm{\bigl(\grad H\bigr)\bigl(z\bigr)}
        =\norm{\bigl(\grad(\varphi\circ f)\bigr)\bigl(z\bigr)}\geq c
      \dott
    \end{eqn}
  }
  In the following let $\psi\colon[0,\eta)\to[0,\infty)$ satisfy
  $\psi(0)=0$ and for all $u\in(0,\eta)$ that
  \begin{eqn}
    \llabel{setting6}
    \psi(u)=\frac{\varphi(u)-\inf\varphi((0,\eta))}{c}
    \dott
  \end{eqn}
  \startnewargseq
  \argument{
    \lref{eq18};
    \lref{setting6};
    the fact that $\lim_{u\to 0}\varphi(u)
      =\inf\varphi((0,\eta))\in[0,\infty)$;
  }{
    that $\psi$ is continuous and strictly increasing on
    $[0,\eta)$, that $\psi|_{(0,\eta)}\in C^1((0,\eta),[0,\infty))$, and that for
    all $z\in W$ with $f(z)<\eta$ it holds that
    \begin{eqn}
      \llabel{eq19}
      \psi'(f(z))\norm{\bigl(\grad H\bigr)\bigl(z\bigr)}\geq1
      \qquad\text{and}\qquad
      \lim_{u\to 0}\psi(u)=0
      \dott
    \end{eqn}
  }
  \argument{
    \lref{eq19};
    \lref{eq13};
    \lref{eq16};
    the assumption that for all $t\in[\tau,T)$ it holds that $e(t)>\ell$;
  }{
    that there exists $a\in(\tau,T)$ which satisfies
    \begin{eqn}
      \llabel{eq20}
      \norm{\theta(a)-\mathbf{v}}<\frac14,
      \qquad 0<e(a)-\ell<\eta,
      \qquad\text{and}\qquad
      \psi(e(a)-\ell)<\frac14
      \dott
    \end{eqn}
  }
  \startnewargseq
  \argument{
    \lref{setting4};
    \lref{setting5}; \lref{eq20};
    the fact that $e$ is non-increasing;
    the fact that $\theta([\tau,T))\subseteq U$;
    the assumption that for all $t\in[\tau,T)$ it holds that
    $e(t)>\ell$;
    the fact that
    $W=\bigl\{z\in U\colon\bigl[(\norm{z-\mathbf{v}}<1)
      \wedge(H(z)>\ell)\bigr]\bigr\}$;
  }{
    that for all $t\in[a,T)$ with $\norm{\theta(t)-\mathbf{v}}<1$
    it holds that
    \begin{eqn}
      \llabel{eq20.1}
      \theta(t)\in W
      \qquad\text{and}\qquad
      0<f(\theta(t))=e(t)-\ell\leq e(a)-\ell<\eta
      \dott
    \end{eqn}
  }
  \argument{
    \lref{eq20.1};
    \lref{eq19};
  }{
    that for all $t\in[a,T)$ with $\norm{\theta(t)-\mathbf{v}}<1$
    it holds that
    \begin{eqn}
      \llabel{eq20.2}
      \psi'(f(\theta(t)))\norm{\bigl(\grad H\bigr)\bigl(\theta(t)\bigr)}\geq1
      \dott
    \end{eqn}
  }
  \argument{
    \lref{eq20.2};
    \lref{eq15};
    \lref{eq20.1};
    the fact that for all $t\in[a,T)$ it holds that
    $\norm{\theta'(t)}\geq0$;
    the chain rule;
  }{
    that for all $t\in[a,T)$ with $\norm{\theta(t)-\mathbf{v}}<1$ it
    holds that
    \begin{eqn}[a]
      \llabel{eq21}
      -\frac{\ddd}{\ddd t}\psi(e(t)-\ell)
        =-\psi'(e(t)-\ell)e'(t)
        =\psi'(f(\theta(t)))\norm{\bigl(\grad H\bigr)\bigl(\theta(t)\bigr)}
          \norm{\theta'(t)}
       \geq\norm{\theta'(t)}
      \dott
    \end{eqn}
  }
  \argument{
    \lref{eq21};
    \lref{eq20};
    the fact that $\psi$ is nonnegative;
    the fundamental theorem of calculus;
  }{
    that for all $b\in(a,T)$ with
    $\theta([a,b])\subseteq\{z\in\R^n\colon\norm{z-\mathbf{v}}<1\}$
    it holds that
    \begin{eqn}
      \llabel{eq22}
      \int_a^b\norm{\theta'(t)}\dd t
        \leq\psi(e(a)-\ell)-\psi(e(b)-\ell)
        \leq\psi(e(a)-\ell)<\frac14
      \dott
    \end{eqn}
  }
  \argument{
    \lref{eq22};
    \lref{eq20};
    the fact that $\theta$ is continuous;
    the fact that for all $b\in(a,T)$ with
    $b=\inf\{t\in[a,T)\colon
      \norm{\theta(t)-\mathbf{v}}\geq\nicefrac12\}$
    it holds that
    \begin{gather}
      \llabel{eq22.1}
      \norm{\theta(b)-\mathbf{v}}=\nicefrac12,
      \qquad
      \theta([a,b])\subseteq\{z\in\R^n\colon\norm{z-\mathbf{v}}<1\},\\
      \text{and}\qquad
      \frac14<\norm{\theta(b)-\theta(a)}
        \leq\int_a^b\norm{\theta'(t)}\dd t
    \end{gather};
  }{
    that
    \begin{eqn}
      \llabel{eq23}
      \theta([a,T))
        \subseteq\{z\in\R^n\colon\norm{z-\mathbf{v}}<\nicefrac12\}
      \dott
    \end{eqn}
  }
  \argument{
    \lref{eq23};
    \lref{eq22};
    the fact that $\theta|_{[\tau,a]}\in C^1([\tau,a],U)$;
    the monotone convergence theorem;
  }{
    that
    \begin{eqn}
      \llabel{eq24}
      \int_\tau^T\norm{\theta'(t)}\dd t
        =\biggl(\int_\tau^a\norm{\theta'(t)}\dd t\biggr)
          +\bigg(\int_a^T\norm{\theta'(t)}\dd t\biggr)
        \leq\biggl(\int_\tau^a\norm{\theta'(t)}\dd t\biggr)
          +\psi(e(a)-\ell)<\infty
      \dott
    \end{eqn}
  }
  \argument{
    \lref{eq24};
  }[verbs=ep]{
    \lref{eq12.2} in the case that for all
    $t\in[\tau,T)$ it holds that $e(t)>\ell$\dott
  }
This completes our proof of \lref{eq12.2}.

  In the following we prove \lref{item2}\dott
  \startnewargseq
  We thus assume that
  $\liminf_{t\to T}\norm{\theta(t)}<\infty$\dott
  \startnewargseq
  \argument{
    \lref{eq12.2};
    the assumption that
    $\liminf_{t\to T}\norm{\theta(t)}<\infty$;
    the fact that for all $b,d\in[\tau,T)$ with $b<d$ it holds
    that $\norm{\theta(d)-\theta(b)}
      \leq\int_b^d\norm{\theta'(t)}\dd t$;
  }{
    that there exists $\mathbf{v}\in\R^n$ such that
    \begin{eqn}
      \llabel{eq25}
      \lim_{t\to T}\norm{\theta(t)-\mathbf{v}}=0
      \dott
    \end{eqn}
  }
  \argument{
    \lref{eq25};
  }[verbs=ep]{
    \lref{item2}\dott
  }

  In the following we prove \lref{item4}\dott
  \startnewargseq
  We thus assume that
  $\liminf_{t\to T}\norm{\theta(t)-\vartheta}=0$,
  let $\mathbb{O}\subseteq\R^n$ be open with
  $\{\vartheta\}\cap\{p\}\subseteq\mathbb{O}\subseteq U\cup\{\vartheta\}$,
  and let $s\in J$\dott
  \startnewargseq
  \argument{
    \lref{item2};
    \lref{eq12.2};
    the assumption that
    $\liminf_{t\to T}\norm{\theta(t)-\vartheta}=0$;
    the fact that
    $\liminf_{t\to T}\norm{\theta(t)}\leq\norm{\vartheta}<\infty$;
  }{
    that
    \begin{eqn}
      \llabel{eq33}
      \int_\tau^T\norm{\theta'(t)}\dd t<\infty
      \qquad\text{and}\qquad
      \lim_{t\to T}\norm{\theta(t)-\vartheta}=0
      \dott
    \end{eqn}
  }
  In our proof of \lref{item4} we distinguish between the case
  $\vartheta\neq p$, the case ($\vartheta=p$ and $n=1$), and the
  case ($\vartheta=p$ and $n\geq2$)\dott
  \startnewargseq
  We first prove \lref{item4} in the case $\vartheta\neq p$\dott
  \startnewargseq
  \argument{
    \lref{eq2};
    \lref{eq3};
    \lref{eq4};\lref{eq33};
    the assumption that $\vartheta\neq p$;
    the triangle inequality;
  }{
    that there exists
    $\zeta\in(s,T)$ such that
    for all $t\in[\zeta,T)$ it holds that $\norm{\theta(t)-p}\geq\nicefrac12\norm{\vartheta-p}$
     and that
    \begin{eqn}[a]
      \llabel{eq34}
      \int_s^T\norm{\omega'(t)}\dd t
        &=\biggl(\int_s^\zeta\norm{\omega'(t)}\dd t\biggr)
          +\biggl(\int_\zeta^T\norm{\omega'(t)}\dd t\biggr)\leq\biggl(\int_s^\zeta\norm{\omega'(t)}\dd t\biggr)
          +\biggl(\int_\zeta^T
            \frac{\norm{\theta'(t)}}{\norm{\theta(t)-p}}\dd t\biggr)\\
        &\leq\biggl(\int_s^\zeta\norm{\omega'(t)}\dd t\biggr)
          +\frac{2}{\norm{\vartheta-p}}
            \biggl(\int_\zeta^T\norm{\theta'(t)}\dd t\biggr)\\
            &\leq\frac{\int_s^\zeta\norm{\theta'(t)}\dd t}
            {\min_{t\in[s,\zeta]}\norm{\theta(t)-p}}
          +\frac{2}{\norm{\vartheta-p}}
            \biggl(\int_\zeta^T\norm{\theta'(t)}\dd t\biggr)<\infty
      \dott
    \end{eqn}
  }
  \argument{
    \lref{eq34};
    the fact that $\omega|_{[s,T)}\in C^1([s,T),\R^n)$;
    \cref{lem:curve-length-speed}
    (applied with $I\curvearrowleft[s,T)$,
    $\gamma\curvearrowleft\omega|_{[s,T)}$ in the notation of
    \cref{lem:curve-length-speed});
  }{
    that
    \begin{eqn}
      \llabel{eq34.1}
      \curvelength{\omega|_{[s,T)}}
        =\int_{[s,T)}\norm{(\omega|_{[s,T)})'(t)}\dd t
        =\int_s^T\norm{\omega'(t)}\dd t<\infty
      \dott
    \end{eqn}
  }
  \argument{
    \lref{eq34.1};
    the fact that $s\in J$ was arbitrary;
  }[verbs=ep]{
    \lref{item4} in the case $\vartheta\neq p$\dott
  }
  We now prove \lref{item4} in the case $\vartheta=p$ and
  $n=1$\dott
  \startnewargseq
  \argument{
    \lref{eq2};
    the assumption that $n=1$;the fact that  $\omega|_{[s,T)}$ is continuous}{that $\omega|_{[s,T)}$ \llabel{arg1} is constant\dott}
    \argument{\lref{arg1};}{
     that
    \begin{eqn}
      \llabel{eq34.2}
      \curvelength{\omega|_{[s,T)}}=0
      \dott
    \end{eqn}
  }
  \argument{
    \lref{eq34.2};
    the fact that $s\in J$ was arbitrary;
  }[verbs=ep]{
    \lref{item4} in the case $\vartheta=p$ and $n=1$\dott
  }
  We now prove \lref{item4} in the case $\vartheta=p$ and
  $n\geq2$\dott
  \startnewargseq
  \argument{
    \lref{eq2};
    \lref{eq33};
    the assumption that $\mathbb{O}$ is open;
    the fact that $p\in\mathbb{O}\subseteq U\cup\{p\}$;
  }{
    that there exist $\mathfrak R\in(0,\infty)$ and
    $\kappa\in(s,T)$ which satisfy
    \begin{eqn}
      \llabel{eq35}
      \{z\in\R^n\colon0<\norm{z-p}<\mathfrak R\}\subseteq U
      \qquad\text{and}\qquad
      \theta([\kappa,T))
        \subseteq\{z\in\R^n\colon0<\norm{z-p}<\mathfrak R\}
      \dott
    \end{eqn}
  }
  In the following let
  $\mathcal{V}=\{z\in\R^n\colon0<\norm{z}<\mathfrak R\}$,
  and let $G\colon \mathcal{V}\to\R$ and
  $\xi\colon[\kappa,T)\to \mathcal{V}$ satisfy for all
  $z\in \mathcal{V}$, $t\in[\kappa,T)$ that
  \begin{eqn}
    \llabel{setting8}
    G(z)=F(p+z)
    \qquad\text{and}\qquad
    \xi(t)=\theta(t)-p
    \dott
  \end{eqn}
  \startnewargseq
  \argument{
    \lref{ass1};
    \lref{eq2};
    \lref{eq35};
    \lref{setting8};
    the assumption that $F\in C^1(U,\R)$ is definable;
    the assumption that $\theta\in C^1([\tau,T),U)$;
    the chain rule;
  }{
    that $G\in C^1(\mathcal{V},\R)$ is definable, that $\xi\in C^1([\kappa,T),\mathcal{V})$,
    and that for all $t\in[\kappa,T)$ it holds that
    \begin{eqn}
      \llabel{eq36}
      \xi'(t)=-(\grad F)(\theta(t))=-(\grad G)(\xi(t))
      \qquad\text{and}\qquad
      \omega(t)=\frac{\xi(t)}{\norm{\xi(t)}}
      \dott
    \end{eqn}
  }
  \argument{
    \lref{eq36};
    the fact that $G\in C^1(\mathcal{V},\R)$ is definable;
    the fact that $\xi\in C^1([\kappa,T),\mathcal{V})$;
    \cref{prop:annular-length}
    (applied with $F\curvearrowleft G$,
    $I\curvearrowleft[\kappa,T)$, $\theta\curvearrowleft \xi$,
    $\lambda\curvearrowleft([\kappa,T)\ni t\mapsto1\in[0,\infty))$ in the notation of
    \cref{prop:annular-length});
  }{
    that there exist $\delta\in(0,\mathfrak R)$ and
    $\mathbf{L}\in[0,\infty)$ which satisfy
    \begin{eqn}
      \llabel{eq37}
      \int_{\{t\in[\kappa,T)\colon0<\norm{\xi(t)}\leq\delta\}}
        \norm{\omega'(t)}\dd t\leq \mathbf{L}
      \dott
    \end{eqn}
  }
  \startnewargseq
  \argument{
    \lref{eq4};
    \lref{eq33};
    \lref{setting8};
     \lref{eq37};
    the assumption that $\vartheta=p$;
  }{
    that there exists $\mu\in(\kappa,T)$ such that 
    $\forall\,t\in[\mu,T)\colon\norm{\xi(t)}\leq\delta$ and that
    \begin{eqn}
      \llabel{eq38}
      \int_s^T\norm{\omega'(t)}\dd t
        =\biggl(\int_s^\mu\norm{\omega'(t)}\dd t\biggr)
          +\biggl(\int_\mu^T\norm{\omega'(t)}\dd t\biggr)
        \leq\biggl(\int_s^\mu\norm{\omega'(t)}\dd t\biggr)+\mathbf{L}<\infty
      \dott
    \end{eqn}
  }
  \argument{
    \lref{eq38};
    the fact that $\omega|_{[s,T)}\in C^1([s,T),\R^n)$;
    \cref{lem:curve-length-speed}
    (applied with $I\curvearrowleft[s,T)$,
    $\gamma\curvearrowleft\omega|_{[s,T)}$ in the notation of
    \cref{lem:curve-length-speed});
  }{
    that
    \begin{eqn}
      \llabel{eq38.1}
      \curvelength{\omega|_{[s,T)}}
        =\int_{[s,T)}\norm{(\omega|_{[s,T)})'(t)}\dd t
        =\int_s^T\norm{\omega'(t)}\dd t<\infty
      \dott
    \end{eqn}
  }
  \argument{
    \lref{eq38.1};
    the fact that $s\in J$ was arbitrary;
  }[verbs=ep]{
    \lref{item4} in the case $\vartheta=p$ and $n\geq2$\dott
  }
  This completes our proof of \lref{item4}.
\end{aproof}

\subsection{Finite spherical length for general open domains}
\label{subsec:open-domain-length}

\begin{athm}{proposition}{prop:spherical-length-away}
  [Finite spherical length for \GF\ trajectories that avoid the reference point]
  Let $n\in\N$, $\vartheta\in\R^n$, $\tau\in\R$,
  $T\in(\tau,\infty]$, let $U\subseteq\R^n$ be open,
  let $F\in C^1(U,\R)$ be definable,
  let $\theta\in C^1([\tau,T),U\backslash\{\vartheta\})$
  satisfy for all $t\in[\tau,T)$ that
  \begin{eqn}
    \llabel{ass1}
    \theta'(t)=-(\grad F)(\theta(t)),
  \end{eqn}
  and let $P\colon U\backslash\{\vartheta\}\to\R^n$ satisfy for all
  $v\in U\backslash\{\vartheta\}$ that
  $P(v)=|v-\vartheta|^{-1}(v-\vartheta)$\dott
  \cfload[.]
  Then $P\circ\theta$ has finite length\cfadd{definition:curve-length}\dott
  \cfout[.]
\end{athm}
\begin{aproof}
  Throughout this proof let $\mathbb S=\{x\in\R^n\colon\norm{x}=1\}$,
  let $\phi\colon\R\to(-1,1)$ satisfy for all $e\in\R$
  that $\phi(e)=e/\sqrt{1+e^2}$, and let
  $\omega\colon[\tau,T)\to\R^n$ satisfy
  \begin{eqn}
    \llabel{setting0}
    \omega=P\circ\theta
    \dott
  \end{eqn}
  \cfload[.]
  \startnewargseq
  \argument{
    the assumption that for all $v\in U\backslash\{\vartheta\}$
    it holds that $P(v)=|v-\vartheta|^{-1}(v-\vartheta)$;
    \lref{setting0};
    the fact that $\theta([\tau,T))\subseteq U\backslash\{\vartheta\}$;
  }{
    that for all $t\in[\tau,T)$ it holds that
    \begin{eqn}
      \llabel{eq0}
      \omega(t)=\frac{\theta(t)-\vartheta}{\norm{\theta(t)-\vartheta}}
      \dott
    \end{eqn}
  }
  \argument{
    the assumption that $F$ is definable;
  }{
    that $U$ is definable\dott
  }
  In our proof of \cref{prop:spherical-length-away} we distinguish between
  the case that there exists $v\in\R^n\backslash\{\vartheta\}$ such that
  $\limsup_{t\to T}\norm{\theta(t)-v}=0$ and the case that there exists no
  $v\in\R^n\backslash\{\vartheta\}$ such that
  $\limsup_{t\to T}\norm{\theta(t)-v}=0$\dott
  \startnewargseq
  We first prove that $P\circ\theta$ has finite length in the case that
  there exists $v\in\R^n\backslash\{\vartheta\}$ such that
  $\limsup_{t\to T}\norm{\theta(t)-v}=0$\dott
  \startnewargseq
  In the following let $v\in\R^n\backslash\{\vartheta\}$ satisfy
  $\limsup_{t\to T}\norm{\theta(t)-v}=0$\dott
  \startnewargseq
  \argument{
    \itref{maintheorem1}{item4}
    (applied with $\vartheta\curvearrowleft v$, $p\curvearrowleft\vartheta$,
    $O\curvearrowleft\emptyset$ in the notation of \cref{maintheorem1});
    \lref{ass1};
    \lref{setting0};
    \lref{eq0};
    the fact that $\limsup_{t\to T}\norm{\theta(t)-v}=0$;
    the fact that $\{v\}\cap\{\vartheta\}=\emptyset$;
    the fact that $\theta^{-1}(\{\vartheta\})=\emptyset$;
  }[verbs=ep]{
    that $P\circ\theta$ has finite length in the case that
    there exists $v\in\R^n\backslash\{\vartheta\}$ such that
    $\limsup_{t\to T}\norm{\theta(t)-v}=0$\dott
  }
  We now prove that $P\circ\theta$ has finite length in the case that
  there exists no $v\in\R^n\backslash\{\vartheta\}$ such that
  $\limsup_{t\to T}\norm{\theta(t)-v}=0$\dott
  \startnewargseq
  In the following we thus assume that there exists no
  $v\in\R^n\backslash\{\vartheta\}$ such that
  $\limsup_{t\to T}\norm{\theta(t)-v}=0$\dott
  \startnewargseq
  \argument{
    \itref{maintheorem1}{item2}
    (applied with $p\curvearrowleft\vartheta$ in the notation of
    \cref{maintheorem1});
    \lref{ass1};
    \lref{eq0};
    the fact that $U$ is definable and open;
    the fact that $\liminf_{t\to T}\norm{\theta(t)}\in[0,\infty]$;
    the assumption that there exists no $v\in\R^n\backslash\{\vartheta\}$
    such that $\limsup_{t\to T}\norm{\theta(t)-v}=0$;
  }{
    that $\lim_{t\to T}\theta(t)=\vartheta$ or that
    $\lim_{t\to T}\norm{\theta(t)}=\infty$\dott
  }
  In the following let $\nu\in\{-1,1\}$ satisfy
  \begin{eqn}
    \llabel{setting1}
    \nu=
    \begin{cases}
      1 & \colon \lim_{t\to T}\theta(t)=\vartheta,\\
      -1 & \colon \lim_{t\to T}\norm{\theta(t)}=\infty
    \end{cases},
  \end{eqn}
  let $\Gamma,r\colon[\tau,T)\to(0,\infty)$ satisfy for all
  $t\in[\tau,T)$ that
  \begin{eqn}
    \llabel{setting2}
    \Gamma(t)=\norm{\theta(t)-\vartheta}
    \qquad\text{and}\qquad r(t)=(\Gamma(t))^\nu,
  \end{eqn}
  let
  \begin{eqn}
    \llabel{setting3}
    \Omega=\{(s,x)\in(0,1)\times\mathbb S\colon
      \vartheta+s^\nu x\in U\},
  \end{eqn}
  and let $z\colon\Omega\to U$,
  $h=(h(s,x))_{(s,x)\in\Omega}\colon\Omega\to(-1,1)$, and
  $\eta\colon\Omega\to(0,\infty)$ satisfy for all $(s,x)\in\Omega$ that
  \begin{eqn}
    \llabel{setting4}
    z(s,x)=\vartheta+s^\nu x,\qquad
    h(s,x)=\phi(F(z(s,x))),
    \qquad\text{and}\qquad
    \eta(s,x)=\phi'(F(z(s,x)))
    \dott
  \end{eqn}
  \cfload[.]
  \startnewargseq
  \argument{
    \lref{ass1};
    \lref{eq0};
    \lref{setting1};
    \lref{setting2};
    the assumption that $\theta\in C^1([\tau,T),U\backslash\{\vartheta\})$;
    the chain rule;
  }{
    that $\Gamma,r\in C^1([\tau,T),(0,\infty))$,
    that $\omega\in C^1([\tau,T),\R^n)$,
    that $\omega([\tau,T))\subseteq\mathbb S$, that
    $\lim_{t\to T}r(t)=0$, and that there exists $c\in[\tau,T)$
    which satisfies for all $t\in[c,T)$ that
    \begin{eqn}
      \llabel{eq1}
      0<r(t)<1,\qquad (r(t),\omega(t))\in\Omega,
      \qquad\text{and}\qquad z(r(t),\omega(t))=\theta(t)
      \dott
    \end{eqn}
  }
  In the following let
  $\alpha\colon[c,T)\to(0,\infty)$ satisfy for all $t\in[c,T)$ that
  \begin{eqn}
    \llabel{setting5}
    \alpha(t)=\frac{(r(t))^{-2\nu}}{\eta(r(t),\omega(t))}
    \dott
  \end{eqn}
  \startnewargseq
  \argument{
    \lref{setting3};
    \lref{setting4};
    the fact that $U$ is definable and open;
    the assumption that $F\in C^1(U,\R)$;
    the assumption that $F$ is definable;
    the fact that for all $e\in\R$ it holds that
    $\phi'(e)=(1+e^2)^{-3/2}>0$;
  }{
    that $\Omega$ is definable and open in $(0,1)\times\mathbb S$,
    that $h\in \COne{\Omega}{(-1,1)}$ is definable, and that $\alpha$ is continuous\dott
  }
  \argument{
    \lref{setting4};
    the chain rule;
    the fact that for all $x\in\mathbb S$ it holds that
    $T_x\mathbb S=\{w\in\R^n\colon\langle w,x\rangle=0\}$;
  }{
    that for all $(s,x)\in\Omega$ it holds that\cfadd{definition:tangent-space}\cfadd{definition:angular-gradient}
    \begin{eqn}[a]
      \llabel{eq2}
      (\partial_s h)(s,x)
        &=\nu s^{\nu-1}\eta(s,x)\langle(\grad F)(z(s,x)),x\rangle\\
      \qquad\text{and}\qquad
      (\nabla_{\mathbb S}h)(s,x)
        &=s^\nu\eta(s,x)
          ((\grad F)(z(s,x))-\langle(\grad F)(z(s,x)),x\rangle x)
      \dott
    \end{eqn}
  }
  \argument{
    \lref{ass1};
    \lref{eq0};
    \lref{setting2};
    the chain rule;
  }{
    that for all $t\in[\tau,T)$ it holds that
    \begin{eqn}
      \llabel{eq3}
      \Gamma'(t)=\frac{\langle\theta(t)-\vartheta,\theta'(t)\rangle}{\Gamma(t)}
        =\langle\omega(t),\theta'(t)\rangle
        =-\langle\bigl(\grad F\bigr)\bigl(\theta(t)\bigr),\omega(t)\rangle
      \dott
    \end{eqn}
  }
  \argument{
    \lref{eq3};
    \lref{setting2};
    \lref{setting5};
    \lref{eq1};
    \lref{eq2};
    the fact that $\nu^2=1$;
  }{
    that for all $t\in[c,T)$ it holds that
    \begin{eqn}[a]
      \llabel{eq4}
      r'(t)&=\nu(\Gamma(t))^{\nu-1}\Gamma'(t)
        =-\nu(r(t))^{1-\nu}\langle(\grad F)(\theta(t)),\omega(t)\rangle=-\alpha(t)(r(t))^2
          (\partial_s h)(r(t),\omega(t))
      \dott
    \end{eqn}
  }
  \argument{
    \lref{ass1};
    \lref{eq0};
    \lref{setting2};
    \lref{setting5};
    \lref{eq1};
    \lref{eq2};
    \lref{eq3}
  }{
    that for all $t\in[c,T)$ it holds that
    \begin{eqn}[a]
      \llabel{eq5}
      \omega'(t)&=\frac{\theta'(t)}{\Gamma(t)}
        -\frac{\theta(t)-\vartheta}{(\Gamma(t))^2}\Gamma'(t)
        =\frac{\theta'(t)-\Gamma'(t)\omega(t)}{\Gamma(t)}\\
        &=-(r(t))^{-\nu}
          ((\grad F)(\theta(t))
            -\langle(\grad F)(\theta(t)),\omega(t)\rangle\omega(t))
        =-\alpha(t)(\nabla_{\mathbb S}h)(r(t),\omega(t))
      \dott
    \end{eqn}
  }
  \argument{
    \lref{eq5};
      \lref{eq1};
    \lref{eq4};
    the fact that $r\in C^1([\tau,T),(0,\infty))$;
    the fact that $\omega\in C^1([\tau,T),\R^n)$;
    the fact that $\omega([\tau,T))\subseteq\mathbb S$;
    the fact that $\mathbb S$ is a nonempty, compact, and definable \submanifold{2}{n};
    the fact that $\alpha$ is continuous;
    the fact that $\lim_{t\to T}r(t)=0$;
    \itref{cor:angular-length}{item2}
    (applied with $D\curvearrowleft n$, $\mathfrak R\curvearrowleft1$,
    $M\curvearrowleft\mathbb S$, $\mathcal O\curvearrowleft\Omega$,
    $\mathsf a\curvearrowleft1$, $\mathsf b\curvearrowleft1$,
    $I\curvearrowleft[c,T)$, $r\curvearrowleft r|_{[c,T)}$,
    $x\curvearrowleft\omega|_{[c,T)}$,
    $A\curvearrowleft([c,T)\ni t\mapsto\mathrm{I}_n\in\R^{n\times n})$,
    $\tau\curvearrowleft c$ in the notation of
    \cref{cor:angular-length});
  }{
    that
    \begin{eqn}
      \llabel{eq6}
      \int_c^T\norm{\omega'(t)}\dd t<\infty
      \dott
    \end{eqn}
  }
  \argument{
    \lref{eq6};
    \lref{setting0};
    the fact that $\omega\in C^1([\tau,T),\R^n)$;
    the fact that $[\tau,c]$ is compact;
    \cref{lem:curve-length-speed};
  }[verbs=ep]{
    that $P\circ\theta$ has finite length in the case that
    there exists no $v\in\R^n\backslash\{\vartheta\}$ such that
    $\limsup_{t\to T}\norm{\theta(t)-v}=0$\dott
  }
  \finishproofthis
  \startnewargseq
\end{aproof}

In the proof of the finite spherical length result in \cref{maintheorem2} below we
employ the well-known fact, which we recall in \cref{lem:fibre_is_convex}
below, that the fibres/level sets of \GF\ trajectories are convex.
Only for completeness, we include a proof for \cref{lem:fibre_is_convex} here.

\begin{athm}{lemma}{lem:fibre_is_convex}[Convexity of the level sets of \GF\ trajectories]
Let $ n \in \mathbb{N} $,
$ \vartheta \in \R^n $,
let $ U \subseteq \R^n $ be open,
let $ I \subseteq \R $ be convex with nonempty interior,
and
let $ F \in C^1( U, \R ) $, $ \theta \in C^1( I, U ) $
satisfy for all $ t \in I $ that
\begin{equation}
  \theta'(t) = - (\nabla F)(\theta(t))
  .
\end{equation}
Then it holds that
$
  \{ t \in I \colon \theta( t ) = \vartheta \}
$
is convex.
\end{athm}
\begin{aproof}
  Throughout this proof let $A\subseteq I$ satisfy
  \begin{eqn}
    \llabel{setting1}
    A=\{t\in I\colon\theta(t)=\vartheta\}
    \dott
  \end{eqn}
  \cfload[.]
  \startnewargseq
  We claim that for all $s,t\in A$ with $s<t$ it holds that
  \begin{eqn}
    \llabel{claim1}
    [s,t]\subseteq A
    \dott
  \end{eqn}
  \startnewargseq
  To prove \lref{claim1}, let $s,t\in A$ satisfy $s<t$\dott
  \startnewargseq
  \argument{
    \lref{setting1};
    the fact that $s,t\in A$;the fact that $s<t$;
    the assumption that $I$ is convex;
  }{
    that
    \begin{eqn}
      \llabel{eq1}
      [s,t]\subseteq I
      \qquad\text{and}\qquad
      \theta(s)=\theta(t)=\vartheta
      \dott
    \end{eqn}
  }
  \argument{
    \lref{eq1};
    the assumption that $F\in C^1(U,\R)$;
    the assumption that $\theta\in C^1(I,U)$;
    the assumption that for all $u\in I$ it holds that
    $\theta'(u)=-\bigl(\grad F\bigr)\bigl(\theta(u)\bigr)$;
    the chain rule;
  }{
    that $(F\circ\theta)|_{[s,t]}\in C^1([s,t],\R)$ and that
    for all $u\in(s,t)$ it holds that
    \begin{eqn}
      \llabel{eq2}
      \frac{\ddd}{\ddd u}F(\theta(u))
        =\langle\bigl(\grad F\bigr)\bigl(\theta(u)\bigr),\theta'(u)\rangle
        =-\norm{\theta'(u)}^2
      \dott
    \end{eqn}
  }
  \argument{
    \lref{eq2};
    \lref{eq1};
    the fact that $(F\circ\theta)|_{[s,t]}\in C^1([s,t],\R)$;
    the fundamental theorem of calculus;
  }{
    that
    \begin{eqn}
      \llabel{eq3}
      \int_s^t\norm{\theta'(u)}^2\dd u
        =F(\theta(s))-F(\theta(t))=0
      \dott
    \end{eqn}
  }
  \argument{
    \lref{eq3};
    \lref{eq1};
    the assumption that $\theta\in C^1(I,U)$;
    the fact that $s<t$;
  }{
    that for all $u\in[s,t]$ it holds that
    \begin{eqn}
      \llabel{eq4}
     \theta'(u)=0
      \dott
    \end{eqn}
  }
  \argument{
    \lref{eq4};
    \lref{eq1};
    the assumption that $\theta\in C^1(I,U)$;
    the fundamental theorem of calculus;
  }{
    that for all $u\in[s,t]$ it holds that
    \begin{eqn}
      \llabel{eq5}
      \theta(u)=\theta(s)+\int_s^u\theta'(v)\dd v=\vartheta
      \dott
    \end{eqn}
  }
  \argument{
    \lref{eq5};
    \lref{setting1};
    \lref{eq1};
    the fact that $s,t\in A$ with $s<t$ were arbitrary;
  }[verbs=ep]{
    \lref{claim1}\dott
  }

  \startnewargseq
  \argument{ \lref{setting1};
    \lref{claim1};
  }[verbs=ep]{
    that $\{t\in I\colon\theta(t)=\vartheta\}$ is convex\dott
  }
\end{aproof}

\begin{athm}{corollary}{maintheorem2}[Finite spherical length for \GF\ trajectories on a general open set]
Let $ n \in \mathbb{N} $,
$ \vartheta \in \R^n $,
$ \tau \in \R $, $ T \in ( \tau, \infty ] $,
let $ U \subseteq \R^n $ be open,
let $ F \in C^1( U, \R ) $ be definable,
let
$ \theta \in C^1( [\tau,T), U ) $ satisfy for all $ t \in [ \tau, T ) $ that
\begin{equation}
  \theta'(t) = - (\nabla F)( \theta( t ) )
  ,
\end{equation}
and let
$ P \colon U \to \R^n $
satisfy for all $ v \in U \backslash \{ \vartheta \} $
that
$
  P( v ) = | v - \vartheta |^{ - 1 } ( v - \vartheta )
$.
Then $ P \circ \theta $ has finite length.
\end{athm}
\begin{aproof}
  \argument{
    \cref{prop:spherical-length-away};
  }{
    that if $\vartheta\notin U$, then $P\circ\theta$ has finite length\dott
  }
  In the following assume that $\vartheta\in U$.
  Throughout this proof let
  $A=\{t\in[\tau,T)\colon\theta(t)=\vartheta\}$\dott\cfload[.]
  \startnewargseq
  \argument{
    \cref{lem:fibre_is_convex}
    (applied with $I\curvearrowleft[\tau,T)$ in the notation of
    \cref{lem:fibre_is_convex});
    the fact that $[\tau,T)$ is convex with nonempty interior;
    the fact that $A=\{t\in[\tau,T)\colon\theta(t)=\vartheta\}$;
  }{
    that $A$ is convex\dott
  }
  In the following let $\mathcal J$ satisfy
  \begin{eqn}
    \llabel{setting7}
    \mathcal J=\bigl\{J\subseteq[\tau,T)\backslash A\colon\bigl[
      (J\text{ is convex})\wedge
      \bigl(\forall\,V\text{ convex},\ J\subseteq V\subseteq[\tau,T)\backslash A
        \colon V=J\neq\emptyset\bigr)\bigr]\bigr\}
    \dott
  \end{eqn}
  \startnewargseq
  \argument{
    \lref{setting7};
    the fact that $A$ is convex;
    the fact that $A$ is closed in $[\tau,T)$;
  }{
    that $\mathcal J$ contains at most two intervals\dott
  }
  We claim that for all $J\in\mathcal J$ it holds that
  \begin{eqn}
    \llabel{eq10}
    \curvelength{(P\circ\theta)|_J}<\infty
    \dott
  \end{eqn}
  To prove \lref{eq10}, let $J\in\mathcal J$, $c\in J$, and let
  $a=\inf J$ and $b=\sup J$.
  \startnewargseq
  \argument{
    \cref{prop:spherical-length-away}
    (applied with $\tau\curvearrowleft c$, $T\curvearrowleft b$,
    $P\curvearrowleft P|_{U\backslash\{\vartheta\}}$,
    $\theta\curvearrowleft\theta|_{[c,b)}$
    in the notation of \cref{prop:spherical-length-away});
    the fact that $[c,b)\subseteq J$;
  }{
    that
    \begin{eqn}
      \llabel{eq11}
      \curvelength{(P\circ\theta)|_{[c,b)}}<\infty
      \dott
    \end{eqn}
  }
  \argument{  \cref{lem:curve-length-speed};
    the fact that $(P\circ\theta)|_J\in C^1(J,\R^n)$;
  }{
    that if $a\in J$, then it holds that
    \begin{eqn}
      \llabel{eq12}
      \curvelength{(P\circ\theta)|_{[a,c]}}<\infty
      \dott
    \end{eqn}
  }
  \argument{
    the assumption that $\theta'=-(\grad F)\circ\theta$;
    the chain rule;
  }{
    that if $a\notin J$, then for all $t\in[-c,-a)$ it holds that
    \begin{eqn}
      \llabel{eq13}
      \frac{\ddd}{\ddd t}\theta(-t)
        =-\theta'(-t)=\bigl(\grad F\bigr)\bigl(\theta(-t)\bigr)
        =-\bigl(\grad(-F)\bigr)\bigl(\theta(-t)\bigr)
      \dott
    \end{eqn}
  }
  \argument{
    \lref{eq13};
    \cref{prop:spherical-length-away}
    (applied with $\tau\curvearrowleft-c$, $T\curvearrowleft-a$,
    $F\curvearrowleft-F$,
    $P\curvearrowleft P|_{U\backslash\{\vartheta\}}$,
    $\theta\curvearrowleft([-c,-a)\ni t\mapsto\theta(-t)\in U)$
    in the notation of \cref{prop:spherical-length-away});
  }{
    that if $a\notin J$, then it holds that
    \begin{eqn}
      \llabel{eq14}
      \curvelength{(P\circ\theta)|_{(a,c]}}<\infty
      \dott
    \end{eqn}
  }
  \argument{
    \lref{eq14};
    \lref{eq11};
    \lref{eq12};
    the triangle inequality;
    the fact that $J\in\mathcal J$ was arbitrary;
  }[verbs=ep]{
    \lref{eq10}\dott
  }

  \startnewargseq
  \argument{
    \lref{eq10};
    the fact that $\mathcal J$ contains at most two intervals;
    the fact that for all $t\in A$ it holds that
    $P(\theta(t))=P(\vartheta)$;
    the fact that for all $t\in[\tau,T)\backslash A$ it holds that
    $\norm{P(\theta(t))}=1$;
    the triangle inequality;
  }{
    that for all $k\in\N$, $t_0,t_1,\dots,t_k\in[\tau,T)$ with
    $t_0<t_1<\ldots<t_k$ it holds that
    \begin{eqn}
      \llabel{eq15}
      \sum_{j=1}^k\norm{P(\theta(t_j))-P(\theta(t_{j-1}))}
        \leq\sum_{J\in\mathcal J}\curvelength{(P\circ\theta)|_J}
          +2\bigl(1+\norm{P(\vartheta)}\bigr)<\infty
      \dott
    \end{eqn}
  }
  \argument{
    \lref{eq15};
  }[verbs=ep]{
    that $\curvelength{P\circ\theta}<\infty$\dott
  }
\end{aproof}
\subsection{Finite spherical length for general time intervals}
\label{subsec:general-time-intervals}

In the proof of the finite spherical length result in \cref{maintheorem3} below we
employ the elementary length bound for decomposed curves without the
boundary points in the following \cref{lem:curve-length-splitting}.
Only for completeness, we include a proof for
\cref{lem:curve-length-splitting} here.

\begin{athm}{lemma}{lem:curve-length-splitting}
  Let $n\in\N$, $B\in\R$, let $I\subseteq\R$ be convex,
  let $c\in\operatorname{int}(I)$, let $a=\inf(I)$, let $b=\sup(I)$,
  and let $\gamma\colon I\to\R^n$ satisfy
  for all $t\in I$ that
  \begin{eqn}
    \llabel{ass1}
    \norm{\gamma(t)}\leq B
    \dott
  \end{eqn}
  \cfload[.]
  Then
  \begin{eqn}
    \llabel{concl1}
    \curvelength{\gamma}
      \leq\curvelength{\gamma|_{(a,c]}}
        +\curvelength{\gamma|_{[c,b)}}+4B
    \dott
  \end{eqn}
  \cfout[.]
\end{athm}
\begin{aproof}
  \argument{
    the assumption that $I$ is convex;
    the fact that $a=\inf I$ and $b=\sup I$;
    the fact that $c\in\operatorname{int}(I)$;
  }{
    that
    \begin{eqn}
      \llabel{eq1}
      a<c<b,
      \qquad (a,b)\subseteq I,
      \qquad\text{and}\qquad
      I\backslash(a,b)\subseteq\{a,b\}\cap\R
      \dott
    \end{eqn}
  }
  \argument{
    \lref{ass1};
    the triangle inequality;
  }{
    that for all $s,t\in I$ it holds that
    \begin{eqn}
      \llabel{eq2}
      \norm{\gamma(t)-\gamma(s)}
        \leq\norm{\gamma(t)}+\norm{\gamma(s)}\leq2B
      \dott
    \end{eqn}
  }
  \argument{
    \lref{eq1};
    the triangle inequality;
  }{
    that for all $k\in\N$, $t=(t_0,\dots,t_k)\in(a,b)^{k+1}$
    with $t_0<t_1<\ldots<t_k$ it holds that
    \begin{eqn}
      \llabel{eq3}
      \sum_{j=1}^k\norm{\gamma(t_j)-\gamma(t_{j-1})}
        \leq\curvelength{\gamma|_{(a,c]}}
          +\curvelength{\gamma|_{[c,b)}}
      \dott
    \end{eqn}
  }
  \argument{
    \lref{eq3};
    \lref{eq1};
    \lref{eq2};
    the fact that for all $k\in\N$, $t=(t_0,t_1,\dots,t_k)\in I^{k+1}$
    with $t_0<t_1<\ldots<t_k$ it holds that
    $\{t_0,t_1,\dots,t_k\}\backslash(a,b)\subseteq\{t_0,t_k\}$;
  }{
    that for all $k\in\N$, $t=(t_0,t_1,\dots,t_k)\in I^{k+1}$
    with $t_0<t_1<\ldots<t_k$ it holds that
    \begin{eqn}[a]
      \llabel{eq4}
      \sum_{j=1}^k\norm{\gamma(t_j)-\gamma(t_{j-1})}
        &\leq\biggl(\sum_{j\in\{1,\dots,k\},\,
            (t_{j-1},t_j)\in(a,b)^2}
          \norm{\gamma(t_j)-\gamma(t_{j-1})}\biggr)+4B
      \\
        &\leq\curvelength{\gamma|_{(a,c]}}
          +\curvelength{\gamma|_{[c,b)}}+4B
      \dott
    \end{eqn}
  }
  \argument{
    \lref{eq4};
  }[verbs=ep]{
    that
    \begin{eqn}
      \curvelength{\gamma}
        \leq\curvelength{\gamma|_{(a,c]}}
          +\curvelength{\gamma|_{[c,b)}}+4B
      \dott
    \end{eqn}
  }
\end{aproof}

\begin{athm}{theorem}{maintheorem3}
Let $ n \in \mathbb{N} $,
$ \vartheta \in \R^n $,
let $ U \subseteq \R^n $ be open,
let $ F \in C^1( U, \R ) $ be definable,
let $ I \subseteq \R $ be convex with nonempty interior,
let
$ \theta \in C^1( I, U ) $ satisfy
\begin{eqn}
  \theta' = - (\nabla F)\circ\theta,
\end{eqn}
and let
$ P \colon U \to \R^n $
satisfy for all $ v \in U \backslash \{ \vartheta \} $
that
$
  P( v ) = | v - \vartheta |^{ - 1 } ( v - \vartheta )
$.
Then $ P \circ \theta $ has finite length.
\end{athm}
\begin{aproof}
  Throughout this proof let $a=\inf I$, $b=\sup I$,
  $c\in\operatorname{int}(I)$, let
  $\mathbf{P}\colon\R^n\to\R^n$ satisfy that $\mathbf{P}|_U=P$ and that for all
  $v\in\R^n\backslash\{\vartheta\}$ it holds that
  \begin{eqn}
    \llabel{setting1}
    \mathbf{P}(v)=\frac{v-\vartheta}{\norm{v-\vartheta}},
  \end{eqn}
  let $B\in(0,\infty)$ satisfy $B=1+\norm{\mathbf{P}(\vartheta)}$,
  and let $\gamma\colon[-c,-a)\to U$ satisfy for all
  $t\in[-c,-a)$ that
  \begin{eqn}
    \llabel{setting2}
    \gamma(t)=\theta(-t)
    \dott
  \end{eqn}
  \cfload[.]
  \startnewargseq
  \argument{
    the assumption that $I$ is convex with nonempty interior;
    the fact that $a=\inf I$ and $b=\sup I$;
    the fact that $c\in\operatorname{int}(I)$;
  }{
    that
    \begin{eqn}[a]
      \llabel{eq1}
      -\infty\leq a<c<b\leq\infty,\qquad (a,b)\subseteq I,
    \qquad\text{and}\qquad
      I\backslash(a,b)\subseteq\{a,b\}\cap\R
      \dott
    \end{eqn}
  }
  \argument{
    \lref{setting1};
    the fact that $B=1+\norm{\mathbf{P}(\vartheta)}$;
  }{
    that for all $v\in\R^n$ it holds that
    \begin{eqn}
      \llabel{eq2}
      \norm{\mathbf{P}(v)}\leq B
      \dott
    \end{eqn}
  }
  \argument{
    \lref{setting1};
    \lref{eq1};
    the assumption that $U$ is open;
    the assumption that $F\in C^1(U,\R)$; the fact that $F$ is definable;
    the assumption that $\theta\in C^1(I,U)$;
    the assumption that for all $t\in I$ it holds that
    $\theta'(t)=-\bigl(\grad F\bigr)\bigl(\theta(t)\bigr)$;
    the fact that $\mathbf{P}|_U=P$;
    \cref{maintheorem2}
    (applied with $\tau\curvearrowleft c$, $T\curvearrowleft b$,
    $\theta\curvearrowleft\theta|_{[c,b)}$,
    $P\curvearrowleft\mathbf{P}|_U$ in the notation of \cref{maintheorem2});
  }{
    that
    \begin{eqn}
      \llabel{eq3}
      \curvelength{(P\circ\theta)|_{[c,b)}}<\infty
      \dott
    \end{eqn}
  }
  \argument{
    \lref{setting2};
    \lref{eq1};
    the assumption that $F\in C^1(U,\R)$ ; the fact that $F$ is definable;
    the assumption that $\theta\in C^1(I,U)$;
    the assumption that for all $t\in I$ it holds that
    $\theta'(t)=-\bigl(\grad F\bigr)\bigl(\theta(t)\bigr)$;
    the chain rule;
  }{
    that $-F\in C^1(U,\R)$, that $-F$ is definable,
    that $\gamma\in C^1([-c,-a),U)$,
    and that for all $t\in[-c,-a)$ it holds that
    \begin{eqn}
      \llabel{eq4}
      \gamma'(t)=-\theta'(-t)
        =\bigl(\grad F\bigr)\bigl(\gamma(t)\bigr)=-\bigl(\grad(-F)\bigr)\bigl(\gamma(t)\bigr)
      \dott
    \end{eqn}
  }
  \argument{
    \lref{eq4};
    \lref{setting1};
    \lref{eq1};
    the assumption that $U$ is open;
    the fact that $-F\in C^1(U,\R)$; the fact that $-F$ is definable;
    the fact that $\gamma\in C^1([-c,-a),U)$;
    \cref{maintheorem2}
    (applied with $\tau\curvearrowleft-c$, $T\curvearrowleft-a$,
    $F\curvearrowleft-F$, $\theta\curvearrowleft\gamma$,
    $P\curvearrowleft\mathbf{P}$ in the notation of \cref{maintheorem2});
  }{
    that
    \begin{eqn}
      \llabel{eq5}
      \curvelength{\mathbf{P}\circ\gamma}<\infty
      \dott
    \end{eqn}
  }
  \argument{
    \lref{eq5};
    \lref{setting2};
    the fact that $\mathbf{P}|_U=P$;
  }{
    that
    \begin{eqn}
      \llabel{eq6}
      \curvelength{(P\circ\theta)|_{(a,c]}}
        =\curvelength{\mathbf{P}\circ\gamma}<\infty
      \dott
    \end{eqn}
  }
  \argument{
    \lref{eq6};
    \lref{eq2};
    \lref{eq3};
    the fact that $\mathbf{P}|_U=P$;
    the fact that $\theta(I)\subseteq U$;
    \cref{lem:curve-length-splitting}
    (applied with $\gamma\curvearrowleft P\circ\theta$
    in the notation of \cref{lem:curve-length-splitting});
  }[verbs=ep]{
    that $P\circ\theta$ has finite length\dott
  }
\end{aproof}

\subsection{Finite spherical length for locally definable potentials}
\label{sec:locally-definable-corollary}

In \cref{cor:spherical} below we establish finite spherical length of all \GF\
trajectories with limit points in the domain of a locally definable potential.
In our proof of \cref{cor:spherical} we employ the well-known consequence of the \KL\ inequality, which
we recall in the following \cref{lem:locally-definable-endpoint}, that such
\GF\ trajectories are also convergent.
The proof of \cref{lem:locally-definable-endpoint} is based on an application
of Kurdyka~\cite[Theorem~1]{K98}.
For completeness, we include a proof for
\cref{lem:locally-definable-endpoint} here.

\begin{athm}{proposition}{lem:locally-definable-endpoint}[Convergence of \GF\ trajectories]
Let $n\in\N$, let $U\subseteq\R^n$ be open,
let $F\in C^1(U,\R)$ be locally definable,
let $I\subseteq\R$ be convex with nonempty interior,
let $\theta\in C^1(I,U)$ satisfy for all $t\in I$ that
$\theta'(t)=-\bigl(\nabla F\bigr)\bigl(\theta(t)\bigr)$,
and let $b\in U$ satisfy
$\liminf_{t\to\sup(I)}\norm{\theta(t)-b}=0$\dott\cfload[.]
Then $ \limsup_{ t \to \sup(I) } | \theta(t) - b | = 0 $\dott\cfout[.]
\end{athm}
\begin{aproof}
  Throughout this proof let $T=\sup(I)$,
  let $h\colon I\to\R$ satisfy $h=F\circ\theta$, and let
  $(t_j)_{j\in\N}\subseteq I$ satisfy
  \begin{eqn}
    \llabel{setting1}
    \lim_{j\to\infty}t_j=T
    \qquad\text{and}\qquad
    \lim_{j\to\infty}\theta(t_j)=b
    \dott
  \end{eqn}
  \cfload[.]
  \startnewargseq
  We claim that
  \begin{eqn}
    \llabel{claim0}
    \limsup_{t\to T}\norm{\theta(t)-b}=0
    \dott
  \end{eqn}
  \startnewargseq
  In our proof of \lref{claim0} we distinguish between the case
  $T\in I$ and the case $T\notin I$\dott
  \startnewargseq
  We first prove \lref{claim0} in the case $T\in I$\dott
  \startnewargseq
  \argument{
    the fact that $\theta$ is continuous;
    the assumption that $\liminf_{t\to T}\norm{\theta(t)-b}=0$;
  }{
    that
    \begin{eqn}
      \llabel{eq0}
      \limsup_{t\to T}\norm{\theta(t)-b}
        =\norm{\theta(T)-b}
        =\liminf_{t\to T}\norm{\theta(t)-b}=0
      \dott
    \end{eqn}
  }
  \argument{
    \lref{eq0};
  }[verbs=ep]{
    \lref{claim0} in the case $T\in I$\dott
  }
  We now prove \lref{claim0} in the case $T\notin I$\dott
  \startnewargseq
  \argument{
    the fact that $h=F\circ\theta$;
    the assumption that $F\in C^1(U,\R)$;
    the assumption that $\theta\in C^1(I,U)$;
    the assumption that for all $t\in I$ it holds that
    $\theta'(t)=-\bigl(\grad F\bigr)\bigl(\theta(t)\bigr)$;
    the chain rule;
  }{
    that for all $t\in I$ it holds that
    \begin{eqn}
      \llabel{eq1}
      h'(t)=\langle\bigl(\grad F\bigr)\bigl(\theta(t)\bigr),\theta'(t)\rangle
        =-\norm{\bigl(\grad F\bigr)\bigl(\theta(t)\bigr)}^2
        =-\norm{\theta'(t)}^2\leq0
      \dott
    \end{eqn}
  }
  \argument{
    \lref{eq1};
    \lref{setting1};
    the fact that $I$ is convex;
    the fact that $F$ is continuous;
  }{
    that $h$ is non-increasing and that for all $t\in I$
    it holds that
    \begin{eqn}
      \llabel{eq2}
      h(t)\geq\lim_{s\to T}h(s)
        =\lim_{j\to\infty}h(t_j)=F(b)
      \dott
    \end{eqn}
  }
  In the following we distinguish
  between the case $h^{-1}(\{F(b)\})\neq\emptyset$ and the case
  that for all $t\in I$ it holds that $h(t)>F(b)$\dott
  \startnewargseq
  We first prove \lref{claim0} in the case
  $h^{-1}(\{F(b)\})\neq\emptyset$\dott
  \startnewargseq
  \argument{  \lref{eq1};
    \lref{eq2};
    the fact that $h$ is non-increasing;
    the fact that $I$ is convex;
    the fundamental theorem of calculus;
  }{
    that for all $s,u\in I$ with $s\leq u$ and $h(s)=F(b)$
    it holds that
    \begin{eqn}
      \llabel{eq3}
      \int_s^u\norm{\theta'(t)}^2\dd t
        =h(s)-h(u)=0
      \dott
    \end{eqn}
  }
  \argument{
    \lref{eq3};
    the fact that $\theta'$ is continuous;
    the fundamental theorem of calculus;
  }{
    that for all $s,u\in I$ with $s\leq u$ and $h(s)=F(b)$
    it holds that
    \begin{eqn}
      \llabel{eq4}
      \theta(u)=\theta(s)+\int_s^u\theta'(t)\dd t=\theta(s)
      \dott
    \end{eqn}
  }
  \argument{
    \lref{eq4};
    \lref{setting1};
    the fact that $T=\sup(I)\notin I$;
    the assumption that $h^{-1}(\{F(b)\})\neq\emptyset$;
  }[verbs=ep]{
    \lref{claim0} in the case
    $h^{-1}(\{F(b)\})\neq\emptyset$\dott
  }
  We now prove \lref{claim0} in the case that for all $t\in I$
  it holds that $h(t)>F(b)$\dott
  \startnewargseq
  \argument{
    the fact that $F$ is locally definable;
    the fact that $b\in U$;
    the fact that $U$ is open;
  }{
    that there exist $\rho\in(0,\infty)$, $V\subseteq U$ which satisfy that
    \begin{eqn}[a]
      \llabel{eq5}
      \{x\in\R^n\colon\norm{x-b}\leq\rho\}&\subseteq U,\qquad
      V=\{x\in\R^n\colon\norm{x-b}<\rho\},
      \qquad\text{and}\qquad
      F|_V\text{ is definable}
      \dott
    \end{eqn}
  }
  In the following let $W\subseteq V$ satisfy
  \begin{eqn}
    \llabel{setting2}
    W=\{x\in V\colon F(x)>F(b)\}
  \end{eqn}
  and let $f\colon W\to(0,\infty)$ satisfy that for all $x\in W$
  it holds that
  \begin{eqn}
    \llabel{setting3}
    f(x)=F(x)-F(b)
    \dott
  \end{eqn}
  \startnewargseq
  \argument{
    \lref{eq5};
    \lref{setting2};
    \lref{setting3};
    the assumption that $F\in C^1(U,\R)$;
  }{
    that $W$ is open, bounded, and definable, that $f$ is
    definable, and that $f\in C^1(W,(0,\infty))$\dott
    \llabel{eq6}
  }
  \argument{
    \lref{eq6};
    \lref{setting2};
    \lref{setting3};
    \cite[Theorem~1]{K98};
    the chain rule;
  }{
    that there exist $c,\eta\in(0,\infty)$ and a strictly
    increasing $\psi\in C^1((0,\eta),(0,\infty))$ which
    satisfy for all $x\in V$ with
    $0<F(x)-F(b)<\eta$ that
    \begin{eqn}
      \llabel{eq7}
      \psi'(F(x)-F(b))\norm{\bigl(\grad F\bigr)\bigl(x\bigr)}\geq c
      \dott
    \end{eqn}
  }
  In the following let $\varphi\colon[0,\eta)\to[0,\infty)$ satisfy that
  $\varphi(0)=0$ and that for all $s\in(0,\eta)$ it holds that
  \begin{eqn}
    \llabel{setting4}
    \varphi(s)=\frac{\psi(s)-\inf\psi((0,\eta))}{c}
    \dott
  \end{eqn}
  \startnewargseq
  \argument{
    \lref{setting4};
    the fact that $\psi$ is positive and strictly increasing;
  }{
    that
    \begin{eqn}
      \llabel{eq8}
      \lim_{s\to 0}\psi(s)=\inf\psi((0,\eta))\in[0,\infty)
      \qquad\text{and}\qquad
      \lim_{s\to 0}\varphi(s)=0
      \dott
    \end{eqn}
  }
  \argument{
    \lref{eq8};
    \lref{eq7};
    \lref{setting4};
    the fact that $\psi\in C^1((0,\eta),(0,\infty))$;
    the fact that $c>0$; the fact that $\varphi(0)=0$;
  }{
    that $\varphi$ is continuous, that
    $\varphi|_{(0,\eta)}\in C^1((0,\eta),\R)$, and that for all
    $x\in V$ with $0<F(x)-F(b)<\eta$ it holds that
    \begin{eqn}
      \llabel{eq9}
      \varphi'(F(x)-F(b))\norm{\bigl(\grad F\bigr)\bigl(x\bigr)}\geq1
      \dott
    \end{eqn}
  }
  We claim that for all $s,u\in I$ with $s\leq u$,
  $0<h(s)-F(b)<\eta$, and
  $\norm{\theta(s)-b}+\varphi(h(s)-F(b))<\rho$ it holds that
  \begin{eqn}
    \llabel{claim1}
    \norm{\theta(u)-b}
      \leq\norm{\theta(s)-b}+\varphi(h(s)-F(b))
    \dott
  \end{eqn}
  \startnewargseq
  To prove \lref{claim1}, let $s,u\in I$, assume $s\leq u$,
  assume $0<h(s)-F(b)<\eta$, and assume
  $\norm{\theta(s)-b}+\varphi(h(s)-F(b))<\rho$\dott
  \startnewargseq
  \argument{
    the fact that $h$ is non-increasing;
    the assumption that for all $t\in I$ it holds that $h(t)>F(b)$;
    the fact that $h(s)-F(b)<\eta$;
    the fact that $I$ is convex;
  }{
    that $[s,u]\subseteq I$ and that for all $t\in[s,u]$
    it holds that
    \begin{eqn}
      \llabel{eq10}
      0<h(t)-F(b)\leq h(s)-F(b)<\eta
      \dott
    \end{eqn}
  }
  \argument{
    \lref{eq10};
    \lref{eq1};
    \lref{eq9};
    the fact that $h=F\circ\theta$;
    the fact that $\varphi|_{(0,\eta)}\in C^1((0,\eta),\R)$;
    the chain rule;
  }{
    that for all $t\in[s,u]$ with $\theta(t)\in V$
    it holds that
    \begin{eqn}
      \llabel{eq11}
      -\frac{\ddd}{\ddd t}\varphi(h(t)-F(b))
        =\varphi'(h(t)-F(b))\norm{\theta'(t)}^2
        \geq\norm{\theta'(t)}
      \dott
    \end{eqn}
  }
  \argument{
    \lref{eq11};
    the fact that $\varphi([0,\eta))\subseteq[0,\infty)$;
    the fundamental theorem of calculus;
    the triangle inequality;
  }{
    that for all $v\in[s,u]$ with
    $\theta([s,v])\subseteq V$ it holds that
    \begin{eqn}[a]
      \llabel{eq12}
      \norm{\theta(v)-b}
      &\leq\norm{\theta(s)-b}+\int_s^v\norm{\theta'(t)}\dd t
      \leq\norm{\theta(s)-b}+\varphi(h(s)-F(b))
        -\varphi(h(v)-F(b))
      \\
      &\leq\norm{\theta(s)-b}+\varphi(h(s)-F(b))<\rho
      \dott
    \end{eqn}
  }
  \argument{
    \lref{eq12};
    \lref{eq5};
    the fact that $\norm{\theta(s)-b}<\rho$;
    the fact that $\theta$ is continuous;
  }[verbs=ep]{
    that
    \begin{eqn}
      \llabel{eq15}
      \theta([s,u])\subseteq V
      \dott
    \end{eqn}
  }
  \argument{
    \lref{eq15};
    \lref{eq12};
    the fact that $s,u\in I$ with $s\leq u$,
    $0<h(s)-F(b)<\eta$, and
    $\norm{\theta(s)-b}+\varphi(h(s)-F(b))<\rho$ were arbitrary;
  }[verbs=ep]{
    \lref{claim1}\dott
  }

  \startnewargseq
  \argument{
    \lref{setting1};
    \lref{eq2};
    \lref{eq8};
    the assumption that for all $t\in I$ it holds that $h(t)>F(b)$;
  }{
    that there exists $J\in\N$ which satisfies for all
    $j\in\N\cap[J,\infty)$ that
    \begin{eqn}
      \llabel{eq16}
      0<h(t_j)-F(b)<\eta
      \qquad\text{and}\qquad
      \norm{\theta(t_j)-b}+\varphi(h(t_j)-F(b))<\rho
      \dott
    \end{eqn}
  }
  \startnewargseq
  \argument{  \lref{claim1};
    \lref{eq16};
    the fact that $T=\sup(I)\notin I$;
  }{
    that for all $j\in\N\cap[J,\infty)$ it holds that
    \begin{eqn}
      \llabel{eq17}
      0\leq\limsup_{t\to T}\norm{\theta(t)-b}
        \leq\norm{\theta(t_j)-b}+\varphi(h(t_j)-F(b))
      \dott
    \end{eqn}
  }
  \argument{
    \lref{setting1};
    \lref{eq2};
  }{
    that
    \begin{eqn}[a]
      \llabel{eq17.1}
      \lim_{j\to\infty}\norm{\theta(t_j)-b}
        =0
        \qquad\text{and}
      \qquad
      \lim_{j\to\infty}(h(t_j)-F(b))
        =0
      \dott
    \end{eqn}
  }
  \argument{
    \lref{eq17.1};\lref{eq2};
    \lref{eq8};
  }{
    that
    \begin{eqn}
      \llabel{eq17.2}
      \lim_{j\to\infty}\varphi(h(t_j)-F(b))
        =\lim_{s\to 0}\varphi(s)=0
      \dott
    \end{eqn}
  }
  \argument{
    \lref{eq17.2};
    \lref{eq17.1};
  }{
    that
    \begin{eqn}
      \llabel{eq18}
      \lim_{j\to\infty}
        \bigl(\norm{\theta(t_j)-b}+\varphi(h(t_j)-F(b))\bigr)=0
      \dott
    \end{eqn}
  }
  \argument{
    \lref{eq18};
    \lref{eq17};
  }[verbs=ep]{
    \lref{claim0} in the case that for all $t\in I$ it holds that
    $h(t)>F(b)$\dott
  }
  This completes the proof of \lref{claim0}\dott
\end{aproof}

\begin{athm}{corollary}{cor:spherical}[Finite spherical length for \GF\ trajectories of locally definable potentials]
Let $ n \in \mathbb{N} $,
$ \vartheta \in \R^n $,
let $ U \subseteq \R^n $ be open,
let $ F \in C^1( U, \R ) $ be locally definable,
let $ I \subseteq \R $ be convex with nonempty interior,
let
$ \theta \in C^1( I, U ) $ satisfy
\begin{eqn}
\label{eq:first_theorem_ODE_B}
  \theta' = - (\nabla F) \circ \theta
  ,
\end{eqn}
let $ a, b \in U $ satisfy
$
  \liminf_{ t \to \inf( I ) } | \theta(t) - a |
  =
  \liminf_{ t \to \sup( I ) } | \theta(t) - b |
  = 0
$,
and let
$ P \colon U \to \R^n $
satisfy for all $ v \in U \backslash \{ \vartheta \} $
that
$
  P( v ) = | v - \vartheta |^{ - 1 } ( v - \vartheta )
$.
Then $ P \circ \theta $ has finite length.
\end{athm}

\begin{aproof}
  Throughout this proof let $G\colon U\to\R$ satisfy $G=-F$,
  let $J\subseteq\R$, $K\subseteq U$ satisfy
  \begin{eqn}
    \llabel{setting1}
    J=\{-t\colon t\in I\}
    \qquad\text{and}\qquad K=\theta(I)\cup\{a,b\},
  \end{eqn}
  and let $\gamma\colon J\to U$ satisfy for all $u\in J$ that
  \begin{eqn}
    \llabel{setting2}
    \gamma(u)=\theta(-u)
    \dott
  \end{eqn}
  \cfload[.]
  \startnewargseq
  \argument{
    \cref{lem:locally-definable-endpoint};
    the assumption that $U$ is open;
    the assumption that $F\in C^1(U,\R)$;
    the assumption that $F$ is locally definable;
    the assumption that $I$ is convex with nonempty interior;
    the assumption that $\theta\in C^1(I,U)$;
    the assumption that for all $t\in I$ it holds that
    $\theta'(t)=-\bigl(\grad F\bigr)\bigl(\theta(t)\bigr)$;
    the fact that $b\in U$;
    the assumption that
    $\liminf_{t\to\sup(I)}\norm{\theta(t)-b}=0$;
  }{
    that
    \begin{eqn}
      \llabel{eq1}
      \limsup_{t\to\sup(I)}\abs{\theta(t)-b}=0
      \dott
    \end{eqn}
  }
  \argument{
    \lref{setting1};
    \lref{setting2};
    the fact that $G=-F$;
    the assumption that $I$ is convex with nonempty interior;
    the assumption that $F\in C^1(U,\R)$;
    the assumption that $F$ is locally definable;
    the assumption that $\theta\in C^1(I,U)$;
    the assumption that for all $t\in I$ it holds that
    $\theta'(t)=-\bigl(\grad F\bigr)\bigl(\theta(t)\bigr)$;
  }{
    that $J$ is convex with nonempty interior, that
    $G\in C^1(U,\R)$, that $G$ is locally definable, that
    $\gamma\in C^1(J,U)$, and that for all $u\in J$ it holds that
    \begin{eqn}
      \llabel{eq2}
      \gamma'(u)=-\theta'(-u)
        =\bigl(\grad F\bigr)\bigl(\theta(-u)\bigr)=-\bigl(\grad G\bigr)\bigl(\gamma(u)\bigr)
      \dott
    \end{eqn}
  }
  \argument{
    \lref{setting1};
    \lref{setting2};
    the assumption that
    $\liminf_{t\to\inf(I)}\norm{\theta(t)-a}=0$;
  }{
    that
    \begin{eqn}[a]
      \llabel{eq3}
      \sup(J)=-\inf(I)
      \qquad\text{and}
      \qquad
      \liminf_{u\to\sup(J)}\norm{\gamma(u)-a}
        =\liminf_{t\to\inf(I)}\norm{\theta(t)-a}=0
      \dott
    \end{eqn}
  }
  \argument{
    \lref{eq3};
    \lref{setting1};
    \lref{setting2};
    \lref{eq2};
    the fact that $J$ is convex with nonempty interior;
    the fact that $G\in C^1(U,\R)$;
    the fact that $G$ is locally definable;
    the fact that $\gamma\in C^1(J,U)$;
    the fact that $a\in U$;
    \cref{lem:locally-definable-endpoint}
    (applied with $F\curvearrowleft G$, $I\curvearrowleft J$,
    $\theta\curvearrowleft\gamma$, $b\curvearrowleft a$
    in the notation of \cref{lem:locally-definable-endpoint});
  }{
    that
    \begin{eqn}
      \llabel{eq4}
      \limsup_{t\to\inf(I)}\abs{\theta(t)-a}
        =\limsup_{u\to\sup(J)}\abs{\gamma(u)-a}=0
      \dott
    \end{eqn}
  }
  \argument{
    \lref{eq4};
    \lref{setting1};
    \lref{eq1};
    the fact that $\theta$ is continuous;
  }{
    that $K$ is compact\dott
    \llabel{arg1}
  }
  \argument{
    \lref{arg1};
    the fact that $K\subseteq U$;
    the assumption that $U$ is open;
    the assumption that $F$ is locally definable;
  }{
    that there exist $N\in\N$, $x=(x_1,\dots,x_N)\in K^N$,
    $\rho=(\rho_1,\dots,\rho_N)\in(0,\infty)^N$, and open 
    $B_j\subseteq U$, $j\in\{1,\dots,N\}$, which satisfy that
    \begin{enumerate}[label=(\roman*)]
      \item for all $j\in\{1,\dots,N\}$ it holds that
      $B_j=\{y\in\R^n\colon\norm{y-x_j}<\rho_j\}$
      and that $F|_{B_j}$ is definable and
      \item it holds that
      \begin{eqn}
        \llabel{eq6}
        K\subseteq\midcup_{j=1}^N B_j\subseteq U
        \dott
      \end{eqn}
    \end{enumerate}
  }
  In the following let $V\subseteq U$ satisfy
  \begin{eqn}
    \llabel{setting3}
    V=\midcup_{j=1}^N B_j
    \dott
  \end{eqn}
  \startnewargseq
  \argument{
    \lref{setting3};
  }{
    that
    \begin{eqn}
      \llabel{eq7}
      \operatorname{Graph}(F|_V)
        =\midcup_{j=1}^N\operatorname{Graph}(F|_{B_j})
      \dott
    \end{eqn}
  }
  \argument{
    \lref{eq7};
    \lref{setting1};
    \lref{eq6};
    \lref{setting3};
    the fact that for all $j\in\{1,\dots,N\}$ it holds that
    $B_j$ is open and that $F|_{B_j}$ is definable;
    the assumption that $F\in C^1(U,\R)$;
    the assumption that $\theta\in C^1(I,U)$;
  }{
    that $V$ is open, that $F|_V\in C^1(V,\R)$,
    that $F|_V$ is definable, and that $\theta\in C^1(I,V)$\dott
    \llabel{arg2}
  }
  \argument{
    \lref{arg2};
    \lref{setting3};
    the assumption that $I$ is convex with nonempty interior;
    the assumption that for all $t\in I$ it holds that
    $\theta'(t)=-\bigl(\grad F\bigr)\bigl(\theta(t)\bigr)$;
    the assumption that for all $v\in U\backslash\{\vartheta\}$
    it holds that $P(v)=(v-\vartheta)/\norm{v-\vartheta}$;
    \cref{maintheorem3}
    (applied with $U\curvearrowleft V$, $F\curvearrowleft F|_V$,
    $P\curvearrowleft P|_V$ in the notation of \cref{maintheorem3});
  }[verbs=ep]{
    that
    \begin{eqn}
      \curvelength{P\circ\theta}
        =\curvelength{(P|_V)\circ\theta}<\infty
      \dott
    \end{eqn}
  }
\end{aproof}

\subsection{Infinite spherical length for analytic potentials}
\label{sec:analytic-counterexample}

\cref{maintheorem} above, in particular, proves for every definable $C^1$-potential $ F \colon \R^n \to \R $
that the conclusion of the gradient conjecture at infinity applies.
In the following elementary results, \cref{cor:analytic-graph-counterexample}
and \cref{prop:analytic-counterexample} below,
we provide explicit examples of analytic potentials $ F \colon \R^2 \to \R $
for which the conclusion of the gradient conjecture at infinity fails to hold (cf., \eg, also \cite[Theorem~1.1, page 592]{DHL22_MR4453693}).
In particular, for every o-minimal structure $ \mathcal{S} $ we have that
$ F \colon \R^2 \to \R $ in \cref{cor:analytic-graph-counterexample}
and $ F \colon \R^2 \to \R $ in \cref{prop:analytic-counterexample}
fail to be $ \mathcal{S} $-definable. However, we have that
$ F $ in \cref{cor:analytic-graph-counterexample}
and $ F $ in \cref{prop:analytic-counterexample}
are locally $ \Ran $-definable (cf.\ \cref{maincor} above).

\begin{athm}{lemma}{lem:graph-gradient-trajectory}
  [Graphs of curves as gradient trajectories]
  Let $n\in\N$, $g\in C^2(\R,\R^n)$, let
  $F\colon\R^{n+1}\to\R$ satisfy for all $x\in\R$, $y\in\R^n$ that
  \begin{eqn}
    \llabel{ass1}
    F(x,y)=\langle g'(x),y\rangle
      +\int_0^x\bigl(1-\langle g''(s),g(s)\rangle\bigr)\dd s,
  \end{eqn}
  and let $\theta\colon\R\to\R^{n+1}$ satisfy for all $t\in\R$ that
  $\theta(t)=(t,g(t))$.
  Then
  \begin{enumerate}[label=(\roman*)]
    \item \llabel{item1} it holds that $F\in C^1(\R^{n+1},\R)$,
    \item \llabel{item2} it holds for all $t\in\R$ that
    $\theta'(t)=(\nabla F)(\theta(t))$, and
    \item \llabel{item3} it holds for all $x\in\R$, $y\in\R^n$ that
    $F(x,y)=x+\int_0^x\norm{g'(s)}^2\dd s
      +\langle g'(x),y-g(x)\rangle+\langle g'(0),g(0)\rangle$.
  \end{enumerate}
\end{athm}
\begin{aproof}
  \argument{
    \lref{ass1};
    the fact that $g\in C^2(\R,\R^n)$;
    the fundamental theorem of calculus;
    the product rule;
  }[verbs=i]{
    that $F\in C^1(\R^{n+1},\R)$ and that for all $x\in\R$, $y\in\R^n$
    it holds that
    \begin{eqn}
      \llabel{eq1}
      (\partial_xF)(x,y)
        =\langle g''(x),y\rangle+1-\langle g''(x),g(x)\rangle
        =1+\langle g''(x),y-g(x)\rangle
    \end{eqn}
    \begin{eqn}
      \llabel{eq2}
      \text{and}\qquad (\nabla_yF)(x,y)=g'(x).
    \end{eqn}
  }
  \argument{
    \lref{eq1};
    \lref{eq2};
  }[verbs=s]{
    that for all $t\in\R$ it holds that
    \begin{eqn}
      (\nabla F)(\theta(t))=(1,g'(t))=\theta'(t).
    \end{eqn}
  }
  \argument{
    \lref{ass1};
    integration by parts;
  }[verbs=ep]{
    that for all $x\in\R$, $y\in\R^n$ it holds that
    \begin{eqn}
      F(x,y)=\langle g'(x),y\rangle+x
        -\bigl[\langle g'(s),g(s)\rangle\bigr]_{s=0}^{s=x}
        +\int_0^x\langle g'(s),g'(s)\rangle\dd s.
    \end{eqn}
  }
\end{aproof}

\begin{athm}{corollary}{cor:analytic-graph-counterexample}
  Let $g\colon\R\to\R$ and $F\colon\R^2\to\R$ satisfy for all
  $x,y\in\R$ that $g(x)=x\sin(x)$ and
  \begin{eqn}
    \llabel{ass1}
    F(x,y)=-g'(x)y-\int_0^x\bigl(1-g''(s)g(s)\bigr)\dd s,
  \end{eqn}
  let $\theta\colon[0,\infty)\to\R^2$ satisfy for all
  $t\in[0,\infty)$ that
  $\theta(t)=(t,t\sin(t))$,
  and let $\omega\colon(0,\infty)\to\R^2$ satisfy for all
  $t\in(0,\infty)$ that $\omega(t)=\norm{\theta(t)}^{-1}\theta(t)$.
  Then
  \begin{enumerate}[label=(\roman*)]
    \item \llabel{item1} it holds that $F$ is real analytic,
    \item \llabel{item2} it holds for all $t\in[0,\infty)$ that
    $\theta'(t)=-(\nabla F)(\theta(t))$, and
    \item \llabel{item3} it does not hold that there exists $w\in\R^2$ such that
    $\limsup_{t\to\infty}\norm{w-\omega(t)}=0$.
  \end{enumerate}
\end{athm}
\begin{aproof}
  \argument{
    \lref{ass1};
    \cref{lem:graph-gradient-trajectory}
    (applied with $n\curvearrowleft1$, $F\curvearrowleft-F$
    in the notation of \cref{lem:graph-gradient-trajectory});
    the fact that $g$ is real analytic;
  }{
    that $F$ is real analytic and that for all $t\in[0,\infty)$ it holds that
    $\theta'(t)=-(\nabla F)(\theta(t))$.
  }
  \startnewargseq
  \setnote
  \argument{
    the fact that for all $t\in[0,\infty)$ it holds that
    $\theta(t)=(t,t\sin(t))$;
    the fact that for all $t\in(0,\infty)$ it holds that
    $\omega(t)=\norm{\theta(t)}^{-1}\theta(t)$;
  }[verbs=d]{
    that for all $t\in(0,\infty)$ it holds that
    \begin{eqn}
      \llabel{eq1}
      \omega(t)=\frac{(1,\sin(t))}{\sqrt{1+\sin^2(t)}}.
    \end{eqn}
  }
  \argument{
    \lref{eq1};
  }{
    that for all $k\in\N$ it holds that
    \begin{eqn}
      \llabel{eq2}
      \omega\bigl(2\pi k+\tfrac{\pi}{2}\bigr)=\frac{(1,1)}{\sqrt2}
      \qquad\text{and}\qquad
      \omega\bigl(2\pi k+\tfrac{3\pi}{2}\bigr)=\frac{(1,-1)}{\sqrt2}.
    \end{eqn}
  }
  \argument{
    \lref{eq2};
    the fact that $\lim_{k\to\infty}2\pi k=\infty$;
  }{
    that it does not hold that there exists $w\in\R^2$ such that
    $\limsup_{t\to\infty}\norm{w-\omega(t)}=0$.
  }
\end{aproof}

\begin{athm}{proposition}{prop:analytic-counterexample}
  Let $F\colon\R^2\to\R$ satisfy for all $x,y\in\R$ that
  \begin{eqn}[a]
    \llabel{ass1}
    F(x,y)&=-(x^2+y^2)-\frac{1}{3}(x^2+y^2)^3
    \\
    &\quad +(x^2+y^2)
      \bigl[(x^2-y^2)\sin(x^2+y^2)-2xy\cos(x^2+y^2)\bigr],
  \end{eqn}
  let $\theta\colon[0,\infty)\to\R^2$ satisfy for all
  $t\in[0,\infty)$ that
  \begin{eqn}
    \llabel{ass2}
    \theta(t)=e^{2t}
      \bigl(\cos(\tfrac12e^{4t}),\sin(\tfrac12e^{4t})\bigr),
  \end{eqn}
  and let $P\colon\R^2\to\R^2$ satisfy for all
  $v\in\R^2\backslash\{0\}$ that
  $P(v)=\norm{v}^{-1}v$\dott\cfload[.]
  Then
  \begin{enumerate}[label=(\roman*)]
    \item \llabel{item1} it holds that $F$ and $\theta$ are analytic,
    \item \llabel{item2} it holds for all $t\in[0,\infty)$ that
    $\norm{\theta(t)}=e^{2t}$,
    \item \llabel{item3} it holds for all $t\in[0,\infty)$ that
    $(P\circ\theta)(t)
      =\bigl(\cos(\tfrac12e^{4t}),\sin(\tfrac12e^{4t})\bigr)$,
    \item \llabel{item4} it holds for all $v\in\R^2$ that
    $\limsup_{t\to\infty}\norm{(P\circ\theta)(t)-v}>0$,
    \item \llabel{item5} it holds that
    $\curvelength{P\circ\theta}=\infty$, and
    \item \llabel{item6} it holds for all $t\in[0,\infty)$ that
    $\theta'(t)=-(\grad F)(\theta(t))$\dott
  \end{enumerate}
  \cfout[.]
\end{athm}
\begin{aproof}
  Throughout this proof let $E,T\colon\R\to\R^2$ satisfy for all
  $t\in\R$ that
  \begin{eqn}
    \llabel{setting1}
    E(t)=(\cos(t),\sin(t))
    \qquad\text{and}\qquad
    T(t)=(-\sin(t),\cos(t)),
  \end{eqn}
  let $r,\alpha\colon[0,\infty)\to(0,\infty)$ satisfy for all
  $t\in[0,\infty)$ that
  \begin{eqn}
    \llabel{setting2}
    r(t)=e^{2t}
    \qquad\text{and}\qquad \alpha(t)=\frac12e^{4t},
  \end{eqn}
  let $H=(H(\rho,t))_{(\rho,t)\in(0,\infty)\times\R}
  \colon(0,\infty)\times\R\to\R$ satisfy for all
  $\rho\in(0,\infty)$, $t\in\R$ that
  \begin{eqn}
    \llabel{setting3}
    H(\rho,t)=F(\rho E(t)),
  \end{eqn}
  and let $(s_j)_{j\in\N}\subseteq[0,\infty)$ satisfy for all
  $j\in\N$ that
  \begin{eqn}
    \llabel{setting4}
    s_j=\frac14\log(2\pi j)
    \dott
  \end{eqn}
  \cfload[.]
  \startnewargseq
  \argument{
    \lref{ass1};
    \lref{ass2};
  }[verbs=ep]{
    \lref{item1}\dott
  }
  \startnewargseq
  \argument{
    \lref{ass2};\lref{setting2};
    the fact that for all $t\in\R$ it holds that
    $(\cos(t))^2+(\sin(t))^2=1$;
  }{
    that for all $t\in[0,\infty)$ it holds that
    \begin{eqn}
      \llabel{eq1}
      \norm{\theta(t)}
        =r(t)\sqrt{(\cos(\alpha(t)))^2+(\sin(\alpha(t)))^2}
        =r(t)=e^{2t}>0
      \dott
    \end{eqn}
  }
  \argument{
    \lref{eq1};
  }[verbs=ep]{
    \lref{item2}\dott
  }
  \startnewargseq
  \argument{
    \lref{item2};
    \lref{ass2};
    the assumption that for all $v\in\R^2\backslash\{0\}$
    it holds that $P(v)=\norm{v}^{-1}v$;
  }[verbs=ep]{
    \lref{item3}\dott
  }
  \startnewargseq
  \argument{
    \lref{item3};
    \lref{setting2};
    \lref{setting4};
    the fact that for all $j\in\N$ it holds that
    $\cos(\pi j)=(-1)^j$ and $\sin(\pi j)=0$;
  }{
    that for all $j\in\N$ it holds that
    \begin{eqn}
      \llabel{eq2}
      0<s_j<s_{j+1},
      \qquad \alpha(s_j)=\pi j,
      \qquad (P\circ\theta)(s_j)=((-1)^j,0),
      \qquad\text{and}\qquad
      \lim_{k\to\infty}s_k=\infty
      \dott
    \end{eqn}
  }
  \argument{
    \lref{eq2};
    the triangle inequality;
  }{
    that for all $v\in\R^2$ it holds that
    \begin{eqn}[a]
      \llabel{eq4}
      \limsup_{t\to\infty}\norm{(P\circ\theta)(t)-v}
        \geq\max\{\norm{(1,0)-v},\norm{(-1,0)-v}\}
      \geq\frac12\norm{(1,0)-(-1,0)}=1>0
      \dott
    \end{eqn}
  }
  \argument{
    \lref{eq4};
  }[verbs=ep]{
    \lref{item4}\dott
  }
  \startnewargseq
  \argument{
    \lref{eq2};
  }{
    that for all $m\in\N$ it holds that
    \begin{eqn}
      \llabel{eq5}
      \curvelength{P\circ\theta}
        \geq\sum_{j=1}^m
          \norm{(P\circ\theta)(s_{j+1})-(P\circ\theta)(s_j)}
        =\sum_{j=1}^m2=2m
      \dott
    \end{eqn}
  }
  \argument{
    \lref{eq5};
  }[verbs=ep]{
    \lref{item5}\dott
  }
  \startnewargseq
  \argument{
    \lref{setting1};
    the fact that for all $t\in\R$ it holds that
    $(\cos(t))^2+(\sin(t))^2=1$;
  }{
    that for all $t\in\R$ it holds that
    \begin{eqn}
      \llabel{eq6}
      E'(t)=T(t),
      \qquad \norm{E(t)}=\norm{T(t)}=1,
      \qquad\text{and}\qquad
      \langle E(t),T(t)\rangle=0
      \dott
    \end{eqn}
  }
  \argument{
    \lref{ass2};
    \lref{setting1};
    \lref{setting2};
    the chain rule;
  }{
    that for all $t\in[0,\infty)$ it holds that
    \begin{eqn}[a]
      \llabel{eq7}
      \theta(t)=r(t)E(\alpha(t)),
      \qquad 2\alpha(t)=(r(t))^2,\qquad
      r'(t)=2r(t),
      \qquad\text{and}\qquad
      \alpha'(t)=2e^{4t}=2(r(t))^2
      \dott
    \end{eqn}
  }
  \argument{
    \lref{ass1};
    \lref{setting1};
    \lref{setting3};
  }{
    that for all $\rho\in(0,\infty)$, $t\in\R$ it holds that
    \begin{eqn}[a]
      \llabel{eq8}
      H(\rho,t)
        &=-\rho^2-\frac{\rho^6}{3}
          +\rho^4\bigl(\cos(2t)\sin(\rho^2)
            -\sin(2t)\cos(\rho^2)\bigr)
      \\
      &=-\rho^2-\frac{\rho^6}{3}
          -\rho^4\sin(2t-\rho^2)
      \dott
    \end{eqn}
  }
  \argument{
    \lref{eq8};
    the chain rule;
  }{
    that for all $\rho\in(0,\infty)$, $t\in\R$ it holds that
    \begin{equation}
      \llabel{eq9}
      \bigl(\partial_\rho H\bigr)\bigl(\rho,t\bigr)
        =-2\rho-2\rho^5-4\rho^3\sin(2t-\rho^2)
          +2\rho^5\cos(2t-\rho^2)\quad
      \text{and}\quad\bigl(\partial_t H\bigr)\bigl(\rho,t\bigr)
        =-2\rho^4\cos(2t-\rho^2)
      \dott
    \end{equation}
  }
  \argument{
    \lref{eq9};
    \lref{eq7};
  }{
    that for all $u\in[0,\infty)$ it holds that
    \begin{eqn}[a]
      \llabel{eq10}
    \textstyle  \bigl(\partial_\rho H\bigr)\bigl(r(u),\alpha(u)\bigr)
        =-2r(u)-2(r(u))^5+2(r(u))^5=-2r(u)\quad\text{and}\quad
      \bigl(\partial_t H\bigr)\bigl(r(u),\alpha(u)\bigr)
        =-2(r(u))^4
      \dott
    \end{eqn}
  }
  \argument{
    \lref{item1};
    \lref{setting3};
    \lref{eq6};
    the chain rule;
  }{
    that for all $\rho\in(0,\infty)$, $t\in\R$ it holds that
    \begin{eqn}[a]
      \llabel{eq11}
      \bigl(\partial_\rho H\bigr)\bigl(\rho,t\bigr)
        =\langle\bigl(\grad F\bigr)\bigl(\rho E(t)\bigr),E(t)\rangle\qquad\text{and}\qquad
      \bigl(\partial_t H\bigr)\bigl(\rho,t\bigr)
        =\rho\langle\bigl(\grad F\bigr)\bigl(\rho E(t)\bigr),T(t)\rangle
      \dott
    \end{eqn}
  }
  \argument{
    \lref{eq11};
    \lref{eq6};
    the fact that for all $t\in\R$ it holds that
    $E(t),T(t)$ form an orthonormal basis of $\R^2$;
  }{
    that for all $\rho\in(0,\infty)$, $t\in\R$ it holds that
    \begin{eqn}
      \llabel{eq12}
      \bigl(\grad F\bigr)\bigl(\rho E(t)\bigr)
        =\bigl(\partial_\rho H\bigr)\bigl(\rho,t\bigr)E(t)
          +\frac1\rho\bigl(\partial_t H\bigr)\bigl(\rho,t\bigr)
            T(t)
      \dott
    \end{eqn}
  }
  \argument{
    \lref{eq12};
    \lref{eq7};
    \lref{eq10};
  }{
    that for all $u\in[0,\infty)$ it holds that
    \begin{eqn}
      \llabel{eq13}
      -\bigl(\grad F\bigr)\bigl(\theta(u)\bigr)
        =2r(u)E(\alpha(u))+2(r(u))^3T(\alpha(u))
      \dott
    \end{eqn}
  }
  \argument{
    \lref{eq6};
    \lref{eq7};
    the chain rule;
  }{
    that for all $u\in[0,\infty)$ it holds that
    \begin{eqn}[a]
      \llabel{eq14}
      \theta'(u)
        =r'(u)E(\alpha(u))+r(u)\alpha'(u)T(\alpha(u))
      =2r(u)E(\alpha(u))+2(r(u))^3T(\alpha(u))
      \dott
    \end{eqn}
  }
  \argument{
    \lref{eq14};
    \lref{eq13};
  }[verbs=ep]{
    \lref{item6}\dott
  }
\end{aproof}

\subsection{Finite spherical length and directional convergence for DNN training flows}
\label{sec:directional-convergence}

In \cref{cor:neural-network-training-full} below we establish finite spherical
length and directional convergence of \DNN\ training flows.
Our proof of \cref{cor:neural-network-training-full} uses the elementary and
well-known fact, which we recall in the following
\cref{lem:finite-length-endpoint-limit}, that every curve with finite
length is also convergent.
Only for completeness, we include a proof for
\cref{lem:finite-length-endpoint-limit} here.

\begin{athm}{lemma}{lem:finite-length-endpoint-limit}[Convergence of curves with finite length]
Let $n\in\N$,
let $U\subseteq\R^n$ be a set,
let $I\subseteq\R$ be convex with nonempty interior,
and let $\gamma\colon I\to U$ have finite length\cfadd{definition:curve-length}\dott\cfload[.]
Then there exists a unique $v\in\overline{U}$ such that
\begin{eqn}[d]
  \llabel{concl1}
  \limsup_{t\to\sup(I)}\norm{\gamma(t)-v}=0
  \dott
\end{eqn}
\cfout[.]
\end{athm}

\begin{aproof}
  Throughout this proof let $a\in\operatorname{int}(I)$,
  $b\in\R\cup\{\infty\}$, $J\subseteq I$ satisfy
  \begin{eqn}
    \llabel{setting1}
    b=\sup(I)
    \qquad\text{and}\qquad J=I\cap[a,b),
  \end{eqn}
  let $\ell\colon J\to[0,\infty)$ satisfy for all $t\in J$ that
  \begin{eqn}
    \llabel{setting2}
    \ell(t)=\curvelength{\gamma|_{[a,t]}},
  \end{eqn}
  and let $L\in[0,\infty]$ satisfy
  \begin{eqn}
    \llabel{setting3}
    L=\sup\{\ell(t)\colon t\in J\}
    \dott
  \end{eqn}
  \cfload[.]\startnewargseq
  \argument{
    \lref{setting1};
    the fact that $a\in\operatorname{int}(I)$;
    the fact that $I$ is convex;
  }{
    that
    \begin{eqn}
      \llabel{eq1}
      a<b,
      \qquad J=[a,b),
      \qquad\text{and}\qquad \sup(J)=b
      \dott
    \end{eqn}
  }
  \argument{
    \lref{eq1};
    \lref{setting2};
    the assumption that $\gamma$ has finite length;
  }{
    that for all $u,t\in J$ with $u\leq t$ it holds that
    \begin{eqn}
      \llabel{eq2}
      0=\ell(a)\leq\ell(u)\leq\ell(t)\leq\curvelength{\gamma}<\infty
      \dott
    \end{eqn}
  }
  \argument{
    \lref{eq2};
    \lref{setting2};
  }{
    that for all $u,t\in J$ with $u\leq t$ it holds that
    \begin{eqn}
      \llabel{eq3}
      \norm{\gamma(t)-\gamma(u)}\leq\ell(t)-\ell(u)
      \dott
    \end{eqn}
  }
  \argument{
    \lref{setting3};
    \lref{eq1};
    \lref{eq2};
  }{
    that
    \begin{eqn}
      \llabel{eq4}
      L\in[0,\curvelength{\gamma}]
      \qquad\text{and}\qquad
      \lim_{t\to b}\ell(t)=L
      \dott
    \end{eqn}
  }
  \argument{
    \lref{eq4};
    \lref{eq1};
  }{
    that there exists $(s_j)_{j\in\N}\subseteq J$ which satisfies
    for all $j\in\N$ that
    \begin{eqn}
      \llabel{eq5}
      s_j<s_{j+1},
      \qquad \lim_{k\to\infty}s_k=b,
      \qquad\text{and}\qquad
      \lim_{k\to\infty}\ell(s_k)=L
      \dott
    \end{eqn}
  }
  \startnewargseq
  \argument{
    \lref{eq2};
    \lref{eq3};
    \lref{eq4};  \lref{eq5};
  }{
    that for all $\varepsilon\in(0,\infty)$ there exists $K\in\N$
    such that for all $j,k\in\N\cap[K,\infty)$ it holds that
    \begin{eqn}[a]
      \llabel{eq6}
      \norm{\gamma(s_j)-\gamma(s_k)}
        \leq\ell(\max\{s_j,s_k\})-\ell(\min\{s_j,s_k\})
      \leq L-\ell(s_K)<\varepsilon
      \dott
    \end{eqn}
  }
  \argument{
    \lref{eq6};
    the fact that for all $j\in\N$ it holds that $\gamma(s_j)\in U$;
  }{
    that there exists $v\in\overline{U}$ which satisfies
    \begin{eqn}
      \llabel{eq7}
      \lim_{j\to\infty}\gamma(s_j)=v
      \dott
    \end{eqn}
  }
  \startnewargseq
  \argument{
    \lref{eq1};
    \lref{eq3};
    \lref{eq5};\lref{eq7};
  }{
    that for all $t\in J$ it holds that
    \begin{eqn}[a]
      \llabel{eq8}
      \norm{\gamma(t)-v}
        =\lim_{j\to\infty}\norm{\gamma(t)-\gamma(s_j)}
        \leq\lim_{j\to\infty}\bigl(\ell(s_j)-\ell(t)\bigr)
        =L-\ell(t)
      \dott
    \end{eqn}
  }
  \argument{
    \lref{eq8};
    \lref{setting1};
    \lref{eq1};
    \lref{eq4};
  }{
    that
    \begin{eqn}
      \llabel{eq9}
      0\leq\limsup_{t\to b}\norm{\gamma(t)-v}
        \leq\lim_{t\to b}\bigl(L-\ell(t)\bigr)=0
      \dott
    \end{eqn}
  }
  \argument{
    \lref{eq9};
    \lref{setting1};
  }[verbs=ep]{
    that 
    \begin{eqn}
      \llabel{eq9.1}
      \limsup_{t\to\sup(I)}\norm{\gamma(t)-v}=0
      \dott
    \end{eqn}
  }
  \argument{
    \lref{eq9.1};
    \lref{setting1};
    the triangle inequality;
  }{
    that for all $w\in\overline{U}$ with
    $\limsup_{t\to b}\norm{\gamma(t)-w}=0$ it holds that
    \begin{eqn}
      \llabel{eq10}
      0\leq\norm{v-w}
        \leq\limsup_{t\to b}
          \bigl(\norm{v-\gamma(t)}+\norm{\gamma(t)-w}\bigr)=0
      \dott
    \end{eqn}
  }
  \argument{
    \lref{eq10};
    \lref{setting1};
  }[verbs=ep]{
    that $v$ is the unique element of $\overline{U}$ which satisfies
    \begin{eqn}
      \limsup_{t\to\sup(I)}\norm{\gamma(t)-v}=0
      \dott
    \end{eqn}
  }
\end{aproof}

\begin{athm}{corollary}{cor:neural-network-training-full}[Finite spherical length for \DNN\ training flows]
Let $ m, L \in \N $, $ \ell = ( \ell_0, \ell_1, \dots, \ell_L ) \in \N^{ L + 1 } $,
let
$
  n =
  \sum_{ i = 1 }^L
  \ell_i ( \ell_{ i - 1 } + 1 )
$,
let $ A \in C^1( \R, \R ) $, $ \mathscr{L} \in C^1( \R^{ \ell_L } \times \R^{ \ell_L } \times \R^n , [ 0, \infty ) ) $ be definable,
for every $ \theta = ( \theta_1, \dots, \theta_n ) \in \R^n $ let
$
  \rel{ v }{ \theta }
  = ( \rell{ v }{ \theta }{ 1 }, \dots, \rell{ v }{ \theta }{ \ell_v } )
  \colon \R^{ \ell_0 } \to \R^{ \ell_v }
$,
$
  v \in \{ 0, 1, \dots, L \}
$,
satisfy for all $ v \in \{ 0, 1, \dots, L - 1 \} $,
$ x = ( x_1, \dots, x_{ \ell_0 } ) \in \R^{ \ell_0 } $,
$ i \in \{ 1, 2, \dots, \ell_{ v + 1 } \} $
that
\begin{equation}
\label{eq:neural_network_def-full}
\begin{split}
  \rell{ v + 1 }{ \theta }{ i }( x )
&
  =
  \theta_{ \ell_{ v + 1 } \ell_{ v } + i + \sum_{ h = 1 }^v \ell_h ( \ell_{ h - 1 } + 1 ) }
\\ & \quad
  + \textstyle\sum\limits_{ j = 1 }^{ \ell_v }
  \theta_{ ( i - 1 ) \ell_v + j + \sum_{ h = 1 }^v \ell_h ( \ell_{ h - 1 } + 1 ) }
  \Bigl(
    x_{ \min\{ j, \ell_0 \} }
    \indicator{ \{ 0 \} }( v )
    +
    A( \rell{ v }{ \theta }{ j }( x ) )
    \indicator{ \N }( v )
  \Bigr)
  ,
\end{split}
\end{equation}
let
$ \vartheta \in \R^n $,
$ x_1, x_2, \dots, x_m \in \R^{ \ell_0 } $,
$ y_1, y_2, \dots, y_m \in \R^{ \ell_L } $,
let $ F \colon \R^n \to \R $ satisfy for all $ \theta \in \R^n $ that
\begin{equation}
  F( \theta ) =
  \frac{ 1 }{ m }
  \sum_{ i = 1 }^m \mathscr{L}\bigl( \rel{ L }{ \theta }( x_i ) , y_i, \theta \bigr)
  ,
\end{equation}
and let
$
  P \colon \R^n \to \R^n
$
satisfy for all
$
  \theta \in \R^n \backslash \{ \vartheta \}
$
that
$ P( \theta ) = \norm{ \theta - \vartheta }^{ - 1 } ( \theta - \vartheta ) $.
Then
\begin{enumerate}[label=(\roman*)]
\item \llabel{item1} it holds that $ F \in C^1( \R^n, \R ) $ is definable,
\item \llabel{item2} it holds for all $ p \in \R^n $ that there exists $ \theta \in C^1( [0,\infty), \R^n ) $
such that for all $ t \in [0,\infty ) $ it holds that
$ \theta( 0 ) = p $
and
$
  \theta'(t) = - \bigl(\nabla F\bigr)\bigl( \theta(t) \bigr)
$,
\item \llabel{item3} it holds for all $ \theta \in C^1( [0,\infty), \R^n ) $
with $ \forall \, t \in [0,\infty) \colon \theta'(t) = - \bigl(\nabla F\bigr)\bigl( \theta(t) \bigr) $
that
$P\circ\theta$ has finite length, and
\item \llabel{item4} it holds for all $ \theta \in C^1( [0,\infty), \R^n ) $
with $ \forall \, t \in [0,\infty) \colon \theta'(t) = - \bigl(\nabla F\bigr)\bigl( \theta(t) \bigr) $
that there exists a unique $ v \in \R^n $ such that
$
  \limsup_{ t \to \infty } | ( P \circ \theta )( t ) - v | = 0
$.
\end{enumerate}
\end{athm}

\begin{aproof}
  Throughout this proof let
  $Q_v=(Q_v^1,\dots,Q_v^{\ell_v})
  \colon\R^n\times\R^{\ell_0}\to\R^{\ell_v}$,
  $v\in\{1,\dots,L\}$, satisfy for all $v\in\{1,\dots,L\}$,
  $w\in\R^n$, $z\in\R^{\ell_0}$ that
  \begin{eqn}
    \llabel{setting1}
    Q_v(w,z)=\rel{v}{w}(z)
    \dott
  \end{eqn}
  \cfload[.]\startnewargseq
  \argument{
    \lref{setting1};
    the assumption that for all $w\in\R^n$ it holds that
    $F(w)=\frac1m\sum_{i=1}^m
      \mathscr L(\rel{L}{w}(x_i),y_i,w)$;
    the assumption that $\mathscr L$ takes values in $[0,\infty)$;
  }{
    that for all $w\in\R^n$ it holds that
    \begin{eqn}
      \llabel{eq1}
      F(w)=\frac1m\sum_{i=1}^m
        \mathscr L(Q_L(w,x_i),y_i,w)\geq0
      \dott
    \end{eqn}
  }
  \argument{
    \cref{eq:neural_network_def-full};
    \lref{setting1};
  }{
    that for all $w=(w_1,\dots,w_n)\in\R^n$,
    $z=(z_1,\dots,z_{\ell_0})\in\R^{\ell_0}$,
    $i\in\{1,\dots,\ell_1\}$ it holds that
    \begin{eqn}
      \llabel{eq2}
      Q_1^i(w,z)=w_{\ell_1\ell_0+i}
        +\sum_{j=1}^{\ell_0}w_{(i-1)\ell_0+j}z_j
      \dott
    \end{eqn}
  }
  \argument{
    \lref{eq2};
    \cref{eq:neural_network_def-full};
    \lref{setting1};
    the fact that $A\in C^1(\R,\R)$ is definable;
    induction;
  }{
    that for all $v\in\{1,\dots,L\}$ it holds that
    $Q_v\in C^1(\R^n\times\R^{\ell_0},\R^{\ell_v})$ is
    definable\dott
    \llabel{arg1}
  }
  \argument{
    \lref{arg1};
    \lref{eq1};
    the assumption that
    $\mathscr L\in C^1(\R^{\ell_L}\times\R^{\ell_L}\times\R^n,
      [0,\infty))$ is definable;
  }[verbs=ep]{
    \lref{item1}\dott
  }

  In the following we prove \lref{item2}\dott\startnewargseq
  To prove \lref{item2}, let $p\in\R^n$\dott\startnewargseq
  \argument{
    the fact that $F\in C^1(\R^n,\R)$;
    Peano's existence theorem;
  }{
    that there exist $T\in(0,\infty]$,
    $\xi\in C^1([0,T),\R^n)$ which satisfy that
    \begin{enumerate}[label=(\Roman*)]
      \item \llabel{itemI} it holds that $\xi(0)=p$,
      \item \llabel{itemII} for all $t\in[0,T)$ it holds that
      $\xi'(t)=-\bigl(\grad F\bigr)\bigl(\xi(t)\bigr)$, and
      \item \llabel{itemIII} there exist no $S\in(T,\infty]$,
      $\chi\in C^1([0,S),\R^n)$ such that
      $\chi|_{[0,T)}=\xi$ and that for all $u\in[0,S)$ it holds that
      $\chi'(u)=-\bigl(\grad F\bigr)\bigl(\chi(u)\bigr)$\dott
    \end{enumerate}
  }
  \startnewargseq
  \argument{
    \lref{itemII};
    the fact that $F\in C^1(\R^n,\R)$;
    the chain rule;
  }{
    that for all $t\in[0,T)$ it holds that
    \begin{eqn}
      \llabel{eq4}
      (F\circ\xi)'(t)
        =\langle\bigl(\grad F\bigr)\bigl(\xi(t)\bigr),\xi'(t)\rangle
        =-\norm{\xi'(t)}^2
      \dott
    \end{eqn}
  }
  \argument{
    \lref{eq4};
    \lref{eq1};
    \lref{itemI};
    the fundamental theorem of calculus;
  }{
    that for all $s,t\in[0,T)$ with $s\leq t$ it holds that
    \begin{eqn}
      \llabel{eq5}
      0\leq\int_s^t\norm{\xi'(u)}^2\dd u
        \leq\int_0^t\norm{\xi'(u)}^2\dd u
        =F(p)-F(\xi(t))\leq F(p)
      \dott
    \end{eqn}
  }
  \argument{
    \lref{eq5};
    the fact that $\xi\in C^1([0,T),\R^n)$;
    the fundamental theorem of calculus;
    the Cauchy--Schwarz inequality;
  }{
    that for all $s,t\in[0,T)$ with $s\leq t$ it holds that
    \begin{eqn}[a]
      \llabel{eq6}
      \norm{\xi(t)-\xi(s)}
        \leq\int_s^t\norm{\xi'(u)}\dd u
        \leq\sqrt{t-s}
          \biggl(\int_s^t\norm{\xi'(u)}^2\dd u\biggr)^{1/2}
      \leq\sqrt{F(p)}\sqrt{t-s}
      \dott
    \end{eqn}
  }
  We claim that $T=\infty$\dott\startnewargseq
  We prove this claim by contradiction\dott\startnewargseq
  In the following we thus assume that $T<\infty$\dott
  \startnewargseq
  \argument{
    \lref{eq6};
    the assumption that $T<\infty$;
  }{
    that there exists $q\in\R^n$ which satisfies
    \begin{eqn}
      \llabel{eq7}
      \lim_{t\nearrow T}\xi(t)=q
      \dott
    \end{eqn}
  }
  \startnewargseq
  \argument{
    the fact that $q\in\R^n$;
    the fact that $F\in C^1(\R^n,\R)$;
    Peano's existence theorem;
  }{
    that there exist $\varepsilon\in(0,\infty)$,
    $\eta\in C^1([T,T+\varepsilon),\R^n)$ which satisfy for all
    $t\in[T,T+\varepsilon)$ that
    \begin{eqn}
      \llabel{eq8}
      \eta(T)=q
      \qquad\text{and}\qquad
      \eta'(t)=-\bigl(\grad F\bigr)\bigl(\eta(t)\bigr)
      \dott
    \end{eqn}
  }
  \argument{
    \lref{eq8};
    \lref{itemII};
    \lref{eq7};
    the fact that $\grad F$ is continuous;
  }{
    that
    \begin{eqn}
      \llabel{eq9}
      \lim_{t\nearrow T}\xi'(t)
        =-\bigl(\grad F\bigr)\bigl(q\bigr)=\eta'(T)
      \dott
    \end{eqn}
  }
  In the following let $\zeta\colon[0,T+\varepsilon)\to\R^n$
  satisfy for all $t\in[0,T+\varepsilon)$ that
  \begin{eqn}
    \llabel{setting2}
    \zeta(t)=
    \begin{cases}
      \xi(t) & \colon t\in[0,T),\\
      \eta(t) & \colon t\in[T,T+\varepsilon)
    \end{cases}
    \dott
  \end{eqn}
  \startnewargseq
  \argument{
    \lref{itemII};
    \lref{eq7};
    \lref{eq8};
    \lref{eq9};
    \lref{setting2};
    the fact that $\xi\in C^1([0,T),\R^n)$ and
    $\eta\in C^1([T,T+\varepsilon),\R^n)$;
  }{
    that $\zeta\in C^1([0,T+\varepsilon),\R^n)$ and that for all
    $t\in[0,T+\varepsilon)$ it holds that
    \begin{eqn}
      \llabel{eq10}
      \zeta|_{[0,T)}=\xi
      \qquad\text{and}\qquad
      \zeta'(t)=-\bigl(\grad F\bigr)\bigl(\zeta(t)\bigr)
      \dott
    \end{eqn}
  }
  \argument{
    \lref{eq10};
    \lref{itemIII};
    the fact that $\zeta\in C^1([0,T+\varepsilon),\R^n)$;}{that $T\geq T+\varepsilon>T$\dott}
   This contradiction proves that $T=\infty$\dott
    \startnewargseq
  \argument{
    the fact that $T=\infty$;
    \lref{itemI};
    \lref{itemII};
    the fact that $p\in\R^n$ was arbitrary;
  }[verbs=ep]{
    \lref{item2}\dott
  }

  In the following we prove \lref{item3} and \lref{item4}\dott
  \startnewargseq
  To prove \lref{item3} and \lref{item4}, let
  $\theta\in C^1([0,\infty),\R^n)$ and assume for all
  $t\in[0,\infty)$ that
  \begin{eqn}
    \llabel{setting3}
    \theta'(t)=-\bigl(\grad F\bigr)\bigl(\theta(t)\bigr)
    \dott
  \end{eqn}
  \startnewargseq
  \argument{
    \lref{setting3};
    the fact that $F\in C^1(\R^n,\R)$ is definable;
    the assumption that for all $w\in\R^n\backslash\{\vartheta\}$
    it holds that $P(w)=\norm{w-\vartheta}^{-1}(w-\vartheta)$;
    \cref{maintheorem3}
    (applied with $U\curvearrowleft\R^n$,
    $I\curvearrowleft[0,\infty)$ in the notation of
    \cref{maintheorem3});
  }{
    that
    \begin{eqn}
      \llabel{eq11}
      \curvelength{P\circ\theta}<\infty
      \dott
    \end{eqn}
  }
  \argument{
    \lref{eq11};
    the fact that $\theta$ was arbitrary;
  }[verbs=ep]{
    \lref{item3}\dott
  }
  \startnewargseq
  \argument{
    \lref{eq11};
    \cref{lem:finite-length-endpoint-limit}
    (applied with $U\curvearrowleft\R^n$,
    $I\curvearrowleft[0,\infty)$,
    $\gamma\curvearrowleft P\circ\theta$ in the notation of
    \cref{lem:finite-length-endpoint-limit});
  }{
    that there exists a unique $v\in\R^n$ such that
    \begin{eqn}
      \llabel{eq12}
      \limsup_{t\to\infty}\norm{(P\circ\theta)(t)-v}=0
      \dott
    \end{eqn}
  }
  \argument{
    \lref{eq12};
    the fact that $\theta$ was arbitrary;
  }[verbs=ep]{
    \lref{item4}\dott
  }
\end{aproof}

\section{Finite spherical length for ambient metrics}
\label{sec:ambient-metrics}

In this section we prove \cref{thm:ambient-conditional-reference} and
\cref{thm:ambient-uniform-general-interval} in the introduction and thereby
extend and generalize the conclusions of \cref{maintheorem3}
(\cref{maintheorem} in the introduction) and \cref{cor:spherical}
(\cref{maincor} in the introduction) from the standard Euclidean metric
to a large class of ambient metrics.
\Cref{thm:ambient-conditional-reference} is proved as
\cref{maintheorem:ambient1} in \cref{subsec:ambient-finite-length} below and
\cref{thm:ambient-uniform-general-interval} is proved as
\cref{maintheorem:ambient2} in \cref{subsec:ambient-uniform-length} below.

\subsection{A criterion for absolute continuity}

In this subsection we provide in \cref{lem:AC-gluing} below for every
$a\in\R$, $b\in[a,\infty)$ sufficient conditions to ensure that a continuous
function $u\colon[a,b]\to\R$ is absolutely continuous.
For our proof of \cref{lem:AC-gluing} we first recall the concept of a
(signed) Radon measure (cf.\ \cref{def:radon-measure,def:signed-radon-measure}
below) and we recall some basic properties for (signed) Radon measures;
cf.\ \cref{lem:radon-decomposition-uniqueness} below.
Only for completeness we include here a proof for
\cref{lem:radon-decomposition-uniqueness}.

\begin{athm}{definition}{def:radon-measure}[Radon measure]
  Let $J\subseteq\R$ be an open interval\dott \cfload[.]
  Then we say that $\lambda$ is a Radon measure on $J$ if and only if we have that
  \begin{enumerate}[label=(\roman*)]
    \item \llabel{item:measure} it holds that $\lambda\colon\mathcal B(J)\to[0,\infty]$
    is a measure on $J$,
    \item \llabel{item1} it holds for all compact $A\subseteq J$ that
    $\lambda(A)<\infty$, and
    \item \llabel{item2} it holds for all $A\in\mathcal B(J)$ that
    $\lambda(A)=\sup_{B\subseteq A\text{ compact}}\lambda(B)
      =\inf_{B\subseteq J\text{ open},\,B\supseteq A}\lambda(B)$\dott
  \end{enumerate}
  \cfout[.]
\end{athm}

\begin{athm}{lemma}{lem:radon-decomposition-uniqueness}[Uniqueness of the Jordan decomposition]
  Let $J\subseteq\R$ be an open interval, let $P_0,P_1\in\mathcal B(J)$,
  for every $i,j\in\{0,1\}$ let $\mu_{i,j}$ be a Radon measure on
  $J$\cfadd{def:radon-measure}, let $\mathcal D\subseteq\mathcal B(J)$ satisfy
  for all $i\in\{0,1\}$ that
  $\mathcal D=\{A\in\mathcal B(J)\colon\textstyle\sum_{j=0}^{1}\mu_{i,j}(A)<\infty\}$,
  let $\lambda\colon\mathcal D\to\R$ satisfy for all $i\in\{0,1\}$, $A\in\mathcal D$ that
  \begin{eqn}
    \llabel{ass:decomposition}
    \lambda(A)=\mu_{i,0}(A)-\mu_{i,1}(A),
  \end{eqn}
  and assume $\sum_{i=0}^{1}(\mu_{i,0}(J\backslash P_i)+\mu_{i,1}(P_i))=0$\dott
  \cfload[.]
  Then
  \begin{enumerate}[label=(\roman*)]
    \item \llabel{item:uniqueness} it holds for all $j\in\{0,1\}$ that
    $\mu_{0,j}=\mu_{1,j}$ and
    \item \llabel{item:separating-set} it holds for all $i,j\in\{0,1\}$ that
    $\mu_{i,j}((P_0\backslash P_1)\cup(P_1\backslash P_0))=0$\dott
  \end{enumerate}
  \cfout[.]
\end{athm}
\begin{aproof}
  \argument{
    the assumption that
    $\mathcal D=\{A\in\mathcal B(J)\colon\sum_{j=0}^{1}\mu_{0,j}(A)<\infty\}$;
    the fact that $\mu_{0,0}$ and $\mu_{0,1}$ are Radon measures on $J$;
  }{
    that for all compact $K\subseteq J$ and all $B\in\mathcal B(J)$
    with $B\subseteq K$ it holds that
    \begin{eqn}
      \llabel{eq:finite-subsets}
      \mu_{0,0}(B)+\mu_{0,1}(B)\leq\mu_{0,0}(K)+\mu_{0,1}(K)<\infty\qquad\text{and}\qquad
     B\in\mathcal D
      \dott
    \end{eqn}
  }
  \argument{
    \lref{eq:finite-subsets};
    \lref{ass:decomposition};
    the assumption that $\sum_{i=0}^{1}(\mu_{i,0}(J\backslash P_i)+\mu_{i,1}(P_i))=0$;
  }{
    that for all compact $K\subseteq J$ it holds that
    \begin{gather}
      \llabel{eq:compact-parts}
      \mu_{0,0}(K)=\lambda(K\cap P_0)\leq\mu_{1,0}(K)
        =\lambda(K\cap P_1)\leq\mu_{0,0}(K)\\
      \text{and}\qquad
      \mu_{0,1}(K)=-\lambda(K\backslash P_0)\leq\mu_{1,1}(K)
        =-\lambda(K\backslash P_1)\leq\mu_{0,1}(K)
      \dott
    \end{gather}
  }
  \argument{
    \lref{eq:compact-parts};
    the fact that $\mu_{0,0}$, $\mu_{0,1}$, $\mu_{1,0}$, and $\mu_{1,1}$ are Radon measures on $J$;
  }{
    that for all $A\in\mathcal B(J)$ it holds that
    \begin{gather}
      \llabel{eq:parts}
      \mu_{0,0}(A)=\sup_{K\subseteq A\text{ compact}}\mu_{0,0}(K)
        =\sup_{K\subseteq A\text{ compact}}\mu_{1,0}(K)=\mu_{1,0}(A)\\
      \text{and}\qquad
      \mu_{0,1}(A)=\sup_{K\subseteq A\text{ compact}}\mu_{0,1}(K)
        =\sup_{K\subseteq A\text{ compact}}\mu_{1,1}(K)=\mu_{1,1}(A)
      \dott
    \end{gather}
  }
  \argument{
    \lref{eq:parts};
  }[verbs=ep]{
    \lref{item:uniqueness}\dott
  }
  \startnewargseq
  \argument{
    \lref{item:uniqueness};
    the assumption that $\sum_{i=0}^{1}(\mu_{i,0}(J\backslash P_i)+\mu_{i,1}(P_i))=0$;
  }{
    that
    \begin{gather}
      \llabel{eq:null-differences}
      0\leq\mu_{0,0}(P_0\backslash P_1)+\mu_{0,1}(P_0\backslash P_1)
        \leq\mu_{1,0}(J\backslash P_1)+\mu_{0,1}(P_0)=0\\
      \text{and}\qquad
      0\leq\mu_{0,0}(P_1\backslash P_0)+\mu_{0,1}(P_1\backslash P_0)
        \leq\mu_{0,0}(J\backslash P_0)+\mu_{1,1}(P_1)=0
      \dott
    \end{gather}
  }
  \argument{
    \lref{eq:null-differences};
    \lref{item:uniqueness};
  }[verbs=ep]{
    \lref{item:separating-set}\dott
  }
\end{aproof}

\begin{athm}{definition}{def:signed-radon-measure}
  [Signed Radon measure]
  Let $J\subseteq\R$ be an open interval\dott \cfload[.]
  Then
  \begin{enumerate}[label=(\roman*)]
    \item \llabel{item:criterion} \llabel{item1} we say that $\lambda$ is a signed Radon measure
    on $J$ if and only if there exist $P\in\mathcal B(J)$ and Radon measures
    $\mu_0$ and $\mu_1$ on $J$\cfadd{def:radon-measure} such that
    \begin{enumerate}[label=(\alph*)]
      \item \llabel{item:domain} it holds that
      $\lambda\colon\midcap_{i=0}^{1}((\mu_i)^{-1}([0,\infty)))\to\R$ is a function from
      $\midcap_{i=0}^{1}((\mu_i)^{-1}([0,\infty)))$ to $\R$,
      \item \llabel{item:singularity} it holds that
      $\mu_0(J\backslash P)=\mu_1(P)=0$,
      and
      \item \llabel{item:decomposition} it holds for all $A\in\lambda^{-1}(\R)$ that
      $\lambda(A)=\mu_0(A)-\mu_1(A)$,
    \end{enumerate}
    \cfconsiderloaded{def:signed-radon-measure}
    \item \llabel{item:total-variation} for every signed Radon measure $\lambda$ on $J$
    we denote by $\TVmeasure{\lambda}\colon\mathcal B(J)\to[0,\infty]$
    the function which satisfies for all $P\in\mathcal B(J)$ and all Radon measures
    $\mu_0$ and $\mu_1$ on $J$ with
    $\lambda^{-1}(\R)=\midcap_{i=0}^{1}((\mu_i)^{-1}([0,\infty)))$,
    $\mu_0(J\backslash P)=\mu_1(P)=0$, and
    $\forall\, A\in\lambda^{-1}(\R)\colon\lambda(A)=\mu_0(A)-\mu_1(A)$
    that $\TVmeasure{\lambda}=\mu_0+\mu_1$,
    \item \llabel{item2} for every signed Radon measure $\lambda$ on $J$
    we say that $\lambda$ is finite if and only if
    $\TVmeasure{\lambda}(J)<\infty$,
    and
    \item \llabel{item3} for every signed Radon measure $\lambda$ on $J$
    and every measurable
    $g\colon J\to\R$ with $\int_J\lvert g(t)\rvert\,\TVmeasure{\lambda}(\ddd t)<\infty$
    we denote by $\signedint_J g(t)\,\lambda(\ddd t)\in\R$
    the real number which satisfies for all $P\in\mathcal B(J)$ and all Radon measures
    $\mu_0$ and $\mu_1$ on $J$ with
    $\lambda^{-1}(\R)=\midcap_{i=0}^{1}((\mu_i)^{-1}([0,\infty)))$,
    $\mu_0(J\backslash P)=\mu_1(P)=0$, and
    $\forall\, A\in\lambda^{-1}(\R)\colon\lambda(A)=\mu_0(A)-\mu_1(A)$
    that
    \begin{eqn}
      \llabel{eq2}
      \signedint_J g(t)\,\lambda(\ddd t)
        =\int_J g(t)\,\mu_0(\ddd t)-\int_J g(t)\,\mu_1(\ddd t)
      \dott
    \end{eqn}
  \end{enumerate}
  \cfout[.]
\end{athm}

\begin{athm}{lemma}{lem:BV-locality}
  Let $a\in\R$, $b\in[a,\infty)$,
  let $u\in C([a,b],\R)$ have finite length\cfadd{definition:curve-length},
  let $E\subseteq [a,b]$ be closed, assume $u|_E=0$, and let
  $U\colon\R\to\R$ satisfy for all $t\in\R$ that
  \begin{equation}
    \llabel{ass1}
    U(t)=
    \begin{cases}
      u(a) & \colon t<a\\
      u(t) & \colon t\in [a,b]\\
      u(b) & \colon t>b.
    \end{cases}
  \end{equation}
  \cfload[.]
  Then there exists a finite signed Radon measure $\mu$ on
  $\R$\cfadd{def:signed-radon-measure} such that
  \begin{enumerate}[label=(\roman*)]
    \item \llabel{item:derivative-measure} it holds for all
    $\varphi\in C^\infty(\R,\R)$ with compact support that
    $\signedint_{\R}\varphi(t)\,\mu(\ddd t)
      =-\int_{\R}U(t)\varphi'(t)\dd t$,
    \item \llabel{item1} it holds for all $A\in\mathcal B(\R)$ that
    $\mu(A\cap E)=0$,
    \item \llabel{item2} it holds that $\TVmeasure{\mu}(E)=0$, and
    \item \llabel{item3} it holds for all $s,t\in[a,b]$ with $s\leq t$ that
    $\mu((s,t])=u(t)-u(s)$ and $\mu(\{s\})=0$\dott
  \end{enumerate}
  \cfout[.]
\end{athm}
\begin{aproof}
  Throughout this proof let $\pi\colon\R\to [a,b]$ satisfy for all $t\in\R$ that
  \begin{eqn}
    \llabel{setting0}
    \pi(t)=\min\{b,\max\{a,t\}\}
  \end{eqn}
  and
  for all $s,t\in E$ with $s\leq t$ let $\mathcal G_{s,t}\subseteq\mathcal B(\R)$ satisfy
  \begin{eqn}
    \llabel{setting1}
    \mathcal G_{s,t}
      =\{(p,q)\colon p,q\in E,\ s\leq p<q\leq t,
        \ (p,q)\cap E=\emptyset\}
    \dott
  \end{eqn}
  \cfload[.]
  \startnewargseq
  \argument{
    \lref{ass1};
    \lref{setting0};
    the assumption that $u$ is continuous;
  }{
    that $U=u\circ\pi$ is \llabel{arg1} continuous\dott
  }
  \argument{
    the assumption that $u$ has finite length\cfadd{definition:curve-length};
\lref{arg1};
    the fact that $\pi$ is non-decreasing;
  }{
    that
    \begin{eqn}[a]
      \llabel{eq0.1}
      &\sup_{\substack{n\in\N,\ (\tau_0,\tau_1\ldots,\tau_n)\in\R^{n+1}\\
        \tau_0\leq\tau_1\leq\ldots\leq\tau_n}}
        \sum_{j=1}^{n}\lvert U(\tau_j)-U(\tau_{j-1})\rvert
      \leq
      \sup_{\substack{n\in\N,\ (\tau_0,\tau_1\ldots,\tau_n)\in [a,b]^{n+1}\\
        a=\tau_0\leq \tau_1\leq\ldots\leq\tau_n=b}}
        \sum_{j=1}^{n}\lvert u(\tau_j)-u(\tau_{j-1})\rvert
      <\infty
      \dott
    \end{eqn}
  }
  \argument{
    \lref{eq0.1};
    \cite[Theorem~6.13]{Apostol74} (applied with $a\curvearrowleft a-1$,
    $b\curvearrowleft b+1$, $f\curvearrowleft U|_{[a-1,b+1]}$
    in the notation of \cite[Theorem~6.13]{Apostol74});
    the fact that for all $t\in(-\infty,a-1]$
    it holds that $U(t)=U(a-1)=u(a)$;
    the fact that for all $t\in[b+1,\infty)$
    it holds that $U(t)=U(b+1)=u(b)$;
  }{
    that there exist bounded and non-decreasing
    $P,Q\colon\R\to\R$ which satisfy that
    \begin{eqn}
      \llabel{eq0.1a}
      U=P-Q
      \dott
    \end{eqn}
  }
  In the following let $F,H\colon\R\to\R$ satisfy for all $t\in\R$
  that
  \begin{eqn}
    \llabel{setting2}
    F(t)=\inf_{s\in (t,\infty)}P(s)
    \qquad \text{and}\qquad H(t)=\inf_{s\in (t,\infty)}Q(s)
    \dott
  \end{eqn}
  \cfload[.]
  \startnewargseq
  \argument{
    \lref{setting2};
    the fact that $P$ and $Q$ are bounded and non-decreasing;
  }{
    that $F$ and $H$ are bounded and non-decreasing and that for all
    $t\in\R$ it holds that
    \begin{eqn}
      \llabel{eq0.1b}
      \limsup_{s\searrow t}\abs{P(s)-F(t)}=0
      \qquad \text{and}\qquad \limsup_{s\searrow t}\abs{Q(s)-H(t)}=0
      \dott
    \end{eqn}
  }
  \argument{
    \lref{eq0.1b};
    the fact that for all $t,x,s\in\R$
    with $t<x<s$ it holds that $F(t)\leq F(x)\leq P(s)$
    and $H(t)\leq H(x)\leq Q(s)$;
    the fact that $F$ and $H$ are bounded and non-decreasing
  }{
    that $F$ and $H$ are right-continuous\dott
  }
  \argument{
    \lref{eq0.1a};
    \lref{eq0.1b};
    the fact that $U$ is continuous;
  }{
    that for all $t\in\R$ it holds that
    \begin{eqn}[a]
      \llabel{eq0.1.1}
      F(t)-H(t)
        =\lim_{s\searrow t}\bigl(P(s)-Q(s)\bigr)
      =\lim_{s\searrow t}U(s)=U(t)
      \dott
    \end{eqn}
  }
  \argument{
    the fact that $F$ and $H$ are non-decreasing and right-continuous;
    \cite[Theorem~12.3]{Billingsley95};
    \cite[Theorem~12.4]{Billingsley95};
  }{
    that there exist Radon measures $\lambda,\eta$ on
    $\R$\cfadd{def:radon-measure} which satisfy for all $s,t\in\R$
    with $s<t$ that
    \begin{eqn}
      \llabel{setting3}
      \lambda((s,t])=F(t)-F(s)
      \qquad \text{and}\qquad \eta((s,t])=H(t)-H(s)
      \dott
    \end{eqn}
  }
  \startnewargseq
  \argument{
    \lref{setting3};
    the fact that $F$ and $H$ are bounded;
  }{
    that
    \begin{eqn}[a]
      \llabel{eq0.2}
      \lambda(\R)
        &=\lim_{n\to\infty}\bigl(F(n)-F(-n)\bigr)<\infty
      \qquad\text{and}\qquad
      \eta(\R)
        &=\lim_{n\to\infty}\bigl(H(n)-H(-n)\bigr)<\infty
      \dott
    \end{eqn}
  }
  \startnewargseq
  In the following let $\mu\colon\mathcal B(\R)\to\R$ satisfy
  for all $A\in\mathcal B(\R)$ that
  \begin{eqn}
    \llabel{setting4}
    \mu(A)=\lambda(A)-\eta(A)
    \dott
  \end{eqn}
  \cfload[.]
  \startnewargseq
  \argument{
    \lref{eq0.2};
    \lref{setting4};
    the fact that $\lambda$ and $\eta$ are countably additive;
  }{
    that $\mu$ is real-valued and countably additive and that $\mu(\emptyset)=0$
    \llabel{arg2}\dott
  }
  \argument{
    \lref{arg2};
    the Hahn decomposition theorem;
  }{
    that there exist Borel measurable $S,N\subseteq\R$ which satisfy
    \begin{eqn}
      \llabel{eq0.2a.1}
      S\cup N=\R
      \qquad\text{and}\qquad
      S\cap N=\emptyset
    \end{eqn}
    and which satisfy for all $B,C\in\mathcal B(\R)$ with
    $B\subseteq S$ and $C\subseteq N$ that
    \begin{eqn}
      \llabel{eq0.2a.2}
      \mu(B)\geq0
      \qquad\text{and}\qquad
      \mu(C)\leq0
      \dott
    \end{eqn}
  }
  \startnewargseq
  \argument{
    \lref{eq0.2a.1}; \lref{eq0.2a.2};
  }{
    that for all $A\in\mathcal B(\R)$ it holds that
    \begin{eqn}
      \llabel{eq0.2a}
      \mu(A\cap S)\geq0
      \qquad\text{and}\qquad
      \mu(A\backslash S)\leq0
      \dott
    \end{eqn}
  }
  In the following let $\boldsymbol{\mu}_1,\boldsymbol{\mu}_2\colon\mathcal B(\R)\to[0,\infty)$
  satisfy for all $A\in\mathcal B(\R)$ that
  \begin{eqn}
    \llabel{setting4.1}
    \boldsymbol{\mu}_1(A)=\mu(A\cap S)
    \qquad\text{and}\qquad
    \boldsymbol{\mu}_2(A)=-\mu(A\backslash S)
    \dott
  \end{eqn}
  \cfload[.]
  \startnewargseq
  \argument{
    \lref{eq0.2};
    \lref{setting4};
    \lref{eq0.2a};
    \lref{setting4.1};
    the fact that $\mu$ is countably additive;
  }{
    that $\boldsymbol{\mu}_1$ and $\boldsymbol{\mu}_2$ are finite nonnegative Borel measures
    and that for all $A\in\mathcal B(\R)$ it holds that
    \begin{eqn}[a]
      \llabel{eq0.2b}
      0\leq\boldsymbol{\mu}_1(A)
        =\lambda(A\cap S)-\eta(A\cap S)\leq\lambda(A)\quad\text{and}\quad
      0\leq\boldsymbol{\mu}_2(A)
        =\eta(A\backslash S)-\lambda(A\backslash S)\leq\eta(A)
      \dott
    \end{eqn}
  }
  \argument{
    \lref{setting4.1};
    the fact that $\mu$ is additive;
    the fact that $\mu(\emptyset)=0$;
  }{
    that
    \begin{eqn}
      \llabel{eq0.2c}
      \mu=\boldsymbol{\mu}_1-\boldsymbol{\mu}_2
      \qquad\text{and}\qquad
      \boldsymbol{\mu}_1(\R\backslash S)=\boldsymbol{\mu}_2(S)=0
      \dott
    \end{eqn}
  }
  \argument{
    \lref{eq0.2};
    \lref{eq0.2b};
    the fact that $\lambda$ and $\eta$ are nonnegative measures;
  }{
    that for every compact $C\subseteq\R$ it holds that
    \begin{eqn}
      \llabel{eq0.2c.1}
      \boldsymbol{\mu}_1(C)\leq\lambda(C)\leq\lambda(\R)<\infty
      \qquad \text{and}\qquad
      \boldsymbol{\mu}_2(C)\leq\eta(C)\leq\eta(\R)<\infty
      \dott
    \end{eqn}
  }
  \argument{\itref{def:radon-measure}{item2};
    \lref{eq0.2};
    the fact that $\lambda$ and $\eta$ are Radon measures on $\R$;
  }{
    that for all $A\in\mathcal B(\R)$,
    $\varepsilon\in(0,\infty)$ there exist compact
    $K_1,K_2\subseteq A$ such that
    \begin{eqn}
      \llabel{eq0.2c.2}
      0\leq\lambda(A)-\lambda(K_1)<\varepsilon
      \qquad \text{and}\qquad
      0\leq\eta(A)-\eta(K_2)<\varepsilon
      \dott
    \end{eqn}
  }
  \argument{
    \lref{eq0.2c.2};
    \lref{eq0.2b};
    the fact that $\boldsymbol{\mu}_1$, $\boldsymbol{\mu}_2$,
    $\lambda$, and $\eta$ are finite nonnegative measures;
  }{
    that for all $A\in\mathcal B(\R)$,
    $\varepsilon\in(0,\infty)$ there exist compact
    $K_1,K_2\subseteq A$ such that
    \begin{eqn}[a]
      \llabel{eq0.2c.3}
      &0\leq\boldsymbol{\mu}_1(A)-\boldsymbol{\mu}_1(K_1)
        =\boldsymbol{\mu}_1(A\backslash K_1)
        \leq\lambda(A\backslash K_1)
      =\lambda(A)-\lambda(K_1)<\varepsilon
      \\
     \text{and}\qquad &0\leq\boldsymbol{\mu}_2(A)-\boldsymbol{\mu}_2(K_2)
        =\boldsymbol{\mu}_2(A\backslash K_2)
        \leq\eta(A\backslash K_2)
      =\eta(A)-\eta(K_2)<\varepsilon
      \dott
    \end{eqn}
  }
  \argument{
    \lref{eq0.2c.3};
    the fact that for every $A\in\mathcal B(\R)$, $j\in\{1,2\}$ and every compact
    $C\subseteq A$ it holds that
    $\boldsymbol{\mu}_j(C)\leq\boldsymbol{\mu}_j(A)$;
  }{
    that for all $A\in\mathcal B(\R)$, $j\in\{1,2\}$
    it holds that
    \begin{eqn}
      \llabel{eq0.2c.4}
      \boldsymbol{\mu}_j(A)
        =\sup\{\boldsymbol{\mu}_j(C)\colon
          C\subseteq A\text{ is compact}\}
      \dott
    \end{eqn}
  }
  \argument{\itref{def:radon-measure}{item2};
    \lref{eq0.2};
    the fact that $\lambda$ and $\eta$ are Radon measures on $\R$;
  }{
    that for all $A\in\mathcal B(\R)$,
    $\varepsilon\in(0,\infty)$ there exist open
    $V_1,V_2\subseteq\R$ with $A\subseteq V_1\cap V_2$ such that
    \begin{eqn}
      \llabel{eq0.2c.5}
      0\leq\lambda(V_1)-\lambda(A)<\varepsilon
      \qquad \text{and}\qquad
      0\leq\eta(V_2)-\eta(A)<\varepsilon
      \dott
    \end{eqn}
  }
  \argument{
    \lref{eq0.2c.5};
    \lref{eq0.2b};
    the fact that $\boldsymbol{\mu}_1$, $\boldsymbol{\mu}_2$,
    $\lambda$, and $\eta$ are finite nonnegative measures
  }{
    that for all $A\in\mathcal B(\R)$,
    $\varepsilon\in(0,\infty)$ there exist open
    $V_1,V_2\subseteq\R$ with $A\subseteq V_1\cap V_2$ such that
    \begin{gather}
      \llabel{eq0.2c.6}
      0\leq\boldsymbol{\mu}_1(V_1)-\boldsymbol{\mu}_1(A)
        =\boldsymbol{\mu}_1(V_1\backslash A)
        \leq\lambda(V_1\backslash A)
      =\lambda(V_1)-\lambda(A)<\varepsilon
      \\
     \text{and}\qquad 0\leq\boldsymbol{\mu}_2(V_2)-\boldsymbol{\mu}_2(A)
        =\boldsymbol{\mu}_2(V_2\backslash A)
        \leq\eta(V_2\backslash A)
      =\eta(V_2)-\eta(A)<\varepsilon
      \dott
    \end{gather}
  }
  \argument{
    \lref{eq0.2c.6};
    the fact that for every $A\in\mathcal B(\R)$, $j\in\{1,2\}$ and every open
    $W\subseteq\R$ with $A\subseteq W$
    it holds that $\boldsymbol{\mu}_j(A)\leq\boldsymbol{\mu}_j(W)$;
  }{
    that for all $A\in\mathcal B(\R)$, $j\in\{1,2\}$
    it holds that
    \begin{eqn}
      \llabel{eq0.2c.7}
      \boldsymbol{\mu}_j(A)
        =\inf\{\boldsymbol{\mu}_j(W)\colon
          A\subseteq W\subseteq\R,\ W\text{ is open}\}
      \dott
    \end{eqn}
  }
  \argument{
    the fact that $\boldsymbol{\mu}_1$ and $\boldsymbol{\mu}_2$
    are finite nonnegative Borel measures;
  }{
    that
    \begin{eqn}
      \llabel{eq0.2c.8}
      \{A\in\mathcal B(\R)\colon
        \boldsymbol{\mu}_1(A)+\boldsymbol{\mu}_2(A)<\infty\}
        =\mathcal B(\R)
      \dott
    \end{eqn}
  }
  \argument{
    \lref{eq0.2c.1};
    \lref{eq0.2c.4};
    \lref{eq0.2c.7};
    the fact that $\boldsymbol{\mu}_1$ and $\boldsymbol{\mu}_2$
    are finite nonnegative Borel measures;
  }{
    that $\boldsymbol{\mu}_1$ and $\boldsymbol{\mu}_2$ are finite
    Radon measures on $\R$\cfadd{def:radon-measure}\llabel{arg3}\dott
  }
  \argument{
    \lref{arg3};
    \lref{eq0.2c};
    \lref{eq0.2c.8};
    the fact that $S$ is Borel measurable;
  }{
    that $\mu$ is a signed Radon measure on $\R$\cfadd{def:signed-radon-measure}
    and that
    \begin{eqn}
      \llabel{eq0.2c.9}
      \TVmeasure{\mu}=\boldsymbol{\mu}_1+\boldsymbol{\mu}_2
      \dott
    \end{eqn}
  }
  \argument{
    \lref{eq0.2c.9};
    \lref{eq0.2};
    \lref{eq0.2b};
  }{
    that $\mu$ is a finite signed Radon measure on $\R$ and that
    \begin{eqn}
      \llabel{eq0.3}
      \TVmeasure{\mu}(\R)=\boldsymbol{\mu}_1(\R)+\boldsymbol{\mu}_2(\R)
        \leq\lambda(\R)+\eta(\R)<\infty
      \dott
    \end{eqn}
  }
  In the following let $\nu\colon\mathcal B(\R)\to\R$
  satisfy for all $A\in\mathcal B(\R)$ that
  \begin{eqn}
    \llabel{setting4.2}
    \nu(A)=\mu(A\cap E)
    \dott
  \end{eqn}
  \cfload[.]
  \startnewargseq
  \argument{
    \lref{eq0.1.1};
    \lref{setting3};
    \lref{setting4};
  }{
    that for all $s,t\in\R$ with $s<t$ it holds that
    \begin{eqn}[a]
      \llabel{eq0.4}
      \mu((s,t])
        =F(t)-F(s)-H(t)+H(s)
      =U(t)-U(s)
      \dott
    \end{eqn}
  }
  \argument{
    \lref{eq0.4};
    \lref{eq0.1};
    the integration by parts formula in \cite[Proposition~2.1]{DatchevShapiro23};
    the fact that for all compactly supported $\varphi\in C^\infty(\R,\R)$
    it holds that $\varphi$ is absolutely continuous;
  }{
    that for all compactly supported $\varphi\in C^\infty(\R,\R)$
    it holds that
    \begin{eqn}
      \llabel{eq0.5}
      \signedint_{\R}\varphi(t)\,\mu(\ddd t)
        =-\int_{\R}U(t)\varphi'(t)\dd t
      \dott
    \end{eqn}
  }
  \argument{
    \lref{eq0.3};
    \lref{eq0.5};
  }[verbs=ep]{
    \lref{item:derivative-measure}\dott
  }
  \argument{
    the fact that $\boldsymbol{\mu}_1$ and $\boldsymbol{\mu}_2$
    are finite nonnegative countably additive measures;
  }{
    that for every non-increasing
    $(A_n)_{n\in\N}\subseteq\mathcal B(\R)$, every $m\in\N$,
    and every $j\in\{1,2\}$ it holds that
    \begin{eqn}
      \llabel{eq0.6.1}
     \textstyle 0\leq\boldsymbol{\mu}_j(A_m)
        -\boldsymbol{\mu}_j\bigl(\bigcap_{n\in\N}A_n\bigr)
        =\sum_{k=m}^{\infty}\boldsymbol{\mu}_j(A_k\backslash A_{k+1})
        \leq\boldsymbol{\mu}_j(\R)<\infty
      \dott
    \end{eqn}
  }
  \argument{
    \lref{eq0.6.1};
  }{
    that for every non-increasing
    $(A_n)_{n\in\N}\subseteq\mathcal B(\R)$ and every
    $j\in\{1,2\}$ it holds that
    \begin{eqn}
      \llabel{eq0.6.2}
     \textstyle \limsup_{m\to\infty}\bigl|\boldsymbol{\mu}_j(A_m)
       -\boldsymbol{\mu}_j\bigl(\bigcap_{n\in\N}A_n\bigr)\bigr|=0
      \dott
    \end{eqn}
  }
  \argument{
    \lref{eq0.6.2};
    \lref{eq0.2c};
    the fact that $\boldsymbol{\mu}_1$ and $\boldsymbol{\mu}_2$
    are finite measures;
  }{
    that for every non-increasing 
    $(A_n)_{n\in\N}\subseteq\mathcal B(\R)$ it holds that
    \begin{eqn}[a]
      \llabel{eq0.6.3}
      \lim_{m\to\infty}\mu(A_m)
        &=\lim_{m\to\infty}
          \bigl(\boldsymbol{\mu}_1(A_m)-\boldsymbol{\mu}_2(A_m)\bigr)
     =\textstyle\boldsymbol{\mu}_1\bigl(\bigcap_{n\in\N}A_n\bigr)
          -\boldsymbol{\mu}_2\Bigl(\bigcap_{n\in\N}A_n\Bigr)
        =\mu\bigl(\bigcap_{n\in\N}A_n\bigr)
      \dott
    \end{eqn}
  }
  \argument{
    \lref{eq0.6.3};
    \lref{eq0.4};
    the fact that $U$ is continuous;
  }{
    that for all $x\in\R$ it holds that
    \begin{eqn}[a]
      \llabel{eq0.7}
      \mu(\{x\})
        =\lim_{n\to\infty}\mu((x-1/n,x])
        =\lim_{n\to\infty}\bigl(U(x)-U(x-1/n)\bigr)=0
      \dott
    \end{eqn}
  }
  \argument{
    \lref{eq0.7};
    \lref{ass1};
    \lref{eq0.4};
  }{
    that
    for all $s,t\in [a,b]$ with $s\leq t$ it holds that
    \begin{eqn}
      \llabel{eq1}
      \mu([s,t])=\mu((s,t))=u(t)-u(s)
      \dott
    \end{eqn}
  }
  \argument{
    the assumption that $E$ is closed in $[a,b]$;
  }{
    that \llabel{argt1} $E$ is compact and Borel measurable\dott
  }
  \argument{\lref{argt1};
    \lref{eq0.2c};
    \lref{setting4.2};
    the fact that $\boldsymbol{\mu}_1$ and $\boldsymbol{\mu}_2$
    are finite nonnegative Borel measures;
  }{
    that $\nu$ is real-valued and countably additive and that
    for all $A\in\mathcal B(\R)$ it holds that
    \begin{eqn}
      \llabel{eq1.1}
      \nu(A)=\boldsymbol{\mu}_1(A\cap E)-\boldsymbol{\mu}_2(A\cap E)
      \dott
    \end{eqn}
  }
  \argument{
    the fact that $E$ is closed;
    the fact that every open subset of $\R$ is the union of at most
    countably many pairwise disjoint open intervals;
    \lref{setting1};
  }{
    that for all $s,t\in E$ with $s\leq t$ it holds that
    $\mathcal G_{s,t}$ is at most countable, that its elements are
    pairwise disjoint bounded open intervals with endpoints in $E$,
    and that
    \begin{equation}
      \llabel{eq2}
     \textstyle (s,t)\backslash E=\bigcup_{G\in\mathcal G_{s,t}}G.
    \end{equation}
  }
  \argument{\lref{eq2};
    \lref{eq1};
    the assumption that $u|_E=0$;
  }{
    that for all $s,t\in E$ with $s\leq t$ it holds that
    $\mu([s,t])=0$ and that for all $G\in\mathcal G_{s,t}$ it holds that
    \begin{eqn}
      \llabel{eq3}
      \mu(G)=0
      \dott
    \end{eqn}
  }
  \argument{
    \lref{eq3};
    \lref{eq2};
  }{
    that for all $s,t\in E$ with $s\leq t$ it holds that
    \begin{eqn}
      \llabel{eq4}
      \mu(E\cap[s,t])
        =\mu([s,t])-\sum_{G\in\mathcal G_{s,t}}\mu(G)=0
      \dott
    \end{eqn}
  }
  \argument{
    the assumption that $E$ is compact;
  }{
    that for all $r\in\R$ with $E\cap(-\infty,r]\neq\emptyset$
    there exist $c,d\in E$ with $c\leq d$ such that
    \begin{eqn}
      \llabel{eq5}
      c=\min\{x\colon x\in E\},
      \qquad d=\max\{x \colon x\in E\cap(-\infty,r]\},
      \qquad\text{and}\qquad E\cap(-\infty,r]=E\cap[c,d]
      \dott
    \end{eqn}
  }
  \argument{  \lref{setting4.2};
    \lref{eq5};
    \lref{eq4};
    the fact that $\mu(\emptyset)=0$;
  }{
    that for all $r\in\R$ it holds that
    \begin{eqn}
      \llabel{eq6}
      \nu((-\infty,r])=\mu(E\cap(-\infty,r])=0
      \dott
    \end{eqn}
  }
  \argument{
    \lref{eq6};
    \lref{eq1.1};
    the fact that for all $j\in\{1,2\}$ it holds that
    $\mathcal B(\R)\ni A\mapsto\boldsymbol{\mu}_j(A\cap E)$
    is a finite nonnegative Borel measure;
  }{
    that for all $A\in\mathcal B(\R)$ it holds that
    \begin{eqn}
      \llabel{eq6.1}
      \boldsymbol{\mu}_1(A\cap E)=\boldsymbol{\mu}_2(A\cap E)
      \dott
    \end{eqn}
  }
  \argument{
    \lref{eq6.1};
    \lref{eq1.1};
  }{
    that
    \begin{eqn}
      \llabel{eq7}
      \nu=0
      \dott
    \end{eqn}
  }
  \argument{
    \lref{eq7};
    \lref{setting4.2};
  }[verbs=ep]{
    \lref{item1}\dott
  }
  \startnewargseq
  \argument{
    \lref{setting4.1}, \lref{setting4.2}; \lref{eq7};
    the fact that $S$ and $E$ are Borel measurable;
  }{
    that
    \begin{eqn}
      \llabel{eq7.1}
      \boldsymbol{\mu}_1(E)=\mu(E\cap S)=\nu(S)=0
      \qquad \text{and}\qquad
      \boldsymbol{\mu}_2(E)=-\mu(E\backslash S)=-\nu(\R\backslash S)=0
      \dott
    \end{eqn}
  }
  \argument{
    \lref{eq7.1};
    \lref{eq0.2c.9};
  }{
    that
    \begin{eqn}
      \llabel{eq8}
      \TVmeasure{\mu}(E)=\boldsymbol{\mu}_1(E)+\boldsymbol{\mu}_2(E)=0
      \dott
    \end{eqn}
  }
  \argument{
    \lref{eq8};
  }[verbs=ep]{
    \lref{item2}\dott
  }
  \startnewargseq
  \argument{
    \lref{ass1}; \lref{eq0.4}; \lref{eq0.7};
  }[verbs=ep]{
    \lref{item3}\dott
  }
\end{aproof}

\begin{athm}{lemma}{lem:AC-gluing}
  Let $N\in\N$, $a\in\R$, $b\in[a,\infty)$,
  let $u\in C([a,b],\R)$ have finite length\cfadd{definition:curve-length},
  for every $j\in\{1,2,\dots,N\}$ let $E_j\subseteq[a,b]$ be closed,
  let $\mathfrak a_j\colon[a,b]\to\R$ be absolutely continuous,
  and assume $\mathfrak a_j|_{E_j}=u|_{E_j}$,
  and assume for all $c\in[a,b)$, $d\in(c,b]$ with
  $[c,d]\cap(\midcup_{j=1}^{N}E_j)=\emptyset$
  that $u|_{[c,d]}$ is absolutely continuous\dott
  \cfload[.]
  Then $u$ is absolutely continuous\dott \cfout[.]
\end{athm}
\begin{aproof}
  Throughout this proof assume without loss of generality that $a<b$\dott
  In the following let $\lambda$ be the Lebesgue measure on $\R$,
  let $U\colon\R\to\R$ and $W_j\colon\R\to\R$,
  $j\in\{1,2,\dots,N\}$, satisfy for all $t\in\R$,
  $j\in\{1,2,\dots,N\}$ that
  \begin{eqn}
    \llabel{setting1}
    U(t)=u(\min\{b,\max\{a,t\}\})
    \qquad \text{and}\qquad
    W_j(t)=\mathfrak a_j(\min\{b,\max\{a,t\}\}),
  \end{eqn}
  let $G=(a,b)\backslash(\midcup_{j=1}^{N}E_j)$, and let
  \begin{eqn}
    \llabel{setting3}
    \mathcal J=\{(p,q)\in\mathbb Q^2\colon p<q,\ [p,q]\subseteq G\}
    \dott
  \end{eqn}
  \cfload[.]
  \startnewargseq
  \argument{
    the assumption that $u$ has finite length\cfadd{definition:curve-length};
    \lref{setting1};
    the assumption that $u$ is continuous;
    \itref{lem:BV-locality}{item:derivative-measure,item3}
    (applied with $E\curvearrowleft\emptyset$
    in the notation of \cref{lem:BV-locality});
  }{
    that there exists a finite signed Radon measure $\mu$ on
    $\R$\cfadd{def:signed-radon-measure} which satisfies for all
    $\varphi\in C^\infty(\R,\R)$ with compact support that
    \begin{eqn}
      \llabel{eq:measure-tests}
      \signedint_{\R}\varphi(t)\,\mu(\ddd t)
        =-\int_{\R}U(t)\varphi'(t)\dd t
    \end{eqn}
    and which satisfies for all $s,t\in[a,b]$ with $s\leq t$ that
    \begin{eqn}
      \llabel{eq1}
      \mu((s,t])=u(t)-u(s)
      \qquad \text{and}\qquad
      \mu(\{s\})=0
      \dott
    \end{eqn}
  }
  \startnewargseq
  \argument{
    \lref{eq1};
    the fact that for every signed Radon measure $\nu$ on $\R$
    and every compact $K\subseteq\R$ it holds that
    $\TVmeasure{\nu}(K)<\infty$;
    the fact that for all $t\in\R$ it holds that
    $\TVmeasure{\mu}(\{t\})=\lvert\mu(\{t\})\rvert$;
  }{
    that
    \begin{eqn}
      \llabel{eq2}
      \TVmeasure{\mu}([a,b])<\infty
      \qquad \text{and}\qquad
      \TVmeasure{\mu}(\{a,b\})=0
      \dott
    \end{eqn}
  }
  \argument{
    the assumption that for all $j\in\{1,2,\dots,N\}$ it holds that
    $\mathfrak a_j$ is absolutely continuous;
    the fundamental theorem of calculus for absolutely continuous functions;
  }{
    that for all $j\in\{1,2,\dots,N\}$ it holds that
    $\mathfrak a_j$ is differentiable $\lambda$-almost everywhere on $(a,b)$
    and that for all $t\in[a,b]$ it holds that
    \begin{eqn}
      \llabel{eq3.1}
      \int_a^b\lvert(\mathfrak a_j)'(s)\rvert\dd s<\infty
      \qquad \text{and}\qquad
      \mathfrak a_j(t)=\mathfrak a_j(a)
        +\int_a^t(\mathfrak a_j)'(s)\dd s
      \dott
    \end{eqn}
  }
  In the following for every $j\in\{1,2,\dots,N\}$ let
  $g_j\colon\R\to\R$ be Borel measurable and satisfy for all
  $x\in\R\backslash(a,b)$, $y\in(a,b)$ for which
  $\mathfrak a_j$ is differentiable at $y$ that
  \begin{eqn}
    \llabel{setting4}
    g_j(x)=0
    \qquad \text{and}\qquad
    g_j(y)=(\mathfrak a_j)'(y)
  \end{eqn}
  and let
  $\sigma_j^+,\sigma_j^-\colon\mathcal B(\R)\to[0,\infty)$,
  $\sigma_j\colon\mathcal B(\R)\to\R$, and $P_j\subseteq\R$ satisfy for all
  $B\in\mathcal B(\R)$ that 
  \begin{eqn}
 \textstyle   \llabel{setting5}
     \sigma_j^+(B)=\int_B\max\{g_j(t),0\}\dd t,\qquad\sigma_j^-(B)=\int_B\max\{-g_j(t),0\}\dd t,\qquad  \sigma_j(B)=\sigma_j^+(B)-\sigma_j^-(B),
  \end{eqn}
  and
  \begin{eqn}
    \llabel{setting6}
    P_j=\{t\in\R\colon g_j(t)\geq0\}
    \dott
  \end{eqn}
  \cfload[.]
  \startnewargseq
  \argument{
    \lref{eq3.1};
    \lref{setting4};
    \lref{setting5};
    the fact that for all finite nonnegative Borel measures $\zeta$
    on $\R$ it holds that $\zeta$ is a Radon measure on $\R$;
  }{
    that for all $j\in\{1,2,\dots,N\}$ it holds that $\sigma_j^+$ and
    $\sigma_j^-$ are finite nonnegative Radon measures on
    $\R$\cfadd{def:radon-measure} and that
    \begin{eqn}
      \llabel{eq3.2}
      \sigma_j^+(\R)+\sigma_j^-(\R)
        =\int_{\R}\lvert g_j(t)\rvert\dd t
        =\int_a^b\lvert(\mathfrak a_j)'(t)\rvert\dd t<\infty
      \dott
    \end{eqn}
  }
  \argument{
    \lref{setting5};
    \lref{setting6};
  }{
    that for all $j\in\{1,2,\dots,N\}$ it holds that
    \begin{eqn}
      \llabel{eq3.3}
      \sigma_j^+(\R\backslash P_j)=\sigma_j^-(P_j)=0
      \dott
    \end{eqn}
  }
  \argument{
    \lref{eq3.3};
    \lref{setting5};
    \lref{eq3.2};
    the fact that for all $j\in\{1,2,\dots,N\}$ it holds that
    $\sigma_j^+$ and $\sigma_j^-$ are Radon measures on $\R$;
    the fact that for all $j\in\{1,2,\dots,N\}$ it holds that $P_j$ is Borel measurable;
  }{
    that for all $j\in\{1,2,\dots,N\}$ it holds that $\sigma_j$ is a
    finite signed Radon measure on $\R$ and that for all
    $B\in\mathcal B(\R)$ it holds that
    \begin{eqn}
      \llabel{eq3.4}
      \TVmeasure{\sigma_j}(B)
        =\sigma_j^+(B)+\sigma_j^-(B)
        =\int_B\lvert g_j(t)\rvert\dd t
      \dott
    \end{eqn}
  }
  \argument{
    \lref{setting1};
    \lref{eq3.1};
    \lref{setting4};
  }{
    that for all $j\in\{1,2,\dots,N\}$ it holds that $W_j$ is bounded
    and continuous and that for all $t\in\R$ it holds that
    \begin{eqn}
      \llabel{eq3.5}
      W_j(t)=\mathfrak a_j(a)+\int_{(-\infty,t)}g_j(s)\dd s
      \dott
    \end{eqn}
  }
  \argument{
    the assumption that for all $j\in\{1,2,\dots,N\}$ it holds that
    $\mathfrak a_j$ is absolutely continuous;
    the integration by parts formula for absolutely continuous functions;
  }{
    that for all $j\in\{1,2,\dots,N\}$ and all compactly supported
    $\varphi\in C^\infty(\R,\R)$ it holds that
    \begin{eqn}[a]
      \llabel{eq3.5a}
      \int_a^b\mathfrak a_j(t)\varphi'(t)\dd t
        =\mathfrak a_j(b)\varphi(b)-\mathfrak a_j(a)\varphi(a)
    -\int_a^b(\mathfrak a_j)'(t)\varphi(t)\dd t
      \dott
    \end{eqn}
  }
  \argument{
    \lref{eq3.5a};
    \lref{setting1};
    \lref{setting4};
    \lref{setting5};
    \lref{eq3.5};
    the fundamental theorem of calculus
  }{
    that for all $j\in\{1,2,\dots,N\}$ and all compactly supported
    $\varphi\in C^\infty(\R,\R)$ it holds that
    \begin{eqn}[a]
      \llabel{eq3.6}
      -\int_{\R}W_j(t)\varphi'(t)\dd t
        &=-\mathfrak a_j(a)\int_{-\infty}^a\varphi'(t)\dd t
          -\int_a^b\mathfrak a_j(t)\varphi'(t)\dd t
  -\mathfrak a_j(b)\int_b^{\infty}\varphi'(t)\dd t
      \\
        &=-\mathfrak a_j(a)\varphi(a)
          -\biggl(\int_a^b\mathfrak a_j(t)\varphi'(t)\dd t\biggr)
          +\mathfrak a_j(b)\varphi(b)
      \\
        &=\int_a^b(\mathfrak a_j)'(t)\varphi(t)\dd t
      =\int_{\R}\varphi(t)g_j(t)\dd t
          =\signedint_{\R}\varphi(t)\,\sigma_j(\ddd t)
      \dott
    \end{eqn}
  }
  \argument{
    \lref{setting4};
    \lref{eq3.4};
  }{
    that for all $j\in\{1,2,\dots,N\}$,
    $B\in\mathcal B(\R)$ it holds that
    \begin{eqn}
      \llabel{eq3}
      \TVmeasure{\sigma_j}(B)
        =\int_{B\cap(a,b)}\lvert(\mathfrak a_j)'(t)\rvert\dd t
      \dott
    \end{eqn}
  }
  \argument{
    \lref{eq3.1};
    the triangle inequality;
  }{
    that for all $j\in\{1,2,\dots,N\}$, $n\in\N$,
    $\tau_0,\tau_1\ldots,\tau_n\in[a,b]$ with
    $a=\tau_0\leq\tau_1\leq\ldots\leq\tau_n=b$ it holds that
    \begin{eqn}[a]
      \llabel{eq4}
      &\sum_{k=1}^{n}
        \lvert(u-\mathfrak a_j)(\tau_k)-(u-\mathfrak a_j)(\tau_{k-1})\rvert
    \leq\sum_{k=1}^{n}\lvert u(\tau_k)-u(\tau_{k-1})\rvert
        +\int_a^b\lvert(\mathfrak a_j)'(t)\rvert\dd t
      \dott
    \end{eqn}
  }
  \argument{
    \lref{eq4};
    the assumption that $u$ has finite length\cfadd{definition:curve-length};
    \lref{setting1};
    the fact that for all $j\in\{1,2,\dots,N\}$ it holds that
    $u-\mathfrak a_j$ is continuous;
    the fact that for all $j\in\{1,2,\dots,N\}$ it holds that
    $(u-\mathfrak a_j)|_{E_j}=0$;
    \itref{lem:BV-locality}{item:derivative-measure,item2}
    (applied for every $j\in\{1,2,\dots,N\}$ with $u\curvearrowleft u-\mathfrak a_j$,
    $E\curvearrowleft E_j$, $U\curvearrowleft U-W_j$ in the notation of \cref{lem:BV-locality});
  }{
    that for all $j\in\{1,2,\dots,N\}$ there exists a finite signed
    Radon measure $\rho_j$ on $\R$ with $\TVmeasure{\rho_j}(E_j)=0$
    which satisfies for all $\varphi\in C^\infty(\R,\R)$
    with compact support that
    \begin{eqn}
      \llabel{eq5}
      \signedint_{\R}\varphi(t)\,\rho_j(\ddd t)
        =-\int_{\R}\bigl(U(t)-W_j(t)\bigr)\varphi'(t)\dd t
      \dott
    \end{eqn}
  }
  \startnewargseq
  \argument{
    \lref{eq:measure-tests};
    \lref{eq3.6};
    \lref{eq5};
  }{
    that for all $j\in\{1,2,\dots,N\}$,
    $\varphi\in C^\infty(\R,\R)$ with compact support it holds that
    \begin{eqn}
      \llabel{eq5.1}
      \signedint_{\R}\varphi(t)\,\mu(\ddd t)
        =\signedint_{\R}\varphi(t)\,\rho_j(\ddd t)
          +\signedint_{\R}\varphi(t)\,\sigma_j(\ddd t)
      \dott
    \end{eqn}
  }
  \argument{
    \lref{eq5.1};
    the fact that $\mu$ is a finite signed Radon measure on $\R$;
    the fact that for all $j\in\{1,2,\dots,N\}$ it holds that
    $\rho_j$ and $\sigma_j$ are finite signed Radon measures on $\R$;
  }{
    that for all $j\in\{1,2,\dots,N\}$,
    $B\in\mathcal B(\R)$ it holds that
    \begin{eqn}
      \llabel{eq5.2}
      \mu(B)=\rho_j(B)+\sigma_j(B)
      \dott
    \end{eqn}
  }
  \argument{
    \lref{eq5.2};\lref{eq3};
    \lref{eq5};
  }{
    that for all $B\in\mathcal B(\R)$,
    $j\in\{1,2,\dots,N\}$ with $\lambda(B)=0$ it holds that
    \begin{eqn}
      \llabel{eq6}
      0\leq\TVmeasure{\mu}(B\cap E_j)
        \leq\TVmeasure{\rho_j}(E_j)
          +\TVmeasure{\sigma_j}(B\cap E_j)=0
      \dott
    \end{eqn}
  }
  \argument{
    the fact that $G=(a,b)\backslash(\midcup_{j=1}^{N}E_j)$;
    \lref{setting3};
    the fact that for all $j\in\{1,2,\dots,N\}$ it holds that $E_j$
    is closed;
    the fact that $\mathbb Q$ is countable and dense in $\R$;
  }{
    that $G$ is open, that $\mathcal J$ is at most countable, and that
    \begin{eqn}
      \llabel{eq7}
      G=\bigcup_{(p,q)\in\mathcal J}(p,q)
      \dott
    \end{eqn}
  }
  \argument{
    the fact that $G=(a,b)\backslash(\midcup_{j=1}^{N}E_j)$;
    \lref{setting3};
    the fact that $u$ is absolutely continuous on every
    compact subinterval of $[a,b]\backslash(\midcup_{j=1}^{N}E_j)$;
  }{
    that for all $(p,q)\in\mathcal J$ it holds that
    $u|_{[p,q]}$ is absolutely continuous\dott
    \llabel{eq7.1}
  }
  \argument{
    \lref{eq7.1};
    the fundamental theorem of calculus for absolutely continuous functions;
  }{
    that for all $(p,q)\in\mathcal J$ it holds that $u$ is differentiable
    $\lambda$-almost everywhere on $(p,q)$ and that
    \begin{eqn}
      \llabel{eq7.2}
      \int_p^q\lvert u'(t)\rvert\dd t<\infty
      \dott
    \end{eqn}
  }
  In the following for every $(p,q)\in\mathcal J$ let
  $h_{p,q}\colon(p,q)\to\R$ and
  $\eta_{p,q}\colon\mathcal B(\R)\to\R$ satisfy for all
  $B\in\mathcal B(\R)$, $t\in(p,q)$ for which $u$ is differentiable at $t$ that
  \begin{eqn}
    \llabel{setting7}
    h_{p,q}(t)=u'(t)
    \qquad\text{and}\qquad
    \eta_{p,q}(B)=\int_{B\cap(p,q)}h_{p,q}(s)\dd s
    \dott
  \end{eqn}
  \cfload[.]
  \startnewargseq
  \argument{
    \lref{eq7.2};
    \lref{setting7};
    the fact that every finite nonnegative Borel measure on $\R$
    is a Radon measure on $\R$;
  }{
    that for all $(p,q)\in\mathcal J$ it holds that $\eta_{p,q}$
    is a finite signed Radon measure on $\R$ and that for all
    $B\in\mathcal B(\R)$ it holds that
    \begin{eqn}[a]
      \llabel{eq7.3}
      \TVmeasure{\eta_{p,q}}(B)
        =\int_{B\cap(p,q)}\max\{h_{p,q}(t),0\}\dd t
          +\int_{B\cap(p,q)}\max\{-h_{p,q}(t),0\}\dd t
    =\int_{B\cap(p,q)}\lvert h_{p,q}(t)\rvert\dd t
      \dott
    \end{eqn}
  }
  \argument{
    \lref{setting1};
    \lref{setting3};
    \lref{eq:measure-tests};
    the fact that for all $(p,q)\in\mathcal J$ it holds that
    $u|_{[p,q]}$ is absolutely continuous;
    the integration by parts formula for absolutely continuous functions;
  }{
    that for all $(p,q)\in\mathcal J$ and all $\varphi\in C^\infty(\R,\R)$
    with compact support in $(p,q)$ it holds that
    \begin{eqn}[a]
      \llabel{eq7.4}
      \signedint_{\R}\varphi(t)\,\mu(\ddd t)
        &=-\int_{\R}U(t)\varphi'(t)\dd t
          =-\int_p^q u(t)\varphi'(t)\dd t
      \\
        &=u(p)\varphi(p)-u(q)\varphi(q)
          +\int_p^q u'(t)\varphi(t)\dd t=\int_p^q u'(t)\varphi(t)\dd t
      \dott
    \end{eqn}
  }
  \argument{
    \lref{eq7.4};
    \lref{setting7};
    \lref{eq7.3};
  }{
    that for all $(p,q)\in\mathcal J$ and all $\varphi\in C^\infty(\R,\R)$
    with compact support in $(p,q)$ it holds that
    \begin{eqn}
      \llabel{eq7.5}
      \signedint_{\R}\varphi(t)\,\mu(\ddd t)
        =\int_p^q\varphi(t)h_{p,q}(t)\dd t
        =\signedint_{\R}\varphi(t)\,\eta_{p,q}(\ddd t)
      \dott
    \end{eqn}
  }
  \argument{
    \lref{eq7.5};
    \lref{eq7.3};
    the fact that $\mu$ is a signed Radon measure on $\R$;
  }{
    that for all $(p,q)\in\mathcal J$, $B\in\mathcal B(\R)$
    it holds that
    \begin{eqn}
      \llabel{eq7.6}
      \mu(B\cap(p,q))=\eta_{p,q}(B\cap(p,q))
      \dott
    \end{eqn}
  }
  \argument{
    \lref{eq7.6};
    \lref{setting7};
    \lref{eq7.3};
  }{
    that for all $(p,q)\in\mathcal J$, $B\in\mathcal B(\R)$
    it holds that
    \begin{eqn}[a]
      \llabel{eq7.7}
      \TVmeasure{\mu}(B\cap(p,q))
        =\TVmeasure{\eta_{p,q}}(B\cap(p,q))
      =\int_{B\cap(p,q)}\lvert h_{p,q}(t)\rvert\dd t
          =\int_{B\cap(p,q)}\lvert u'(t)\rvert\dd t
      \dott
    \end{eqn}
  }
  \argument{
    \lref{eq7.7};
    \lref{eq7.2};
  }{
    that for all $(p,q)\in\mathcal J$,
    $B\in\mathcal B(\R)$ with $\lambda(B)=0$ it holds that
    \begin{eqn}
      \llabel{eq8}
      \TVmeasure{\mu}(B\cap(p,q))
        =\int_{B\cap(p,q)}\lvert u'(t)\rvert\dd t=0
      \dott
    \end{eqn}
  }
  \argument{
    \lref{eq8};
    \lref{eq7};
    the fact that $\TVmeasure{\mu}$ is countably subadditive;
  }{
    that for all $B\in\mathcal B(\R)$ with $\lambda(B)=0$ it holds that
    \begin{eqn}
      \llabel{eq9}
      0\leq\TVmeasure{\mu}(B\cap G)
        \leq\sum_{(p,q)\in\mathcal J}\TVmeasure{\mu}(B\cap(p,q))=0
      \dott
    \end{eqn}
  }
  \argument{
    \lref{eq9};
    the fact that $G=(a,b)\backslash(\midcup_{j=1}^{N}E_j)$;
    \lref{eq2};
    \lref{eq6};
    the fact that $\TVmeasure{\mu}$ is subadditive;
  }{
    that for all $B\in\mathcal B(\R)$ with $B\subseteq[a,b]$
    and $\lambda(B)=0$ it holds that
    \begin{eqn}
      \llabel{eq10}
      0\leq\TVmeasure{\mu}(B)
        \leq\TVmeasure{\mu}(B\cap G)
          +\sum_{j=1}^{N}\TVmeasure{\mu}(B\cap E_j)
          +\TVmeasure{\mu}(\{a,b\})=0
      \dott
    \end{eqn}
  }
  \argument{
    \lref{eq10};
    \lref{eq2};
    \cite[Theorem~6.11]{Rudin87}
  }{
    that for all $\varepsilon\in(0,\infty)$ there exists
    $\delta\in(0,\infty)$ such that for all $B\in\mathcal B(\R)$
    with $B\subseteq[a,b]$ and $\lambda(B)<\delta$ it holds that
    \begin{eqn}
      \llabel{eq11}
      \TVmeasure{\mu}(B)<\varepsilon
      \dott
    \end{eqn}
  }
  \argument{
    \lref{eq11};
    \lref{eq1};
  }{
    that for all $\varepsilon\in(0,\infty)$ there exists
    $\delta\in(0,\infty)$ which satisfies for all $m\in\N$,
    $s_1,t_1,s_2,t_2,\ldots,s_m,t_m\in[a,b]$ with
    $s_1<t_1\leq s_2<t_2\leq\ldots\leq s_m<t_m$ and
    $\sum_{k=1}^{m}(t_k-s_k)<\delta$ that
    \begin{eqn}[a]
      \llabel{eq12}
      \sum_{k=1}^{m}\lvert u(t_k)-u(s_k)\rvert
        =\sum_{k=1}^{m}\lvert\mu((s_k,t_k])\rvert
     \leq\sum_{k=1}^{m}\TVmeasure{\mu}((s_k,t_k])
        =\textstyle\TVmeasure{\mu}\bigl(\midcup_{k=1}^{m}(s_k,t_k]\bigr)
        <\varepsilon
      \dott
    \end{eqn}
  }
  \argument{
    \lref{eq12};
  }[verbs=ep]{
    that $u$ is absolutely continuous\dott
  }
\end{aproof}
\subsection{Uniform upper bounds for the length of GF trajectories close to the origin}

\begin{athm}{lemma}{lem:ambient-annular}
  Let $n\in\N\backslash\{1\}$, $K,R\in(0,\infty)$,
  let $V=\{z\in\R^n\colon0<\norm{z}<R\}$,
  let $W\subseteq V$ be definable and open, and let
  $g\in C^1(W,\R)$ be definable\cfadd{definition:definable-function}\dott
  \cfload[.]
  Then there exist $\varrho\in(0,R)$, $\mathfrak{C}\in(0,\infty)$ such that
  for every $\rho\in(0,\varrho)$,
  every open $U\subseteq\R^n$,
  every $F\in C^1(U,\R)$,
  every interval $I\subseteq\R$,
  every locally absolutely continuous $y\colon I\to U\backslash\{0\}$,
  every Lebesgue measurable $\lambda\colon I\to[0,\infty)$,
  and every Lebesgue measurable $B\colon I\to\R^{n\times n}$
  with the property that $U\cap V\subseteq W$ and
  $F|_{U\cap V}=g|_{U\cap V}$,
  with the property that for Lebesgue-almost all $t\in I$ it holds
  that $y$ is differentiable at $t$, that $B(t)$ is symmetric and positive-definite,
  and that
  \begin{eqn}
    \llabel{ass1}
    y'(t)=-\lambda(t)B(t)\grad F(y(t)),
  \end{eqn}
  and with the property that for Lebesgue-almost all
  $t\in\{s\in I\colon\norm{y(s)}<R\}$ it holds that
  \begin{eqn}
    \llabel{ass2}
    \max\nolimits_{v\in\R^n,\,\norm{v}=1}
      \langle B(t)v,v\rangle
    \leq K\min\nolimits_{v\in\R^n,\,\norm{v}=1}
      \langle B(t)v,v\rangle
  \end{eqn}
  we have that
  \begin{eqn}
    \llabel{concl1}
    \int_{\{t\in I\colon\rho/2\leq\norm{y(t)}\leq\rho\}}
      \norm{y'(t)}\dd t\leq \mathfrak{C}\rho\qquad\text{and}\qquad
    \int_{\{t\in I\colon\norm{y(t)}\leq\rho\}}
      \norm{y'(t)}\dd t\leq2\mathfrak{C}\rho
    \dott
  \end{eqn}
  \cfout[.]
\end{athm}
\begin{aproof}
  Throughout this proof let $\varphi\colon\R\to(-1,1)$ and
  $\mathfrak g\colon W\to(-1,1)$
  satisfy for all $e\in\R$ that
  \begin{eqn}
    \llabel{normalization-map}
    \varphi(e)=\frac{e}{\sqrt{1+e^2}}
    \qquad\text{and}\qquad
    \mathfrak g=\varphi\circ g,
  \end{eqn}
  let $\mathbb{V}\subseteq\R^n$ and
  $\Omega\subseteq(0,R)\times\R^n$ satisfy
  \begin{eqn}
    \llabel{setting1}
    \mathbb{V}=\{z\in\R^n\colon1/2\leq\norm{z}\leq1\}
    \qquad\text{and}\qquad
    \Omega=\{(\rho,z)\in(0,R)\times\R^n\colon\rho z\in W\},
  \end{eqn}
  let $f=(f(\rho,z))_{(\rho,z)\in\Omega}\colon\Omega\to\R$ and
  $H\colon\Omega\to[0,\infty)$ satisfy
  for all $(\rho,z)\in\Omega$ that
  \begin{eqn}
    \llabel{setting2}
    f(\rho,z)=\mathfrak g(\rho z)
    \qquad\text{and}\qquad
    H(\rho,z)=\rho\norm{\grad\mathfrak g(\rho z)},
  \end{eqn}
  and let $\mathcal E\subseteq(0,R)\times\R$ satisfy
  \begin{eqn}
    \llabel{setting3}
    \mathcal E=\{(\rho,e)\in(0,R)\times\R
      \colon(\exists\,z\in\mathbb{V}\colon
      (\rho,z)\in\Omega,\ f(\rho,z)=e)\}
    \dott
  \end{eqn}
  \cfload[.]
  \startnewargseq
  \argument{
    \lref{normalization-map};
    the assumption that $g\in C^1(W,\R)$ is definable;
    the chain rule;
  }{
    that $\mathfrak g\in C^1(W,(-1,1))$ is definable and that
    for all $e\in\R$ it holds that
    \begin{eqn}
      \llabel{normalization-derivative}
      \varphi'(e)=(1+e^2)^{-3/2}>0
      \dott
    \end{eqn}
  }
  \argument{
    \lref{setting1};
    \lref{setting2};
    the fact that $\mathfrak g\in C^1(W,(-1,1))$ is definable;
    the assumption that $W$ is open and definable;
    the chain rule;
  }{
    that $\mathbb{V}$ is compact and definable, that $\Omega$ is open and definable, that
    $f\in C^1(\Omega,\R)$ is definable, that $H$ is continuous and definable, and that for all
    $(\rho,z)\in\Omega$ it holds that
    \begin{eqn}
      \llabel{eq1}
      \nabla_z f(\rho,z)=\rho\grad\mathfrak g(\rho z)
      \qquad\text{and}\qquad
      \norm{\nabla_z f(\rho,z)}=H(\rho,z)
      \dott
    \end{eqn}
  }
  \argument{
    \lref{setting3};
    the fact that $\Omega$, $f$, and $\mathbb{V}$ are definable;
    \cite[Theorem~1.13]{Coste};
  }{
    that $\mathcal E$ is definable and that for all
    $(\rho,e)\in\mathcal E$ it holds that
    \begin{eqn}
      \llabel{eq:nonempty-fiber}
      \{z\in\mathbb{V}\colon(\rho,z)\in\Omega,\ f(\rho,z)=e\}
      \neq\emptyset
      \dott
    \end{eqn}
  }
  \argument{
    \lref{eq:nonempty-fiber};
    the fact that $H(\Omega)\subseteq[0,\infty)$;
  }{
    that for all $(\rho,e)\in\mathcal E$ it holds that
    \begin{eqn}
      \llabel{eq:infimum-values}
      \emptyset\neq\{H(\rho,z)\colon
        z\in\mathbb{V},\ (\rho,z)\in\Omega,\ f(\rho,z)=e\}
      \subseteq[0,\infty)
      \dott
    \end{eqn}
  }
  \argument{
    \lref{eq:infimum-values};
  }{
    that there exists
    $\mathfrak h\colon\mathcal E\to[0,\infty)$
    which satisfies for all $(\rho,e)\in\mathcal E$ that
    \begin{eqn}
      \llabel{eq:infimum}
      \mathfrak h(\rho,e)=\inf\{H(\rho,z)\colon
        z\in\mathbb{V},\ (\rho,z)\in\Omega,\ f(\rho,z)=e\}
      \dott
    \end{eqn}
  }
  \startnewargseq
  \argument{
    \lref{eq:infimum};
  }{
    that
    \begin{eqn}[a]
      \llabel{eq:infimum-characterization}
      &\{(\rho,e,\mathfrak h(\rho,e))\colon(\rho,e)\in\mathcal E\}\\
      &=\Bigl\{(\rho,e,c)\in\mathcal E\times\R\colon\Bigl[
        \bigl(\forall\,z\in\mathbb{V},\ (\rho,z)\in\Omega,\ f(\rho,z)=e\colon
          c\leq H(\rho,z)\bigr)\\
      &\qquad\wedge\Bigl(\forall\,b\in(c,\infty)\ \exists\,z\in\mathbb{V}\colon
        \bigl[((\rho,z)\in\Omega)\wedge(f(\rho,z)=e)\wedge(H(\rho,z)<b)\bigr]\Bigr)
        \Bigr]\Bigr\}
      \dott
    \end{eqn}
  }
  \argument{
    \lref{eq:infimum-characterization};
    the fact that $\mathcal E$, $\mathbb{V}$, $\Omega$, $f$, and $H$
    are definable;
    \cite[Theorem~1.13]{Coste};
  }{
    that $\mathfrak h$ is definable\dott
  }
  In the following let $\mathcal Z\subseteq\mathcal E$ satisfy
  \begin{eqn}
    \llabel{eq:zero-set}
    \mathcal Z=\{(\rho,e)\in\mathcal E\colon\mathfrak h(\rho,e)=0\}
    \dott
  \end{eqn}
  \startnewargseq
  \argument{
    \lref{eq:zero-set};
    the fact that $\mathfrak h$ is definable;
    \cite[Theorem~1.13]{Coste};
  }{
    that $\mathcal Z$ is definable\dott
  }
  \startnewargseq
  \argument{
    \lref{eq:infimum};
    \lref{eq:zero-set};
    the fact that $\mathcal Z$ is definable;
    \cite[Theorem~3.1]{Coste};
  }{
    that there exists a definable
    $\bfX=\bigl(\bfX(\rho,\eta,e)\bigr)_{(\rho,e)\in\mathcal Z,\,\eta\in(0,\infty)}
      =\bigl(\bfX_1(\rho,\eta,e),\dots,\bfX_n(\rho,\eta,e)\bigr)_{(\rho,e)\in\mathcal Z,\,\eta\in(0,\infty)}
      \colon\{(\rho,\eta,e)\in(0,R)\times(0,\infty)\times\R
      \colon(\rho,e)\in\mathcal Z\}\to\mathbb V$
    which satisfies for all $(\rho,e)\in\mathcal Z$,
    $\eta\in(0,\infty)$ that
    \begin{eqn}
      \llabel{eq:small-selector}
      (\rho,\bfX(\rho,\eta,e))\in\Omega,\qquad
      f(\rho,\bfX(\rho,\eta,e))=e,\qquad\text{and}\qquad
      H(\rho,\bfX(\rho,\eta,e))<\eta
      \dott
    \end{eqn}
  }
  \startnewargseq
  \argument{
    the fact that $\bfX$ is definable;
    \cite[Theorem~6.7]{Coste};
    \cite[Theorem~6.6]{Coste};
  }{
    that there exists a cylindrical definable $C^1$ cell decomposition
    $\mathscr Q$ of $\R^3$,
    which is compatible with
    $\{(\rho,\eta,e)\in\R^3\colon(\rho,e)\in\mathcal Z,\ \eta>0\}$
    and which satisfies for all $S\in\mathscr Q$ with
    $S\subseteq\{(\rho,\eta,e)\in\R^3\colon(\rho,e)\in\mathcal Z,\ \eta>0\}$ that
    \begin{eqn}
      \llabel{eq:selector-smooth-cells}
      \bfX|_S\in\COne{S}{\R^n}
      \cfadd{definition:planar-cell-decomposition}\dott
    \end{eqn}
  }
  In the following let $\mathcal U\subseteq\R^3$ satisfy
  \begin{eqn}[a]
    \llabel{eq:selector-derivative-domain}
    \mathcal U=\bigl\{(\rho,\eta,e)\in\R^3\colon
      (\rho,e)\in\mathcal Z,\ \eta>0,\
      \bigl[\exists\,S\in\mathscr Q,\ \delta\in(0,\infty)\colon
        \{(\rho,\eta)\}\times(e-\delta,e+\delta)\subseteq S\bigr]\bigr\}
    \dott
  \end{eqn}
  \startnewargseq
  \argument{
    \lref{eq:selector-smooth-cells};
    \lref{eq:selector-derivative-domain};
    the chain rule;
    \cite[Theorem~1.13]{Coste};
  }{
    that $\mathcal U$ is definable and that for all
    $(\rho,\eta,e)\in\mathcal U$ there exists
    $\delta\in(0,\infty)$ such that
    \begin{eqn}
      \llabel{eq:selector-local-regularity}
      \bigl((e-\delta,e+\delta)\ni t\mapsto\bfX(\rho,\eta,t)\in\R^n\bigr)
        \in C^1\bigl((e-\delta,e+\delta),\R^n\bigr)
      \dott
    \end{eqn}
  }
  \argument{
    \lref{eq:selector-local-regularity};
    the fact that $\bfX$ is definable;
    \cite[Theorem~1.13]{Coste};
  }{
    that there exists a definable
    $\bfY=(\bfY_1,\dots,\bfY_n)\colon\mathcal U\to\R^n$
    which satisfies for all $(\rho,\eta,e)\in\mathcal U$ that
    \begin{eqn}
      \llabel{eq:selector-restricted-derivative}
      \bfY(\rho,\eta,e)=(\partial_e\bfX)(\rho,\eta,e)
      \dott
    \end{eqn}
  }
  \startnewargseq
  \argument{
    \lref{eq:selector-smooth-cells};
    \lref{eq:selector-derivative-domain};
    \lref{eq:selector-restricted-derivative};
    \cite[Theorem~6.6]{Coste}
  }{
    that there exists a cylindrical definable $C^1$ cell decomposition
    $\mathscr F$ of $\R^3$ which is a refinement of
    $\mathscr Q$\cfadd{definition:partition-refinement},
    which is compatible with
    $\mathcal U$\cfadd{definition:partition-compatibility},
    and which satisfies for every open interval $I\subseteq\R$
    and every $S\in\mathscr F$, $\rho\in(0,R)$,
    $\eta\in(0,\infty)$, $k\in\{1,\dots,n\}$, $s,t\in I$ with
    $S\subseteq\{(\rho,\eta,e)\in\R^3\colon(\rho,e)\in\mathcal Z,\ \eta>0\}$
    and $I=\{u\in\R\colon(\rho,\eta,u)\in S\}$ that
    \begin{eqn}
      \llabel{eq:selector-derivative-signs}
      \operatorname{sign}\bigl(\bfY_k(\rho,\eta,s)\bigr)
        =\operatorname{sign}\bigl(\bfY_k(\rho,\eta,t)\bigr)
      \dott
    \end{eqn}
  }
  \startnewargseq
  \argument{
    \lref{eq:selector-smooth-cells};
    \lref{eq:selector-restricted-derivative};
    \lref{eq:selector-derivative-signs};
    the mean value theorem;
  }{
    that there exist definable partitions
    $\mathscr L_{\rho,\eta}$ of
    $\{e\in\R\colon(\rho,e)\in\mathcal Z\}$,
    $\rho\in(0,R)$, $\eta\in(0,\infty)$,
    which satisfy for all $\rho\in(0,R)$, $\eta\in(0,\infty)$ that
    \begin{enumerate}[label=(\roman*)]
      \item \llabel{item:selector-pieces} for all
      $J\in\mathscr L_{\rho,\eta}$ it holds that $J$ is a singleton
      or an open interval,
      \item \llabel{item:selector-regularity} for all open intervals
      $J\in\mathscr L_{\rho,\eta}$ it holds that
      $\bigl(J\ni e\mapsto\bfX(\rho,\eta,e)\in\R^n\bigr)\in C^1(J,\R^n)$,
      \item \llabel{item:selector-monotonicity} for all open intervals
      $J\in\mathscr L_{\rho,\eta}$, $k\in\{1,2,\dots,n\}$ it holds that
      $J\ni e\mapsto\bfX_k(\rho,\eta,e)\in\R$
      is non-decreasing or non-increasing, and
      \item \llabel{item:selector-bound} it holds that
      $\#\mathscr L_{\rho,\eta}\leq\#\mathscr F$\dott
    \end{enumerate}
  }
  \startnewargseq
  \argument{
    \itref{lem:ambient-annular}{item:selector-regularity,item:selector-monotonicity};
    the fact that
    $\bfX\bigl(\{(\rho,\eta,e)\in\R^3
      \colon(\rho,e)\in\mathcal Z,\ \eta>0\}\bigr)\subseteq[-1,1]^n$;
    the fundamental theorem of calculus;
  }{
    that for every $\rho\in(0,R)$, $\eta\in(0,\infty)$, every 
    open interval $J\in\mathscr L_{\rho,\eta}$,
    and every $k\in\{1,2,\dots,n\}$, $s,t\in J$ with $s<t$ it holds that
    \begin{eqn}
      \llabel{eq:selector-coordinate-variation}
      \int_s^t\bigl|(\partial_e\bfX_k)(\rho,\eta,e)\bigr|\dd e
        =\bigl|\bfX_k(\rho,\eta,t)-\bfX_k(\rho,\eta,s)\bigr|\leq2
      \dott
    \end{eqn}
  }
  \argument{
    \lref{eq:selector-coordinate-variation};
    the monotone convergence theorem;
    the fact that for all $v=(v_1,\dots,v_n)\in\R^n$ it holds that
    $\norm{v}\leq\sum_{k=1}^n|v_k|$;
  }{
    that for every $\rho\in(0,R)$, $\eta\in(0,\infty)$
    and every interval $J\in\mathscr L_{\rho,\eta}$ it holds that
    \begin{eqn}
      \llabel{eq:selector-piece-variation}
      \int_J\norm{(\partial_e\bfX)(\rho,\eta,e)}\dd e
        \leq\sum_{k=1}^n\int_J
          \bigl|(\partial_e\bfX_k)(\rho,\eta,e)\bigr|\dd e
        \leq2n
      \dott
    \end{eqn}
  }
  \argument{
    \itref{lem:ambient-annular}{item:selector-pieces,item:selector-bound};
    \lref{eq:selector-piece-variation};
  }{
    that for all
    $\rho\in(0,R)$, $\eta\in(0,\infty)$, $a,b\in\R$
    with $a<b$ and $\{\rho\}\times(a,b)\subseteq\mathcal Z$
    it holds that
    \begin{eqn}[a]
      \llabel{eq:selector-variation}
      \int_a^b\norm{(\partial_e\bfX)(\rho,\eta,e)}\dd e
        =\sum_{\substack{J\in\mathscr L_{\rho,\eta}\\J\text{ open interval}}}
          \int_{J\cap(a,b)}\norm{(\partial_e\bfX)(\rho,\eta,e)}\dd e
        \leq2n(\#\mathscr L_{\rho,\eta})\leq2n(\#\mathscr F)
      \dott
    \end{eqn}
  }
  \argument{
    \lref{eq1};
    \lref{eq:small-selector};
    \lref{item:selector-regularity};
    the chain rule;
    the Cauchy--Schwarz inequality;
  }{
    that for every open interval $J\subseteq\R$ and every
    $\rho\in(0,R)$, $\eta\in(0,\infty)$, $e\in J$ with
    $J\in\mathscr L_{\rho,\eta}$ it holds that
    \begin{eqn}[a]
      \llabel{eq:selector-pointwise-estimate}
      1=\langle(\nabla_z f)(\rho,\bfX(\rho,\eta,e)),
        (\partial_e\bfX)(\rho,\eta,e)\rangle\leq H(\rho,\bfX(\rho,\eta,e))
          \norm{(\partial_e\bfX)(\rho,\eta,e)}
        \leq\eta\norm{(\partial_e\bfX)(\rho,\eta,e)}
      \dott
    \end{eqn}
  }
  \argument{
    \lref{eq:selector-pointwise-estimate};
    \itref{lem:ambient-annular}{item:selector-pieces,item:selector-bound};
  }{
    that for all $\rho\in(0,R)$, $\eta\in(0,\infty)$,
    $a,b\in\R$ with $a<b$ and
    $\{\rho\}\times(a,b)\subseteq\mathcal Z$ it holds that
    for Lebesgue-almost all $e\in(a,b)$ it holds that
    \begin{eqn}
      \llabel{eq:selector-ae-estimate}
      1\leq\eta\norm{(\partial_e\bfX)(\rho,\eta,e)}
      \dott
    \end{eqn}
  }
  \argument{
    \lref{eq:selector-ae-estimate};
    \lref{eq:selector-variation};
  }{
    that for all $\rho\in(0,R)$, $\eta\in(0,\infty)$,
    $a,b\in\R$ with $a<b$ and
    $\{\rho\}\times(a,b)\subseteq\mathcal Z$ it holds that
    \begin{eqn}
      \llabel{eq:zero-interval}
      b-a=\int_a^b1\dd e\leq\eta\int_a^b
        \norm{(\partial_e\bfX)(\rho,\eta,e)}\dd e
        \leq2n\eta(\#\mathscr F)
      \dott
    \end{eqn}
  }
  \argument{
    \lref{eq:zero-interval};
    the fact that $\mathcal Z$ is definable;
    the fact that every definable subset of $\R$ is a finite union
    of points and intervals;
  }{
    that for all $\rho\in(0,R)$ it holds that
    \begin{eqn}
      \llabel{eq:zero-fibers}
      \#\{e\in\R\colon(\rho,e)\in\mathcal Z\}<\infty
      \dott
    \end{eqn}
  }
  \argument{
    \lref{eq:infimum};
    \lref{eq:zero-set};
    the fact that $\mathfrak h$ is definable;
    \cite[Theorem~3.1]{Coste};
  }{
    that there exists a definable
    $\zeta=(\zeta(\rho,e))_{(\rho,e)\in\mathcal E}
      =(\zeta_1(\rho,e),\dots,\zeta_n(\rho,e))_{(\rho,e)\in\mathcal E}\colon\mathcal E\to\mathbb{V}$
    which satisfies for all $(\rho,e)\in\mathcal E$ that
    \begin{eqn}
      \llabel{eq2}
      (\rho,\zeta(\rho,e))\in\Omega
      \qquad\text{and}\qquad f(\rho,\zeta(\rho,e))=e
    \end{eqn}
    and which satisfies for all $(\rho,e)\in\mathcal E\backslash\mathcal Z$ that
    \begin{eqn}
      \llabel{eq:approximate-selector}
      H(\rho,\zeta(\rho,e))<2\mathfrak h(\rho,e)
      \dott
    \end{eqn}
  }
  \startnewargseq
  \argument{
    the fact that $\mathcal E$ is definable; the fact that $\zeta$ is definable;
    the fact that $\mathcal Z$ is definable;
    \cite[Theorem~6.7]{Coste};
    \cite[Theorem~6.6]{Coste};
  }{
    that there exists a finite cylindrical $C^2$ cell decomposition
    $\mathscr P$ of $\R^2$, which is compatible with
    $(0,R)\times\R$, $\mathcal E$, and $\mathcal Z$, which satisfies
    for all open $S\in\mathscr P$ with $S\subseteq\mathcal E$ that
    \begin{eqn}
      \llabel{eq2.1}
      \zeta|_S\in C^1(S,\R^n)
      \cfadd{definition:planar-cell-decomposition}\dott
    \end{eqn}
  }
  \startnewargseq
  \argument{
    \lref{eq2.1};
    the fact that for all open $S\in\mathscr P$, $k\in\{1,2,\dots,n\}$
    with $S\subseteq\mathcal E$ it holds that
    $\partial_e(\zeta_k|_S)$ is definable;
    \cite[Theorem~6.6]{Coste}
  }{
    that there exists a finite cylindrical $C^2$ refinement $\mathscr D$
    of $\mathscr P$ which satisfies for all open
    $S\in\mathscr D$ and all $k\in\{1,2,\dots,n\}$, $\xi,\eta\in S$
    with $S\subseteq\mathcal E$ that
    $\zeta|_S\in C^1(S,\R^n)$ and 
    \begin{eqn}
      \llabel{eq2.2}
      \operatorname{sign}\bigl(\partial_e(\zeta_k|_S)(\xi)\bigr)
        =\operatorname{sign}\bigl(\partial_e(\zeta_k|_S)(\eta)\bigr)
      \dott
    \end{eqn}
  }
  In the following let $\pi\colon\R^2\to\R$ satisfy for all
  $(s,e)\in\R^2$ that $\pi(s,e)=s$,
  let $\mathscr X$ satisfy
  \begin{eqn}
    \llabel{setting3.0a}
    \mathscr X=\{\pi(S)\colon S\in\mathscr D\},
  \end{eqn}
  and for every $S\in\mathscr D$, $\rho\in\R$ let
  $\mathbf{T}_{\rho,S}\subseteq\R$ satisfy
  \begin{eqn}
    \llabel{setting3.0b}
    \mathbf{T}_{\rho,S}=\{e\in\R\colon(\rho,e)\in S\}
    \dott
  \end{eqn}
  \startnewargseq
  \argument{
    \lref{setting3.0a};
    the fact that $\mathscr D$ is a finite cylindrical $C^2$ cell
    decomposition;
  }{
    that every element of $\mathscr X$ is a singleton or a nonempty
    open interval, that
    \begin{eqn}
      \llabel{eq2.3}
      \#\mathscr X<\infty
      \qquad\text{and}\qquad
      \midcup_{\Lambda\in\mathscr X}\Lambda=\R,
    \end{eqn}
    and that for all $\Lambda,\Upsilon\in\mathscr X$ with $\Lambda\neq\Upsilon$ it holds that
    \begin{eqn}
      \llabel{eq2.4}
      \Lambda\cap\Upsilon=\emptyset
      \dott
    \end{eqn}
  }
  \argument{
    \lref{setting3};
    \lref{setting3.0a};
    the fact that $\mathscr P$ is compatible with $\mathcal E$;
    the fact that $\mathscr D$ is a refinement of $\mathscr P$;
  }{
    that
    \begin{eqn}[a]
      \llabel{eq2.5}
      \mathcal E
        =\midcup_{S\in\mathscr D,\,S\subseteq\mathcal E}S
       \qquad\text{and}\qquad\pi(\mathcal E)
       =\midcup_{S\in\mathscr D,\,S\subseteq\mathcal E}\pi(S)
      \dott
    \end{eqn}
  }
  In the following let $\varrho\in\R$ satisfy
  \begin{eqn}
    \llabel{setting3.0c}
    \varrho=\frac12\min\bigl(
      \{R\}\cup\{a\in(0,R)\colon\{a\}\in\mathscr X\}
    \bigr)
    \dott
  \end{eqn}
  \startnewargseq
  \argument{
    \lref{eq2.3};
    \lref{setting3.0c};
    the assumption that $R>0$;
  }{
    that
    \begin{eqn}
      \llabel{eq2.6}
      0<\varrho<R
      \qquad\text{and}\qquad
      \{a\in(0,\varrho)\colon\{a\}\in\mathscr X\}=\emptyset
      \dott
    \end{eqn}
  }
  \argument{
    \lref{eq2.6};
    \lref{eq2.3};
    \lref{eq2.4};
    the fact that $\mathscr D$ is compatible with $(0,R)\times\R$;
    the fact that every element of $\mathscr X$ is a singleton or
    an open interval;
  }{
    that there exists an open interval $\mathcal O\in\mathscr X$
    which satisfies that
    \begin{eqn}
      \llabel{eq2.7}
      (0,\varrho)\subseteq\mathcal O\subseteq(0,R)
      \dott
    \end{eqn}
  }
  In the following let $\mathscr G,\mathscr W\subseteq\mathscr D$
  satisfy $ \mathscr G
      =\{S\in\mathscr D\colon S\subseteq\mathcal E,\
          \pi(S)=\mathcal O,\ S\text{ is not open}\}$ and
  \begin{eqn}[a]
    \llabel{setting3.0d}
    \mathscr W
      =\{S\in\mathscr D\colon S\subseteq\mathcal E,\
          \pi(S)=\mathcal O,\ S\text{ is open}\}
    \dott
  \end{eqn}
  \startnewargseq
  \argument{
    \lref{setting3};
    \lref{setting3.0b};
    \lref{eq2.7};
    \lref{setting3.0d};
    the fact that $\mathfrak g(W)\subseteq(-1,1)$;
  }{
    that for all $S\in\mathscr G\cup\mathscr W$, 
    $\rho\in\mathcal O$ it holds that
    \begin{eqn}[a]
      \llabel{eq2.8}
      \emptyset\neq \mathbf{T}_{\rho,S}
        &\subseteq\{f(\rho,z)\colon z\in\mathbb{V},\ (\rho,z)\in\Omega\}
        \subseteq(-1,1)
      \dott
    \end{eqn}
  }
  \argument{
    \lref{eq2.8};
    the fact that $\mathscr D$ is a cylindrical $C^2$ cell
    decomposition of $\mathbb R^2$;
    the fact that $\mathcal O$ is an open interval;
  }{
    that there exist definable 
    $\chi_S\in C^2(\mathcal O,\R)$, $S\in\mathscr G$, and definable
    $\theta_Q,\psi_Q\in C^2(\mathcal O,\R)$, $Q\in\mathscr W$,
    which satisfy for all $S\in\mathscr G$, $p\in\mathcal O$ that
    $S=\{(s,\chi_S(s))\colon s\in\mathcal O\}$ and
    $\mathbf{T}_{p,S}=\{\chi_S(p)\}$, and which satisfy for all
    $Q\in\mathscr W$, $p\in\mathcal O$ that
    \begin{gather}
      \llabel{eq2.9}
      \theta_Q(p)<\psi_Q(p),\qquad
      Q=\{(s,e)\in\mathcal O\times\R
        \colon\theta_Q(s)<e<\psi_Q(s)\},\\
      \text{and}\qquad
      \mathbf{T}_{p,Q}=(\theta_Q(p),\psi_Q(p))
      \dott
    \end{gather}
  }
  \startnewargseq
  \argument{
    \lref{setting3.0a};
    \lref{setting3.0b};
    \lref{eq2.4};
    \lref{eq2.5};
    \lref{eq2.7};
    \lref{setting3.0d};
  }{
    that for all $\rho\in(0,\varrho)$ it holds that
    \begin{equation}
      \llabel{eq2.11}
      \{e\in\R\colon(\rho,e)\in\mathcal E\}
        =\midcup_{S\in\mathscr G\cup\mathscr W}\mathbf{T}_{\rho,S}
    \end{equation}
    and that for all $S,T\in\mathscr G\cup\mathscr W$, $\rho\in(0,\varrho)$ with $S\neq T$ it holds that
    \begin{eqn}
      \llabel{eq2.12}
      \mathbf{T}_{\rho,S}\cap \mathbf{T}_{\rho,T}=\emptyset
      \dott
    \end{eqn}
  }
  \argument{ \lref{eq:zero-fibers};
    \lref{eq2.2};
    \lref{eq2.7};
    \lref{setting3.0d};
    \lref{eq2.9};
    the fact that $\mathscr D$ is compatible with $\mathcal Z$;
    the fact that for every $S\in\mathscr W$ it holds that
    $\zeta|_S\in C^1(S,\R^n)$;
  }{
    that for all $S\in\mathscr W$,
    $\rho\in(0,\varrho)$ it holds that
    \begin{eqn}
      \llabel{eq2.13}
      \bigl(\mathbf{T}_{\rho,S}\ni e\mapsto\zeta(\rho,e)\in\R^n\bigr)
        \in C^1(\mathbf{T}_{\rho,S},\R^n)
      \qquad\text{and}\qquad S\cap\mathcal Z=\emptyset,
    \end{eqn}
    and that for all $S\in\mathscr W$, $\rho\in(0,\varrho)$,
    $k\in\{1,2,\dots,n\}$, $u,v\in\mathbf{T}_{\rho,S}$ it holds that
    \begin{eqn}
      \llabel{eq2.14}
      \operatorname{sign}\bigl(\partial_e\zeta_k(\rho,u)\bigr)
        =\operatorname{sign}\bigl(\partial_e\zeta_k(\rho,v)\bigr)
      \dott
    \end{eqn}
  }
  In the following let $N\in\N$ satisfy
  \begin{eqn}
    \llabel{setting3.0e}
    N=\max\{1,\#\mathscr W\}
  \end{eqn}
  and for every $p\in(0,\varrho)$ let
  $\mathscr J_p\subseteq\{(a,b)\colon a,b\in\R,\ a<b\}$ and
  $A_p\subseteq\R$ satisfy
  \begin{eqn}
    \llabel{eq2.15}
    \mathscr J_p=\{(\theta_S(p),\psi_S(p))\colon S\in\mathscr W\}
    \qquad\text{and}\qquad
    A_p=\{\chi_S(p)\colon S\in\mathscr G\}
    \dott
  \end{eqn}
  \startnewargseq
  \argument{
    \lref{eq2.9};
    \lref{eq2.11};
    \lref{eq2.12};
    \lref{eq2.13};
    \lref{eq2.14};
    \lref{setting3.0e};
    \lref{eq2.15};
    the fact that $\mathscr G$ and $\mathscr W$ are finite;
    the mean value theorem;
  }{that
    \begin{enumerate}[label=(\Roman*)]
      \item \llabel{item1} for all $p\in(0,\varrho)$ it holds that
      \begin{eqn}
        \llabel{eq3}
        \{e\in\R\colon(p,e)\in\mathcal E\}
          =A_p\cup\bigl(\midcup_{J\in\mathscr J_p}J\bigr)
        \qquad\text{and}\qquad
        \#\mathscr J_p\leq N
      \end{eqn}
      and
      \item \llabel{item2} for all $p\in(0,\varrho)$, $J\in\mathscr J_p$,
      $k\in\{1,\dots,n\}$ it holds that
      $\bigl(J\ni e\mapsto\zeta(p,e)\in\R^n\bigr)\in C^1(J,\R^n)$ and that
      $J\ni e\mapsto\zeta_k(p,e)\in\R$ is non-decreasing or non-increasing\dott
    \end{enumerate}
  }
  In the following let $\mathfrak{C}=4nNK$\dott
  \startnewargseq

  We claim that for every $\rho\in(0,\varrho)$,
  every open $U\subseteq\R^n$,
  every $F\in C^1(U,\R)$,
  every interval $I\subseteq\R$,
  every locally absolutely continuous $y\colon I\to U\backslash\{0\}$,
  every Lebesgue measurable $\lambda\colon I\to[0,\infty)$,
  and every Lebesgue measurable $B\colon I\to\R^{n\times n}$
  with the property that $U\cap V\subseteq W$ and
  $F|_{U\cap V}=\mathfrak g|_{U\cap V}$,
  with the property that for Lebesgue-almost all $t\in I$ it holds
  that $y$ is differentiable at $t$, that $B(t)$ is symmetric and positive-definite,
  and that
  \begin{eqn}
    y'(t)=-\lambda(t)B(t)\grad F(y(t)),
  \end{eqn}
  and with the property that for Lebesgue-almost all
  $t\in\{s\in I\colon\norm{y(s)}<R\}$ it holds that
  \begin{eqn}
    \max_{v\in\R^n,\,\norm{v}=1}\langle B(t)v,v\rangle
      \leq K\min_{v\in\R^n,\,\norm{v}=1}\langle B(t)v,v\rangle,
  \end{eqn}
  we have that
  \begin{eqn}
    \llabel{eq3.1}
    \int_{\{t\in I\colon\rho/2\leq\norm{y(t)}\leq\rho\}}
      \norm{y'(t)}\dd t\leq\mathfrak{C}\rho
    \dott
  \end{eqn}
  \startnewargseq
  To prove this claim, let $\rho\in(0,\varrho)$,
  let $U\subseteq\R^n$ be open,
  let $F\in C^1(U,\R)$,
  let $I\subseteq\R$ be an interval,
  let $y\colon I\to U\backslash\{0\}$ be locally absolutely continuous,
  let $\lambda\colon I\to[0,\infty)$
  and $B\colon I\to\R^{n\times n}$ be Lebesgue measurable,
  assume that
  \begin{eqn}
    \llabel{setting3.1}
    U\cap V\subseteq W
    \qquad\text{and}\qquad
    F|_{U\cap V}=\mathfrak g|_{U\cap V},
  \end{eqn}
  assume that for Lebesgue-almost all $t\in I$ it holds that
  $y$ is differentiable at $t$, that $B(t)$ is symmetric and positive-definite,
  and that
  \begin{eqn}
    \llabel{setting3.2}
    y'(t)=-\lambda(t)B(t)\grad F(y(t)),
  \end{eqn}
  and assume that for Lebesgue-almost all
  $t\in\{s\in I\colon\norm{y(s)}<R\}$ it holds that
  \begin{eqn}
    \llabel{setting3.3}
    \max_{v\in\R^n,\,\norm{v}=1}\langle B(t)v,v\rangle
      \leq K\min_{v\in\R^n,\,\norm{v}=1}\langle B(t)v,v\rangle
    \dott
  \end{eqn}
  \startnewargseq
  In the following let $\mathfrak e\colon I\to\R$ satisfy
  $\mathfrak e=F\circ y$ and let
  $\mu,\nu\colon I\to\R$
  satisfy for all $t\in I$ that
  \begin{eqn}
    \llabel{setting4}
    \mu(t)=\min_{v\in\R^n,\,\norm{v}=1}
      \langle B(t)v,v\rangle
    \qquad\text{and}\qquad
    \nu(t)=\max_{v\in\R^n,\,\norm{v}=1}
      \langle B(t)v,v\rangle
    \dott
  \end{eqn}
  \cfload[.]
  \startnewargseq
  \argument{\lref{setting3.2};
    the assumption that $F\in C^1(U,\R)$;
    the assumption that $y$ is locally absolutely continuous;
    the chain rule;
    the fact that for Lebesgue-almost all $t\in I$ it holds that
    $\lambda(t)\geq0$; the fact that for Lebesgue-almost all $t\in I$ it holds that $B(t)$ is positive-definite;
  }{
    that $\mathfrak e$ is locally absolutely continuous and non-increasing,
    and that for Lebesgue-almost all $t\in I$ it holds that
    \begin{eqn}
      \llabel{eq4}
      -\mathfrak e'(t)
        =\lambda(t)\langle B(t)\grad F(y(t)),\grad F(y(t))\rangle
        \geq0
      \dott
    \end{eqn}
  }
  \argument{
    \lref{eq4};
    \lref{setting3.2};
    \lref{setting3.3};
    \lref{setting4};
  }{
    that for Lebesgue-almost all
    $t\in\{s\in I\colon\norm{y(s)}<R\}$ it holds that
    $0<\mu(t)\leq\nu(t)\leq K\mu(t)$, $\norm{y'(t)}
        \leq\lambda(t)\nu(t)\norm{\grad F(y(t))}$, and 
    \begin{eqn}[a]
      \llabel{eq5}
      -\mathfrak e'(t)
        \geq\lambda(t)\mu(t)\norm{\grad F(y(t))}^2
         \geq K^{-1}\norm{\grad F(y(t))}\norm{y'(t)}
      \dott
    \end{eqn}
  }
  In the following for every $J\in\mathscr J_\rho$ let
  $\mathbf{a}_J,\mathbf{b}_J\in\R$ satisfy
  $J=(\mathbf{a}_J,\mathbf{b}_J)$\dott
  \startnewargseq
  \argument{
    \lref{item2};
    the fact that for all $J\in\mathscr J_\rho$ it holds that
    $J$ is a nonempty, bounded open interval;
    the fact that $\zeta(\mathcal E)\subseteq\mathbb{V}\subseteq[-1,1]^n$;
  }{
    that for all $J\in\mathscr J_\rho$ there exists
    $v_J=(v_{J,1},\dots,v_{J,n})\in[-1,1]^n$ which satisfies
    \begin{eqn}
      \llabel{eq6}
      \limsup_{e\searrow\mathbf{a}_J}\norm{\zeta(\rho,e)-v_J}=0
      \dott
    \end{eqn}
  }
  In the following for every $J\in\mathscr J_\rho$ let
  $q_J\colon J\to\R$ satisfy for all $e\in J$ that
  \begin{eqn}
    \llabel{setting5}
    q_J(e)=\sum_{k=1}^n\abs{\zeta_k(\rho,e)-v_{J,k}}
    \dott
  \end{eqn}
  \startnewargseq
  \argument{
    \lref{item2};
    \lref{eq6};
    \lref{setting5};
    the fact that $\zeta(\mathcal E)\subseteq[-1,1]^n$;
    the fact that for all $J\in\mathscr J_\rho$ it holds that
    $v_J\in[-1,1]^n$;
  }{
    that for all $J\in\mathscr J_\rho$, $e\in J$ it holds that
    $q_J\in C^1(J,\R)$ is non-decreasing and that
    \begin{eqn}
      \llabel{eq7}
      0\leq q_J(e)\leq2n
      \qquad\text{and}\qquad
      (q_J)'(e)=\sum_{k=1}^n
        \abs{\partial_e\zeta_k(\rho,e)}
        \geq\norm{\partial_e\zeta(\rho,e)}
      \dott
    \end{eqn}
  }
  \argument{
    \lref{eq2.7};
    \lref{eq2.9};
    \lref{eq2.15};
  }{
    that for all $J\in\mathscr J_\rho$ there exists
    $S\in\mathscr W$ such that
    \begin{eqn}
      \llabel{eq:annular-interval-section}
      J=\bigl(\theta_S(\rho),\psi_S(\rho)\bigr)
      \qquad\text{and}\qquad
      \{\rho\}\times J\subseteq S
      \dott
    \end{eqn}
  }
  \argument{
    \lref{eq:annular-interval-section};
    \lref{setting3.0d};
    \lref{eq2.13};
  }{
    that for all $S\in\mathscr W$ it holds that
    $S\subseteq\mathcal E\backslash\mathcal Z$ and that for all
    $J\in\mathscr J_\rho$ it holds that
    \begin{eqn}
      \llabel{eq:annular-selector-domain}
      \{\rho\}\times J\subseteq\mathcal E\backslash\mathcal Z
      \dott
    \end{eqn}
  }
  \argument{
    \lref{eq:annular-selector-domain};
    \lref{eq1};
    \lref{eq2};
    \lref{item2};
    the chain rule;
    the Cauchy--Schwarz inequality;
  }{
    that for all $J\in\mathscr J_\rho$, $e\in J$ it holds that
    \begin{eqn}[a]
      \llabel{eq:annular-selector-chain-rule}
      1
        =\langle(\nabla_z f)(\rho,\zeta(\rho,e)),
                  (\partial_e\zeta)(\rho,e)\rangle
        \leq\norm{(\nabla_z f)(\rho,\zeta(\rho,e))}
                  \norm{(\partial_e\zeta)(\rho,e)}
        =H(\rho,\zeta(\rho,e))
                  \norm{(\partial_e\zeta)(\rho,e)}
      \dott
    \end{eqn}
  }
  \argument{
    \lref{eq:annular-selector-chain-rule};
    \lref{eq:approximate-selector};
    \lref{eq7};\lref{eq:annular-selector-domain};
  }{
    that for all $J\in\mathscr J_\rho$, $e\in J$ it holds that
    \begin{eqn}
      \llabel{eq:annular-selector-slope}
      1\leq H(\rho,\zeta(\rho,e))
              \norm{(\partial_e\zeta)(\rho,e)}
        \leq2\mathfrak h(\rho,e)(q_J)'(e)
      \dott
    \end{eqn}
  }
  \argument{
    \lref{setting1};
    \lref{setting2};
    \lref{eq:infimum};
    \lref{setting3.1};
    the fact that $U\cap V$ is open;
  }{
    that for all $J\in\mathscr J_\rho$, $e\in J$,
    $z\in\mathbb{V}$ with $\rho z\in U\cap V$ and
    $f(\rho,z)=e$ it holds that $(\rho,z)\in\Omega$ and that
    \begin{eqn}
      \llabel{eq:annular-gradient-comparison}
      \mathfrak h(\rho,e)\leq H(\rho,z)
        =\rho\norm{\grad\mathfrak g(\rho z)}
        =\rho\norm{\grad F(\rho z)}
      \dott
    \end{eqn}
  }
  \argument{\lref{eq7};
    \lref{eq:annular-selector-slope};
    \lref{eq:annular-gradient-comparison};
  }{
    that for all $J\in\mathscr J_\rho$, $e\in J$,
    $z\in\mathbb{V}$ with $\rho z\in U\cap V$ and
    $f(\rho,z)=e$ it holds that
    \begin{eqn}
      \llabel{eq8}
      1\leq2\mathfrak h(\rho,e)(q_J)'(e)
        \leq2\rho\norm{\grad F(\rho z)}(q_J)'(e)
      \dott
    \end{eqn}
  }
  In the following for every $J\in\mathscr J_\rho$,
  $\delta\in(0,1/2)$ let $\alpha_{J,\delta},\beta_{J,\delta}\in\R$
  and $Q_{J,\delta}\colon\R\to\R$ satisfy for all $e\in\R$ that  $\alpha_{J,\delta}=\mathbf{a}_J+\delta(\mathbf{b}_J-\mathbf{a}_J)$,
  \begin{eqn}[a]
    \llabel{setting6}
    \beta_{J,\delta}=\mathbf{b}_J-\delta(\mathbf{b}_J-\mathbf{a}_J),\quad\text{and}\quad
    Q_{J,\delta}(e)
      =q_J\bigl(\max\{\alpha_{J,\delta},
                       \min\{e,\beta_{J,\delta}\}\}\bigr)
  \end{eqn}
  and let $\mathcal T_{J,\delta}\subseteq I$ satisfy
  \begin{eqn}
    \llabel{setting7}
    \mathcal T_{J,\delta}
      =\{t\in I\colon\rho/2\leq\norm{y(t)}\leq\rho,\
          \alpha_{J,\delta}<\mathfrak e(t)<\beta_{J,\delta}\}
    \dott
  \end{eqn}
  \startnewargseq
  \argument{
    \lref{eq7};
    \lref{setting6};
    the fact that for all $J\in\mathscr J_\rho$,
    $\delta\in(0,1/2)$ it holds that
    $\mathbf{a}_J<\alpha_{J,\delta}<\beta_{J,\delta}<\mathbf{b}_J$;
    the fact that for all $J\in\mathscr J_\rho$ it holds that
    $q_J\in C^1(J,\R)$;
  }{
    that for all $J\in\mathscr J_\rho$, $\delta\in(0,1/2)$
    it holds that $Q_{J,\delta}$ is globally Lipschitz and
    non-decreasing and that
    \begin{eqn}
      \llabel{eq9}
      Q_{J,\delta}(\R)\subseteq[0,2n]
      \qquad\text{and}\qquad
      Q_{J,\delta}|_{(\alpha_{J,\delta},\beta_{J,\delta})}
        =q_J|_{(\alpha_{J,\delta},\beta_{J,\delta})}
      \dott
    \end{eqn}
  }
  \argument{\lref{eq9};
    the fact that $\mathfrak e$ is non-increasing;
    the fact that $\mathfrak e$ is locally absolutely continuous;
  }{
    that for all $J\in\mathscr J_\rho$, $\delta\in(0,1/2)$ it holds that
    $Q_{J,\delta}\circ\mathfrak e$ is locally absolutely continuous
    and that for all $s,t\in I$ with $s<t$ it holds that
    \begin{eqn}
      \llabel{eq9.1}
      Q_{J,\delta}(\mathfrak e(s))
        \geq Q_{J,\delta}(\mathfrak e(t))
      \dott
    \end{eqn}
  }
  \argument{
    \lref{setting1};
    \lref{setting2};
    \lref{setting3.1};
    the fact that $0<\rho<R$;
    the fact that $\mathfrak e=F\circ y$;
  }{
    that for all $t\in I$ with
    $\rho/2\leq\norm{y(t)}\leq\rho$ it holds that
    \begin{eqn}[a]
      \llabel{eq9.2}
      &0<\rho/2\leq\norm{y(t)}\leq\rho<R,
        \qquad y(t)\in U\cap V\subseteq W,\qquad
        y(t)/\rho\in\mathbb{V},\\
      &\qquad\text{and}\qquad
        f(\rho,y(t)/\rho)=\mathfrak g(y(t))=F(y(t))=\mathfrak e(t)
      \dott
    \end{eqn}
  }
  \argument{
    \lref{eq9.2};
    \lref{eq7};
    \lref{eq8} (applied with $z\curvearrowleft y(t)/\rho$, $e\curvearrowleft\mathfrak e(t)$ in the notation of \lref{eq8});
    \lref{setting7};
    the fact that for all $J\in\mathscr J_\rho$, $\delta\in(0,1/2)$
    it holds that $(\alpha_{J,\delta},\beta_{J,\delta})\subseteq J$;
    the fact that $\rho>0$;
  }{
    that for all $J\in\mathscr J_\rho$, $\delta\in(0,1/2)$,
    $t\in\mathcal T_{J,\delta}$ it holds that
    \begin{eqn}
      \llabel{eq9.3}
      \mathfrak e(t)\in(\alpha_{J,\delta},\beta_{J,\delta})
        \subseteq J,
        \qquad
        (q_J)'(\mathfrak e(t))\geq0,
        \qquad\text{and}\qquad
        (q_J)'(\mathfrak e(t))\norm{\grad F(y(t))}
          \geq\frac1{2\rho}
      \dott
    \end{eqn}
  }
  \argument{
    \lref{eq9.3};
    \lref{eq9};
    the fact that for all $J\in\mathscr J_\rho$ it holds that
    $q_J\in C^1(J,\R)$;
    the fact that $\mathfrak e$ is locally absolutely continuous;
    the chain rule; the fundamental theorem of calculus for absolutely continuous functions;
  }{
    that for all $J\in\mathscr J_\rho$, $\delta\in(0,1/2)$ it holds that
    for Lebesgue-almost all $t\in\mathcal T_{J,\delta}$ it holds that
    $Q_{J,\delta}$ is differentiable at $\mathfrak e(t)$,
    that $Q_{J,\delta}\circ\mathfrak e$ is differentiable at $t$,
    and that
    \begin{eqn}
      \llabel{eq9.4}
      -\frac{\dd}{\dd t}Q_{J,\delta}(\mathfrak e(t))
        =-(Q_{J,\delta})'(\mathfrak e(t))\mathfrak e'(t)
        =(q_J)'(\mathfrak e(t))\bigl(-\mathfrak e'(t)\bigr)
      \dott
    \end{eqn}
  }
  \argument{
    \lref{eq9.4};
    \lref{eq5};
    \lref{eq9.2};
    \lref{eq9.3};
  }{
    that for all $J\in\mathscr J_\rho$, $\delta\in(0,1/2)$ it holds that
    for Lebesgue-almost all $t\in\mathcal T_{J,\delta}$ it holds that
    \begin{eqn}[a]
      \llabel{eq10}
      -\frac{\dd}{\dd t}Q_{J,\delta}(\mathfrak e(t))
        =(q_J)'(\mathfrak e(t))\bigl(-(\mathfrak e)'(t)\bigr)\geq K^{-1}(q_J)'(\mathfrak e(t))
          \norm{\grad F(y(t))}\norm{y'(t)}\geq\frac{\norm{y'(t)}}{2K\rho}
      \dott
    \end{eqn}
  }
  \argument{
    \lref{eq9.1};
    the fact that for all $J\in\mathscr J_\rho$, $\delta\in(0,1/2)$
    it holds that $Q_{J,\delta}\circ\mathfrak e$ is locally absolutely continuous;
    the fundamental theorem of calculus for absolutely continuous functions;
  }{
    that for all $J\in\mathscr J_\rho$, $\delta\in(0,1/2)$,
    $s,t\in I$ with $s<t$ it holds that
    for Lebesgue-almost all $u\in[s,t]$ it holds that
    $Q_{J,\delta}\circ\mathfrak e$ is differentiable
    at $u$ and that
    \begin{eqn}
      \llabel{eq10.1}
      -\frac{\dd}{\dd u}Q_{J,\delta}(\mathfrak e(u))\geq0
      \dott
    \end{eqn}
  }
  \argument{
    \lref{eq10.1};
    \lref{item1};
    \lref{eq9};
    \lref{eq9.1};
    \lref{eq10};
    the fact that for all $J\in\mathscr J_\rho$, $\delta\in(0,1/2)$
    it holds that $Q_{J,\delta}\circ\mathfrak e$ is locally absolutely continuous;
      the fact that $\mathfrak{C}=4nNK$;
    the fundamental theorem of calculus for absolutely continuous functions;
  }{
    that for all $s,t\in I$, $\delta\in(0,1/2)$ with $s<t$ 
    it holds that
    \begin{eqn}[a]
      \llabel{eq11}
      &\int_{[s,t]\cap(\midcup_{J\in\mathscr J_\rho}\mathcal T_{J,\delta})}
        \norm{y'(u)}\dd u
        \leq\sum_{J\in\mathscr J_\rho}
          \int_{[s,t]\cap\mathcal T_{J,\delta}}
            \norm{y'(u)}\dd u\\
      &\leq2K\rho\sum_{J\in\mathscr J_\rho}
          \int_{[s,t]\cap\mathcal T_{J,\delta}}
            -\frac{\dd}{\dd u}Q_{J,\delta}(\mathfrak e(u))\dd u
        \leq2K\rho\sum_{J\in\mathscr J_\rho}
          \int_s^t-\frac{\dd}{\dd u}
            Q_{J,\delta}(\mathfrak e(u))\dd u\\
      &=2K\rho\sum_{J\in\mathscr J_\rho}
          \bigl(Q_{J,\delta}(\mathfrak e(s))
                -Q_{J,\delta}(\mathfrak e(t))\bigr)
        \leq2K\rho\sum_{J\in\mathscr J_\rho}
          Q_{J,\delta}(\mathfrak e(s))\\
      &\leq2K\rho\sum_{J\in\mathscr J_\rho}2n
        =4nK\rho(\#\mathscr J_\rho)
        \leq4nK\rho N=\mathfrak{C}\rho
      \dott
    \end{eqn}
  }
  \argument{
    \lref{eq11};
    \lref{setting6};
    \lref{setting7};
    the monotone convergence theorem;
  }{
    that for all $s,t\in I$ with $s<t$ it holds that
    \begin{eqn}
      \llabel{eq12}
      \int_{\{u\in[s,t]\colon\rho/2\leq\norm{y(u)}\leq\rho,\
                         \mathfrak e(u)\in\bigcup_{J\in\mathscr J_\rho}J\}}
        \norm{y'(u)}\dd u\leq \mathfrak{C}\rho
      \dott
    \end{eqn}
  }
  \argument{
    the fact that $A_\rho$ is finite;
    the fact that $\mathfrak e$ is locally absolutely continuous;
    \cref{lem:AC-level-set}
  }{
    that for all $s,t\in I$ with $s<t$ it holds that
    for Lebesgue-almost all
    $u\in\{v\in[s,t]\colon\mathfrak e(v)\in A_\rho\}$ it holds that
    \begin{eqn}
      \llabel{eq12.1}
      \mathfrak e'(u)=0
      \dott
    \end{eqn}
  }
  \argument{
    \lref{eq12.1};
    \lref{eq4};
    the fact that for Lebesgue-almost all $u\in I$ it holds that
    $B(u)$ is positive-definite; the fact that for Lebesgue-almost all $u\in I$ it holds that $\lambda(u)\geq0$;
  }{
    that for all $s,t\in I$ with $s<t$ it holds that
    for Lebesgue-almost all
    $u\in\{v\in[s,t]\colon\mathfrak e(v)\in A_\rho\}$ it holds that
    \begin{eqn}
      \llabel{eq12.2}
      \lambda(u)\langle B(u)\grad F(y(u)),\grad F(y(u))\rangle=0
      \qquad\text{and}\qquad
      \lambda(u)\grad F(y(u))=0
      \dott
    \end{eqn}
  }
  \argument{
    \lref{eq12.2};
    \lref{setting3.2};
  }{
    that for all $s,t\in I$ with $s<t$ it holds that
    for Lebesgue-almost all
    $u\in\{v\in[s,t]\colon\mathfrak e(v)\in A_\rho\}$ it holds that
    \begin{eqn}
      \llabel{eq13}
      y'(u)=-B(u)\bigl(\lambda(u)\grad F(y(u))\bigr)=0
      \dott
    \end{eqn}
  }
  \argument{
    \lref{eq13};
    the fact that $A_\rho$ is finite;
    the fact that for all $s,t\in I$ with $s<t$ it holds that
    $\mathfrak e|_{[s,t]}$ is continuous;
  }{
    that for all $s,t\in I$ with $s<t$ it holds that
    \begin{eqn}
      \llabel{eq13.1}
      \int_{\{u\in[s,t]\colon\mathfrak e(u)\in A_\rho\}}
        \norm{y'(u)}\dd u=0
      \dott
    \end{eqn}
  }
  \argument{\lref{item1};
    \lref{setting3};
    \lref{eq9.2};
    the fact that $\rho\in(0,\varrho)$;
  }{
    that for all $u\in I$ with
    $\rho/2\leq\norm{y(u)}\leq\rho$ it holds that
    \begin{eqn}
      \llabel{eq13.2}
      (\rho,\mathfrak e(u))\in\mathcal E
      \qquad\text{and}\qquad
      \mathfrak e(u)\in A_\rho\cup\bigl(\midcup_{J\in\mathscr J_\rho}J\bigr)
      \dott
    \end{eqn}
  }
  \argument{
    \lref{eq13.2};
    \lref{eq12};
    \lref{eq13.1};
  }{
    that for all $s,t\in I$ with $s<t$ it holds that
    \begin{eqn}[a]
      \llabel{eq14.0}
      &\int_{\{u\in[s,t]\colon\rho/2\leq\norm{y(u)}\leq\rho\}}
        \norm{y'(u)}\dd u\\
        &
        \leq\int_{\{u\in[s,t]\colon\mathfrak e(u)\in A_\rho\}}
          \norm{y'(u)}\dd u+
      \int_{\{u\in[s,t]\colon\rho/2\leq\norm{y(u)}\leq\rho,\
                         \mathfrak e(u)\in\bigcup_{J\in\mathscr J_\rho}J\}}
          \norm{y'(u)}\dd u\\
          &
        \leq0+\mathfrak{C}\rho=\mathfrak{C}\rho
      \dott
    \end{eqn}
  }
  \argument{\lref{eq14.0};the monotone convergence theorem}{that 
  \begin{equation}\llabel{eq14}
    \int_{\{u\in I\colon\rho/2\leq\norm{y(u)}\leq\rho\}}
        \norm{y'(u)}\dd u\leq \mathfrak C\rho.
  \end{equation}
  }
  \argument{
    \lref{eq14};
    the arbitrariness of $\rho,U,F,I,y,\lambda,B$ subject to
    the stated assumptions;
  }[verbs=ep]{
    \lref{eq3.1}\dott
  }

  \startnewargseq
  \argument{
    \lref{eq3.1} (applied for every $j\in\N_0$ with $\rho\curvearrowleft2^{-j}\rho$);
  }{
    that for every $\rho\in(0,\varrho)$,
    every open $U\subseteq\R^n$,
    every $F\in C^1(U,\R)$,
    every interval $I\subseteq\R$,
    every locally absolutely continuous $y\colon I\to U\backslash\{0\}$,
    every Lebesgue measurable $\lambda\colon I\to[0,\infty)$,
    and every Lebesgue measurable $B\colon I\to\R^{n\times n}$
    with the property that $U\cap V\subseteq W$ and
    $F|_{U\cap V}=\mathfrak g|_{U\cap V}$,
    with the property that for Lebesgue-almost all $t\in I$ it holds
    that $y$ is differentiable at $t$, that $B(t)$ is symmetric and positive-definite,
    and that
    \[
      y'(t)=-\lambda(t)B(t)\grad F(y(t)),
    \]
    and with the property that for Lebesgue-almost all
    $t\in\{s\in I\colon\norm{y(s)}<R\}$ it holds that
    \[
      \max_{v\in\R^n,\,\norm{v}=1}\langle B(t)v,v\rangle
        \leq K\min_{v\in\R^n,\,\norm{v}=1}\langle B(t)v,v\rangle,
    \]
    we have that
    \begin{eqn}[a]
      \llabel{eq14.1}
      \int_{\{t\in I\colon\norm{y(t)}\leq\rho\}}\norm{y'(t)}\dd t
        =\sum_{j=0}^{\infty}
          \int_{\{t\in I\colon2^{-j-1}\rho<\norm{y(t)}\leq2^{-j}\rho\}}
            \norm{y'(t)}\dd t
        \leq\sum_{j=0}^{\infty}\mathfrak{C}2^{-j}\rho=2\mathfrak{C}\rho
      \dott
    \end{eqn}
  }
  \startnewargseq
  To prove \lref{concl1}, let $\rho\in(0,\varrho)$,
  let $U\subseteq\R^n$ be open,
  let $F\in C^1(U,\R)$,
  let $I\subseteq\R$ be an interval,
  let $y\colon I\to U\backslash\{0\}$ be locally absolutely continuous,
  let $\lambda\colon I\to[0,\infty)$
  and $B\colon I\to\R^{n\times n}$ be Lebesgue measurable,
  assume that
  \begin{eqn}
    \llabel{eq:original-agreement}
    U\cap V\subseteq W
    \qquad\text{and}\qquad
    F|_{U\cap V}=g|_{U\cap V},
  \end{eqn}
  assume that for Lebesgue-almost all $t\in I$ it holds that
  $y$ is differentiable at $t$, that $B(t)$ is symmetric and positive-definite,
  and that
  \begin{eqn}
    \llabel{eq:original-flow}
    y'(t)=-\lambda(t)B(t)\grad F(y(t)),
  \end{eqn}
  and assume that for Lebesgue-almost all
  $t\in\{s\in I\colon\norm{y(s)}<R\}$ it holds that
  \begin{eqn}
    \llabel{eq:original-matrix-bound}
    \max_{v\in\R^n,\,\norm{v}=1}\langle B(t)v,v\rangle
      \leq K\min_{v\in\R^n,\,\norm{v}=1}\langle B(t)v,v\rangle
    \dott
  \end{eqn}
  In the following let $G\colon U\to\R$ satisfy for all $z\in U$ that
  \begin{eqn}
    \llabel{eq:normalized-extension}
    G(z)=\varphi(F(z))
  \end{eqn}
  and let $\kappa\colon I\to[0,\infty)$ satisfy for all $t\in I$ that
  \begin{eqn}
    \llabel{eq:normalization-coefficient}
    \kappa(t)=\frac{\lambda(t)}{\varphi'(F(y(t)))}
    \dott
  \end{eqn}
  \startnewargseq
  \argument{
    \lref{normalization-map};
    \lref{normalization-derivative};
    \lref{eq:original-agreement};
    \lref{eq:normalized-extension};
    the assumption that $F\in C^1(U,\R)$;
    the chain rule;
  }{
    that $G\in C^1(U,\R)$, that $G(U)\subseteq(-1,1)$, that
    $G|_{U\cap V}=\mathfrak g|_{U\cap V}$, and that
    for all $z\in U$ it holds that
    \begin{eqn}
      \llabel{eq:normalized-gradient}
      \grad G(z)=\varphi'(F(z))\grad F(z)
      \dott
    \end{eqn}
  }
  \argument{
    \lref{normalization-derivative};
    \lref{eq:normalization-coefficient};
    the assumption that $F\in C^1(U,\R)$;
    the assumption that $y$ is locally absolutely continuous;
    the assumption that $\lambda$ is Lebesgue measurable;
  }{
    that $I\ni t\mapsto\varphi'(F(y(t)))\in(0,\infty)$ is continuous,
    that $\kappa$ is Lebesgue measurable, and that for all $t\in I$
    it holds that
    \begin{eqn}
      \llabel{eq:normalization-coefficient-properties}
      \varphi'(F(y(t)))>0
      \qquad\text{and}\qquad \kappa(t)\geq0
      \dott
    \end{eqn}
  }
  \argument{
    \lref{eq:original-flow};
    \lref{eq:normalized-gradient};
    \lref{eq:normalization-coefficient-properties};    \lref{eq:normalization-coefficient};
  }{
    that for Lebesgue-almost all $t\in I$ it holds that
    \begin{eqn}[a]
      \llabel{eq:normalized-flow}
      y'(t)&=-\lambda(t)B(t)\grad F(y(t))
        =-\kappa(t)\varphi'(F(y(t)))B(t)\grad F(y(t))
         =-\kappa(t)B(t)\grad G(y(t))
      \dott
    \end{eqn}
  }
  \argument{
    \lref{eq3.1} (applied with $F\curvearrowleft G$,
    $\lambda\curvearrowleft\kappa$);
    \lref{eq:original-agreement};
    \lref{eq:normalized-flow};
    \lref{eq:original-flow};
    \lref{eq:original-matrix-bound};
    \lref{eq:normalized-gradient};
    \lref{eq:normalization-coefficient-properties};
  }{
    that
    \begin{eqn}
      \llabel{eq:original-annular-estimate}
      \int_{\{t\in I\colon\rho/2\leq\norm{y(t)}\leq\rho\}}
        \norm{y'(t)}\dd t\leq\mathfrak C\rho
      \dott
    \end{eqn}
  }
  \argument{
    \lref{eq14.1} (applied with $F\curvearrowleft G$,
    $\lambda\curvearrowleft\kappa$);
    \lref{eq:original-agreement};
    \lref{eq:original-flow};
    \lref{eq:original-matrix-bound};\lref{eq:normalized-gradient};
    \lref{eq:normalization-coefficient-properties};\lref{eq:normalized-flow};
  }{
    that
    \begin{eqn}
      \llabel{eq:original-ball-estimate}
      \int_{\{t\in I\colon\norm{y(t)}\leq\rho\}}
        \norm{y'(t)}\dd t\leq2\mathfrak C\rho
      \dott
    \end{eqn}
  }
  \argument{
    \lref{eq:original-annular-estimate};
    \lref{eq:original-ball-estimate};
    the arbitrariness of $\rho,U,F,I,y,\lambda,B$ subject to
    the stated assumptions;
  }[verbs=ep]{
    \lref{concl1}\dott
  }
\end{aproof}

\subsection{Uniform upper bounds for spherical lengths under integrable radial-tangential mixing}

The spherical projection $\omega(t)$
(see \cref{prop:ambient-state.ass1} in \cref{prop:ambient-state} and
\cref{prop:ambient-mixing.ass:polar} in \cref{prop:ambient-mixing} below)
points purely in the radial direction and its tangential projection
$\omega(t)-\langle\omega(t),\omega(t)\rangle\omega(t)=\omega(t)-\omega(t)=0$
vanishes.
Applying the (inverse) ambient metric matrix $B(t)$ to $\omega(t)$
need not preserve this property, since the tangential projection
$B(t)\omega(t)-\langle B(t)\omega(t),\omega(t)\rangle\omega(t)$
of $B(t)\omega(t)$
(see \cref{prop:ambient-state.ass3} in \cref{prop:ambient-state} and
\cref{eq:ambient-mixing} in \cref{prop:ambient-mixing} below)
need not vanish.
In \cref{prop:ambient-mixing} below we establish upper bounds for the
spherical length of all \GF\ trajectories close to the origin under, among
other hypotheses, the assumption that the ratio
between the norm
$\norm{B(t)\omega(t)-\langle B(t)\omega(t),\omega(t)\rangle\omega(t)}$
of the tangential projection of $B(t)\omega(t)$ and the norm
$\norm{\langle B(t)\omega(t),\omega(t)\rangle\omega(t)}
  =\langle B(t)\omega(t),\omega(t)\rangle$
of the radial projection of $B(t)\omega(t)$ is controlled by an integrable
radial bound in the sense of \cref{eq:ambient-mixing} in
\cref{prop:ambient-mixing} (see also \cref{prop:ambient-state.ass3} in
\cref{prop:ambient-state}).

\begin{athm}{proposition}{prop:ambient-state}
  Let $n\in\N\backslash\{1\}$, $K,R\in(0,\infty)$,
  let $\mathbb S=\{z\in\R^n\colon\norm{z}=1\}$,
  let $W\subseteq\{z\in\R^n\colon0<\norm{z}<R\}$ be definable
  and open, let
  $F\in C^1(W,\R)$ be definable\cfadd{definition:definable-function},
   let $\mathbf{k}\colon(0,R)\to[0,\infty)$ be locally bounded and
  definable, and assume
  \begin{eqn}
    \llabel{ass0}
    \int_0^R \mathbf{k}(s)\dd s<\infty
    \dott
  \end{eqn}
  \cfload[.]
  Then there exist $\mathfrak{R}\in(0,R)$, a bounded
  $Q\in C(\{(s,x)\in(0,\mathfrak{R})\times\mathbb S\colon sx\in W\},[0,\infty))$,
  and a definable $m\in C((0,\mathfrak{R}),[0,\infty))$
  with $\int_0^{\mathfrak{R}}m(s)\dd s<\infty$
  such that for every interval $I\subseteq\R$,
  every locally absolutely continuous $y\colon I\to W$,
  every Lebesgue measurable $\lambda\colon I\to[0,\infty)$,
  every Lebesgue measurable $B\colon I\to\R^{n\times n}$,
  every $r\colon I\to(0,\mathfrak{R})$,
  every $\omega\colon I\to\mathbb S$,
  and every $q\colon I\to\R$
  with the property that for all $t\in I$ it holds that
  \begin{eqn}
    \llabel{ass1}
    r(t)=\norm{y(t)},\qquad
    \omega(t)=y(t)/r(t),\qquad\text{and}\qquad
    q(t)=Q(r(t),\omega(t))
  \end{eqn}
  and with the property that for Lebesgue-almost all $t\in I$ it holds
  that $y$ is differentiable at $t$, that $B(t)$ is symmetric and
  positive-definite, and that
  \begin{equation}
    \llabel{ass2}
    y'(t)=-\lambda(t)B(t)\grad F(y(t)),
    \qquad
    \max\nolimits_{z\in\R^n,\,\norm{z}=1}\langle B(t)z,z\rangle
      \leq K\min\nolimits_{z\in\R^n,\,\norm{z}=1}\langle B(t)z,z\rangle,
  \end{equation}
  \begin{equation}
    \llabel{ass3}
    \text{and}\qquad
    [\langle B(t)\omega(t),\omega(t)\rangle]^{-1}
      \norm{B(t)\omega(t)
        -\langle B(t)\omega(t),\omega(t)\rangle\omega(t)}
      \leq r(t)\mathbf{k}(r(t))
  \end{equation}
  we have that
  \begin{enumerate}[label=(\roman*)]
    \item \llabel{item1} it holds that $q$ is locally absolutely continuous,
    \item \llabel{item2} it holds for Lebesgue-almost all $t\in I$ that
    $q$, $r$, and $\omega$ are differentiable at $t$, and
    \item \llabel{item3} it holds for Lebesgue-almost all $t\in I$ that
    $q'(t)\leq-\norm{\omega'(t)}+m(r(t))\abs{r'(t)}$\dott
  \end{enumerate}
  \cfout[.]
\end{athm}
\begin{aproof}
  Throughout this proof let
  $\Omega=\{(s,x)\in(0,R)\times\mathbb S\colon sx\in W\}$,
  let $f\colon W\to(-1,1)$ satisfy for all $z\in W$ that
  \begin{eqn}
    \llabel{bounded-transform}
    f(z)=\frac{F(z)}{\sqrt{1+(F(z))^2}},
  \end{eqn}
  and let $h=(h(s,x))_{(s,x)\in\Omega}\colon\Omega\to\R$,
  $u\colon\Omega\to\R^n$, and
  $G\colon\Omega\to[0,\infty)$ satisfy
  for all $(s,x)\in\Omega$ that
  \begin{eqn}
    \llabel{setting1}
    h(s,x)=f(sx),\qquad
    u(s,x)=\nabla_{\mathbb S}h(s,x)\cfadd{definition:angular-gradient},\qquad\text{and}\qquad
    G(s,x)=\norm{u(s,x)}
    \dott
  \end{eqn}
  \cfload[.]
  \startnewargseq
  \argument{
    \lref{bounded-transform};
    the assumption that $F\in C^1(W,\R)$ is definable;
    the chain rule;
  }{
    that $f\in C^1(W,(-1,1))$ is definable and that for all $z\in W$
    it holds that
    \begin{eqn}
      \llabel{bounded-potential}
      -1<f(z)<1
      \dott
    \end{eqn}
  }
  \argument{
    \lref{bounded-transform};
    the assumption that $F\in C^1(W,\R)$;
    the chain rule;
  }{
    that for all $z\in W$ it holds that
    \begin{eqn}
      \llabel{bounded-gradient}
      \grad f(z)=\frac{\grad F(z)}{(1+(F(z))^2)^{3/2}}
      \dott
    \end{eqn}
  }
  \argument{
    \lref{setting1};
    the assumption that $W$ is open and definable;
    the fact that $f\in C^1(W,(-1,1))$ is definable;
    the chain rule;
  }{
    that $\Omega$ is definable and open in $(0,R)\times\mathbb S$,
    that $h\in\COne{\Omega}{\R}$ is definable, and that for all
    $(s,x)\in\Omega$ it holds that
    \begin{eqn}
      \llabel{eq:polar-radial-derivative}
      \partial_s h(s,x)
        =\langle\grad f(sx),x\rangle
      \dott
    \end{eqn}
  }
  \argument{
    \lref{setting1};
    the fact that $f\in C^1(W,(-1,1))$;
    the chain rule;
  }{
    that for all $(s,x)\in\Omega$, $\delta\in(0,\infty)$,
    $\gamma\in C^1((-\delta,\delta),\mathbb S)$ with
    $\gamma(0)=x$ and $\{s\}\times\gamma((-\delta,\delta))\subseteq\Omega$
    it holds that
    $\bigl((-\delta,\delta)\ni t\mapsto h(s,\gamma(t))\in\R\bigr)
      \in C^1((-\delta,\delta),\R)$ and that
    \begin{eqn}
      \llabel{eq:angular-curve-derivative}
      0=\frac{\dd}{\dd t}\norm{\gamma(t)}^2\Big|_{t=0}
        =2\langle x,\gamma'(0)\rangle
      \qquad\text{and}\qquad
      \frac{\dd}{\dd t}h(s,\gamma(t))\Big|_{t=0}
        =s\langle\grad f(sx),\gamma'(0)\rangle
      \dott
    \end{eqn}
  }
  \argument{
    \lref{eq:angular-curve-derivative};
    the fact that for all $x\in\mathbb S$ it holds that
    $T_x\mathbb S=\{v\in\R^n\colon\langle v,x\rangle=0\}$;
  }{
    that for all $(s,x)\in\Omega$ it holds that
    \begin{gather}
      \bigl\langle s\bigl(\grad f(sx)
        -\langle\grad f(sx),x\rangle x\bigr),x\bigr\rangle
        =s\langle\grad f(sx),x\rangle(1-\norm{x}^2)=0\notag\\
      \llabel{eq:angular-projection-tangent}
      \text{and}\qquad
      s\bigl(\grad f(sx)-\langle\grad f(sx),x\rangle x\bigr)
        \in T_x\mathbb S
      \dott
    \end{gather}
    and that for all $(s,x)\in\Omega$, $\delta\in(0,\infty)$,
    $\gamma\in C^1((-\delta,\delta),\mathbb S)$ with
    $\gamma(0)=x$ and $\{s\}\times\gamma((-\delta,\delta))\subseteq\Omega$
    it holds that
    \begin{eqn}[a]
      \llabel{eq:angular-projection-pairing}
      &\bigl\langle s\bigl(\grad f(sx)
        -\langle\grad f(sx),x\rangle x\bigr),\gamma'(0)\bigr\rangle\\
      &=s\langle\grad f(sx),\gamma'(0)\rangle
        -s\langle\grad f(sx),x\rangle\langle x,\gamma'(0)\rangle
        =\frac{\dd}{\dd t}h(s,\gamma(t))\Big|_{t=0}
      \dott
    \end{eqn}
  }
  \argument{
    \lref{setting1};
    \lref{eq:polar-radial-derivative};
    \lref{eq:angular-projection-tangent};
    \lref{eq:angular-projection-pairing};
    the fact that $\Omega$ is definable and open in $(0,R)\times\mathbb S$;
    the fact that $h\in\COne{\Omega}{\R}$ is definable;
    the fact that $f\in C^1(W,(-1,1))$ is definable;
    \cref{lem:angular-gradient-existence};
    \cref{lem:angular-gradient-well-defined};
  }{
    that $u$ and $G$ are continuous and definable and that for all
    $(s,x)\in\Omega$ it holds that
    \begin{gather}
      \partial_s h(s,x)
        =\langle\grad f(sx),x\rangle,\qquad
        u(s,x)=s\bigl(\grad f(sx)-\langle\grad f(sx),x\rangle x\bigr),\notag\\
      \llabel{eq0}
      \grad f(sx)=\partial_s h(s,x)x+\frac{u(s,x)}s,
        \qquad\text{and}\qquad\langle u(s,x),x\rangle=0
      \dott
    \end{gather}
  }
  \argument{
    the fact that $\mathbf{k}$ is definable;
    $C^1$ cell decomposition \cite[Theorem 6.6 and Theorem 6.7]{Coste};
  }{
    that there exists $\varepsilon\in(0,R)$ which satisfies that
    \begin{eqn}
      \llabel{eq:mixing-regularity}
      \mathbf{k}|_{(0,\varepsilon)}\in C^1((0,\varepsilon),[0,\infty))
      \dott
    \end{eqn}
  }
  \startnewargseq
  \argument{
    \lref{setting1};
    \lref{bounded-potential};
    the fact that $\mathbb S$ is a nonempty, compact, and definable \submanifold{2}{n};
    the fact that $\Omega$ is definable and open in $(0,R)\times\mathbb S$;
    the fact that $h\in\COne{\Omega}{\R}$ is definable;
    \cref{lem:divide}
    (applied with $D\curvearrowleft n$, $M\curvearrowleft\mathbb S$,
    $R\curvearrowleft\varepsilon$,
    $\Omega\curvearrowleft\Omega\cap((0,\varepsilon)\times\mathbb S)$,
    $h\curvearrowleft h|_{\Omega\cap((0,\varepsilon)\times\mathbb S)}$,
    $G\curvearrowleft G|_{\Omega\cap((0,\varepsilon)\times\mathbb S)}$,
    $a\curvearrowleft((0,\varepsilon)\ni s\mapsto-1\in\R)$,
    $b\curvearrowleft((0,\varepsilon)\ni s\mapsto1\in\R)$,
    $c\curvearrowleft K^{-1/2}$ in the notation of \cref{lem:divide});
  }{
    that there exist $\rho\in(0,\varepsilon)$, $N\in\N$, $H\in(0,\infty)$
    and definable
    $\mathfrak e_0,\mathfrak e_1,\dots,\mathfrak e_N\in C^2((0,\rho),\R)$,
    $\mathcal F_0,\mathcal F_1,\dots,\mathcal F_N
      \in C((0,\rho)\times\R,[0,H])$,
    $E_1=(E_1(s,z))_{(s,z)\in(0,\rho)\times(0,H)}$,
    $E_2=(E_2(s,z))_{(s,z)\in(0,\rho)\times(0,H)}$, $\dots$,
    $E_N=(E_N(s,z))_{(s,z)\in(0,\rho)\times(0,H)}
      \in C^2((0,\rho)\times(0,H),\R)$, and
    $T_1,T_2,\dots,T_N\in C^1((0,\rho)\times(0,H),(0,\infty))$
    which satisfy that
    \begin{enumerate}[label=(\Roman*)]
      \item \llabel{atlas1} for all $s\in(0,\rho)$ it holds that
      \begin{eqn}
        \llabel{eq2}
        -1=\mathfrak e_0(s)<\cdots<\mathfrak e_N(s)=1,
      \end{eqn}
      \item \llabel{atlas2} for all $s\in(0,\rho)$, $e\in\R$ with
      $\inf\bigl(\{\infty\}\cup
        \{z\in\R\colon\exists\,x\in\mathbb S\colon
          (s,x)\in\Omega\cap h^{-1}(\{e\})\cap G^{-1}(\{z\})\}\bigr)=0$
      it holds that
      \begin{eqn}
        \llabel{eq2.1}
        e\in\{\mathfrak e_0(s),\mathfrak e_1(s),\dots,\mathfrak e_N(s)\},
      \end{eqn}
      \item \llabel{atlas3} for all $i\in\{0,1,\dots,N-1\}$, $s\in(0,\rho)$,
      $v\in(-\infty,\mathfrak e_i(s)]$,
      $w\in[\mathfrak e_{i+1}(s),\infty)$, and
      $z,\zeta\in(\mathfrak e_i(s),\mathfrak e_{i+1}(s))$ with $z<\zeta$
      it holds that
      \begin{eqn}
        \llabel{eq2.2}
        0=\mathcal F_i(s,v)<\mathcal F_i(s,z)
          <\mathcal F_i(s,\zeta)<\mathcal F_i(s,w)=H,
      \end{eqn}
      \item \llabel{atlas4} for all $i\in\{0,1,\dots,N-1\}$ it holds that
      $\mathcal F_i|_{\{(s,v)\in(0,\rho)\times\R\colon
        \mathfrak e_i(s)<v<\mathfrak e_{i+1}(s)\}}$
      is twice continuously differentiable,
      \item \llabel{atlas5} for all $i\in\{1,2,\dots,N\}$,
      $s\in(0,\rho)$, $z\in(0,H)$ it holds that
      \begin{eqn}
        \llabel{eq2.4}
        \mathfrak e_{i-1}(s)<E_i(s,z)<\mathfrak e_i(s),\quad
        \mathcal F_{i-1}(s,E_i(s,z))=z,\quad\text{and}\quad
        T_i(s,z)=\partial_z E_i(s,z)>0,
      \end{eqn}
      and
      \item \llabel{atlas6} for all $i\in\{1,2,\dots,N\}$,
      $(s,x)\in\Omega\cap((0,\rho)\times\mathbb S)$ with
      $\mathfrak e_{i-1}(s)<h(s,x)<\mathfrak e_i(s)$ it holds that
      \begin{eqn}
        \llabel{eq2.5}
        T_i(s,\mathcal F_{i-1}(s,h(s,x)))\leq K^{-1/2}G(s,x)
        \dott
      \end{eqn}
    \end{enumerate}
  }
  In the following for every $i\in\{0,1,\dots,N-1\}$ let
  $\mathbf{\Omega}_i\subseteq\Omega\cap((0,\rho)\times\mathbb S)$ satisfy
  \begin{eqn}
    \llabel{setting2}
    \mathbf{\Omega}_i=\{(s,x)\in\Omega\cap((0,\rho)\times\mathbb S)\colon
      \mathfrak e_i(s)<h(s,x)<\mathfrak e_{i+1}(s)\},
  \end{eqn}
  let $Z_i\colon\Omega\cap((0,\rho)\times\mathbb S)\to[0,H]$
  satisfy for all $(s,x)\in\Omega\cap((0,\rho)\times\mathbb S)$ that
  \begin{eqn}
    \llabel{setting2.1}
    Z_i(s,x)=\mathcal F_i(s,h(s,x)),
  \end{eqn}
  and let $c_i\colon(0,\rho)\times(0,H)\to\R$ satisfy
  for all $(s,z)\in(0,\rho)\times(0,H)$ that
  \begin{eqn}
    \llabel{setting2.2}
    c_i(s,z)=\partial_s E_{i+1}(s,z)
    \dott
  \end{eqn}
  \startnewargseq
  \argument{
    \lref{atlas1};
    \lref{atlas3};
    \lref{atlas4};
    \lref{setting2};
    \lref{setting2.1};
    \lref{setting2.2};
    the fact that $\Omega$ is definable and open in $(0,R)\times\mathbb S$;
    the fact that $h\in\COne{\Omega}{\R}$ is definable;
    the fact that for all $j\in\{0,1,\dots,N\}$ it holds that
    $\mathfrak e_j\in C^2((0,\rho),\R)$ is definable;
    the fact that for all $i\in\{0,1,\dots,N-1\}$ it holds that
    $\mathcal F_i$ is continuous and definable; the fact that for all
    $i\in\{0,1,\dots,N-1\}$ it holds that $E_{i+1}\in C^2((0,\rho)\times(0,H),\R)$;
  }{
    that for all $i\in\{0,1,\dots,N-1\}$ it holds that $\mathbf{\Omega}_i$ is
    definable and open in $(0,\rho)\times\mathbb S$, that $Z_i$ is
    continuous and definable, that $Z_i|_{\mathbf{\Omega}_i}\in\COne{\mathbf{\Omega}_i}{[0,H]}$
    and $c_i\in C^1((0,\rho)\times(0,H),\R)$,
    and that for all $(s,x)\in\Omega\cap((0,\rho)\times\mathbb S)$ it holds that
    \begin{gather}
      \llabel{eq3}
      \mathbf{\Omega}_i=\{(v,\vartheta)\in\Omega\cap((0,\rho)\times\mathbb S)\colon
        0<Z_i(v,\vartheta)<H\}\\
      \llabel{eq:coordinate-saturation}
      \text{and}\qquad Z_i(s,x)=\begin{cases}
        0&\colon h(s,x)\leq\mathfrak e_i(s)\\
        H&\colon h(s,x)\geq\mathfrak e_{i+1}(s)\\
        \mathcal F_i(s,h(s,x))&\colon(s,x)\in\mathbf{\Omega}_i
      \end{cases}
      \dott
    \end{gather}
  }
  \argument{
    \lref{setting2};
    \lref{eq3};
  }{
    that for all $i\in\{0,1,\dots,N-1\}$,
    $(s,x)\in\mathbf{\Omega}_i$ it holds that
    \begin{eqn}
      \llabel{eq:inverse-domain}
      0<Z_i(s,x)<H\qquad\text{and}\qquad
      \mathfrak e_i(s)<h(s,x)<\mathfrak e_{i+1}(s)
      \dott
    \end{eqn}
  }
  \argument{
    \lref{atlas5};
    \lref{setting2.1};
    \lref{eq:inverse-domain};
  }{
    that for all $i\in\{0,1,\dots,N-1\}$,
    $(s,x)\in\mathbf{\Omega}_i$ it holds that
    \begin{eqn}[a]
      \llabel{eq:inverse-values}
      \mathfrak e_i(s)<E_{i+1}(s,Z_i(s,x))<\mathfrak e_{i+1}(s)\quad
      \text{and}\quad
      \mathcal F_i\bigl(s,E_{i+1}(s,Z_i(s,x))\bigr)
        =Z_i(s,x)=\mathcal F_i(s,h(s,x))
      \dott
    \end{eqn}
  }
  \argument{
    \lref{atlas3};
    \lref{atlas5};
    \lref{atlas6};
    \lref{setting2};
    \lref{eq:inverse-domain};
    \lref{eq:inverse-values};
  }{
    that for all $i\in\{0,1,\dots,N-1\}$,
    $(s,z)\in(0,\rho)\times(0,H)$ it holds that
    $\partial_z E_{i+1}(s,z)=T_{i+1}(s,z)$ and that for all
    $i\in\{0,1,\dots,N-1\}$, $(s,x)\in\mathbf{\Omega}_i$ it holds that
    \begin{eqn}
      \llabel{eq4}
      h(s,x)=E_{i+1}(s,Z_i(s,x))\qquad\text{and}\qquad
      0<T_{i+1}(s,Z_i(s,x))\leq K^{-1/2}G(s,x)
      \dott
    \end{eqn}
  }
  \argument{
    \lref{atlas1};
    \lref{atlas3};
    \lref{atlas5};
    the fact that $\Omega$ is definable and open in $(0,R)\times\mathbb S$;
    the fact that $h\in\COne{\Omega}{\R}$ is definable;
    the fact that for all $j\in\{0,1,\dots,N\}$ it holds that
    $\mathfrak e_j\in C^2((0,\rho),\R)$ is definable;
    the fact that for all $i\in\{0,1,\dots,N-1\}$ it holds that
    $\mathcal F_i$ is continuous and definable and that
    $E_{i+1}\in C^2((0,\rho)\times(0,H),\R)$ is definable;
    \cref{lem:residual} (applied with $D\curvearrowleft n$, $M\curvearrowleft\mathbb S$,
    $\mathfrak r_0\curvearrowleft\varepsilon$,
    $\mathfrak r_1\curvearrowleft\rho$,
    $\mathcal O\curvearrowleft\Omega\cap((0,\varepsilon)\times\mathbb S)$,
    $h\curvearrowleft h|_{\Omega\cap((0,\varepsilon)\times\mathbb S)}$,
    $\mathbf c\curvearrowleft K^{-1/2}$,
    $(\mathcal F_1,\mathcal F_2,\dots,\mathcal F_N)\curvearrowleft
      (\mathcal F_0,\mathcal F_1,\dots,\mathcal F_{N-1})$,
    $(\mathbf m_1,\mathbf m_2,\dots,\mathbf m_N)\curvearrowleft
      (\eta_1,\eta_2,\dots,\eta_N)$
    in the notation of \cref{lem:residual});
  }{
    that there exist $\sigma\in(0,\rho)$ and integrable, definable
    $\eta_1,\eta_2,\dots,\eta_N\in C^1((0,\sigma),[1,\infty))$
    which satisfy for all $i\in\{1,2,\dots,N\}$,
    $(s,x)\in\Omega\cap((0,\sigma)\times\mathbb S)$ with
    $\mathfrak e_{i-1}(s)<h(s,x)<\mathfrak e_i(s)$ that
    \begin{eqn}
      \llabel{eq:residual-raw}
      |\partial_s h(s,x)
        -(\partial_s E_i)(s,\mathcal F_{i-1}(s,h(s,x)))|
      \leq\eta_i(s)\bigl(
        (\partial_z E_i)(s,\mathcal F_{i-1}(s,h(s,x)))
        +K^{-1/2}\norm{(\nabla_{\mathbb S}h)(s,x)}\bigr)
      \dott
    \end{eqn}
  }
  \startnewargseq
  \argument{
    \lref{setting1};
    \lref{atlas5};
    \lref{setting2};
    \lref{setting2.1};
    \lref{setting2.2};
    \lref{eq:inverse-domain};
    \lref{eq:residual-raw};
  }{
    that for all $i\in\{0,1,\dots,N-1\}$,
    $(s,x)\in\mathbf{\Omega}_i\cap((0,\sigma)\times\mathbb S)$ it holds that
    \begin{eqn}
      \llabel{eq5}
      |\partial_s h(s,x)-c_i(s,Z_i(s,x))|
        \leq\eta_{i+1}(s)
          \bigl(T_{i+1}(s,Z_i(s,x))+K^{-1/2}G(s,x)\bigr)
      \dott
    \end{eqn}
  }
  In the following for every $i\in\{0,1,\dots,N-1\}$ let
  $\mathfrak n_i\colon(0,\sigma)\to[1,\infty)$,
  $\mathbf{w}_i\colon(0,\sigma)\times(0,H)\to(0,1]$, and
  $P_i\colon(0,\sigma)\times[0,H]\to\R$ satisfy for all
  $s\in(0,\sigma)$, $z\in(0,H)$, $v\in[0,H]$ that
  $\mathfrak n_i(s)=\eta_{i+1}(s)+\sqrt K\,\mathbf{k}(s)$,
  \begin{gather}
    \llabel{setting3}
    \textstyle \mathbf{w}_i(s,z)=\frac{T_{i+1}(s,z)}
      {T_{i+1}(s,z)+8s^2\mathfrak n_i(s)
        \sqrt{(c_i(s,z))^2+(T_{i+1}(s,z))^2}},
    \qquad\text{and}\qquad P_i(s,v)=\int_{(0,v)}\mathbf{w}_i(s,z)\dd z
    \dott
  \end{gather}
  \startnewargseq
  \argument{
    \lref{ass0};
    \lref{eq:mixing-regularity};
    \lref{setting3};
    the assumption that $\mathbf{k}$ is definable;
    the fact that for all $i\in\{0,1,\dots,N-1\}$ it holds that
    $\eta_{i+1}\in C^1((0,\sigma),[1,\infty))$ is definable and integrable;
  }{
    that for all $i\in\{0,1,\dots,N-1\}$ it holds that
    $\mathfrak n_i\in C^1((0,\sigma),[1,\infty))$ is definable and integrable\dott
    \llabel{arg:radial-weight-regularity}
  }
  \argument{
    \lref{arg:radial-weight-regularity};
    \cref{lem:radial-decay};
  }{
    that there exists
    $\nu\in(0,\sigma)$ which satisfies for all
    $i\in\{0,1,\dots,N-1\}$, $s\in(0,\nu)$ that
    \begin{eqn}
      \llabel{eq6}
      s^2(\mathfrak n_i(s))^2\leq\frac1{64}
      \dott
    \end{eqn}
  }
  \startnewargseq
  \argument{
    \lref{setting3};
    the fact that for all $i\in\{0,1,\dots,N-1\}$ it holds that
    $\mathbf{w}_i\in C^1((0,\sigma)\times(0,H),(0,1])$ is definable;
    \cref{lem:primitive};
    the fact that for all $i\in\{0,1,\dots,N-1\}$,
    $s\in(0,\sigma)$ it holds that
    $[0,H]\ni v\mapsto P_i(s,v)\in\R$ is continuous;
  }{
    that there exist $\mu\in(0,\nu)$ and definable and integrable
    $\mathfrak b_0,\mathfrak b_1,\dots,\mathfrak b_{N-1}\in C((0,\mu),[0,\infty))$
    which satisfy for all
    $i\in\{0,1,\dots,N-1\}$ that
    $P_i|_{(0,\mu)\times[0,H]}$ is locally Lipschitz
    and which satisfy for all $i\in\{0,1,\dots,N-1\}$, $s,t\in(0,\mu)$,
     $z,\zeta\in[0,H]$ that
    \begin{gather}
      P_i(s,0)=0,\qquad 0\leq P_i(s,z)\leq z\leq H,\qquad
        \abs{P_i(s,\zeta)-P_i(s,z)}\leq\abs{\zeta-z},\notag\\
      \llabel{eq7}
      \text{and}\qquad\abs{P_i(t,z)-P_i(s,z)}
        \leq\int_{\min\{s,t\}}^{\max\{s,t\}}\mathfrak b_i(v)\dd v
      \dott
    \end{gather}
  }
  \startnewargseq
  \argument{
    \lref{setting1};
    the fact that $\Omega$ is definable and open in $(0,R)\times\mathbb S$;
    the fact that $h\in\COne{\Omega}{\R}$ is definable;
    the fact that for all $j\in\{0,1,\dots,N\}$ it holds that
    $\mathfrak e_j\in C^2((0,\rho),\R)$ is definable;
    \cref{lem:contact-bound}
    (applied for every $j\in\{0,1,\dots,N\}$ with $M\curvearrowleft\mathbb S$,
    $\mathfrak r_0\curvearrowleft\mu$,
    $\mathcal O\curvearrowleft\Omega\cap((0,\mu)\times\mathbb S)$,
    $h\curvearrowleft h|_{\Omega\cap((0,\mu)\times\mathbb S)}$,
    $a\curvearrowleft\mathfrak e_j|_{(0,\mu)}$
    in the notation of \cref{lem:contact-bound});
  }{
    that there exist $\mathfrak{R}\in(0,\mu)$ and integrable, definable
    $\beta_0,\beta_1,\dots,\beta_N\in C^1((0,\mathfrak{R}),[0,\infty))$ which satisfy
    for all $j\in\{0,1,\dots,N\}$,
    $(s,x)\in\Omega\cap((0,\mathfrak{R})\times\mathbb S)$ with
    $h(s,x)=\mathfrak e_j(s)$ and $G(s,x)>0$ that
    \begin{eqn}
      \llabel{eq8}
      \abs{\partial_s h(s,x)-(\mathfrak e_j)'(s)}
        \leq\beta_j(s)G(s,x)
      \dott
    \end{eqn}
  }
  In the following let
  $Q\colon\Omega\cap((0,\mathfrak{R})\times\mathbb S)\to[0,\infty)$
  satisfy for all $(s,x)\in\Omega\cap((0,\mathfrak{R})\times\mathbb S)$ that
  \begin{eqn}
    \llabel{eq:state-definition}
    Q(s,x)=4\sqrt K\sum_{i=0}^{N-1}P_i(s,Z_i(s,x))
  \end{eqn}
  and let
  $\mathfrak m,m\colon(0,\mathfrak{R})\to[0,\infty)$ satisfy
  for all $s\in(0,\mathfrak{R})$ that
  \begin{eqn}[a]
    \llabel{setting4}
  \textstyle  \mathfrak m(s)=4\sqrt K\sum\limits_{i=0}^{N-1}
      \bigl(5\mathfrak n_i(s)+\mathfrak b_i(s)\bigr)+\mathbf{k}(s)
    \quad\text{and}\quad m(s)=\mathfrak m(s)+K\biggl(\sum\limits_{j=0}^N\beta_j(s)\biggr)+(K+1)\mathbf{k}(s)
    \dott
  \end{eqn}
  \startnewargseq
  \argument{
    \lref{eq3};
    \lref{eq7};
    \lref{eq:state-definition};
    the fact that for all $i\in\{0,1,\dots,N-1\}$ it holds that
    $Z_i$ is continuous and that $P_i|_{(0,\mu)\times[0,H]}$
    is locally Lipschitz;
  }{
    that $Q$ is continuous and that for all
    $(s,x)\in\Omega\cap((0,\mathfrak{R})\times\mathbb S)$ it holds that
    \begin{eqn}
      \llabel{eq9}
      0\leq Q(s,x)\leq4\sqrt K\,NH
      \dott
    \end{eqn}
  }
  \argument{
    \lref{ass0};
    \lref{eq:mixing-regularity};
    \lref{setting4};
    the fact that for all $i\in\{0,1,\dots,N-1\}$ and
    $j\in\{0,1,\dots,N\}$ it holds that $\mathfrak n_i$, $\mathfrak b_i$,
    and $\beta_j$ are continuous, definable, nonnegative, and integrable;
  }{
    that $\mathfrak m$ and $m$ are continuous and definable and that
    \begin{eqn}
      \llabel{eq10}
      0\leq\int_0^{\mathfrak{R}}\mathfrak m(s)\dd s
        \leq\int_0^{\mathfrak{R}}m(s)\dd s<\infty
      \dott
    \end{eqn}
  }

  We claim that for every $a\in\R$, $b\in(a,\infty)$,
  every absolutely continuous $y\colon[a,b]\to W$,
  every Lebesgue measurable $\lambda\colon[a,b]\to[0,\infty)$,
  every Lebesgue measurable $B\colon[a,b]\to\R^{n\times n}$,
  every $r\colon[a,b]\to(0,\mathfrak{R})$,
  every $\omega\colon[a,b]\to\mathbb S$,
  and every $q\colon[a,b]\to\R$
  with the property that for all $t\in[a,b]$ it holds that
  \begin{eqn}
    r(t)=\norm{y(t)},\qquad
    \omega(t)=y(t)/r(t),\qquad\text{and}\qquad
    q(t)=Q(r(t),\omega(t))
  \end{eqn}
  and with the property that for Lebesgue-almost all $t\in[a,b]$ it holds
  that $y$ is differentiable at $t$, that $B(t)$ is symmetric and
  positive-definite, and that
  \begin{gather}
    y'(t)=-\lambda(t)B(t)\grad f(y(t)),\qquad
    \max_{z\in\R^n,\,\norm{z}=1}\langle B(t)z,z\rangle
      \leq K\min_{z\in\R^n,\,\norm{z}=1}\langle B(t)z,z\rangle,\\
    \text{and}\qquad\frac{\norm{B(t)\omega(t)
       -\langle B(t)\omega(t),\omega(t)\rangle\omega(t)}}
      {\langle B(t)\omega(t),\omega(t)\rangle}
      \leq r(t)\mathbf{k}(r(t))
  \end{gather}
  we have that
  \begin{enumerate}[label=\alph*)]
    \item \llabel{item:compact-absolute-continuity} $q$ is absolutely continuous and
    \item \llabel{item:compact-differential-inequality}
    for Lebesgue-almost all $t\in[a,b]$ it holds that
    $q$, $r$, and $\omega$ are differentiable at $t$ and that
    \begin{eqn}
      \llabel{claim1}
      q'(t)\leq-\norm{\omega'(t)}+m(r(t))\abs{r'(t)}
      \dott
    \end{eqn}
  \end{enumerate}
  \startnewargseq
  To prove this claim, let $a\in\R$, $b\in(a,\infty)$,
  let $y\colon [a,b]\to W$ be absolutely continuous,
  let $\lambda\colon [a,b]\to[0,\infty)$
  and $B\colon [a,b]\to\R^{n\times n}$ be Lebesgue measurable,
  let $r\colon [a,b]\to(0,\mathfrak{R})$,
  $\omega\colon [a,b]\to\mathbb S$, and $q\colon [a,b]\to\R$,
  assume for all $t\in [a,b]$ that
  \begin{eqn}
    \llabel{setting4.1}
    r(t)=\norm{y(t)},\qquad
    \omega(t)=y(t)/r(t),\qquad\text{and}\qquad
    q(t)=Q(r(t),\omega(t)),
  \end{eqn}
  and assume that for Lebesgue-almost all $t\in [a,b]$ it holds that
  $y$ is differentiable at $t$, that $B(t)$ is symmetric and
  positive-definite, and that
  \begin{gather}
    \llabel{setting4.2}
    y'(t)=-\lambda(t)B(t)\grad f(y(t)),\qquad
    \max_{z\in\R^n,\,\norm{z}=1}\langle B(t)z,z\rangle
      \leq K\min_{z\in\R^n,\,\norm{z}=1}\langle B(t)z,z\rangle,\\
    \llabel{setting4.3}
    \text{and}\qquad\frac{\norm{B(t)\omega(t)
       -\langle B(t)\omega(t),\omega(t)\rangle\omega(t)}}
      {\langle B(t)\omega(t),\omega(t)\rangle}
      \leq r(t)\mathbf{k}(r(t))
    \dott
  \end{gather}
  \startnewargseq
  In the following let
  $\mathbf{h}\colon [a,b]\to\R$ and $U\colon [a,b]\to\R^n$
  satisfy for all $t\in [a,b]$ that
  \begin{eqn}
    \llabel{setting5}
    \mathbf{h}(t)=h(r(t),\omega(t))\qquad\text{and}\qquad
    U(t)=u(r(t),\omega(t)),
  \end{eqn}
  let $\mathcal D=\{t\in [a,b]\colon B(t)\text{ is symmetric
  and positive-definite}\}$, and let
  $\mathsf A\colon [a,b]\to\R^{n\times n}$,
  $\alpha\colon [a,b]\to[0,\infty)$,
  $\xi\colon [a,b]\to\R^n$, and $\mathbf{v}\colon [a,b]\to\R$ satisfy for all
  $t\in\mathcal D$ that
  \begin{gather}
    \textstyle \xi(t)=\frac{B(t)\omega(t)}{\langle B(t)\omega(t),\omega(t)\rangle}
      -\omega(t),\qquad
    \alpha(t)=\frac{\lambda(t)\langle B(t)\omega(t),\omega(t)\rangle}{(r(t))^2},\notag\\
    \llabel{setting6}
   \textstyle \mathsf A(t)=\frac{B(t)}{\langle B(t)\omega(t),\omega(t)\rangle}
      -\bigl(\omega(t)+\xi(t)\bigr)\bigl(\omega(t)+\xi(t)\bigr)^{\mathsf T},
      \quad\text{and}\quad \mathbf{v}(t)=(\partial_s h)(r(t),\omega(t))+\frac{\langle\xi(t),U(t)\rangle}{r(t)}
  \end{gather}
  and for all $t\in [a,b]\backslash\mathcal D$ that
  \begin{eqn}
    \llabel{eq:zero-extensions}
    \mathsf A(t)=0,\qquad \alpha(t)=0,\qquad
    \xi(t)=0,\qquad\text{and}\qquad \mathbf{v}(t)=0
    \dott
  \end{eqn}
  \startnewargseq
  \argument{
    \lref{setting1};
    \lref{setting4.1};
    \lref{setting5};
    the fact that $f\in C^1(W,(-1,1))$ and $h\in\COne{\Omega}{\R}$;
    the fact that $u$ is continuous;
    the fact that $y$ is absolutely continuous;
    the chain rule;
  }{
    that $r$, $\omega$, and $\mathbf{h}=f\circ y$ are absolutely
    continuous, that $U$ is continuous, and that $r([a,b])$ is a compact
    subset of $(0,\mathfrak{R})$\dott
  }
  \argument{
    \lref{setting6};
    \lref{eq:zero-extensions};
  }{
    that for all $t\in [a,b]$ it holds that $\mathsf A(t)$ is symmetric
    and $\mathsf A(t)\omega(t)=0$ and that for all
    $t\in\mathcal D$, $z\in\R^n$ with
    $\langle z,\omega(t)\rangle=0$ it holds that
    \begin{eqn}[a]
      \llabel{eq11}
      \langle\mathsf A(t)z,z\rangle
        =\frac{\langle B(t)z,z\rangle}{\langle B(t)\omega(t),\omega(t)\rangle}
          -\langle\xi(t),z\rangle^2=\min_{\kappa\in\R}
          \frac{\langle B(t)(\kappa\omega(t)+z),\kappa\omega(t)+z\rangle}
            {\langle B(t)\omega(t),\omega(t)\rangle}
          \geq0
      \dott
    \end{eqn}
  }
  \argument{
    \lref{setting4.2};
    \lref{setting4.3};
    \lref{eq11};
  }{
    that for Lebesgue-almost all $t\in [a,b]$ it holds that for all
    $z\in\R^n$ with
    $\langle z,\omega(t)\rangle=0$ it holds that
    \begin{eqn}
      \llabel{eq12}
      K^{-1}\norm{z}^2\leq\langle\mathsf A(t)z,z\rangle
        \leq K\norm{z}^2\qquad\text{and}\qquad
      \frac{\norm{\xi(t)}}{r(t)}\leq\mathbf{k}(r(t))
      \dott
    \end{eqn}
  }
  In the following let $\gamma,\ell\colon [a,b]\to[0,\infty)$ satisfy
  for all $t\in [a,b]$ that
  \begin{eqn}
    \llabel{setting7}
    \gamma(t)=\sqrt{\langle\mathsf A(t)U(t),U(t)\rangle}\qquad\text{and}\qquad
    \ell(t)=\begin{cases}
      \norm{\omega'(t)}&\colon\omega\text{ is differentiable at }t\\
      0&\colon\text{otherwise}
    \end{cases}
    \dott
  \end{eqn}
  \startnewargseq
  \argument{
    \lref{setting5};
    \lref{setting6};
    \lref{eq:zero-extensions};
    \lref{eq11};
    \lref{setting7};
    the fact that $B$ and $\lambda$ are Lebesgue measurable;
    the fact that $r$, $\omega$, and $U$ are continuous;
    the fact that $\omega$ is absolutely continuous;
  }{
    that $\alpha,\gamma,\ell$ are Lebesgue measurable and that $\ell$
    is integrable\dott
  }
  \argument{
    \lref{eq12};
    \lref{setting7};
    the fact that for all $t\in [a,b]$ it holds that
    $\mathsf A(t)$ is symmetric, $\mathsf A(t)\omega(t)=0$,
    and $\langle U(t),\omega(t)\rangle=0$;
  }{
    that for Lebesgue-almost all $t\in [a,b]$ it holds that
    \begin{eqn}
      \llabel{eq:tangent-norm-bound}
      \langle\mathsf A(t)U(t),\omega(t)\rangle
        =\langle U(t),\mathsf A(t)\omega(t)\rangle=0
      \qquad\text{and}\qquad
      K^{-1}\norm{U(t)}^2\leq(\gamma(t))^2
      \dott
    \end{eqn}
  }
  \argument{
    \lref{eq12};
    \lref{setting7};
    \lref{eq:tangent-norm-bound};
    the fact that for all $t\in [a,b]$ it holds that
    $\mathsf A(t)$ is symmetric and $\langle U(t),\omega(t)\rangle=0$;
    the Cauchy--Schwarz inequality;
  }{
    that for Lebesgue-almost all $t\in [a,b]$ it holds that
    \begin{eqn}
      \llabel{eq:quadratic-norm-lower}
      (\gamma(t))^4
        =\langle\mathsf A(t)U(t),U(t)\rangle^2
        \leq\norm{\mathsf A(t)U(t)}^2\norm{U(t)}^2
        \leq K(\gamma(t))^2\norm{\mathsf A(t)U(t)}^2
    \end{eqn}
    and
    \begin{eqn}[a]
      \llabel{eq:quadratic-norm-comparison}
      0&\leq\bigl\langle
        \mathsf A(t)\bigl(U(t)-K^{-1}\mathsf A(t)U(t)\bigr),
        U(t)-K^{-1}\mathsf A(t)U(t)\bigr\rangle\\
      &=(\gamma(t))^2-2K^{-1}\norm{\mathsf A(t)U(t)}^2
        +K^{-2}\langle(\mathsf A(t))^2U(t),\mathsf A(t)U(t)\rangle
        \leq(\gamma(t))^2-K^{-1}\norm{\mathsf A(t)U(t)}^2
      \dott
    \end{eqn}
  }
  \argument{
    \lref{eq:tangent-norm-bound};
    \lref{eq:quadratic-norm-lower};
    \lref{eq:quadratic-norm-comparison};
  }{
    that for Lebesgue-almost all $t\in [a,b]$ it holds that
    \begin{eqn}
      \llabel{eq13}
      K^{-1/2}\norm{U(t)}\leq\gamma(t)\qquad\text{and}\qquad
      K^{-1}(\gamma(t))^2\leq\norm{\mathsf A(t)U(t)}^2
        \leq K(\gamma(t))^2
      \dott
    \end{eqn}
  }
  \argument{
    \lref{eq0};
    \lref{setting4.1};
    \lref{setting4.2};
    \lref{setting5};
    \lref{setting6};
    the chain rule;
  }{
    that for Lebesgue-almost all $t\in [a,b]$ it holds that
    $r$ and $\omega$ are differentiable at $t$, that
    \begin{eqn}[a]
      \llabel{eq14}
      r'(t)&=\langle y'(t),\omega(t)\rangle
        =-\lambda(t)\langle B(t)\grad f(y(t)),\omega(t)\rangle\\
        &=-\lambda(t)\biggl(
          (\partial_s h)(r(t),\omega(t))\langle B(t)\omega(t),\omega(t)\rangle
          +\frac{\langle B(t)U(t),\omega(t)\rangle}{r(t)}\biggr)\\
        &=-\lambda(t)\langle B(t)\omega(t),\omega(t)\rangle
          \biggl((\partial_s h)(r(t),\omega(t))
          +\frac{\langle\xi(t),U(t)\rangle}{r(t)}\biggr)
          =-\alpha(t)(r(t))^2\mathbf{v}(t),
    \end{eqn}
    and that
    \begin{eqn}[a]
      \llabel{eq:angular-evolution}
      \omega'(t)&=\frac{y'(t)-r'(t)\omega(t)}{r(t)}
        =-\frac{\lambda(t)}{r(t)}B(t)
          \biggl((\partial_s h)(r(t),\omega(t))\omega(t)
          +\frac{U(t)}{r(t)}\biggr)-\frac{r'(t)}{r(t)}\omega(t)\\
        &=-\alpha(t)\biggl(
          \frac{B(t)U(t)}{\langle B(t)\omega(t),\omega(t)\rangle}
          +r(t)(\partial_s h)(r(t),\omega(t))
            \bigl(\omega(t)+\xi(t)\bigr)\biggr)
          -\frac{r'(t)}{r(t)}\omega(t)\\
        &=-\alpha(t)\mathsf A(t)U(t)
          -\alpha(t)\bigl(r(t)(\partial_s h)(r(t),\omega(t))
            +\langle\xi(t),U(t)\rangle\bigr)
            \bigl(\omega(t)+\xi(t)\bigr)
          -\frac{r'(t)}{r(t)}\omega(t)\\
        &=-\alpha(t)\mathsf A(t)U(t)
          -\alpha(t)r(t)\mathbf{v}(t)\bigl(\omega(t)+\xi(t)\bigr)
          -\frac{r'(t)}{r(t)}\omega(t)\\
        &=-\alpha(t)\mathsf A(t)U(t)-\alpha(t)r(t)\mathbf{v}(t)\xi(t)
          =-\alpha(t)\mathsf A(t)U(t)+\frac{\xi(t)}{r(t)}r'(t)
      \dott
    \end{eqn}
  }
  \argument{
    \lref{setting5};
    \lref{setting6};
    \lref{setting7};
    \lref{eq14};
    \lref{eq:angular-evolution};
    the fact that $h\in\COne{\Omega}{\R}$;
    the chain rule;
  }{
    that for Lebesgue-almost all $t\in [a,b]$ it holds that
    $\mathbf{h}$ is differentiable at $t$ and that
    \begin{eqn}[a]
      \llabel{eq15}
      \mathbf{h}'(t)&=(\partial_s h)(r(t),\omega(t))r'(t)+\langle U(t),\omega'(t)\rangle\\
        &=\mathbf{v}(t)r'(t)-\alpha(t)(\gamma(t))^2
        =-\alpha(t)\bigl((\gamma(t))^2+(r(t))^2(\mathbf{v}(t))^2\bigr)
      \dott
    \end{eqn}
  }
  \argument{
    \lref{setting6};
      \lref{eq12};
    \lref{setting7};
    \lref{eq13};
    \lref{eq:angular-evolution};
    the triangle inequality;
  }{
    that for Lebesgue-almost all $t\in [a,b]$ it holds that
    \begin{equation}
      \llabel{eq16}
      \ell(t)\leq\sqrt K\,\alpha(t)\gamma(t)
          +\mathbf{k}(r(t))\abs{r'(t)}
    \end{equation}
    and
    \begin{eqn}[a]
      \llabel{eq:reverse-angular-speed}
      \alpha(t)\gamma(t)
        &\leq\sqrt K\,\alpha(t)\norm{\mathsf A(t)U(t)}
          =\sqrt K\biggl|-\omega'(t)+\frac{\xi(t)}{r(t)}r'(t)\biggr|\\
        &\leq\sqrt K\biggl(\norm{\omega'(t)}
          +\frac{\norm{\xi(t)}}{r(t)}\abs{r'(t)}\biggr)
          \leq\sqrt K\bigl(\ell(t)+\mathbf{k}(r(t))\abs{r'(t)}\bigr)
      \dott
    \end{eqn}
  }
  \argument{
    \lref{eq:reverse-angular-speed};
    the fact that $r$ is absolutely continuous;
    the fact that $\ell$ is integrable;
    the fact that $\mathbf{k}$ is bounded on $r([a,b])$;
  }{
    that $\alpha\gamma$ is integrable\dott
  }
  In the following for every $i\in\{0,1,\dots,N-1\}$ let
  $\theta_i\colon [a,b]\to[0,H]$
  satisfy for all $t\in [a,b]$ that
  \begin{eqn}
    \llabel{setting8}
    \theta_i(t)=Z_i(r(t),\omega(t))
  \end{eqn}
  and let $\mathcal O_i\subseteq [a,b]$ satisfy
  \begin{eqn}
    \llabel{setting8.1}
    \mathcal O_i=\{s\in [a,b]\colon0<\theta_i(s)<H\}
    \dott
  \end{eqn}
  \startnewargseq
  \argument{
    \lref{setting4.1};
    the fact that $0<\mathfrak{R}<\rho<R$;
    the fact that $\Omega=\{(s,x)\in(0,R)\times\mathbb S\colon sx\in W\}$;
  }{
    that for all $t\in [a,b]$ it holds that
    \begin{eqn}
      \llabel{eq:time-domain}
      r(t)\omega(t)=y(t)\in W
      \qquad\text{and}\qquad
      (r(t),\omega(t))\in\Omega\cap((0,\rho)\times\mathbb S)
      \dott
    \end{eqn}
  }
  \argument{
    \lref{eq3};
    \lref{setting8};
    \lref{setting8.1};
    \lref{eq:time-domain};
  }{
    that for all $i\in\{0,1,\dots,N-1\}$ it holds that
    \begin{eqn}[a]
      \llabel{eq:time-preimage}
      \mathcal O_i
        =\{t\in [a,b]\colon0<Z_i(r(t),\omega(t))<H\}=\{t\in [a,b]\colon(r(t),\omega(t))\in\mathbf{\Omega}_i\}
      \dott
    \end{eqn}
  }
  \argument{
    \lref{setting2};
    \lref{setting5};
    \lref{eq:time-domain};
    \lref{eq:time-preimage};
  }{
    that for all $i\in\{0,1,\dots,N-1\}$ it holds that
    \begin{eqn}[a]
      \llabel{eq:time-strip}
      \mathcal O_i
        &=\{t\in [a,b]\colon
          \mathfrak e_i(r(t))<h(r(t),\omega(t))<\mathfrak e_{i+1}(r(t))\}\\
        &=\{t\in [a,b]\colon
          \mathfrak e_i(r(t))<\mathbf{h}(t)<\mathfrak e_{i+1}(r(t))\}
      \dott
    \end{eqn}
  }
  \argument{
    \lref{eq2};
    \lref{eq:time-strip};
  }{
    that for all $i,j\in\{0,1,\dots,N-1\}$, $t\in\mathcal O_i$ with $i<j$
    it holds that
    $\mathbf{h}(t)<\mathfrak e_{i+1}(r(t))\leq\mathfrak e_j(r(t))$
    and $t\notin\mathcal O_j$\dott
    \llabel{arg:ordered-time-strips}
  }
  \argument{
    \lref{arg:ordered-time-strips};
  }{
    that for all
    $i,j\in\{0,1,\dots,N-1\}$ with $i\neq j$ it holds that
    \begin{eqn}
      \llabel{eq:time-strips-disjoint}
      \mathcal O_i\cap\mathcal O_j=\emptyset
      \dott
    \end{eqn}
  }
  \argument{
    \lref{eq3};
    \lref{eq5};
    \lref{setting3};
    \lref{setting4.2};
    \lref{setting6};
    \lref{eq12};
    \lref{eq13};
    the triangle inequality;
    the Cauchy--Schwarz inequality;
  }{
    that for all $i\in\{0,1,\dots,N-1\}$ it holds that for Lebesgue-almost all
    $t\in\mathcal O_i$ it holds that
    \begin{eqn}[a]
      \llabel{eq17}
      &\abs{\mathbf{v}(t)-c_i(r(t),\theta_i(t))}
        \leq\abs{(\partial_s h)(r(t),\omega(t))-c_i(r(t),\theta_i(t))}
          +\frac{\abs{\langle\xi(t),U(t)\rangle}}{r(t)}\\
      &\leq\abs{(\partial_s h)(r(t),\omega(t))-c_i(r(t),\theta_i(t))}
          +\frac{\norm{\xi(t)}}{r(t)}\norm{U(t)}\\
      &\leq\eta_{i+1}(r(t))\bigl(T_{i+1}(r(t),\theta_i(t))+K^{-1/2}\norm{U(t)}\bigr)
          +\mathbf{k}(r(t))\norm{U(t)}\\
      &\leq\eta_{i+1}(r(t))T_{i+1}(r(t),\theta_i(t))
          +\bigl(\eta_{i+1}(r(t))+\sqrt K\,\mathbf{k}(r(t))\bigr)\gamma(t)\\
      &\leq\mathfrak n_i(r(t))T_{i+1}(r(t),\theta_i(t))
          +\mathfrak n_i(r(t))\gamma(t)
          =\mathfrak n_i(r(t))\bigl(T_{i+1}(r(t),\theta_i(t))+\gamma(t)\bigr)
      \dott
    \end{eqn}
  }
  \argument{
    \lref{atlas5};
    \lref{atlas6};
    \lref{setting2};
    \lref{eq3};
    \lref{eq13};
  }{
    that for all $i\in\{0,1,\dots,N-1\}$ it holds that for Lebesgue-almost all
    $t\in\mathcal O_i$ it holds that
    $(r(t),\omega(t))\in\mathbf{\Omega}_i$, that
    $\mathfrak e_i(r(t))<h(r(t),\omega(t))<\mathfrak e_{i+1}(r(t))$,
    and that
    \begin{eqn}[a]
      \llabel{eq:time-inverse-bound}
      0&<T_{i+1}(r(t),\theta_i(t))
        =T_{i+1}\bigl(r(t),Z_i(r(t),\omega(t))\bigr)\\
      &=T_{i+1}\bigl(r(t),\mathcal F_i\bigl(r(t),h(r(t),\omega(t))\bigr)\bigr)
        \leq K^{-1/2}G(r(t),\omega(t))=K^{-1/2}\norm{U(t)}\leq\gamma(t)
      \dott
    \end{eqn}
  }
  \argument{
    \lref{eq3};
    \lref{eq4};
    \lref{setting5};
    \lref{setting8};
    \lref{setting8.1};
  }{
    that for all $i\in\{0,1,\dots,N-1\}$, $t\in\mathcal O_i$ it holds that
    \begin{eqn}
      \llabel{eq:time-energy-identity}
      \mathbf{h}(t)=h(r(t),\omega(t))
        =E_{i+1}\bigl(r(t),Z_i(r(t),\omega(t))\bigr)
        =E_{i+1}(r(t),\theta_i(t))
      \dott
    \end{eqn}
  }
  \argument{
    \lref{eq:time-energy-identity};
    \lref{atlas5};
    \lref{setting2.2};
    \lref{eq14};
    \lref{eq15};
    the fact that for all $i\in\{0,1,\dots,N-1\}$ it holds that
    $Z_i|_{\mathbf{\Omega}_i}\in\COne{\mathbf{\Omega}_i}{[0,H]}$ and
    $E_{i+1}\in C^2((0,\rho)\times(0,H),\R)$;
    the fact that $r$ and $\omega$ are absolutely continuous;
    the chain rule;
  }{
    that for all $i\in\{0,1,\dots,N-1\}$ it holds that for Lebesgue-almost all
    $t\in\mathcal O_i$ it holds that $\theta_i$ is differentiable at $t$ and that
    \begin{eqn}[a]
      \llabel{eq18}
      &T_{i+1}(r(t),\theta_i(t))(\theta_i)'(t)
        =(\partial_z E_{i+1})(r(t),\theta_i(t))(\theta_i)'(t)\\
      &=\frac{\dd}{\dd t}E_{i+1}(r(t),\theta_i(t))
        -(\partial_s E_{i+1})(r(t),\theta_i(t))r'(t)\\
      &=\mathbf{h}'(t)-c_i(r(t),\theta_i(t))r'(t)
        =\bigl[\mathbf{v}(t)-c_i(r(t),\theta_i(t))\bigr]r'(t)
          -\alpha(t)(\gamma(t))^2\\
      &=-\alpha(t)(r(t))^2\mathbf{v}(t)
          \bigl[\mathbf{v}(t)-c_i(r(t),\theta_i(t))\bigr]
          -\alpha(t)(\gamma(t))^2\\
      &=-\alpha(t)\bigl[(\gamma(t))^2
          +(r(t))^2\mathbf{v}(t)[\mathbf{v}(t)-c_i(r(t),\theta_i(t))]\bigr]
      \dott
    \end{eqn}
  }
  \argument{
    \lref{eq14};
    the fact that for all $t\in [a,b]$ it holds that $\alpha(t)\geq0$;
  }{
    that for Lebesgue-almost all $t\in [a,b]$ it holds that
    \begin{eqn}
      \llabel{eq:radial-absolute-speed}
      \abs{r'(t)}=\abs{-\alpha(t)(r(t))^2\mathbf{v}(t)}
        =\alpha(t)(r(t))^2\abs{\mathbf{v}(t)}
      \dott
    \end{eqn}
  }
  \argument{
    \lref{setting3};
    \lref{eq6};
    \lref{eq:radial-absolute-speed};
    \lref{eq17};
    \lref{eq:time-inverse-bound};
    \lref{eq18};
    the fact that for all $t\in [a,b]$ it holds that $\alpha(t)\geq0$;
    \cref{lem:attenuation} (applied with $r\curvearrowleft r(t)$,
    $T\curvearrowleft T_{i+1}(r(t),\theta_i(t))$, $g\curvearrowleft\gamma(t)$,
    $n\curvearrowleft\mathfrak n_i(r(t))$, $v\curvearrowleft \mathbf{v}(t)$,
    $c\curvearrowleft c_i(r(t),\theta_i(t))$
    in the notation of \cref{lem:attenuation});
  }{
    that for all $i\in\{0,1,\dots,N-1\}$ it holds that for Lebesgue-almost all
    $t\in\mathcal O_i$ it holds that
    \begin{eqn}[a]
      \llabel{eq19}
      &\mathbf{w}_i(r(t),\theta_i(t))(\theta_i)'(t)
        =\frac{\mathbf{w}_i(r(t),\theta_i(t))}{T_{i+1}(r(t),\theta_i(t))}
          T_{i+1}(r(t),\theta_i(t))(\theta_i)'(t)\\
      &=-\alpha(t)\frac{(\gamma(t))^2
          +(r(t))^2\mathbf{v}(t)\bigl[\mathbf{v}(t)-c_i(r(t),\theta_i(t))\bigr]}
        {T_{i+1}(r(t),\theta_i(t))+8(r(t))^2\mathfrak n_i(r(t))
          \sqrt{(c_i(r(t),\theta_i(t)))^2+(T_{i+1}(r(t),\theta_i(t)))^2}}\\
      &\leq\alpha(t)\min\Bigl\{2\mathfrak n_i(r(t))(r(t))^2\abs{\mathbf{v}(t)},
        -\frac{\gamma(t)}4
          +5\mathfrak n_i(r(t))(r(t))^2\abs{\mathbf{v}(t)}\Bigr\}\\
      &=\min\Bigl\{2\mathfrak n_i(r(t))\abs{r'(t)},
        -\frac{\alpha(t)\gamma(t)}4
          +5\mathfrak n_i(r(t))\abs{r'(t)}\Bigr\}
      \dott
    \end{eqn}
  }
  In the following for every $i\in\{0,1,\dots,N-1\}$, $\delta\in(0,H/2)$ let
  $\mathbf{P}_i^\delta\colon [a,b]\to\R$ satisfy for all $t\in [a,b]$ that
  \begin{eqn}
    \llabel{setting10}
    \mathbf{P}_i^\delta(t)
      =P_i\bigl(r(t),\max\{\delta,\min\{\theta_i(t),H-\delta\}\}\bigr),
  \end{eqn}
  for every $\delta\in(0,H/2)$ let $\mathbf{q}^\delta\colon [a,b]\to\R$ and
  $\mathcal U^\delta\subseteq [a,b]$ satisfy for all $t\in [a,b]$ that
  \begin{eqn}
    \llabel{setting10.2}
    \mathbf{q}^\delta(t)=4\sqrt K\sum_{i=0}^{N-1}\mathbf{P}_i^\delta(t)
    \qquad\text{and}\qquad
    \mathcal U^\delta=\bigcup_{i=0}^{N-1}
      \{s\in [a,b]\colon\delta<\theta_i(s)<H-\delta\},
  \end{eqn}
  for every $j\in\{0,1,\dots,N\}$ let $\Gamma_j\subseteq [a,b]$ satisfy
  \begin{eqn}
    \llabel{setting10.4}
    \Gamma_j=\{s\in [a,b]\colon \mathbf{h}(s)=\mathfrak e_j(r(s))\},
  \end{eqn}
  and let $\mathcal C\subseteq [a,b]$ satisfy
  \begin{eqn}
    \llabel{setting10.5}
    \mathcal C=\bigcup_{j=0}^N\Gamma_j
    \dott
  \end{eqn}
  \startnewargseq
  \argument{
    \lref{eq3};
    \lref{setting3};
    \lref{eq7};
    \lref{setting8};
    \lref{setting8.1};
    \lref{eq19};
    \lref{setting10};
    the fact that for all $i\in\{0,1,\dots,N-1\}$ it holds that
    $Z_i$, $\mathbf{w}_i$, $\mathfrak n_i$, and $\mathfrak b_i$ are continuous
    and that $Z_i|_{\mathbf{\Omega}_i}\in\COne{\mathbf{\Omega}_i}{[0,H]}$;
    the fact that $r$ and $\omega$ are absolutely continuous;
    the fact that for all $t\in [a,b]$ it holds that $(r(t),\omega(t))\in\Omega$;
    the fact that $\alpha$ and $\gamma$ are Lebesgue measurable;
    the fact that $\alpha\gamma$ is integrable;
    \itref{lem:truncated}{item2,item3,item4}
    (applied for every $i\in\{0,1,\dots,N-1\}$,
    $\delta\in(0,H/2)$ with
    $d\curvearrowleft n$, $M\curvearrowleft\mathbb S$,
    $\varepsilon\curvearrowleft\delta$,
    $\mathcal O\curvearrowleft\Omega\cap((0,\mathfrak{R})\times\mathbb S)$,
    $Z\curvearrowleft Z_i|_{\Omega\cap((0,\mathfrak{R})\times\mathbb S)}$,
    $\Omega\curvearrowleft\mathbf{\Omega}_i\cap((0,\mathfrak{R})\times\mathbb S)$,
    $\omega\curvearrowleft\mathbf{w}_i|_{(0,\mathfrak{R})\times(0,H)}$,
    $\mathbf m\curvearrowleft\mathfrak n_i|_{(0,\mathfrak{R})}$,
    $\mathfrak b\curvearrowleft\mathfrak b_i|_{(0,\mathfrak{R})}$,
    $P\curvearrowleft P_i|_{(0,\mathfrak{R})\times[0,H]}$,
    $J\curvearrowleft[a,b]$,
    $\varrho\curvearrowleft r$, $\xi\curvearrowleft\omega$,
    $v\curvearrowleft\theta_i$,
    $Q\curvearrowleft\mathbf{P}_i^\delta$);
  }{
    that for all $i\in\{0,1,\dots,N-1\}$, $\delta\in(0,H/2)$ we have that
    \begin{enumerate}[label=(\Alph*)]
      \item \llabel{item:truncated-local-absolute-continuity}
      it holds that $\mathbf{P}_i^\delta$ is locally absolutely continuous on $[a,b]$ and
      \item it holds for Lebesgue-almost every $t\in [a,b]$ that
      \begin{eqn}
        \llabel{eq20}
        (\mathbf{P}_i^\delta)'(t)
          \leq-\frac{\alpha(t)\gamma(t)}4
            \indicator{\{s\in [a,b]\colon\delta<\theta_i(s)<H-\delta\}}(t)
            +\bigl(5\mathfrak n_i(r(t))+\mathfrak b_i(r(t))\bigr)\abs{r'(t)}
        \dott
      \end{eqn}
    \end{enumerate}
  }
  \argument{
    \lref{setting10};
    \lref{item:truncated-local-absolute-continuity};
  }{
    that for all $i\in\{0,1,\dots,N-1\}$, $\delta\in(0,H/2)$
    it holds that $\mathbf{P}_i^\delta$ is absolutely continuous\dott
    \llabel{arg:truncated-absolute-continuity}
  }
  \argument{
    \lref{setting10.2};
    \lref{arg:truncated-absolute-continuity};
  }{
    that for all $\delta\in(0,H/2)$ it holds that $\mathbf{q}^\delta$ is absolutely
    continuous and that for Lebesgue-almost all $t\in [a,b]$ it holds that
    \begin{eqn}
      \llabel{eq:truncated-sum-derivative}
      (\mathbf{q}^\delta)'(t)=4\sqrt K\sum_{i=0}^{N-1}(\mathbf{P}_i^\delta)'(t)
      \dott
    \end{eqn}
  }
  \argument{
    \lref{setting8.1};
     \lref{eq:time-strips-disjoint};
    \lref{setting10.2};
  }{
    that for all $\delta\in(0,H/2)$, $t\in [a,b]$ it holds that
    \begin{eqn}
      \llabel{eq:truncated-indicators}
      \sum_{i=0}^{N-1}
        \indicator{\{s\in [a,b]\colon\delta<\theta_i(s)<H-\delta\}}(t)
        =\indicator{\mathcal U^\delta}(t)
      \dott
    \end{eqn}
  }
  \argument{
    \lref{setting4};
    \lref{eq16};
    \lref{eq20};
    \lref{eq:truncated-sum-derivative};
    \lref{eq:truncated-indicators};
  }{
    that for all $\delta\in(0,H/2)$ it holds that for
    Lebesgue-almost all $t\in [a,b]$ it holds that
    \begin{eqn}[a]
      \llabel{eq21}
      &(\mathbf{q}^\delta)'(t)\leq4\sqrt K\sum_{i=0}^{N-1}\biggl[
          -\frac{\alpha(t)\gamma(t)}4
            \indicator{\{s\in [a,b]\colon\delta<\theta_i(s)<H-\delta\}}(t)
          +\bigl(5\mathfrak n_i(r(t))+\mathfrak b_i(r(t))\bigr)\abs{r'(t)}\biggr]\\
        &=-\sqrt K\,\alpha(t)\gamma(t)\indicator{\mathcal U^\delta}(t)
          +\bigl(\mathfrak m(r(t))-\mathbf{k}(r(t))\bigr)\abs{r'(t)}\\
        &\leq\bigl(-\ell(t)+\mathbf{k}(r(t))\abs{r'(t)}\bigr)
            \indicator{\mathcal U^\delta}(t)
          +\bigl(\mathfrak m(r(t))-\mathbf{k}(r(t))\bigr)\abs{r'(t)}\\
        &=-\ell(t)\indicator{\mathcal U^\delta}(t)
          +\mathfrak m(r(t))\abs{r'(t)}
          -\mathbf{k}(r(t))\bigl(1-\indicator{\mathcal U^\delta}(t)\bigr)\abs{r'(t)}\leq-\ell(t)\indicator{\mathcal U^\delta}(t)
          +\mathfrak m(r(t))\abs{r'(t)}
      \dott
    \end{eqn}
  }
  \argument{
    \lref{setting10.2};
  }{
    that for all $\kappa,\delta\in(0,H/2)$ with $\kappa\leq\delta$ it holds that
    \begin{eqn}
      \llabel{eq:truncated-set-monotonicity}
      \mathcal U^\delta\subseteq\mathcal U^\kappa
      \dott
    \end{eqn}
  }
  \argument{
    \lref{setting1};
    \lref{bounded-potential};
    \lref{setting5};
    \lref{eq:time-domain};
  }{
    that for all $t\in [a,b]$ it holds that
    \begin{eqn}
      \llabel{eq:time-potential-range}
      \mathbf{h}(t)=h(r(t),\omega(t))=f(r(t)\omega(t))=f(y(t))\in(-1,1)
      \dott
    \end{eqn}
  }
  \argument{
    \lref{setting10.4};
    \lref{setting10.5};
  }{
    that
    \begin{eqn}
      \llabel{eq:contact-times}
      \mathcal C=\{t\in [a,b]\colon
        \mathbf{h}(t)\in\{\mathfrak e_0(r(t)),\dots,\mathfrak e_N(r(t))\}\}
      \dott
    \end{eqn}
  }
  \argument{
    \lref{eq2};
    \lref{eq:time-strip};
    \lref{eq:time-potential-range};
    \lref{eq:contact-times};
  }{
    that
    \begin{eqn}[a]
      \llabel{eq:time-noncontact-union}
      \bigcup_{i=0}^{N-1}\mathcal O_i
        &=\Bigl\{t\in [a,b]\colon\mathbf{h}(t)\in
          \bigcup_{i=0}^{N-1}(\mathfrak e_i(r(t)),\mathfrak e_{i+1}(r(t)))\Bigr\}\\
        &=\{t\in [a,b]\colon
          \mathbf{h}(t)\notin\{\mathfrak e_0(r(t)),\dots,\mathfrak e_N(r(t))\}\}
          =[a,b]\backslash\mathcal C
      \dott
    \end{eqn}
  }
  \argument{
    \lref{setting8.1};
    \lref{setting10.2};
    \lref{eq:time-noncontact-union};
  }{
    that
    \begin{eqn}[a]
      \llabel{eq:truncated-set-exhaustion}
      \bigcup_{\delta\in(0,H/2)}\mathcal U^\delta
        &=\bigcup_{i=0}^{N-1}\{t\in [a,b]\colon0<\theta_i(t)<H\}=\bigcup_{i=0}^{N-1}\mathcal O_i=[a,b]\backslash\mathcal C
      \dott
    \end{eqn}
  }
  \argument{
    \lref{eq:state-definition};
    \lref{setting4.1};
    \lref{setting8};
  }{
    that for all $t\in [a,b]$ it holds that
    \begin{eqn}[a]
      \llabel{eq:time-state-sum}
      q(t)&=Q(r(t),\omega(t))
        =4\sqrt K\sum_{i=0}^{N-1}P_i\bigl(r(t),Z_i(r(t),\omega(t))\bigr)=4\sqrt K\sum_{i=0}^{N-1}P_i(r(t),\theta_i(t))
      \dott
    \end{eqn}
  }
  \argument{
    \lref{eq7};
    \lref{setting8};
    \lref{setting10};
  }{
    that for all $i\in\{0,1,\dots,N-1\}$, $\delta\in(0,H/2)$, $t\in [a,b]$
    it holds that
    \begin{eqn}[a]
      \llabel{eq:truncated-primitive-error}
      \abs{\mathbf{P}_i^\delta(t)-P_i(r(t),\theta_i(t))}
        &=\abs[\Big]{P_i\bigl(r(t),\max\{\delta,\min\{\theta_i(t),H-\delta\}\}\bigr)
          -P_i(r(t),\theta_i(t))}\\
        &\leq\abs[\big]{\max\{\delta,\min\{\theta_i(t),H-\delta\}\}-\theta_i(t)}
          \leq\delta
      \dott
    \end{eqn}
  }
  \argument{
    \lref{setting10.2};
    \lref{eq:time-state-sum};
    \lref{eq:truncated-primitive-error};
    the triangle inequality;
  }{
    that for all $\delta\in(0,H/2)$, $t\in [a,b]$ it holds that
    \begin{eqn}[a]
      \llabel{eq22}
      \abs{\mathbf{q}^\delta(t)-q(t)}
        &=4\sqrt K\abs[\bigg]{\sum_{i=0}^{N-1}
          \bigl(\mathbf{P}_i^\delta(t)-P_i(r(t),\theta_i(t))\bigr)}\\
        &\leq4\sqrt K\sum_{i=0}^{N-1}
          \abs{\mathbf{P}_i^\delta(t)-P_i(r(t),\theta_i(t))}
          \leq4\sqrt K\sum_{i=0}^{N-1}\delta=4\sqrt K\,N\delta
      \dott
    \end{eqn}
  }
  \argument{\lref{eq:truncated-sum-derivative};
    \lref{eq21};
    the fundamental theorem of calculus;
  }{
    that for all $\delta\in(0,H/2)$, $\zeta,\chi\in [a,b]$ with $\zeta<\chi$
    it holds that
    \begin{eqn}[a]
      \llabel{eq:truncated-integrated-bound}
      \mathbf{q}^\delta(\chi)-\mathbf{q}^\delta(\zeta)
        &=\int_\zeta^\chi(\mathbf{q}^\delta)'(t)\dd t\leq-\int_\zeta^\chi\ell(t)\indicator{\mathcal U^\delta}(t)\dd t
          +\int_\zeta^\chi\mathfrak m(r(t))\abs{r'(t)}\dd t\\
        &=-\int_{[\zeta,\chi]\cap\mathcal U^\delta}\ell(t)\dd t
          +\int_\zeta^\chi\mathfrak m(r(t))\abs{r'(t)}\dd t
      \dott
    \end{eqn}
  }
  \argument{
    \lref{eq:truncated-set-monotonicity};
    \lref{eq:truncated-set-exhaustion};
    the fact that $\ell$ is nonnegative and integrable;
    the monotone convergence theorem;
  }{
    that for all $\zeta,\chi\in [a,b]$ with $\zeta<\chi$ it holds that
    \begin{eqn}
      \llabel{eq:truncated-integral-limit}
      \limsup_{\delta\searrow0}\biggl|
        \int_{[\zeta,\chi]\cap\mathcal U^\delta}\ell(t)\dd t
        -\int_{[\zeta,\chi]\backslash\mathcal C}\ell(t)\dd t\biggr|=0
      \dott
    \end{eqn}
  }
  \argument{
    \lref{eq22};
    \lref{eq:truncated-integrated-bound};
    \lref{eq:truncated-integral-limit};
  }{
    that for all $\zeta,\chi\in [a,b]$ with $\zeta<\chi$ it holds that
    \begin{eqn}[a]
      \llabel{eq23}
      q(\chi)-q(\zeta)
        =\lim_{\delta\searrow0}
          \bigl[\mathbf{q}^\delta(\chi)-\mathbf{q}^\delta(\zeta)\bigr]\leq-\int_{[\zeta,\chi]\backslash\mathcal C}\ell(t)\dd t
          +\int_\zeta^\chi\mathfrak m(r(t))\abs{r'(t)}\dd t
      \dott
    \end{eqn}
  }
  \argument{
    \lref{eq23};
    the fact that $Q$, $r$, and $\omega$ are continuous;
    the fact that $r$ is absolutely continuous;
    the fact that $\mathfrak m$ is bounded on $r([a,b])$;
  }{
    that
    \begin{eqn}[]
      [a,b]\ni t\longmapsto q(t)-\int_a^t\mathfrak m(r(s))\abs{r'(s)}\dd s
    \end{eqn}
    is continuous and non-increasing and that $q$ has finite length\cfadd{definition:curve-length}\dott
  }
  In the following for every $j\in\{0,1,\dots,N\}$ let $\Psi_j\colon(0,\mathfrak{R})\to\R$ satisfy for all $s\in(0,\mathfrak{R})$ that
  \begin{eqn}
    \llabel{setting11}
    \Psi_j(s)=4\sqrt K\sum_{i=0}^{j-1}P_i(s,H)
    \dott
  \end{eqn}
  \startnewargseq
  \argument{
    \lref{setting10.4};
    the fact that $\mathbf{h}$ and $r$ are continuous;
    the fact that for all $j\in\{0,1,\dots,N\}$ it holds that
    $\mathfrak e_j$ is continuous;
  }{
    that for all $j\in\{0,1,\dots,N\}$ it holds that
    $\Gamma_j$ is closed in $[a,b]$\dott
  }
  \argument{
    \lref{eq2};
    \lref{setting10.4};
  }{
    that for all $j,k\in\{0,1,\dots,N\}$ with $j\neq k$ it holds that
    \begin{eqn}
      \llabel{eq:contact-sets-disjoint}
      \Gamma_j\cap\Gamma_k=\emptyset
      \dott
    \end{eqn}
  }
  \argument{
    \lref{setting5};
    \lref{setting10.4};
  }{
    that for all $j\in\{0,1,\dots,N\}$, $t\in\Gamma_j$ it holds that
    \begin{eqn}
      \llabel{eq:contact-level}
      h(r(t),\omega(t))=\mathbf{h}(t)=\mathfrak e_j(r(t))
      \dott
    \end{eqn}
  }
  \argument{
    \lref{eq2};
    \lref{eq:contact-level};
  }{
    that for all $j\in\{0,1,\dots,N\}$, $i\in\{0,1,\dots,N-1\}$,
    $t\in\Gamma_j$ with $i<j$ it holds that
    \begin{eqn}
      \llabel{eq:contact-order-upper}
      h(r(t),\omega(t))=\mathfrak e_j(r(t))\geq\mathfrak e_{i+1}(r(t))
      \dott
    \end{eqn}
  }
  \argument{
    \lref{eq2};
    \lref{eq:contact-level};
  }{
    that for all $j\in\{0,1,\dots,N\}$, $i\in\{0,1,\dots,N-1\}$,
    $t\in\Gamma_j$ with $i\geq j$ it holds that
    \begin{eqn}
      \llabel{eq:contact-order-lower}
      h(r(t),\omega(t))=\mathfrak e_j(r(t))\leq\mathfrak e_i(r(t))
      \dott
    \end{eqn}
  }
  \argument{
    \lref{eq:coordinate-saturation};
    \lref{setting8};
    \lref{eq:contact-order-upper};
    \lref{eq:contact-order-lower};
  }{
    that for all $j\in\{0,1,\dots,N\}$, $i\in\{0,1,\dots,N-1\}$,
    $t\in\Gamma_j$ it holds that
    \begin{eqn}
      \llabel{eq:contact-coordinates}
      \theta_i(t)=Z_i(r(t),\omega(t))=\begin{cases}
        H&\colon i<j,\\
        0&\colon i\geq j
      \end{cases}
      \dott
    \end{eqn}
  }
  \argument{
    \lref{eq7}; \lref{eq:time-state-sum};
    \lref{setting11};
    \lref{eq:contact-coordinates};
  }{
    that for all $j\in\{0,1,\dots,N\}$, $t\in\Gamma_j$ it holds that
    \begin{eqn}[a]
      \llabel{eq24}
      q(t)&=4\sqrt K\sum_{i=0}^{N-1}P_i(r(t),\theta_i(t))=4\sqrt K\Biggl[\sum_{i=0}^{j-1}P_i(r(t),H)
          +\sum_{i=j}^{N-1}P_i(r(t),0)\Biggr]\\
        &=4\sqrt K\sum_{i=0}^{j-1}P_i(r(t),H)=\Psi_j(r(t))
      \dott
    \end{eqn}
  }
  \argument{
    \lref{eq2};
    \lref{setting2};
    \lref{eq:coordinate-saturation};
  }{
    that for all $i,k\in\{0,1,\dots,N-1\}$,
    $(s,x)\in\mathbf{\Omega}_i\cap((0,\mathfrak{R})\times\mathbb S)$ with $k\neq i$
    it holds that
    \begin{eqn}
      \llabel{eq:strip-saturated-coordinates}
      Z_k(s,x)=\begin{cases}
        H&\colon k<i,\\
        0&\colon k>i
      \end{cases}
      \dott
    \end{eqn}
  }
  \argument{
    \lref{eq7};
    \lref{eq:state-definition};
    \lref{setting11};
    \lref{eq:strip-saturated-coordinates};
  }{
    that for all $i\in\{0,1,\dots,N-1\}$,
    $(s,x)\in\mathbf{\Omega}_i\cap((0,\mathfrak{R})\times\mathbb S)$ it holds that
    \begin{eqn}[a]
      \llabel{eq25}
      Q(s,x)&=4\sqrt K\sum_{k=0}^{N-1}P_k(s,Z_k(s,x))=4\sqrt K\Biggl[\sum_{k=0}^{i-1}P_k(s,H)
          +P_i(s,Z_i(s,x))+\sum_{k=i+1}^{N-1}P_k(s,0)\Biggr]\\
        &=4\sqrt K\sum_{k=0}^{i-1}P_k(s,H)+4\sqrt K\,P_i(s,Z_i(s,x))=\Psi_i(s)+4\sqrt K\,P_i(s,Z_i(s,x))
      \dott
    \end{eqn}
  }
  \argument{
    \lref{setting4.1};
    \lref{setting8};
    \lref{eq:time-preimage};
    \lref{eq25};
  }{
    that for all $i\in\{0,1,\dots,N-1\}$, $t\in\mathcal O_i$
    it holds that
    \begin{eqn}[a]
      \llabel{eq:time-strip-state}
      q(t)&=Q(r(t),\omega(t))
        =\Psi_i(r(t))+4\sqrt K\,P_i\bigl(r(t),Z_i(r(t),\omega(t))\bigr)\\
        &=\Psi_i(r(t))+4\sqrt K\,P_i(r(t),\theta_i(t))
      \dott
    \end{eqn}
  }
  \argument{
    \lref{setting8.1};
    the fact that for all $i\in\{0,1,\dots,N-1\}$ it holds that
    $\theta_i$ is continuous;
  }{
    that for all $i\in\{0,1,\dots,N-1\}$, $\zeta,\chi\in [a,b]$
    with $\zeta<\chi$ and $[\zeta,\chi]\subseteq\mathcal O_i$
 there exists $\delta\in(0,H/2)$ such that for all
    $t\in[\zeta,\chi]$ it holds that
    \begin{eqn}
      \llabel{eq:strip-truncation-margin}
      \delta<\theta_i(t)<H-\delta
      \dott
    \end{eqn}
  }
  \argument{
    \lref{setting10};
    \lref{eq:time-strip-state};
    \lref{eq:strip-truncation-margin};
  }{
    that for all $i\in\{0,1,\dots,N-1\}$, $\zeta,\chi\in [a,b]$
    with $\zeta<\chi$ and $[\zeta,\chi]\subseteq\mathcal O_i$
    it holds that there exists $\delta\in(0,H/2)$ such that for all
    $t\in[\zeta,\chi]$ it holds that
    \begin{eqn}
      \llabel{eq:time-strip-truncated-state}
      q(t)=\Psi_i(r(t))+4\sqrt K\,P_i(r(t),\theta_i(t))
        =\Psi_i(r(t))+4\sqrt K\,\mathbf{P}_i^\delta(t)
      \dott
    \end{eqn}
  }
  \argument{
    \lref{eq:time-strips-disjoint};
    \lref{eq:time-noncontact-union};
    the fact that for all $i\in\{0,1,\dots,N-1\}$ it holds that
    $\mathcal O_i$ is open in $[a,b]$;
  }{
    that for all $\zeta,\chi\in [a,b]$ with $\zeta<\chi$ and
    $[\zeta,\chi]\cap\mathcal C=\emptyset$
    there exists $i\in\{0,1,\dots,N-1\}$ such that
    \begin{eqn}
      \llabel{eq:noncontact-interval-strip}
      [\zeta,\chi]\subseteq\mathcal O_i
      \dott
    \end{eqn}
  }
  \argument{
    \lref{eq7};
    \lref{setting11};
    the fact that for all $i\in\{0,1,\dots,N-1\}$ it holds that
    $P_i|_{(0,\mathfrak{R})\times[0,H]}$ is locally Lipschitz;
    the fact that $r$ is absolutely continuous;
  }{
    that for all $j\in\{0,1,\dots,N\}$ it holds that $\Psi_j$ is
    locally Lipschitz and $\Psi_j\circ r$ is absolutely continuous\dott
  }
  \argument{
    \lref{eq:time-strip-truncated-state};
    \lref{eq:noncontact-interval-strip};
    the fact that for all $i\in\{0,1,\dots,N-1\}$, $\delta\in(0,H/2)$
    it holds that $\mathbf{P}_i^\delta$ is absolutely continuous;
    the fact that for all $i\in\{0,1,\dots,N-1\}$ it holds that
    $\Psi_i\circ r$ is absolutely continuous;
  }{
    that for all $\zeta,\chi\in [a,b]$ with $\zeta<\chi$ and
    $[\zeta,\chi]\cap\mathcal C=\emptyset$ it holds that
    \begin{eqn}
      \llabel{arg1}
      q|_{[\zeta,\chi]}\text{ is absolutely continuous}
      \dott
    \end{eqn}
  }
  \argument{
    \lref{eq24};
    \lref{arg1};
    the fact that $q$ is continuous and has finite length\cfadd{definition:curve-length};
    the fact that for all $j\in\{0,1,\dots,N\}$ it holds that
    $\Gamma_j$ is closed in $[a,b]$;
    the fact that for all $j\in\{0,1,\dots,N\}$ it holds that
    $\Psi_j\circ r$ is absolutely continuous;
    \cref{lem:AC-gluing} (applied with $N\curvearrowleft N+1$,
    $u\curvearrowleft q$, $E_j\curvearrowleft\Gamma_{j-1}$,
    $\mathfrak a_j\curvearrowleft\Psi_{j-1}\circ r$
    in the notation of \cref{lem:AC-gluing});
  }[verbs=ep]{
    \lref{item:compact-absolute-continuity}\dott
  }
  \startnewargseq
  \argument{
    \lref{eq15};
    \lref{setting10.4};
    the fact that for all $j\in\{0,1,\dots,N\}$ it holds that
    $\mathbf{h}-\mathfrak e_j\circ r$ is absolutely continuous;
    \cref{lem:AC-level-set};
  }{
    that for all $j\in\{0,1,\dots,N\}$ it holds that for Lebesgue-almost all
    $t\in\Gamma_j$ it holds that
    \begin{eqn}
      \llabel{eq26}
      0=\mathbf{h}'(t)-(\mathfrak e_j)'(r(t))r'(t)
        =-\alpha(t)(\gamma(t))^2
          +\bigl(\mathbf{v}(t)-(\mathfrak e_j)'(r(t))\bigr)r'(t)
      \dott
    \end{eqn}
  }
  \argument{
    \lref{eq8};
    \lref{setting6};
    \lref{eq12};
    \lref{eq13};
    \lref{eq26};
  }{
    that for all $j\in\{0,1,\dots,N\}$ it holds that for Lebesgue-almost all
    $t\in\Gamma_j$ with $U(t)\neq0$ it holds that
    \begin{eqn}[a]
      \llabel{eq27}
      \alpha(t)(\gamma(t))^2
        &=\bigl(\mathbf{v}(t)-(\mathfrak e_j)'(r(t))\bigr)r'(t)
          =\abs{\mathbf{v}(t)-(\mathfrak e_j)'(r(t))}\abs{r'(t)}\\
        &=\abs[\bigg]{(\partial_s h)(r(t),\omega(t))-(\mathfrak e_j)'(r(t))
          +\frac{\langle\xi(t),U(t)\rangle}{r(t)}}\abs{r'(t)}\\
        &\leq\bigl(\beta_j(r(t))+\mathbf{k}(r(t))\bigr)\norm{U(t)}\abs{r'(t)}\leq\sqrt K\bigl(\beta_j(r(t))+\mathbf{k}(r(t))\bigr)
          \gamma(t)\abs{r'(t)}
      \dott
    \end{eqn}
  }
  \argument{
    \lref{eq13};
    \lref{eq27};
    the fact that for all $t\in [a,b]$ with $U(t)=0$ it holds that
    $\gamma(t)=0$;
  }{
    that for all $j\in\{0,1,\dots,N\}$ it holds that for Lebesgue-almost all
    $t\in\Gamma_j$ it holds that
    \begin{eqn}
      \llabel{eq28}
      \alpha(t)\gamma(t)
        \leq\sqrt K\bigl(\beta_j(r(t))+\mathbf{k}(r(t))\bigr)\abs{r'(t)}
      \dott
    \end{eqn}
  }
  \argument{
    \lref{eq16};
    \lref{eq28};
  }{
    that for all $j\in\{0,1,\dots,N\}$ it holds that for Lebesgue-almost all
    $t\in\Gamma_j$ it holds that
    \begin{eqn}
      \llabel{eq29}
      \ell(t)\leq
        \bigl(K\beta_j(r(t))+(K+1)\mathbf{k}(r(t))\bigr)\abs{r'(t)}
      \dott
    \end{eqn}
  }
  \argument{
    \lref{setting10.5};
    \lref{eq:contact-sets-disjoint};
    \lref{eq29};
  }{
    that for all $\zeta,\chi\in [a,b]$ with $\zeta<\chi$ it holds that
    \begin{eqn}[a]
      \llabel{eq:contact-integral-bound}
      \displaystyle\int_{[\zeta,\chi]\cap\mathcal C}\textstyle\ell(t)\dd t
        &\textstyle=\sum\limits_{j=0}^N\displaystyle\int_{[\zeta,\chi]\cap\Gamma_j}\textstyle\ell(t)\dd t\\
        &\textstyle\leq\sum\limits_{j=0}^N\displaystyle\int_{[\zeta,\chi]\cap\Gamma_j}\textstyle
          \bigl(K\beta_j(r(t))+(K+1)\mathbf{k}(r(t))\bigr)\abs{r'(t)}\dd t\\
        &\textstyle=K\sum\limits_{j=0}^N\displaystyle\int_{[\zeta,\chi]\cap\Gamma_j}\textstyle
          \beta_j(r(t))\abs{r'(t)}\dd t
          +(K+1)\displaystyle\int_{[\zeta,\chi]\cap\mathcal C}\textstyle
            \mathbf{k}(r(t))\abs{r'(t)}\dd t\\
        &\textstyle\leq\displaystyle\int_\zeta^\chi\textstyle\biggl[K\sum\limits_{j=0}^N\beta_j(r(t))
          +(K+1)\mathbf{k}(r(t))\biggr]\abs{r'(t)}\dd t
      \dott
    \end{eqn}
  }
  \argument{
    \lref{setting4};
    \lref{eq23};
    \lref{eq:contact-integral-bound};
  }{
    that for all $\zeta,\chi\in [a,b]$ with $\zeta<\chi$ it holds that
    \begin{eqn}[a]
      \llabel{eq30}
      \textstyle q(\chi)-q(\zeta)
        &\textstyle\leq-\displaystyle\int_{[\zeta,\chi]\backslash\mathcal C}\textstyle\ell(t)\dd t
          +\displaystyle\int_\zeta^\chi\textstyle\mathfrak m(r(t))\abs{r'(t)}\dd t\\
        &\textstyle=-\displaystyle\int_\zeta^\chi\textstyle\ell(t)\dd t
          +\displaystyle\int_{[\zeta,\chi]\cap\mathcal C}\textstyle\ell(t)\dd t
          +\displaystyle\int_\zeta^\chi\textstyle\mathfrak m(r(t))\abs{r'(t)}\dd t\\
        &\textstyle\leq-\displaystyle\int_\zeta^\chi\textstyle\ell(t)\dd t
          +\displaystyle\int_\zeta^\chi\textstyle\biggl[\mathfrak m(r(t))
            +K\sum\limits_{j=0}^N\beta_j(r(t))+(K+1)\mathbf{k}(r(t))\biggr]
            \abs{r'(t)}\dd t\\
        &\textstyle=-\displaystyle\int_\zeta^\chi\textstyle\ell(t)\dd t
          +\displaystyle\int_\zeta^\chi\textstyle m(r(t))\abs{r'(t)}\dd t
      \dott
    \end{eqn}
  }
  \argument{
    \lref{setting7};
    \lref{eq30};
    the fact that $q$, $r$, and $\omega$ are absolutely continuous;
    the fact that $\ell$ is integrable;
    the fact that $m$ is continuous and bounded on $r([a,b])$;
    the Lebesgue differentiation theorem;
  }[verbs=ep]{
    \lref{item:compact-differential-inequality}\dott
  }
  This and the arbitrariness of $a,b,y,\lambda,B,r,\omega,q$
  complete our proof of
  \lref{item:compact-absolute-continuity,item:compact-differential-inequality}\dott
  \startnewargseq

  To prove \lref{item1,item2,item3}, let $\mathbb{I}\subseteq\R$ be an interval,
  let $\mathbf{y}\colon \mathbb{I}\to W$ be locally absolutely continuous,
  let $\boldsymbol{\lambda}\colon \mathbb{I}\to[0,\infty)$ and
  $\mathbf{B}\colon \mathbb{I}\to\R^{n\times n}$ be Lebesgue measurable, and let
  $\mathbf{r}\colon \mathbb{I}\to(0,\mathfrak{R})$, $\boldsymbol{\omega}\colon \mathbb{I}\to\mathbb S$, and
  $\mathfrak q\colon \mathbb{I}\to\R$ satisfy for all $t\in \mathbb{I}$ that
  \begin{eqn}
    \llabel{eq:original-state}
    \mathbf{r}(t)=\norm{\mathbf{y}(t)},\qquad
    \boldsymbol{\omega}(t)=\mathbf{y}(t)/\mathbf{r}(t),\qquad\text{and}\qquad
    \mathfrak q(t)=Q(\mathbf{r}(t),\boldsymbol{\omega}(t)),
  \end{eqn}
  and assume that for Lebesgue-almost all $t\in \mathbb{I}$ it holds that
  $\mathbf{y}$ is differentiable at $t$, that $\mathbf{B}(t)$ is symmetric and positive-definite, and that
  \begin{equation}
    \llabel{eq:original-dynamics}
    \begin{gathered}
      \mathbf{y}'(t)=-\boldsymbol{\lambda}(t)\mathbf{B}(t)\grad F(\mathbf{y}(t)),\qquad
      \max_{z\in\R^n,\,\norm{z}=1}\langle \mathbf{B}(t)z,z\rangle
        \leq K\min_{z\in\R^n,\,\norm{z}=1}\langle \mathbf{B}(t)z,z\rangle,\\
      \text{and}\qquad\frac{\norm{\mathbf{B}(t)\boldsymbol{\omega}(t)
         -\langle \mathbf{B}(t)\boldsymbol{\omega}(t),\boldsymbol{\omega}(t)\rangle\boldsymbol{\omega}(t)}}
        {\langle \mathbf{B}(t)\boldsymbol{\omega}(t),\boldsymbol{\omega}(t)\rangle}
        \leq \mathbf{r}(t)\mathbf{k}(\mathbf{r}(t))
      \dott
    \end{gathered}
  \end{equation}
  In the following let $\kappa\colon \mathbb{I}\to[0,\infty)$ satisfy for all
  $t\in \mathbb{I}$ that
  \begin{eqn}
    \llabel{eq:normalization-coefficient}
    \kappa(t)=\boldsymbol{\lambda}(t)(1+(F(\mathbf{y}(t)))^2)^{3/2}
    \dott
  \end{eqn}
  \startnewargseq
  \argument{
    \lref{eq:normalization-coefficient};
    the fact that $F\circ \mathbf{y}$ is continuous;
    the assumption that $\boldsymbol{\lambda}$ is Lebesgue measurable;
  }{
    that $\kappa$ is Lebesgue measurable\dott
  }
  \argument{
    \lref{bounded-gradient};
    \lref{eq:normalization-coefficient};
  }{
    that for all $t\in \mathbb{I}$ it holds that
    \begin{eqn}
      \llabel{eq:normalization-gradient}
      \kappa(t)\grad f(\mathbf{y}(t))
        =\boldsymbol{\lambda}(t)(1+(F(\mathbf{y}(t)))^2)^{3/2}
          \frac{\grad F(\mathbf{y}(t))}{(1+(F(\mathbf{y}(t)))^2)^{3/2}}
        =\boldsymbol{\lambda}(t)\grad F(\mathbf{y}(t))
      \dott
    \end{eqn}
  }
  \argument{
    \lref{eq:original-dynamics};
    \lref{eq:normalization-gradient};
  }{
    that for Lebesgue-almost all $t\in \mathbb{I}$ it holds that
    \begin{eqn}
      \llabel{eq:normalized-flow}
      \mathbf{y}'(t)=-\boldsymbol{\lambda}(t)\mathbf{B}(t)\grad F(\mathbf{y}(t))
        =-\kappa(t)\mathbf{B}(t)\grad f(\mathbf{y}(t))
      \dott
    \end{eqn}
  }
  \argument{
    \itref{prop:ambient-state}{item:compact-absolute-continuity,item:compact-differential-inequality};
    \lref{eq:original-dynamics};
    \lref{eq:normalized-flow};
  }{
    that for all $a,b\in \mathbb{I}$ with $a<b$ it holds that
    $\mathfrak q|_{[a,b]}$ is absolutely continuous and that for Lebesgue-almost
    all $t\in[a,b]$ it holds that $\mathfrak q$, $\mathbf{r}$, and $\boldsymbol{\omega}$ are
    differentiable at $t$ and that
    \begin{eqn}
      \llabel{eq:original-compact-state-estimate}
      \mathfrak q'(t)\leq-\norm{\boldsymbol{\omega}'(t)}+m(\mathbf{r}(t))\abs{\mathbf{r}'(t)}
      \dott
    \end{eqn}
  }
  \argument{
    \lref{eq:original-compact-state-estimate};
    the arbitrariness of $\mathbb{I},\mathbf{y},\boldsymbol{\lambda},\mathbf{B},\mathbf{r},\boldsymbol{\omega},\mathfrak q$ subject to
    the stated assumptions;
  }[verbs=ep]{
    \lref{item1,item2,item3}\dott
  }
  \startnewargseq
\end{aproof}

\begin{athm}{proposition}{prop:ambient-mixing}
  Let $n\in\N\backslash\{1\}$, $K,R\in(0,\infty)$,
  let
  $V=\{z\in\R^n\colon0<\norm{z}<R\}$,
  let $\mathbb S=\{z\in\R^n\colon\norm{z}=1\}$,
  let $W\subseteq V$ be definable and open,
  let $g\in C^1(W,\R)$ be definable,
  and let $\mathbf{k}\colon(0,R)\to[0,\infty)$ be locally bounded,
  definable, and integrable\dott
  \cfload[.]
  Then there exist $\varrho\in(0,R)$, $\varepsilon\in(0,\varrho)$,
  $\mathfrak{C},L\in(0,\infty)$ such that
  for every open $U\subseteq\R^n$,
  every $F\in C^1(U,\R)$,
  every nonempty interval $I\subseteq\R$,
  every locally absolutely continuous $y\colon I\to U\backslash\{0\}$,
  every Lebesgue measurable $\lambda\colon I\to[0,\infty)$,
  every Lebesgue measurable $B\colon I\to\R^{n\times n}$,
  and every $\omega\colon I\to\mathbb S$
  with the property that $U\cap V\subseteq W$ and
  $F|_{U\cap V}=g|_{U\cap V}$,
  with the property that for all $t\in I$ it holds that
  \begin{eqn}
    \llabel{ass:polar}
    \omega(t)=\norm{y(t)}^{-1}y(t),
  \end{eqn}
  with the property that for Lebesgue-almost all $t\in I$ it holds
  that $y$ is differentiable at $t$, that $B(t)$ is symmetric and
  positive-definite, and that
  \begin{eqn}
    \label{eq:ambient-flow}
    y'(t)=-\lambda(t)B(t)\grad F(y(t)),
  \end{eqn}
  and with the property that for Lebesgue-almost all
  $t\in\{s\in I\colon \norm{y(s)}<R\}$ it holds that
  \begin{equation}
    \label{eq:ambient-mixing}
    [\langle B(t)\omega(t),\omega(t)\rangle]^{-1}
      \norm{B(t)\omega(t)
        -\langle B(t)\omega(t),\omega(t)\rangle\omega(t)}
      \leq \norm{y(t)}\mathbf{k}(\norm{y(t)})
  \end{equation}
  \begin{equation}
    \llabel{ass:condition}
    \text{and}\qquad
    \max\nolimits_{z\in\R^n,\,\norm{z}=1}\langle B(t)z,z\rangle
      \leq K\min\nolimits_{z\in\R^n,\,\norm{z}=1}\langle B(t)z,z\rangle
  \end{equation}
  we have that
  \begin{enumerate}[label=(\roman*)]
    \item \llabel{item:annular} it holds for all $\rho\in(0,\varrho)$ that
    \begin{eqn}
      \label{eq:ambient-annular}
      \int_{\{t\in I\colon\rho/2\leq \norm{y(t)}\leq\rho\}}
        \norm{y'(t)}\dd t\leq\mathfrak{C}\rho
      \qquad\text{and}\qquad
      \int_{\{t\in I\colon \norm{y(t)}\leq\rho\}}
        \norm{y'(t)}\dd t\leq2\mathfrak{C}\rho,
    \end{eqn}
    \item \llabel{item:radial-limit} it holds for all
    $\tau\in\{\inf I,\sup I\}\backslash I$ with
    $\liminf_{t\to\tau}\norm{y(t)}=0$ that
    $\limsup_{t\to\tau}\norm{y(t)}=0$,
    \item \llabel{item:angular} it holds that
    $\int_{\{t\in I\colon\norm{y(t)}\leq\varepsilon\}}
      \norm{\omega'(t)}\dd t\leq L$, and
    \item \llabel{item:direction} it holds for all $c\in I$,
    $\tau\in\{\inf I,\sup I\}\backslash I$ with
    $\liminf_{t\to\tau}\norm{y(t)}=0$ there exists
    $v\in\mathbb S$ such that
    \begin{eqn}
      \llabel{concl:direction}
      \textstyle\int_{I\cap(\min\{c,\tau\},\max\{c,\tau\})}
        \norm{\omega'(t)}\dd t<\infty
      \qquad\text{and}\qquad
      \limsup_{t\to\tau}\norm{\omega(t)-v}=0
      \dott
    \end{eqn}
  \end{enumerate}
  \cfout[.]
\end{athm}
\begin{aproof}
  \startnewargseq
  \argument{
    the fact that $n\in\N\backslash\{1\}$;
    the fact that $K,R\in(0,\infty)$;
    the fact that $V=\{z\in\R^n\colon0<\norm{z}<R\}$;
    the fact that $W\subseteq V$ is definable and open;
    the fact that $g\in C^1(W,\R)$ is definable;
    \cref{lem:ambient-annular};
  }{
    that there exist $\varrho\in(0,R)$, $\mathfrak{C}\in(0,\infty)$
    which satisfy that for every $\rho\in(0,\varrho)$, every open
    $\mathcal O\subseteq\R^n$, every $G\in C^1(\mathcal O,\R)$,
    every nonempty interval $J\subseteq\R$, every locally absolutely
    continuous $z\colon J\to\mathcal O\backslash\{0\}$, every
    Lebesgue measurable $\mu\colon J\to[0,\infty)$,
    and every Lebesgue measurable $A\colon J\to\R^{n\times n}$
    with the property that
    \begin{eqn}
      \llabel{ass:uniform-annular-extension}
      \mathcal O\cap V\subseteq W
      \qquad\text{and}\qquad
      G|_{\mathcal O\cap V}=g|_{\mathcal O\cap V},
    \end{eqn}
    with the property that for Lebesgue-almost all $t\in J$ it holds
    that $z$ is differentiable at $t$, that $A(t)$ is symmetric and
    positive-definite, and that
    \begin{eqn}
      \llabel{ass:uniform-annular-flow}
      z'(t)=-\mu(t)A(t)\grad G(z(t))
    \end{eqn}
    and with the property that for Lebesgue-almost all
    $t\in\{s\in J\colon\norm{z(s)}<R\}$ it holds that
    \begin{eqn}
      \llabel{ass:uniform-annular-condition}
      \max_{u\in\R^n,\,\norm{u}=1}\langle A(t)u,u\rangle
        \leq K\min_{u\in\R^n,\,\norm{u}=1}\langle A(t)u,u\rangle
    \end{eqn}
    we have that
    \begin{eqn}
      \llabel{eq:uniform-annular}
      \int_{\{t\in J\colon\rho/2\leq\norm{z(t)}\leq\rho\}}
        \norm{z'(t)}\dd t\leq\mathfrak{C}\rho
      \qquad\text{and}\qquad
      \int_{\{t\in J\colon\norm{z(t)}\leq\rho\}}
        \norm{z'(t)}\dd t\leq2\mathfrak{C}\rho
      \dott
    \end{eqn}
  }
  \startnewargseq
  \argument{
    the fact that $n\in\N\backslash\{1\}$ and $K,R\in(0,\infty)$;
    the fact that $\mathbb S=\{z\in\R^n\colon\norm{z}=1\}$;
    the fact that $W\subseteq\{z\in\R^n\colon0<\norm{z}<R\}$
    is definable and open;
    the fact that $g\in C^1(W,\R)$ is definable;
    the fact that $\mathbf k$ is nonnegative, locally bounded, definable,
    and integrable;
    \cref{prop:ambient-state} (applied with $F\curvearrowleft g$
    in the notation of \cref{prop:ambient-state});
  }{
    that there exist $\mathfrak R\in(0,R)$, bounded continuous
    $Q\colon\{(s,x)\in(0,\mathfrak R)\times\mathbb S\colon sx\in W\}
      \to[0,\infty)$, and
    continuous definable integrable $m\colon(0,\mathfrak R)\to[0,\infty)$
    which satisfy that for every nonempty interval $J\subseteq\R$, every
    locally absolutely continuous
    $z\colon J\to W$, every
    Lebesgue measurable $\mu\colon J\to[0,\infty)$, every
    Lebesgue measurable $A\colon J\to\R^{n\times n}$, every
    $r\colon J\to(0,\mathfrak R)$, every $\xi\colon J\to\mathbb S$,
    and every $q\colon J\to\R$
    with the property that for all $t\in J$ it holds that
    \begin{eqn}
      \llabel{ass:uniform-state-coordinates}
      r(t)=\norm{z(t)},\qquad
      \xi(t)=\frac{z(t)}{r(t)},
      \qquad\text{and}\qquad
      q(t)=Q(r(t),\xi(t)),
    \end{eqn}
    and with the property that for Lebesgue-almost all $t\in J$ it holds
    that $z$ is differentiable at $t$, that $A(t)$ is symmetric and
    positive-definite, and that
    \begin{gather}
      \llabel{ass:uniform-state-flow}
      z'(t)=-\mu(t)A(t)\grad g(z(t)),\qquad
      \max_{u\in\R^n,\,\norm{u}=1}\langle A(t)u,u\rangle
        \leq K\min_{u\in\R^n,\,\norm{u}=1}\langle A(t)u,u\rangle,\\
      \llabel{ass:uniform-state-mixing}
      \text{and}\qquad
      \frac{\norm{A(t)\xi(t)
        -\langle A(t)\xi(t),\xi(t)\rangle\xi(t)}}
        {\langle A(t)\xi(t),\xi(t)\rangle}
        \leq r(t)\mathbf{k}(r(t))
    \end{gather}
    we have that $q$ is locally absolutely continuous and that for
    Lebesgue-almost all $t\in J$ it holds that $q$, $r$, and $\xi$ are
    differentiable at $t$ and that
    \begin{eqn}
      \llabel{eq:uniform-state}
      q'(t)\leq-\norm{\xi'(t)}+m(r(t))\abs{r'(t)}
      \dott
    \end{eqn}
  }
  \startnewargseq
  \argument{
    \itref{lem:radial-integrability}{item2}
    (applied with $R\curvearrowleft\mathfrak R/2$,
    $m\curvearrowleft m|_{(0,\mathfrak R/2]}$ in the notation of
    \cref{lem:radial-integrability});
  }{
    that there exist $\eta\in(0,\min\{\varrho,\mathfrak R/2\})$
    and a definable, non-increasing, and integrable
    $\psi\colon(0,\eta]\to[0,\infty)$ which satisfy for all
    $s\in(0,\eta]$ that
    \begin{eqn}
      \llabel{eq:majorant}
      m(s)\leq\psi(s)
      \dott
    \end{eqn}
  }
  \startnewargseq
  In the following let $M\in[1,\infty)$ satisfy for all
  $(s,x)\in(0,\mathfrak R)\times\mathbb S$ with $sx\in W$ that
  $Q(s,x)\leq M$,
  let $L\in(0,\infty)$ satisfy
  \begin{eqn}
    \llabel{setting:constants}
    L=M(1+2\mathfrak{C})+2\mathfrak{C}\int_0^\eta\psi(s)\dd s,
  \end{eqn}
  and let $\chi\colon(0,\infty)\to[0,1]$ and
  $\mathfrak m\colon(0,\infty)\to[0,\infty)$ satisfy for all
  $s\in(0,\infty)$ that
  \begin{eqn}
    \llabel{setting:cutoff}
    \chi(s)=\begin{cases}
      1&\colon0<s\leq\eta/2\\
      2-2s/\eta&\colon\eta/2<s<\eta\\
      0&\colon s\geq\eta
    \end{cases}
    \qquad\text{and}\qquad
    \mathfrak m(s)=\begin{cases}
      m(s)&\colon0<s\leq\eta\\
      0&\colon s>\eta
    \end{cases}
    \dott
  \end{eqn}
  \startnewargseq
  We claim that for every open $U\subseteq\R^n$,
  every $F\in C^1(U,\R)$,
  every nonempty interval $I\subseteq\R$,
  every locally absolutely continuous $y\colon I\to U\backslash\{0\}$,
  every Lebesgue measurable $\lambda\colon I\to[0,\infty)$,
  every Lebesgue measurable $B\colon I\to\R^{n\times n}$,
  and every $\omega\colon I\to\mathbb S$
  with the property that $U\cap V\subseteq W$ and
  $F|_{U\cap V}=g|_{U\cap V}$,
  with the property that for all $t\in I$ it holds that
  \begin{eqn}
    \llabel{ass:claim-polar}
    \omega(t)=\frac{y(t)}{\norm{y(t)}},
  \end{eqn}
  with the property that for Lebesgue-almost all $t\in I$ it holds
  that $y$ is differentiable at $t$, that $B(t)$ is symmetric and
  positive-definite, and that
  \begin{eqn}
    \llabel{ass:claim-flow}
    y'(t)=-\lambda(t)B(t)\grad F(y(t)),
  \end{eqn}
  and with the property that for Lebesgue-almost all
  $t\in\{s\in I\colon\norm{y(s)}<R\}$ it holds that
  \begin{gather}
    \llabel{ass:claim-condition}
    \max_{z\in\R^n,\,\norm{z}=1}\langle B(t)z,z\rangle
      \leq K\min_{z\in\R^n,\,\norm{z}=1}\langle B(t)z,z\rangle\\
    \llabel{ass:claim-mixing}
    \text{and}\qquad
    \frac{\norm{B(t)\omega(t)
      -\langle B(t)\omega(t),\omega(t)\rangle\omega(t)}}
      {\langle B(t)\omega(t),\omega(t)\rangle}
      \leq\norm{y(t)}\mathbf{k}\bigl(\norm{y(t)}\bigr)
  \end{gather}
  we have that
  \begin{enumerate}[label=(\Roman*)]
    \item \llabel{claim:annular} for all $\rho\in(0,\varrho)$ it holds that
    \begin{eqn}
      \llabel{eq:claim-annular}
      \int_{\{t\in I\colon\rho/2\leq\norm{y(t)}\leq\rho\}}
        \norm{y'(t)}\dd t\leq\mathfrak{C}\rho
      \qquad\text{and}\qquad
      \int_{\{t\in I\colon\norm{y(t)}\leq\rho\}}
        \norm{y'(t)}\dd t\leq2\mathfrak{C}\rho,
    \end{eqn}
    \item \llabel{claim:radial-limit} for all
    $\tau\in\{\inf I,\sup I\}\backslash I$ with
    $\liminf_{t\to\tau}\norm{y(t)}=0$ it holds that
    $\limsup_{t\to\tau}\norm{y(t)}=0$,
    \item \llabel{claim:angular} it holds that
    $\int_{\{t\in I\colon\norm{y(t)}\leq\eta/2\}}
      \norm{\omega'(t)}\dd t\leq L$,
    and
    \item \llabel{claim:direction} for all $c\in I$,
    $\tau\in\{\inf I,\sup I\}\backslash I$ with
    $\liminf_{t\to\tau}\norm{y(t)}=0$ there exists
    $v\in\mathbb S$ such that
    \begin{eqn}
      \llabel{eq:claim-direction}
      \int_{I\cap(\min\{c,\tau\},\max\{c,\tau\})}
        \norm{\omega'(t)}\dd t<\infty
      \qquad\text{and}\qquad
      \limsup_{t\to\tau}\norm{\omega(t)-v}=0
      \dott
    \end{eqn}
  \end{enumerate}
  \startnewargseq
  To prove \lref{claim:annular,claim:radial-limit,claim:angular,claim:direction}, let
  $U\subseteq\R^n$ be open, let $F\in C^1(U,\R)$,
  assume that
  \begin{eqn}
    \llabel{setting:extension}
    U\cap V\subseteq W
    \qquad\text{and}\qquad
    F|_{U\cap V}=g|_{U\cap V},
  \end{eqn}
  let $I\subseteq\R$ be a nonempty interval,
  let $y\colon I\to U\backslash\{0\}$ be locally absolutely continuous,
  let $\lambda\colon I\to[0,\infty)$ and
  $B\colon I\to\R^{n\times n}$ be Lebesgue measurable,
  let $\omega\colon I\to\mathbb S$ satisfy for all $t\in I$ that
  \begin{eqn}
    \llabel{setting:polar}
    \omega(t)=\frac{y(t)}{\norm{y(t)}},
  \end{eqn}
  assume that for Lebesgue-almost all $t\in I$ it holds that
  $y$ is differentiable at $t$, that $B(t)$ is symmetric and positive-definite, and that
  \begin{eqn}
    \llabel{setting:flow}
    y'(t)=-\lambda(t)B(t)\grad F(y(t)),
  \end{eqn}
  assume that for Lebesgue-almost all
  $t\in\{s\in I\colon\norm{y(s)}<R\}$ it holds that
  \begin{gather}
    \llabel{setting:condition}
    \max_{z\in\R^n,\,\norm{z}=1}\langle B(t)z,z\rangle
      \leq K\min_{z\in\R^n,\,\norm{z}=1}\langle B(t)z,z\rangle\\
    \llabel{setting:mixing}
    \text{and}\qquad
    \frac{\norm{B(t)\omega(t)
      -\langle B(t)\omega(t),\omega(t)\rangle\omega(t)}}
      {\langle B(t)\omega(t),\omega(t)\rangle}
      \leq\norm{y(t)}\mathbf{k}\bigl(\norm{y(t)}\bigr),
  \end{gather}
  and let $r\colon I\to(0,\infty)$ satisfy for all $t\in I$
  that
  \begin{eqn}
    \llabel{setting:radius}
    r(t)=\norm{y(t)}
    \dott
  \end{eqn}
  \startnewargseq
  \argument{
    \lref{setting:polar};
    \lref{setting:radius};
    the fact that $y$ is locally absolutely continuous;
    the chain rule;
  }{
    that $r$ and $\omega$ are locally absolutely continuous and that
    for Lebesgue-almost all $t\in I$ it holds that
    \begin{eqn}
      \llabel{eq:radial-derivative}
      r'(t)=\langle\omega(t),y'(t)\rangle
      \qquad\text{and}\qquad
      \abs{r'(t)}\leq\norm{y'(t)}
      \dott
    \end{eqn}
  }
  \argument{
    \lref{eq:uniform-annular} (applied with $J\curvearrowleft I$,
    $\mathcal O\curvearrowleft U$, $G\curvearrowleft F$,
    $z\curvearrowleft y$, $\mu\curvearrowleft\lambda$,
    $A\curvearrowleft B$ in the notation of \lref{eq:uniform-annular});
    \lref{setting:extension};
    \lref{setting:flow};
    \lref{setting:condition};
  }[verbs=ep]{
    \lref{claim:annular}\dott
  }
  \startnewargseq
  \argument{
    \lref{claim:annular};
    \lref{setting:radius};
    \lref{eq:radial-derivative};
  }{
    that for all $\rho\in(0,\varrho)$ it holds that
    \begin{eqn}
      \llabel{eq:radial-annular}
      \int_{\{t\in I\colon\rho/2\leq r(t)\leq\rho\}}
        \abs{r'(t)}\dd t
        \leq\int_{\{t\in I\colon\rho/2\leq r(t)\leq\rho\}}
          \norm{y'(t)}\dd t\leq\mathfrak{C}\rho
      \dott
    \end{eqn}
  }

  To prove \lref{claim:radial-limit}, let
  $\tau\in\{\inf I,\sup I\}\backslash I$ satisfy
  \begin{eqn}
    \llabel{setting:terminal}
    \liminf_{t\to\tau}\norm{y(t)}=0
  \end{eqn}
  and let
  $\mathbf{p}\colon I\to(0,\varrho/2]$ satisfy for all $t\in I$ that
  \begin{eqn}
    \llabel{setting:clipped-radius}
    \mathbf{p}(t)=\min\{r(t),\varrho/2\}
    \dott
  \end{eqn}
  \startnewargseq
  \argument{
    \lref{setting:clipped-radius};
    the fact that $r$ is locally absolutely continuous;
    \cref{lem:AC-level-set};
    the chain rule;
  }{
    that $\mathbf{p}$ is locally absolutely continuous and that for
    Lebesgue-almost all $t\in I$ it holds that
    \begin{eqn}
      \llabel{eq:clipped-derivative}
      \mathbf{p}'(t)=\indicator{(0,\varrho/2)}(r(t))r'(t)
      \dott
    \end{eqn}
  }
  \argument{
    \lref{claim:annular};
    \lref{setting:radius};
    \lref{eq:radial-derivative};
    \lref{eq:clipped-derivative};
  }{
    that
    \begin{eqn}
      \llabel{eq:clipped-integral}
      \int_I\abs{\mathbf{p}'(t)}\dd t
        =\int_{\{t\in I\colon r(t)<\varrho/2\}}\abs{r'(t)}\dd t
        \leq\int_{\{t\in I\colon r(t)\leq\varrho/2\}}\norm{y'(t)}\dd t
        \leq\mathfrak{C}\varrho<\infty
      \dott
    \end{eqn}
  }
  \argument{
    \lref{eq:clipped-integral};
    the fact that $\mathbf{p}$ is locally absolutely continuous;
    the fundamental theorem of calculus;
  }{
    that there exists $\kappa\in[0,\varrho/2]$ which satisfies that
    \begin{eqn}
      \llabel{eq:clipped-limit}
      \limsup_{t\to\tau}\abs{\mathbf{p}(t)-\kappa}=0
      \dott
    \end{eqn}
  }
  \startnewargseq
  \argument{
    \lref{setting:radius};
    \lref{setting:terminal};
    \lref{setting:clipped-radius};
    \lref{eq:clipped-limit};
  }{
    that
    \begin{eqn}
      \llabel{eq:zero-clipped-limit}
      \kappa=\liminf_{t\to\tau}\min\{r(t),\varrho/2\}
        =\min\Bigl\{\liminf_{t\to\tau}r(t),\varrho/2\Bigr\}=0
      \dott
    \end{eqn}
  }
  \argument{
    \lref{setting:radius};
    \lref{setting:clipped-radius};
    \lref{eq:clipped-limit};
    \lref{eq:zero-clipped-limit};
  }{that
    \begin{eqn}
      \llabel{eq:radial-limit}
      \lim_{t\to\tau}\norm{y(t)}
        =\lim_{t\to\tau}r(t)
        =\lim_{t\to\tau}\mathbf{p}(t)=0
      \dott
    \end{eqn}
  }
  This and the arbitrariness of $\tau$ complete our proof of
  \lref{claim:radial-limit}.
  \startnewargseq
  We next prove \lref{claim:angular}.
  \startnewargseq
  \argument{
    \lref{eq:uniform-state};
    \lref{setting:extension};
    \lref{setting:polar};
    \lref{setting:flow};
    \lref{setting:condition};
    \lref{setting:mixing};
    \lref{setting:radius};
    the fact that $\mathfrak R<R$;
    the fact that $y$ is locally absolutely continuous;
    the fact that $\lambda$ and $B$ are Lebesgue measurable;
    the fact that for all $t\in I$ with $r(t)<\mathfrak R$
    it holds that $y(t)\in U\cap V\subseteq W$;
    the fact that for all $z\in U\cap V$ it holds that
    $\nabla F(z)=\nabla g(z)$;
  }{
    that for all nonempty intervals $J\subseteq I$ with
    $r(J)\subseteq(0,\mathfrak R)$ it holds that
    $J\ni t\mapsto Q(r(t),\omega(t))\in[0,M]$ is locally absolutely
    continuous and that for Lebesgue-almost all $t\in J$ it holds that
    \begin{eqn}
      \llabel{eq:local-state}
      \frac{\dd}{\dd t}Q(r(t),\omega(t))
        \leq-\norm{\omega'(t)}+m(r(t))\abs{r'(t)}
      \dott
    \end{eqn}
  }
  In the following let $S\colon I\to[0,M]$ satisfy for all $t\in I$ that
  \begin{eqn}
    \llabel{setting:state}
    S(t)=\begin{cases}
      \chi(r(t))Q(r(t),\omega(t))&\colon r(t)<\mathfrak R,\\
      0&\colon r(t)\geq\mathfrak R
    \end{cases}
    \dott
  \end{eqn}
  \startnewargseq
  \argument{
    \lref{setting:cutoff};
    \lref{eq:local-state};
    \lref{setting:state};
    the fact that $r$ is locally absolutely continuous;
  }{
    that $S$ is locally absolutely continuous on
    $\{t\in I\colon r(t)<\mathfrak R\}$ and that
    \begin{eqn}
      \llabel{eq:state-zero}
      S|_{\{t\in I\colon r(t)>\eta\}}=0
      \dott
    \end{eqn}
  }
  \argument{
    \lref{eq:state-zero};
    the fact that $\eta<\mathfrak R$;
    the fact that $r$ is continuous;
  }{
    that $S$ is locally absolutely continuous on $I$\dott
  }
  \argument{
    \lref{setting:cutoff};
    the fact that $r$ is locally absolutely continuous;
    \cref{lem:AC-level-set};
    the chain rule;
  }{
    that $\chi\circ r$ is locally absolutely continuous and that for
    Lebesgue-almost all $t\in I$ it holds that
    \begin{eqn}
      \llabel{eq:cutoff-derivative}
      \abs{(\chi\circ r)'(t)}
        =\frac2\eta\indicator{(\eta/2,\eta)}(r(t))\abs{r'(t)}
      \dott
    \end{eqn}
  }
  \argument{
    \lref{setting:cutoff};
    \lref{eq:local-state};
    \lref{setting:state};
    \lref{eq:cutoff-derivative};
    the fact that for all $(s,x)\in(0,\mathfrak R)\times\mathbb S$
    with $sx\in W$
    it holds that $0\leq Q(s,x)\leq M$;
  }{
    that for Lebesgue-almost all $t\in I$ it holds that
    \begin{eqn}
      \llabel{eq:state-derivative}
      S'(t)\leq-\chi(r(t))\norm{\omega'(t)}
        +\mathfrak m(r(t))\abs{r'(t)}
        +M\abs{(\chi\circ r)'(t)}
      \dott
    \end{eqn}
  }
  \argument{
    \lref{eq:majorant};
    \lref{setting:cutoff};
    \lref{eq:radial-annular};
    \itref{lem:radial-integrability}{item4}
    (applied with $R\curvearrowleft\mathfrak R/2$,
    $m\curvearrowleft m|_{(0,\mathfrak R/2]}$,
    $\delta\curvearrowleft\eta$,
    $S\curvearrowleft\eta$ in the notation of
    \cref{lem:radial-integrability});
  }{
    that
    \begin{eqn}
      \llabel{eq:radial-error}
      \int_I\mathfrak m(r(t))\abs{r'(t)}\dd t
        =\int_{\{t\in I\colon r(t)\leq\eta\}}
          m(r(t))\abs{r'(t)}\dd t
        \leq2\mathfrak{C}\int_0^\eta\psi(s)\dd s<\infty
      \dott
    \end{eqn}
  }
  \argument{
    \lref{eq:radial-annular};
    \lref{eq:cutoff-derivative};
    the fact that $\eta\in(0,\varrho)$;
  }{
    that
    \begin{eqn}
      \llabel{eq:cutoff-error}
      \int_I\abs{(\chi\circ r)'(t)}\dd t
        =\frac2\eta\int_{\{t\in I\colon\eta/2<r(t)<\eta\}}
          \abs{r'(t)}\dd t
        \leq\frac2\eta\mathfrak{C}\eta=2\mathfrak{C}
      \dott
    \end{eqn}
  }
  \argument{
    \lref{setting:cutoff};
    \lref{eq:state-derivative};
    the fact that $S$ is locally absolutely continuous;
    the fact that for all $t\in I$ it holds that $0\leq S(t)\leq M$;
    the fundamental theorem of calculus;
  }{
    that for all $a,b\in I$ with $a<b$ it holds that
    \begin{eqn}[a]
      \llabel{eq:compact-angular}
      \int_{\{t\in[a,b]\colon r(t)\leq\eta/2\}}
        \norm{\omega'(t)}\dd t
      &\leq\int_a^b\chi(r(t))\norm{\omega'(t)}\dd t\\
      &\leq S(a)-S(b)+\int_a^b\mathfrak m(r(t))\abs{r'(t)}\dd t
        +M\int_a^b\abs{(\chi\circ r)'(t)}\dd t\\
      &\leq M+\int_a^b\mathfrak m(r(t))\abs{r'(t)}\dd t
        +M\int_a^b\abs{(\chi\circ r)'(t)}\dd t
      \dott
    \end{eqn}
  }
  \argument{
    \lref{setting:constants};
    \lref{eq:radial-error};
    \lref{eq:cutoff-error};
    \lref{eq:compact-angular};
  }{
    that for all $a,b\in I$ with $a<b$ it holds that
    \begin{eqn}
      \llabel{eq:compact-angular-bound}
      \int_{\{t\in[a,b]\colon r(t)\leq\eta/2\}}
        \norm{\omega'(t)}\dd t
        \leq M+2\mathfrak{C}\int_0^\eta\psi(s)\dd s+2M\mathfrak{C}=L
      \dott
    \end{eqn}
  }
  \argument{
    \lref{setting:radius};
    \lref{eq:compact-angular-bound};
    the monotone convergence theorem;
  }[verbs=ep]{
    \lref{claim:angular}\dott
  }

  To prove \lref{claim:direction}, let $c\in I$ and
  $\boldsymbol{\tau}\in\{\inf I,\sup I\}\backslash I$
  satisfy $\liminf_{t\to\boldsymbol{\tau}}\norm{y(t)}=0$ and let
  \begin{eqn}
    \llabel{setting:terminal-interval}
    J=I\cap(\min\{c,\boldsymbol{\tau}\},\max\{c,\boldsymbol{\tau}\})
    \dott
  \end{eqn}
  \startnewargseq
  \argument{
    \lref{claim:radial-limit};
    \lref{setting:terminal-interval};
    the fact that $\eta>0$;
  }{
    that there exists $d\in J$ which satisfies for all
    $t\in I\cap(\min\{d,\boldsymbol{\tau}\},\max\{d,\boldsymbol{\tau}\})$ that
    \begin{eqn}
      \llabel{eq:terminal-radius}
      \norm{y(t)}\leq\frac{\eta}{2}
      \dott
    \end{eqn}
  }
  \startnewargseq
  \argument{
    \lref{claim:angular};
    \lref{setting:terminal-interval};
    \lref{eq:terminal-radius};
    the fact that $\omega$ is locally absolutely continuous;
  }{
    that
    \begin{eqn}
      \llabel{eq:terminal-length}
      \int_J\norm{\omega'(t)}\dd t
        \leq\int_{\min\{c,d\}}^{\max\{c,d\}}\norm{\omega'(t)}\dd t
          +\int_{\{t\in I\colon\norm{y(t)}\leq\eta/2\}}
            \norm{\omega'(t)}\dd t<\infty
      \dott
    \end{eqn}
  }
  \argument{
    the fact that $\omega$ is locally absolutely continuous;
    the fundamental theorem of calculus;
  }{
    that for all $u,w\in J$ with $u<w$ it holds that
    \begin{eqn}
      \llabel{eq:angular-cauchy}
      \norm{\omega(w)-\omega(u)}
        =\biggl|\int_u^w\omega'(t)\dd t\biggr|
        \leq\int_u^w\norm{\omega'(t)}\dd t
      \dott
    \end{eqn}
  }
  \argument{
    \lref{setting:terminal-interval};
    \lref{eq:terminal-length};
    \lref{eq:angular-cauchy};
    the fact that $\mathbb S$ is complete;
  }{
    that there exists $v\in\mathbb S$ such that
    \begin{eqn}
      \llabel{eq:direction}
      \limsup_{t\to\boldsymbol{\tau}}\norm{\omega(t)-v}=0
      \dott
    \end{eqn}
  }
  \argument{
    \lref{eq:terminal-length};
    \lref{eq:direction};
    the fact that $c,\boldsymbol{\tau}$ were arbitrary subject to
    the stated assumptions;
  }[verbs=p]{
    \lref{claim:direction}\dott
  }
  This and the arbitrariness of $U,F,I,y,\lambda,B,\omega$
  complete our proof of
  \lref{claim:annular,claim:radial-limit,claim:angular,claim:direction}.

  \startnewargseq
  \argument{
    \lref{eq:majorant};
    \itref{prop:ambient-mixing}{claim:annular,claim:radial-limit,claim:angular,claim:direction};
  }[verbs=ep]{
    \lref{item:annular,item:radial-limit,item:angular,item:direction}\dott
  }
\end{aproof}

\subsection{Finite spherical length for GF trajectories converging to finite points}
\label{subsec:ambient-finite-length}

\begin{athm}{theorem}{maintheorem:ambient1}
Let $n\in\N$, $\vartheta_0,\vartheta_1,\vartheta_2\in\R^n$, $c,\varepsilon\in(0,\infty)$,
let $I\subseteq\R$ be an interval with nonempty interior,
let $U\subseteq\R^n$ be open,
let $F\in C^1(U,\R)$,
let $V\subseteq U$ be a set, assume that $F|_V$ is definable,
let $B\colon\R^n\to\R^{n\times n}$ and
$\lambda\colon\R\to[0,\infty)$ be measurable,
assume for all $v\in\R^n$ that $B(v)$ is symmetric positive-definite,
let $\theta\colon I\to U$ be locally absolutely continuous,
assume for Lebesgue-almost all $t\in I$ that
\begin{equation}
  \label{maintheorem:ambient1.flow}
  \EMPH{\theta'(t)=-\lambda(t)B(\theta(t))(\nabla F)(\theta(t))},
\end{equation}
assume for all $i\in\{0,1,2\}$, $v\in\R^n$ with
$0<\norm{v-\vartheta_i}
  (\indicator{\theta(I)}(\vartheta_1)+|i-1|)
  <\varepsilon$ that
$v\in V$ and $\opnorm{B(v)-B(\vartheta_i)}\leq c\norm{v-\vartheta_i}^{\varepsilon}$,
assume $\sum_{z\in\{-1,1\}}\liminf_{t\to z\sup(zI)}
  \bigl([1-\indicator{I}(z\sup(zI))]\,
    \norm{\theta(t)-\vartheta_{z+1}}\bigr)=0$,
and let $P\colon U\to\R^n$ satisfy for all
$v\in U\backslash\{\vartheta_1\}$ that
$P(v)=\norm{v-\vartheta_1}^{-1}(v-\vartheta_1)$.
Then \EMPH{$P\circ\theta$} has finite length.
\end{athm}

\begin{aproof}
  Throughout this proof let $a=\inf I$, $b=\sup I$, let
  $E\colon I\to\R$ satisfy for all $t\in I$ that $E(t)=F(\theta(t))$,
  and let $H\subseteq I$, $\mathcal C\subseteq\R^n$ satisfy
  \begin{eqn}
    \llabel{setting1}
    H=\theta^{-1}(\{\vartheta_1\})
    \qquad\text{and}\qquad
    \mathcal C=\{\vartheta_0,\vartheta_2\}
      \cup\bigl(\{\vartheta_1\}\cap\theta(I)\bigr)
    \dott
  \end{eqn}
  \cfload[.]
  \startnewargseq
  \argument{
    \cref{maintheorem:ambient1.flow};
    the fact that $\theta$ is locally absolutely continuous;
    the fact that $F\in C^1(U,\R)$;
    the fact that $E=F\circ\theta$;
    the chain rule;
    the fact that for all $t\in\R$ it holds that $\lambda(t)\geq0$;
    the fact that for all $v\in\R^n$ it holds that $B(v)$ is positive-definite;
  }{
    that $E$ is locally absolutely continuous and non-increasing
    and that for Lebesgue-almost all $t\in I$ it holds that
    $E$ and $\theta$ are differentiable at $t$ and that
    \begin{eqn}
      \llabel{eq1}
      E'(t)=-\lambda(t)\langle\nabla F(\theta(t)),
        B(\theta(t))\nabla F(\theta(t))\rangle\leq0
      \dott
    \end{eqn}
  }
  \argument{
    \lref{eq1};
    the fact that $E=F\circ\theta$;
  }{
    that for all $u,v\in I$ with $u<v$ and $\theta(u)=\theta(v)$
    it holds that $E|_{[u,v]}$ is constant and that for
    Lebesgue-almost all $t\in[u,v]$ it holds that
    \begin{eqn}
      \llabel{eq2}
      \lambda(t)\langle\nabla F(\theta(t)),
        B(\theta(t))\nabla F(\theta(t))\rangle=0
      \dott
    \end{eqn}
  }
  \argument{
    \lref{eq2};
    \cref{maintheorem:ambient1.flow};
    the fact that for all $v\in\R^n$ it holds that $B(v)$ is positive-definite;
    the fact that $\theta$ is locally absolutely continuous;
    the fundamental theorem of calculus;
  }{
    that for all $u,v,t\in I$ with $u\leq t\leq v$ and
    $\theta(u)=\theta(v)$ it holds that
    \begin{eqn}
      \llabel{eq3}
      \theta(t)=\theta(u)
      \dott
    \end{eqn}
  }
  \argument{
    \lref{eq3};
    \lref{setting1};
    the fact that $\theta$ is continuous;
  }{
    that $H$ is closed in $I$ and that $H$ is an interval\dott
  }
  In the following let $\mathcal J\subseteq\mathcal P(I)$ satisfy
  \begin{eqn}
    \llabel{setting1.1}
    \mathcal J=
    \begin{cases}
      \{I\}&\colon H=\emptyset,\\
      \bigl\{\{t\in I\colon t<\inf H\},
        \{t\in I\colon t>\sup H\}\bigr\}\backslash\{\emptyset\}
        &\colon H\neq\emptyset
    \end{cases}
    \dott
  \end{eqn}
  \startnewargseq
  \argument{
    \lref{setting1.1};
  }{
    that for all $J\in\mathcal J$ it holds that
    $J$ is a nonempty interval and that $J$ is open in $I$\dott
  }

  In our proof of $\curvelength{P\circ\theta}<\infty$ we distinguish between
  the case $n=1$ and the case $n\geq2$\dott\cfload[.]
  \startnewargseq
  We first prove $\curvelength{P\circ\theta}<\infty$ in the case
  \begin{eqn}
    \llabel{case1}
    n=1
    \dott
  \end{eqn}
  \startnewargseq
  \argument{
    \lref{setting1.1};
    the fact that for all $J\in\mathcal J$ it holds that $J$ is an interval;
    the fact that $\theta$ is continuous;
    the fact that for all $v\in U\backslash\{\vartheta_1\}$ it holds that
    $P(v)=\norm{v-\vartheta_1}^{-1}(v-\vartheta_1)$;
    the intermediate value theorem;
  }{
    that for all $J\in\mathcal J$ it holds that
    $(P\circ\theta)|_J$ is constant with value in $\{-1,1\}$\dott
    \llabel{arg1}
  }
  \argument{
    \lref{arg1};
    \lref{setting1};
    \lref{setting1.1};
    the fact that $H$ is convex;
    the fact that $H$ is closed in $I$;
    the fact that for all $J\in\mathcal J$ it holds that $J$ is an interval;
    the fact that $(P\circ\theta)|_H$ is constant;
    the triangle inequality;
  }[verbs=ep]{
    that $P\circ\theta$ has finite length in the case $n=1$\dott
  }

  We now prove $\curvelength{P\circ\theta}<\infty$ in the case
  \begin{eqn}
    \llabel{case2}
    n\geq2
    \dott
  \end{eqn}
  \startnewargseq
  \argument{
    \lref{setting1};
    the fact that for all $i\in\{0,1,2\}$, $v\in\R^n$ with
    $0<\norm{v-\vartheta_i}
      \bigl(\indicator{\theta(I)}(\vartheta_1)+\abs{i-1}\bigr)
      <\varepsilon$ it holds that $v\in V$ and that
    $\opnorm{B(v)-B(\vartheta_i)}\leq c\norm{v-\vartheta_i}^{\varepsilon}$;
  }{
    that for all $p\in\mathcal C$, $v\in\R^n$ with
    $0<\norm{v-p}<\varepsilon/2$ it holds that
    \begin{eqn}
      \llabel{eq5}
      v\in V
      \qquad\text{and}\qquad
      \opnorm{B(v)-B(p)}\leq c\norm{v-p}^{\varepsilon}
      \dott
    \end{eqn}
  }
  \argument{
    the fact that for all $p\in\mathcal C$ it holds that
    $B(p)$ is symmetric and positive-definite;
    the spectral theorem;
  }{
    that for all $p\in\mathcal C$ it holds that there exists an invertible
    $L_p\in\R^{n\times n}$ which satisfies that
    \begin{eqn}
      \llabel{eq5.1}
      L_p(L_p)^{\mathsf T}=B(p)
      \dott
    \end{eqn}
  }
  In the following let $q\in\mathbb Q\cap(0,\varepsilon]$,
  for every $p\in\mathcal C$ let $W_p\subseteq\R^n$ satisfy
  $W_p=\{z\in\R^n\colon p+L_pz\in U\}$ and let $f_p\colon W_p\to\R$ and
  $\mathbf{B}_p\colon\R^n\to\R^{n\times n}$
  satisfy for all $z\in W_p$, $w\in\R^n$ that
  \begin{eqn}
    \llabel{setting2}
    f_p(z)=F(p+L_pz)
    \qquad\text{and}\qquad
    \mathbf{B}_p(w)=(L_p)^{-1}B(p+L_pw)((L_p)^{-1})^{\mathsf T},
  \end{eqn}
  let $\kappa=\max_{p\in\mathcal C}
    c\opnorm{(L_p)^{-1}}^2\opnorm{L_p}^{\varepsilon}$,
  let $R\in(0,1)$ satisfy
  \begin{eqn}
    \llabel{setting3}
    R\max_{p\in\mathcal C}\opnorm{L_p}<\varepsilon/2
    \qquad\text{and}\qquad
    \kappa R^{\varepsilon}\leq1/2,
  \end{eqn}
  and let $\mathbf{k}\colon(0,R)\to[0,\infty)$ satisfy for all
  $r\in(0,R)$ that $\mathbf{k}(r)=2\kappa r^{q-1}$\dott
  \startnewargseq
  \argument{
    \lref{eq5};
    \lref{eq5.1};
    \lref{setting2};
    \lref{setting3};
    the fact that $V\subseteq U$;
    the fact that $F|_V$ is definable;
    the fact that $F\in C^1(U,\R)$;
  }{
    that for all $p\in\mathcal C$ it holds that $W_p$ is open,
    that $f_p\in C^1(W_p,\R)$, that
    $\{z\in\R^n\colon0<\norm{z}<R\}\subseteq W_p$,
    and that $f_p|_{\{z\in\R^n\colon0<\norm{z}<R\}}$ is definable\dott
  }
  \argument{
    \lref{eq5};
    \lref{eq5.1};
    \lref{setting2};
    \lref{setting3};
  }{
    that for all $p\in\mathcal C$, $z\in\R^n$ with
    $0<\norm{z}<R$ it holds that
    \begin{eqn}
      \llabel{eq6}
      \opnorm{\mathbf{B}_p(z)-\Id}
      =\opnorm[\big]{(L_p)^{-1}\bigl(B(p+L_pz)-B(p)\bigr)((L_p)^{-1})^{\mathsf T}}
      \leq c\opnorm{(L_p)^{-1}}^2\norm{L_pz}^{\varepsilon}
      \leq\kappa\norm{z}^{\varepsilon}\leq\frac12
      \dott
    \end{eqn}
  }
  \argument{
    \lref{eq5.1};
    \lref{setting2};
    the fact that $B$ is measurable;
    the fact that for all $v\in\R^n$ it holds that
    $B(v)$ is symmetric and positive-definite;
  }{
    that for all $p\in\mathcal C$ it holds that $\mathbf{B}_p$ is measurable
    and that for all $p\in\mathcal C$, $z\in\R^n$ it holds that
    $\mathbf{B}_p(z)$ is symmetric and positive-definite\dott
    \llabel{arg2}
  }
  \argument{
    \lref{arg2};
    \lref{eq6};
  }{
    that for all $p\in\mathcal C$, $r\in(0,R)$, $\omega\in\R^n$,
    $v\in\R^n\backslash\{0\}$, $\mu\in\R$ with
    $\mathbf{B}_p(r\omega)v=\mu v$ and $\norm{\omega}=1$ it holds that
    \begin{eqn}
      \llabel{eq6.1}
      \frac12\leq\mu\leq\frac32
      \dott
    \end{eqn}
  }
  \argument{
    \lref{eq6.1};
    \lref{setting3};
    \lref{eq6};
    the fact that for all $r\in(0,R)$ it holds that
    $\mathbf{k}(r)=2\kappa r^{q-1}$;
    the fact that for all $\omega\in\R^n$ with $\norm{\omega}=1$
    it holds that $\opnorm{\Id-\omega\omega^{\mathsf T}}\leq1$;
  }{
    that for all $p\in\mathcal C$, $r\in(0,R)$, $\omega\in\R^n$ with
    $\norm{\omega}=1$ it holds that
    \begin{eqn}[a]
      \llabel{eq7}
      \frac{\norm{(\Id-\omega\omega^{\mathsf T})\mathbf{B}_p(r\omega)\omega}}
        {\langle\mathbf{B}_p(r\omega)\omega,\omega\rangle}\leq2\opnorm{\mathbf{B}_p(r\omega)-\Id}
      \leq2\kappa r^{\varepsilon}
       \leq2\kappa r^q=r\mathbf{k}(r)
      \dott
    \end{eqn}
  }
  \argument{
    the fact that $q\in\mathbb Q\cap(0,\infty)$;
    the fact that for all $r\in(0,R)$ it holds that
    $\mathbf{k}(r)=2\kappa r^{q-1}$;
  }{
    that $\mathbf{k}$ is nonnegative, locally bounded, and definable,
    and that
    \begin{eqn}
      \llabel{eq8}
      \int_0^R\mathbf{k}(r)\dd r=\frac{2\kappa R^q}{q}<\infty
      \dott
    \end{eqn}
  }
  \argument{
    \lref{eq8};
    \lref{case2};
    \lref{setting3};
    \lref{eq6.1};
    \lref{eq7};
    the fact that for all $p\in\mathcal C$ it holds that $W_p$ is open;
    the fact that for all $p\in\mathcal C$ it holds that
    $\{z\in\R^n\colon0<\norm{z}<R\}\subseteq W_p$;
    the fact that for all $p\in\mathcal C$ it holds that
    $f_p\in C^1(W_p,\R)$;
    the fact that for all $p\in\mathcal C$ it holds that
    $f_p|_{\{z\in\R^n\colon0<\norm{z}<R\}}$ is definable;
    the fact that $\lambda$ is measurable with values in $[0,\infty)$;
    the fact that for all $p\in\mathcal C$ it holds that
    $\mathbf{B}_p$ is measurable;
    the fact that for all $p\in\mathcal C$, $z\in\R^n$ it holds that
    $\mathbf{B}_p(z)$ is symmetric and positive-definite;
    \itref{prop:ambient-mixing}{item:annular,item:angular}
    (applied for every $p\in\mathcal C$ with $K\curvearrowleft3$,
    $V\curvearrowleft\{z\in\R^n\colon0<\norm{z}<R\}$,
    $W\curvearrowleft\{z\in\R^n\colon0<\norm{z}<R\}$,
    $g\curvearrowleft f_p|_{\{z\in\R^n\colon0<\norm{z}<R\}}$,
    $U\curvearrowleft W_p$, $F\curvearrowleft f_p$
    in the notation of \cref{prop:ambient-mixing});
  }{
    that for all $p\in\mathcal C$ it holds that there exist
    $\varrho_p,\eta_p,\mathfrak{C}_p\in(0,\infty)$ which satisfy that
    \begin{eqn}
      \llabel{eq8.1}
      0<\eta_p<\varrho_p<R
    \end{eqn}
    and which satisfy that for every nonempty interval $J\subseteq I$,
    every locally absolutely continuous $y\colon J\to W_p\backslash\{0\}$,
    and every $h\in(0,\varrho_p)$
    with the property that for Lebesgue-almost all $t\in J$ it holds
    that $y$ is differentiable at $t$ and that
    $y'(t)=-\lambda(t)\mathbf{B}_p(y(t))\nabla f_p(y(t))$
    we have that
    \begin{eqn}
      \llabel{eq10}
      \int_{\{t\in J\colon\norm{y(t)}\leq h\}}\norm{y'(t)}\dd t
        \leq2\mathfrak{C}_p h
      \qquad\text{and}\qquad
      \int_{\{t\in J\colon\norm{y(t)}\leq\eta_p\}}
        \norm[\bigg]{\frac{\ddd}{\ddd t}
          \biggl(\frac{y(t)}{\norm{y(t)}}\biggr)}\dd t\leq \mathfrak{C}_p
      \dott
    \end{eqn}
  }
  In the following let $\rho,D,T\in\R$ satisfy
  \begin{eqn}
    \llabel{setting3.1}
    \rho=\min_{p\in\mathcal C}\frac{\eta_p}{2\opnorm{(L_p)^{-1}}},
    \qquad
    D=\max_{p\in\mathcal C}2\mathfrak{C}_p\opnorm{L_p}\opnorm{(L_p)^{-1}},
    \qquad\text{and}\qquad
    T=\max_{p\in\mathcal C}\mathfrak{C}_p\opnorm{L_p}\opnorm{(L_p)^{-1}}
    \dott
  \end{eqn}
  \startnewargseq
  \argument{
    \lref{setting1};
    \lref{eq5.1};
    \lref{setting2};
    \lref{eq8.1};
    \lref{eq10};
    \lref{setting3.1};
    the fact that $\mathcal C$ is finite and nonempty;
  }{
    that $\rho,D,T\in(0,\infty)$ and that for all $p\in\mathcal C$
    it holds that
    \begin{eqn}
      \llabel{eq8.2}
      \opnorm{(L_p)^{-1}}\rho\leq\eta_p/2<\eta_p<\varrho_p,
      \qquad
      2\mathfrak{C}_p\opnorm{L_p}\opnorm{(L_p)^{-1}}\leq D,
      \qquad\text{and}\qquad
      \mathfrak{C}_p\opnorm{L_p}\opnorm{(L_p)^{-1}}\leq T
      \dott
    \end{eqn}
  }

  We claim that for every $p\in\mathcal C$,
  every nonempty interval $J\subseteq I$,
  and every $h\in(0,\rho)$ with $p\notin\theta(J)$ we have that
  \begin{eqn}
    \llabel{eq15}
    \int_{\{t\in J\colon\norm{\theta(t)-p}\leq h\}}
      \norm{\theta'(t)}\dd t\leq Dh
    \qquad\text{and}\qquad
    \int_{\{t\in J\colon\norm{\theta(t)-p}\leq\rho\}}
      \norm[\bigg]{\frac{\ddd}{\ddd t}
        \biggl(\frac{\theta(t)-p}{\norm{\theta(t)-p}}\biggr)}\dd t
      \leq T
    \dott
  \end{eqn}
  \startnewargseq
  To prove \lref{eq15}, let $p\in\mathcal C$, $h\in(0,\rho)$,
  let $J\subseteq I$ be a nonempty interval with $p\notin\theta(J)$,
  and let $z\colon J\to W_p\backslash\{0\}$,
  $\omega\colon J\to\R^n$, and $\mathbf{w}\colon J\to\R^n$ satisfy for all
  $t\in J$ that
  \begin{eqn}
    \llabel{setting4}
    z(t)=(L_p)^{-1}(\theta(t)-p),
    \qquad \omega(t)=\frac{z(t)}{\norm{z(t)}},
    \qquad\text{and}\qquad \mathbf{w}(t)=\frac{\theta(t)-p}{\norm{\theta(t)-p}}
    \dott
  \end{eqn}
  \startnewargseq
  \argument{
    \cref{maintheorem:ambient1.flow};
    \lref{eq5.1};
    \lref{setting2};
    \lref{setting4};
    the fact that $\theta$ is locally absolutely continuous;
    the chain rule;
  }{
    that $z$, $\omega$, and $\mathbf{w}$ are locally absolutely continuous
    and that for Lebesgue-almost all $t\in J$ it holds that
    $z$ is differentiable at $t$ and that
    \begin{eqn}[a]
      \llabel{eq9}
      z'(t)
      =-\lambda(t)(L_p)^{-1}B(\theta(t))\nabla F(\theta(t))=-\lambda(t)\mathbf{B}_p(z(t))\nabla f_p(z(t))
      \dott
    \end{eqn}
  }
  \argument{
    \lref{setting4};
  }{
    that for all $t\in J$ it holds that
    \begin{eqn}
      \llabel{eq11}
      \norm{z(t)}\leq\opnorm{(L_p)^{-1}}\norm{\theta(t)-p}
      \qquad\text{and}\qquad
      \mathbf{w}(t)=\frac{L_p\omega(t)}{\norm{L_p\omega(t)}}
      \dott
    \end{eqn}
  }
  \argument{
    \lref{eq11};
    \lref{setting4};
    the fact that $\theta$, $z$, $\omega$, and $\mathbf{w}$
    are locally absolutely continuous;
    the chain rule;
    the fact that for all $v\in\R^n$ it holds that
    $\norm{L_pv}\geq\opnorm{(L_p)^{-1}}^{-1}\norm{v}$;
  }{
    that for Lebesgue-almost all $t\in J$ it holds that
    $\theta$, $z$, $\omega$, and $\mathbf{w}$ are differentiable at $t$ and that
    \begin{eqn}
      \llabel{eq12}
      \norm{\theta'(t)}\leq\opnorm{L_p}\norm{z'(t)}
      \qquad\text{and}\qquad
      \norm{\mathbf{w}'(t)}
      =\norm[\bigg]{\frac{(\Id-\mathbf{w}(t)(\mathbf{w}(t))^{\mathsf T})L_p\omega'(t)}
        {\norm{L_p\omega(t)}}}
        \leq\opnorm{L_p}\opnorm{(L_p)^{-1}}\norm{\omega'(t)}
      \dott
    \end{eqn}
  }
  \argument{
    \lref{eq12};
    \lref{eq10};
    \lref{eq8.2};
    \lref{setting4};
    \lref{eq9};
    \lref{eq11};
    the chain rule;
  }{
    that
    \begin{eqn}[a]
      \llabel{eq13}
      &\int_{\{t\in J\colon\norm{\theta(t)-p}\leq h\}}
        \norm{\theta'(t)}\dd t
        =\int_{\{t\in J\colon\norm{L_pz(t)}\leq h\}}
          \norm{L_pz'(t)}\dd t\\
        &\leq\opnorm{L_p}
          \int_{\{t\in J\colon\norm{z(t)}\leq\opnorm{(L_p)^{-1}}h\}}
            \norm{z'(t)}\dd t\leq2\mathfrak{C}_p\opnorm{L_p}\opnorm{(L_p)^{-1}}h
    \end{eqn}
    and 
    \begin{eqn}[a]
      \llabel{eq13.1}
      &\int_{\{t\in J\colon\norm{\theta(t)-p}\leq\rho\}}
        \norm{\mathbf{w}'(t)}\dd t
        =\int_{\{t\in J\colon\norm{L_pz(t)}\leq\rho\}}
          \norm{\mathbf{w}'(t)}\dd t\\
        &\leq\opnorm{L_p}\opnorm{(L_p)^{-1}}
          \int_{\{t\in J\colon\norm{z(t)}\leq\opnorm{(L_p)^{-1}}\rho\}}
            \norm{\omega'(t)}\dd t\leq\opnorm{L_p}\opnorm{(L_p)^{-1}}
          \int_{\{t\in J\colon\norm{z(t)}\leq\eta_p\}}
            \norm{\omega'(t)}\dd t\\
        &\leq\mathfrak{C}_p\opnorm{L_p}\opnorm{(L_p)^{-1}}
      \dott
    \end{eqn}
  }
  \argument{
    \lref{eq13.1};
    \lref{eq13};
    \lref{eq8.2};
    \lref{setting4};
    the fact that $p$, $J$, and $h$ are arbitrary;
  }[verbs=ep]{
    \lref{eq15}\dott
  }
  This completes our proof of \lref{eq15}\dott
  \startnewargseq

  We claim that for every $\mathbf{p}\in\mathcal C$,
  every nonempty interval $\mathbf{J}\subseteq I$,
  and every $\boldsymbol{\tau}\in\{\inf \mathbf{J},\sup \mathbf{J}\}\backslash \mathbf{J}$
  with $\mathbf{p}\notin\theta(\mathbf{J})$ and
  $\liminf_{t\to\boldsymbol{\tau},\,t\in \mathbf{J}}\norm{\theta(t)-\mathbf{p}}=0$
  it holds that
  \begin{eqn}
    \llabel{eq14}
    \limsup_{t\to\boldsymbol{\tau},\,t\in \mathbf{J}}\norm{\theta(t)-\mathbf{p}}=0
    \dott
  \end{eqn}
  \startnewargseq
  To prove \lref{eq14}, let $\mathbf{p}\in\mathcal C$,
  let $\mathbf{J}\subseteq I$ be a nonempty interval with $\mathbf{p}\notin\theta(\mathbf{J})$,
  let $\boldsymbol{\tau}\in\{\inf \mathbf{J},\sup \mathbf{J}\}\backslash \mathbf{J}$ satisfy
  $\liminf_{t\to\boldsymbol{\tau},\,t\in \mathbf{J}}\norm{\theta(t)-\mathbf{p}}=0$,
  and let $r\colon \mathbf{J}\to(0,\infty)$ and
  $\psi\colon \mathbf{J}\to(0,\rho/2]$ satisfy for all $t\in \mathbf{J}$ that
  \begin{eqn}
    \llabel{setting4.1}
    r(t)=\norm{\theta(t)-\mathbf{p}}
    \qquad\text{and}\qquad
    \psi(t)=\min\{r(t),\rho/2\}
    \dott
  \end{eqn}
  \startnewargseq
  \argument{
    \lref{setting4.1};
    the fact that $\theta$ is locally absolutely continuous;
  }{
    that $r$ and $\psi$ are locally absolutely continuous\dott
  }
  \argument{
    \lref{setting4.1};
    the chain rule;
  }{
    that for Lebesgue-almost all $t\in \mathbf{J}$ it holds that
    $\theta$, $r$, and $\psi$ are differentiable at $t$ and that
    \begin{eqn}
      \llabel{eq14.1}
      \abs{\psi'(t)}
        \leq\indicator{(0,\rho/2]}(r(t))\abs{r'(t)}
        \leq\indicator{(0,\rho/2]}(r(t))\norm{\theta'(t)}
      \dott
    \end{eqn}
  }
  \argument{
    \lref{eq14.1};
    \lref{eq15};
    \lref{setting4.1};
  }{
    that
    \begin{eqn}
      \llabel{eq14.2}
      \int_{\mathbf{J}}\abs{\psi'(t)}\dd t
        \leq\int_{\{t\in \mathbf{J}\colon\norm{\theta(t)-\mathbf{p}}\leq\rho/2\}}
          \norm{\theta'(t)}\dd t
        \leq D\rho/2<\infty
      \dott
    \end{eqn}
  }
  \argument{
    \lref{eq14.2};
    \lref{setting4.1};
    the fact that $\psi$ is locally absolutely continuous;
    the fact that $\liminf_{t\to\boldsymbol{\tau},\,t\in \mathbf{J}}\norm{\theta(t)-\mathbf{p}}=0$;
    the fundamental theorem of calculus;
  }{
    that
    \begin{eqn}
      \llabel{eq14.3}
      \limsup_{t\to\boldsymbol{\tau}}\abs{\psi(t)}=0
      \dott
    \end{eqn}
  }
  \argument{
    \lref{eq14.3};
    \lref{setting4.1};
    the fact that $\rho>0$;
    the fact that for all $t\in \mathbf{J}$ with $\psi(t)<\rho/2$ it holds that
    $r(t)=\psi(t)$;
    the fact that $\mathbf{p}$, $\mathbf{J}$, and $\boldsymbol{\tau}$ are arbitrary;
  }[verbs=ep]{
    \lref{eq14}\dott
  }
  This completes our proof of \lref{eq14}\dott
  \startnewargseq

  We claim that for all $\delta\in\{a,b\}\backslash I$ 
  there exists $s\in\operatorname{int}(I)$ such that
  \begin{eqn}
    \llabel{claim:endpoints}
    \int_{I\cap(\min\{s,\delta\},\max\{s,\delta\})}
      \norm{\theta'(t)}\dd t<\infty
    \qquad\text{and}\qquad
    \lim_{t\to\delta}\theta(t)=
    \begin{cases}
      \vartheta_0&\colon\delta=a\\
      \vartheta_2&\colon\delta=b
    \end{cases}
    \dott
  \end{eqn}
  \startnewargseq
  To prove \lref{claim:endpoints}, let $\delta\in\{a,b\}\backslash I$
  and let $\mathfrak{p}\in\mathcal C$ satisfy
  \begin{eqn}
    \llabel{setting5}
    \mathfrak{p}=
    \begin{cases}
      \vartheta_0&\colon\delta=a\\
      \vartheta_2&\colon\delta=b
    \end{cases}
    \dott
  \end{eqn}
  \startnewargseq
  \argument{
    \lref{setting5};
    the fact that
    $\sum_{\sigma\in\{-1,1\}}\liminf_{t\to\sigma\sup(\sigma I)}
      \bigl([1-\indicator{I}(\sigma\sup(\sigma I))]\,
        \norm{\theta(t)-\vartheta_{\sigma+1}}\bigr)=0$;
  }{
    that
    \begin{eqn}
      \llabel{eq16}
      \liminf_{t\to\delta}\norm{\theta(t)-\mathfrak{p}}=0
      \dott
    \end{eqn}
  }
  \argument{
    \lref{eq3};
    the fact that $\theta$ is continuous;
    the fact that $\delta$ is an endpoint of $I$;
  }{
    that there exists $d\in\operatorname{int}(I)$ which satisfies that
    \begin{eqn}
      \llabel{eq17}
      \theta|_{I\cap(\min\{d,\delta\},\max\{d,\delta\})}\equiv \mathfrak{p}
      \quad\text{or}\quad
      \mathfrak{p}\notin\theta\bigl(I\cap(\min\{d,\delta\},\max\{d,\delta\})\bigr)
      \dott
    \end{eqn}
  }
  In the following let
  $\mathbb{J}=I\cap(\min\{d,\delta\},\max\{d,\delta\})$\dott
  \startnewargseq
  \argument{
    the fact that $\delta\in\{\inf \mathbb{J},\sup \mathbb{J}\}\backslash \mathbb{J}$;
  }{
    that if $\theta|_{\mathbb{J}}\equiv \mathfrak{p}$, then
    \begin{eqn}
      \llabel{eq18}
      \limsup_{t\to\delta}\norm{\theta(t)-\mathfrak{p}}=0
      \qquad\text{and}\qquad
      \int_{\mathbb{J}}\norm{\theta'(t)}\dd t=0
      \dott
    \end{eqn}
  }
  \argument{
    \lref{eq14};
    \lref{eq16};
    the fact that
    $I\cap(\min\{d,\delta\},\max\{d,\delta\})\subseteq\mathbb{J}$;
  }{
    that if $\mathfrak{p}\notin\theta(\mathbb{J})$, then
    \begin{eqn}
      \llabel{eq19}
      \limsup_{t\to\delta}\norm{\theta(t)-\mathfrak{p}}=0
      \dott
    \end{eqn}
  }
  \argument{
    \lref{eq19};
    the fact that $\rho>0$;
    the fact that $\mathbb{J}\subseteq I$ is a nonempty interval and that
    $\delta\in\{\inf\mathbb{J},\sup\mathbb{J}\}\backslash\mathbb{J}$;
  }{
    that if $\mathfrak{p}\notin\theta(\mathbb{J})$, then there exists $e\in \mathbb{J}$
    such that for all
    $t\in I\cap(\min\{e,\delta\},\max\{e,\delta\})$ it holds that
    $\norm{\theta(t)-\mathfrak{p}}\leq\rho/2$\dott
    \llabel{arg3}
  }
  \argument{
    \lref{arg3};
    \lref{eq15};
    \lref{setting5};
  }{
    that if $\mathfrak{p}\notin\theta(\mathbb{J})$, then there exists $e\in\mathbb{J}$
    such that
    \begin{eqn}
      \llabel{eq20}
      \int_{I\cap(\min\{e,\delta\},\max\{e,\delta\})}
        \norm{\theta'(t)}\dd t
        \leq\int_{\{t\in\mathbb{J}\colon\norm{\theta(t)-\mathfrak{p}}\leq\rho/2\}}
          \norm{\theta'(t)}\dd t\leq D\rho/2<\infty
      \dott
    \end{eqn}
  }
  \argument{
    \lref{eq20};
    \lref{setting5};
    \lref{eq17};
    \lref{eq18};
    \lref{eq19};
    the fact that $\mathbb{J}\subseteq\operatorname{int}(I)$;
    the fact that $\delta\in\{a,b\}\backslash I$ is arbitrary;
  }[verbs=ep]{
    \lref{claim:endpoints}\dott
  }
  This completes our proof of \lref{claim:endpoints}\dott

  \startnewargseq
  \argument{
    \lref{claim:endpoints};
    the fact that $a=\inf I$ and $b=\sup I$;
    the fact that $I$ is an interval with nonempty interior;
  }{
    that there exist $u,v\in I$ with $u<v$ which satisfy that
    if $a\in I$, then $u=a$, that if $b\in I$, then $v=b$, and that
    \begin{eqn}
      \llabel{eq20.1}
      \int_{I\cap(-\infty,u)}\norm{\theta'(t)}\dd t
        +\int_{I\cap(v,\infty)}\norm{\theta'(t)}\dd t<\infty
      \dott
    \end{eqn}
  }
  \argument{
    \lref{eq20.1};
    the fact that $[u,v]\subseteq I$ is compact;
    the fact that $\theta$ is locally absolutely continuous;
    \cref{lem:curve-length-speed};
  }{
    that
    \begin{eqn}[a]
      \llabel{eq21}
      \curvelength{\theta}
        &=\int_I\norm{\theta'(t)}\dd t\\
        &=\int_{I\cap(-\infty,u)}\norm{\theta'(t)}\dd t
          +\int_u^v\norm{\theta'(t)}\dd t
          +\int_{I\cap(v,\infty)}\norm{\theta'(t)}\dd t<\infty
      \dott
    \end{eqn}
  }
  \argument{
    \lref{setting1};
    \lref{setting1.1};
    the fact that $\theta$ is locally absolutely continuous;
    the fact that for all $x\in U\backslash\{\vartheta_1\}$ it holds that
    $P(x)=\norm{x-\vartheta_1}^{-1}(x-\vartheta_1)$;
    the chain rule;
  }{
    that for all $J\in\mathcal J$ it holds that $(P\circ\theta)|_J$
    is locally absolutely continuous and that for Lebesgue-almost
    all $t\in J$ it holds that $P\circ\theta$ and $\theta$ are
    differentiable at $t$ and that
    \begin{eqn}
      \llabel{eq22}
      \norm{(P\circ\theta)'(t)}
        \leq\frac{\norm{\theta'(t)}}{\norm{\theta(t)-\vartheta_1}}
      \dott
    \end{eqn}
  }
  In our proof of $\curvelength{P\circ\theta}<\infty$ in the case $n\geq2$
  we distinguish between the case $\vartheta_1\notin\mathcal C$
  and the case $\vartheta_1\in\mathcal C$\dott
  \startnewargseq
  We first prove $\curvelength{P\circ\theta}<\infty$ in the case
  \begin{eqn}
    \llabel{case3}
    n\geq2
    \qquad\text{and}\qquad
    \vartheta_1\notin\mathcal C
    \dott
  \end{eqn}
  \startnewargseq
  \argument{
    \lref{case3};
    \lref{setting1};
    \lref{setting1.1};
    \lref{claim:endpoints};
    the fact that $\theta$ is continuous;
  }{
    that $H=\emptyset$, that $\mathcal J=\{I\}$,
    and that
    \begin{eqn}
      \llabel{eq23}
      \inf_{t\in I}\norm{\theta(t)-\vartheta_1}>0
      \dott
    \end{eqn}
  }
  \argument{
    \lref{eq23};
    \lref{eq21};
    \lref{eq22};
    \cref{lem:curve-length-speed};
  }[verbs=ep]{
    that
    \begin{eqn}
      \curvelength{P\circ\theta}
        \leq\frac{\curvelength{\theta}}
          {\inf_{t\in I}\norm{\theta(t)-\vartheta_1}}<\infty
      \dott
    \end{eqn}
  }
  This proves $\curvelength{P\circ\theta}<\infty$ in the case
  $n\geq2$ and $\vartheta_1\notin\mathcal C$\dott
  \startnewargseq
  We now prove $\curvelength{P\circ\theta}<\infty$ in the case
  \begin{eqn}
    \llabel{case4}
    n\geq2
    \qquad\text{and}\qquad
    \vartheta_1\in\mathcal C
    \dott
  \end{eqn}
  \startnewargseq
  \argument{
    \lref{eq15};
    \lref{eq21};
    \lref{eq22};
    \cref{lem:curve-length-speed};
  }{
    that for all $J\in\mathcal J$ it holds that
    \begin{eqn}[a]
      \llabel{eq24}
      \curvelength{(P\circ\theta)|_J}
      &=\int_J\norm{(P\circ\theta)'(t)}\dd t\\
      &=\int_{\{t\in J\colon\norm{\theta(t)-\vartheta_1}\leq\rho\}}
          \norm{(P\circ\theta)'(t)}\dd t+\int_{\{t\in J\colon\norm{\theta(t)-\vartheta_1}>\rho\}}
          \norm{(P\circ\theta)'(t)}\dd t\\
      &\leq T+\int_{\{t\in J\colon\norm{\theta(t)-\vartheta_1}>\rho\}}
          \frac{\norm{\theta'(t)}}{\norm{\theta(t)-\vartheta_1}}\dd t\\
      &\leq T+\rho^{-1}\int_{\{t\in J\colon\norm{\theta(t)-\vartheta_1}>\rho\}}
          \norm{\theta'(t)}\dd t\\
      &\leq T+\rho^{-1}\int_I\norm{\theta'(t)}\dd t
        =T+\rho^{-1}\curvelength{\theta}<\infty
      \dott
    \end{eqn}
  }
  \argument{
    \lref{setting1};
    \lref{setting1.1};
    \lref{eq24};
    the fact that $H$ is convex;
    the fact that $H$ is closed in $I$;
    the fact that for all $J\in\mathcal J$ it holds that $J$ is an interval;
    the fact that for all $t\in I\backslash H$ it holds that
    $\norm{P(\theta(t))}=1$;
    the triangle inequality;
  }[verbs=ep]{
    that 
    \begin{eqn}
      \curvelength{P\circ\theta}
        \leq\sum_{J\in\mathcal J}\curvelength{(P\circ\theta)|_J}
          +2\sum_{x\in\{\vartheta_1\}\cap\theta(I)}
            \bigl(1+\norm{P(x)}\bigr)<\infty
      \dott
    \end{eqn}
  }
  This completes our proof of $\curvelength{P\circ\theta}<\infty$ in the case
  $n\geq2$\dott
  \startnewargseq
\end{aproof}

\Cref{maintheorem:ambient1} includes \cref{maincor} in the introduction
as a special case.
More specifically, note that \cref{maintheorem:ambient1}
(applied with $\vartheta_0\curvearrowleft a$, $\vartheta_1\curvearrowleft\vartheta$,
$\vartheta_2\curvearrowleft b$,
$B\curvearrowleft(\R^n\ni v\mapsto\mathrm{I}_n\in\R^{n\times n})$,
$\lambda\curvearrowleft(\R\ni t\mapsto1\in[0,\infty))$
in the notation of \cref{maintheorem:ambient1}) establishes \cref{maincor}.

\subsection{Uniform upper bounds for spherical lengths under H\"{o}lder estimates}
\label{subsec:ambient-uniform-length}

\begin{athm}{theorem}{maintheorem:ambient2}
Let $n\in\N$, $\vartheta\in\R^n$, $c,\varepsilon\in(0,\infty)$,
let $V=\{v\in\R^n\colon0<|v-\vartheta|<\varepsilon\}$,
let $\mathbb S=\{A\in\R^{n\times n}\colon A\text{ is symmetric positive-definite}\}$,
let $\mathbb B\in\mathbb S$,
and let $f\in C(V,\R)$ be definable.
Then there exist $\delta,L\in(0,\infty)$ such that
for every interval with nonempty interior $I\subseteq\R$,
every open $U\subseteq\R^n$,
every $F\in C^1(U,\R)$,
every locally absolutely continuous $\theta\colon I\to U$,
and all measurable $P\colon U\to\R^n$, $\lambda\colon\R\to[0,\infty)$, and
$B\colon\R^n\to\mathbb S$
with $F|_{V\cap U}=f|_{V\cap U}$,
with the property that for all $v\in U\backslash\{\vartheta\}$ it holds that
$P(v)=|v-\vartheta|^{-1}(v-\vartheta)$ and
$\opnorm{B(v)-\mathbb B}\indicator{V}(v)\leq c|v-\vartheta|^\varepsilon$,
and with the property that for Lebesgue-almost all $t\in I$ it holds that
\begin{equation}
  \label{maintheorem:ambient2.flow}
  \EMPH{\theta'(t)=-\lambda(t)B(\theta(t))(\nabla F)(\theta(t))}
\end{equation}
we have that
\begin{enumerate}[label=(\roman*)]
\item \label{maintheorem:ambient2.item1}
there exists a finite $\mathbb{A}\subseteq I$ such that
$I\backslash \mathbb{A}\ni t\mapsto(P\circ\theta)(t)\in\R^n$
is locally absolutely continuous,
\item \label{maintheorem:ambient2.item2}
it holds that $P\circ\theta$ is Lebesgue-almost everywhere differentiable, and
\item \label{maintheorem:ambient2.item3} it holds that
\begin{equation}
  \label{maintheorem:ambient2.length}
  \int_I |(P\circ\theta)'(s)|\,
    \indicator{[0,\delta]}(|\theta(s)-\vartheta|)\,\mathrm{d}s\leq L.
\end{equation}
\end{enumerate}
\end{athm}

\begin{aproof}
  Let $I\subseteq\R$ be an interval with nonempty interior,
  let $U\subseteq\R^n$ be open,
  let $F\in C^1(U,\R)$,
  let $\theta\colon I\to U$ be locally absolutely continuous,
  let $P\colon U\to\R^n$, $\lambda\colon\R\to[0,\infty)$, and
  $B\colon\R^n\to\mathbb S$ be measurable,
  assume that
  \begin{eqn}
    \llabel{setting:admissible-potential}
    F|_{V\cap U}=f|_{V\cap U},
  \end{eqn}
  assume for all $v\in U\backslash\{\vartheta\}$ that
  \begin{eqn}
    \llabel{setting:admissible-projection-metric}
    P(v)=\frac{v-\vartheta}{\norm{v-\vartheta}}
    \qquad\text{and}\qquad
    \opnorm{B(v)-\mathbb B}\indicator{V}(v)
      \leq c\norm{v-\vartheta}^{\varepsilon},
  \end{eqn}
  and assume that for Lebesgue-almost all $t\in I$ it holds that
  $\theta$ is differentiable at $t$ and that
  \begin{eqn}
    \llabel{setting:admissible-flow}
    \theta'(t)=-\lambda(t)B(\theta(t))\nabla F(\theta(t))
    \dott
  \end{eqn}
  \startnewargseq
  In the following let $E\colon I\to\R$ satisfy for all $t\in I$
  that $E(t)=F(\theta(t))$, and let $H\subseteq I$ satisfy
  \begin{eqn}
    \llabel{setting:fibre}
    H=\theta^{-1}(\{\vartheta\})
    \dott
  \end{eqn}
  \startnewargseq
  \argument{
    \lref{setting:admissible-flow};
    the fact that $\theta$ is locally absolutely continuous;
    the fact that $F\in C^1(U,\R)$;
    the fact that $E=F\circ\theta$;
    the fact that for all $t\in\R$ it holds that $\lambda(t)\geq0$;
    the fact that for all $v\in\R^n$ it holds that $B(v)$ is positive-definite;
    the chain rule;
  }{
    that $E$ is locally absolutely continuous and non-increasing
    and that for Lebesgue-almost all $t\in I$ it holds that
    $E$ and $\theta$ are differentiable at $t$ and that
    \begin{eqn}
      \llabel{eq:energy}
      E'(t)=-\lambda(t)\langle\nabla F(\theta(t)),
        B(\theta(t))\nabla F(\theta(t))\rangle\leq0
      \dott
    \end{eqn}
  }
  \argument{
    \lref{eq:energy};
    the fact that $E=F\circ\theta$;
  }{
    that for all $u,v\in I$ with $u<v$ and $\theta(u)=\theta(v)$
    it holds that $E|_{[u,v]}$ is constant and that for
    Lebesgue-almost all $t\in[u,v]$ it holds that
    \begin{eqn}
      \llabel{eq:zero-dissipation}
      \lambda(t)\langle\nabla F(\theta(t)),
        B(\theta(t))\nabla F(\theta(t))\rangle=0
      \dott
    \end{eqn}
  }
  \argument{
    \lref{eq:zero-dissipation};
    \lref{setting:admissible-flow};
    the fact that for all $v\in\R^n$ it holds that $B(v)$ is positive-definite;
  }{
    that for all $u,v\in I$ with $u<v$ and $\theta(u)=\theta(v)$
    it holds that for Lebesgue-almost all $t\in[u,v]$ it holds that
    $\theta$ is differentiable at $t$ and that
    \begin{eqn}
      \llabel{eq:zero-speed}
      \theta'(t)=-\lambda(t)B(\theta(t))\nabla F(\theta(t))=0
      \dott
    \end{eqn}
  }
  \argument{
    \lref{eq:zero-speed};
    the fact that $\theta$ is locally absolutely continuous;
    the fundamental theorem of calculus;
  }{
    that for all $u,v,t\in I$ with $u\leq t\leq v$ and
    $\theta(u)=\theta(v)$ it holds that
    \begin{eqn}
      \llabel{eq:no-return}
      \theta(t)=\theta(u)+\int_u^t\theta'(s)\dd s=\theta(u)
      \dott
    \end{eqn}
  }
  \argument{
    \lref{eq:no-return};
    \lref{setting:fibre};
    the fact that $\theta$ is continuous;
  }{
    that $H$ is an interval and that $H$ is closed in $I$\dott
  }
  In the following let $\mathcal J\subseteq\mathcal P(I)$ and $\mathbb{A}\subseteq I$
  satisfy
  \begin{equation}
    \llabel{setting:pieces}
    \mathcal J=\begin{cases}
      \{I\}&\colon H=\emptyset\\
      \bigl\{\{t\in I\colon t<\inf H\},
        \{t\in I\colon t>\sup H\}\bigr\}\backslash\{\emptyset\}
        &\colon H\neq\emptyset
    \end{cases}
  \end{equation}
  and
  \begin{eqn}
    \llabel{setting:exceptional-set}
    \mathbb{A}=\begin{cases}
      \emptyset&\colon H=\emptyset\\
      \{\inf H,\sup H\}\cap I&\colon H\neq\emptyset
    \end{cases}
    \dott
  \end{eqn}
  \startnewargseq
  \argument{ \lref{setting:pieces};\lref{setting:exceptional-set};
    the fact that $H$ is an interval and that $H$ is closed in $I$;
  }{
    that
    \begin{eqn}
      \llabel{eq:pieces}
      \bigcup_{J\in\mathcal J}J=I\backslash H,\qquad
      \#\mathcal J\leq2,\qquad
      \mathbb{A}\subseteq H,\qquad\text{and}\qquad\#\mathbb{A}\leq2,
    \end{eqn}
    that $H\backslash \mathbb{A}$ is open in $I$, that for all $J\in\mathcal J$
    it holds that $J$ is a nonempty interval and that $J$ is open in $I$,
    and that for all $J,K\in\mathcal J$ with $J\neq K$ it holds that
    $J\cap K=\emptyset$\dott
  }
  \argument{
    \lref{eq:pieces};
    \lref{setting:admissible-projection-metric};
    \lref{setting:fibre};
    the fact that $\theta$ is locally absolutely continuous;
    the chain rule;
  }{
    that $(P\circ\theta)|_{I\backslash \mathbb{A}}$ is locally absolutely continuous\dott
    \llabel{arg:local-absolute-continuity}
  }
  \argument{
    \lref{arg:local-absolute-continuity};
    \lref{eq:pieces};
  }[verbs=ep]{
    \lref{item1}\dott
  }
  \startnewargseq
  \argument{
    \lref{eq:pieces};
    the fact that $(P\circ\theta)|_{I\backslash \mathbb{A}}$ is locally absolutely continuous;
  }[verbs=ep]{
    \lref{item2}\dott
  }
  \startnewargseq
  \argument{
    \lref{eq:pieces};
    the fact that $(P\circ\theta)|_H$ is constant;
  }{
    that for Lebesgue-almost all $t\in H$
    it holds that $P\circ\theta$ is differentiable at $t$ and that
    \begin{eqn}
      \llabel{eq:zero-fibre-derivative}
      (P\circ\theta)'(t)=0
      \dott
    \end{eqn}
  }
  In our proof of \lref{item3}
  we distinguish between the case $n=1$ and the case $n\geq2$\dott
  \startnewargseq
  We first prove \lref{item3} in the case that
  \begin{eqn}
    \llabel{one:case}
    n=1
    \dott
  \end{eqn}
  \startnewargseq
  In the following let $\delta,L\in(0,\infty)$ satisfy
  \begin{eqn}
    \llabel{one:setting:constants}
    \delta=\varepsilon/2
    \qquad\text{and}\qquad
    L=1
    \dott
  \end{eqn}
  \startnewargseq
  \argument{
    \lref{setting:admissible-projection-metric};
    \lref{setting:fibre};
    \lref{setting:pieces};
    \lref{one:case};
    the fact that for all $J\in\mathcal J$ it holds that $J$ is an interval;
    the fact that $\theta$ is continuous;
    the intermediate value theorem;
  }{
    that for all $J\in\mathcal J$ it holds that
    $(P\circ\theta)|_J$ is constant with value in $\{-1,1\}$\dott
    \llabel{one:arg:constancy}
  }
  \argument{
    \lref{one:arg:constancy};
    \lref{eq:pieces};
    \lref{eq:zero-fibre-derivative};
    the fact that for all $J\in\mathcal J$ it holds that $J$ is open in $I$;
  }{
    that for Lebesgue-almost all $t\in I$ it holds that
    $P\circ\theta$ is differentiable at $t$ and that
    \begin{eqn}
      \llabel{one:eq:zero-derivative}
      (P\circ\theta)'(t)=0
      \dott
    \end{eqn}
  }
  \argument{
    \lref{one:eq:zero-derivative};
    \lref{one:setting:constants};
  }{
    that
    \begin{eqn}
      \llabel{one:eq:integral}
      \int_I\norm{(P\circ\theta)'(t)}
        \indicator{[0,\delta]}(\norm{\theta(t)-\vartheta})\dd t=0\leq L
      \dott
    \end{eqn}
  }
  This completes our proof of \lref{item3}
  in the case $n=1$\dott
  \startnewargseq
  We now prove \lref{item3} in the case that
  \begin{eqn}
    \llabel{case:constants-higher}
    n\geq2
    \dott
  \end{eqn}
  \startnewargseq
  \argument{
    the fact that $\mathbb B$ is symmetric and positive-definite;
    the spectral theorem for symmetric matrices;
  }{
    that there exists an invertible $Y\in\R^{n\times n}$ which satisfies that
    \begin{eqn}
      \llabel{eq:factorization}
      YY^{\mathsf T}=\mathbb B
      \qquad\text{and}\qquad
      Y^{-1}\mathbb B(Y^{-1})^{\mathsf T}=\Id
      \dott
    \end{eqn}
  }
  In the following let $q\in\mathbb Q\cap(0,\varepsilon]$,
  $\kappa\in(0,\infty)$, $R\in(0,1)$ satisfy
  \begin{eqn}
    \llabel{setting:normalization-radius}
    \kappa=c\opnorm{Y^{-1}}^2\opnorm{Y}^{\varepsilon},
    \qquad
    \opnorm{Y}R<\varepsilon,
    \qquad\text{and}\qquad
    \kappa R^{\varepsilon}\leq\frac12,
  \end{eqn}
  let $\mathbb{V},\mathcal S\subseteq\R^n$ satisfy
  \begin{eqn}
    \llabel{setting:normalized-domains}
    \mathbb{V}=\{z\in\R^n\colon0<\norm{z}<R\}
    \qquad\text{and}\qquad
    \mathcal S=\{\omega\in\R^n\colon\norm{\omega}=1\},
  \end{eqn}
  let $g\colon\mathbb{V}\to\R$ and $W\subseteq\mathbb{V}$
  satisfy for all $z\in\mathbb{V}$ that
  \begin{equation}
    \llabel{setting:normalized-potential}
    g(z)=f(\vartheta+Yz)
    \qquad\text{and}\qquad W=\bigcup\bigl\{\mathcal U\subseteq\mathbb{V}\colon
      \mathcal U\text{ is open and }g|_{\mathcal U}\in C^1(\mathcal U,\R)\bigr\},
  \end{equation}
  and let $\mathbf k\colon(0,R)\to[0,\infty)$ satisfy for all
  $s\in(0,R)$ that
  \begin{eqn}
    \llabel{setting:mixing-density}
    \mathbf k(s)=2\kappa s^{q-1}
    \dott
  \end{eqn}
  \startnewargseq
  \argument{
    \lref{eq:factorization};
    \lref{setting:normalization-radius};
    \lref{setting:normalized-domains};
    \lref{setting:normalized-potential};
    the fact that $f\in C(V,\R)$ is definable;
  }{
    that $g$ is continuous and definable, that $W$ is definable and open,
    that $g|_W$ is definable, and that
    \begin{eqn}
      \llabel{eq:normalized-regularity}
      \{\vartheta+Yz\colon z\in\mathbb{V}\}\subseteq V
      \qquad\text{and}\qquad
      g|_W\in C^1(W,\R)
      \dott
    \end{eqn}
  }
  \argument{
    \lref{setting:normalization-radius};
    \lref{setting:mixing-density};
  }{
    that $\mathbf k$ is nonnegative, locally bounded, and definable
    and that
    \begin{eqn}
      \llabel{eq:mixing-density}
      \int_0^R\mathbf k(s)\dd s=\frac{2\kappa R^q}{q}<\infty
      \dott
    \end{eqn}
  }
  \argument{
    \lref{setting:normalization-radius};
    \lref{setting:normalized-domains};
  }{
    that for all $s\in(0,R)$, $\mathbf A\in\mathbb S$ with
    $\opnorm{\mathbf A-\Id}\leq\kappa s^{\varepsilon}$ it holds that
    \begin{eqn}
      \llabel{eq:normalized-condition}
      \max_{x\in\mathcal S}\langle\mathbf A x,x\rangle
        \leq1+\opnorm{\mathbf A-\Id}\leq\frac32
        \leq3(1-\opnorm{\mathbf A-\Id})
        \leq3\min_{x\in\mathcal S}\langle\mathbf A x,x\rangle
      \dott
    \end{eqn}
  }
  \argument{
    \lref{eq:normalized-condition};
    \lref{setting:normalization-radius};
    \lref{setting:mixing-density};
    the fact that for all $x\in\mathcal S$ it holds that
    $\opnorm{\Id-xx^{\mathsf T}}\leq1$;
  }{
    that for all $s\in(0,R)$, $\mathbf A\in\mathbb S$, $x\in\mathcal S$
    with $\opnorm{\mathbf A-\Id}\leq\kappa s^{\varepsilon}$ it holds that
    \begin{eqn}
      \llabel{eq:normalized-mixing}
      \frac{\norm{\mathbf A x-\langle\mathbf A x,x\rangle x}}
        {\langle\mathbf A x,x\rangle}
        =\frac{\norm{(\Id-xx^{\mathsf T})(\mathbf A-\Id)x}}
          {\langle\mathbf A x,x\rangle}
        \leq2\opnorm{\mathbf A-\Id}\leq2\kappa s^{\varepsilon}
        \leq2\kappa s^q=s\mathbf k(s)
      \dott
    \end{eqn}
  }
  \argument{
    \lref{eq:normalized-mixing};
    \lref{case:constants-higher};
    \lref{setting:normalization-radius};
    \lref{setting:normalized-domains};
    \lref{eq:normalized-regularity};
    \lref{eq:mixing-density};
    \lref{eq:normalized-condition};
    \itref{prop:ambient-mixing}{item:angular}
    (applied with $K\curvearrowleft3$, $V\curvearrowleft\mathbb{V}$,
    $\mathbb S\curvearrowleft\mathcal S$, $g\curvearrowleft g|_W$
    in the notation of \cref{prop:ambient-mixing});
  }{
    that there exist $\eta\in(0,R)$, $T\in(0,\infty)$ which satisfy that
    for every open $\mathcal O\subseteq\R^n$,
    every $G\in C^1(\mathcal O,\R)$, every nonempty interval $J\subseteq\R$,
    every locally absolutely continuous $y\colon J\to\mathcal O\backslash\{0\}$,
    every Lebesgue measurable $\mu\colon J\to[0,\infty)$,
    every Lebesgue measurable $\mathbf A\colon J\to\R^{n\times n}$,
    and every $\omega\colon J\to\mathcal S$ with the property that
    \begin{eqn}
      \llabel{ass:angular-extension}
      \mathcal O\cap\mathbb{V}\subseteq W
      \qquad\text{and}\qquad
      G|_{\mathcal O\cap\mathbb{V}}=g|_{\mathcal O\cap\mathbb{V}},
    \end{eqn}
    with the property that for all $t\in J$ it holds that
    $\omega(t)=y(t)/\norm{y(t)}$, with the property that for
    Lebesgue-almost all $t\in J$ it holds that $y$ is differentiable
    at $t$, that $\mathbf A(t)$ is symmetric and positive-definite, and that
    \begin{eqn}
      \llabel{ass:angular-flow}
      y'(t)=-\mu(t)\mathbf A(t)\nabla G(y(t)),
    \end{eqn}
    and with the property that for Lebesgue-almost all
    $t\in\{s\in J\colon\norm{y(s)}<R\}$ it holds that
    $\opnorm{\mathbf A(t)-\Id}\leq\kappa\norm{y(t)}^{\varepsilon}$
    we have that
    \begin{eqn}
      \llabel{eq:uniform-angular}
      \int_{\{t\in J\colon\norm{y(t)}\leq\eta\}}
        \norm{\omega'(t)}\dd t\leq T
      \dott
    \end{eqn}
  }
  In the following let $\delta,L\in(0,\infty)$ satisfy
  \begin{eqn}
    \llabel{setting:constants}
    \delta=\frac{\eta}{\opnorm{Y^{-1}}}
    \qquad\text{and}\qquad
    L=2\opnorm{Y}\opnorm{Y^{-1}}T,
  \end{eqn}
  let $\mathcal O=\{z\in\R^n\colon\vartheta+Yz\in U\}$,
  and let $G\colon\mathcal O\to\R$ satisfy for all $z\in\mathcal O$ that
  \begin{eqn}
    \llabel{setting:extension}
    G(z)=F(\vartheta+Yz)
    \dott
  \end{eqn}
  \startnewargseq
  \argument{
    \lref{setting:extension};
    the fact that $U$ is open;
    the fact that $F\in C^1(U,\R)$;
    the chain rule;
  }{
    that $\mathcal O$ is open, that $G\in C^1(\mathcal O,\R)$, and that
    for all $z\in\mathcal O$ it holds that
    \begin{eqn}
      \llabel{eq:extension-gradient}
      \nabla G(z)=Y^{\mathsf T}\nabla F(\vartheta+Yz)
      \dott
    \end{eqn}
  }
  \argument{
    \lref{eq:extension-gradient};
    \lref{setting:admissible-potential};
    \lref{setting:normalized-potential};
    \lref{eq:normalized-regularity};
    \lref{setting:extension};
  }{
    that
    \begin{eqn}
      \llabel{eq:compatibility}
      \mathcal O\cap\mathbb{V}\subseteq W
      \qquad\text{and}\qquad
      G|_{\mathcal O\cap\mathbb{V}}=g|_{\mathcal O\cap\mathbb{V}}
      \dott
    \end{eqn}
  }

  We claim that for all $J\in\mathcal J$ it holds that
  \begin{eqn}
    \llabel{claim:component-bound}
    \int_J\norm{(P\circ\theta)'(t)}
      \indicator{[0,\delta]}(\norm{\theta(t)-\vartheta})\dd t
      \leq\opnorm{Y}\opnorm{Y^{-1}}T
    \dott
  \end{eqn}
  \startnewargseq
  To prove \lref{claim:component-bound}, let $J\in\mathcal J$
  and let $z\colon J\to\mathcal O\backslash\{0\}$,
  $\omega,\mathbf w\colon J\to\mathcal S$, and
  $\mathbf B\colon J\to\R^{n\times n}$ satisfy for all $t\in J$ that
  \begin{gather}
    \llabel{setting:curve}
    z(t)=Y^{-1}(\theta(t)-\vartheta),\qquad
    \omega(t)=\frac{z(t)}{\norm{z(t)}},\\
    \llabel{setting:projected-curve-metric}
    \mathbf w(t)=P(\theta(t))
    \qquad\text{and}\qquad
    \mathbf B(t)=Y^{-1}B(\theta(t))(Y^{-1})^{\mathsf T}
    \dott
  \end{gather}
  \startnewargseq
  \argument{
    \lref{setting:admissible-projection-metric};
    \lref{setting:curve};
    \lref{setting:projected-curve-metric};
    the fact that $\theta$ is locally absolutely continuous;
    the fact that $B\colon\R^n\to\mathbb S$ is measurable;
    the chain rule;
  }{
    that $z$, $\omega$, and $\mathbf w$ are locally absolutely continuous,
    that $\mathbf B$ is measurable, and that for all $t\in J$
    it holds that $\mathbf B(t)$ is symmetric and positive-definite\dott
  }
  \argument{
    \lref{setting:admissible-flow};
    \lref{eq:extension-gradient};
    \lref{setting:curve};
    \lref{setting:projected-curve-metric};
    the chain rule;
  }{
    that for Lebesgue-almost all $t\in J$ it holds that $\theta$ and $z$
    are differentiable at $t$ and that
    \begin{eqn}
      \llabel{eq:normalized-flow}
      z'(t)=Y^{-1}\theta'(t)
        =-\lambda(t)Y^{-1}B(\theta(t))\nabla F(\theta(t))
        =-\lambda(t)\mathbf B(t)\nabla G(z(t))
      \dott
    \end{eqn}
  }
  \argument{
    \lref{setting:admissible-projection-metric};
    \lref{eq:factorization};
    \lref{setting:normalization-radius};
    \lref{setting:curve};
    \lref{setting:projected-curve-metric};
  }{
    that for all $t\in J$ with $\norm{z(t)}<R$ it holds that
    $\theta(t)\in V$ and that
    \begin{eqn}[a]
      \llabel{eq:metric}
      \opnorm{\mathbf B(t)-\Id}
      =\opnorm{Y^{-1}(B(\theta(t))-\mathbb B)(Y^{-1})^{\mathsf T}}\leq c\opnorm{Y^{-1}}^2\norm{\theta(t)-\vartheta}^{\varepsilon}
        \leq\kappa\norm{z(t)}^{\varepsilon}
      \dott
    \end{eqn}
  }
  \argument{
    \lref{eq:metric};
    \lref{eq:uniform-angular};
    \lref{eq:extension-gradient};
    \lref{eq:compatibility};
    \lref{setting:curve};
    \lref{eq:normalized-flow};
    the fact that $z$ is locally absolutely continuous;
    the fact that $\lambda|_J$ and $\mathbf B$ are measurable;
    the fact that for all $t\in J$ it holds that
    $\mathbf B(t)$ is symmetric and positive-definite;
    the chain rule;
  }{
    that
    \begin{eqn}[a]
      \llabel{eq:angular}
      \int_{\{t\in J\colon\norm{z(t)}\leq\eta\}}\norm{\omega'(t)}\dd t
        &=\int_{\{t\in J\colon\norm{z(t)}\leq\eta\}}
          \norm[\bigg]{\frac{\ddd}{\ddd t}
            \biggl(\frac{z(t)}{\norm{z(t)}}\biggr)}\dd t\\
        &=\int_{\{t\in J\colon\norm{z(t)}\leq\eta\}}
          \frac{\norm{(\Id-\omega(t)(\omega(t))^{\mathsf T})z'(t)}}
            {\norm{z(t)}}\dd t\leq T
      \dott
    \end{eqn}
  }
  \argument{
    \lref{setting:admissible-projection-metric};
    \lref{setting:curve};
    \lref{setting:projected-curve-metric};
  }{
    that for all $t\in J$ it holds that
    \begin{eqn}
      \llabel{eq:transfer-coordinates}
      \norm{z(t)}\leq\opnorm{Y^{-1}}\norm{\theta(t)-\vartheta}
      \qquad\text{and}\qquad
      \mathbf w(t)=\frac{Y\omega(t)}{\norm{Y\omega(t)}}
      \dott
    \end{eqn}
  }
  \argument{
    \lref{eq:transfer-coordinates};
    the fact that $\omega$ and $\mathbf w$ are locally absolutely continuous;
    the fact that for all $x\in\mathcal S$ it holds that
    $\norm{Yx}\geq\opnorm{Y^{-1}}^{-1}$ and that
    $\opnorm{\Id-xx^{\mathsf T}}\leq1$;
    the chain rule;
  }{
    that for Lebesgue-almost all $t\in J$ it holds that $\omega$
    and $\mathbf w$ are differentiable at $t$ and that
    \begin{eqn}
      \llabel{eq:transfer-derivative}
      \norm{\mathbf w'(t)}
        =\norm[\bigg]{\frac{(\Id-\mathbf w(t)(\mathbf w(t))^{\mathsf T})
            Y\omega'(t)}{\norm{Y\omega(t)}}}
        \leq\opnorm{Y}\opnorm{Y^{-1}}\norm{\omega'(t)}
      \dott
    \end{eqn}
  }
  \argument{
    \lref{eq:transfer-derivative};
    \lref{setting:constants};
    \lref{setting:projected-curve-metric};
    \lref{eq:angular};
    \lref{eq:transfer-coordinates};
  }{
    that
    \begin{eqn}[a]
      \llabel{eq:component-bound}
      \int_J\norm{(P\circ\theta)'(t)}
        \indicator{[0,\delta]}(\norm{\theta(t)-\vartheta})\dd t
      &=\int_{\{t\in J\colon\norm{\theta(t)-\vartheta}\leq\delta\}}
          \norm{\mathbf w'(t)}\dd t\\
      &\leq\opnorm{Y}\opnorm{Y^{-1}}
          \int_{\{t\in J\colon\norm{z(t)}\leq\opnorm{Y^{-1}}\delta\}}
            \norm{\omega'(t)}\dd t\\
      &=\opnorm{Y}\opnorm{Y^{-1}}
          \int_{\{t\in J\colon\norm{z(t)}\leq\eta\}}\norm{\omega'(t)}\dd t
        \leq\opnorm{Y}\opnorm{Y^{-1}}T
      \dott
    \end{eqn}
  }
  This and the arbitrariness of $J\in\mathcal J$ complete our proof
  of \lref{claim:component-bound}\dott

  \startnewargseq
  \argument{
    \lref{eq:pieces};
    \lref{eq:zero-fibre-derivative};
    \lref{setting:constants};
    \lref{claim:component-bound};
  }{
    that
    \begin{eqn}[a]
      \llabel{eq:total-bound}
      \int_I\norm{(P\circ\theta)'(t)}
        \indicator{[0,\delta]}(\norm{\theta(t)-\vartheta})\dd t
      &=\sum_{J\in\mathcal J}\int_J\norm{(P\circ\theta)'(t)}
          \indicator{[0,\delta]}(\norm{\theta(t)-\vartheta})\dd t\\
      &\leq(\#\mathcal J)\opnorm{Y}\opnorm{Y^{-1}}T
        \leq2\opnorm{Y}\opnorm{Y^{-1}}T=L
      \dott
    \end{eqn}
  }
  This completes our proof of \lref{item3}
  in the case $n\geq2$\dott
  \startnewargseq
\end{aproof}

\subsection*{Acknowledgements}
This work has been partially supported by the National Natural Science Foundation of China (NSFC) under grant number W2531010.
This work has also been supported by the Ministry of Culture and Science NRW as part of the Lamarr Fellow Network.
We also gratefully acknowledge the Cluster of Excellence EXC 2044/2-390685587, Mathematics M{\"u}nster: Dynamics-Geometry-Structure funded by the Deutsche Forschungsgemeinschaft (DFG, German Research Foundation).

\subsection*{Use of large language models}
The authors gratefully acknowledge extensive assistance from {\sc OpenAI}'s {\sc GPT-5.6 Sol} and {\sc GPT-6 Astra/Sol} in the development of this work. These large language models were used via {\sc ChatGPT}/{\sc Codex} as research collaborators for developing mathematical arguments and ideas, exploring proof strategies, testing arguments, identifying gaps and possible counterexamples, refining proofs and exposition, and reviewing the literature. All mathematical statements, claims, and proofs were independently checked by the authors and all of the sentences, statements, and arguments in this article are the sole responsibility of the authors.

\bibliographystyle{abbrv-doi}

\input{main.bbl}
\end{document}

%% file: commands.tex
\ExplSyntaxOn

\bool_new:N \g_kuckuck_debug_colors_bool
\bool_new:N \g_kuckuck_debug_logs_bool

\seq_new:N \g_cflist_loaded
\seq_new:N \g_cflist_pending

\NewDocumentCommand{\cfadd} { m } {
  \seq_if_in:NnF \g_cflist_loaded { #1 } {
    \seq_if_in:NnF \g_cflist_pending { #1 } {
      \seq_gput_right:Nn \g_cflist_pending { #1 }
    }
  }
}

\NewDocumentCommand{\cfconsiderloaded} { m } {
  \seq_gput_right:Nn \g_cflist_loaded {#1}
}

\NewDocumentCommand{\cfremove} { m } {
  \seq_gremove_all:Nn \g_cflist_pending { #1 }
}

\NewDocumentCommand{\cfload} { o } {
  \seq_if_empty:NTF \g_cflist_pending {
    \IfValueTF{#1}{\ignorespaces}{\unskip}
  } {
    (cf.\ \cref{\seq_use:Nn \g_cflist_pending {,}})\IfValueTF{#1}{#1~}{\unskip}
    \seq_gconcat:NNN \g_cflist_loaded \g_cflist_loaded \g_cflist_pending
    \seq_gclear:N \g_cflist_pending
    \IfValueT{#1}{\ignorespaces}
  }
}

\NewDocumentCommand{\cfclear} {} {
  \seq_gclear:N \g_cflist_loaded
  \seq_gclear:N \g_cflist_pending
}

\NewDocumentCommand{\cfout} { o } {
  \seq_if_empty:NTF \g_cflist_pending {\unskip\IfValueT{#1}{\ignorespaces}} {
    (cf.\ \cref{\seq_use:Nn \g_cflist_pending {,}})\IfValueTF{#1}{#1~}{\unskip}
    \seq_gclear:N \g_cflist_pending
    \IfValueT{#1}{\ignorespaces}
  }
}

\NewDocumentCommand{\ifnocf} { m } {
  \seq_if_empty:NT \g_cflist_pending { #1 }
}

\ExplSyntaxOff

\ExplSyntaxOn

\bool_new:N \g_noteobserve

\NewDocumentCommand{\setnote}{}{
  \bool_gset_true:N \g_noteobserve
}

\NewDocumentCommand{\setobserve}{}{
  \bool_gset_false:N \g_noteobserve
}

\NewDocumentCommand{\nobs}{ o }{
  \IfValueT{#1}{
    \str_if_eq:noTF {note} {#1} {
      \bool_gset_true:N \g_noteobserve
    } {
      \str_if_eq:noTF {Note} {#1} {
        \bool_gset_true:N \g_noteobserve
      } {
        \bool_gset_false:N \g_noteobserve
      }
    }
  }
  \bool_if:nTF { \g_noteobserve } {
    \bool_gset_false:N \g_noteobserve
    note
  } {
    \bool_gset_true:N \g_noteobserve
    observe
  }
  \IfValueF{#1}{~}
}

\NewDocumentCommand{\Nobs}{ o }{
  \IfValueT{#1}{
    \str_if_eq:noTF {note} {#1} {
      \bool_gset_true:N \g_noteobserve
    } {
      \str_if_eq:noTF {Note} {#1} {
        \bool_gset_true:N \g_noteobserve
      } {
        \bool_gset_false:N \g_noteobserve
      }
    }
  }
  \bool_if:nTF { \g_noteobserve } {
    \bool_gset_false:N \g_noteobserve
    Note
  } {
    \bool_gset_true:N \g_noteobserve
    Observe
  }
  \IfValueF{#1}{~}
}

\ExplSyntaxOff

\ExplSyntaxOn

\bool_new:N \g_hencetherefore

\NewDocumentCommand{\hence}{ o }{
  \IfValueT{#1}{
    \str_if_eq:noTF {hence} {#1} {
      \bool_gset_true:N \g_hencetherefore
    } {
      \str_if_eq:noTF {Hence} {#1} {
        \bool_gset_true:N \g_hencetherefore
      } {
        \bool_gset_false:N \g_hencetherefore
      }
    }
  }
  \bool_if:nTF { \g_hencetherefore } {
    \bool_gset_false:N \g_hencetherefore
    hence
  } {
    \bool_gset_true:N \g_hencetherefore
    therefore
  }
  \IfValueF{#1}{~}
}

\NewDocumentCommand{\Hence}{ o }{
  \IfValueT{#1}{
    \str_if_eq:noTF {hence} {#1} {
      \bool_gset_true:N \g_hencetherefore
    } {
      \str_if_eq:noTF {Hence} {#1} {
        \bool_gset_true:N \g_hencetherefore
      } {
        \bool_gset_false:N \g_hencetherefore
      }
    }
  }
  \bool_if:nTF { \g_hencetherefore } {
    \bool_gset_false:N \g_hencetherefore
    Hence,~we~obtain
  } {
    \bool_gset_true:N \g_hencetherefore
    Therefore,~we~obtain
  }
  \IfValueF{#1}{~}
}

\ExplSyntaxOff

\ExplSyntaxOn

\int_new:N \g_furthermore

\NewDocumentCommand{\Moreover}{ o o }{
  \IfValueT{#1}{
    \str_case:nn {#1} {
      {Next} {\int_gset:Nn {\g_furthermore} {0}}
      {Furthermore} {\int_gset:Nn {\g_furthermore} {1}}
      {Moreover} {\int_gset:Nn {\g_furthermore} {2}}
      {In~addition} {\int_gset:Nn {\g_furthermore} {3}}
      {note} {\bool_gset_true:N \g_noteobserve}
      {observe} {\bool_gset_false:N \g_noteobserve}
    }
    \IfValueT{#2}{
      \str_case:nn {#2} {
        {Next} {\int_gset:Nn {\g_furthermore} {0}}
        {Furthermore} {\int_gset:Nn {\g_furthermore} {1}}
        {Moreover} {\int_gset:Nn {\g_furthermore} {2}}
        {In~addition} {\int_gset:Nn {\g_furthermore} {3}}
        {note} {\bool_gset_true:N \g_noteobserve}
        {observe} {\bool_gset_false:N \g_noteobserve}
      }
    }
  }
  \int_case:nn { \int_mod:nn {\g_furthermore} {4} } {
    { 0 } { Next~\nobs that}
    { 1 } { Furthermore,~\nobs that}
    { 2 } { Moreover,~\nobs that}
    { 3 } { In~addition,~\nobs that}
  }
  \int_incr:N \g_furthermore
  \peek_charcode:NTF , {  } { ~ }
}

\ExplSyntaxOff

\ExplSyntaxOn

\seq_const_from_clist:Nn \g_prove_mru {
  establish,
  demonstrate,
  prove,
  show,
  imply,
  ensure
}

\prop_new:N \l__verbs
\prop_put:Nnn \l__verbs {show} {shows}
\prop_put:Nnn \l__verbs {imply} {implies}
\prop_put:Nnn \l__verbs {demonstrate} {demonstrates}
\prop_put:Nnn \l__verbs {prove} {proves}
\prop_put:Nnn \l__verbs {establish} {establishes}
\prop_put:Nnn \l__verbs {ensure} {ensures}
\prop_put:Nnn \l__verbs {assure} {assures}

\tl_new:N \g_wordtmp
\seq_new:N \l_mytmps

\cs_generate_variant:Nn \str_if_in:nnTF { nVTF }
\cs_generate_variant:Nn \str_if_in:nnTF { xVTF }

\NewDocumentCommand{\prove}{ o }{
  \IfValueTF{#1}{
    \seq_clear:N \l_mytmps
    \seq_map_inline:Nn \g_prove_mru {
      \str_if_eq:nnTF {##1} {ensure} {
        \str_set:Nn \l_temps {n}
      } {
        \str_set:Nx \l_temps {\str_head_ignore_spaces:n {##1}}
      }
      \str_if_in:xVTF {#1} \l_temps {
        \seq_put_right:Nn \l_mytmps {##1}
      } { }
    }
    \seq_get_right:NN \l_mytmps \g_wordtmp
  } {
    \seq_get_right:NN \g_prove_mru \g_wordtmp
  }
  \tl_use:N \g_wordtmp
  \IfValueTF{#1}{}{~}
  \seq_gput_left:NV \g_prove_mru \g_wordtmp
  \seq_gremove_duplicates:N \g_prove_mru
}

\NewDocumentCommand{\proves}{ o }{
  \IfValueTF{#1}{
    \seq_clear:N \l_mytmps
    \seq_map_inline:Nn \g_prove_mru {
      \str_if_eq:nnTF {##1} {ensure} {
        \str_set:Nn \l_temps {n}
      } {
        \str_set:Nx \l_temps {\str_head_ignore_spaces:n {##1}}
      }
      \str_if_in:xVTF {#1} \l_temps {
        \seq_put_right:Nn \l_mytmps {##1}
      } { }
    }
    \seq_get_right:NN \l_mytmps \g_wordtmp
  } {
    \seq_get_right:NN \g_prove_mru \g_wordtmp
  }
  \str_set:NV \l_tmpa_str \g_wordtmp
  \prop_get:NVN \l__verbs \l_tmpa_str \l_tmpa_tl
  \tl_use:N \l_tmpa_tl
  \IfValueTF{#1}{}{~}
  \seq_gput_left:NV \g_prove_mru \g_wordtmp
  \seq_gremove_duplicates:N \g_prove_mru
}

\ExplSyntaxOff

\ExplSyntaxOn

\newcommand{\llabel}[1]{\savelabel{#1}\label{\loc.#1}\ignorespaces}

\clist_new:N \l_localreflist
\clist_new:N \l_reflist

\NewDocumentCommand{\lref} { m } {
  \clist_set:No \l_localreflist {#1}
  \clist_clear:N \l_reflist
  \clist_map_inline:Nn \l_localreflist { \clist_put_right:Nn \l_reflist {\loc.##1} }
  \cref{\l_reflist}
}

\NewDocumentCommand{\Lref} { m } {
  \clist_set:No \l_localreflist {#1}
  \clist_clear:N \l_reflist
  \clist_map_inline:Nn \l_localreflist { \clist_put_right:Nn \l_reflist {\loc.##1} }
  \Cref{\l_reflist}
}

\NewDocumentCommand{\itref}{ m m }{
  \clist_set:No \l_localreflist {#2}
  \clist_clear:N \l_reflist
  \clist_map_inline:Nn \l_localreflist { \clist_put_right:Nn \l_reflist {#1.##1} }
  \cref{\l_reflist}~in~\cref{#1}
}

\ExplSyntaxOff

\ExplSyntaxOn

\seq_new:N \l_enum_seq
\int_new:N \l_num_items

\bool_new:N \g_commaused_bool

\providecommand{\comma}{}

\cs_new:Nn \enum_it:nn {
  \int_case:nnF {\l_num_items - #1} {
    {0} {
      \renewcommand{\comma}{}
      #2\space
    }
    {1} {
      \bool_gset_false:N \g_commaused_bool
      \renewcommand{\comma}{,~\bool_gset_true:N \g_commaused_bool}
      #2
      \bool_if:NTF \g_commaused_bool {} {,~}
      and~
    }
  } {
    \bool_gset_false:N \g_commaused_bool
    \renewcommand{\comma}{,~\bool_gset_true:N \g_commaused_bool}
    #2
    \bool_if:NTF \g_commaused_bool {} {,~}
  }
}

\cs_new:Nn \enum_it_U:nn {
  \int_case:nnF {\l_num_items - #1} {
    {0} {
      \renewcommand{\comma}{}
      #2
      \space
    }
    {1} {
      \bool_gset_false:N \g_commaused_bool
      \renewcommand{\comma}{,~\bool_gset_true:N \g_commaused_bool}
      #2
      \bool_if:NTF \g_commaused_bool {} {,~}
      and~
    }
  } {
    \bool_gset_false:N \g_commaused_bool
    \renewcommand{\comma}{,~\bool_gset_true:N \g_commaused_bool}
    \int_compare:nTF {#1=1} {
      \text_titlecase_first:n {#2}
    } {
      #2
    }
    \bool_if:NTF \g_commaused_bool {} {,~}
  }
}

\cs_generate_variant:Nn \tl_if_eq:nnTF {onTF}

\cs_new:Nn \enum:nnnn {
  \seq_set_split:Nnn \l_enum_seq ; {#1}
  \seq_remove_all:Nn \l_enum_seq { }
  \seq_remove_all:Nn \l_enum_seq {#2}
  \bool_if:NT \g_kuckuck_debug_logs_bool { \seq_log:N \l_enum_seq }
  \int_set:Nn \l_num_items {\seq_count:N \l_enum_seq}
  \bool_if:NT \g_kuckuck_debug_logs_bool { \int_log:N \l_num_items }
  \int_case:nnF {\l_num_items} {
    { 0 } { 0 }
    { 1 } {
      \IfBooleanTF{#4} {
        \tl_set:Nn \l_text_case_exclude_arg_tl {\cref}
        \bool_set_false:N \l_text_titlecase_check_letter_bool
        \text_titlecase_first:n {\seq_use:Nn \l_enum_seq {}}
      } {
        \seq_use:Nn \l_enum_seq {}
      }
      \space
      \tl_if_eq:onTF{#3}{-}{}{
        \bool_if:NTF \l_plural_bool {
          \prove[#3]~
        } {
          \proves[#3]~
        }
      }
    }
    { 2 } {
      \IfBooleanTF{#4} {
        \tl_set:Nn \l_text_case_exclude_arg_tl {\cref}
        \bool_set_false:N \l_text_titlecase_check_letter_bool
        \text_titlecase_first:n {\seq_item:Nn \l_enum_seq {1}}
        {} ~and~
        \seq_item:Nn \l_enum_seq {2}
      } {
        \seq_use:Nn \l_enum_seq {~and~}
      }
      \space
      \tl_if_eq:onTF{#3}{-}{}{
        \prove[#3]~
      }
    }
  } {
    \IfBooleanTF{#4} {
      \tl_set:Nn \l_text_case_exclude_arg_tl {\cref}
      \seq_indexed_map_function:NN \l_enum_seq \enum_it_U:nn
    } {
      \seq_indexed_map_function:NN \l_enum_seq \enum_it:nn
    }
    \tl_if_eq:onTF{#3}{-}{}{
      \prove[#3]~
    }
  }
}

\cs_generate_variant:Nn \enum:nnnn {nxnn}
\cs_generate_variant:Nn \enum:nnnn {nxxn}

\NewDocumentCommand{\enum}{O{} m O{-} s}{
  \IfBooleanTF{#4}{
    \enum:nxnn {#2} {#1} {sindep} \BooleanFalse
  } {
    \enum:nxxn {#2} {#1} {#3} \BooleanFalse
  }
}

\ExplSyntaxOff

\ExplSyntaxOn

\NewDocumentCommand{\dott}{}{\ifnocf{.}\space}

\bool_new:N \g_arg_start_bool
\bool_gset_true:N \g_arg_start_bool

\NewDocumentCommand{\startnewargseq}{}{\bool_gset_true:N \g_arg_start_bool \tl_set:Nn \g_label_tl {}}

\cs_generate_variant:Nn \seq_if_in:NnTF {NxTF}
\cs_generate_variant:Nn \seq_remove_all:Nn {Nx}

\int_new:N \l_random_int

\cs_generate_variant:Nn \tl_if_head_eq_catcode:nNTF {oNTF}
\cs_generate_variant:Nn \tl_if_head_eq_catcode:nNTF {VNTF}

\cs_generate_variant:Nn \tl_if_head_eq_catcode:nNTF {eNTF}
\cs_generate_variant:Nn \tl_log:n {o}
\cs_generate_variant:Nn \tl_log:n {f}
\cs_generate_variant:Nn \tl_log:n {x}
\cs_generate_variant:Nn \tl_log:n {e}

\cs_generate_variant:Nn \tl_if_in:nnTF {onTF}
\cs_generate_variant:Nn \tl_if_in:nnTF {nVTF}
\cs_generate_variant:Nn \tl_if_in:NnTF {NeTF}

\cs_generate_variant:Nn \tl_if_head_eq_meaning:nNTF {VNTF}

\seq_const_from_clist:Nn \g_arg_mru_this {
  Ahpr,
  Tapr,
  Ctapr,
  H
}

\seq_const_from_clist:Nn \g_arg_mru_nothis {
  Ia,
  Nwc,
  N,
  Itns,
  Fm,
  Itnswc,
  Mo
}

\seq_new:N \l_arg_seq
\tl_new:N \l_cons_tl
\tl_new:N \l_dummy_tl

\bool_new:N \g_debug_bool
\bool_gset_false:N \g_debug_bool

\bool_new:N \l_insidearg_bool

\bool_new:N \g_firstargletter_bool

\sys_gset_rand_seed:n {0903}

\bool_new:N \l_plural_bool
\tl_new:N \l_arg_verbs_tl

\NewDocumentCommand{\argument}{mom}{
  \bool_if:NT \g_kuckuck_debug_colors_bool { \color{blue} }
  \bool_set_false:N \l_plural_bool
  \tl_set:Nn \l_arg_verbs_tl {sindep}
  \keys_define:nn { benno/argument } {
    plural .value_forbidden:n = true,
    plural .code:n = {\bool_set_true:N \l_plural_bool},
    verbs .value_required:n = false,
    verbs .tl_set:N = \l_arg_verbs_tl,
  }
  \IfValueT{#2}{
    \keys_set:nn { benno/argument } {#2}
  }
  \bool_if:NT \g_kuckuck_debug_logs_bool { \bool_log:N \l_plural_bool }
  \bool_gset_true:N \l_insidearg_bool
  \seq_set_split:Nnn \l_arg_seq ; {#1}
  \seq_remove_all:Nn \l_arg_seq { }
  \bool_if:NT \g_kuckuck_debug_logs_bool { \seq_log:N \l_arg_seq }
  \tl_set:Nn \l_cons_tl {#3}
  \tl_trim_spaces:N \l_cons_tl
  \seq_if_in:NxTF \l_arg_seq {\lref{\g_label_tl}} {
    \seq_remove_all:Nx \l_arg_seq {\lref{\g_label_tl}}
    \seq_get_left:NNTF \l_arg_seq \l_dummy_tl {
      \tl_trim_spaces:N \l_dummy_tl
      \bool_gset_false:N \g_firstargletter_bool
      \tl_if_head_eq_catcode:nNTF \l_dummy_tl a {
        \bool_gset_true:N \g_firstargletter_bool
      } {
        \tl_if_head_eq_meaning:VNTF \l_dummy_tl {\cref} {
          \tl_set:Nx \l_tmpa_tl {\tl_tail:N \l_dummy_tl}
          \tl_set:Nx \l_tmpb_tl {\tl_head:N \l_tmpa_tl}
          \bool_gset_true:N \g_firstargletter_bool
          \tl_if_in:NeTF \l_tmpb_tl {lem\c_colon_str} {} {
            \tl_if_in:NeTF \l_tmpb_tl {thm\c_colon_str} {} {
              \tl_if_in:NeTF \l_tmpb_tl {prop\c_colon_str} {} {
                \tl_if_in:NeTF \l_tmpb_tl {cor\c_colon_str} {} {
                  \bool_gset_false:N \g_firstargletter_bool
                }
              }
            }
          }
        } {
        }
      }
      \bool_if:NTF \g_firstargletter_bool {
        \seq_set_eq:NN \l_tmpa_seq \g_arg_mru_this
        \seq_remove_all:Nn \l_tmpa_seq {H}
        \seq_get_right:NN \l_tmpa_seq \l_tmpa_tl
        \int_case:nnF {\seq_count:N \l_arg_seq} {
          {1} {
            \str_case:VnF {\l_tmpa_tl} {
              {Ahpr} {
                \bool_if:NT \g_debug_bool {C1.1}
                \seq_gput_left:Nn \g_arg_mru_this {Ahpr}
                \seq_gremove_duplicates:N \g_arg_mru_this
                \enum:nxnn {#1} {\lref{\g_label_tl}} {-} {\BooleanTrue}
                \hence~
                \bool_if:NTF \l_plural_bool {
                  \prove[\l_arg_verbs_tl]~\ignorespaces #3
                } {
                  \proves[\l_arg_verbs_tl]~\ignorespaces #3
                }
              }
              {Tapr} {
                \bool_if:NT \g_debug_bool {C1.2}
                \seq_gput_left:Nn \g_arg_mru_this {Tapr}
                \seq_gremove_duplicates:N \g_arg_mru_this
                \enum[\lref{\g_label_tl}]{
                  This;
                  #1
                }[\l_arg_verbs_tl]\ignorespaces #3
              }
              {Ctapr} {
                \bool_if:NT \g_debug_bool {C1.3}
                \seq_gput_left:Nn \g_arg_mru_this {Ctapr}
                \seq_gremove_duplicates:N \g_arg_mru_this
                Combining~
                \enum[\lref{\g_label_tl}]{
                  this;
                  #1
                } \proves[\l_arg_verbs_tl]~\ignorespaces #3
              }
            } {}
          }
        } {
          \str_case:VnF {\l_tmpa_tl} {
            {Ahpr} {
              \bool_if:NT \g_debug_bool {C2.1}
              \seq_gput_left:Nn \g_arg_mru_this {Ahpr}
              \seq_gremove_duplicates:N \g_arg_mru_this
              \enum:nxnn {#1} {\lref{\g_label_tl}} {-} {\BooleanTrue}
              \hence~
              \prove[\l_arg_verbs_tl]~\ignorespaces #3
            }
            {Tapr} {
              \bool_if:NT \g_debug_bool {C2.2}
              \seq_gput_left:Nn \g_arg_mru_this {Tapr}
              \seq_gremove_duplicates:N \g_arg_mru_this
              \enum[\lref{\g_label_tl}]{
                This;
                #1
              }[\l_arg_verbs_tl]\ignorespaces #3
            }
            {Ctapr} {
              \int_case:nn {\int_rand:nn {0} {1}} {
                {0} {
                  \bool_if:NT \g_debug_bool {C2.3}
                  \seq_gput_left:Nn \g_arg_mru_this {Ctapr}
                  \seq_gremove_duplicates:N \g_arg_mru_this
                  Combining~
                  \enum[\lref{\g_label_tl}]{
                    this;
                    #1
                  } \proves[\l_arg_verbs_tl]~\ignorespaces #3
                }
                {1} {
                  \bool_if:NT \g_debug_bool {C2.4}
                  \seq_gput_left:Nn \g_arg_mru_this {Ctapr}
                  \seq_gremove_duplicates:N \g_arg_mru_this
                  Combining~
                  \enum:nxnn {#1} {\lref{\g_label_tl}} {-} {\BooleanFalse}
                  \hence~
                  \proves[\l_arg_verbs_tl]~\ignorespaces #3
                }
              }
            }
          } {}
        }
      } {
        \seq_set_eq:NN \l_tmpa_seq \g_arg_mru_this
        \seq_remove_all:Nn \l_tmpa_seq {H}
        \seq_remove_all:Nn \l_tmpa_seq {Ahpr}
        \seq_get_right:NN \l_tmpa_seq \l_tmpa_tl
        \int_case:nnF {\seq_count:N \l_arg_seq} {
          {1} {
            \str_case:VnF {\l_tmpa_tl} {
              {Tapr} {
                \bool_if:NT \g_debug_bool {C3.1}
                \seq_gput_left:Nn \g_arg_mru_this {Tapr}
                \seq_gremove_duplicates:N \g_arg_mru_this
                \enum[\lref{\g_label_tl}]{
                  This;
                  #1
                }[\l_arg_verbs_tl]\ignorespaces #3
              }
              {Ctapr} {
                \bool_if:NT \g_debug_bool {C3.2}
                \seq_gput_left:Nn \g_arg_mru_this {Ctapr}
                \seq_gremove_duplicates:N \g_arg_mru_this
                Combining~
                \enum[\lref{\g_label_tl}]{
                  this;
                  #1
                } \proves[\l_arg_verbs_tl]~\ignorespaces #3
              }
            } {}
          }
        } {
          \str_case:VnF {\l_tmpa_tl} {
            {Tapr} {
              \bool_if:NT \g_debug_bool {C4.1}
              \seq_gput_left:Nn \g_arg_mru_this {Tapr}
              \seq_gremove_duplicates:N \g_arg_mru_this
              \enum[\lref{\g_label_tl}]{
                This;
                #1
              }[\l_arg_verbs_tl]\ignorespaces #3
            }
            {Ctapr} {
              \int_case:nn {\int_rand:nn {0} {1}} {
                {0} {
                  \bool_if:NT \g_debug_bool {C4.2}
                  \seq_gput_left:Nn \g_arg_mru_this {Ctapr}
                  \seq_gremove_duplicates:N \g_arg_mru_this
                  Combining~
                  \enum[\lref{\g_label_tl}]{
                    this;
                    #1
                  } \proves[\l_arg_verbs_tl]~\ignorespaces #3
                }
                {1} {
                  \bool_if:NT \g_debug_bool {C4.3}
                  \seq_gput_left:Nn \g_arg_mru_this {Ctapr}
                  \seq_gremove_duplicates:N \g_arg_mru_this
                  Combining~
                  \enum:nxnn {#1} {\lref{\g_label_tl}} {-} {\BooleanFalse}
                  \hence~
                  \proves[\l_arg_verbs_tl]~\ignorespaces #3
                }
              }
            }
          } {}
        }
      }
    } {
      \tl_if_head_eq_catcode:oNTF \l_cons_tl a {
        \seq_set_eq:NN \l_tmpa_seq \g_arg_mru_this
        \seq_remove_all:Nn \l_tmpa_seq {Ctapr}
        \seq_remove_all:Nn \l_tmpa_seq {Ahpr}
        \seq_get_right:NN \l_tmpa_seq \l_tmpa_tl
        \str_case:VnF {\l_tmpa_tl} {
          {H} {
            \bool_if:NT \g_debug_bool {C5.1}
            \seq_gput_left:Nn \g_arg_mru_this {H}
            \seq_gremove_duplicates:N \g_arg_mru_this
            Hence,~we~obtain~\ignorespaces #3
          }
          {Tapr} {
            \bool_if:NT \g_debug_bool {C5.2}
            \seq_gput_left:Nn \g_arg_mru_this {Tapr}
            \seq_gremove_duplicates:N \g_arg_mru_this
            This~\proves[\l_arg_verbs_tl]~\ignorespaces #3
          }
        } {}
      } {
        \bool_if:NT \g_debug_bool {C6.1}
        \seq_gput_left:Nn \g_arg_mru_this {Tapr}
        \seq_gremove_duplicates:N \g_arg_mru_this
        This~\proves[\l_arg_verbs_tl]~\ignorespaces #3
      }
    }
  } {
    \int_compare:nNnTF {\seq_count:N \l_arg_seq} = {0} {
      \bool_if:NTF \g_arg_start_bool {
        \bool_if:NT \g_debug_bool {C7.1}
        \Nobs\unskip
        #3
      } {
        \bool_if:NT \g_debug_bool {C7.2}
        \Moreover~
        #3
      }
    } {
      \bool_if:NTF \g_arg_start_bool {
        \bool_if:NT \g_debug_bool {C8.1}
        \bool_if:NT \g_kuckuck_debug_logs_bool { \tl_log:N \l_arg_verbs_tl }
        \Nobs~that~
        \enum{
          #1
        }[\l_arg_verbs_tl]\ignorespaces #3
      } {
        \int_compare:nNnTF {\seq_count:N \l_arg_seq} = {1} {
          \seq_set_eq:NN \l_tmpa_seq \g_arg_mru_nothis
          \seq_remove_all:Nn \l_tmpa_seq {Nwc}
          \seq_remove_all:Nn \l_tmpa_seq {Itnswc}
          \seq_get_right:NN \l_tmpa_seq \l_tmpa_tl
        } {
          \seq_get_right:NN \g_arg_mru_nothis \l_tmpa_tl
        }
        \str_case:VnF {\l_tmpa_tl} {
          {Mo} {
            \bool_if:NT \g_debug_bool {C9.1}
            \seq_gput_left:Nn \g_arg_mru_nothis {Mo}
            \seq_gremove_duplicates:N \g_arg_mru_nothis
            Moreover,~\nobs~that~
            \enum{
              #1
            }[\l_arg_verbs_tl]\ignorespaces #3
          }
          {Fm} {
            \bool_if:NT \g_debug_bool {C9.2}
            \seq_gput_left:Nn \g_arg_mru_nothis {Fm}
            \seq_gremove_duplicates:N \g_arg_mru_nothis
            Furthermore,~\nobs~that~
            \enum{
              #1
            }[\l_arg_verbs_tl]\ignorespaces #3
          }
          {Ia} {
            \bool_if:NT \g_debug_bool {C9.3}
            \seq_gput_left:Nn \g_arg_mru_nothis {Ia}
            \seq_gremove_duplicates:N \g_arg_mru_nothis
            In~addition,~\nobs~that~
            \enum{
              #1
            }[\l_arg_verbs_tl]\ignorespaces #3
          }
          {N} {
            \bool_if:NT \g_debug_bool {C9.4}
            \seq_gput_left:Nn \g_arg_mru_nothis {N}
            \seq_gremove_duplicates:N \g_arg_mru_nothis
            Next,~\nobs~that~
            \enum{
              #1
            }[\l_arg_verbs_tl]\ignorespaces #3
          }
          {Itns} {
            \bool_if:NT \g_debug_bool {C9.5}
            \seq_gput_left:Nn \g_arg_mru_nothis {Itnswc}
            \seq_gput_left:Nn \g_arg_mru_nothis {Itns}
            \seq_gremove_duplicates:N \g_arg_mru_nothis
            In~the~next~step~we~\nobs~that~
            \enum{
              #1
            }[\l_arg_verbs_tl]\ignorespaces #3
          }
          {Nwc} {
            \bool_if:NT \g_debug_bool {C9.6}
            \seq_gput_left:Nn \g_arg_mru_nothis {Nwc}
            \seq_gremove_duplicates:N \g_arg_mru_nothis
            Next~we~combine~
            \enum{
              #1
            }to~obtain~\ignorespaces #3
          }
          {Itnswc} {
            \bool_if:NT \g_debug_bool {C9.7}
            \seq_gput_left:Nn \g_arg_mru_nothis {Itns}
            \seq_gput_left:Nn \g_arg_mru_nothis {Itnswc}
            \seq_gremove_duplicates:N \g_arg_mru_nothis
            In~the~next~step~we~combine~
            \enum{
              #1
            }to~obtain~\ignorespaces #3
          }
        } {}
      }
    }
  }
  \bool_gset_false:N \g_arg_start_bool
  \bool_gset_false:N \l_insidearg_bool
  \cfload[.]%
  \bool_if:NT \g_kuckuck_debug_colors_bool { \color{black} }
}

\tl_new:N \g_label_tl
\tl_gset:Nn \g_label_tl { }

\NewDocumentCommand{\savelabel}{m}{
  \bool_if:NTF \l_insidearg_bool {
    \tl_gset:Nn \g_label_tl {#1}
  } {
    \tl_gset:Nn \g_label_tl { }
  }
}

\ExplSyntaxOff

\ExplSyntaxOn

\NewDocumentEnvironment {athm} {m m o} {
\str_if_eq:noTF {example} {#1} {
  \bool_gset_true:N \g_example_bool
} {
  \bool_gset_false:N \g_example_bool
}
\cfclear
\IfNoValueTF{#3}{
\begin{#1}\label{#2}\global\def\loc{#2}
}{
\begin{#1}[#3]\label{#2}\global\def\loc{#2}
}
}{
\end{#1}
}

\NewDocumentEnvironment {adef} {m} {
\begin{definition}\label{#1}\global\def\loc{#1}
}{
\end{definition}
}

\NewDocumentEnvironment{aproof} {} {
\bool_if:NTF \g_example_bool {
  \bool_gset_true:N \g_arg_start_bool
  \begin{proof}[Proof~for~\cref{\loc}]
} {
  \bool_gset_true:N \g_arg_start_bool
  \begin{proof}[Proof~of~\cref{\loc}]
}
\bool_gset_false:N \g_finishproof_bool
}{
\bool_if:NTF \g_finishproof_bool {}
{\finishproofthus}
\end{proof}
}

\NewDocumentCommand{\finishproofthus} {} {
  \bool_gset_true:N \g_finishproof_bool
  \bool_if:NTF \g_example_bool {
    The~proof~for~\cref{\loc}~is~thus~complete.
  } {
    The~proof~of~\cref{\loc}~is~thus~complete.
  }
}
\NewDocumentCommand{\finishproofthis} {} {
  \bool_gset_true:N \g_finishproof_bool
  \bool_if:NTF \g_example_bool {
    This~completes~the~proof~for~\cref{\loc}.
  } {
    This~completes~the~proof~of~\cref{\loc}.
  }
}

\ExplSyntaxOff

\usepackage{mleftright}

\ExplSyntaxOn

\cs_generate_variant:Nn \seq_set_split:Nnn { NnV }

\tl_new:N \l_content_tl
\tl_new:N \l_last_tl

\tl_new:N \l_eqn_type

\NewDocumentEnvironment{eqn}{ob}{
  \tl_set:Nn \l_tmpa_tl {#2}
  \tl_trim_spaces:N \l_tmpa_tl
  \tl_set:Nx \l_last_tl  { \tl_range:Nnn { \l_tmpa_tl } { -1 } { -1 } }
  \tl_if_in:nVTF {.,} \l_last_tl {
    \tl_set:Nx \l_content_tl  { \tl_range:Nnn { \l_tmpa_tl } { 1 }  { -2 } }
  } {
    \tl_set:NV \l_content_tl  { \l_tmpa_tl }
    \tl_set:Nn \l_last_tl {}
  }
  \tl_set:Nn \l_eqn_type {d}
  \keys_define:nn { benno/argument } {
    inline .value_forbidden:n = true,
    inline .code:n = {\tl_set:Nn \l_eqn_type {i}},
    i .value_forbidden:n = true,
    i .code:n = {\tl_set:Nn \l_eqn_type {i}},
    align .value_forbidden:n = true,
    align .code:n = {\tl_set:Nn \l_eqn_type {a}},
    a .value_forbidden:n = true,
    a .code:n = {\tl_set:Nn \l_eqn_type {a}},
    display .value_forbidden:n = true,
    display .code:n = {\tl_set:Nn \l_eqn_type {d}},
    d .value_forbidden:n = true,
    d .code:n = {\tl_set:Nn \l_eqn_type {d}},
    oneline .value_forbidden:n = true,
    oneline .code:n = {\tl_remove_all:Nn \l_content_tl {&}
                        \tl_remove_all:Nn \l_content_tl {\\}},
    autosize .value_forbidden:n = true,
    autosize .code:n = {
                        \tl_replace_all:Nnn \l_content_tl {\lvert}{\mleft\lvert}
                        \tl_replace_all:Nnn \l_content_tl {\rvert}{\mright\rvert}
                        \tl_replace_all:Nnn \l_content_tl {\lVert}{\mleft\lVert}
                        \tl_replace_all:Nnn \l_content_tl {\rVert}{\mright\rVert}
                        \tl_replace_all:Nnn \l_content_tl {\langle}{\mleft\langle}
                        \tl_replace_all:Nnn \l_content_tl {\rangle}{\mright\rangle}
                        \tl_replace_all:Nnn \l_content_tl {(|}{\mleft\lvert}
                        \tl_replace_all:Nnn \l_content_tl {|)}{\mright\rvert}
                        \tl_replace_all:Nnn \l_content_tl {(||}{\mleft\lVert}
                        \tl_replace_all:Nnn \l_content_tl {||)}{\mright\rVert}
                        \tl_replace_all:Nnn \l_content_tl {(<}{\mleft\langle}
                        \tl_replace_all:Nnn \l_content_tl {>)}{\mright\rangle}
                        \tl_replace_all:Nnn \l_content_tl {(}{\mleft(}
                        \tl_replace_all:Nnn \l_content_tl {)}{\mright)}
                        \tl_replace_all:Nnn \l_content_tl {[}{\mleft[}
                        \tl_replace_all:Nnn \l_content_tl {]}{\mright]}
                        \tl_replace_all:Nnn \l_content_tl {\{}{\mleft\{}
                        \tl_replace_all:Nnn \l_content_tl {\}}{\mright\}}
                        \tl_replace_all:Nnn \l_content_tl {\lceil}{\mleft\lceil}
                        \tl_replace_all:Nnn \l_content_tl {\rceil}{\mright\rceil}
                        \tl_replace_all:Nnn \l_content_tl {\lfloor}{\mleft\lfloor}
                        \tl_replace_all:Nnn \l_content_tl {\rfloor}{\mright\rfloor}
                        },
  }
  \IfValueT{#1}{
    \keys_set:nn { benno/argument } {#1}
  }
  \bool_if:NT \g_kuckuck_debug_logs_bool { \tl_log:N \l_content_tl }
  \str_case:Vn \l_eqn_type {
    {i} {
      \math
      \l_content_tl
    }
    {a} {
      \equation
      \aligned
      \l_content_tl
    }
    {d} {
      \equation
      \l_content_tl
    }
  }
} {
  \str_case:Vn \l_eqn_type {
    {i} {
      \tl_if_eq:NnTF \l_last_tl {.} {
        \endmath \ifnocf{.}
      } {
        \l_last_tl\endmath
      }
    }
    {a} {
      \tl_if_eq:NnTF \l_last_tl {.}
        {\ifnocf{.}}
        {\l_last_tl}
      \endaligned
      \endequation
    }
    {d} {
      \tl_if_eq:NnTF \l_last_tl {.} {
        \ifnocf{.}\endequation
      } {
        \l_last_tl\endequation
      }
    }
  }
}

\ExplSyntaxOff

\makeatletter
\ExplSyntaxOn
\seq_new:N \g_abbrs
\prop_new:N \g_abbr_counts
\tl_new:N \l_abbr_count_tl

\ExplSyntaxOn

\bool_new:N \g_forexample

\NewDocumentCommand{\eg}{ o }{
	\IfValueT{#1}{
		\str_if_eq:noTF {fe} {#1} {
			\bool_gset_true:N \g_forexample
		} {\bool_gset_false:N \g_forexample}
	}
	\bool_if:nTF { \g_forexample } {
		\bool_gset_false:N \g_forexample
		for~example
	}{
		\bool_gset_true:N \g_forexample
		for~instance
	}
}

\NewDocumentCommand{\abbr}{m m O{#1} m m O{#4} m}{
	\expandafter\newcommand\csname#3\endcsname[1][]{
		\seq_if_in:NnTF \g_abbrs {#1} {
			\prop_get:NnN \g_abbr_counts {#1} \l_abbr_count_tl
			\prop_gput:Nnx \g_abbr_counts {#1} {\int_eval:n {\l_abbr_count_tl + 1}}
			\hyperref[#1]{#7}
		} {
			\seq_gput_left:Nn \g_abbrs {#1}
			\prop_gput:Nnn \g_abbr_counts {#1} {1}
			\expandafter\gdef\csname#1@def\endcsname{#2}
			\phantomsection\label{#1}
			\str_if_eq:nnTF{##1}{}{\emph{#2}}{##1}~(\hyperref[#1]{#7})
		}
	}
	\expandafter\newcommand\csname#6\endcsname[1][]{
		\seq_if_in:NnTF \g_abbrs {#1} {
			\prop_get:NnN \g_abbr_counts {#1} \l_abbr_count_tl
			\prop_gput:Nnx \g_abbr_counts {#1} {\int_eval:n {\l_abbr_count_tl + 1}}
			\hyperref[#1]{#4}
		} {
			\expandafter\gdef\csname#1@def\endcsname{#5}
			\seq_gput_left:Nn \g_abbrs {#1}
			\prop_gput:Nnn \g_abbr_counts {#1} {1}
			\phantomsection\label{#1}
			\str_if_eq:nnTF{##1}{}{\emph{#5}}{##1}~(\hyperref[#1]{#4})
		}
	}
}

\ExplSyntaxOff

%% file: main.bbl
\begin{thebibliography}{10}

\bibitem{Apostol74}
T.~M. Apostol.
\newblock {\em Mathematical analysis}.
\newblock Addison-Wesley Publishing Co., Reading, Mass.-London-Don Mills, Ont.,
  second edition, 1974.

\bibitem{MR713046}
V.~I. Arnold.
\newblock Some open problems in the theory of singularities.
\newblock In {\em Singularities, {P}art 1 ({A}rcata, {C}alif., 1981)},
  volume~40 of {\em Proc. Sympos. Pure Math.}, 57--69. Amer. Math. Soc.,
  Providence, RI, 1983.
\newblock Translated from the Russian.

\bibitem{MR972342}
E.~Bierstone and P.~D. Milman.
\newblock Semianalytic and subanalytic sets.
\newblock {\em Inst. Hautes Etudes Sci. Publ. Math.} 67, 5--42, 1988.
\newblock \url{http://www.numdam.org/item?id=PMIHES_1988__67__5_0}.

\bibitem{Billingsley95}
P.~Billingsley.
\newblock {\em Probability and measure}.
\newblock Wiley Series in Probability and Mathematical Statistics. John Wiley
  \& Sons, Inc., New York, third edition, 1995.
\newblock A Wiley-Interscience Publication.

\bibitem{MR1659509}
J.~Bochnak, M.~Coste, and M.-F. Roy.
\newblock {\em Real algebraic geometry}, volume~36 of {\em Ergebnisse der
  Mathematik und ihrer Grenzgebiete (3) [Results in Mathematics and Related
  Areas (3)]}.
\newblock Springer-Verlag, Berlin, 1998.
\newblock Translated from the 1987 French original, Revised by the authors.
\newblock
  \href{https://doi.org/10.1007/978-3-662-03718-8}{\nolinkurl{doi:10.1007/978-3-662-03718-8}}.

\bibitem{MR2274510}
J.~Bolte, A.~Daniilidis, and A.~Lewis.
\newblock The {{\L}ojasiewicz} inequality for nonsmooth subanalytic functions
  with applications to subgradient dynamical systems.
\newblock {\em SIAM J. Optim.} 17(4), 1205--1223, 2007.
\newblock
  \href{https://doi.org/10.1137/050644641}{\nolinkurl{doi:10.1137/050644641}}.

\bibitem{Boumal23}
N.~Boumal.
\newblock {\em An introduction to optimization on smooth manifolds}.
\newblock Cambridge University Press, Cambridge, 2023.
\newblock
  \href{https://doi.org/10.1017/9781009166164}{\nolinkurl{doi:10.1017/9781009166164}}.

\bibitem{BBI01}
D.~Burago, Y.~Burago, and S.~Ivanov.
\newblock {\em A Course in Metric Geometry}, volume~33 of {\em Graduate Studies
  in Mathematics}.
\newblock American Mathematical Society, Providence, RI, 2001.
\newblock
  \href{https://doi.org/10.1090/gsm/033}{\nolinkurl{doi:10.1090/gsm/033}}.

\bibitem{Cai25}
Y.~Cai, K.~Zhou, J.~Wu, S.~Mei, M.~Lindsey, and P.~L. Bartlett.
\newblock Implicit bias of gradient descent for non-homogeneous deep networks.
\newblock In {\em Proceedings of the 42nd International Conference on Machine
  Learning}, volume 267 of {\em Proceedings of Machine Learning Research},
  6426--6504. PMLR, 2025.
\newblock \url{https://proceedings.mlr.press/v267/cai25m.html}.

\bibitem{ChizatBach20}
L.~Chizat and F.~Bach.
\newblock Implicit bias of gradient descent for wide two-layer neural networks
  trained with the logistic loss.
\newblock In {\em Proceedings of the 33rd Conference on Learning Theory},
  volume 125 of {\em Proceedings of Machine Learning Research}, 1305--1338.
  PMLR, 2020.
\newblock \url{https://proceedings.mlr.press/v125/chizat20a.html}.

\bibitem{Coste}
M.~Coste.
\newblock An introduction to {O}-minimal geometry.
\newblock Lecture notes, November 1999.
\newblock \url{https://univ-rennes.hal.science/hal-05413940/document}.

\bibitem{MR2166522}
D.~D'Acunto and V.~Grandjean.
\newblock A gradient inequality at infinity for tame functions.
\newblock {\em Rev. Mat. Complut.} 18(2), 493--501, 2005.
\newblock
  \href{https://doi.org/10.5209/rev_rema.2005.v18.n2.16699}{\nolinkurl{doi:10.5209/rev_rema.2005.v18.n2.16699}}.

\bibitem{MR2208534}
D.~D'Acunto and V.~Grandjean.
\newblock On gradient at infinity of semialgebraic functions.
\newblock {\em Ann. Polon. Math.} 87, 39--49, 2005.
\newblock
  \href{https://doi.org/10.4064/ap87-0-4}{\nolinkurl{doi:10.4064/ap87-0-4}}.

\bibitem{DHL22_MR4453693}
A.~Daniilidis, M.~Haddou, and O.~Ley.
\newblock A convex function satisfying the {{\L}ojasiewicz} inequality but
  failing the gradient conjecture both at zero and infinity.
\newblock {\em Bull. Lond. Math. Soc.} 54(2), 590--608, 2022.
\newblock
  \href{https://doi.org/10.1112/blms.12586}{\nolinkurl{doi:10.1112/blms.12586}}.

\bibitem{DatchevShapiro23}
K.~Datchev and J.~Shapiro.
\newblock Exponential time-decay for a one-dimensional wave equation with
  coefficients of bounded variation.
\newblock {\em Math. Nachr.} 296(11), 4978--4994, 2023.
\newblock
  \href{https://doi.org/10.1002/mana.202200459}{\nolinkurl{doi:10.1002/mana.202200459}}.

\bibitem{DJR25}
S.~Dereich, A.~Jentzen, and A.~Riekert.
\newblock Averaged {Adam} accelerates stochastic optimization in the training
  of deep neural network approximations for partial differential equation and
  optimal control problems.
\newblock arXiv:2501.06081, 2025.
\newblock \url{https://arxiv.org/abs/2501.06081}.

\bibitem{DoJR25}
T.~Do, A.~Jentzen, and A.~Riekert.
\newblock Non-convergence to the optimal risk for {Adam} and stochastic
  gradient descent optimization in the training of deep neural networks.
\newblock arXiv:2503.01660, 2025.
\newblock \url{https://arxiv.org/abs/2503.01660}.

\bibitem{Fischer2010Peano}
A.~Fischer.
\newblock Peano differentiable extensions in o-minimal structures.
\newblock {\em Adv. Geom.} 10(1), 31--49, 2010.
\newblock
  \href{https://doi.org/10.1515/ADVGEOM.2009.026}{\nolinkurl{doi:10.1515/ADVGEOM.2009.026}}.

\bibitem{FS01_MR1844525}
P.~Fortuny and F.~Sanz.
\newblock Gradient vector fields do not generate twister dynamics.
\newblock {\em J. Differential Equations} 174(1), 91--100, 2001.
\newblock
  \href{https://doi.org/10.1006/jdeq.2000.3926}{\nolinkurl{doi:10.1006/jdeq.2000.3926}}.

\bibitem{Gabrielov1968}
A.~M. Gabrielov.
\newblock Projections of semi-analytic sets.
\newblock {\em Funct. Anal. Appl.} 2(4), 282--291, 1968.
\newblock
  \href{https://doi.org/10.1007/BF01075680}{\nolinkurl{doi:10.1007/BF01075680}}.

\bibitem{GJL22}
D.~Gallon, A.~Jentzen, and F.~Lindner.
\newblock Blow up phenomena for gradient descent optimization methods in the
  training of artificial neural networks.
\newblock arXiv:2211.15641, 2022.
\newblock \url{https://arxiv.org/abs/2211.15641}.

\bibitem{Goldstein2009_MR2568063}
P.~Goldstein.
\newblock Gradient flow of a harmonic function in {$\mathbb R^3$}.
\newblock {\em J. Differential Equations} 247(9), 2517--2557, 2009.
\newblock
  \href{https://doi.org/10.1016/j.jde.2009.07.027}{\nolinkurl{doi:10.1016/j.jde.2009.07.027}}.

\bibitem{G07_MR2290270}
V.~Grandjean.
\newblock On the limit set at infinity of a gradient trajectory of a
  semialgebraic function.
\newblock {\em J. Differential Equations} 233(1), 22--41, 2007.
\newblock
  \href{https://doi.org/10.1016/j.jde.2006.10.009}{\nolinkurl{doi:10.1016/j.jde.2006.10.009}}.

\bibitem{G14MR3200186}
V.~Grandjean.
\newblock Gradient trajectories for plane singular metrics {I}: {O}scillating
  trajectories.
\newblock {\em Demonstr. Math.} 47(1), 69--78, 2014.
\newblock
  \href{https://doi.org/10.2478/dema-2014-0006}{\nolinkurl{doi:10.2478/dema-2014-0006}}.

\bibitem{GS13_MR3072668}
V.~Grandjean and F.~Sanz.
\newblock On restricted analytic gradients on analytic isolated surface
  singularities.
\newblock {\em J. Differential Equations} 255(7), 1684--1708, 2013.
\newblock
  \href{https://doi.org/10.1016/j.jde.2013.05.020}{\nolinkurl{doi:10.1016/j.jde.2013.05.020}}.

\bibitem{Ichikawa92_MR1156895}
F.~Ichikawa.
\newblock Thom's conjecture on singularities of gradient vector fields.
\newblock {\em Kodai Math. J.} 15(1), 134--140, 1992.
\newblock
  \href{https://doi.org/10.2996/kmj/1138039532}{\nolinkurl{doi:10.2996/kmj/1138039532}}.

\bibitem{JKR25}
A.~Jentzen, J.~Kranz, and A.~Riekert.
\newblock {PADAM}: Parallel averaged {Adam} reduces the error for stochastic
  optimization in scientific machine learning.
\newblock arXiv:2505.22085, 2025.
\newblock Accepted for publication in \emph{J. Comput. Math.}.
  \url{https://arxiv.org/abs/2505.22085}.

\bibitem{JKW25Book}
A.~Jentzen, B.~Kuckuck, and P.~von Wurstemberger.
\newblock Mathematical introduction to deep learning: Methods, implementations,
  and theory.
\newblock arXiv:2310.20360v3, 2025.
\newblock \url{https://arxiv.org/abs/2310.20360}.

\bibitem{JT20}
Z.~Ji and M.~Telgarsky.
\newblock Directional convergence and alignment in deep learning.
\newblock In {\em Advances in Neural Information Processing Systems},
  volume~33, 17176--17186, 2020.
\newblock
  \url{https://proceedings.neurips.cc/paper_files/paper/2020/hash/c76e4b2fa54f8506719a5c0dc14c2eb9-Abstract.html}.

\bibitem{MR3569073}
G.~Jones, J.~Kirby, and T.~Servi.
\newblock Local interdefinability of {W}eierstrass elliptic functions.
\newblock {\em J. Inst. Math. Jussieu} 15(4), 673--691, 2016.
\newblock
  \href{https://doi.org/10.1017/S1474748014000425}{\nolinkurl{doi:10.1017/S1474748014000425}}.

\bibitem{JordanEtAl24Muon}
K.~Jordan, Y.~Jin, V.~Boza, J.~You, F.~Cesista, L.~Newhouse, and J.~Bernstein.
\newblock {Muon}: An optimizer for hidden layers in neural networks.
\newblock Blog post, December 2024.
\newblock \url{https://kellerjordan.github.io/posts/muon/}.

\bibitem{KingmaBa15}
D.~P. Kingma and J.~Ba.
\newblock {Adam}: A method for stochastic optimization.
\newblock In {\em Proceedings of the 3rd International Conference on Learning
  Representations}, 2015.
\newblock \url{https://arxiv.org/abs/1412.6980}.

\bibitem{KGDJ25}
J.~Kranz, D.~Gallon, S.~Dereich, and A.~Jentzen.
\newblock {SAD} neural networks: Divergent gradient flows and asymptotic
  optimality via {o}-minimal structures.
\newblock In {\em Advances in Neural Information Processing Systems},
  volume~38, 4428--4464, 2025.
\newblock
  \href{https://doi.org/10.52202/085713-0138}{\nolinkurl{doi:10.52202/085713-0138}}.

\bibitem{Kuiper1991}
N.~H. Kuiper.
\newblock On gradient curves of an analytic function near a critical point.
\newblock IHES Preprint M/91/36, 1991.
\newblock
  \url{https://repo-archives.ihes.fr/FONDS_IHES/I_Prepublications/KUIPER/1972-1991/M_91_36/M_91_36.pdf}.

\bibitem{K98}
K.~Kurdyka.
\newblock On gradients of functions definable in o-minimal structures.
\newblock {\em Ann. Inst. Fourier (Grenoble)} 48(3), 769--783, 1998.
\newblock
  \href{https://doi.org/10.5802/aif.1638}{\nolinkurl{doi:10.5802/aif.1638}}.

\bibitem{KMP00}
K.~Kurdyka, T.~Mostowski, and A.~Parusi\'nski.
\newblock Proof of the gradient conjecture of {R}. {T}hom.
\newblock {\em Ann. of Math. (2)} 152(3), 763--792, 2000.
\newblock
  \href{https://doi.org/10.2307/2661354}{\nolinkurl{doi:10.2307/2661354}}.

\bibitem{KP06_MR2325137}
K.~Kurdyka and A.~Parusi\'nski.
\newblock Quasi-convex decomposition in o-minimal structures. {A}pplication to
  the gradient conjecture.
\newblock In {\em Singularity theory and its applications}, volume~43 of {\em
  Adv. Stud. Pure Math.}, 137--177. Math. Soc. Japan, Tokyo, 2006.
\newblock
  \href{https://doi.org/10.2969/aspm/04310137}{\nolinkurl{doi:10.2969/aspm/04310137}}.

\bibitem{Lee13}
J.~M. Lee.
\newblock {\em Introduction to smooth manifolds}, volume 218 of {\em Graduate
  Texts in Mathematics}.
\newblock Springer, New York, second edition, 2012.
\newblock
  \href{https://doi.org/10.1007/978-1-4419-9982-5}{\nolinkurl{doi:10.1007/978-1-4419-9982-5}}.

\bibitem{Loi_MR3483169}
T.~L. Loi.
\newblock {\L}ojasiewicz inequalities in o-minimal structures.
\newblock {\em Manuscripta Math.} 150(1--2), 59--72, 2016.
\newblock
  \href{https://doi.org/10.1007/s00229-015-0806-y}{\nolinkurl{doi:10.1007/s00229-015-0806-y}}.

\bibitem{Loj84_MR771152}
S.~{\L}ojasiewicz.
\newblock Sur les trajectoires du gradient d'une fonction analytique.
\newblock In {\em Geometry seminars, 1982--1983 ({B}ologna, 1982/1983)},
  115--117. Univ. Stud. Bologna, Bologna, 1984.

\bibitem{MR1275210}
S.~{\L}ojasiewicz.
\newblock Sur la g\'eom\'etrie semi- et sous-analytique.
\newblock {\em Ann. Inst. Fourier (Grenoble)} 43(5), 1575--1595, 1993.
\newblock
  \href{https://doi.org/10.5802/aif.1384}{\nolinkurl{doi:10.5802/aif.1384}}.

\bibitem{LyuLi20}
K.~Lyu and J.~Li.
\newblock Gradient descent maximizes the margin of homogeneous neural networks.
\newblock In {\em Proceedings of the 8th International Conference on Learning
  Representations}, 2020.
\newblock \url{https://arxiv.org/abs/1906.05890}.

\bibitem{Moussu1997_MR1437048}
R.~Moussu.
\newblock Sur la dynamique des gradients. {E}xistence de vari\'et\'es
  invariantes.
\newblock {\em Math. Ann.} 307(3), 445--460, 1997.
\newblock
  \href{https://doi.org/10.1007/s002080050043}{\nolinkurl{doi:10.1007/s002080050043}}.

\bibitem{P05}
A.~Parusi\'nski.
\newblock On the gradient conjecture for definable functions.
\newblock In M.~J. Edmundo, D.~Richardson, and A.~Wilkie, editors, {\em
  O-minimal Structures, Lisbon 2003}, 143--158. Cuvillier Verlag, G\"ottingen,
  2005.
\newblock \url{https://cuvillier.de/get/ebook/3812/9783865375575.pdf}.

\bibitem{Rudin87}
W.~Rudin.
\newblock {\em Real and complex analysis}.
\newblock McGraw-Hill Book Co., New York, third edition, 1987.

\bibitem{Sanz1998_MR1656007}
F.~Sanz.
\newblock Non-oscillating solutions of analytic gradient vector fields.
\newblock {\em Ann. Inst. Fourier (Grenoble)} 48(4), 1045--1067, 1998.
\newblock
  \href{https://doi.org/10.5802/aif.1648}{\nolinkurl{doi:10.5802/aif.1648}}.

\bibitem{Sanz2024_MR4789467}
F.~Sanz~S\'anchez.
\newblock Half a century with the problem of the gradient of an analytic
  function.
\newblock In {\em Handbook of geometry and topology of singularities {V}.
  {F}oliations}, 335--395. Springer, Cham, 2024.
\newblock
  \href{https://doi.org/10.1007/978-3-031-52481-3_7}{\nolinkurl{doi:10.1007/978-3-031-52481-3_7}}.

\bibitem{Soudry18}
D.~Soudry, E.~Hoffer, M.~S. Nacson, S.~Gunasekar, and N.~Srebro.
\newblock The implicit bias of gradient descent on separable data.
\newblock {\em J. Mach. Learn. Res.} 19, Paper No. 70, 57 pages, 2018.
\newblock \url{http://jmlr.org/papers/v19/18-188.html}.

\bibitem{Thom89}
R.~Thom.
\newblock Probl\`emes rencontr\'es dans mon parcours math\'ematique: un bilan.
\newblock {\em Inst. Hautes \'Etudes Sci. Publ. Math.} 70, 199--214, 1989.
\newblock \url{http://www.numdam.org/item?id=PMIHES_1989__70__199_0}.

\bibitem{vanDenDries1986}
L.~van~den Dries.
\newblock A generalization of the {Tarski--Seidenberg} theorem, and some
  nondefinability results.
\newblock {\em Bull. Amer. Math. Soc. (N.S.)} 15(2), 189--193, 1986.
\newblock
  \href{https://doi.org/10.1090/S0273-0979-1986-15468-6}{\nolinkurl{doi:10.1090/S0273-0979-1986-15468-6}}.

\bibitem{vdD_MR1633348}
L.~van~den Dries.
\newblock {\em Tame topology and o-minimal structures}, volume 248 of {\em
  London Mathematical Society Lecture Note Series}.
\newblock Cambridge University Press, Cambridge, 1998.
\newblock
  \href{https://doi.org/10.1017/CBO9780511525919}{\nolinkurl{doi:10.1017/CBO9780511525919}}.

\bibitem{Wilkie2009_MR2605320}
A.~J. Wilkie.
\newblock {O}-minimal structures.
\newblock In {\em S\'eminaire Bourbaki. Volume 2007/2008. Expos\'es 982--996},
  number 326 in Ast\'erisque, 131--142. Soci\'et\'e math\'ematique de France,
  2009.
\newblock Talk no. 985.
\newblock \url{https://www.numdam.org/item/AST_2009__326__131_0/}.

\end{thebibliography}
